\documentclass[10pt]{article}
\usepackage{mypackages}
\usepackage{makeidx}

\usepackage[T1]{fontenc}
\usepackage{titlesec}
\title{\textbf{\textsc{A wavelet-inspired $L^3$-based convex integration framework for the Euler equations}}}

\author{{\textsc{vikram giri, hyunju kwon, and matthew novack}}}

\date{}
\makeindex

\allowdisplaybreaks

\titleformat{\subsection}[runin]
      {\normalfont\bfseries}
      {\thesubsection}
      {0.5em}
      {}
      [.]
\titleformat{\subsubsection}[runin]
      {\normalfont\bfseries}
      {\thesubsubsection}
      {0.5em}
      {}
      [.]

\begin{document}

\maketitle

\begin{abstract}
In this work, we develop a wavelet-inspired, $L^3$-based convex integration framework for constructing weak solutions to the three-dimensional incompressible Euler equations.  The main innovations include a new multi-scale building block, which we call an intermittent Mikado bundle; a wavelet-inspired inductive set-up which includes assumptions on spatial and temporal support, in addition to $L^p$ and pointwise estimates for Eulerian and Lagrangian derivatives; and sharp decoupling lemmas, inverse divergence estimates, and space-frequency localization technology which is well-adapted to functions satisfying $L^p$ estimates for $p$ other than $1$, $2$, or $\infty$. We develop these tools in the context of the Euler-Reynolds system, enabling us to give both a new proof of the intermittent Onsager theorem from~\cite{NV22} in this paper, and a proof of the $L^3$-based strong Onsager conjecture in the companion paper~\cite{GKN23}.
\end{abstract}

\setcounter{tocdepth}{1}
\tableofcontents

\section{Introduction}

\subsection{The \texorpdfstring{$L^3$-based}{llllll} strong Onsager conjecture}

We consider the three-dimensional incompressible Euler equations on $[0,T]\times\mathbb{T}^3$, which are given by 
\begin{equation}\label{eqn:Euler}
    \begin{cases}
    \pa_t u + (u\cdot \na) u + \na p = 0 \\
    \div \, u =0 \, .
    \end{cases}
\end{equation}
Smooth solutions of these equations satisfy a pointwise energy balance obtained by taking the dot product of the first equation in \eqref{eqn:Euler} with $u$.  Integration of this balance in time and space then implies that smooth solutions conserve the total kinetic energy $\sfrac 12 \| u(t,\cdot) \|_{L^2(\T^3)}^2$.  However, there is significant mathematical and physical motivation behind the study of \emph{weak} solutions of \eqref{eqn:Euler} which allow for the dissipation of kinetic energy.  These dissipative weak solutions of \eqref{eqn:Euler} will satisfy the local energy identity 
\begin{align}
    \partial_t \left( \frac 12 |u|^2 \right) + \div \left( \left( \frac{1}{2} |u|^2 + p \right) u  \right) = - D[u] \label{thurzday:evening}
\end{align} 
in the sense of distributions, where the Duchon-Robert measure $D[u]$ captures the dissipation due to possible singularities \cite{DuchonRobert00}. For Euler flows arising as vanishing-viscosity limits of suitable Navier-Stokes flows, this measure is non-negative \cite{DuchonRobert00}, and the resulting inequality in \eqref{thurzday:evening} is referred to as the \emph{local energy inequality}.

The well-known Onsager conjecture \cite{Onsager49} postulates that $L^\infty_tC^{\sfrac 13}_{x}$ serves as a threshold, below which weak solutions of the Euler equations \eqref{eqn:Euler} may dissipate the total kinetic energy \cite{Onsager49}, and above which solutions must conserve the kinetic energy.  Recent years have seen remarkable success in the validation of Onsager's conjecture. The conservation of kinetic energy for solutions in $L^3_t B^{\alpha}_{3,\infty}$ for $\alpha>\sfrac 13$ has been proven by Constantin, E, and Titi in \cite{ConstantinETiti94} (see also \cite{Eyink75, CCFS08,DE,DuchonRobert00,DI,DDI}), and the flexibility statement was proven by Isett in \cite{Isett2018} and extended by Buckmaster, De Lellis, Sz\'ekelyhidi, and Vicol \cite{BDLSV17}.  The proofs in \cite{Isett2018, BDLSV17} utilize the convex integration framework initiated by De Lellis and Sz\'ekelyhidi in \cite{DLS09, DeLellisSzekelyhidi13}, inspired by Nash's work \cite{Nash} and following work of Scheffer~\cite{Scheffer93} and Shnirelman~\cite{Shnirelman00}; we refer the reader to the survey papers \cite{BVReview, DSReview} for further history of the Onsager program. 

The regularity threshold $C^{\sfrac 13}$ is also intimately connected to Kolmogorov's 1941 (K41) phenomenological theory of turbulence \cite{K2, K3, K1}, which may be interpreted as suggesting that turbulent fluids enjoy uniform $L^p_t B^{\sfrac 13}_{p,\infty,x}$ regularity in the vanishing viscosity limit for $p\in[1,\infty)$.  Here we define the inhomogeneous Besov norms for $s\in(0,1)$ and $p\in[1,\infty]$ by 
\begin{align*}
    \norm{v}_{B_{p,\infty}^s(\T^3)}
    \sim \norm{v}_{L^p(\T^3)} + \sup_{|z|>0} \frac{\norm{v(\cdot +z)-v}_{L^p(\T^3)}}{|z|^s}\, .
\end{align*}
Such uniform regularity bounds would then imply that dissipative solutions of Euler obtained as vanishing viscosity limits enjoy $L^\infty_t C^{\sfrac 13}_{x}$ regularity, or the \emph{maximum} amount of regularity identified by Onsager as allowing for the dissipation of kinetic energy. In the case $p=3$, K41 scaling is strongly supported by experimental evidence \cite[Figure 8.8]{Frisch95}, \cite[Figure5]{ChenEtAl05}, \cite[Figure 3]{IshiharaEtAl09}, \cite[Figure 1]{ISY20}, indicating that $B^{\sfrac 13}_{3,\infty}$ may indeed be a natural function space for turbulent flows.  However, it is well known that turbulent fluids exhibit deviations from the K41 scaling $S_p(\ell)\sim \sfrac p3$ when $p \neq 3$. When $p<3$, one typically observes that $\sfrac{\zeta_p}p> 1/3$, while for $p>3$, one typically observes that $\sfrac{\zeta_p}p<1/3$; see \cite[Figure 8.8]{Frisch95}, or \cite[Figure 6]{ISY20} for a recent numerical simulation. These observations suggest that the H\"older space $C^{\sfrac13}$ in which Onsager's theorem has been proven may not be the most reasonable space for turbulent flows. In this direction, the third author and Vicol recently proved an intermittent Onsager theorem \cite{NV22} for non-conservative solutions in $C^0_t(H^{\sfrac 12 -}\cap L^{\infty-})\subset C^0_tB_{3,\infty}^{\sfrac 13-}$; see Theorem~\ref{thm:main:ER} below.

With the significance of the local energy inequality, the $L^3$-based Besov space $B^{\sfrac 13}_{3,\infty}$, and intermittency in mind, we can now introduce the $L^3$-based Onsager conjecture. 
\begin{theorem-non}[\bf $L^3$-based strong Onsager conjecture] Let $\upbeta\in (0,1)$ and $T\in (0,\infty)$. 
\begin{enumerate}[(a)]
    \item{\bf (Conservation and local energy equality)} For any $\upbeta>\sfrac13$, if a weak solution to the Euler equations belongs to $C^0([0,T]; B^{\upbeta}_{3,\infty}(\T^3))$, then it satisfies the local energy identity \eqref{thurzday:evening} in the sense of distributions with $D[u] \equiv 0$.
    \item{\bf (Dissipation and local energy inequality)} For any $\upbeta<\sfrac13$, there exist weak solutions $u$ to the Euler equations belonging to $C^0([0,T]; B^{\upbeta}_{3,\infty}(\T^3))$ which satisfy the local energy balance \eqref{thurzday:evening} in the sense of distributions, where $D[u]$ is non-negative and does not identically vanish. 
\end{enumerate}
\end{theorem-non}

The rigidity part has been established by Duchon-Robert \cite{DuchonRobert00}. For the flexibility part, on the other hand, some partial results are known. The current best result is due to the second author and De Lellis~\cite{DK22}, who showed the existence of H\"older continuous weak solutions to the Euler equations in $C^{\upbeta}_{t,x}$ for any $\upbeta<\sfrac17$ which also satisfy the strict local energy inequality \eqref{thurzday:evening}; we also refer to earlier results of De Lellis and Sz\'ekelyhidi \cite{DLeSz2010} and Isett \cite{Is22}, the latter of which formulated the strong $C^0$ Onsager conjecture. In the companion paper~\cite{GKN23}, we give a proof of the flexible side for $\upbeta \in [\sfrac17, \sfrac13)$, thus resolving the $L^3$-based strong Onsager conjecture.

\begin{theorem}[\bf Dissipation and local energy inequality \cite{GKN23}]\label{thm:main}
For any fixed $\upbeta\in (0,\sfrac{1}{3})$ and $T>0$, we can find a weak solution $u$ in $C^0_t(B^{\upbeta}_{3,\infty}\cap L^{\frac{1}{1-3\upbeta}})$ to the Euler equations \eqref{eqn:Euler} which dissipates the total kinetic energy and satisfies the local energy inequality \eqref{thurzday:evening} with $D[u]$ non-negative.
\end{theorem}

The proof of this theorem is lengthy and technical, and it is the main motivation for the present work.  However, not all of the technology which we have developed in the course of the proof of Theorem~\ref{thm:main} is specific to the construction of solutions satisfying the local energy inequality. Indeed a significant portion of this technology should be applicable in \emph{any} construction of intermittent,\footnote{Here, ``intermittent'' means that different $L^p_x$ norms satisfy very different bounds}\index{intermittent} high-regularity weak solutions to a variety of fluid equations.  For example, the technology developed here provides an improved proof of the intermittent Onsager theorem from \cite{NV22}.
\begin{theorem}[\bf Dissipation, but no local energy inequality \cite{NV22}]\label{thm:main:ER}
For any fixed $\upbeta\in (0,\sfrac{1}{3})$ and $T>0$, there exist weak solutions $u$ to \eqref{eqn:Euler} belonging to $C^0_t(B^{\upbeta}_{3,\infty}\cap L^{\frac{1}{1-3\upbeta}})$ which dissipate the total kinetic energy.
\end{theorem}

We have structured this manuscript around the following goals.  First, we develop the tools which will be used in \cite{GKN23} to prove Theorem~\ref{thm:main}.  Second, we apply these tools to give a new proof of Theorem~\ref{thm:main:ER}.  We however ensure that our application of these tools \emph{coincides} with the construction required for the proof of Theorem~\ref{thm:main}, so that we may freely reference results from this manuscript in~\cite{GKN23}.  One consequence of this goal is that certain portions of this paper are formulated in a way which is convenient for the proof of Theorem~\ref{thm:main}, but not strictly necessary for the proof of Theorem~\ref{thm:main:ER}. However, we shall always isolate and explain the results which are necessary \emph{only} for the proof of Theorem~\ref{thm:main}, so that the reader who wishes to ignore them can safely do so.  We notate these results with an asterisk; for a first example of this notation, we refer to Sections~\ref{opsection:pressure} and \ref{opsection:vel:inc:pot} from the table of contents.
\begin{remark}[\bf * Notation]\label{not.astnot}
   Throughout this article, any section, lemma, theorem, etc. which is amended with an asterisk * is only essential for the proof of Theorem~\ref{thm:main} given in~\cite{GKN23}, and not essential for our proof of Theorem~\ref{thm:main:ER} in this paper.
\end{remark}

In subsections~\ref{ss:cs}--\ref{ss:tk} of the introduction, we outline the contents of this paper, focusing respectively on the novel aspects of our wavelet-inspired scheme, the role of the intermittent pressure in pointwise estimates, and the technical tools we have developed.  Then in subsection~\ref{ss:gp}, we give two guides to the rest of the paper; one aimed at understanding the proof of Theorem~\ref{thm:main}, and the other aimed at understanding the proof of Theorem~\ref{thm:main:ER}.

\subsection{The wavelet-inspired scheme}\label{ss:cs}

As is customary in convex integration constructions of weak solutions to the Euler equations, the solution $u$ satisfying the conclusions of Theorem~\ref{thm:main:ER} will be constructed as a limit of a sequence of approximate solutions $u_q$.  In all existing iterations for the Euler equations, $u_q$ is equal (up to negligible errors) to the frequency truncation $\mathbb{P}_{\leq \la_q}$ of the limiting solution $u$, where $\la_q\to \infty$ at a slightly super-exponential rate as $q\to \infty$. As a consequence of such a construction, velocity increments\index{velocity increment} $w_q=u_q - u_{q-1}$ and $w_{q'} = u_{q'} - u_{q'-1}$ for $q'\neq q$ have no significant overlap in their active frequencies. In our new proof of Theorem~\ref{thm:main:ER}, however, $u_q$ functions as a partial wavelet decomposition of the limiting solution $u$, in the sense that $w_{q}$ and $w_{q'}$ may have frequency overlap even if $q\neq q'$.  We use the parameter $\bn$\index{$\bn$} to quantify the number of velocity increments which have frequency overlap; that is, $w_{q}$ and $w_{q'}$ have non-trivial frequency overlap if and only if $|q-q'| \leq \sfrac \bn 2$. Furthermore, the frequency support of $w_{q+1}$ in our setting is not contained in between $\la_{q}$ and $\la_{q+1}$, but rather $\la_{q+\half}$ and $\la_{\qbn}$.  To highlight this distinction, we often use the notations $\hat w_\qbn := w_{q+1}$ to emphasize that the maximum frequency of $w_{q+1}$ is $\la_\qbn$, and $\hat u_{\qbn-1} = u_q$ to emphasize that the maximum frequency present in $u_q$ is $\lambda_{\qbn-1}$.\index{$\hat u_q$}\index{$\hat w_\qbn$}

This perspective greatly affects the structure of the Euler-Reynolds system at stage $q$, which is the system satisfied by $u_q$.  In our wavelet-inspired setting, $u_q$ satisfies
\begin{equation}\label{eqn:ER:intro:new}
    \begin{cases}
    \pa_t u_q + \div \left( u_q \otimes u_q \right) + \na p_q = \div \left( R_q - \pi_q \Id \right) \\
    \div \, u_q =0  \, , 
    \end{cases}
\end{equation}
where $\kappa_q:=\sfrac 12\tr\left(R_q-\pi_q\Id\right)$.  The Reynolds stress $R_q$ and intermittent pressure $\pi_q$ can be decomposed into components
$$  R_q = \sum_{k=q}^{\qbn-1} R_q^k \, , \qquad \pi_q = \sum_{k=q}^{\qbn-1} \pi_q^k  \, . $$
The superscript $k$ indicates that the stress or pressure oscillates at frequencies no larger than $\la_k$.  The velocity increment $w_{q+1}$ is then designed to cancel out $R_q^q - \pi_q^q\Id$, leaving $R_q^k - \pi_q^k\Id$ untouched for $q+1\leq k \leq \qbn-1$.  This stands in contrast to all existing schemes, in which the entire Reynolds stress is cancelled.  

In order to replace the lack of frequency separation between various velocity increments, we instead impose that $w_{q}$ and $w_{q'}$ have disjoint spatial support if $|q-q'| < \bn$.  Therefore spatial support information is a key component of our inductive assumptions.  In order to successfully propagate the spatial support information we require, we utilize a new stationary solution to the Euler equations as our main building block, which we call an ``intermittent Mikado bundle.''  Intermittent Mikado bundles $\BB_{q+1}$ are multi-scale shear flows consisting of a product of a high frequency, highly-intermittent shear (Mikado, following~\cite{DaneriSzekelyhidi17}) flow $\WW_{\qbn}$, and an essentially homogeneous\footnote{{Homogeneous here means the opposite of intermittent.}}\index{homogeneous} shear (Mikado) flow $\rhob_{q+1}$.  The frequency support of $\WW_\qbn$ is contained in the set $[r_q \la_\qbn, \la_\qbn]$, where $r_q = \la_{q+\half}\la_\qbn^{-1}$, whereas the frequency support of $\rhob_{q+1}$ is highly concentrated around $\la_{q+1}$.  We point out that the intermittency ratio $r_q\approx (\la_q \la_\qbn^{-1})^{\sfrac 12}$ has been identified as the ``Goldilocks ratio'' in \cite{NV22} for producing solutions to the 3D Euler equations in $C^0_t B^{\sfrac 13-}_{3,\infty}$. The second key component of our spatial support toolkit is a synthetic Littlewood-Paley projector $\tP_{\leq \la_q}$, which replaces the kernel corresponding to the usual Fourier projector $\mathbb{P}_{\leq \la_q}$ onto frequencies no larger than $\la_q$ with a kernel which is compactly supported in a ball of size $\approx \la_q^{-1}$. As one would expect, the synthetic Littlewood-Paley projector obeys the usual derivative estimates with cost $\la_q$, but produces outputs supported in the fattened (by $\la_q^{-1}$) support of the input.

The flexibility afforded by the wavelet-inspired scheme and the multi-scale intermittent Mikado bundles allows us to rectify one of the seemingly unnatural components of the construction in \cite{NV22} of solutions satisfying Theorem~\ref{thm:main:ER}.  In \cite{NV22}, the velocity increment $w_{q+1}$ consisted of a collection of \emph{sub-}increments $w_{q+1,k}$, all with varying degrees of intermittency (i.e. scaling between $L^2$ and $L^\infty$ norms).  These sub-increments were designed to cancel a collection of \emph{sub-}stresses produced at a fixed stage $q\mapsto q+1$.  Our wavelet-inspired scheme instead produces a sequence of perfectly self-similar velocity increments, which obey uniform intermittent scaling laws in terms of the Goldilocks intermittency ratio.  Furthermore, there is no longer a need for the sub-stresses or sub-increments which complicated the scheme in \cite{NV22}. In this sense, our wavelet-inspired scheme is a natural generalization of the classical Fourier-inspired convex integration (Nash iteration) schemes.

\subsection{Pointwise estimates}\label{ss:ps}
One of the difficulties of an intermittent scheme, such as those in \cite{BV19, BCV, CheskidovLuo, MS} is the lack of homogeneity in estimates.  For example, inductive assumptions on $\nabla u_q$ in \cite{NV22} are propagated in $L^2$, meaning that the local $L^\infty$ norm of $\nabla u_q$ may vary greatly across different space-time regions.  This affects the stability of solutions to the transport equation with velocity $u_q$, which is used to flow the intermittent bundles (\`a la Taylor's frozen turbulence hypothesis).  Similarly, the size of the Reynolds stress $R_q$ (or $R_q^q$ in our case) will vary greatly across different space-time regions, forcing us to normalize $w_{q+1}$ as roughly $|R_q^q|^{\sfrac 12} \mathbb{B}_{q+1}$ so as to enact a quadratic cancellation between $w_{q+1}\otimes w_{q+1}$ and $R_q^q$.  One role of the intermittent pressure $\pi_q$ is to streamline these estimates by building into $\pi_q$ information regarding the local size of $R_q$, $\nabla u_q$, and their derivatives.  Our inductive estimates assert that
$$  |R_q^q| \leq \pi_q^q \, , \qquad |\nabla \hat u_{q}|^2 \leq r_{q-\bn}^{-2} \la_q ^2 \pi_q^q \, , $$
with similar bounds holding for $R_q^k$ and $\nabla \hat u_{k}$ for $k\neq q$.  Using $\Dtq$ to denote the material derivative $\pa_t + \hat u_q\cdot \nabla$, we are in fact able to show the much stronger estimates (which we refer to as ``pointwise estimates'')\index{pointwise estimates}
\begin{subequations}\label{eq:wednezday}
\begin{align}
    \left| D^N D_{t,q}^M R_q^q \right|  &\leq \pi_q^q \la_q^N \left(r_{q-\bn}^{-1} \la_q (\pi_q^q)^{\sfrac 12}\right)^{M} \,, \\
    \left| D^N \Dtq^M \nabla \hat u_q \right|  &\leq r_{q-\bn}^{-1} \la_q (\pi_q^q)^{\sfrac 12} \la_q^N \left(r_{q-\bn}^{-1} \la_q (\pi_q^q)^{\sfrac 12}\right)^{M} \\
    \left| D^N D_{t,q}^M \pi_q^q \right|  &\leq \pi_q^q \la_q^N \left(r_{q-\bn}^{-1} \la_q (\pi_q^q)^{\sfrac 12}\right)^{M} \,. 
\end{align}
\end{subequations}
These estimates show that we can use $\pi_q^q$ in conjunction with the parameters $\la_q$ and $r_q$ as multiplicative factors controlling the pointwise size of both spatial and material derivatives on $R_q^q$, $\nabla \hat u_q$, and $\pi_q^q$. While we still choose to formulate estimates in terms of carefully constructed cutoff functions as in \cite{NV22}, the intermittent pressure serves to centralize all the necessary size and frequency information needed throughout the iteration.

\subsection{Toolkit}\label{ss:tk}

At a technical level, this manuscript contains generalizations of a number of the tools from \cite{BMNV21} and \cite{NV22}.  First among these is a sharp $L^p$ decoupling estimate for products $fg$, where $f$ has maximum effective frequency $\la$ and $g$ is periodic to scale $\Lambda^{-1} \ll \lambda^{-1}$.  Estimates for such a product in $L^1$ and $L^2$ were first shown by Buckmaster and Vicol in \cite{BV19}.  We generalize this estimate to any $p\in[1,\infty]$.   With a sharp $L^p$ decoupling estimate in hand, we construct an inverse divergence operator inspired by \cite{BMNV21} which is well-adapted to error terms of the form $f g \circ \Phi$, where $f$ and $g$ satisfy the same properties as above, $\Phi$ is a low frequency flow map, and $g$ can be written as the iterated divergence of a tensor potential $\div^\dpot G = g$.  Our inverse divergence operator can produce estimates in any Lebesgue space, propagates arbitrarily large numbers of spatial and material derivatives, preserves the spatial support of the inputs $f$ and $G$, and can be iterated an arbitrarily large number of times.  Finally, we have generalized the cutoff machinery developed in \cite{BMNV21} for intermittent functions with $L^1$ or $L^2$ estimates to intermittent functions with $L^p$ estimates for any $p\in[1,\infty)$; for the sake of convenience and concreteness, we specify to the cases ($L^1$, $L^{\sfrac 32}$, and $L^3$) which are used to measure current errors, stress and pressure errors, and velocity fields in the proofs of Theorems~\ref{thm:main} and \ref{thm:main:ER}.

\subsection{Guides to the paper}\label{ss:gp}
We present guides to Theorem~\ref{thm:main} and Theorem~\ref{thm:main:ER}.

\subsubsection{Guide to Theorem~\ref{thm:main}}

The reader interested in the proof of Theorem~\ref{thm:main} can proceed as follows.  While the inductive assumptions which are included in Proposition~\ref{prop:main} can be found in section~\ref{ss:is} of this paper, the presentation in \cite[section~2]{GKN23} is better adapted to Theorem~\ref{thm:main}.  The reader can also read the proof of Theorem~\ref{thm:main}, assuming the main inductive proposition, from \cite[section~2]{GKN23}.  The next step is to understand the construction of the velocity increment, which is carefully carried out in sections~\ref{ss:bundles}--\ref{s.corr-ec-tor} of this paper and \cite[section~4]{GKN23}.  Specifically, section~\ref{ss:bundles} introduces the definitions and estimates for intermittent Mikado bundles, section~\ref{ss:nic} presents the required non-inductive cutoffs, and section~\ref{s.corr-ec-tor} constructs and estimates the velocity increment, \emph{except} for the placement procedure for the bundles, which is the content of \cite[section~4]{GKN23}.  General readers will at this point be mostly prepared to read the estimates of the primary error terms in both the Euler-Reynolds system and the relaxed local energy inequality, which are contained in section~\ref{sec:ER:main} and \cite[section~5]{GKN23}.  The experienced reader may instead choose to skip some more familiar details, such as the Euler-Reynolds estimates, by reading \cite[sections~3,4]{GKN23}, which contain a streamlined version, mostly without proofs, of the construction of the velocity increment and the main estimates related to Euler-Reynolds errors.  One caveat for all readers is that most error terms, even the more familiar Euler-Reynolds errors, require the construction of a pressure increment.  The abstract construction of intermittent pressure increments is contained in section~\ref{opsection:pressure}, and is rather technical, so the reader may read one or two of the proofs and then take the rest of the statements for granted.  Sections~\ref{sec:inductive.cutoffs} and \ref{opsection:vel:inc:pot} are also rather technical in nature, and the reader need not understand the details of these sections in order to follow the rest of the proof.  The former section constructs the inductive (velocity) cutoffs, while the latter shows that the velocity increment can be written as an iterated Laplacian, roughly speaking, and constructs the pressure increment associated with the velocity.  Readers wishing to skip this step may find summaries of these tools in \cite[section~3 and Appendix~A]{GKN23}.  The last crucial step in the proof of Theorem~\ref{thm:main} is the construction of the intermittent and anticipated pressure in \cite[section~6]{GKN23}.  Finally, the reader should consult section~\ref{sec:parm} for parameter choices and inequalities and appendix~\ref{sec:app:tools} for technical tools as needed.

\subsubsection{Guide to Theorem~\ref{thm:main:ER}}
The reader interested in the proof of Theorem~\ref{thm:main:ER} need not consult \cite{GKN23} at any point. The reader has two options, the first of which is to follow the outline to the proof given in the proof of Proposition~\ref{prop:main:ER}, which includes a treatment of the intermittent pressure $\pi_q$.  Alternatively, the reader who prefers to ignore the intermittent pressure can follow the outline given in Remark~\ref{rem:toozday}, which replaces the intermittent pressure with methodology more similar to that of \cite{NV22}.

\subsection*{Acknowledgements} 
The authors acknowledge the hospitality and working environment at the Institute for Advanced Study during
the special year on the h-principle, when they first started working on this project. VG was supported by the NSF under Grants DMS-FRG-1854344 and DMS-1946175 while at Princeton University. HK was supported by the NSF under Grant DMS-1926686 while a member at the IAS and would like to extend gratitude to the employer, the Forschungsinstitut f\"ur Mathematik (FIM) at ETH Z\"urich. MN was supported by the NSF under Grant DMS-1926686 while a member at the IAS.  The authors thank Camillo De Lellis and Vlad Vicol for their commentary on a draft of this manuscript.

\section{Inductive propositions and proofs of main theorems}\label{ss:is}

In this section, we present the main inductive assumptions and propositions required for Theorems~\ref{thm:main} and \ref{thm:main:ER}.  The inductive assumptions which are required for Theorem~\ref{thm:main} but not \ref{thm:main:ER} are sorted into subsection~\ref{ss:inductive:LEI}.  Then in subsection~\ref{ss:sunday}, we present the inductive propositions required for both the main theorems and outline how the contents of this paper contribute to the proofs of Theorems~\ref{thm:main} and \ref{thm:main:ER}.

\subsection{General notations and parameters}\label{sec:not.general}

We first introduce the primary parameters 
$$ \beta\, , \bn\, , b\, , \la_q\, , \delta_q\, , r_q \, , \Ga_q  \, , \varepsilon_\Gamma $$
which appear in the inductive hypotheses. First, we choose an $L^3$-based regularity index $\beta\in[\sfrac17,\sfrac 13)$. Since $\beta<\sfrac 13$, we can choose $\bn\in 6\N$ such that\index{$\beta$} \index{$\bn$}
\begin{equation}
        \beta < \frac 13 \cdot \frac{\sfrac \bn 3}{\sfrac \bn 3 + 2} - \frac{2}{\sfrac \bn 3 +2} \, , \qquad \beta< \frac {2}{3} \cdot \frac{\sfrac{\bn}{2}-1}{\bn} \, . \label{eq:choice:of:bn}
    \end{equation}
This in turn enables a choice of $b\in (1,\sfrac{25}{24})$ close to $1$ such that\index{$\beta$}\index{$\bn$}\index{$b$}
 \begin{subequations}\label{ineq:b}
    \begin{align}
        \beta < \frac{1}{3b^{\bn}} \cdot \frac{1+b+\dots+b^{\sfrac \bn 3 -1}}{1+b+\dots+b^{\sfrac \bn 3 +1}}& - \frac{2\left(1+(b-1)(1+\cdots+b^{\sfrac \bn 2 -1})^2\right)}{1+b+\dots+b^{\sfrac \bn 3 +1}} \, , \quad  \frac{2}{3b^{\sfrac{\bn}{2}}} \cdot \frac{1+\dots+b^{\sfrac \bn 2 - 2}}{1+\dots+b^{\bn-1}} \label{ineq:b:first} \\
        b^{\bn} < 2& \, , \qquad \frac{(b^{\sfrac \bn 2 - 1}+\dots+b+1)^2}{b^{\sfrac \bn 2 -1}+\dots+b+1} (b-1) < (b-1)^{\sfrac 12} \, . \label{ineq:b:second}
    \end{align}
    \end{subequations}
Indeed the inequalities in \eqref{ineq:b:first} are possible since \eqref{eq:choice:of:bn} is just \eqref{ineq:b:first} evaluated at $b=1$, and both expressions in \eqref{ineq:b:first} are continuous in $b$ in a neighborhood of $b=1$. The first inequality in \eqref{ineq:b:second} is trivial, and the second is possible since the fraction in the expression is continuous at $b=1$ and equal to $\sfrac \bn 2$ if $b=1$. It is clear that as $\beta\rightarrow \sfrac 13$, we are forced to choose $\bn \rightarrow \infty$ and $b\rightarrow 1$.
    
We now define the frequency parameter $\la_q$, the amplitude parameter $\de_q$, the intermittency parameter $r_q$, and the multi-purpose parameter $\Gamma_q$ by \index{$\la_q$} \index{$\de_q$} \index{$r_q$} \index{$\Ga_q$}
\begin{align}
    &\la_q = 2^{\left \lceil (b^q) \log_2 a \right \rceil}\approx a^{(b^q)} \, , \qquad  \de_q = \la_q^{-2\beta} \, ,
    \label{eq:def:la:de}\\
    r_q = \frac{\la_{q+\half}\Ga_q}{\la_\qbn}&\, , \qquad
    \Ga_q = 2^{\left \lceil \varepsilon_\Gamma \log_2 \left( \frac{\la_{q+1}}{\la_q} \right) \right \rceil}  
    \approx \left( \frac{\la_{q+1}}{\la_q} \right)^{\varepsilon_\Gamma} \approx \la_{q}^{(b-1)\varepsilon_\Gamma} \, .   \label{eq:deffy:of:gamma:intro}
    \end{align}
The large positive integer $a$ and the small positive number $0<\varepsilon_\Gamma\ll (b-1)^2<1$ are defined in \eqref{i:choice:of:a} and \eqref{i:par:4} of subsection \ref{sec:para.q.ind}, respectively.  Note that the intermittency parameter $r_q$ is determined by the ``$\sfrac 12$ rule'' as in \cite{NV22}.

We now introduce further parameters
$$ \tau_q \, ,  \Lambda_q\, , \Tau_q\, , \badshaq\, . $$
We shall often decompose $u_q=\hat u_q + (u_q-\hat u_q)$, and heuristically speaking, the gradient of velocity $\nabla\hat u_{q'}$ will have spatial derivative cost $\approx\la_{q'}$ and $L^3$ norm $ \approx \tau_{q'}^{-1}\approx \de_{q'}^{\sfrac12} r_{q'-\bn}^{-\sfrac13}\la_{q'}$.  We in fact adjust the definition of $\tau_q^{-1}$ using the parameter $\La_q$ (slightly larger than $\la_q$), which accounts for small spatial frequency losses due to mollification, and introduce the parameter $\Tau_q^{-1}$ (much larger than $\tau_q^{-1})$, which accounts for temporal frequency losses due to mollification.  We set
\begin{align}\label{eq:defn:tau}
\la_q < \La_q = {\la_q \Ga_q^{10} } \, , \qquad 
\tau_q^{-1} &= \delta_q^{\sfrac 12} \la_q r_{q-\bn}^{-\sfrac 13} \Ga_q^{35} \ll \Tau_q^{-1} \, ,
\end{align}
and refer to \eqref{v:global:par:ineq} for the precise definition of $\Tau_q$. For the $L^\infty$ norm of $R_q^q$ (and other inductive objects), we use the parameter $\badshaq$, which will satisfy (we refer to \eqref{eq:badshaq:choice} for the precise choice of $\badshaq$)
 \begin{align*}
    \la_q^{\frac 1{\bn}}\les\Ga_q^{\badshaq} \les \la_q^{\frac{12}{\bn}}\, . 
\end{align*}
 
Finally, we will inductively propagate spatial and material derivative estimates, where we use the notation and parameters
\begin{align*}
    D_{t,q} = \pa_t  + (\hat u_q\cdot \na)\, , \qquad \NcutSmall\, , \Nindt\, , \Nind\, , \Nfin \, .
\end{align*}
The integers $\mathsf{N}_{\bullet}$ above quantify the number of spatial and material derivative estimates propagated inductively and satisfy the ordering (see subsection~\ref{sec:para.q.ind} for the precise choices)\index{$\NcutLarge$}\index{$\NcutLarge$}\index{$\Nind$}\index{$\Nfin$}
\begin{align*}
1\ll \NcutSmall\ll \Nindt \ll  \Nind \ll \Nfin\, . 
\end{align*}
In particular, $\Nindt$ helps us keep track of both sharp and lossy material derivative estimates.  For this purpose, we use the following notation, which roughly says that ``the first $N_*$ material derivatives cost $\tau^{-1}$, while additional derivatives cost $\Tau^{-1}$.'' We also list a few other notations in the subsequent two remarks.\index{$\MM{n,N_*,\tau^{-1},\Tau^{-1}}$}
\begin{remark}[\bf Geometric upper bounds with two bases]\label{not:M}
    For all $n\geq 0$, we define
    \[\MM{n,N_*,\tau^{-1},\Tau^{-1}} := \tau^{-\min\{n,N_*\}} \Tau^{-\max\{n-N_*,0\}}\,.\]
\end{remark}

\begin{remark}[\bf Space-time norms]\label{rem:notation:space:time}
In the remainder of the paper, we shall always measure objects using uniform-in-time norms $\sup_{t\in[T_1,T_2]}\|\cdot (t)\|$, where $\| \cdot (t)\|$ is any of a variety of norms used to measure functions defined on $\T^3\times[T_1,T_2]$ but restricted to time $t$. In a slight abuse of notation, we shall always abbreviate these space-time norms with simply $\| \cdot \|$.
\end{remark}

\begin{remark}[\bf Space-time balls]\label{rem:notation:space:time:balls}\index{$B(\Omega, \la^{-1})$}\index{$B(\Omega, \la^{-1}, \tau)$}
For any set $\Omega\subseteq \mathbb{T}^3\times \mathbb{R}$, we shall use the notations
\begin{subequations}\label{eq:space:time:balls}
\begin{align}
B(\Omega,\lambda^{-1}) &:= \left\{ (x,t) \, : \, \exists \, (x_0,t) \in \Omega \textnormal{ with } |x-x_0| \leq \lambda^{-1} \right\}\\
B(\Omega,\lambda^{-1},\tau) &:= \left\{ (x,t) \, : \, \exists \, (x_0,t_0) \in \Omega \textnormal{ with } |x-x_0| \leq \lambda^{-1} \, , |t-t_0| \leq \tau \right\}
\end{align}
\end{subequations}
for space and space-time neighborhoods of $\Omega$ of radius $\lambda^{-1}$ in space and $\tau$ in time, respectively. 
\end{remark}

\subsection{Relaxed equations}\label{ss:relaxed}
We assume that there exists an \emph{approximate solution} $(u_q, p_q, R_q, -\pi_q)$ at the $q^{\textnormal{th}}$ step, $q\geq 0$, where $u_q:\T^3\times[-\tau_{q-1},T+\tau_{q-1}]\footnote{We adopt the convention that $\tau_{-1}:=1$.}\rightarrow \R^3$ is the velocity field, $p_q:\T^3\times[-\tau_{q-1},T+\tau_{q-1}]\rightarrow \R$ is the pressure, $R_q:\T^3\times[-\tau_{q-1},T+\tau_{q-1}]\rightarrow \R^{3\times 3}_{\textnormal{symm}}$ is the symmetric stress tensor, and $\pi_q:\T^3\times[-\tau_{q-1},T+\tau_{q-1}]\rightarrow \R$ is a scalar field which we shall refer to as the \emph{intermittent pressure}\index{intermittent pressure}.  We assume that the approximate solution satisfies the Euler-Reynolds system\index{Euler-Reynolds system}
\begin{align}\label{eqn:ER}
\begin{cases}
\partial_t u_q + \div (u_q \otimes u_q) + {\nabla p_q} = \div(R_q -\pi_q \Id) \\
\div\, u_q = 0 \, .
\end{cases}
\end{align}
We use the decomposition and notations
\begin{equation}\label{eq:hat:no:hat}
u_q = \underbrace{\hat u_{q-1} + \hat w_q}_{=: \hat u_q} + \hat w_{q+1} + \dots + \hat w_{q+\bn-1} =: \hat u_{q+\bn-1} \, 
\end{equation}
for the velocity field; one purpose of the notation $\hat u_{\qbn-1}$ is to emphasize that $u_q$ has effective maximum spatial frequency $\lambda_{\qbn-1}$.  The stress error \index{stress error} $R_q$\index{$R_q$} has a decomposition
\begin{align}
R_q &= \sum_{k=q}^{q+\bn-1} R_q^k \, , \label{eq:ER:decomp:basic} 
\end{align}
where each $R_q^k$\index{$R_q^k$} is a symmetric stress tensor. The intermittent pressure $\pi_q$ has a decomposition\index{$\pi_q^k$}
\begin{align}
\pi_q &= \sum_{k=q}^{{\infty}} \pi_q^k \, . \label{eq:pi:decomp:basic}
\end{align}

In our wavelet-inspired scheme, the Reynolds stress $R_q$ will have a wide band of frequency support in between $\lambda_q$ and $\lambda_{q+\bn-1}$ (effectively speaking). We correct the portion of it which lives at frequencies no higher than $\lambda_q$. We denote this portion by $R^q_q$. More generally, we denote the portions of $R_q$ with spatial derivative cost $\lambda_{k}$ by $R_q^k$. 

\subsection{Inductive assumptions for velocity cutoff functions}
\label{sec:cutoff:inductive}
The inductively-defined velocity cutoff functions $\psi_{i,q'}$ partition space-time into distinct level sets of the gradient of velocity. We first record here the key properties which will be required throughout the inductive assumptions, and the local $L^\infty$ estimates for velocity increments $\hat w_{q'}$ and velocity $\hat u_{q'}$, obtained as a consequence of the definition of $\psi_{i,q'}$, can be found in subsection \ref{sec:inductive:secondary:velocity}. The concrete construction of $\psi_{i,q+\bn}$ and the verification of \eqref{eq:inductive:partition}--\eqref{eq:inductive:timescales} for $q\mapsto q+1$ (i.e., $q'=q+\bn$) will be given in Section \ref{sec:inductive.cutoffs}.\index{velocity cutoffs}

All assumptions in subsection~\ref{sec:cutoff:inductive} are assumed to hold for all $0\leq q-1\leq q'\leq q+\bn-1$.  First, we assume that the velocity cutoff functions form a partition of unity:
\begin{align}\label{eq:inductive:partition}
    \sum_{i\geq 0} \psi_{i,q'}^6 \equiv 1, \qquad \mbox{and} \qquad \psi_{i,q'}\psi_{i',q'}=0 \quad \textnormal{for}\quad|i-i'| \geq 2 \, .
\end{align}
Second, we assume that there exists an $\imax = \imax(q') \geq 0$, which is bounded uniformly in $q'$ by
\begin{align}
\imax(q') \leq \frac{\badshaq+12}{(b-1)\varepsilon_\Ga} \, ,
\label{eq:imax:upper:lower}
\end{align}
such that
\begin{align}
\psi_{i,q'} &\equiv 0 \quad \mbox{for all} \quad i > \imax(q')\,,
\qquad \mbox{and} \qquad
\Gamma_{q'}^{\imax(q')} \leq 
\Ga_{q'-\bn}^{\sfrac{\badshaq}{2}+18}
 \delta_{q'}^{-\sfrac 12}r_{q'-\bn}^{-\sfrac {2}3} \, .
\label{eq:imax:old}
\end{align}
For all $0 \leq i \leq \imax$, we assume the following pointwise derivative bounds for the cutoff functions $\psi_{i,q'}$. First, for mixed space and material derivatives and multi-indices $\alpha,\beta \in {\mathbb N}^k$, $k \geq 0$, $0 \leq |\alpha| + |\beta| \leq \Nfin$, we assume that
\begin{align}
&\frac{{\bf 1}_{\supp \psi_{i,q'}}}{\psi_{i,q'}^{1- (K+M)/\Nfin}} \left|\left(\prod_{l=1}^k D^{\alpha_l} D_{t,q'-1}^{\beta_l}\right) \psi_{i,q'}\right| \leq \Gamma_{q'} (\Gamma_{q'}  \lambda_{q'})^{|\alpha|} 
\MM{|\beta|,\NindSmall - \NcutSmall,  \Gamma_{q'}^{i+3}  \tau_{q'-1}^{-1}, \Gamma_{q'-1} \Tau_{q'-1}^{-1}} \, .
\label{eq:sharp:Dt:psi:i:q:old}
\end{align}
Next, with $\alpha, \beta,k$ as above, $N\geq 0$ and $D_{q'}:=\hat w_{q'}\cdot\nabla$, we assume that
\begin{align}
&\frac{{\bf 1}_{\supp \psi_{i,q'}}}{\psi_{i,q'}^{1- (N+K+M)/\Nfin}} \left| D^N \left( \prod_{l=1}^k D_{q'}^{\alpha_l} D_{t,q'-1}^{\beta_l}\right)  \psi_{i,q'} \right| \notag\\
&\qquad \leq \Gamma_{q'} ( \Gamma_{q'}  \lambda_{q'})^N
(\Gamma_{q'}^{i-5} \tau_{q'}^{-1})^{|\alpha|}
\MM{|\beta|,\Nindt-\NcutSmall,  \Gamma_{q'}^{i+3}  \tau_{q'-1}^{-1}, \Gamma_{q'-1}  \Tau_{q'-1}^{-1}}
\label{eq:sharp:Dt:psi:i:q:mixed:old}
\end{align}
for $0 \leq N+ |\alpha| + |\beta| \leq \Nfin$. Finally, for $0\leq i \leq \imax(q')$, we assume the $L^1$ bound\index{$\CLebesgue$}
\begin{align}
\norm{\psi_{i,q'}}_{1}  \leq \Gamma_{q'}^{-3i+\CLebesgue} \qquad \mbox{where} \qquad \CLebesgue = \frac{6+b}{b-1} \, .
\label{eq:psi:i:q:support:old}
\end{align}
Lastly, we assume that local timescales dictated by velocity cutoffs at a fixed point in space-time are decreasing in $q$. More precisely, for all $q' \leq q+\bn-1$ and all $q''\leq q'-1$, we assume 
\begin{equation}
    \psi_{i',q'}  \psi_{i'',q''} \not \equiv 0 \quad \implies \quad 
    \tau_{q'} \Gamma_{q'}^{-i'} \leq \tau_{q''} \Gamma_{q''}^{-i'' -25}   \, . \label{eq:inductive:timescales}
\end{equation}
This will be useful when we upgrade material derivative from $D_{t,q''}$ to $D_{t,q'}$.

\subsection{Inductive bounds on the intermittent pressure \texorpdfstring{$\pi_q$}{piqq}}\label{sec:pi:inductive}

The intermittent pressure $\pi_q$ is designed to majorize derivatives of errors and velocity increments pointwise. In this subsection, we introduce estimates for $\pi_q$ which are part of the proof of Theorem~\ref{thm:main:ER}, and establish precise relations between the intermittent pressure and errors/velocity increments. The reader who is interested in the proof of Theorem~\ref{thm:main} should refer to \cite[subsection~2.4]{GKN23} for a complete listing of the inductive assumptions related to the intermittent pressure.  On the other hand, the reader only interested in the proof of Theorem~\ref{thm:main:ER} can refer to the proof of Proposition~\ref{prop:main:ER} for an outline of how to verify the inductive assumptions from this subsection.  Alternatively, it is possible to prove Theorem~\ref{thm:main:ER} by treating the more familiar $L^p$ bounds on the Reynolds stress in Remark~\ref{sec:stress:inductive} as the main inductive assumptions and ignoring the rest of the content of this subsection. This approach is completely analogous to that of \cite{NV22}, and we discuss this further in Remark~\ref{rem:toozday}.

\subsubsection{\texorpdfstring{$L^{\sfrac 32}$}{tpdfs1}, \texorpdfstring{$L^\infty$}{tpdfs2}, and pointwise bounds for \texorpdfstring{$\pi_q^k$}{tpdfs3}}

We assume that for $q\leq k \leq q+\bn-1$ and $N+M\leq 2\Nind$, $\pi_q^k$ satisfies 
\begin{subequations}\label{eq:pressure:inductive}
\begin{align}
\norm{ \psi_{i,k-1} D^N D_{t,k-1}^M  \pi_q^k }_{\sfrac 32}  
 &\leq \Ga_q\Ga_k \delta_{k+\bn} \Lambda_k^N \MM{M, \NindRt, \Gamma_{k-1}^{i} \tau_{k-1}^{-1} ,  \Tau_{k-1}^{-1} } \, .
\label{eq:pressure:inductive:dtq} \\
\norm{ \psi_{i,k-1} D^N D_{t,k-1}^M \pi_q^k }_{\infty}  
 &\leq \Ga_q {\Gamma_k^{\badshaq+1}} \Lambda_k^N \MM{M, \NindRt, \Gamma_{k-1}^{i} \tau_{k-1}^{-1} ,  \Tau_{k-1}^{-1} } \, , \label{eq:pressure:inductive:dtq:uniform} \\
 \label{eq:ind:pi:by:pi}
    \left|\psi_{i,k-1} D^N D_{t,k-1}^M \pi_q^{k}\right| &\leq \Ga_q {\Gamma_k}\pi_q^k  \Lambda_k^N \MM{M, \NindRt, \Gamma_{k-1}^{i} \tau_{k-1}^{-1} , \Tau_{k-1}^{-1} } \, .
\end{align}
\end{subequations}
Throughout the paper, we shall use the phrase ``pointwise estimates'' to refer to bounds on stress errors, current errors, or velocities in terms of various $\pi$'s which resemble the third bound in either of the above displays.

\subsubsection{Pointwise bounds for errors, velocities, and velocity cutoffs}
We assume that we have the pointwise estimates\index{pointwise estimates}
\begin{subequations}\label{eq:inductive:pointwise}
\begin{align}
    \label{eq:ind:stress:by:pi}
    \left|\psi_{i,k-1} D^N D_{t,k-1}^M R^k_q\right| &< \Ga_q \Ga_k^{-8} \pi^k_q \Lambda_k^N \MM{M, \NindRt, \Gamma_{k-1}^{i  +20} \tau_{k-1}^{-1} ,  \Tau_{k-1}^{-1}\Ga_{k-1}^{10} } \, ,\\
    \label{eq:ind:velocity:by:pi}
    \left|\psi_{i,k-1} D^N D_{t,k-1}^M \hat w_k \right| &< \Ga_q r_{k-\bn}^{-1} (\pi_q^k)^{\sfrac12} {\Lambda}_k^{N} \MM{M, \NindRt, \Gamma_{k-1}^{i} \tau_{k-1}^{-1} ,  \Tau_{k-1}^{-1} {\Ga_{k-1}^2}} \, ,
\end{align}
\end{subequations}
where the first bound holds for $q\leq k\leq q+\bn-1$ and $N+M \leq 2\Nind$, and the second bound holds for $N+M\leq \sfrac{3\Nfin}{2}$.  While the main $L^p$ estimates on the Reynolds stress will follow from the pointwise estimates in terms of the pressure (see Remark~\ref{sec:stress:inductive}), we are forced to assume that $R_q^k$ has a decomposition $R_q^k = R_q^{k,l}+R_q^{k,*}$ \index{$R_q^{k,l},\, R_q^{k,*}$}, where $R_q^{k,*}$ satisfies the stronger bound 
\begin{align}
\norm{ D^N D_{t,k-1}^M  R_q^{k,*} }_{\infty}
\leq \Ga_q^2 \Tau_k^{2\Nindt}\delta_{k+2\bn} \Lambda_k^N \MM{M,\Nindt, \tau_{k-1}^{-1},\Tau_{k-1}^{-1}}
\label{eq:Rnnl:inductive:dtq} 
\end{align}
for all $N+M \leq 2\Nind$. The extra superscript $l$ stands for ``local," in the sense that $R_q^{k,l}$ is a stress error over which we maintain control of the spatial support, whereas $\ast$ refers to non-local terms which are negligibly small.  The reader can safely ignore such non-local error terms.

Finally, we assume that for all $q\leq q'\leq q+\bn-1$,
\begin{align}\label{eq:psi:q:q'}
    \sum_{i=0}^{\imax} \psi_{i,q'}^2 \delta_{q'} r_{q'-\bn}^{-\sfrac23} \Gamma_{q'}^{2i} &\leq {2^{q-q'}} \Ga_{q'} {r_{q'-\bn}^{-2}} \pi_{q}^{q'} \, .
\end{align}
Combining this bound with~\eqref{eq:nasty:D:vq:old} and~\eqref{eq:defn:tau} shows that for $N+M\leq \sfrac{3\Nfin}{2}$,
\begin{equation}\notag
\left|  D^N \Dtq^M \nabla \hat u_q \right| \leq \Ga_q^{50} r_{q-\bn}^{-1} \La_q (\pi_q^q)^{\sfrac 12} \La_q^n \left( \Ga_q^{50} r_{q-\bn}^{-1} \La_q (\pi_q^q)^{\sfrac 12} \right)^N \, .
\end{equation}
\begin{remark}[\bf Velocity cutoffs, timescales, and intermittent pressure]\label{rem.summing.psi}
Using the timescale parameter $\tau_q^{-1}\approx \delta_q^{\sfrac 12}\la_q r_{q-\bn}^{-\sfrac 13}$ defined precisely in subsection~\ref{sec:para.q.ind}, item~\eqref{i:par:4}, we may now record the following version of \eqref{eq:psi:q:q'} for $q'=q$;
\begin{align}
    \psi_{i,q} \tau_q^{-1} \Ga_q^i \leq \la_q \Ga_q \left(\pi_q^q\right)^{\sfrac12} r_q^{-1} \, .
\end{align}
\end{remark}
\begin{remark}[\bf $L^p$ estimates on Reynolds errors from pointwise estimates]\label{sec:stress:inductive}
The estimates on $R^k_q$ in \eqref{eq:ind:stress:by:pi} and the estimates on $\pi_q^k$ in \eqref{eq:pressure:inductive} imply that for $q\leq k \leq q+\bn-1$ and $N+M\leq 2\Nind$, $R_q^k$ satisfies
\begin{subequations}\label{eq:Rn:inductive}
\begin{align}
\norm{ \psi_{i,k-1} D^N D_{t,k-1}^M  R_q^k }_{\sfrac 32}  
 &\leq \Ga_q^2 \Ga_k^{-7} \delta_{k+\bn} \Lambda_k^N \MM{M, \NindRt, \Gamma_{k-1}^{i+20} \tau_{k-1}^{-1} ,  \Tau_{k-1}^{-1}\Gamma_q^{10} } \, ,
\label{eq:Rn:inductive:dtq} \\
\norm{ \psi_{i,k-1} D^N D_{t,k-1}^M R_q^k }_{\infty}  
 &\leq \Ga_q^2 \Ga_k^{-7} {\Gamma_k^{\badshaq}} \Lambda_k^N \MM{M, \NindRt, \Gamma_{k-1}^{i+20} \tau_{k-1}^{-1} ,  \Tau_{k-1}^{-1}\Gamma_q^{10} }\, . \label{eq:Rn:inductive:dtq:uniform}
\end{align}
\end{subequations}
\end{remark}

\subsection{Dodging principle ingredients}\label{ss:dodging}
As discussed in the introduction, one of the crucial elements for the wavelet-inspired scheme is dodging between velocity increments, which is elaborated upon in Hypothesis~\ref{hyp:dodging1}. To construct a new velocity increment with such dodging, it is necessary to keep a record of the density of previous velocity increments as stated in Hypothesis~\ref{hyp:dodging2}. 
These two hypotheses can be seen as improved and inductive versions of the ``pipe dodging'' technique used in \cite{BMNV21} or \cite{NV22}, and will be verified rigorously for $q\mapsto q+1$ in \cite[section~4]{GKN23}.  We however outline the main heuristics behind the proof following the statement of Lemma~\ref{lem:dodging}.  
\begin{hypothesis}[\bf Effective dodging]\label{hyp:dodging1}
For $q',q''\leq q+\bn-1$ that satisfy $0<|q''-q'|\leq \bn-1$, we have that\footnote{Here we are considering the support of $\hat w_{q}$ in time and space, then expanding to a ball of radius $\la_q^{-1}\Ga_{q+1}$ in space only; see \eqref{eq:space:time:balls}.}
\begin{equation}\label{eq:ind:dodging}
    B\left(\supp \hat w_{q'} , \lambda_{q'}^{-1}\Gamma_{q'+1} \right) \cap  B\left( \supp \hat w_{q''} , \lambda_{q''}^{-1}\Gamma_{q''+1} \right)= \emptyset \, .
\end{equation}
\end{hypothesis}

\begin{hypothesis}[\textbf{Density of old pipe bundles}]\label{hyp:dodging2}
There exists a $q$-independent constant $\const_D$ such that the following holds.  Let $\bar q', \bar q''$ satisfy $q \leq \bar q'' < \bar q'\leq q+\bn-1$, and set\footnote{The reasoning behind the choice of $d(\bar q', \bar q'')$ is as follows. The set should be small enough that it can be contained in the support of a single $\bar q''$ velocity cutoff.  Since these functions oscillate at frequencies no larger than $\approx\lambda_{q''}$, the first number inside the minimum ensures that this is the case.  The set should also be no larger than the size of a periodic cell for pipes of thickness $\bar q'$, which is ensured by the second number inside the minimum.}
\begin{equation} \label{eq:diam:def}
d(\bar q', \bar q'') :=  \min\left[ (\lambda_{\bar q''}\Gamma_{\bar q''}^7)^{-1} , (\lambda_{\bar q' - \half} \Gamma_{\bar q' - \bn})^{-1} \right] \, .
\end{equation}
Let $t_0\in\R$ be any time and $\Omega\subset\T^3$ be a convex set of diameter at most $d(\bar q', \bar q'')$.  Let $i$ be such that $\Omega \times \{t_0\}\cap \supp \psi_{i,\bar q''} \neq \emptyset$. 
Let $\Phi_{\bar q''}$ be the flow map such that
\begin{align*}
\begin{cases}
\pa_t \Phi_{\bar q''} + \left(\hat u_{\bar q''} \cdot \na \right) \Phi_{\bar q''} = 0\\
\Phi_{\bar q''}(t_0,x) = x \, .
\end{cases}
\end{align*}
We define $\Omega(t)=\Phi_{\bar q''}(t)^{-1}(\Omega)$.\footnote{For any set $\Omega'\subset \T^3$, $\Phi_{\bar q''}(t)^{-1}(\Omega')=\{x\in \T^3: \Phi_{\bar q''}(t,x) \in \Omega'\}$. We shall also sometimes use the notation $\Omega\circ \Phi_{\bar q''}(t)$.} 
Then there exists a set\footnote{Heuristically this set is $\cup_t \suppp_x \hat w_{\bar q'}(\cdot,t) \cap \Omega(t)$, but in order to ensure that $(\partial_t+\hat u_{\bar q''} \cdot \nabla)\mathbf{1}_L\equiv 0$, $L$ does not include any ``time cutoffs" which turn pipes on and off.} $L=L(\bar q',\bar q'', \Omega, {t_0})\subseteq \T^3\times \R$ such that for all $t\in(t_0-\tau_{\bar q''}\Gamma_{\bar q''}^{-i+2},t_0+\tau_{\bar q''}\Gamma_{\bar q''}^{-i+2})$,
\begin{align}\label{eq:ind:dodging2}
    (\partial_t + \hat u_{\bar q''}\cdot\nabla) \mathbf{1}_{L}(t,\cdot) \equiv 0 \qquad \textnormal{and} \qquad \supp_x \hat w_{\bar q'}(x,t) \cap \Omega(t) \subseteq L \cap \{t\} \, .
\end{align}
Here, the first identity holds in distribution sense. Furthermore, there exists a finite family of Lipschitz curves $\{\ell_{j,L}\}_{j=1}^{\const_D}$ of length at most $2d(\bar q', \bar q'')$ which satisfy 
\begin{equation}\label{eq:concentrazion}
    L \cap \{t=t_0\} \subseteq \bigcup_{j=1}^{\const_D} B\left( \ell_{j,L} , 3\la_{\bar q'}^{-1} \right) \, .
\end{equation}
\end{hypothesis} 
\begin{remark}[\bf Segments of deformed pipes of thickness $\la_{\bar q'}^{-1}$]\label{rem:deformed:pipes}
We will sometimes refer to a $3\la_{\bar q'}^{-1}$ neighborhood of a Lipschitz curve of length at most $2(\la_{\bar q'-\sfrac{\bn}{2}}\Ga_{\bar q'-\bn})^{-1}$ as a ``segment of deformed pipe" - see Definition~\ref{def:sunday:sunday}. Since $(\la_{\bar q'-\sfrac{\bn}{2}}\Ga_{\bar q'-\bn})^{-1}$ will be the scale to which our high-frequency pipes will be periodized, Hypothesis~\ref{hyp:dodging2} then asserts that at each step of the iteration, our algorithm can use at most a finite number of high-frequency pipe segments inside any single periodic cell.
\end{remark}

\subsection{Inductive velocity bounds}
\label{sec:inductive:secondary:velocity}
In this subsection, we present inductive $L^\infty$-bounds for velocity increments and velocity, which are derived from the construction of velocity cutoffs.  All inductive assumptions in subsection~\ref{sec:inductive:secondary:velocity} except for \eqref{est.upsilon.ptwise} at $q\mapsto q+1$ will be verified in Section \ref{sec:inductive.cutoffs}.

Assume that $0\leq q'\leq q+\bn-1$.  First, for $0 \leq i \leq \imax$, $k\geq 1$, $\alpha, \beta \in \N^k$, we assume that
\begin{align}
&\norm{\left( \prod_{l=1}^k D^{\alpha_l} D_{t,q'-1}^{\beta_l} \right) \hat w_{q'} }_{L^\infty(\supp \psi_{i,q'})} \leq \Gamma_{q'}^{i+2}\de_{q'}^{\sfrac12} r_{q'-\bn}^{-\sfrac13} (\la_{q'}\Ga_{q'})^{|\alpha|} \MM{|\beta|,\Nindt, \Gamma_{{q'}}^{i+3}  \tau_{q'-1}^{-1},  \Gamma_{{q'-1}} \Tau_{q'-1}^{-1}}
\label{eq:nasty:D:wq:old}
\end{align}
for $|\alpha|+|\beta|  \leq {\sfrac{3\Nfin}{2}+1}$. We also assume that for $N\geq 0$,
\begin{subequations}\label{eq:nasty}
\begin{align}
&\norm{ D^N \Big( \prod_{l=1}^k D_{q'}^{\alpha_l} D_{t,q'-1}^{\beta_l} \Big) \hat w_{q'}}_{L^\infty(\supp \psi_{i,q'})} \notag\\
&\qquad \leq
(\Gamma_{{q'}}^{i+2}\delta_{q'
}^{\sfrac 12}r_{q'-\bn}^{-\sfrac 13})^{|\alpha|+1} (\la_{q'}\Ga_{q'})^{N+|\alpha|} \MM{|\beta|,\Nindt,  \Gamma_{q'}^{i+3}  \tau_{q'-1}^{-1},  \Gamma_{q'-1} \Tau_{q'-1}^{-1}}   \label{eq:nasty:Dt:wq:old} \\
&\qquad \leq
\Gamma_{q'}^{i+2} \delta_{q'}^{\sfrac 12} r_{q'-\bn}^{-\sfrac 13} (\lambda_{q'}\Ga_{q'})^N (\Gamma_{q'}^{i-5}  \tau_{q'}^{-1})^{|\alpha|}  \MM{|\beta|,\Nindt, \Gamma_{q'}^{i+3}  \tau_{q'-1}^{-1},  \Gamma_{q'-1} \Tau_{q'-1}^{-1}}
\label{eq:nasty:Dt:wq:WEAK:old}
\end{align}
\end{subequations}
whenever $N+|\alpha|+|\beta|\leq  {\sfrac{3\Nfin}{2}+1}$. Next, we assume
\begin{align}
&\norm{\left( \prod_{l=1}^k D^{\alpha_l} D_{t,q'}^{\beta_l} \right) D \hat u_{q'} }_{L^\infty(\supp \psi_{i,q'})}  \leq \tau_{q'}^{-1}\Ga_{q'}^{i-4} (\lambda_{q'}\Ga_{q'})^{|\alpha|} \MM{|\beta|,\Nindt,\Gamma_{q'}^{i-5} \tau_{q'}^{-1},   \Gamma_{q'-1} \Tau_{q'-1}^{-1}}
\label{eq:nasty:D:vq:old}
\end{align}
for $|\alpha|+|\beta| \leq {\sfrac{3\Nfin}{2}}$.  In addition, we assume the lossy bounds
\begin{subequations}\label{eq:bob:old}
\begin{align}
\norm{\left( \prod_{l=1}^k D^{\alpha_l} D_{t,q'}^{\beta_l} \right)  \hat u_{q'}}_{L^\infty(\supp \psi_{i,q'})} &\leq \tau_{q'}^{-1}\Ga_{q'}^{i+2} \la_{q'} (\lambda_{q'}\Ga_{q'})^{|\alpha|} \MM{|\beta|,\Nindt,\Gamma_{q'}^{i-5} \tau_{q'}^{-1},   \Gamma_{q'-1} \Tau_{q'-1}^{-1}}
\label{eq:bob:Dq':old} \\
\left\| D^{|\alpha|} \partial_t^{|\beta|} \hat{u}_{q'} \right\|_{L^\infty} & \leq \Lambda_{q'}^{\sfrac 12} \Lambda_q^{|\alpha|} \Tau_{q'}^{-|\beta|} \, , \label{eq:bobby:old}
\end{align}
\end{subequations}
hold, where the first bounds holds for $|\alpha|+|\beta| \leq \sfrac{3\Nfin}{2}+1$, and the second bound holds for $|\alpha|+|\beta|\leq 2\Nfin$.

\begin{remark}[\bf Upgrading material derivatives]
\label{rem:D:t:q':orangutan}
By applying Lemma~\ref{lem:cooper:1} and \eqref{eq:nasty:Dt:wq:WEAK:old}, we have the bound
\begin{align}
&\norm{ D^N  D_{t,q'}^{M}  \hat w_{q'} }_{L^\infty(\supp \psi_{i,q'})} \les \Gamma_{q'}^{i+2} \delta_{q'}^{\sfrac 12} r_{q'-\bn}^{-\sfrac 13} (\la_{q'}\Ga_{q'})^N 
\MM{M,\Nindt, \Gamma_{q'}^{i-5}  \tau_{q'}^{-1},  \Gamma_{q'-1} \Tau_{q'-1}^{-1}}
\label{eq:nasty:Dt:uq:orangutan}
\end{align}
for all $N+M \leq {\sfrac{3\Nfin}{2}+1}$. Specifically, we set $B = D_{t,q'-1}$ and $A = D_{q'}$, so that $A+B = D_{t,q'}$.  Then the estimate \eqref{eq:nasty:Dt:uq:orangutan} follows from the aforementioned Lemma and \eqref{ineq:tau:q}. We similarly have that \eqref{eq:sharp:Dt:psi:i:q:mixed:old} and \eqref{condi.Nindt} imply that for all $N+M \leq \Nfin$,
\begin{align}
\frac{{\bf 1}_{\supp \psi_{i,q'}}}{\psi_{i,q'}^{1- (N+M)/\Nfin}} \left| D^N  D_{t,q'}^{M}  \psi_{i,q'} \right| &\leq \Gamma_{q'} (\lambda_{q'}\Ga_{q'})^N
\MM{M,\Nindt-\NcutSmall, \Gamma_{q'}^{i-5} \tau_{q'}^{-1}, \Gamma_{q'-1}  \Tau_{q'-1}^{-1}} \notag\\
&\lesssim \Gamma_{q'} (\lambda_{q'}\Ga_{q'})^N \MM{M,\Nindt,\Gamma_{q'}^{i-4}\tau_{q'}^{-1},\Gamma_{q'-1}^2\Tau_{q'-1}^{-1}} \, .
\label{eq:nasty:Dt:psi:i:q:orangutan}
\end{align}
\end{remark}

\opsubsection{Inductive assumptions for the local energy inequality}\label{ss:inductive:LEI}
In this subsection, we record several extra inductive assumptions which are only used in the proof of Theorem~\ref{thm:main}, but not in the proof of Theorem~\ref{thm:main:ER}. All assumptions in this subsection will be verified for $q\mapsto q+1$ in the companion paper~\cite{GKN23}, and we refer to \cite[section~2]{GKN23} for a presentation of these inductive assumptions which is integrated with the rest of the inductive assumptions required for the proof of Theorem~\ref{thm:main}.

\opsubsubsection{Approximate solution}\label{ss:relaxed:LEI}
First, we assume that the approximate solution now includes a scalar field $\varphi_q:\T^3\times[-\tau_{q-1},T+\tau_{q-1}]$, which is called the \emph{current error}\index{current error}.  The current error plays the role of the Reynolds stress in the relaxation of the local energy inequality, given by\index{$\kappa_q$}\index{relaxed local energy identity}
\begin{align}\label{ineq:relaxed.LEI}
\pa_t \left( \frac 12 |u_q|^2 \right)
    + \div\left( \left(\frac 12 |u_q|^2 + p_q\right) u_q\right)
    = (\pa_t + \hat u_q \cdot \na ) \ka_q  
    + \div((R_q-\pi_q\Id) \hat u_q) + \div \ph_q - E(x,t)\, .
\end{align}
We use the notation $\ka_q = \sfrac{\tr(R_q -\pi_q \Id)}2$, and $E(x,t)$ is given continuous function which independent of $q$ and will become the Duchon-Robert measure of the limiting solution. The current error $\varphi_q$ has a decomposition
\begin{align}
\ph_q &= \sum_{k=q}^{q+\bn-1} \ph_q^k \,. \label{eq:LEI:decomp:basic} 
\end{align}
Analogous to $R_q^k$, the portions $\ph_q^k$ of $\ph_q$ have spatial derivative cost $\la_k$ in an effective sense.

\opsubsubsection{Bounds for intermittent pressure $\pi_q^k$ for $k\geq \qbn$}
For $q+\bn\leq k \leq q+\Npr-1$\index{$\Npr$} (where $\Npr$ is defined in subsection~\ref{sec:para.q.ind}, item~\ref{i:par:4.5}) and $N+M\leq 2\Nind$, we assume that $\pi_q^k$ satisfies
\begin{subequations}\label{eq:pressure:inductive:largek}
\begin{align}
\norm{ \psi_{i,q{+\bn-1}} D^N D_{t,q{+\bn-1}}^M  \pi_q^k }_{\sfrac 32}  
 &\leq \Ga_q \Ga_{{k}} \delta_{k+\bn} \Lambda_{q+\bn-1}^N \MM{M, \NindRt, \Gamma_{{q+\bn-1}}^{i} \tau_{{q+\bn-1}}^{-1} ,  \Tau_{q+\bn-1}^{-1} }
\label{eq:pressure:inductive:dtq:largek} \\
\norm{ \psi_{i,q{+\bn-1}} D^N D_{t,q{+\bn-1}}^M \pi_q^k }_{\infty}  
 &\leq \Ga_q \Gamma_{{q+\bn-1}}^{\badshaq+1} \Lambda_{q+\bn-1}^N \MM{M, \NindRt, \Gamma_{q+\bn-1}^{i} \tau_{q+\bn-1}^{-1} ,  \Tau_{q+\bn-1}^{-1} } \, , \label{eq:pressure:inductive:dtq:uniform:largek} \\
 \label{eq:ind:pi:by:pi:largek}
    \left|\psi_{i,q+\bn-1} D^N D_{t,q+\bn-1}^M \pi_q^{k}\right| &\leq {\Gamma_q}\pi_q^k  \Lambda_{q+\bn-1}^N \MM{M, \NindRt, \Gamma_{q+\bn-1}^{i} \tau_{q+\bn-1}^{-1} ,  \Tau_{q+\bn-1}^{-1} } \, .
\end{align}
\end{subequations}
\opsubsubsection{Lower and upper bounds for \texorpdfstring{$\pi_q^k$}{tpdfs4}}
For $ k\geq q$, we assume that $\pi_q^{k}$ has the lower bound
\begin{align}\label{low.bdd.pi}
    \pi_q^{k} \geq \de_{k+\bn} \, .
\end{align}
{For all $q+\bn-1 \leq k' < k \leq q+\Npr-1$, we assume that $\pi_q^k$ has the upper bound
\begin{align}\label{ind:pi:upper}
    \pi_q^k \leq \pi_q^{k'} \, .
\end{align}}
For all $k\geq q+\Npr$, we assume that
\begin{equation}\label{defn:pikq.large.k}
    \pi_q^k \equiv \Ga_k \de_{k+\bn} \, .
\end{equation}
We finally assume that for all $q\leq q' < q''<\infty$,
\begin{subequations}\label{eq:ind.pr.anticipated}
\begin{align}
    \frac{\de_{q''+\bn}}{\de_{q'+\bn}} \pi_q^{q'} &< 2^{q'-q''}\pi_q^{q''} \,, \qquad \text{if}\, q+\half \leq q'' \label{eq:ind.pr.anticipated.1}\\
    \frac{\de_{q''+\bn}}{\de_{q'+\bn}} \pi_q^{q'} &< \pi_q^{q''} \,,\qquad \text{otherwise} \, . \label{eq:ind.pr.anticipated.2}
\end{align}
\end{subequations}
This final bound says that the $\pi_q^k$'s obey a scaling law which may be roughly translated as ``any $\pi_q^{k+m}$ for $m>0$ can be bounded from below by an appropriately rescaled $\pi_q^k$.''

\opsubsubsection{Pointwise bounds for current error}\label{ind:ptbdd:current}
We assume that we have the pointwise estimate
\begin{equation}    \label{eq:ind:current:by:pi}
    \left|\psi_{i,k-1} D^N D_{t,k-1}^M 
    {\ph_q^k}\right|
    < \Ga_q \Ga_k^{-12}(\pi^k_q)^\frac32 r_k^{-1} \Lambda_k^N \MM{M, \NindRt, \Gamma_{k-1}^{i  +20} \tau_{k-1}^{-1} ,  \Tau_{k-1}^{-1} \Ga_{k-1}^{10} }
\end{equation}
for $N+M\leq \sfrac{\Nind}{4}$. 

\opsubsubsection{More dodging hypotheses}\label{ss:dodging:current}
In order to treat several current errors related to the term $(R_q-\pi_q\Id) \hat u_q$ appearing in \eqref{ineq:relaxed.LEI}, we require the following two additional dodging assumptions, which state that certain velocity increments are either disjoint from pressures and stresses, or may be controlled pointwise via already existing intermittent pressure.

\begin{hypothesis*}[\bf Stress dodging]\label{hyp:dodging4}
For all $k,q''$ such that $q\leq q'' \leq k-1$ and $q\leq k \leq q+\bn-1$, we assume that
\begin{equation}\label{eq:ind:stressdodging:equiv}
    B\left(\supp \hat w_{q''} , \lambda_{q''}^{-1}\Gamma_{q''+1} \right) \cap  \supp R_{q}^{k,l}= \emptyset \, .
\end{equation}
\end{hypothesis*}

\begin{hypothesis*}[\bf Pressure dodging]\label{hyp:dodging5}
We assume that for all $q<k\leq q+\bn-1$, $k\leq k'$, and $N+M\leq 2\Nind$,
\begin{subequations}
\begin{align}
\label{eq:pinl:inductive:dtq} 
    \left|\psi_{i,k-1} D^N D_{t,k-1}^M \left(\hat w_{k}\pi_q^{k'}\right)\right| &<
    \Ga_q 
    {\Ga_k^{-100}} \left(\pi_q^k\right)^{\sfrac32} r_k^{-1} \Lambda_k^{N} \MM{M, \NindRt, \Gamma_{k-1}^{i+{1}} \tau_{k-1}^{-1} , \Gamma_{k}^{-1} \Tau_{k}^{-1} } \, .
\end{align}
\end{subequations}
\end{hypothesis*}

\opsubsubsection{Velocity increment potentials}\label{sec:ind.vel.inc.pot}
We assume that for all $q-1 < q' \leq q+\bn -1$ and $\hat w_{q'}$ as in \eqref{eq:hat:no:hat}, there exists a velocity increment potential $\hat\upsilon_{q'}$ and an error $\hat e_{q'}$ such that $\hat w_{q'}$ can be decomposed as
\begin{align}\label{exp.w.q'}
\hat w_{q'} = \div^{\dpot} \hat \upsilon_{q'} +\hat e_{q'}  \,, 
\end{align}
which written component-wise gives $
\hat w_{q'}^\bullet = 
\pa_{i_1}\cdots \pa_{i_\dpot} \hat\upsilon_{q'}^{(\bullet, i_1, \cdots, i_\dpot)} + \hat e_{q'}^\bullet
$. Next, we assume that $\hat\upsilon_{q'}$ and $\hat e_{q'}$ satisfy
\begin{align}\label{supp.upsilon.e.ind}
B\left(\supp(\hat w_{q''}), \frac14 \la_{q''}\Ga_{q''}^2\right) \cap \left(\supp(\hat \upsilon_{q'}) \cup \supp(\hat e_{q'})\right) = \emptyset 
\end{align}
for any $q+1 \leq q'' < q'$.  In addition, we assume that $\hat\upsilon_{q',k}^\bullet :=\la_{q'}^{\dpot-k}\partial_{i_1}\cdots \partial_{i_k} \hat\upsilon_{q'}^{(\bullet,i_1,\dots,i_\dpot)}$, $0\leq k \leq \dpot$, satisfies the estimates
\begin{align}
    &\left|\psi_{i,q'-1} D^ND_{t,q'-1}^{M} 
   \hat \upsilon_{q',k}
    \right|< \Ga_q \Ga_{q'}\left(\pi_q^{q'}\right)^{\sfrac12} r_{q'-\bn}^{-1}
     (\la_{q'}{\Ga_{q'}})^N
     \MM{M, \Nindt, \Ga_{q'-1}^i \tau_{q'-1}^{-1}, \Tau_{q'-1}^{-1}\Ga_{q'-1}^2}
    \label{est.upsilon.ptwise}
\end{align}
for $N+M \leq \sfrac{3\Nfin}{2}$. Finally, we assume that $\hat e_{q'}$ satisfies the estimates
\begin{align} 
    \norm{D^ND_{t,q'-1}^{M} \hat e_{q'}}_{\infty} 
    \leq \de_{q'+2\bn}^3 \Tau_{q'}^{5\Nindt}\la_{q'}^{-10} 
    (\la_{q'} {\Ga_{q'}})^N
     \MM{M, \Nindt, \tau_{q'-1}^{-1}, \Tau_{q'-1}^{-1}\Ga_{q'-1}^2}
    \label{est.e.inf}\, .
\end{align}
for $N+M \leq \sfrac{3\Nfin}{2}$. The velocity increment potential is used in \cite[section~5.3]{GKN23} to help invert the divergence on a product of a velocity increment with stresses and intermittent pressures.

\subsection{Inductive propositions}\label{ss:sunday}
In this section, we first introduce the inductive proposition required for Theorem~\ref{thm:main}, and point out the inductive assumptions for $q\mapsto q+1$ which are verified in this article. The proof of Theorem~\ref{thm:main} is contained in \cite[subsection~2.7]{GKN23}, and \cite[Section~3]{GKN23} includes a discussion of the portion of the proposition which is verified in this article. Next, we present a simplified inductive proposition which is sufficient for flexibility statements analogous to that contained in Theorem~\ref{thm:main:ER}.

\begin{proposition*}[\bf Inductive proposition for Theorem~\ref{thm:main}]\label{prop:main}
Fix $\be\in (0, \sfrac13)$, and choose $\bn$ satisfying \eqref{eq:choice:of:bn}, $b\in (1, \sfrac{25}{24})$ satisfying \eqref{ineq:b}, $T>0$, and a continuous positive function $E(x,t)\geq 0$. Then there exist parameters $\varepsilon_\Gamma$, $\badshaq$, $\Npr$, $\NcutSmall$, $\Nindt$, $\Nind$, $\Nfin$, depending only on $\be$, $b$, and $\bn$ (see section~\ref{sec:para.q.ind} and subsection~\ref{sec:not.general}) such that we can find sufficiently large $a_*=a_*(b,\be, \bn, T)$ such that for $a\geq a_*(b,\be, \bn, T)$, the following statements hold for any $q\geq 0$. Suppose that an approximate solution $(u_q, p_q, R_q, \ph_q, -\pi_q)$ of the Euler-Reynolds system \eqref{eqn:ER} and the relaxed local energy identity \eqref{ineq:relaxed.LEI} with dissipation measure $E$ on the time interval $[-\tau_{q-1},T+\tau_{q-1}]$ is given, and suppose that there exist partitions of unity $\{\psi_{i,q'}^6\}_{i\geq 0}$ of $[-\tau_{q-1},T+\tau_{q-1}] \times \T^3$ for $q-1\leq q'\leq q+\bn-1$ such that 
\begin{itemize}
        \item $\psi_{i,q'}$ satisfies \eqref{eq:inductive:partition}--\eqref{eq:inductive:timescales}, and
        \item the velocity $u_q$ and the errors $R_q$, $\ph_q$, and $\pi_q$ may be decomposed as in \eqref{eq:hat:no:hat}--\eqref{eq:pi:decomp:basic} and~\eqref{eq:LEI:decomp:basic} so that \eqref{eq:pressure:inductive}--\eqref{eq:psi:q:q'}, \eqref{eq:pressure:inductive:largek}--\eqref{eq:ind:current:by:pi}, Hypotheses~\ref{hyp:dodging1}--\ref{hyp:dodging2} and~\ref{hyp:dodging4}--\ref{hyp:dodging5}, \eqref{eq:nasty:D:wq:old}--\eqref{eq:bob:old}, and \eqref{exp.w.q'}--\eqref{est.e.inf} hold. 
    \end{itemize}
     Then there exist a new partition of unity $\{\psi_{i, q+\bn}^6\}_{i\geq 0}$ of $[-\tau_q,T+\tau_q] \times \T^3$ satisfying \eqref{eq:inductive:partition}--\eqref{eq:inductive:timescales} for $q'=q+\bn$, and a new approximate solution $(u_{q+1}, p_{q+1}, R_{q+1}, \ph_{q+1}, -\pi_{q+1})$ satisfying \eqref{eqn:ER} and \eqref{ineq:relaxed.LEI} on $[-\tau_q,T+\tau_q]$ with dissipation measure $E$ and also the following conditions.  The approximate solution may be decomposed as in \eqref{eq:hat:no:hat}--\eqref{eq:pi:decomp:basic} and \eqref{eq:LEI:decomp:basic} for $q\mapsto q+1$ so that \eqref{eq:pressure:inductive}--\eqref{eq:psi:q:q'}, \eqref{eq:pressure:inductive:largek}--\eqref{eq:ind:current:by:pi}, Hypotheses~\ref{hyp:dodging1}--\ref{hyp:dodging2} and~\ref{hyp:dodging4}--\ref{hyp:dodging5}, \eqref{eq:nasty:D:wq:old}--\eqref{eq:bob:old}, and \eqref{exp.w.q'}--\eqref{est.e.inf} hold for $q\mapsto q+1$.
\end{proposition*}
\begin{proof}[Partial proof of Proposition~\ref{prop:main}]

In section~\ref{s.corr-ec-tor}, we construct a new velocity $u_{q+1} = u_q + \hat w_{q+\bn}$, and in section~\ref{opsection:vel:inc:pot}, we construct the associated velocity increment potential.  In section~\ref{sec:ER:main}, we construct a stress error $\overline R_{q+1}$ defined on $\T^3 \times [-\tau_q, T+\tau_q]$.  Finally, in section~\ref{sec:inductive.cutoffs}, we construct a new partition of unity $\{\psi_{i, q+\bn}^6\}_{i\geq 0}$ of $\T^3\times[-\tau_{q},T+\tau_{q}]$.  From the results in the aforementioned sections, the new velocity, stress error, and partition of unity satisfy the following conditions.
\begin{itemize}
\item $\psi_{i,q+\bn}$ satisfies \eqref{eq:inductive:partition}--\eqref{eq:inductive:timescales} for $q'=q+\bn$. 
\item The pair $(u_{q+1}, p_q, \overline R_{q+1}, -(\pi_q-\pi_q^q))$ solves
\begin{align*}
    \pa_t u_{q+1} + \div (u_{q+1}\otimes u_{q+1}) + \na p_q = \div (-(\pi_q-\pi_q^q)\Id + \overline R_{q+1}), \quad \div u_{q+1}=0\,  , 
\end{align*}
analogous to \eqref{eqn:ER}. 
\item The new velocity $u_{q+1}$ can be decomposed as in \eqref{eq:hat:no:hat}, and the stress $\overline R_{q+1}$ can be decomposed as 
\begin{align*}
    \overline R_{q+1} = \sum_{k=q+1}^{q+\bn}\overline R_{q+1}^k, \quad \overline R_{q+1}^k = \overline R_{q+1}^{k,l} + \overline R_{q+1}^{k,*} \, ,
\end{align*}
analogous to \eqref{eq:ER:decomp:basic}.  Furthermore, we have that 
\eqref{eq:nasty:D:wq:old}--\eqref{eq:bob:old}, \eqref{exp.w.q'}, \eqref{supp.upsilon.e.ind}, \eqref{est.e.inf} hold for $q'=q+\bn$, and $R_{q+1}^{k,l}:=\overline R_{q+1}^{k,l}$ verifies Hypothesis \ref{hyp:dodging4} for $q\mapsto q+1$. 

\item Hypotheses~\ref{hyp:dodging1}--\ref{hyp:dodging2} hold, provided that Lemma \ref{lem:dodging} holds true. This lemma will be verified in \cite[section~4]{GKN23}. 
\end{itemize}

For the full proof of this proposition, we refer to \cite{GKN23}.  In particular, \cite[section~3]{GKN23} recalls the set-up of the proof of the inductive proposition and contains a summary of the specific results from this paper which the proof requires. 
\end{proof}

For the purpose of proving Theorem~\ref{thm:main:ER}, it is enough to propagate the following subset of the inductive assumptions.

\begin{proposition}[\bf Inductive proposition for Theorem~\ref{thm:main:ER}]\label{prop:main:ER}
Fix $\be\in (0, \sfrac13)$, and choose $\bn$ satisfying \eqref{eq:choice:of:bn}, $b\in (1, \sfrac{25}{24})$ satisfying \eqref{ineq:b}, and $T>0$. There exist parameters $\varepsilon_\Gamma$, $\badshaq$, $\dpot$, $\Npr$, $\NcutSmall$, $\Nindt$, $\Nind$, $\Nfin$, depending only on $\be$, $b$, and $\bn$ (see section~\ref{sec:para.q.ind} and subsection~\ref{sec:not.general}) such that we can find sufficiently large $a_*=a_*(b,\be, \bn, T)$ such that for $a\geq a_*(b,\be, \bn, T)$, the following statements hold for any $q\geq 0$. Suppose that we have an approximate solution $(u_q, p_q, R_q, -\pi_q)$ which satisfies the Euler-Reynolds system \eqref{eqn:ER} on the time interval $[-\tau_{q-1},T+\tau_{q-1}]$, and suppose there exist partitions of unity $\{\psi_{i,q'}^6\}_{i\geq 0}$ of $\T^3\times[-\tau_{q-1},T+\tau_{q-1}]$ for $q-1\leq q'\leq q+\bn-1$ such that 
    \begin{itemize}
        \item $\psi_{i,q'}$ satisfies \eqref{eq:inductive:partition}--\eqref{eq:inductive:timescales}.
        \item The velocity $u_q$, the error $R_q$, and the intermittent pressure $\pi_q$ may be decomposed as in \eqref{eq:hat:no:hat}--\eqref{eq:pi:decomp:basic} so that
        \eqref{eq:pressure:inductive}--\eqref{eq:psi:q:q'}, Hypotheses~\ref{hyp:dodging1} and~\ref{hyp:dodging2}, and \eqref{eq:nasty:D:wq:old}--\eqref{eq:bob:old} hold.
    \end{itemize}
     Then there exist a new partition of unity $\{\psi_{i, q+\bn}^6\}_{i\geq 0}$ of $\T^3\times[-\tau_{q},T+\tau_{q}]$ satisfying \eqref{eq:inductive:partition}--\eqref{eq:inductive:timescales} for $q'=q+\bn$ and a new approximate solution $(u_{q+1}, p_{q+1}, R_{q+1}, -\pi_{q+1})$ satisfying \eqref{eqn:ER} for $q\mapsto q+1$ on $\T^3\times[-\tau_{q},T+\tau_{q}]$ as well as the following.  The approximate solution may be decomposed as in \eqref{eq:hat:no:hat}--\eqref{eq:pi:decomp:basic} for $q\mapsto q+1$ so that \eqref{eq:pressure:inductive}--\eqref{eq:psi:q:q'}, Hypothesis \ref{hyp:dodging1} and~\ref{hyp:dodging2}, and \eqref{eq:nasty:D:wq:old}--\eqref{eq:bob:old} hold for $q\mapsto q+1$.
\end{proposition}

\begin{proof}[Outline of the proof of Proposition~\ref{prop:main:ER}]
Throughout this proof, we restrict our attention to the Euler-Reynolds system.  The main components of the proof, drawing from the rest of the article, are as follows. 
\begin{itemize}
    \item First, we construct the new premollified velocity increment $w_{q+1}$ in subsection~\ref{ss:corr-ec-tor} by setting $w_{q+1,\ph} =0$, and hence $w_{q+1}= w_{q+1,R}$. In the definition of $w_{q+1,R}$, furthermore, we set $R_{q,i,k} = -\na\Phi_{(i,k)}(R_\ell - \pi_\ell\Id)\na\Phi_{(i,k)}^T$ in~\eqref{eq:rqnpj}. The velocity increment $\hat w_{q+\bn}$ is then defined in~\eqref{def.w.mollified}.
    \item A new partition of unity $\{\psi_{i,q+\bn}^6\}_{i\geq 0}$ is defined
    on $\T^3\times[-\tau_{q},T+\tau_{q}]$
    as in Definition~\ref{def:psi:i:q:def}. Then, under the restricted inductive assumptions listed in Proposition~\ref{prop:main:ER}, \eqref{eq:inductive:partition}--\eqref{eq:inductive:timescales}, \eqref{eq:nasty:D:wq:old}--\eqref{eq:bob:old}, and \eqref{eq:hat:no:hat} for $q\mapsto q+1$ are verified, by the arguments given in section~\ref{sec:inductive.cutoffs}.
    \item Hypotheses~\ref{hyp:dodging1}-- \ref{hyp:dodging2} are verified in~\cite[section~4]{GKN23}, and we refer to the discussion following the statement of Lemma~\ref{lem:dodging} for an outline of the proof.
    \item Referring to Definition~\ref{def:of:new:stresses}, we set $R_{q+1}=\overline R_{q+1}$ and define $R_{q+1}^k$, $R_{q+1}^{k,l}$ and $R_{q+1}^{k,*}$ in a similar fashion. Then by definition, $R_{q+1}$ satisfies the decomposition~\eqref{eq:ER:decomp:basic} at level $q+1$ from \eqref{defn:primitive.stress}--\eqref{defn:local.stress}. We now have from~\eqref{ER:new:equation} that the triple $(u_{q+1}, p_q, \overline R_{q+1}, -(\pi_q-\pi_q^q))$ solves
\begin{align}
    \pa_t u_{q+1} + \div (u_{q+1}\otimes u_{q+1}) + \na p_q = \div (-(\pi_q-\pi_q^q)\Id + \overline R_{q+1}), \quad \div u_{q+1}=0\,  . \label{heatsie:3}
\end{align}
\item Lastly, we define $\pi_{q+1}=\pi_q -\pi_q^q +\si_{q+1}$ and $p_{q+1} = p_q - \si_{q+1}$, where $\si_{q+1}=\sum_{k=q+\half +1}^{q+\bn} \si_{q+1}^k$ and $\si_{q+1}^k$ are defined by
    \begin{align*}
        \si_{q+1}^k = \si_{S^k_O}^+ + \si_{S^k_C}^+ 
        + \mathbf{1}_{m=q+\bn}(\si_{S^k_{TN}}^+ + \si_{\upsilon}^+) + \de_{q+3\bn} \, ,
    \end{align*}
using the pressure increments associated to stress errors which are defined in Section~\ref{sec:ER:main}. Combined with~\eqref{heatsie:3}, this shows that~\eqref{eqn:ER} and \eqref{eq:pi:decomp:basic} are satisfied at level $q+1$.
\item In order to verify~\eqref{eq:pressure:inductive}, we appeal to the definition of $\pi_{q+1}$ above, the inductive assumptions in~\eqref{eq:pressure:inductive} for $\pi_q$, and Lemmas~\ref{lem:oscillation:pressure}, \ref{lem:transport:pressure}, \ref{lem:corrector:pressure}, and~\ref{lem:pr.inc.vel.inc.pot}. In order to verify \eqref{eq:ind:stress:by:pi}, we refer again to Lemmas~\ref{lem:oscillation:pressure}, \ref{lem:transport:pressure}, and~\ref{lem:corrector:pressure}, while for \eqref{eq:ind:velocity:by:pi} we refer to Lemma~\ref{lem:pr.inc.vel.inc.pot}.  The nonlocal estimate in \eqref{eq:Rnnl:inductive:dtq} follows by the same estimate at level $q$, the definition of $R_{q+1}$ above, and Lemmas~\ref{lem:oscillation:general:estimate}, \ref{l:transport:error}, and~\ref{l:divergence:corrector:error}.  Finally, \eqref{eq:psi:q:q'} at level $q+1$ follows from the same estimate at level $q$, the above definition of $\pi_{q+1}$, and Lemma~\ref{lem:pr.vel.dom.cutoff}.
\end{itemize}
\end{proof}

\begin{remark}[\bf Inductive proposition without intermittent pressure]\label{rem:toozday}
It is worth pointing out that for the purpose of proving Theorem~\ref{thm:main:ER}, we do not need to propagate pointwise estimates for $\pi_{q+1}^k$ and $R_{q+1}^k$. As in \cite{NV22}, it actually suffices to remove $\pi_q$ from the inductive assumptions entirely and propagate the $L^p$-estimates given in Remark~\ref{sec:stress:inductive}. Upon doing so, \eqref{eqn:ER} no longer contains $\pi_q$ and \eqref{eq:pi:decomp:basic} and \eqref{eq:pressure:inductive}--\eqref{eq:psi:q:q'} are no longer needed.  Then in order to prove the iterative step, one may proceed as follows.
\begin{itemize}
    \item Define cutoffs for $R_\ell$, analogous for those of $\pi_\ell$ in Definition~\ref{def:pressure:cutoff}, by
    \begin{align*}
        g_{i,q}^2(x,t) &= 1 + \sum_{k=0}^{\NcutLarge} \sum_{m=0}^{\NcutSmall} \delta_{\qbn}^{-2} (\Ga_q \La_q)^{-2k} (\Ga_q^{i}\tau_q^{-1}) \left| D^k \Dtq^m R_\ell(x,t) \right|^2 \, ,\\
        \omega_{i,j,q}(x,t) &= \gamma_{0,q}  \left( \Ga_q^{-2j} g_{i,q}(x,t) \right) \, , \qquad j\geq 1 \, , \\
        \omega_{i,0,q}(x,t) &= \tilde\gamma_{0,q}  \left( \Ga_q^{-2j} g_{i,q}(x,t) \right) \, 
    \end{align*}
    where $\gamma_{0,q}$ and $\tilde\gamma_{0,q}$ are defined as in Lemma~\ref{lem:cutoff:construction:first:statement}.
    This definition is completely analogous to that of \cite[(5.24)--(5.26)]{NV22}. Then following the method of~\cite[section~6.7]{BMNV21}, one can obtain estimates for $\omega_{i,j,q}$ on the support of $\psi_{i,q}$ exactly analogous to those obtained for $\omega_{j,q}$ in subsection~\ref{ss:nic}.  
    \item Define
    \begin{align}\label{heatsie:cheatsie}
        R_{q,j,i,k} = \nabla\Phiik \left( \delta_\qbn \Ga_q^{2j} \Id - R_\ell \right) \nabla \Phiik^T \, ,
    \end{align}
    substituting for the definition of $R_{q,i,k}$ in \eqref{eq:rqnpj}.  Then define the velocity increment exactly as in \eqref{eq:a:xi:def}--\eqref{pi:top:bottom}, except choosing $K=1$ in Proposition~\ref{p:split}.
    \item At this point, no modifications are needed to the rest of the argument - only omissions. Specifically, one may skip sections~\ref{opsection:pressure} and~\ref{opsection:vel:inc:pot}, and simply go through the portions of sections~\ref{sec:ER:main} without asterisks, and all of section~\ref{sec:inductive.cutoffs}. This will suffice to prove a reduced inductive proposition which is sufficient for the construction of weak solutions to Euler which however do not satisfy the local energy inequality.
\end{itemize}
\end{remark}

\begin{remark}[\bf Theorem~\ref{thm:main:ER} and different flavors of flexibility results]
    With the above inductive proposition in hand, the proof of any flexibility result, such as that contained in Theorem~\ref{thm:main:ER}, may be carried out in a manner essentially identical to that of \cite{NV22} or \cite{BMNV21}.  Achieving a decreasing kinetic energy profile will require an inductive assumption measuring the difference between the energy profile of $u_q$ and the desired energy profile.  This can be done in the same manner, for example, as in~\cite{BV19}.  We refer the reader to these references for further details.
\end{remark}

\section{Mollification and upgrading material derivatives}\label{ss:mollification}

In this section, we introduce suitable mollifications of $\pi_q^k$, $R_q^k$, $\ka_q^q$, and $\ph_q^q$ in preparation of later analysis; we have opted to include the mollification of the current error $\ph_q^q$ in this section since the method of proof is identical as for the stress or pressure. The following lemma says that the mollified functions satisfy the same estimates essentially as the unmollified ones, ignoring extra $\Gamma_k$ costs. The difference between the mollified function and the original function, on the other hand, can be made small.

\begin{lemma}[\bf Mollification and upgrading material derivative estimates]\label{lem:upgrading}
Assume that \emph{all} inductive assumptions listed in subsections~\ref{ss:relaxed}-\ref{sec:inductive:secondary:velocity} hold. Let $\Pqxt$ be a space-time mollifier for which the kernel is a product of $\mathcal{P}_{q,x}(x)$, which is compactly supported in space at scale $\Lambda_q^{-1}\Gamma_{q-1}^{-\sfrac 12}$, and $\mathcal{P}_{q,t}(t)$, which is compactly supported in time at scale $\Tau_{q-1}\Gamma_{q-1}^{\sfrac 12}$; we further assume that both kernels have vanishing moments up to $10\Nfin$ and are $C^{10\Nfin}$-differentiable. Define\index{$\Pqxt$}\index{$R_\ell$}\index{$\pi_\ell$}\index{$\varphi_\ell$}
\begin{align}
    &R_\ell = \Pqxt R_q^q \, , \qquad 
    \pi_\ell = \Pqxt \pi_q^q, 
    \label{def:mollified:stuff}
\end{align}
on the space-time domain $[-\sfrac{\tau_{q-1}}2,T+\sfrac{\tau_{q-1}}2]\times \T^3$. 
For $q'$ such that $q<q'\leq q+\bn-1$, we define $\mathcal{P}_{q',x,t}$ in an analogous way after making the appropriate parameter substitutions, and we set $R_\ell^{q'}=\mathcal{P}_{q',x,t}R_q^{q'}$ and $\pi_\ell^{q'}=\mathcal{P}_{q',x,t}\pi_q^{q'}$. For $q'$ with $q+\bn\leq q' <q+\Npr$, we define $\overline{\mathcal{P}}_{q+\bn-1,x,t}$ analogously at the spatial scale $\Lambda_{q+\bn-1}^{-1}\Ga_{q+\bn-1}^{-\sfrac12}$ and temporal scale $\Tau_{q+\bn-1}\Ga_{q+\bn-1}^{-\sfrac12}$ and set $\pi_\ell^{q'} = \overline{\mathcal{P}}_{q+\bn-1,x,t} \pi_q^{q'}$. Then the following hold. 
\begin{enumerate}[(i)]
\item\label{item:moll:one} The following relaxed equation (replacing \eqref{eqn:ER}
) is satisfied:
\begin{align}\label{eqn:ER:LEI:new}
&\partial_t u_q + \div (u_q \otimes u_q) + {\nabla p_q} \nonumber\\
&\qquad = \div\left(R_\ell + \sum\limits_{k=q+1}^{q+\bn-1} R_q^k - \left(\pi_\ell + \sum\limits_{k=q+1}^{q+\Npr-1}\pi_q^k \right) \Id\right) + \div \left( R_q^q - R_\ell + \left( \pi_\ell - \pi_q^q \right) \Id \right) \, . 
\end{align}
\item\label{item:moll:two} The inductive assumptions for $\pi_q^q$ in \eqref{eq:pressure:inductive} are replaced with the following upgraded bounds for $\pi_\ell$ for all $N+M\leq \Nfin$:
\begin{subequations}\label{eq:pressure:upgraded}
\begin{align}
\norm{ \psi_{i,q} D^N D_{t,q}^M  \pi_\ell }_{\sfrac 32}  
&\les \Gamma_q^2 \delta_{q+\bn}\left(\Lambda_q\Gamma_q\right)^N \MM{M, \NindRt, \Gamma_{q}^{i} \tau_q^{-1}, \Tau_{q}^{-1} } \, ,
\label{eq:pressure:inductive:dtq:upgraded} \\
\norm{ \psi_{i,q} D^N D_{t,q}^M \pi_\ell }_{\infty} 
&\les \Gamma_q^{2+\badshaq} \left(\Lambda_q\Gamma_q\right)^N \MM{M, \NindRt, \Gamma_{q}^{i} \tau_q^{-1}, \Tau_{q}^{-1} } \, , \label{eq:pressure:inductive:dtq:uniform:upgraded} \\
\left|\psi_{i,q} D^N D_{t,q}^M \pi_\ell\right| &\les \Gamma_q^3 \pi_\ell \left(\Lambda_q\Gamma_q\right)^N \MM{M, \Nindt, \Gamma_{q}^{i} \tau_q^{-1}, \Tau_{q}^{-1} }  \, . \label{eq:pressure:inductive:dtq:pointwise}
\end{align}
\end{subequations}
While we do not replace the inductive bounds in \eqref{eq:pressure:inductive} and \eqref{eq:pressure:inductive:largek} for $k\neq q$, we do record the following additional bounds for $\pi_\ell^k$ with $q<k\leq q+\bn-1$ and $N+M\leq \Nfin$,
\begin{subequations}\label{eq:pressure:upgraded:higher}
\begin{align}
\norm{ \psi_{i,k-1} D^N D_{t,k-1}^M  \pi_\ell^k }_{\sfrac 32}  
&\les \Ga_k^2 \delta_{k+\bn}\left(\Lambda_k\Gamma_{k-1}\right)^N \MM{M, \NindRt, \Gamma_{k-1}^{i+2} \tau_{k-1}^{-1}, \Tau_{k-1}^{-1}\Ga_{k-1} } \, ,
\label{eq:pressure:inductive:dtq:upgraded:higher} \\
\norm{\psi_{i,k-1} D^N D_{t,k-1}^M \pi_\ell^k }_{\infty} 
&\les \Ga_k^{2+\badshaq} \left(\Lambda_k\Gamma_{k-1}\right)^N \MM{M, \NindRt, \Gamma_{k-1}^{i+2} \tau_{k-1}^{-1}, \Tau_{k-1}^{-1}\Ga_{k-1} } \, , \label{eq:pressure:inductive:dtq:uniform:upgraded:higher} \\
\left|\psi_{i,k-1} D^N D_{t,k-1}^M \pi_\ell^k \right| &\leq 2\Ga_k^3 \pi_\ell^k \left(\Lambda_k\Gamma_k\right)^N \MM{M, \Nindt, \Gamma_{k-1}^{i+3} \tau_{k-1}^{-1}, \Tau_{k-1}^{-1} \Ga_{k-1}^2 }  \, . \label{eq:pressure:inductive:dtq:pointwise:higher}
\end{align}
\end{subequations}
and for $\pi_\ell^k$ with $q+\bn\leq k<q+\Npr$ and $N+M\leq \Nfin$,
\begin{subequations}\label{eq:pressure:upgraded:higher:much}
\begin{align}
\norm{ \psi_{i,q+\bn-1} D^N D_{t,q+\bn-1}^M  \pi_\ell^k }_{\sfrac 32}  
&\les \Ga_k^2 \delta_{k+\bn}\left(\Lambda_{q+\bn-1}\Gamma_{q+\bn-1}\right)^N \notag\\
&\qquad\times\MM{M, \NindRt, \Gamma_{q+\bn-1}^{i+2} \tau_{q+\bn-1}^{-1}, \Tau_{q+\bn-1}^{-1}\Ga_{q+\bn-1} } \, ,
\label{eq:pressure:inductive:dtq:upgraded:higher:much} \\
\norm{ \psi_{i,q+\bn-1} D^N D_{t,q+\bn-1}^M  \pi_\ell^k }_{\infty} 
&\les \Ga_{{ q+\bn-1}}^{2+\badshaq} \left(\Lambda_{q+\bn-1}\Gamma_{q+\bn-1}\right)^N\notag\\
&\qquad\times \MM{M, \NindRt, \Gamma_{q+\bn-1}^{i+2} \tau_{q+\bn-1}^{-1}, \Tau_{q+\bn-1}^{-1}\Ga_{q+\bn-1} } \, , \label{eq:pressure:inductive:dtq:uniform:upgraded:higher:much} \\
\left|\psi_{i,q+\bn-1} D^N D_{t,q+\bn-1}^M  \pi_\ell^k \right| &\leq 2\Ga_k^3 \pi_\ell^k \left(\Lambda_{q+\bn-1}\Gamma_{q+\bn-1}^2\right)^N \notag\\
&\qquad\times\MM{M, \Nindt, \Gamma_{q+\bn-1}^{i+3} \tau_{q+\bn-1}^{-1}, \Tau_{q+\bn-1}^{-1}\Ga_{q+\bn-1}^2 }  \, . \label{eq:pressure:inductive:dtq:pointwise:higher:much}
\end{align}
\end{subequations}
The inductive assumptions \eqref{defn:pikq.large.k} and subsection~\ref{sec:cutoff:inductive} remain unchanged. While we do not discard the estimate in \eqref{low.bdd.pi}, we however record the additional estimate
\begin{equation}\label{ind:pi:lower}
   \frac12 \delta_{q+\bn} \leq \pi_\ell \leq 2\pi_q^q \leq 4\pi_\ell \, , \qquad \frac12 \delta_{k+\bn} \leq \pi_\ell^k \leq 2\pi_q^k \leq 4\pi_\ell^k \, .
\end{equation}

\item\label{item:moll:three} The inductive assumptions in \eqref{eq:ind:stress:by:pi}--\eqref{eq:ind:velocity:by:pi} for $k=q$ are replaced with the following upgraded bounds for all $N+M\leq \Nfin$ in the first two inequalities, and $N+M\leq\sfrac{3\Nfin}{2}$ in the third: 
\begin{subequations}\label{eq:inductive:pointwise:upgraded}
\begin{align}
    \left|\psi_{i,q} D^N D_{t,q}^M R_\ell\right| &\les \Gamma_q^{-7} \pi_\ell \left(\Lambda_q\Gamma_q\right)^N \MM{M, \Nindt, \Gamma_{q}^{i} \tau_q^{-1}, \Tau_{q}^{-1} } \label{eq:inductive:pointwise:upgraded:1} \\
    \left|\psi_{i,q} D^N D_{t,q}^M \hat w_k \right| &\les r_{k-\bn}^{-1} \pi_\ell^{\sfrac 12} \left(\Lambda_q\Gamma_q\right)^N \MM{M, \NindRt, \Gamma_{q}^{i} \tau_q^{-1}, \Tau_{q}^{-1} } \, . \label{eq:inductive:pointwise:upgraded:3}
\end{align}
\end{subequations}
For $k$ such that $q<k\leq q+\bn-1$, we have for $N+M\leq \Nfin$ the additional bound
\begin{align}
     \left|\psi_{i,k-1} D^N D_{t,k-1}^M R_\ell^k\right| &\les \Gamma_q^{-7} \pi_\ell^k \left(\Lambda_k\Gamma_k\right)^N \MM{M, \Nindt, \Gamma_{k-1}^{i+{23}} \tau_{k-1}^{-1}, \Tau_{k-1}^{-1} \Ga_{k-1}^{{12}} } \, . \label{eq:inductive:pointwise:upgraded:1:higher}
\end{align}
\item\label{item:moll:four} The symmetric tensor $R_{\ell}-R_q^q$ and the pressure $\pi_q^q-\pi_\ell$ satisfy
\begin{align}\label{eq:Rcomm:bounds}
     &\left\| D^N \Dtq^M \left( \pi_{\ell}-\pi_q^q \right) \right\|_{\infty} + \left\| D^N \Dtq^M \left( R_{\ell}-R_q^q \right) \right\|_{\infty} \notag\\
     &\qquad \qquad \les \Gamma_{q+1} \Tau_{q+1}^{ 4\Nindt} \delta_{q+3\bn}^2 \lambda_{q+1}^N \MM{M,\Nindt, \tau_{q}^{-1},\Gamma_{q}^{-1}\Tau_{q}^{-1}}
\end{align}
for all $N+M\leq 2\Nind$.  For $k$ such that $q<k\leq q+\bn-1$ and $N+M\leq 2\Nind$, we have that
\begin{align}\label{eq:diff:moll:higher:statement}
    &\norm{D^N D_{t,k-1}^M \left(\pi_q^k - \pi_\ell^k\right)}_\infty +  \norm{D^N D_{t,k-1}^M \left(R_q^k - R_\ell^k\right)}_\infty \notag \\ 
    &\qquad \qquad \lec \Gamma_{k+1} \Tau_{k+1}^{4\Nindt} \delta_{k+3\bn}^2 (\Lambda_k \Ga_{k-1})^N \MM{M,\Nindt,\tau_{k-1}^{-1}\Ga_{k-1},\Tau_{k-1}^{-1}\Gamma_{k-1}^{ 11}} \, . 
\end{align}    
and for $k$ with $q+\bn\leq k <q+\Npr$ and $N+M\leq 2\Nind$, 
\begin{align}
    \norm{D^N D_{t,q+\bn-1}^M \left(\pi_q^k - \pi_\ell^k\right)}_\infty 
    &\lec \Gamma_{q+\bn+1} \Tau_{q+\bn+1}^{{4}\Nindt} \delta_{q+4\bn}^2 (\Lambda_{q+\bn-1} \Ga_{q+\bn-1})^N\notag\\
    &\quad\times
    \MM{M,\Nindt,\tau_{q+\bn-1}^{-1}\Ga_{q+\bn-1},\Tau_{q+\bn-1}^{-1}\Gamma_{q+\bn-1}} \, . \label{eq:diff:moll:higher:statement2}
\end{align}
\end{enumerate}
\end{lemma}

\begin{lemma*}[\bf Mollification and upgrading material derivative estimates]\label{lem:upgrading:current}
Assume that \emph{all} inductive assumptions listed in subsections~~\ref{ss:relaxed}--\ref{ss:inductive:LEI} hold. Let $\Pqxt$ and $\mathcal{P}_{q',x,t}$ be defined as in Lemma \ref{lem:upgrading}. 
Define
\begin{align}
    \varphi_\ell = \Pqxt \varphi_q^q \,
    \label{def:mollified:stuff:current}
\end{align}
on the space-time domain $[-\sfrac{\tau_{q-1}}2,T+\sfrac{\tau_{q-1}}2]\times \T^3$. For $q'$ such that $q<q'\leq q+\bn-1$, set $\ph_\ell^{q'}=\mathcal{P}_{q',x,t}\ph_q^{q'}$.
\begin{enumerate}[(i)]
\item\label{item:moll:current:one} The following relaxed equation (replacing 
\eqref{ineq:relaxed.LEI}) is satisfied:
\begin{align}\label{eqn:ER:LEI:new:current}
&\pa_t \left( \frac 12 |u_q|^2 \right)
+ \div\left( \left(\frac 12 |u_q|^2 + {p_q}\right) u_q\right) \nonumber\\
&\qquad = (\pa_t + \hat u_q \cdot \na ) \ka_q
+ \div\left(\left(R_\ell + \sum\limits_{k=q+1}^{q+\bn-1} R_q^k - \left(\pi_\ell + \sum\limits_{k=q+1}^{q+\Npr-1}\pi_q^k \right) \Id \right) \hat u_q \right) \nonumber\\
&\qquad\,\,   + \,\div\left( \left( R_q^q - R_\ell + (\pi_\ell - \pi_q^q) \Id \right) \hat u_q \right)  + \div \left( \varphi_\ell + \sum\limits_{k=q+1}^{q+\bn-1} \varphi_q^k \right) + \div \left( \varphi_q^q - \varphi_\ell \right)  -E(t)\, .
\end{align}
\item\label{item:moll:current:three} The inductive assumptions in \eqref{eq:ind:current:by:pi} for $k=q$ are replaced with the following upgraded bounds for all $N+M\leq \Nfin$: 
\begin{subequations}\label{eq:inductive:pointwise:upgraded:current}
\begin{align}
    \left|\psi_{i,q} D^N D_{t,q}^M \varphi_\ell \right| &\les \Gamma_q^{-11} \pi_\ell^{\sfrac 32} r_q^{-1} \left(\Lambda_q\Gamma_q\right)^N \MM{M, \Nindt, \Gamma_{q}^{i} \tau_q^{-1}, \Tau_{q}^{-1} } \label{eq:inductive:pointwise:upgraded:2} \, . 
\end{align}
\end{subequations}
The difference $\varphi_{\ell}-\vp_q^q$ satisfies
\begin{align}\label{eq:phicomm:bounds}
    \left\| D^N \Dtq^M \left( \varphi_{\ell}-\vp_q^q \right) \right\|_{\infty} \leq \delta_{q+3\bn}^{\sfrac 32} \lambda_{q+1}^N \MM{M,\Nindt,\tau_{q}^{-1},\Gamma_{q}^{-1}\Tau_{q}^{-1}}
\end{align}
for all $N+M\leq \sfrac{\Nind}{4}$. 
\end{enumerate}
\end{lemma*}

\begin{proof}[Proof of Lemmas \ref{lem:upgrading} and \ref{lem:upgrading:current}]
We first note that \eqref{eqn:ER:LEI:new} and \eqref{eqn:ER:LEI:new:current} are immediate from \eqref{eqn:ER}, \eqref{ineq:relaxed.LEI} and the definitions in \eqref{def:mollified:stuff}.  At this point, we split the proof into steps, in which we first carry out the mollifications, and then upgrade the material derivatives. 
\smallskip

\noindent\texttt{Step 1: Mollifying the pressure $\pi_q^k$.} We first consider the case $k=q$ and apply the abstract mollification Proposition~\ref{lem:mollification:general} with the following choices:
\begin{align*}
    &p=\sfrac 32, \infty \, , \quad N_{\rm g}, N_{\rm c} \textnormal{ as in \eqref{i:par:12}} \, , \quad M_t = \Nindt\,  , \quad N_* = 2\Nind \, , \\
    &N_\gamma = \Nfin \, , \quad  \Omega = \supp \psi_{i,q-1}\, , \quad v = \hat u_{q-1}\, , \quad i=i \, , \\
    &\lambda = \Lambda_{q-1}\, , \quad \Lambda = \Lambda_q \Ga_{q-1}\, , \quad \Ga = \Ga_{q-1}\, , \quad \tau = \tau_{q-1}\Gamma_{q-1} \, , \quad \Tau = \Tau_{q-1},\\
    &f = \pi_q^q, \quad \const_{f,\sfrac32} = \Ga_q^2 \de_{q+\bn}, \quad \const_{f,\infty} = \tilde \const_f = \Ga_q^{\badshaq + 2}\, \quad \const_v = \Lambda_{q-1}^{\sfrac 12} \,  .
\end{align*}
First, we have that the assumptions on the parameters in \eqref{eq:moll:assumps:1:1} are satisfied by \eqref{eq:darnit:1}, \eqref{eq:Nindsy:1},\eqref{condi.Nfin0}, \eqref{v:global:par:ineq} and \eqref{eq:imax:old}. The assumptions in \eqref{eq:moll:assumps:1:3} are satisfied from \eqref{eq:darnit:3}, and the assumptions in \eqref{moll.assum.v.est} are satisfied from \eqref{eq:bobby:old}. Next, the assumptions in \eqref{eq:moll:f:1} are satisfied from \eqref{eq:pressure:inductive} (where we apply the bound with $\psi_{i\pm,q-1}$ in order to obtain a bound for $L^p(\supp \psi_{i,q-1})$).  Finally, in order to verify \eqref{eq:moll:f:2}, we apply Remark~\ref{rem:upgrade.material.derivative.end} with the following choices. We set $p=\infty$, $N_x=N_t=\infty$, $N_*=2\Nind$, $\Omega=\T^3\times\R$, $v=-w=\hat u_{q-1}$, $\const_w = \Gamma_{q-1}^{\imax+2}\delta_{q-1}^{\sfrac 12}\lambda_{q-1}^2$, $\lambda_w=\tilde\lambda_w=\Lambda_{q-1}$, $\mu_w=\tilde \mu_w = \Gamma_{q-1}^{-1}\Tau_{q-1}^{-1}$ in \eqref{eq:cooper:w}, while in \eqref{eq:cooper:2:v} and \eqref{eq:cooper:2:f} we set $v=\hat u_{q-1}$, $\const_v=\const_w$, $\lambda_v=\tilde\lambda_v=\Lambda_{q-1}$, $\mu_v=\tilde\mu_v = \Gamma_{q-1}^{-1}\Tau_{q-1}^{-1}$, $f=\pi_q^q$, $\const_f = \Gamma_q^{\badshaq+2}$, $\lambda_f=\tilde\lambda_f=\Lambda_q$, $\mu_f=\tilde \mu_f= \Tau_{q-1}^{-1}$.  Then \eqref{eq:cooper:2:v} and \eqref{eq:cooper:2:f} are satisfied from \eqref{eq:nasty:D:vq:old} at level $q-1$, \eqref{eq:pressure:inductive}, \eqref{eq:imax:old}, and \eqref{v:global:par:ineq}. Next, \eqref{eq:cooper:w} is satisfied from \eqref{eq:bob:Dq':old} at level $q-1$. Thus from \eqref{eq:cooper:f:mat} and \eqref{v:global:par:ineq}, we obtain that 
\begin{align}\label{eq:pi:f'ed:up}
\left\| D^N \partial_t^M \pi_q^q \right\|_\infty \lesssim \Gamma_q^{\badshaq+2} \Lambda_q^N \Tau_{q-1}^{-M}
\end{align}
for $N+M \leq 2\Nind$, thus verifying the final assumption \eqref{eq:moll:f:2} from Lemma~\ref{lem:mollification:general}.

We first apply \eqref{eq:moll:conc:1} to conclude that for $N+M\leq \Nfin$,
\begin{subequations}\label{eq:pressure:upgraded:copied}
\begin{align}
\norm{ \psi_{i,q-1} D^N D_{t,q-1}^M  \pi_\ell }_{\sfrac 32}  
&\les \Gamma_q^2 \delta_{q+\bn}\left(\Lambda_q\Gamma_{q-1}\right)^N \MM{M, \NindRt, \Gamma_{q-1}^{i+2} \tau_{q-1}^{-1}, \Tau_{q-1}^{-1}\Gamma_{q-1} }
\label{eq:pressure:inductive:dtq-1:upgraded} \\
\norm{ \psi_{i,q-1} D^N D_{t,q-1}^M \pi_\ell }_{\infty} 
&\les \Gamma_q^{\badshaq+2} \left(\Lambda_q\Gamma_{q-1}\right)^N \MM{M, \NindRt, \Gamma_{q-1}^{i+2} \tau_{q-1}^{-1}, \Tau_{q-1}^{-1}\Gamma_{q-1} } \, . \label{eq:pressure:inductive:dtq-1:uniform:upgraded}
\end{align}
\end{subequations}
Next, we have from \eqref{eq:moll:conc:2} and \eqref{eq:darnit:2} that the difference $\pi_q^q - \pi_\ell$ satisfies
\begin{align}\label{eq:diff:moll}
    \norm{D^N D_{t,q-1}^M \left(\pi_q^q - \pi_\ell\right)}_\infty \lec \Gamma_{q+1} \Tau_{q+1}^{4\Nindt} \delta_{q+3\bn}^2 (\Lambda_q \Ga_{q-1})^N \MM{M,\Nindt,\tau_{q-1}^{-1}\Ga_{q-1},\Tau_{q-1}^{-1}\Gamma_{q-1}}
\end{align}
for $N+M\leq {2\Nind}$. Note also that since we have a lower bound on $\pi_q^q$ given by \eqref{low.bdd.pi}, the above estimate implies that (after a sufficiently large choice of $\lambda_0$ so that the implicit constant is absorbed)
\begin{align}
    \pi_\ell \geq \pi_q^q - \de_{q+2\bn} \geq \frac12 \de_{q+\bn} \, , \notag 
\end{align}
which is the first inequality for $\pi_\ell$ and $\pi_q^q$ in \eqref{ind:pi:lower}. The other two inequalities there follow similarly.  Finally, we note that by \eqref{eq:ind:pi:by:pi} and \eqref{ind:pi:lower},
\begin{align}
    \left| \psi_{i,q-1} D^N D_{t,q-1}^M \pi_\ell \right| &\leq \left| \psi_{i,q-1} D^N D_{t,q-1}^M \pi_q^q \right| + \left| D^N D_{t,q-1}^M \left(\pi_q^q - \pi_\ell \right) \right| \notag\\
    &\leq \Gamma_q^2 \pi_q^q \Lambda_{q}^N \MM{M, \NindRt, \Gamma_{q-1}^{i} \tau_{q-1}^{-1} ,  \Tau_{q-1}^{-1} } \notag\\
    &\qquad\qquad + \de^2_{q+3\bn} (\Lambda_q \Ga_{q-1})^N \MM{M,\Nindt,\tau_{q-1}^{-1},\Tau_{q-1}^{-1}\Gamma_{q-1}} \notag\\
    &\leq \Gamma_q^3 \pi_\ell (\Lambda_q \Ga_{q-1})^N \MM{M, \NindRt, \Gamma_{q-1}^{i} \tau_{q-1}^{-1} ,  \Tau_{q-1}^{-1} \Ga_{q-1}} \notag 
\end{align}
for $N+M\leq 2\Nind$. For $2\Nind < N+M \leq \Nfin$, we have from \eqref{eq:pressure:inductive:dtq-1:uniform:upgraded} and \eqref{eq:Nind:darnit} that 
\begin{align}
    \left| D^N D_{t,q-1}^M \pi_\ell \right| \leq \delta_{q+\bn}^2 (\Lambda_q \Ga_{q-1}^{\sfrac 12}\Ga_q^{\sfrac 12})^N \MM{M, \NindRt, \Gamma_{q-1}^{i+3} \tau_{q-1}^{-1} ,  \Tau_{q-1}^{-1} \Ga_{q-1}^2} \, . \notag 
\end{align}

In the case $k\neq q$, we may obtain the bounds \eqref{eq:pressure:inductive:dtq:upgraded:higher}, \eqref{eq:pressure:inductive:dtq:uniform:upgraded:higher}, \eqref{eq:pressure:inductive:dtq:upgraded:higher:much},  \eqref{eq:pressure:inductive:dtq:uniform:upgraded:higher:much}, and the second inequality of \eqref{ind:pi:lower},
via an argument identical to the proof of \eqref{eq:pressure:upgraded} and the first inequality of \eqref{ind:pi:lower}. We additionally have the pointwise bound for $q+1\leq k\leq q+\bn-1$ and $N+M\leq \Nfin$
\begin{align}
     \left|\psi_{i,k-1} D^N D_{t,k-1}^M \pi_\ell^k \right| &\leq (\Ga_k^3 \pi_q^k + \delta_{k+\bn}^2) (\Lambda_k \Ga_{k-1}^{\sfrac 12}\Ga_k^{\sfrac 12})^N \MM{M, \NindRt, \Gamma_{k-1}^{i+3} \tau_{k-1}^{-1} ,  \Tau_{k-1}^{-1} \Ga_{k-1}^2} \notag\\
     &\leq 2\Ga_k^3 \pi_\ell^k  (\Lambda_k \Ga_{k-1}^{\sfrac 12}\Ga_k^{\sfrac 12})^N \MM{M, \NindRt, \Gamma_{k-1}^{i+3} \tau_{k-1}^{-1} ,  \Tau_{k-1}^{-1} \Ga_{k-1}^2} \, , \, 
\end{align}
and for $q+\bn\leq k <q+\Npr$ and $N+M\leq \Nfin$
\begin{align}
     \left|\psi_{i,q+\bn-1} D^N D_{t,q+\bn-1}^M \pi_\ell^k \right| &\leq (\Ga_k^3 \pi_q^k + \delta_{k+\bn}^2) (\Lambda_{q+\bn-1} \Ga_{q+\bn-1}^2)^N \MM{M, \NindRt, \Gamma_{q+\bn-1}^{i+3} \tau_{q+\bn-1}^{-1} ,  \Tau_{q+\bn-1}^{-1} \Ga_{q+\bn-1}^2} \notag\\
     &\leq 2\Ga_k^3 \pi_\ell^k  (\Lambda_{q+\bn-1} \Ga_{q+\bn-1}^2)^N \MM{M, \NindRt, \Gamma_{q+\bn-1}^{i+3} \tau_{q+\bn-1}^{-1} ,  \Tau_{q+\bn-1}^{-1} \Ga_{q+\bn-1}^2} \, ,
\end{align}
which again follows from a similar argument as in the proof of the corresponding bounds for $q=k$ and \eqref{ind:pi:lower}. Furthermore, we have that the difference $\pi_q^k - \pi_\ell^k$ satisfies \eqref{eq:diff:moll:higher:statement} and \eqref{eq:diff:moll:higher:statement2}, 
which follows directly from the mollification lemma and \eqref{eq:darnit:2} with $q$ replaced by $k-1$ or $q+\bn$, as in the case $k=q$. Finally, the bounds in \eqref{ind:pi:lower} for $\pi_\ell^m$ follow similarly as before. At this point we have completed the proofs of the required estimates in \eqref{eq:pressure:upgraded:higher}--\eqref{ind:pi:lower} and \eqref{eq:diff:moll:higher:statement}--\eqref{eq:diff:moll:higher:statement2} for $\pi_\ell^k$.
\smallskip

\noindent\texttt{Step 2: Mollifying the stress and current errors.} We apply the abstract mollification Proposition~\ref{lem:mollification:general} with the same choices as before, except for the stress error we choose
\begin{align*}
    &f = R_q^k \,, \quad  q\leq k \leq q+\bn-1 \, , \quad p=\infty \,, \quad \const_{f,\infty} = \Ga_k^{\badshaq + 2} \, , \quad
    \tau = \tau_{k-1} \, , \quad
    c=20 \, , 
    \quad
    \Tau = \Tau_{k-1}\Ga_q^{-10}
    \,.
\end{align*}
We then have that \eqref{eq:moll:assumps:1:1}--\eqref{eq:moll:assumps:1:3} are satisfied as in the previous step, as is \eqref{moll.assum.v.est}. In order to verify \eqref{eq:moll:f:1}, we appeal to \eqref{eq:ind:stress:by:pi} and \eqref{eq:pressure:inductive:dtq:uniform}. In order to verify \eqref{eq:moll:f:2}, we use Remark~\ref{rem:upgrade.material.derivative.end} exactly as in the previous step, but with $R_q^k$ replacing $\pi_q^k$.  Thus from \eqref{eq:moll:conc:1}--\eqref{eq:moll:conc:2} and \eqref{eq:darnit:2}, we have that for $q\leq k\leq q+\bn-1$ (we denote $R_\ell$ by $R_\ell^q$ for concision here)
\begin{subequations}
\begin{align}
    \left| \psi_{i,k-1} D^N D_{t,k-1} R_\ell^{k} \right| &\lesssim \Gamma_k^{\badshaq+2} (\Lambda_k \Gamma_{k-1})^N \MM{M,\Nindt, \Ga_{k-1}^{i+{22}}\tau_{k-1}^{-1},\Tau_{k-1}^{-1}\Gamma_{k-1}^{{11}}} \label{eq:R:ell:new:one}\\
    \left| D^N D_{t,k-1}^M \left( R_\ell^k - R_q^k \right) \right| &\lesssim \Gamma_{k+1} \Tau_{k+1}^{4\Nindt} \delta_{k+3\bn}^2 (\Lambda_k \Ga_{k-1})^N \MM{M,\Nindt,\tau_{k-1}^{-1},\Tau_{k-1}^{-1}\Gamma_{k-1}^{11}} \, , \label{eq:R:ell:new:two}
\end{align}
\end{subequations}
where the first bound holds for $N+M\leq \Nfin$, and the second bound holds for $N+M\leq 2\Nind$. The second bound verifies \eqref{eq:diff:moll:higher:statement} for the difference $R_q^k - R_\ell^k$. Appealing to \eqref{eq:ind:stress:by:pi}, \eqref{eq:R:ell:new:two}, and \eqref{ind:pi:lower}, we then may write that in the case $k=q$,
\begin{align}
    \left| \psi_{i,q-1} D^N D_{t,q-1}^M R_\ell \right| &\leq \left| \psi_{i,q-1} D^N D_{t,q-1}^M R_q^q \right| + \left| D^N D_{t,q-1}^M \left(R_q^q - R_\ell \right) \right| \notag\\
    &\leq \Gamma_q^{-7} \pi_q^q \Lambda_{q}^N \MM{M, \NindRt, \Gamma_{q-1}^{i+20} \tau_{q-1}^{-1} ,  \Tau_{q-1}^{-1}\Ga_{q}^{11} } \notag\\
    &\qquad\qquad + \de_{q+3\bn}^2 (\Lambda_q \Ga_{q-1})^N \MM{M,\Nindt,\tau_{q-1}^{-1},\Tau_{q-1}^{-1}\Gamma_{q-1}^{  {11}}} \notag\\
    &\lesssim \Gamma_q^{-7} \pi_\ell (\Lambda_q \Ga_{q-1})^N \MM{M, \NindRt, \Gamma_{q-1}^{i} \tau_{q-1}^{-1} ,  \Tau_{q-1}^{-1} \Ga_{q-1}^{11}} \notag 
\end{align}
for $N+M\leq 2\Nind$. For $2\Nind < N+M \leq \Nfin$, we have from \eqref{eq:R:ell:new:one} and \eqref{eq:Nind:darnit} that 
\begin{align}
    \left| D^N D_{t,q-1}^M R_\ell \right| \leq \delta_{q+\bn}^2 (\Lambda_q \Ga_{q-1}^{\sfrac 12}\Ga_q^{\sfrac 12})^N \MM{M, \NindRt, \Gamma_{q-1}^{i+ {23}} \tau_{q-1}^{-1} ,  \Tau_{q-1}^{-1} \Ga_{q-1}^{ {12}}} \, . \notag
\end{align}
In the case $q\neq k$, we have that for $N+M\leq \Nfin$,
\begin{align}
     \left|\psi_{i,k-1} D^N D_{t,k-1}^M R_\ell^k \right| &\les (\Ga_k^{-7} \pi_\ell^k + \delta_{k+\bn}^2) (\Lambda_k \Ga_{k-1}^{\sfrac 12}\Ga_k^{\sfrac 12})^N \MM{M, \NindRt, \Gamma_{k-1}^{i+ {23}} \tau_{k-1}^{-1} ,  \Tau_{k-1}^{-1} \Ga_{k-1}^{  {12}}} \, , \notag 
\end{align}
giving the desired bound in \eqref{eq:inductive:pointwise:upgraded:1:higher} after using \eqref{eq:ind.pr.anticipated.1} again.

In the case of the current error, we again apply Proposition~\ref{lem:mollification:general} with the same choices as in the first portion of this step, except we choose
\begin{align*}
    &f = \ph_q^q \,, \quad \const_{f,\infty} = \Ga_q^{\frac{3\badshaq}{2} + 3} r_q^{-1}\, \quad
    c=20\,, \quad 
    \Tau= \Tau_{q-1}\Ga_q^{10}\, ,\quad 
    N_*=\sfrac{\Nind}{4}.
\end{align*}
We then have that \eqref{eq:moll:assumps:1:1}--\eqref{eq:moll:assumps:1:3} are satisfied exactly as in the previous step, as is \eqref{moll.assum.v.est}. In order to verify \eqref{eq:moll:f:1}, we appeal to \eqref{eq:ind:current:by:pi} and \eqref{eq:pressure:inductive:dtq:uniform}. In order to verify \eqref{eq:moll:f:2}, we use Remark~\ref{rem:upgrade.material.derivative.end} exactly as in the first part of this step, but with $\ph_q^q$ replacing $R_q^q$. We conclude that \eqref{eq:moll:f:2} is satisfied with $\tilde \const_f = \const_{f,\infty}$. Thus from \eqref{eq:moll:conc:1}--\eqref{eq:moll:conc:2}, we have that
\begin{subequations}
\begin{align}
    \left| \psi_{i,q-1} D^N \Dtqminus \ph_\ell \right| &\lesssim\Ga_q^{\frac{3\badshaq}{2} + 3} r_q^{-1} (\Lambda_q \Gamma_{q-1})^N \MM{M,\Nindt, \Ga_{q-1}^{i+ {22}}\tau_{q-1}^{-1},\Tau_{q-1}^{-1}\Gamma_{q-1}^{  {11}}} \label{eq:ph:ell:new:one}\\
    \left| D^N \Dtqminus \left( \ph_\ell - \ph_q^q \right) \right| &\lesssim \Gamma_{q+1} \Tau_{q+1}^{4\Nindt} \delta_{q+3\bn}^2 (\Lambda_q \Ga_{q-1})^N \MM{M,\Nindt,\tau_{q-1}^{-1},\Tau_{q-1}^{-1}\Gamma_{q-1}^{ {11}}} \, , \label{eq:ph:ell:new:two}
\end{align}
\end{subequations}
where the first bound holds for $N+M\leq \Nfin$, and the second bound holds for $N+M\leq \sfrac{\Nind}{4}$. Appealing to \eqref{eq:ind:current:by:pi}, \eqref{eq:ph:ell:new:two}, and \eqref{ind:pi:lower}, we then may write that
\begin{align}
    \left| \psi_{i,q-1} D^N D_{t,q-1}^M \ph_\ell \right| &\leq \left| \psi_{i,q-1} D^N D_{t,q-1}^M \ph_q^q \right| + \left| D^N D_{t,q-1}^M \left(\ph_q^q - \ph_\ell \right) \right| \notag\\
    &\leq \Gamma_q^{-11} (\pi_q^q)^{\sfrac 32} r_q^{-1} \Lambda_{q}^N \MM{M, \NindRt, \Gamma_{q-1}^{i+20} \tau_{q-1}^{-1} ,  \Tau_{q-1}^{-1}\Gamma_{q}^{10} } \notag\\
    &\qquad\qquad + \de_{q+2\bn}^2 (\Lambda_q \Ga_{q-1})^N \MM{M,\Nindt,\tau_{q-1}^{-1},\Tau_{q-1}^{-1}\Gamma_{q-1}^{11}} \notag\\
    &\lesssim \Gamma_q^{-11} \pi_\ell^{\sfrac 32} r_q^{-1} (\Lambda_q \Ga_{q-1})^N \MM{M, \NindRt, \Gamma_{q-1}^{i+20} \tau_{q-1}^{-1} ,  \Tau_{q-1}^{-1} \Ga_{q-1}^{11}} \notag 
\end{align}
for $N+M\leq \sfrac{\Nind}{4}$. For $\sfrac{\Nind}{4} < N+M \leq \Nfin$, we have from \eqref{eq:ph:ell:new:one} and \eqref{eq:Nind:darnit} that 
\begin{align}
    \left| D^N D_{t,q-1}^M \ph_\ell \right| \leq \delta_{q+\bn}^2 (\Lambda_q \Ga_{q-1}^{\sfrac 12}\Ga_q^{\sfrac 12})^N \MM{M, \NindRt, \Gamma_{q-1}^{i+23} \tau_{q-1}^{-1} ,  \Tau_{q-1}^{-1} \Ga_{q-1}^{12}} \, . \notag
\end{align}
\smallskip

\noindent\texttt{Step 3: Upgrading material derivatives for $k=q$.}
We begin with the pointwise bounds for $\pi_\ell$. Combining the bounds from Step 1 with \eqref{eq:inductive:timescales} with $q'=q$ and $q''=q-1$, we have that for $N+M\leq \Nfin$,
\begin{align}\label{eq:moll:proving:one}
    \left| \psi_{i,q} D^N \Dtqminus^M \pi_\ell \right| \leq 2\Gamma_q^3 \pi_\ell \left( \Lambda_q\Ga_{q-1}^{\sfrac 12}\Ga_q^{\sfrac 12}\right)^N \MM{M,\Nindt, \tau_q^{-1}\Gamma_q^{i-2},\Tau_{q-1}^{-1} \Gamma_{q-1}^2} \, .
\end{align}
We shall apply Remark~\ref{rem:upgrade.material.derivative.end} (with the adjustment in Remark~\ref{rem:cooper:2:sum} for derivative bounds) with the following choices, at a point $(t,x)\in \textnormal{int}\left( \supp \psi_{i,q} \right)$ for which the neighborhood $\Omega_{t,x} \subset \supp \psi_{i,q}$:
\begin{align}
    \eqref{eq:cooper:w}\textnormal{ choices: } &p=\infty \, , \quad N_x = \infty \, , \quad N_t = \Nindt \, , \quad N_* = \Nfin \, , \quad w=\hat w_q \, ,\notag \\
    &\Omega = \Omega_{t,x} \, , \quad v = \hat u_{q-1} \, , \quad \const_w = \Gamma_q^{i+2} \delta_q^{\sfrac 12}r_{q-\bn}^{-\sfrac 13} \, , \notag\\
    & \lambda_w = \tilde \lambda_w = \Lambda_q \, , \quad \mu_w = \Gamma_{q-1}^{i+3}\tau_{q-1}^{-1} \, , \quad \tilde \mu_w = \Ga_q^{-1} \Tau_q^{-1} \, , \notag \\
    \eqref{eq:cooper:2:v}\textnormal{ choices: } &\const_v = \Gamma_q^{i+2} \delta_q^{\sfrac 12} r_{q-\bn}^{-\sfrac 13} \, , \quad \lambda_v = \tilde\la_v = \Lambda_q \, , \quad \mu_v = \Gamma_q^i \tau_q^{-1} \, , \quad \tilde\mu_v = \Tau_q^{-1} \Ga_q^{-1} \, , \quad \Omega=\Omega_{t,x} \, ,  \notag \\
    \eqref{eq:cooper:2:f}\textnormal{ choices: } &f=\pi_\ell \, , \quad \const_f = \displaystyle\sup_{\Omega_{t,x}} \pi_\ell \, , \quad \lambda_f = \tilde \lambda_f = \Lambda_q (\Ga_{q-1}\Ga_q)^{\sfrac 12} \, , \quad \mu_f = \mu_v \, , \quad \tilde \mu_f = \tilde \mu_v \, , \quad \Omega=\Omega_{t,x} \, . \notag
\end{align}
Then we have that \eqref{eq:cooper:w} holds from \eqref{eq:nasty:D:wq:old} at level $q$, \eqref{eq:cooper:2:v} holds from \eqref{eq:nasty:D:vq:old} at level $q$, and \eqref{eq:cooper:2:f} holds from \eqref{eq:moll:proving:one}. Taking $\Omega_{t,x}$ to be arbitrary and using the continuity of $\pi_\ell$, we thus have from \eqref{eq:cooper:f:mat} that for $N+M\leq \Nfin$,
\begin{align}\notag
    \left| \psi_{i,q} D^N \Dtq^M \pi_\ell \right| \lesssim \Gamma_q^3 \pi_\ell \left( \La_q (\Ga_{q-1}\Ga_q)^{\sfrac 12} \right)^N \MM{M,\Nindt, \tau_q^{-1}\Ga_q^i,\Tau_q^{-1}\Ga_q^{-1}} \, ,
\end{align}
matching \eqref{eq:pressure:inductive:dtq:pointwise}. In order to obtain \eqref{eq:pressure:inductive:dtq:upgraded} and \eqref{eq:pressure:inductive:dtq:uniform:upgraded}, we use the $L^{\sfrac 32}$ and $L^\infty$ bounds on $\pi_\ell$ shown in \eqref{eq:pressure:upgraded}. Combined with Step 1, this concludes the proof of \eqref{item:moll:two}.

In order to prove \eqref{eq:inductive:pointwise:upgraded:1}, we argue in a manner very similar to the proof of \eqref{eq:pressure:inductive:dtq:pointwise} carried out just previously. The only difference is that from Step 2, we have the bound
\begin{align}
    \left| D^N \Dtqminus R_\ell \right| \les \Gamma_q^{-7} \pi_\ell \left(\Lambda_q(\Ga_{q-1}\Ga_q)^{\sfrac 12} \right)^N \MM{M,\Nindt, \Gamma_{q-1}^{i+23}\tau_{q-1}^{-1},\Tau_{q-1}^{-1}\Gamma_{q-1}^{12}} \, .
\end{align}
Carrying out the same steps with the obvious modifications, we deduce that \eqref{eq:inductive:pointwise:upgraded:1} holds as desired.  The proof of \eqref{eq:inductive:pointwise:upgraded:2} is again quite similar, and we omit the details.  To conclude the proof of \eqref{item:moll:three}, we must show \eqref{eq:inductive:pointwise:upgraded:3}. Following the exact same steps as before but beginning instead with the bound \eqref{eq:ind:velocity:by:pi} and appealing to \eqref{ind:pi:lower}, we obtain the desired estimate, concluding the proof of item~\eqref{item:moll:three}.

Finally, we must upgrade the material derivatives to $\Dtq$ on the differences in order to conclude the proofs of \eqref{eq:Rcomm:bounds}--\eqref{eq:phicomm:bounds} from item~\eqref{item:moll:four}. Arguing in a similar fashion as in the first part of this step but applying Remark~\ref{rem:upgrade.material.derivative.end} to the differences, choosing $\const_w=\mu_w=\tilde\mu_w=\const_v=\mu_v=\tilde\mu_v=\Tau_{q+1}^{-1}$ and using the extra prefactors from $\Tau_{q+1}^{4\Nindt}$ to absorb the lossy material derivative cost yields the desired estimates in \eqref{eq:Rcomm:bounds}--\eqref{eq:phicomm:bounds}.

\end{proof}

\section{Intermittent Mikado bundles and synthetic Littlewood-Paley decompositions}\label{ss:bundles}
In this section, we recall the geometric lemmas which enact the cubic and quadratic cancellations and the basic definitions of intermittent Mikado flows in subsection~\ref{ss:im}. Then in subsection~\ref{ss:ib}, we introduce intermittent Mikado bundles.  Finally, in subsectcion~\ref{sec:LP}, we introduce the synthetic Littlewood-Paley decomposition.

\subsection{Definition of intermittent Mikado flows and basic properties}\label{ss:im}

We shall require the following lemmas regarding decompositions of symmetric positive definite tensor fields. Typically such lemmas are stated and applied for tensors in a neighborhood of the identity.  Since it will be convenient for us to decompose tensors for which some rescaling of the original tensors belongs to a neighborhood of the identity, and later estimates (see Lemma~\ref{lem:a_master_est_p}) will depend on the rescaling factor, we include a slightly altered statement with full proof.

\begin{proposition}[\bf Geometric lemma I]\label{p:split}
Let $\Xi\subset \Q^3\cap\mathbb{S}^2$ denote the set $\left\{\sfrac 35 e_i \pm \sfrac 45 e_j\right\}_{1\leq i < j\leq 3}$, and for every $\xi$ in $\Xi$. Then there exists $\epsilon>0$ such that every symmetric 2-tensor in $B(\Id, \epsilon)$ can be written as a unique, positive linear combination of $\xi \otimes \xi$ for $\xi \in \Xi$.\index{$\Xi$} Furthermore, for a given large number $K > 1$, let $C_K$ denote the set 
\begin{equation}\label{defn:CK}
C_K := \bigcup_{1\leq k\leq K} B(k\Id, k\epsilon) \, ,
\end{equation}
which we note is contained in the set of positive definite, symmetric 2-tensors for $\epsilon$ sufficiently small. Then there exist functions $\gamma_{\xi,K}$ for $\xi\in\Xi$ such that every element $R \in C_K$ can also be written as a unique, positive linear combination
\begin{equation}\label{e:split}
R = \sum_{\xi\in\Xi} \left(\gamma_{\xi,K}(R)\right)^2  \xi\otimes \xi \, .
\end{equation}
Additionally, we have that for all $1\leq N\leq 3\Nfin$,
\begin{equation}\label{eq:gamma:xi:derivative:bounds}
    1\lesssim \left|\gamma_{\xi,K} \right| \lesssim K^{\sfrac 12}, 
    \quad \left| D^N \gamma_{\xi,K} \right| \lec 1
    \, , {\qquad \text{on } C_K}
\end{equation}
where the implicit constants above depend on $\Xi$ and $\Nfin$ but not $K$.
\end{proposition}

\begin{proof}
By direct computation, we have that the identity matrix can be written as a strictly positive linear combination of $\xi\otimes \xi$ for $\xi\in\Xi$, and that the set of simple tensors $\{\xi\otimes\xi\}_{\xi\in\Xi}$ is linearly independent in the set of symmetric matrices.  Therefore, there exists $\epsilon<1 $ and linear functions $(\gamma_\xi)^2$ for $\xi\in\Xi$ such that for all $R\in B(\Id,\epsilon)$, 
\begin{align*}
    R = \sum_{\xi\in\Xi} \gamma_\xi^2(R) \xi \otimes \xi \, ,
\end{align*}
and there exist implicit constants depending only on $\Xi$ such that for all $R\in B(\Id,\epsilon)$,
\begin{align}\label{imp:are:good}
    1 \lesssim \gamma_\xi(R) \lesssim 1 \, , \quad \left| D [\gamma_\xi^2(R)] \right| \lesssim 1 \, , \qquad D^N [\gamma_\xi^2(R)] \equiv 0 \quad \forall N \geq 2 \, .
\end{align}
Now let $K$ be given. We define $\gamma_{\xi,K}:C_K \to \R$ by
\begin{equation}\label{def:gamma:xi:R}
    \gamma_{\xi,K}^2 (R) := \gamma^2_\xi(R) = k \gamma_{\xi}^2 \left( \frac R k \right)  \, .
\end{equation}
In the last identity, $1\leq k\leq K$ is chosen to satisfy $\sfrac{R}{k} \in B(\Id,\epsilon)$ (cf. \eqref{defn:CK}), and the identity holds because of linearity of $(\gamma_\xi)^2$. Then, we have
\begin{equation}\notag
    \sum_{\xi\in\Xi} \gamma_{\xi,K}^2(R) \xi \otimes \xi = \sum_{\xi\in\Xi} \gamma_\xi^2\left( \frac R k \right) k \xi \otimes \xi = R \, ,
\end{equation}
and \eqref{e:split} is satisfied. Also, we have that for all $R\in C_K$,
\begin{align}
    1 \lesssim \gamma_{\xi,K}(R) \lesssim K^{\sfrac 12} \, , \quad \left| D [\gamma_{\xi,K}^2(R)] \right| \lesssim 1 \, , \qquad D^N [\gamma_{\xi,K}^2(R)] \equiv 0 \quad \forall N \geq 2 \, , \notag
\end{align}
where the implicit constants are those from \eqref{imp:are:good} and depend only on $\Xi$. We immediately deduce from the lower bound for $\gamma_{\xi,K}(R)$ that
\begin{align*}
    \left| D \gamma_{\xi,K} (R) \right| \leq \frac{\left| D[\gamma_{\xi,K}^2(R)] \right|}{\left| \gamma_{\xi,K}(R) \right|} \lesssim 1 \, .
\end{align*}
Now for $N\geq 1$, we may write that
\begin{align}
    2\gamma_{\xi,K}(R) D^{N+1} \gamma_{\xi,K}(R) = D^{N+1} \left( \gamma_{\xi,K}^2(R) \right) + \sum_{0<N'<N+1} c_{N,N'} D^{N'} \left( \gamma_{\xi,K}(R) \right) D^{N+1-N'} \left( \gamma_{\xi,K}(R) \right) \, . \notag
\end{align}
Assuming by induction that $|D^{N''} \gamma_{\xi,K}(R)| \lesssim 1$ for $1\leq N'' \leq N$, we use the lower bound for $\gamma_{\xi,K}(R)$ to divide both sides by $\gamma_{\xi,K}(R)$ and deduce that $|D^{N+1} \gamma_{\xi,K}(R)| \lesssim 1$, concluding the proof of \eqref{eq:gamma:xi:derivative:bounds}.
\end{proof}

We now recall~\cite[Lemma~3.3]{DK22}.

\begin{proposition*}[\bf Geometric lemma II]\label{prop:geo2} Let $\{\xi_1, \xi_2, \xi_3, \xi_4\}\subset \mathbb{Z}^{3}$ be a set of nonzero vectors satisfying\index{$\Xi'$}
\begin{align*}
    \{\xi_1, \xi_2, \xi_3\} \text{ is an orthogonal basis of }\R^3 \text{ and } \xi_4 := -(\xi_1+\xi_2+\xi_3).
\end{align*}
Fix $C_0>0$ and let $B_{C_0}:=\{\phi\in\R^3 :|\phi|\leq C_0  \}$. Then, there exist positive functions $\{\td\gamma_{\xi_i}\}_{i=1}^4\subset C^\infty(B_{C_0})$ such that for each $\phi\in B_{C_0}$, we have
\begin{align*}
    \phi = \frac12 \sum_{i=1}^4 (\td\gamma_{\xi_i}(\phi))^3 \xi_i \, . 
\end{align*}
In particular, the set $\{e_1, 2e_2, 2e_3, -(e_1 + 2e_2+ 2e_3)\}$ satisfies the assumption. We denote the set of their normalized vectors by $\Xi' := \{e_1, e_2, e_3, -\sfrac13(e_1 + 2e_2+ 2e_3)\}\subset \Q^3\cap \S^2$, and with slight abuse of the notation we redefine $\td \gamma_{\xi}$ to have 
\begin{align}
    2\phi =  \sum_{\xi\in \Xi'} (\td\gamma_{\xi}(\phi))^3 \xi \, . \label{e:split:ii}
\end{align}
\end{proposition*}

\begin{definition}
For any $\xi\in \Xi\cup \Xi'$, we choose $\xi',\xi'' \in \mathbb{Q}^3\cap\mathbb{S}^2$ such that $\{\xi,\xi',\xi''\}$ is an orthonormal basis of $\R^3$. We then denote by $n_\ast$\index{$n_\ast$} the least positive integer such that $n_\ast \xi, n_\ast \xi' n_\ast \xi'' \in \mathbb{Z}^3$ for all $\xi \in \Xi \cup \Xi'$.
\end{definition}

We now recall \cite[Proposition 4.3]{BMNV21}, which details the choices for shifts enjoyed by a function with sparse support.  In our setting, such functions will be pipe densities, or equivalently the densities associated to their potentials. \index{choice of shifts}
\begin{proposition}[\bf Rotating, Shifting, and Periodizing]\label{prop:pipe:shifted}
Fix $\xi\in\Xi$ (or $\in \Xi'$), where $\Xi$ is as in Proposition~\ref{p:split} (or as in Proposition \ref{prop:geo2}). Let ${r^{-1},\lambda\in\mathbb{N}}$ be given such that $\lambda r\in\mathbb{N}$. Let $\varkappa:\mathbb{R}^2\rightarrow\mathbb{R}$ be a smooth function with support contained inside a ball of radius $\sfrac{1}{4}$. Then for $k\in\{0,...,r^{-1}-1\}^2$, there exist functions $\varkappa^k_{\lambda,r,\xi}:\mathbb{R}^3\rightarrow\mathbb{R}$ defined in terms of $\varkappa$, satisfying the following additional properties:
\begin{enumerate}[(1)]
    \item \label{item:point:1} 
    We have that  $\varkappa^k_{\lambda,r,\xi}$ is simultaneously $\left(\frac{\mathbb{T}^3}{\lambda r}  \right)$-periodic and $\left(\frac{\Tthreexi}{\lambda r n_\ast}  \right)$-periodic. Here, by  $\T^3_\xi$ we refer to a rotation of the standard torus such that $\T^3_\xi$ has a face perpendicular to $\xi$.\index{$\T^3_\xi$}
    \item \label{item:point:2}  Let $F_\xi$ be one of the two faces of the cube $\frac{\Tthreexi}{\lambda r n_\ast}$ which is perpendicular to $\xi$. Let $\mathbb{G}_{\lambda,r}\subset F_\xi\cap \twopi \mathbb{Q}^3$ be the grid consisting of $r^{-2}$-many points spaced evenly at distance $\twopi  (\lambda n_\ast  )^{-1}$ on $F_\xi$ and containing the origin.  Then each grid point $g_{k}$ for $k\in\{0,...,r^{-1}-1\}^{2}$ satisfies
    \begin{equation}\label{e:shifty:support}
    \left(\supp\varkappa_{\lambda,r,\xi}^k\cap F_\xi \right) \subset \bigl\{x: |x-g_{k}| \leq \twopi\left(4\lambda n_\ast\right)^{-1} \bigr\}.
    \end{equation}
    
    \item \label{item:point:2a} The support of $\varkappa_{\lambda,r,\xi}^k$ is a pipe (cylinder) centered around a $\left(\frac{\mathbb{T}^3}{\lambda r}  \right)$-periodic and $\left(\frac{\Tthreexi}{\lambda r n_\ast}  \right)$-periodic line parallel to $\xi$, which passes through the point $g_k$. The radius of the cylinder's cross-section is as in \eqref{e:shifty:support}.
    \item We have that $\xi \cdot \nabla \varkappa_{\lambda,r,\xi}^k = 0$.
    \item \label{item:point:3} For $k\neq k'$, $\supp \varkappa_{\lambda,r,\xi}^k \cap \supp \varkappa_{\lambda,r,\xi}^{k'}=\emptyset$.
\end{enumerate}
\end{proposition}

We now state a slightly modified version of \cite[Proposition 4.4]{BMNV21} or equivalently \cite[Proposition~3.3]{NV22}, which rigorously constructs the $L^2$-normalized intermittent pipe flows and enumerates the necessary properties. 

\begin{proposition}[\bf Intermittent pipe flows for Reynolds corrector]
\label{prop:pipeconstruction}
Fix a vector $\xi$ belonging to the set of rational vectors $\Xi \subset\mathbb{Q}^{3} \cap \mathbb{S}^2 $ from Proposition~\ref{p:split}, $r^{-1},\lambda \in \mathbb{N}$ with $\lambda r\in \mathbb{N}$, and large integers $\Nfin$ and $\Dpot$. There exist vector fields $\mathcal{W}^k_{\xi,\lambda,r}:\mathbb{T}^3\rightarrow\mathbb{R}^3$ for $k\in\{0,...,r^{-1}-1\}^2$ and implicit constants depending on $\Nfin$ and $\Dpot$ but not on $\lambda$ or $r$ such that:
\begin{enumerate}[(1)]
    \item\label{item:pipe:1} There exists $\varrho:\mathbb{R}^2\rightarrow\mathbb{R}$ given by the iterated divergence $\div^\Dpot  \vartheta =: \varrho$ of a pairwise symmetric tensor potential $\vartheta:\mathbb{R}^2\rightarrow\mathbb{R}$ with compact support in a ball of radius $\frac{1}{4}$ such that the following holds.  Let $\varrho_{\xi,\lambda,r}^k$ and $\vartheta_{\xi,\lambda,r}^k$ be defined as in Proposition~\ref{prop:pipe:shifted}, in terms of $\varrho$ and $\vartheta$ (instead of $\varkappa$).  Then there exists $\mathcal{U}^k_{\xi,\lambda,r}:\mathbb{T}^3\rightarrow\mathbb{R}^3$ such that
    if $\{\xi,\xi',\xi''\} \subset \mathbb{Q}^3 \cap \mathbb{S}^2$ form an orthonormal basis of $\R^3$ with $\xi\times\xi'=\xi''$, then we have\footnote{The double index $ii$ indicates that $\div^{\Dpot-2} \left(\vartheta_{\xi,\lambda,r}^k \right)$ is a $2$-tensor, and we are summing over the diagonal components. The factor of $\sfrac 13$ appears because each component on the diagonal of this $3\times 3$ matrix is $\Delta^{-1} \varrho_{\xi,\lambda,r}^{k}$. The formula then follows from the identity $\curl \curl = -\Delta$ for divergence-free vector fields.}
    \begin{equation}
    \mathcal{U}_{\xi,\lambda,r}^k
     =  - \frac 13  \xi' \underbrace{\lambda^{-\Dpot} \xi''\cdot \nabla \left(\div^{\Dpot-2} \left(\vartheta_{\xi,\lambda,r}^k \right)\right)^{ii}}_{=:\varphi_{\xi,\lambda,r}^{\prime \prime k}}
     +  \frac 13
     \xi'' \underbrace{\lambda^{-\Dpot} \xi'\cdot \nabla \left(\div^{\Dpot-2} \left(\vartheta_{\xi,\lambda,r}^k \right)\right)^{ii}
     }_{=:\varphi_{\xi,\lambda,r}^{\prime k}}
     \label{eq:UU:explicit}
        \,, 
    \end{equation}
and thus
\begin{equation}
\label{eq:WW:explicit}
\curl \mathcal{U}^k_{\xi,\lambda,r} = \xi \lambda^{-\Dpot }\div^\Dpot  \left(\vartheta^k_{\xi,\lambda,r}\right) = \xi \varrho^k_{\xi,\lambda,r} =: \mathcal{W}^k_{\xi,\lambda,r}
\,,
\end{equation}
and 
\begin{equation}
    \xi \cdot \nabla \vartheta_{\xi,\lambda,r} =  (\xi \cdot \nabla) \mathcal{W}^k_{\xi,\lambda,r} 
    = (\xi \cdot \nabla )\mathcal{U}^k_{\xi,\lambda,r}
    = 0
    \,.
    \label{eq:derivative:along:pipe}
\end{equation}
    \item\label{item:pipe:2} The sets of functions $\{\mathcal{U}_{\xi,\lambda,r}^k\}_{k}$, $\{\varrho_{\xi,\lambda,r}^k\}_{k}$, $\{\vartheta_{\xi,\lambda,r}^k\}_{k}$, and $\{\mathcal{W}_{\xi,\lambda,r}^k\}_{k}$ satisfy items~\ref{item:point:1}--\ref{item:point:3} in Proposition~\ref{prop:pipe:shifted}.
    \item\label{item:pipe:3} $\mathcal{W}^k_{\xi,\lambda,r}$ is a stationary, pressureless solution to the Euler equations.
    \item\label{item:pipe:4} $\displaystyle{\dashint_{\mathbb{T}^3} \mathcal{W}^k_{\xi,\lambda,r} \otimes \mathcal{W}^k_{\xi,\lambda,r} = \xi \otimes \xi }$.
    \item\label{item:pipe:means}
    $\displaystyle{\dashint_{\mathbb{T}^3} |\mathcal{W}^k_{\xi,\lambda,r}|^2 \mathcal{W}^k_{\xi,\lambda,r} = \dashint_{\T^3} (\varrho_{\xi,\lambda,r}^k)^2 \mathcal{U}^k_{\xi,\lambda,r} = \int_{\T^3} \varrho^k_{\xi,\lambda,r} \mathcal{U}_{\xi,\lambda,r}^k = 0 \, . }$
    \item\label{item:pipe:5} For all $n\leq 3 \Nfin$, 
    \begin{equation}\label{e:pipe:estimates:1}
    {\left\| \nabla^n\vartheta^k_{\xi,\lambda,r} \right\|_{L^p(\mathbb{T}^3)} \lesssim \lambda^{n}r^{\left(\frac{2}{p}-1\right)} }, \qquad {\left\| \nabla^n\varrho^k_{\xi,\lambda,r} \right\|_{L^p(\mathbb{T}^3)} \lesssim \lambda^{n}r^{\left(\frac{2}{p}-1\right)} }
    \end{equation}
    and
    \begin{equation}\label{e:pipe:estimates:2}
    {\left\| \nabla^n\mathcal{U}^k_{\xi,\lambda,r} \right\|_{L^p(\mathbb{T}^3)} \lesssim \lambda^{n-1}r^{\left(\frac{2}{p}-1\right)} }, \qquad {\left\| \nabla^n\mathcal{W}^k_{\xi,\lambda,r} \right\|_{L^p(\mathbb{T}^3)} \lesssim \lambda^{n}r^{\left(\frac{2}{p}-1\right)} }.
    \end{equation}
    \item\label{item:pipe:3.5} We have that $\supp \vartheta_{\xi,\lambda,r}^k \subseteq B\left( \supp\varrho_{\xi,\lambda,r} ,2\lambda^{-1}\right)$.
    \item\label{item:pipe:6} Let $\Phi:\mathbb{T}^3\times[0,T]\rightarrow \mathbb{T}^3$ be the periodic solution to the transport equation
\begin{align}
\label{e:phi:transport}
\partial_t \Phi + v\cdot\nabla \Phi =0\,, 
\qquad 
\Phi|_{t=t_0} &= x\, ,
\end{align}
with a smooth, divergence-free, periodic velocity field $v$. Then
\begin{equation}\label{eq:pipes:flowed:1}
\nabla \Phi^{-1} \cdot \left( \mathcal{W}^k_{\xi,\lambda,r} \circ \Phi \right) = \curl \left( \nabla\Phi^T \cdot \left( \mathcal{U}^k_{\xi,\lambda,r} \circ \Phi \right) \right).
\end{equation}
\item\label{item:pipe:7} For any convolution kernel $K$, $\Phi$ as in \eqref{e:phi:transport}, $A=(\nabla\Phi)^{-1}$, and for $i=1,2,3$,
\begin{align}
\bigg{[} \nabla \cdot \bigg{(} A \, K\ast \left(  \mathcal{W}^k_{\xi,\lambda,r} \otimes  \mathcal{W}^k_{\xi,\lambda,r} \right)(\Phi) A^T \bigg{)} \bigg{]}_i 
& = A_{m}^j K\ast \left( (\mathcal{W}^k_{\xi,\lambda,r})^m  (\mathcal{W}_{\xi,\lambda,r}^k)^l (\Phi)\right) \partial_j A_{l}^i\nonumber\\
& = A_m^j \xi^m \xi^l \partial_jA_{l}^i \, K\ast \left( \left( \varrho^k_{\xi,\lambda,r} \right)^2(\Phi) \right) \, .
\label{eq:pipes:flowed:2}
\end{align}
In the above display, $k$ indicates the choice of placement, $i$ is the component of the vector field on either side of the equality, and $m$, $l$, and $j$ are repeated indices over which summation is implicitly encoded.
\end{enumerate}
\end{proposition}
\begin{proof}
The only small changes relative to the cited Propositions are as follows.  First, we write the pipe density $\varrho$ as the iterated \emph{divergence} of a pairwise symmetric vector potential $\div^\Dpot \vartheta= \varrho$ to match the form required for our inverse divergence operator (cf. Proposition~\ref{prop:intermittent:inverse:div}). By ``pairwise symmetric,'' we mean that permuting the $2n-1$ and $2n$ components for $1\leq n \leq \sfrac \Dpot 2$ leaves $\vartheta$ unchanged. Since one can always rewrite the identity $\Delta f = g$ as $\partial_i \partial_j \delta_{ij} f = g$, it is easy to convert the equality $\Delta^{\sfrac \Dpot 2} \tilde \vartheta=\varrho$ into $\div^\Dpot \vartheta = \varrho$ where $\vartheta$ is a pairwise symmetric tensor (see \eqref{eq:kronecker:tensor}). 

Second, \eqref{item:pipe:means} is new.  We will show that the second and third integrals vanish for any radial pipe density, while the first vanishes by choosing a suitable radial pipe density to have $\int_{\T^3} (\varrho_{\xi, \la, r}^k)^3 dx =0$. In order to compute the second and third integrands, we shall assume that $\xi=e_3$ and leave the case for general $\xi\in \Xi, \Xi'$ to the reader.  Since $\mathcal{U}_{e_3,\lambda,r}$ is mean-zero and divergence free, it can be written as the curl of a radial scalar potential $\mathcal{V}(r)$ according to the formula
$$ \mathcal{U}_{e_3,\lambda,r} = (-\partial_y \mathcal{V}_{e_3,\lambda,r}, \partial_x \mathcal{V}_{e_3,\lambda,r},0) \, . $$
Writing out the above expression in axial coordinates $(x,y,z)\mapsto(R,\theta,z)$ centered around the axis of a single cylinder of the pipe, we have
$$ \mathcal{U}_{e_3,\lambda,r}(R) = (-\sin(\theta) \mathcal{V}'_{e_3,\lambda,r}(R) , \cos(\theta) \mathcal{V}'_{e_3,\lambda,r}(R),0) \, . $$
Then since
$$ \int_{z_1}^{z_2} \int_{0}^{2\pi} \int_{R_1}^{R_2}  \sin(\theta) f(R) \, dR \, d\theta \, dz = \int_{z_1}^{z_2} \int_{0}^{2\pi} \int_{R_1}^{R_2}  \cos(\theta) f(R) \, dR \, d\theta \, dz $$
for any $R_1$, $R_2$, $z_1$, $z_2$ and radial function $f(R)$, and both the second and third integrals from~\eqref{item:pipe:means} can be written in this form, we see that the second and third integrals vanish as desired.

Finally, \eqref{item:pipe:3.5} is new, but it follows immediately from definitions and \eqref{e:shifty:support}.
\end{proof}

We shall require a set of intermittent pipe flows which possess nearly the same properties as above, but which are however normalized in $L^3$, and have non-vanishing cubic mean.  
\begin{proposition*}[\bf Intermittent pipe flows for current corrector]
\label{prop:pipe.flow.current}
Fix a vector $\xi$ belonging to the set of rational vectors $\Xi' \subset\mathbb{Z}^3 $ from Proposition~\ref{prop:geo2}.
The statement is same as in Proposition \ref{prop:pipeconstruction}, but item~\ref{item:pipe:4} is not imposed, and items~\ref{item:pipe:means}--\ref{item:pipe:5} are replaced by
\begin{enumerate}[(1)]
\setcounter{enumi}{4}
    \item\label{item:pipe:means:current}
    $\displaystyle{\dashint_{\mathbb{T}^3} |\mathcal{W}^k_{\xi,\lambda,r}|^2 \mathcal{W}^k_{\xi,\lambda,r} = |\xi|^2\xi} \, $, \quad $\displaystyle \dashint_{\T^3} (\varrho_{\xi,\lambda,r}^k)^2 \mathcal{U}^k_{\xi,\lambda,r} = \dashint_{\T^3} \varrho^k_{\xi,\lambda,r} \mathcal{U}^k_{\xi,\lambda,r} = 0$.
    \item\label{item:pipe:5:current} For all $n\leq 3 \Nfin$, 
    \begin{equation}\label{e:pipe:estimates:1:current}
    {\left\| \nabla^n\vartheta^k_{\xi,\lambda,r} \right\|_{L^p(\mathbb{T}^3)} \lesssim \lambda^{n} r^{\left(\frac{2}{p}-\frac23\right)} }, \qquad {\left\| \nabla^n\varrho^k_{\xi,\lambda,r} \right\|_{L^p(\mathbb{T}^3)} \lesssim \lambda^{n}r^{\left(\frac{2}{p}-\frac23\right)} }
    \end{equation}
    and
    \begin{equation}\label{e:pipe:estimates:2:current}
    {\left\| \nabla^n\mathcal{U}^k_{\xi,\lambda,r} \right\|_{L^p(\mathbb{T}^3)} \lesssim \lambda^{n-1}r^{\left(\frac{2}{p}-\frac23\right)} }, \qquad {\left\| \nabla^n\mathcal{W}^k_{\xi,\lambda,r} \right\|_{L^p(\mathbb{T}^3)} \lesssim \lambda^{n}r^{\left(\frac{2}{p}-\frac23\right)} } \, .
    \end{equation}
\end{enumerate}
\end{proposition*}
\begin{proof}
The differences in \eqref{item:pipe:5:current} relative to \eqref{item:pipe:5} from the preceding proposition are simply a result of the $L^3$ normalization and require no further justification.  In order to ensure \eqref{item:pipe:means:current}, it remains to show that one can construct a radial pipe density $\varrho_{\xi,\lambda,r}$ which has non-vanishing cubic mean and is the iterated Laplacian of a scalar potential, and then convert the scalar potential to a pairwise symmetric tensor potential.  As the latter task has already been carried out in the previous proposition, we can focus on the former. One can start with a smooth function $f:(\sfrac 12,1)\rightarrow \R$ for which $\int_0^{2\pi} (f^{(\Dpot)})^3(x) \, dx \neq 0$, and then define $F(r)=f(\lambda_1 r + \lambda_2)$, where $\lambda_1$ and $\lambda_2$ are chosen to ensure that to leading order, $\Delta^{\sfrac \Dpot 2}_r F \approx \lambda_1^\Dpot f^{(\Dpot)}(\lambda_1 r + \lambda_2)$.  Then periodizing concludes the proof.
\end{proof}

In order to control the geometry of pipes which are deformed by a velocity field on a local Lipschitz timescale, we recall \cite[Lemma~3.7]{NV22}.

\begin{lemma}[\bf Control on Axes, Support, and Spacing]
\label{lem:axis:control}
Consider a convex neighborhood  of space $\Omega\subset \mathbb{T}^3$. Let $v$ be an incompressible velocity field, and define the flow $X(x,t)$ and inverse $\Phi(x,t)=X^{-1}(x,t)$, which solves
\begin{align}\notag
\partial_t \Phi + v\cdot\nabla \Phi =0\,, 
\qquad
\Phi|_{t=t_0} &= x\, .
\end{align}
Define $\Omega(t):=\{ x\in\mathbb{T}^3 : \Phi(x,t) \in \Omega \} = X(\Omega,t)$. For an arbitrary $C>0$, let $\tau>0$ be a timescale parameter and $\Gamma > 3$ a large multiplicative prefactor such that the vector field $v$ satisfies the Lipschitz bound
\begin{equation}\notag
\sup_{t\in [t_0 - \tau,t_0+\tau]} \norm{\nabla v(\cdot,t) }_{L^\infty(\Omega(t))} \lesssim \tau^{-1} \Gamma^{-2} \, .
\end{equation}
Let $\mathcal{W}^k_{\xi,\lambda,r}:\mathbb{T}^3\rightarrow\mathbb{R}^3$ be a set of straight pipe flows constructed as in Proposition~\ref{prop:pipe:shifted}, Proposition~\ref{prop:pipeconstruction}, and Proposition~\ref{prop:pipe.flow.current}
which are $(\sfrac{\mathbb{T}}{\lambda r})^3$-periodic and concentrated around axes $\{A_i\}_{i\in\mathcal{I}}$ oriented in the vector direction $\xi$ for $\xi\in\Xi,\Xi'$, passing through the grid-points in item~\ref{item:point:2} of Proposition~\ref{prop:pipe:shifted}.  Then $\mathcal{W}:=\mathcal{W}^k_{\xi,\lambda,r}(\Phi(x,t)):\Omega(t)\times[t_0-\tau,t_0+\tau]$ satisfies the following conditions:
\begin{enumerate}[(1)]
	\item  We have the inequality
	\begin{equation}\label{eq:diameter:inequality}
	\textnormal{diam}(\Omega(t)) \leq \left(1+\Gamma^{-1}\right)\textnormal{diam}(\Omega) \, .
	\end{equation}
    \item If $x$ and $y$ with $x\neq y$ belong to a particular axis $A_i\subset\Omega$, then 
    \begin{equation}\label{e:axis:variation}
    \frac{X(x,t)-X(y,t)}{|X(x,t)-X(y,t)|} = \frac{x-y}{|x-y|} + \delta_i(x,y,t)    
    \end{equation}
    where $|\delta_i(x,y,t)|<\Gamma^{{-1}}$.
    \item Let $x$ and $y$ belong to $A_i\cap\Omega$ for some $i$, where the axes $A_i$ are defined above.  Denote the length of the axis $A_i(t):=X(A_i\cap\Omega,t)$ in between $X(x,t)$ and $X(y,t)$ by $L(x,y,t)$.  Then
    \begin{equation}\label{e:axis:length}
    L(x,y,t) \leq \left(1+\Gamma^{-1}\right)\left| x-y \right| \, .
    \end{equation}
    \item The support of $\mathcal{W}$ is contained in a $\displaystyle\left(1+\Gamma^{-1}\right)\twopi (4n_\ast\lambda)^{-1}$-neighborhood of the set
    \begin{equation}\label{e:axis:union}
       \bigcup_{i} A_i(t) \, .
    \end{equation}
\item $\mathcal{W}$ is ``approximately periodic" in the sense that for distinct axes $A_i,A_j$ with $i\neq j$, we have
\begin{equation}\label{e:axis:periodicity:1}
    \left(1-\Gamma^{-1}\right) \dist(A_i\cap\Omega,A_j\cap\Omega)
    \leq \dist\left(A_i(t),A_j(t)\right)
    \leq \left(1+\Gamma^{-1}\right) \dist(A_i\cap\Omega,A_j\cap\Omega) \, .
\end{equation}
\end{enumerate}
\end{lemma}

A consequence of Lemma~\ref{lem:axis:control} is that a set of $(\sfrac{\T}{\lambda r})^3$-periodic intermittent pipe flows which are flowed by a locally Lipschitz vector field on the Lipschitz timescale can be decomposed into ``segments of deformed pipe" in the sense of Remark~\ref{rem:deformed:pipes}.  Furthermore, any neighborhood of diameter $\approx (\lambda r)^{-1}$ contains at most a finite number of such segments of deformed pipe.\index{segments of deformed pipes}

\begin{definition}[\bf Segments of deformed pipes]\label{def:sunday:sunday}
A single ``segment of deformed pipe with thickness $\la^{-1}$ and spacing $(\la r)^{-1}$" is defined as a $3\lambda^{-1}$ neighborhood of a Lipschitz curve of length at most $2(\la r)^{-1}$.
\end{definition}

\subsection{Intermittent Mikado bundles}\label{ss:ib}
In the continuous scheme, the building block flows are {\it intermittent Mikado bundles}, which are bundles of pipes carefully designed to dodge previously placed intermittent Mikado bundles. To give the idea, suppose that intermittent Mikado bundles comprised of deformed pipes of thickness $\la_{q+1}^{-1}, \cdots \la_{q+\bn}^{-1}$ are given in a rectangular prism $\Omega_0$ of particular dimensions. If certain conditions are satisfied with respect to the spacing of the new bundles and the dimensions of the prism $\Omega_0$, we can successfully place new bundles of thickness $\la_{q+\bn}^{-1}$ that dodge all given bundles. Furthermore, the pipes in each new bundles will be placed to be at least at a distance $\la_{q+i}^{-1}\Gamma_{q+i}$ away from a given deformed pipe of thickness $\la_{q+i}^{-1}$. We call this additional property {\it effective dodging}, and it will play a crucial role throughout our scheme.

The key observation is that the intermittency alone need not dictate the spacing of the pipes in a bundle. For example, consider a set of pipes of thickness $\lambda_{q+\bn}^{-1}$ and spacing $\lambda_{q+\half}^{-1}$ restricted to the support of a set of a small number of pipes of thickness and spacing $\lambda_{q+1}^{-1}$.  An intermittent Mikado bundle is precisely such an object; a \emph{low} frequency, \emph{small} number of nearly \emph{homogeneous} pipes on which \emph{high} frequency, \emph{large} numbers of \emph{intermittent} pipes live. We call the nearly homogeneous pipes \emph{bundling} pipes.

\begin{proposition}[\bf ``Bundling" pipe flows $\rhob_{\xi,\diamond}^k$ for Reynolds and current correctors]\label{prop:bundling}
Fix a vector $\xi$ belonging to either of the sets of rational vectors from Propositions~\ref{p:split} or \ref{prop:geo2}. Then for $k\in\{1,\dots,\Gamma_q^6\}$, there exist master scalar functions $\ov\rhob_{\xi,k}$ and subsidiary bundling pipe flows $\rhob_{\xi,R}^k:=\ov\rhob_{\xi,k}^3$ for Reynolds correctors and $\rhob_{\xi,\varphi}^k:=\ov\rhob_{\xi,k}^2$ for current correctors satisfying the following.
\begin{enumerate}[(i)]
    \item\label{i:bundling:1} $\rhob^k_{\xi,\diamond}$ is $\left( \sfrac{\T}{\lambda_{q+1}\Gamma_q^{-4}}\right)^3$-periodic and satisfies $\xi \cdot \nabla \rhob^k_{\xi,\diamond}\equiv 0$, where either $\diamond=R$ or $\diamond=\varphi$.
    \item\label{i:bundling:2} The set of functions $\{\rhob^k_{\xi,\diamond}\}_{k}$ satisfies the conclusions of Proposition~\ref{prop:pipe:shifted} with $r^{-1}=\Gamma_q^3$, $\lambda=\lambda_{q+1}\Gamma_q^{-1}$.  In particular, $\supp \rhob^k_{\xi,\diamond} \cap \supp \rhob^{k'}_{\xi,\diamond} = \emptyset$ for $k\neq k'$, and there are $\Gamma_q^6$ disjoint choices of placement.
    \item\label{i:bundling:3} $\displaystyle \int_{\T^3} \ov\rhob_{\xi,k}^6=1$.
    \item\label{i:bundling:4} For all $n\leq 3\Nfin$ and $p\in[1,\infty]$,
    \begin{equation}\label{e:fat:pipe:estimates:1}
    \left\| \nabla^n \rhob^k_{\xi,R} \right\|_{L^p(\mathbb{T}^3)} \lesssim \left(\Gamma_q^{-1}\lambda_{q+1}\right)^n \Gamma_q^{-3\left(\frac 2p -1\right)} \, , \qquad \left\| \nabla^n \rhob^k_{\xi,\varphi} \right\|_{L^p(\mathbb{T}^3)} \lesssim \left(\Gamma_q^{-1}\lambda_{q+1}\right)^n \Gamma_q^{-3\left(\frac 2p - \frac 23 \right)} \, .
    \end{equation}
\end{enumerate}
\end{proposition}
\begin{proof}
The proof is a straightforward adaptation of the proofs of Propositions~\ref{prop:pipeconstruction} or \ref{prop:pipe.flow.current} after construction of an $L^6$ normalized master function $\ov\rhob_{\xi,k}$ which satisfies the shift and support properties from Proposition~\ref{prop:pipe:shifted}.  We omit further details.
\end{proof}

Now we further divide the support of the bundling pipes using the following anisotropic cutoffs and assign different pipes on the support of different cutoffs. We remark that these cutoffs have the same dimensions as the analogous objects in \cite[Definition~5.17]{NV22} and correspond to a length just larger than the scale to which the pipes have been periodizied, which is $(\la_{q+\bn}r_q)^{-1}$.

\begin{definition}[Strongly anisotropic cutoffs]\label{def:etab}
    To each $\xi \in \Xi$, we associate a partition of the orthogonal space $\xi^\perp \in \T^3$ into a grid\footnote{We refer to the grid used in Proposition~\ref{prop:pipe:shifted}, as any periodicity issues have been avoided there.} of squares of sidelength $\approx\lambda_{q+\half}^{-1}$. We index the squares $\mathcal{S}$ in this partition by $I_{\xi}$ which we will also denote by simply $I$. To this grid, we associate a partition of unity $ \etab_{\xi}^I$, i.e.,
    \begin{align}\label{eq:sa:summability}
         \etab_\xi^{I} =\begin{cases}
        1 &\textnormal{  on  } \frac34 \mathcal{S}_{I}\\
        0 &\textnormal{  outside  } \frac54 \mathcal{S}_{I}
        \end{cases}, \qquad
        \sum_{I} (\etab_\xi^{I})^{6} = 1 \, ,
    \end{align}
    which in addition satisfies $(\xi\cdot\na) \etab_\xi =0$ and $\left\| \nabla^N \etab_\xi^I \right\|_\infty \lesssim \lambda_{q+\half}^N$ for all $N\leq 3\Nfin$ and all $I$, where the implicit constants depend only on $\Xi$.  
\end{definition}

\begin{remark}\label{rem:strong:cardinality}
We note that the number of grid squares of sidelength $\la_{q+\half}^{-1}$ partitioning the orthogonal space $\xi^\perp \subset \T^3$ is $\lec \la_{q+\half}^2$. Consequently, we bound the cardinality of the index set $I$ as
\begin{align*}
    |\{I \in \mathcal{S} \}| \lec \la_{q+\half}^2 \,.
\end{align*}
\end{remark}

We now introduce \emph{intermittent pipe bundles}. These objects are \emph{multi-scale} and consist of nearly homogeneous bundling pipes at scale $\lambda_{q+1}^{-1}$, upon which various intermittent pipes are placed on the support of the strongly anisotropic cutoffs. 

\begin{definition}[\bf Intermittent pipe bundles]\label{defn:pipe.bundle}
   We define intermittent pipe bundles by
    \begin{align*}
        \BB_{\xi,R} =\chib_{\xi,R} \sum_I    (\etab_{\xi}^I)^3 \WW_{\xi,R}^I \quad
        \text{and}\quad
        \BB_{\xi,\ph} =\chib_{\xi,\ph} \sum_I    (\etab_{\xi}^I)^2 \WW_{\xi,\ph}^I.
    \end{align*}
where $\chib_{\xi,\diamond} = \chib_{\xi,\diamond}^m$ defined as in Proposition \ref{prop:bundling} for some $m=m_{\xi,\diamond}$ and $\WW_{\xi, \diamond}^I:=\mathcal{W}^{m'}_{\xi,\lambda_{q+\bn},\sfrac{\lambda_{q+\half}\Gamma_q}{\lambda_{q+\bn}}}$, constructed as in Propositions~\ref{prop:pipeconstruction} or \ref{prop:pipe.flow.current}, for some $m'=m'_{\xi,\diamond,I}$. \index{$\diamond$} We use $\diamond$ as a stand-in for either $R$ or $\varphi$ in order to streamline notation.
\end{definition}

\begin{remark}[\bf Choice of the placement]\label{rem:sat:sat:sat} \index{choice of placements}
    The placements $m$ and $m'$ will be chosen to have effective dodging with deformed pipes of thickness $\la_{q+1}^{-1}, \cdots, \la_{q+\half}^{-1}$ and that of thickness $\la_{q+\half+1}^{-1}, \cdots, \la_{q+\bn}^{-1}$, respectively. The requisite properties of these pipes are contained in Hypothesis~\ref{hyp:dodging2}.  The specifics of the placement procedure are contained in \cite[section~4]{GKN23}; see also the discussion following Lemma~\ref{lem:dodging}. 
\end{remark}

\begin{remark}[\bf Notational conventions]\label{rem:notational:conventions}
We shall frequently denote the intermittent pipe bundles defined above as follows:\index{$\BB_{\pxi,\diamond}$}\index{$\WW_{\pxi,\diamond}^I$}\index{$\rhob_\pxi^\diamond$}\index{$\pxi$}
\begin{equation}\label{int:pipe:bundle:convention}
    \BB_{\pxi,\diamond} = \rhob_\pxi^\diamond \sum_I \zetab_\xi^{I,\diamond} \WW_{\pxi,\diamond}^I \, .
\end{equation}
The meaning of this notation is as follows:
\begin{enumerate}[(i)]
    \item We assign a different intermittent Mikado bundle (where the difference is in terms of the placement mentioned in Remark~\ref{rem:sat:sat:sat}) to each mildly anisotropic checkerboard cutoff function $\zeta_{q,\diamond,i,k,\xi,\vecl}$ defined in Definition~\ref{def:checkerboard}. Therefore, the choice of placements $m$ for the bundling pipes will depend on all the indices for $\zeta_{q,\diamond,i,k,\xi,\vecl}$ , as well as the index $j$ for the pressure cutoffs defined in Definition~\ref{def:pressure:cutoff}. We will suppress these indices most of the time and simply write $(\xi)$ in parentheses, where the parentheses is a stand-in for the omitted indices $q,i,k,\vecl,j$. As a result, the bundling pipe has dependence on $(\xi), \diamond$, and so does the intermittent Mikado bundle.
    
    
    \item The subscript ``$\diamond$" in $\BB_{\pxi,\diamond}$ will be equal to either $\varphi$ or $R$, corresponding to velocity increments designed to correct current errors or stress errors, respectively.
    \item We abbreviate the bundling pipes $\rhob_{(\xi),\diamond}$ by $\rhob_{(\xi)}^\diamond$.  We write the $\diamond$ in the exponent to emphasize that the only difference between $\diamond=\varphi$ and $\diamond=R$ is the power of the scalar function $\ov \rhob_{\xi,k}$ used to define them.
    \item We abbreviate the very anisotropic cutoff functions by $\zetab_\xi^{I,\diamond}$.  We do \emph{not} write $\xi$ in parentheses, since $\zetab_\xi^{I,\diamond}$ does not depend on anything besides the vector direction $\xi$ and the index $I$ used to index the partition of unity. Also, the only difference between $\diamond=\varphi$ and $\diamond=R$ is the power, so we write $\diamond$ in the exponent.
    \item We write $\WW_{\pxi,\diamond}^I$ for the following reasons: first, the pipe flow depends on more indices than just $\xi$, so we write $\pxi$ to denote the omitted indices; we include the index $I$ to emphasize that the placement of the intermittent pipe flow depends not just on the omitted indices in $\pxi$, but on the index $I$ as well. 
    Finally, we leave $\diamond$ in the subscript since the difference between $\WW_{\pxi,R}^I$ and $\WW_{\pxi,\varphi}^I$ is more than just a power; the former has vanishing cubic mean, while the latter does not. We note that the placement of $\WW_{\pxi,\diamond}^I$ will depend on $(\xi), \diamond, I$.
\end{enumerate}
\end{remark}

\subsection{Synthetic Littlewood-Paley decomposition}\label{sec:LP}

When we estimate material derivatives of oscillation stress errors, we need dodging in order to estimate the application of the differential operator $\left(\hat u_{k-1} - \hat u_q\right)\cdot\nabla$ to the error; this operator appears in the material derivative estimates of the error term. To ensure that the error term enjoys a spatial support property even though it is defined using an inverse divergence operator and a frequency projection operator, we introduce a {\it synthetic Littlewood-Paley projector} $\tP_{(\la_1, \la_2]}$.\index{$\tP_{(\la_1, \la_2]}$}  This operator is defined using convolution with a compactly supported kernel, and thus behaves like the original projection operator $\mathbb P_{(\la_1, \la_2]}$ in estimates but allows control on the spatial support of the output.\index{synthetic Littlewood-Paley projector}

\begin{definition}[\bf Synthetic Littlewood-Paley projector]\label{def:synth:LP}
Let $\bph\in C_c^\infty(\R)$ satisfy
\begin{align*}
    \supp(\bph)\subset (-1/\sqrt 2, 1/\sqrt 2) \, , \qquad \int_{\R} \bph d s =1 \, , \qquad
    \int_{\R} s^n \bph ds =0
\end{align*}
for $n= 1, \dots, 10\Nfin$. Define $\bph_\la(\cdot) = \la \bph (\la \cdot)$, and set $\ph_\la(x) = \bph_\la(x_1)\bph_\la(x_2)$. For $f \in C^\infty(\T^2)$, we define the \emph{synthetic Littlewood-Paley projectors} by
\begin{align}
\tP_{\la}f(x):= \int_{\R^2} \ph_\la (y) f(x-y)  dy \, , \qquad
\tP_{(\la_1, \la_2]}f 
:= (\tP_{\la_2}- \tP_{\la_1})f \, ,  \label{eq:synth:LP} 
\end{align}
where in the convolution we consider $f$ as a periodic function defined on $\R^2$.\index{$\tP_{\la}$}
\end{definition}
From the definition, it is easy to see that $\supp(\ph_{\la_2}-\ph_{\la_1})\subseteq \supp(\ph_{\la_1})$ and hence $\supp(\tP_{(\la_1, \la_2]}f) \subset B(\supp(f), \la_1^{-1})$. With a bit of care, this property persists even after inverting the divergence.

\begin{lemma}[\bf Inverse divergence with spatial support property] \label{lem:LP.supp}
For given $f\in C^\infty(\T^2)$ and $\Dpot\geq 1$,\footnote{{The value of this number will be specified using the parameter $\dpot$ from item~\eqref{i:par:10}}.} there exists a symmetric tensor field $\Theta_f^{\la_1, \la_2}:\T^2 \to \R^{(2^\Dpot)}$ such that 
\begin{align}
    \tP_{(\la_1, \la_2]}(f) = \tP_{(\la_1, \la_2]}(f- \langle f \rangle)
    =\left(\lambda_1^{-1}\div\right)^{(\Dpot)}\Theta_f^{\la_1, \la_2} \, , \quad
    \supp\left(\Theta_f^{\la_1, \la_2}\right) \subset B(\supp(f), \la_1^{-1}) \, . \label{eq:LP:div:support:basic}
\end{align}
\end{lemma}
\begin{proof}
By a simple computation, we have
\begin{align}
\ph_{\la_2}(x) - \ph_{\la_1}(x)
= (\bph_{\la_2}(x_1)-\bph_{\la_1}(x_1))\bph_{\la_2}(x_2) + \bph_{\la_1}(x_1)(\bph_{\la_2}(x_2)-\bph_{\la_1}(x_2))  \, .  \label{eq:tricky:split}
\end{align}
Now define $g_0(z)= \bph_{\la_2}(z)- \bph_{\la_1}(z)$. We first construct a function $g_\Dpot(z) : \R \rightarrow \R$ with zero mean such that upon differentiating $\Dpot$ many times,
\begin{align*}
    g_\Dpot^{(\Dpot)} = g_0 \, , \qquad \supp(g_\Dpot)\subset (-(\sqrt2 \la_1) ^{-1}, (\sqrt2 \la_1) ^{-1}) \, . 
\end{align*}
The construction follows from applying the following claim iteratively: if $g_i\in C_c^\infty(\R)$ for some $i\in \{0, \dots, \Dpot-1\}$ satisfies $\int s^{n} g_i ds =0$ for all $n=0, \cdots, \Dpot-i$, then we can find $g_{i+1}$ such that
\begin{align*}
    g_{i+1}' = g_i \, , \qquad \supp(g_{i+1})\subset (-(\sqrt2 \la_1) ^{-1}, (\sqrt2 \la_1) ^{-1}) \, , \qquad \int_{\R} s^{n} g_{i+1} ds =0 \, \, \text{for } n=0 , \dots, \Dpot -i-1 \, .
\end{align*}
Assuming the claim, then $g_0$ satisfies $\int_{\R} s^{n} g_0(s) ds =0$ for $n=0, \cdots, \Dpot$, so we can find $g_{\Dpot}$ with zero-mean such that 
$$ g_\Dpot^{(\Dpot)} = g_{\Dpot-1}^{(\Dpot-1)} = \cdots =g_0 \, , \qquad \supp(g_\Dpot)  \subset (-(\sqrt2 \la_1) ^{-1}, (\sqrt2 \la_1) ^{-1}) \, . $$ 
To prove the claim, we define $g_{i+1}$ by $g_{i+1}(t) := \int_{-a}^{t} g_i ds$, where $a$ is chosen so that $\supp(g_i)\subset (-a, a)$. Since $g_i$ has zero-mean, we can easily see that $\supp(g_{i+1})\subset (-a, a)$, and $g_{i+1}(a)= g_{i+1}(-a)=0$. Using the latter, the vanishing moment condition follows from
\begin{align*}
    \int_{\R} s^{n} g_{i+1} ds 
    =\frac 1{n+1}\int_{-a}^a (s^{n+1})' g_{i+1} ds
    = - \frac 1{n+1} \int_{-a}^a s^{n+1} g_{i} ds
    = 0 \, .
\end{align*}

Now, we set $\theta_1^{(1,\dots, 1)}(x_1,x_2)= g_\Dpot(x_1) \bph_{\la_2} (x_2)$, and otherwise $\theta_1^{(i_1,\dots, i_\Dpot)}$ is zero, and $\theta_2^{(2,\dots, 2)}(x_1,x_2)=\bph_{\la_1}(x_1)g_\Dpot(x_2)$, and otherwise $\theta_2^{(i_1,\dots, i_\Dpot)}$ is zero. Then
\begin{align}
    &\pa_{i_1 \cdots i_\Dpot}\theta_1^{(i_1,\cdots, i_\Dpot)} = g_0(x_1) \bph_{\la_2}(x_2) \, ,
    \quad
    \supp(\theta_1^{(i_1,\cdots, i_\Dpot)}) \subset B(0,\la_1^{-1})
    \notag \\
    &\pa_{i_1 \cdots i_\Dpot}\theta_2^{(i_1,\cdots, i_\Dpot)} = \bph_{\la_1}(x_1)g_0(x_2) \, ,
    \quad
    \, \,  \supp(\theta_2^{(i_1,\cdots, i_\Dpot)}) \subset B(0,\la_1^{-1}) 
    \, .  \label{eq:support:props}
\end{align}
Lastly, we define the desired tensor function $\Theta_f^{\la_1, \la_2}$ by 
\begin{align}\label{defn.Xi}
    (\Theta_f^{\la_1, \la_2}) ^{(i_1,\dots, i_\Dpot)}(x_1,x_2)
    := \Theta \ast f (x_1,x_2)
    :=\lambda_1^\Dpot [(\theta_1+ \theta_2)^{(i_1,\cdots, i_\Dpot)}]\ast f (x_1,x_2) \, ,
\end{align}
which by \eqref{eq:tricky:split} and direct computation satisfies $\left(\lambda_1^{-1}\div\right)^{(\Dpot)}\Theta_f^{\lambda_1,\lambda_2} = \tP_{(\lambda_1,\lambda_2]} f$. The desired spatial support property follows from \eqref{defn.Xi} and \eqref{eq:support:props}.  We note that since $\varphi_{\lambda_2}-\varphi_{\lambda_1}$ has zero mean, $\tilde{\mathbb{P}}_{(\lambda_1,\lambda_2]}\langle f \rangle = 0$.
\end{proof}


With the previous Lemma in hand, we aim to apply various synthetic Littlewood-Paley projectors to smooth functions (such as squared pipe densities) and derive estimates for the projected function, and its ``inverse divergence potentials."   We shall generally decompose a smooth, $\left( \sfrac{\T}{\lambda r}\right)^3$-periodic function $\rho$ which has derivative cost $\lambda$ as a sum of the form
\begin{equation}
    \tilde{\mathbb{P}}_{\lambda_0} (\rho) + \left( \sum_{k=1}^{K}  \tilde{\mathbb{P}}_{(\lambda_{k-1},\lambda_k]} (\rho) \right) + \left( \Id -  \tilde{\mathbb{P}}_{\lambda_K} \right)(\rho) \, , \label{eq:decomp:showing} 
\end{equation}
where $\lambda_0$ is slightly larger than $\lambda r$, and $\lambda_K$ is slightly larger than $\lambda$. The terms in the sum are precisely of the form to which the previous lemma applies, and we estimate these in Lemma~\ref{lem:LP.est}.  The bottom and top shells which correspond to the two terms not in the summand are slightly unique cases; for these we record the following Lemma. Note that spatial localization is not relevant for these unique cases, as the lowest shell will have no spatial localization properties at all, and the highest shell will be vanishingly small.

\begin{lemma}[\bf Inverse divergence, special cases]\label{lem:special:cases}
Fix $q\in [1,\infty]$. Let $\Nblank$ a positive integer, $N_{**}\leq \sfrac{\Nblank}{2}$ a positive integer, $r,\lambda$ such that $\lambda r, \lambda \in \mathbb{N}$, and $\rho:(\sfrac{\T}{\lambda r})^2\rightarrow \R$ a smooth function such that there exists a constant $\const_{\rho,q}$ with
\begin{equation}\label{eq:moll:1:as}
 \left\| D^N \rho \right\|_{L^q(\T^2)} \lesssim \const_{\rho,q}  \lambda^N \, .
\end{equation}
for $N\leq \Nblank$.  Let $\lambda_0,\lambda_K$ be given with $\lambda r < \lambda_0 < \lambda < \lambda_K$. If the kernel $\overline{\varphi}$ used in Definition~\ref{def:synth:LP} has $N_{**}$ vanishing moments, then for $p\in[q,\infty]$ we have that
\begin{subequations}
\begin{align}
    \left\| D^N \left( \tilde{\mathbb{P}}_{\lambda_0} \rho \right) \right\|_{L^p} &\lesssim \const_{\rho,q} \left( \frac{\lambda_0}{\lambda r} \right)^{\sfrac 2q -\sfrac 2p} \lambda_0^N \qquad \qquad \qquad \forall N \leq \Nblank \, , \label{eq:lowest:shell:estimates}\\
    \left\| D^N \left( \left( \Id - \tilde{\mathbb{P}}_{\lambda_K} \right) \rho \right) \right\|_{L^\infty} &\lesssim \left(\frac{\lambda}{\lambda_K}\right)^{N_{**}} \const_{\rho,q} \lambda^{N+3} \qquad \forall N\leq \Nblank-N_{**}-3 \, .  \label{eq:remainder:estimates}
\end{align}
\end{subequations}
Furthermore, for any chosen positive even integer $\Dpot$ and any small positive number $\alpha$, there exist adjacent-pairwise symmetric\footnote{By ``adjacent-pairwise symmetric," we mean that permuting the $2n-1$ and $2n$ components for $1\leq n \leq \sfrac \Dpot 2$ leaves $\vartheta$ unchanged.} rank-$\Dpot$ tensor potentials $\vartheta_0$ and $\vartheta_K$ such that for $0\leq k \leq \Dpot$ and $N$ in the same range as above,
\begin{subequations}
\begin{align}
    \div^\Dpot \vartheta_0 &= \tilde{\mathbb{P}}_{\lambda_0} \mathbb{P}_{\neq 0}\rho \, , \qquad\qquad  \left\| D^N \div^{k} \vartheta_0 \right\|_{L^p} \lesssim \lambda_0^\alpha \const_{\rho,q} \left( \frac{\lambda_0}{\lambda r} \right)^{\sfrac 2q -\sfrac 2p} (\lambda r)^{k-\Dpot} \MM{N,\Dpot-k,\lambda r,\lambda_0} \, , \label{eq:lowest:shell:inverse} \\
    \div^\Dpot \vartheta_K &= (\Id - \tilde{\mathbb{P}}_{\lambda_K}) \rho \, , \qquad \left\| D^N \div^{k} \vartheta_K \right\|_{L^\infty} \lesssim \left(\frac{\lambda}{\lambda_K}\right)^{N_{**}} \const_{\rho,q}\lambda^{3} (\lambda r)^{k-\Dpot} \MM{N,\Dpot-k,\lambda r,\lambda}  \, .  \label{eq:remainder:inverse}
\end{align}
\end{subequations}
The implicit constants above depend on $\alpha$ but do not depend on $\lambda$, $\lambda_0$, $\lambda_K$, or $r$.
\end{lemma}
\begin{proof}
For the proof of \eqref{eq:lowest:shell:estimates}, we first define $F(x)=(\tilde{\mathbb{P}}_{\lambda r} \rho )(\sfrac{x}{\lambda r})$ to be the $1$-periodic rescaling of $\tilde{\mathbb{P}}_{\lambda r} \rho$. Then we can write that
\begin{align}
    \sup_{x\in\T^2} \left| D^N \left( \tilde{\mathbb{P}}_{\lambda r} \rho \right) \right|(x) &= (\lambda r)^{N} \sup_{x\in\T^2} \left| D^N F \right|(x) \notag\\
    &= (\lambda r)^{N} \sup_{x\in\T^2} \left| D_x^N \int_{\R^2} \rho(\sfrac{x}{\lambda r}-y) {\varphi}_{\lambda_0}( y) \, dy \right| \notag\\
    &= (\lambda r)^{N} \sup_{x\in\T^2} \left| D_x^N \int_{\R^2} \rho\left(\frac{x-z}{\lambda r}\right) \varphi_{\frac{\lambda_0}{\lambda r}} (z) \, dz \right| \notag\\
    &= (\lambda r)^{N} \sup_{x\in\T^2} \left| \int_{\R^2} \rho\left(\frac{x-z}{\lambda r}\right) (D^N_z \varphi_{\frac{\lambda_0}{\lambda r}}) (z) \, dz \right| \notag \\
    &\lesssim (\lambda r)^N \left( \frac{\lambda_0}{\lambda r} \right)^N \left( \frac{\lambda_0}{\lambda r}\right)^{\sfrac 2q} \const_{\rho,q} = \lambda_0^N \left( \frac{\lambda_0}{\lambda r}\right)^{\sfrac 2q} \const_{\rho,q} \notag
\end{align}
for all $N$, and in particular for all $N\leq\Nblank$.  This proves \eqref{eq:lowest:shell:estimates} for $p=\infty$, and the full estimate follows from interpolation with the trivial $L^q$ estimate.  To prove the second estimate, we use the vanishing moments condition to expand $\rho$ as a Taylor series and eliminate the first $N_{**}-1$ terms; in particular, we have that 
\begin{align}
    &\left| D^N \left( \left( \Id - \tilde{\mathbb{P}}_{\lambda_K} \right) \rho \right) \right|(x) \notag\\
    &\qquad = \left| \int_{\R^2} \varphi_{\lambda_K}(x-y) \left( \sum_{|\beta|=N_{**}} \frac{|\beta|(y-x)^\beta}{\beta!} \int_0^1 (1-\eta)^{N_{**}-1} D^\beta D^N \rho(x+\eta(y-x)) \, , d\eta \right) \,dy \right|  \notag \\
    &\qquad \lesssim \left\| D^{N+N_{**}} \varrho \right\|_{L^\infty} (\lambda_K)^{-N_{**}} \notag\\
    &\qquad \lesssim \left( \frac{\lambda}{\lambda_K} \right)^{N_{**}} \lambda^{N+3} \const_{\rho,q} \, . \notag 
\end{align}
The above computation holds for $N+N_{**}+3\leq \Nblank$, concluding the proof of the second estimate.

To prove the estimates for the tensor potentials, for $k=0,K$ we first define
\begin{subequations}\label{eq:kronecker:tensor}
\begin{align}
    \vartheta_0^{i_1 i_2 \dots i_{\Dpot -1} i_{\Dpot}} &= \delta^{i_1 i_2} \cdots \delta^{i_{\Dpot-1}i_{\Dpot}} \Delta^{-\frac \Dpot 2}\tilde{\mathbb{P}}_{\la_0}  \mathbb{P}_{\neq 0} \rho \, , \\
     \vartheta_K^{i_1 i_2 \dots i_{\Dpot -1} i_{\Dpot}}  &= \delta^{i_1 i_2} \cdots \delta^{i_{\Dpot-1}i_{\Dpot}}(\Id-\tilde {\mathbb{P}}_{\la_K}) \Delta^{-\frac \Dpot 2} \mathbb{P}_{\neq 0}\rho
\end{align}
\end{subequations}
where $\delta^{jl}$ is the usual Kronecker delta.  Then by direct computation and standard Littlewood-Paley analysis, \eqref{eq:lowest:shell:inverse} and \eqref{eq:remainder:inverse} hold. The $\alpha$ loss in the first estimate is due to the failure of the Calderon-Zygmund inequality in endpoint cases.
\end{proof}

We now move to the middle cases from \eqref{eq:decomp:showing}, for which the spatial localization will be important.
\begin{lemma}[\bf General localized inverse divergence]\label{lem:LP.est}
Fix $q\in [1,\infty]$. Let $\rho:\T^2\rightarrow \R$ be a smooth function which is $\left( \sfrac{\T}{\lambda r} \right)^2$-periodic and for $N\leq 2\Nfin$ satisfies
\begin{equation}
    \left\| D^N \rho \right\|_{L^q(\T^2)} \lesssim \const_{\rho,q} \lambda^N \, . \label{eq:gen:id:assump}
\end{equation}
For $\lambda r < \lambda_1 < \lambda_2$, define $\Theta_\rho^{\lambda_1,\lambda_2}$ using Lemma~\ref{lem:LP.supp}. Then for $p\in[q,\infty]$, $0\leq k \leq \Dpot$, $0<\alpha\ll 1$, and $N\leq \Nfin$, we have
\begin{subequations}
\begin{align}
    \left(\lambda_1^{-1}\div\right)^{(\Dpot)}\Theta_\rho^{\la_1, \la_2} &= \tP_{(\la_1, \la_2]}(\rho) = \tP_{(\la_1, \la_2]}(\rho- \langle \rho \rangle)
     \label{eq:LP:equality} \\
    \norm{D^N\pa_{i_1\cdots i_{\Dpot-k}}
    (\lambda_1^{-\Dpot}\Theta_{\rho}^{\la_1, \la_2})^{(i_1, \cdots, i_\Dpot)}}_{L^p(\T^2)} &\lec_{\Dpot,\alpha} \const_{\rho,q} \left(\frac{\min\left(\lambda,\la_2\right)}{\la r}\right)^{\frac 2q -\frac 2p+\alpha}  \la_1^{-k} \min\left(\lambda,\la_2\right)^N  \, , \label{eq:LP:div:estimates} \\
    \supp(\Theta_{\rho}^{\la_1, \la_2}) &\subset B(\supp(\rho) , \la_1^{-1}) \, .  \label{eq:LP:div:support}
\end{align}
\end{subequations}
The implicit constants above depend on $\alpha$ but do not depend on $\lambda$, $\lambda_1$, $\lambda_2$, or $r$.
\end{lemma}
\begin{proof}
The spatial property immediately follows from Lemma \ref{lem:LP.supp}. To obtain $L^p$-norm estimates, we will obtain $L^q$ and $L^\infty$ norm estimates and then interpolate them.  We first rescale by setting
\begin{align}
    \tilde\rho (\cdot) = \rho\left( \frac{\cdot}{\lambda r} \right) \, , \qquad \tilde\lambda_1 = \frac{\lambda_1}{\lambda r} \, , \qquad \tilde\lambda_2 = \frac{\lambda_2}{\lambda r} \, , \qquad \tilde \lambda = \frac{\lambda}{\lambda r} = r^{-1} \, , \label{eq:rescaled:choices}
\end{align}
so that $\tilde \rho$ is $\T^2$ periodic and satisfies
\begin{align}
    \left\| D^N \tilde \rho \right\|_{L^q(\T^2)} \lesssim \const_{\rho,q} \tilde\lambda^{N} \, . \notag
\end{align}
Constructing $\theta_1$ and $\theta_2$ as in the previous lemma but for the choices in \eqref{eq:rescaled:choices}, we have
$$ \pa_{1}^{\Dpot-k} \theta_1^{(1,\dots, 1)}(x_1,x_2)= g_{k}(x_1)\bph_{\tilde\la_2}(x_2) \, , \quad \pa_{2}^{\Dpot-k} \theta_2^{(2,\dots,2)}(x_1,x_2)=\bph_{\tilde\la_1}(x_1) g_{k}(x_2)\, . $$
By direction computation, i.e. simply integrating a difference of mollifiers, we have that $\tilde g_{k}$ satisfies
\begin{align*}
\norm{D^N g_{k}}_{L^1(\R)}
&\lesssim_\Dpot \tilde\lambda_1^{-k} \MM{N,k-1,\tilde\lambda_1,\tilde\lambda_2} \, , \qquad
\left\| D^N g_k \right\|_{L^\infty(\R)} \lesssim_\Dpot \tilde\lambda_1^{1-k} \MM{N,k-1,\tilde\lambda_1,\tilde\lambda_2} \, , \qquad k\geq 1 \, , \notag\\
\norm{D^N g_{0}}_{L^1(\R)}
&\lesssim_\Dpot \tilde \lambda_2^N \, , \qquad\qquad\qquad \qquad\qquad \qquad 
\left\| D^N g_0 \right\|_{L^\infty(\R)} \lesssim_\Dpot \tilde \lambda_2^{N+1} \, .
\end{align*}
Then we have the bounds
\begin{align*}
     &\norm{D^N\pa_{1}^{\Dpot-k} \theta_1^{(1,\dots, 1)}}_{L^1(\R^2)} \lec_\Dpot \tilde \lambda_2^N \tilde \lambda_1^{-k} \, , \qquad
     \norm{D^N\pa_{1}^{\Dpot-k} \theta_1^{(1,\dots, 1)}}_{L^\infty(\R^2)} \lec_{\Dpot} \tilde \lambda_2^{N+2} \tilde \lambda_1^{-k} \, , \\
     &\norm{D^N\pa_{2}^{\Dpot-k} \theta_2^{(2,\dots, 2)}}_{L^1(\R^2)} \lec_\Dpot \tilde \lambda_2^N \tilde \lambda_1^{-k} \, , \qquad
     \norm{D^N\pa_{2}^{\Dpot-k} \theta_2^{(2,\dots, 2)}}_{L^\infty(\R^2)} \lec_{\Dpot} \tilde \lambda_2^{N+1} \tilde \lambda_1^{-k+1} \, .
\end{align*}
Thus it follows by interpolation for $1/q' = 1 - 1/q$ that
\begin{align*}
    \norm{D^N\pa_{1}^{\Dpot-k} \theta_1^{(1,\dots, 1)}}_{L^{q'}(\R^2)} \lec_\Dpot \tilde \lambda_2^{N + \sfrac 2q} \tilde \lambda_1^{-k} \, ,
    \qquad
     \norm{D^N\pa_{2}^{\Dpot-k} \theta_2^{(2,\dots, 2)}}_{L^{q'}(\R^2)} \lec_{\Dpot} \tilde \lambda_2^{N + \sfrac 1q} \tilde \lambda_1^{-k+1} \, .
\end{align*}
We therefore have that for $k=0,\dots,\Dpot$,
\begin{align*}
\norm{D^N\pa_{i_1\cdots i_{\Dpot-k}}(\Theta_{\tilde \rho}^{\tilde \la_1, \tilde \la_2})^{(i_1, \cdots, i_\Dpot)}}_{L^q(\T^2)}
&\lec \tilde \la_1^{\Dpot-k}\min\left(\tilde\lambda,\tilde\la_2\right)^N \const_{\rho,{q}} \\
\norm{D^N\pa_{i_1\cdots i_{\Dpot-k}}(\Theta_{\tilde \rho}^{\tilde \la_1, \tilde \la_2})^{(i_1, \cdots, i_\Dpot)}}_{L^\infty(\T^2)} &\lec_{\Dpot} \tilde \la_1^{\Dpot-k}\min\left(\tilde\lambda,\tilde\la_2\right)^{N+ {\sfrac 2q}+\alpha} \const_{\rho,{q}} \, ,
\end{align*}
where if $\tilde\lambda_2\leq\tilde\lambda$, we let the derivatives fall on $\theta_i$, and if $\tilde\lambda_2>\tilde\lambda$, we let the derivatives fall on $\tilde \rho$.
Using the interpolation inequality, we obtain
\begin{align*}
    \norm{D^N\pa_{i_1\cdots i_{\Dpot-k}}(\Theta_{\tilde \rho}^{\tilde \la_1, \tilde \la_2})^{(i_1, \cdots, i_\Dpot)}}_{L^p(\T^2)} \lec_\Dpot \tilde \la_1^{\Dpot-k}\min(\tilde\lambda,\tilde\lambda_2)^{N+ \sfrac 2q -\sfrac 2p +\alpha} \const_{\rho,{q}} \, .
\end{align*}
Undoing our original rescaling, we find that
\begin{align*}
    \norm{D^N\pa_{i_1\cdots i_{\Dpot-k}}(\Theta_{\rho}^{\la_1,\la_2})^{(i_1, \cdots, i_\Dpot)}}_{L^p(\T^2)} &\lec_\Dpot \left(\lambda r\right)^{N+\Dpot-k} \norm{D^{N} \left[ \pa_{i_1\cdots i_{\Dpot-k}}(\Theta_{\tilde \rho}^{\tilde \la_1, \tilde \la_2})^{(i_1, \cdots, i_\Dpot)}\right]}_{L^p(\T^2)} \\
    &\leq \left( \frac{\min(\lambda,\lambda_2)}{\lambda r} \right)^{\frac 2q -\frac 2p+\alpha} \const_{\rho,{q}} \lambda_1^{\Dpot-k} \min(\lambda,\lambda_2)^{N} \, .
\end{align*}
\end{proof}

\section{Non-inductive cutoffs}\label{ss:nic}
In this section, we introduce all the non-inductive cutoffs which will be required throughout the proof.  First, we introduce a collection of time cutoffs in subsection~\ref{sec:cutoff:temporal:definitions}.  Then in subsection~\ref{s:deformation}, we can estimate flow maps related to the flow of $\nabla \hat u_{q'}$ for $q'\leq \qbn-1$ on the support of time and velocity cutoffs.  Then in subsection~\ref{ss:ipc}, we introduce the intermittent pressure cutoffs for $\pi_\ell$.  Subsection~\ref{sec:cutoff:checkerboard:definitions} contains the definitions and estimates for the mildly and strongly anistropic checkerboard cutoffs, whose properties are put to use in the discussion following Lemma~\ref{lem:dodging}.  Finally, in subsection~\ref{sec:cutoff:total:definitions}, we introduce the cumulative cutoff functions given as a product of all previously defined types of cutoffs.  The last subsection of this section then contains a number of ``cutoff aggregation lemmas'' which allow us to turn estimates in localized regions of space-time into global pointwise and $L^p$ bounds.

\subsection{Time cutoffs}
\label{sec:cutoff:temporal:definitions}
Let $\chi:(-1,1)\rightarrow[0,1]$ be a $C^\infty$ function which induces a partition of unity according to
\begin{align}
 \sum_{k \in \Z} \chi^6(\cdot - k) \equiv 1 \, .
\label{eq:chi:cut:partition:unity}
\end{align}
Consider the translated and rescaled function 
\begin{equation*}
    \chi\left(2 t \tau_{q}^{-1}\Gamma^{i+2}_{q} - k\right) \, ,
\end{equation*}
which is supported in the set of times $t$ satisfying
\begin{equation}\label{eq:chi:support}
\left| t - \sfrac 12 \tau_q \Gamma_{q}^{-i-2} k \right| \leq \sfrac 12 \tau_q\Gamma_{q}^{-i-2} \quad \iff t\in \left[ (k-1)\sfrac 12 \tau_q \Gamma_{q}^{-i-2}, (k+1) \sfrac 12 \tau_q \Gamma_{q}^{-i-2} \right]  \, .
\end{equation}
We then define temporal cut-off functions\index{$\chi_{i,k,q}$}
\begin{align}
 \chi_{i,k,q}(t) =  \chi\left(2t \tau_{q}^{-1}\Gamma^{i+2}_{q} - k\right) \, .
 \label{eq:chi:cut:def}
\end{align}
It is then clear that 
\begin{align}
{|\partial_t^m \chi_{i,k,q}| \les (\Gamma_{q}^{i+2} \tau_{q}^{-1})^m}
\label{eq:chi:cut:dt}
\end{align}
for $m\geq 0$ and
\begin{equation}\label{e:chi:overlap}
    \chi_{i,k_1,q}(t)\chi_{i,k_2,q}(t) = 0
\end{equation}
for all $t\in\mathbb{R}$ unless $|k_1-k_2|\leq 1$. In analogy to $\psi_{i\pm,q}$, we define
\begin{equation}\label{e:chi:plus:minus:definition}
    \chi_{i, k \pm, q}(t) := \left( \chi_{i,k-1,q}^6(t) + \chi_{i,k,q}^6(t) + \chi_{i,k+1,q}^6(t)  \right)^\frac{1}{6} \, ,
\end{equation}
which are cutoffs with the property that
\begin{equation}\label{e:chi:overlaps}
    \chi_{i,k\pm,q} \equiv 1 \textnormal{ on } \supp{(\chi_{i,k,q})} \, .
\end{equation}
Next, we define the cutoffs $\tilde\chi_{i,k,q}$\index{$\tilde \chi_{i,k,q}$} by
\begin{equation}\label{eq:chi:tilde:cut:def}
\tilde\chi_{i,k,q}(t) = \chi\left( t \tau_q^{-1}\Gamma_{q}^{i} - k\Gamma_q^{-2} \right) \, .
\end{equation}
For comparison with \eqref{eq:chi:support}, we have that $\tilde\chi_{i,k,q}$ is supported in the set of times $t$ satisfying
\begin{equation}\label{eq:chi:tilde:support}
\left| t-\tau_q \Gamma_{q}^{-i-2} k \right| \leq \tau_q\Gamma_{q}^{-i} \,.
\end{equation}
Let $(i,k)$ and $(\istar,\kstar)$ be such that $\supp \chi_{i,k,q} \cap \supp \chi_{\istar,\kstar,q}\neq\emptyset$ and $\istar\in\{i-1,i,i+1\}$. Then as a consequence of these definitions and a sufficiently large choice of $\lambda_0$,
\begin{equation}\label{eq:tilde:chi:contains}
\supp \chi_{i,k,q} \subset \supp \tilde\chi_{\istar,\kstar,q}\, .
\end{equation}

\subsection{Estimates on flow maps}\label{s:deformation}
\label{sec:cutoff:flow:maps}
We can now make estimates regarding the flows of the vector field $\hat u_{q'}$ for $q'\leq q+\bn-1$ on the support of a velocity and time cutoff function. This section is completely analogous to \cite[Section~6.4]{BMNV21}, and we omit the proofs.
\begin{lemma}[\bf Lagrangian paths don't jump many supports]
\label{lem:dornfelder}
Let $q'\leq q+\bn-1$ and $(x_0,t_0)$ be given. Assume that the index $i$ is such that $\psi_{i,q'}^2(x_0,t_0) \geq \kappa^2$, where $\kappa\in\left[\frac{1}{16},1\right]$. Then the forward flow $(X(t),t) := (X(x_0,t_0;t),t)$ of the velocity field $\hat u_{q'}$ originating at $(x_0,t_0)$ has the property that $\psi_{i,q'}^2(X(t),t) \geq\sfrac{\kappa^2}{2}$ for all $t$ such that $|t - t_0|\leq \tau_{q'} \Gamma_{q'}^{-i+4}$.
\end{lemma}

We note that $\psi_{i,q'}$ for $q'\leq q+\bn-1$ are given inductively. The proof of the lemma uses their properties recorded in subsection \ref{sec:cutoff:inductive} only.

\begin{corollary}[\bf Backwards Lagrangian paths don't jump many supports]
\label{cor:dornfelder}
Suppose $(x_0,t_0)$ is such that $\psi^2_{i,q'}(x_0,t_0)\geq \kappa^2$, where $\kappa\in\left[\sfrac{1}{16},1\right]$. For $\abs{t-t_0}\leq\tau_{q'}\Gamma_{q'}^{-i+3}$, define $x$ to satisfy
\[
x_0=X(x,t;t_0) \, .
\]
That is, the forward flow $X$ of the velocity field $\hat u_{q'}$, originating at $x$ at time $t$, reaches the point $x_0$ at time $t_0$.
Then we have
\begin{equation*}
\psi_{i,q'}(x,t)\neq 0 \,.
\end{equation*}
\end{corollary}

\begin{definition}[\bf Flow maps]\label{def:transport:maps} We define $\Phi_{i,k,q'}(x,t)=\Phi_{(i,k)}(x,t)$ to be the flows induced by $\hat u_{q'}$ with initial datum at time $k {\tau_{q'}}\Gamma_{q}^{-i-2}$ given by the identity, i.e.\index{$\Phi_{i,k,q}$}\index{$\Phiik$}\index{flow maps}
\begin{equation}\label{e:Phi}
\left\{\begin{array}{l}
(\partial_t + \hat u_{q'} \cdot\nabla) \Phi_{i,k,q'} = 0 \\
\Phi_{i,k,q'}(x,k{\tau_{q'}}\Gamma_{q'}^{-i-2})=x\, .
\end{array}\right.
\end{equation}
\end{definition}
\noindent We will use $D\Phi_{(i,k)}$ to denote the gradient of $\Phi_{(i,k)}$ (which is a thus matrix-valued function).  The inverse of the matrix $D\Phi_{(i,k)}$ is denoted by $\left(D\Phi_{(i,k)}\right)^{-1}$, in contrast to $D\Phi_{(i,k)}^{-1}$, which is the gradient of the inverse map $\Phi_{(i,k)}^{-1}$.

\begin{corollary}[\bf Deformation bounds]
\label{cor:deformation}
For $k \in \Z$, $0 \leq i \leq  i_{\rm max}$, $q'\leq q+\bn-1$, and $2 \leq N \leq \sfrac{3\Nfin}{2}+1$, we have the following bounds on the support of $\psi_{i,q'}(x,t){\tilde\chi_{i,k,q'}(t)}$.
\begin{subequations}
\begin{align}
\norm{D\Phi_{(i,k)} - {\rm Id}}_{L^\infty(\supp(\psi_{i,q'} \tilde\chi_{i,k,q'} ))} &\lesssim \Gamma_{q'}^{-1}
\label{eq:Lagrangian:Jacobian:1}\\
\norm{D^N\Phi_{(i,k)} }_{L^\infty(\supp(\psi_{i,q'} \tilde\chi_{i,k,q'} ))} & \lesssim \Gamma_{q'}^{-1} (\lambda_{q'}\Ga_{q'})^{N-1} \label{eq:Lagrangian:Jacobian:2}\\
\norm{(D\Phi_{(i,k)})^{-1} - {\rm Id}}_{L^\infty(\supp(\psi_{i,q'} \tilde\chi_{i,k,q'} ))} & \lesssim \Gamma_{q'}^{-1}\label{eq:Lagrangian:Jacobian:3}\\
\norm{D^{N-1}\left((D\Phi_{(i,k)})^{-1}\right) }_{L^\infty(\supp(\psi_{i,q'} \tilde\chi_{i,k,q'} ))} & \lesssim \Gamma_{q'}^{-1} (\lambda_{q'}\Ga_{q'})^{N-1} \label{eq:Lagrangian:Jacobian:4} \\
\norm{D^N\Phi^{-1}_{(i,k)} }_{L^\infty(\supp(\psi_{i,q'} \tilde\chi_{i,k,q'} ))} & \lesssim \Gamma_{q'}^{-1} (\lambda_{q'}\Ga_{q'})^{N-1} \label{eq:Lagrangian:Jacobian:7}
\end{align}
\end{subequations}
Furthermore, we have the following bounds for $1\leq N+M\leq \sfrac{3\Nfin}{2}$ and $0\leq N'\leq N$:
\begin{subequations}
\begin{align}
\left\| D^{N-N'} D_{t,q'}^M D^{N'+1} \Phi_{(i,k)} \right\|_{L^\infty(\supp(\psi_{i,q'}\tilde\chi_{i,k,q'}))} &\leq  (\lambda_{q'}\Ga_{q'})^{N} \MM{M,\NindSmall,\Gamma_{q'}^{i} \tau_{q'}^{-1},\Tau_{q'-1}^{-1}\Gamma_{q'-1}}\label{eq:Lagrangian:Jacobian:5}\\
\left\| D^{N-N'} D_{t,q'}^M D^{N'} (D \Phi_{(i,k)})^{-1} \right\|_{L^\infty(\supp(\psi_{i,q'}\tilde\chi_{i,k,q'}))} &\leq (\lambda_{q'}\Ga_{q'})^{N} \MM{M,\NindSmall,\Gamma_{q'}^{i} \tau_{q'}^{-1},{\Tau}_{q'-1}^{-1}\Gamma_{q'-1}} \, . \label{eq:Lagrangian:Jacobian:6}
\end{align}
\end{subequations}
\end{corollary}

\subsection{Intermittent pressure cutoffs}\label{ss:ipc}

In this section, we introduce cutoff functions for the level sets of $\pi_\ell$. Estimates for $\pi_\ell$ are provided by \eqref{eq:pressure:inductive:dtq:upgraded}--\eqref{eq:pressure:inductive:dtq:pointwise}.

\subsubsection{Definition of the intermittent pressure cutoffs}
\label{sec:cutoff:stress:definitions}
We first introduce a partition of unity which is slightly more general than is needed at the moment; however, the generality will prove useful in the construction of the velocity cutoffs. The statement is almost identical to \cite[Lemma~6.2]{BMNV21}.  The only slight difference is that \eqref{eq:tilde:partition} holds for the sixth power (the least common multiple of two and three, corresponding to cubic and quadratic error terms, respectively), and the estimates in \eqref{item:cutoff:estimates} hold for arbitrary integer powers of the cutoff functions. The more general bounds follow from the fact that since the cutoff functions are defined by gluing together exponential functions, raising to a power is (locally) equivalent to dilation.
\begin{lemma}\label{lem:cutoff:construction:first:statement}
For all $q\geq 1$ and $0\leq m \leq \NcutSmall$, there exist smooth cutoff functions $\tilde\gamma_{m,q},\gamma_{m,q}:[0,\infty)\rightarrow[0,1]$ which satisfy the following.
\begin{enumerate}[(1)]
    \item\label{item:cutoff:1} The function $\tilde\gamma_{m,q}$ satisfies ${\bf 1}_{[0,\frac{1}{4}\Gamma_{q}^{2(m+1)}]} \leq \tilde\gamma_{m,q} \leq {\bf 1}_{[0,\Gamma_{q}^{2(m+1)}]}$.
    \item\label{item:cutoff:2}  The function $\gamma_{m,q}$ satisfies ${\bf 1}_{[1,\frac{1}{4}\Gamma_{q}^{2(m+1)}]} \leq \gamma_{m,q} \leq {\bf 1}_{[\frac{1}{4},\Gamma_{q}^{2(m+1)}]}$.
    \item For all $y\geq 0$, a partition of unity is formed as
    \begin{align}
    \tilde \gamma_{m,q}^{6}(y) + \sum_{{i\geq 1}} \gamma_{m,q}^{ 6}\bigl(\Gamma_{q}^{-2i(m+1)} y\bigr) = 1 \, .
    \label{eq:tilde:partition}
    \end{align}
    \item $\tilde\gamma_{m,q}$ and $\gamma_{m,q}(\Gamma_{q}^{-2i(m+1)}\cdot)$ satisfy
    \begin{align}
   \supp \tilde\gamma_{m,q}(\cdot) \cap \supp \gamma_{m,q}\bigl(\Gamma_{q}^{-2i(m+1)}\cdot\bigr) &= \emptyset \quad \textnormal{if} \quad i \geq 2,\notag\\
   \supp \gamma_{m,q}\bigl(\Gamma_{q}^{-2i(m+1)}\cdot\bigr) \cap \supp \gamma_{m,q}\bigl(\Gamma_{q}^{-2i'(m+1)}\cdot\bigr) &= \emptyset \quad \textnormal{if} \quad |i-i'|\geq 2 \, . \label{eq:psi:support:base:case}
    \end{align}
    \item\label{item:cutoff:estimates} For $0\leq N \leq \Nfin$, when $0\leq y<\Gamma_{q}^{2(m+1)}$ we have
    \begin{align}
    {|D^N \tilde \gamma_{m,q}(y)|}    &\lesssim  {(\tilde \gamma_{m,q}(y))^{1-N/\Nfin}}
    \Gamma_{q}^{-2N(m+1)}. \label{eq:DN:psi:q:0}
    \end{align}
For $\frac{1}{4}<y<1$ we have
    \begin{align}
     {|D^N  \gamma_{m,q}(y)|} &\lesssim {( \gamma_{m,q}(y))^{1- N / \Nfin}} \, , \label{eq:DN:psi:q}
    \end{align}
    while for $\frac{1}{4}\Gamma_{q}^{2(m+1)}<y<\Gamma_{q}^{2(m+1)}$ we have
    \begin{align}
     {|D^N  \gamma_{m,q}(y)|} &\lesssim \Gamma_{q}^{-2N(m+1)} {( \gamma_{m,q}(y))^{1- N / \Nfin}} \, . \label{eq:DN:psi:q:gain}
    \end{align}
In each of the above inequalities, the implicit constants depend on $N$ but not $m$ or $q$. If $\gamma_{m,q}$ or $\tilde\gamma_{m,q}$ is replaced on the left hand side with $\gamma_{m,q}^p$, respectively $\tilde\gamma_{m,q}^p$ for $p\in \mathbb{N}$, then a similar inequality holds after substituting the same power on the right-hand side and changing implicit constants.
\end{enumerate}
\end{lemma}

We now introduce the intermittent pressure cut-off functions.\index{$\omega_{j,q}$}
\begin{definition}[\bf Intermittent pressure cutoff functions]\label{def:pressure:cutoff}
For $j\geq 1$ the cut-off functions are defined by\index{$j$}\index{$\omega_{j,q}$}
\begin{align}
\omega_{j,q}(x,t) = \gamma_{0} \Big( \Gamma_{q}^{-2 j} \, (\de_{q+\bn})^{-1}\pi_\ell(x,t) \Big)
\,,
\label{eq:omega:cut:def}
\end{align}
while for $j=0$ we let 
\begin{align}
\omega_{0,q}(x,t) = \tilde\gamma_{0}\Big((\de_{q+\bn})^{-1}\pi_\ell(x,t)\Big)
\,,
\label{eq:omega:cut:def:0}
\end{align}
where $\gamma_0 := \gamma_{0,q}$ and $\tilde \gamma_0 := \td \gamma_{0,q}$.
\end{definition}
\noindent An immediate consequence of \eqref{eq:tilde:partition} with $m=0$ is that $\{\omega_{j,q}^{{6}} \}_{j\geq 0}$ satisfies
\begin{align}
\sum_{j\geq 0} \omega_{j,q}^{{6}} = 1 \, , \qquad \omega_{j,q} \omega_{j',q} \equiv 0 \quad \textnormal{if} \quad |j-j'| >1 
\label{eq:omega:cut:partition:unity}
\end{align}
on $\T^3 \times \R$. 

\subsubsection{Estimates for intermittent pressure cutoffs}
\begin{lemma}[\bf Simple derivative bounds]
\label{lem:D:Dt:Rn:sharp}
For all $m+k\leq \Nfin$ and $j\geq 0$, we have that
\begin{subequations}
\begin{align}
{\bf 1}_{\supp (\omega_{j,q}\psi_{i,q})} | D^k D_{t,q}^m \pi_\ell (x,t) | 
&\leq \Ga_q^{2j+6}\de_{q+\bn}
(\Ga_q\La_q)^k \MM{m, \NindRt, \Gamma_{q}^{i} \tau_{q}^{-1} , \Tau_{q}^{-1} } \, ,
\label{pt.est.pi.lem}
\\
\sfrac 14 \delta_{q+\bn} \Ga_q^{2j} &\leq {\bf 1}_{\supp (\omega_{j,q})} \pi_\ell \label{pt.est.pi.lem.lower} \\
\sfrac 18 \sum_j \omega_{j,q} \delta_{q+\bn} \Gamma_q^{2j} &\leq \pi_\ell \, ,
\label{pt.est.pi.cutoff}\\
{\bf 1}_{\supp (\omega_{j,q}\psi_{i,q})} | D^k D_{t,q}^m R_\ell (x,t) | 
&\leq \Gamma_{q}^{2j-4} \delta_{q+\bn} (\Ga_q\La_q)^k \MM{m, \NindRt, \Gamma_{q}^{i} \tau_{q}^{-1} , \Tau_{q}^{-1} } \, .
\label{pt.est.R.lem}
\end{align}
\end{subequations}
\end{lemma}
\begin{proof}
First, observe that by the construction of $\om_{j,q}$, we have that for all $j\geq 0$,
\begin{align}\label{pt.est.pi.0}
    {\bf 1}_{\supp (\omega_{j,q})} |\pi_\ell|
    ={\bf 1}_{\supp (\omega_{j,q})} \pi_\ell
    \leq \Ga_q^{2(j+1)}\de_{q+\bn}\,. 
\end{align}
Then, recalling the pointwise estimate \eqref{eq:pressure:inductive:dtq:pointwise} and using \eqref{pt.est.pi.0}, we have that
\begin{align*}
 {\bf 1}_{\supp (\omega_{j,q})}|\psi_{i,q} D^k D_{t,q}^m \pi_\ell (x,t) | 
&\lesssim {\bf 1}_{\supp (\omega_{j,q})} \Gamma_q^3 \pi_\ell (\Ga_q\La_q)^k \MM{m, \NindRt, \Gamma_{q}^{i} \tau_{q}^{-1} , \Tau_{q}^{-1} }\\
&\leq \Ga_q^{2(j+3)}\de_{q+\bn}
(\Ga_q\La_q)^k \MM{m, \NindRt, \Gamma_{q}^{i} \tau_{q}^{-1} , \Tau_{q}^{-1} } \, .
\end{align*}
To obtain the lower bounds on $\pi_\ell$ on the support of $\omega_{j,q}$, we appeal to \eqref{ind:pi:lower} in the case $j=0$ and the definition of $\omega_{j,q}$ in the case $j\geq 1$. Summing over $j$ and appealing to \eqref{eq:omega:cut:partition:unity} yields \eqref{pt.est.pi.cutoff}. Next, we can obtain the pointwise estimates \eqref{pt.est.R.lem} for $R_q^q$ in a similar way by using \eqref{eq:inductive:pointwise:upgraded:1}. Finally, we obtain \eqref{pt.est.pi.cutoff} from \eqref{ind:pi:lower}, the definition of $\omega_{j,q}$ for $j\geq 0$.
\end{proof}

\begin{corollary}[\bf Higher derivative bounds]
\label{cor:D:Dt:Rn:sharp:new}
For $q \geq 0$, $0\leq i\leq \imax$, and $\alpha,\beta \in \N_0^k$ with $|\alpha| + |\beta| \leq \Nfin$, we have 
\begin{subequations}
\begin{align}
&\norm{\left(\prod_{\ell=1}^{k} D^{\alpha_\ell} D_{t,q}^{\beta_\ell}\right) \pi_\ell}_{L^\infty(\supp (\psi_{i,q} \omega_{j,q}))} 
 \les \Gamma_{q}^{2j+6 } \delta_{q+\bn}  
(\Gamma_{q} \Lambda_q)^{|\alpha|}  
\MM{|\beta|, \NindRt,\Gamma_{q}^{i} \tau_q^{-1}, {\Tau}_q^{-1}}
\label{eq:D:Dt:pi:sharp:new}
\\
&\norm{\left(\prod_{\ell=1}^{k} D^{\alpha_\ell} D_{t,q}^{\beta_\ell}\right) R_\ell}_{L^\infty(\supp (\psi_{i,q} \omega_{j,q}))} 
 \les \Gamma_{q}^{2j-4 } \delta_{q+\bn}  
(\Gamma_{q} \Lambda_q)^{|\alpha|}  
\MM{|\beta|, \NindRt,\Gamma_{q}^{i} \tau_q^{-1}, {\Tau}_q^{-1}} \, .
\label{eq:D:Dt:R:sharp:new}
\end{align}
\end{subequations}
\end{corollary}
\begin{proof}[Proof of Corollary~\ref{cor:D:Dt:Rn:sharp:new}]
We only work on the estimate for $\pi_\ell$ because the estimates for $R_q^q$ can be obtained in a completely analogous way from Lemma~\ref{lem:D:Dt:Rn:sharp} and Lemma~\ref{lem:cooper:2}, Remark~\ref{rem:cooper:2:sum}. We then apply Lemma~\ref{lem:cooper:2} with $v = \hat u_{q}$, $f = \pi_\ell$, $\Omega = \supp \psi_{i,q} \cap \supp \omega_{j,q}$, and $p=\infty$. In view of estimate \eqref{eq:nasty:D:vq:old} at level $q$, the assumption \eqref{eq:cooper:2:v} holds with $\const_v = \tau_q^{-1}\Gamma_q^i\Lambda_q^{-1}$, $\lambda_v = \tilde \lambda_v = \Lambda_q$, $N_x = \infty$, $\mu_v = \Gamma_{q}^{i} \tau_q^{-1}$, $\tilde \mu_v = \Gamma_{q}^{-1} \Tau_q^{-1}$, and $N_t = \Nindt$. On the other hand, the bound \eqref{pt.est.pi.lem} implies assumption~\eqref{eq:cooper:2:f} with $\const_f = \Gamma_{q+1}^{2j+6} \delta_{q+\bn}$, $\lambda_f =  \tilde \lambda_f = \Gamma_q \Lambda_q$,  $\mu_f = \Gamma_{q}^{i} \tau_q^{-1}$, $\tilde \mu_f = \Tau_q^{-1}$, and $N_t = \Nindt$. We then deduce from the bound \eqref{eq:cooper:2:f:2} that \eqref{eq:D:Dt:pi:sharp:new} holds, thereby concluding the proof.
\end{proof}

\begin{lemma*}[\bf Current error estimates]\label{curr:press:lemma}
For all $m+k\leq \Nfin$ and $j\geq 0$, we have that
\begin{equation}
{\bf 1}_{\supp (\omega_{j,q}\psi_{i,q})} | D^k D_{t,q}^m \ph_\ell (x,t) | 
\leq \Gamma_{q}^{3j-7} \delta_{q+\bn}^\frac32 r_q^{-1} (\Ga_q\La_q)^k \MM{m, \NindRt, \Gamma_{q}^{i} \tau_{q}^{-1} , \Tau_{q}^{-1} }
\,.
\label{pt.est.ph.lem}
\end{equation}
For $q\geq 0$, $0\leq i \leq \imax$, and $\alpha,\beta\in \mathbb{N}^k_0$ with $|\alpha|+|\beta|\leq \Nfin$, we have
\begin{equation}
    \\
\norm{\left(\prod_{\ell=1}^{k} D^{\alpha_\ell} D_{t,q}^{\beta_\ell}\right) \ph_\ell}_{L^\infty(\supp (\psi_{i,q} \omega_{j,q}))} 
 \les \Gamma_{q}^{3j-7 }  \delta_{q+\bn}^\frac32r_q^{-1}
(\Gamma_{q} \Lambda_q)^{|\alpha|}  
\MM{|\beta|, \NindRt,\Gamma_{q}^{i} \tau_q^{-1}, {\Tau}_q^{-1}} \, .
\label{eq:D:Dt:ph:sharp:new}
\end{equation}
\end{lemma*}
\begin{proof}
The proof is completely analagous to the proofs of Lemma~\ref{lem:D:Dt:Rn:sharp} and Corollary~\ref{cor:D:Dt:Rn:sharp:new}, and we omit the details.
\end{proof}

\begin{lemma}[\bf Maximal $j$ index]
\label{lem:maximal:j}
Fix $q\geq 0$. There exists a $\jmax = \jmax(q) \geq 1$, determined by the formula
\begin{align}
\jmax = \inf \left\{ j \, : \, \frac 14 \Ga_q^{2j} \delta_{q+\bn} \geq \Ga_q^{3+\badshaq} \right\}
\label{eq:j:max:def}
\end{align}
and which is bounded independently of $q$, such that\index{$\jmax$}
\begin{align}\label{eq:omega:j:is:zero}
\omega_{j,q} \equiv 0 \qquad \mbox{for all} \qquad j > \jmax \, .
\end{align}
Moreover, we have the bound
\begin{align}\label{ineq:jmax:use}
\Gamma_{q}^{2j_{\rm max}} \leq \de_{q+\bn}^{-1} \Gamma_q^{\badshaq+6} \, .
\end{align}
\end{lemma}

\begin{proof}[Proof of Lemma~\ref{lem:maximal:j}]
The proof of \eqref{eq:omega:j:is:zero} follows immediately from the definition in \eqref{eq:j:max:def}, the bound \eqref{pt.est.pi.lem}, and the bound \eqref{eq:pressure:inductive:dtq:uniform:upgraded}, where the extra factor of $\Ga_q$ absorbs the implicit constant in \eqref{eq:pressure:inductive:dtq:uniform:upgraded}. Checking that $\jmax$ is independent of $q$ is a simple calculation, as is the bound in \eqref{ineq:jmax:use}.
\end{proof}

\begin{lemma}[\bf Derivative bounds]
\label{lem:D:Dt:omega:sharp}
For $q\geq 0$, $0 \leq i \leq \imax$, $0 \leq j \leq \jmax$, and $N + M \leq \Nfin$, we have
\begin{align}
\frac{{\bf 1}_{\supp \psi_{i,q}} |D^N D_{t,q}^M \omega_{j,q}|}{\omega_{j,q}^{1-(N+M)/\Nfin}} 
\les (\Gamma_{q}^5 \Lambda_q)^N \MM{M, \Nindt,\Gamma_{q}^{i+4} \tau_{q}^{-1}, {\Tau}_{q}^{-1}} \, .
\label{eq:D:Dt:omega:sharp}
\end{align}
\end{lemma}

\begin{proof}[Proof of Lemma~\ref{lem:D:Dt:omega:sharp}]
We shall apply the mixed-derivative Fa'a di Bruno formula from \cite[Lemma~A.5]{BMNV21} with the following choices, where we use the parameter names from there:
\begin{align*}
    &\psi = \gamma_0 \textnormal{ or } \tilde\gamma_0 \, , \quad \Gamma_\psi = \Ga_q \, , \quad v = \hat u_q \, , \\
    &\Gamma= \delta_{q+\bn}^{\sfrac 12} \Ga_q^{-j} \, , \quad \lambda=\tilde\lambda=\La_q\Ga_q \, , \quad \mu = \tau_q^{-1}\Ga_q^i \, , \quad \tilde \mu = \Tau_q^{-1} \, , \\
    &N_x = \infty \, , \quad N_t = \Nindt \, , \quad h = \pi_\ell \, , \quad \const_h = \delta_{q+\bn}\Ga_q^{2j+6} \, .
\end{align*}
The assumption \cite[A.24]{BMNV21} is verified due to \eqref{eq:DN:psi:q:0}--\eqref{eq:DN:psi:q:gain}, and \cite[(A.25)]{BMNV21} is verified due to \eqref{eq:D:Dt:pi:sharp:new}, which holds on the support of $\omega_{j,q}\psi_{i,q}$. From conclusion \cite[(A.26)]{BMNV21} and the equality $(\Gamma_\psi\Gamma)^{-2}\const_h=\Ga_q^4$, we find that \eqref{eq:D:Dt:omega:sharp} holds; note that for the $N=M=0$ case, we just use the fact that $\omega_{j,q}\leq 1$ rather than incur the loss $\const_h\Ga^{-2}$ from \cite[(A.26)]{BMNV21}.
\end{proof}

\begin{lemma}[\bf Support bounds]\label{lem:omega:support}
For any $r \geq \sfrac 32$ and $0\leq j \leq \jmax$, we have that 
\begin{align}
\norm{\omega_{j,q}}_{L^r}  \lesssim  \Gamma_{q}^{\frac{3(1-j)}{r}} \, .
\label{eq:omega:support}
\end{align}
\end{lemma}
\begin{proof}[Proof of Lemma~\ref{lem:omega:support}]
We prove only the case $r=\sfrac 32$, at which point the remaining estimates follow from Lebesgue interpolation and the fact that $\omega_{j,q}\leq 1$ for all $j,q$. For $j=0,1$ the estimate is trivial from the pointwise bound for $\omega_{j,q}$, and so we consider now $j\geq 2$. Using Chebyshev's inequality, \eqref{eq:pressure:inductive:dtq:upgraded}, and \eqref{pt.est.pi.lem.lower}, we have that
\begin{align*}
\left\| \omega_{j,q} \right\|_{\sfrac 32}^{\sfrac 32} &\leq \sup_{t\in \mathbb{R}} \int_{\T^3} \mathbf{1}_{\{\pi_\ell(t,\cdot) \geq \sfrac 14 \delta_{q+\bn}\Gamma_q^{2j}\}} dx
\lesssim \frac{\left\| \pi_\ell \right\|_{\sfrac 32}^{\sfrac 32}}{\delta_{q+\bn}^{\sfrac 32}\Gamma_q^{3j}} 
\lesssim \Ga_q^{3(1-j)} \, .
\end{align*}
\end{proof}

\subsection{Mildly and strongly anisotropic checkerboard cutoffs}
\label{sec:cutoff:checkerboard:definitions}
We first construct mildly anisotropic\index{$\chi_{q,\xi,l/l^\perp}$} checkerboard cutoff functions which are well-suited for intermittent pipe flows with axes parallel to $e_1$.  The construction for general $\xi\in\Xi$ follows by rotation. We include all the details since the power for which the partition is summable to 1 is absolutely crucial for the definition of the perturbation in \eqref{eq:rqnpj} and its estimates in Lemma~\ref{lem:a_master_est_p}, and the Reynolds oscillation errors in subsections~\ref{ss:ssO}. These summability properties are also crucial in the estimates for the current oscillation errors in \cite[section~5.2]{GKN23}.

\noindent\texttt{Step 1: Partitioning the space perpendicular to $x_1$.}  Consider a partition of $\T^2_{x_2,x_3}$ into the squares defined using the periodized base square
\begin{equation}\label{eq:prism:zero}
    \left\{ (x_2,x_3)\in\T^2 \, : \, 0 \leq x_2,x_3 \leq \frac{\pi}{8} \Gamma_{q}^5\left(\lambda_{q+1} \right)^{-1} \right\}
\end{equation}
and its periodized translations by 
$$\bigl( l_2 \cdot \sfrac{\pi}{8} \cdot  \Gamma_{q}^5(\lambda_{q+1})^{-1}, l_3 \cdot \sfrac \pi 8 \cdot \Gamma_{q}^5(\lambda_{q+1})^{-1} \bigr) \, $$
for
$$ l_2,l_3\in\{0,\dots, 16\Gamma_{q}^{-5}\lambda_{q+1} -1\} \, .
$$
Note that the periodized squares evenly partition $[-\pi,\pi]^2$. We let $l^\perp := (l_2,l_3)$ be an ordered pair using the indices defined above, and choose $\{\mathcal{X}_{q,e_1,l^\perp}\}_{l^\perp}$ to be a $C^\infty$ partition of unity adapted to these periodized squares such that
\begin{subequations}
\begin{align}\label{eq:checkie:squared}
    \sum_{l^\perp} \mathcal{X}_{q,e_1,l^\perp}^2(x_2,x_3) &\equiv 1, \quad \forall (x_2,x_3) \in \T^2_{x_2,x_3} \, , \quad\,  \mathcal{X}_{q,e_1,l^\perp} \mathcal{X}_{q,e_1,\tilde l^\perp} \equiv 0 \quad  \textnormal{if $|l_2-\tilde l_2| >1$ $|l_3 - \tilde l_3| >1$} \, , \\
    \supp\mathcal{X}_{q,e_1,l^\perp_0}
    &= [-\sfrac18 \Ga_q^5 \la_{q+1}^{-1}, \sfrac58 \Ga_q^5 \la_{q+1}^{-1}]^2 \qquad \textnormal{for  $l^\perp_0=(0,0)$}  \, .
    \label{supp.X}
\end{align}
\end{subequations}
We shall later need that
\begin{equation}\label{eq:summy:summ:0}
    \left \langle \sum_{l^\perp} \chi_{q,e_1,l^\perp}^3(x_2,x_3) \right \rangle = c_3 \, ,
\end{equation}
where the constant $c_3$ is geometric and bounded independently of $q$.

\noindent\texttt{Step 2: Partitioning the space parallel to $x_1$. } Next, consider a partition of $\T_{x_1}$ into the line segments defined using the base line segment
\begin{equation}\label{eq:prism:zero:zero}
    \left\{ x_1 \in \T \, : \, 0 \leq x_1 \leq \frac{\pi}{8} \lambda_q^{-1} \Gamma_q^{-8} \right\}
\end{equation}
and its translations by 
$$ l \cdot \sfrac12 \cdot \lambda_q^{-1} \Gamma_q^{-8} \, , \qquad l \in \{0 , \dots , 16 \lambda_q^{-1}\Gamma_q^{-8}-1 \} \, .   $$
Note that the segments evenly partition $[-\pi,\pi]$. Choose $\{\mathcal{X}_{q,e_1,l}\}_{l}$ to be a $C^\infty$ partition of unity adapted to these segments such that for $N\leq 3\Nfin$,
\begin{subequations}
\begin{align}\label{eq:checkie:sixthed}
    \sum_{l} \mathcal{X}_{q,e_1,l}^6(x_1) &\equiv 1 \quad \forall (x_1) \in \T_{x_1} \, , \quad \mathcal{X}_{q,e_1,l} \mathcal{X}_{q,e_1,\tilde l} \equiv 0 \quad  \textnormal{ if } |l-\tilde l|>1 \, , \qquad \left| D^N \mathcal{X}_{q,\xi',l'} \right| \lesssim (\lambda_q \Gamma_q^8)^N \, , \\
    \supp(\mathcal{X}_{q,e_1,0}) &=
    [-\sfrac18 \la_q^{-1}\Ga_q^{-8}  , \sfrac58 \la_q^{-1}\Ga_q^{-8}  ]  \, . \label{supp.X:long}
\end{align}
\end{subequations}

\noindent\texttt{Step 3: Reynolds cutoffs.} Combining $l,l^\perp$ into integer triples $\vec{l}=(l,l_2,l_3)=(l,l^\perp)$, we now have a division of $\T^3$ into rectangular prisms indexed by $\vecl$.  We define 
\begin{align}\notag
    \mathcal{X}_{q,e_1,\vecl,R}(x_1,x_2,x_3) = \mathcal{X}_{q,e_1,l}^3(x_1) \mathcal{X}_{q,e_1,l^\perp}(x_2,x_3) 
\end{align}
and note that
\begin{align}\notag
    \sum_{\vecl} \mathcal{X}_{q,e_1,\vecl,R}^2 (x_1,x_2,x_3) \equiv 1 \qquad \forall \, (x_1,x_2,x_3) \in \T^3 \, .
\end{align}

\noindent\texttt{Step 4: Current cutoffs.} We combine $l,l^\perp$ into integer triples $\vec{l}$ as above but now define
\begin{align}\notag
    \mathcal{X}_{q,e_1,\vecl,\vp}(x_1,x_2,x_3) = \mathcal{X}_{q,e_1,l}^2(x_1) \mathcal{X}_{q,e_1,l^\perp}(x_2,x_3) 
\end{align}
and note that for each fixed value of $l=l_0$,
\begin{align}\notag
    \sum_{\vecl \,:\,l=l_0} \mathcal{X}_{q,e_1,\vecl,\vp}^2 (x_1,x_2,x_3) \equiv \mathcal{X}_{q,e_1,l_0}^4(x_1) \qquad \forall \, (x_1,x_2,x_3) \in \T^3 \, .
\end{align}
Conversely, for each fixed value of $l^\perp=l^\perp_0$, we have that
\begin{align}\notag
    \sum_{\vecl \,:\,l^\perp=l^\perp_0} \mathcal{X}_{q,e_1,\vecl,\vp}^3 (x_1,x_2,x_3) \equiv \mathcal{X}_{q,e_1,l^\perp_0}^3(x_2,x_3) \, .
\end{align}
With the time-independent cutoffs in hand, we define the time-dependent cutoff which is adapted to the flows of the velocity field $\hat u_q$.\index{$\zeta_{q,\diamond,i,k,\xi,\vecl}$}

\begin{definition}[\bf Mildly anisotropic checkerboard cutoff functions]\label{def:checkerboard}
Given $q$, {$\xi\in\Xi$}, $i\leq \imax$, and $k\in\mathbb{Z}$, we define
\begin{equation}\label{eq:checkerboard:definition}
    \zeta_{q,\diamond,i,k,{\xi},\vec{l}}\,(x,t) = \mathcal{X}_{q,\xi,\vec{l},\diamond}\left(\Phi_{i,k,q}(x,t)\right) \, .
\end{equation}
\end{definition}
These cutoff functions satisfy properties which we enumerate in the following lemma.

\begin{lemma}\label{lem:checkerboard:estimates}
The cutoff functions $\{\zeta_{q,\diamond,i,k,{\xi},\vec{l}}\}_{\vec{l}}$ satisfy the following properties.
\begin{enumerate}[(i)]
    \item\label{item:check:1} The material derivative $\Dtq (\zeta_{q,\diamond,i,k,{\xi},\vec{l}})$ vanishes.  
    \item\label{item:check:2} We have the summability properties for all $(x,t)\in \T^3\times \R$; 
    \begin{subequations}\label{eq:checkerboard:partition}
    \begin{align}
    \sum_{\vec{l}} \bigl(\zeta_{q,R,i,k,{\xi},\vec{l}}\,(x,t)\bigr)^{2} &\equiv 1 \, , \label{eq:summy:summ:1} \\
    \sum_{\vecl \, : \, l=l_0} \zeta_{q,\vp,i,k,\xi,\vecl}^2(x,t) &\equiv \mathcal{X}_{q,\xi,l_0}^4(\Phi_{i,k,q}(x,t)) \, , \label{eq:summy:summ:2} \\
    \sum_{\vecl \, : \, l^\perp = l^\perp_0} \zeta_{q,\vp,i,k,\xi,\vecl}^3(x_1,x_2,x_3) &= \mathcal{X}_{q,\xi,l_0^\perp}^3(\Phi_{i,k,q}(x,t)) \, . \label{eq:summy:summ:3}
    \end{align}
    \end{subequations}
    \item\label{item:check:3} 
    Let $A=(\nabla\Phi_{(i,k)})^{-1}$.  Then we have the spatial derivative estimate
    \begin{align}\label{eq:checkerboard:derivatives}
        \bigl\| D^{N_1} \Dtq^M  ({\xi^\ell A_\ell^j \partial_j} )^{N_2} \zeta_{q,\diamond,i,k,\xi,\vec{l}} \bigr\|_{L^\infty\left(\supp \psi_{i,q}\tilde\chi_{i,k,q} \right)} &\lesssim \left(\Gamma_{q}^{-5} \lambda_{q+1} \right)^{N_1} \left(\Gamma_q^8\lambda_{q}\right)^{N_2} \notag\\
        &\qquad \qquad \times \MM{M,\Nindt,\Gamma_{q}^{i}\tau_q^{-1},\Tau_q^{-1}\Gamma_{q}^{-1}} \, .
    \end{align}
    for all $N_1+N_2+M\leq \sfrac{3\Nfin}{2}+1$.
    \item\label{item:check:4} There exists an implicit dimensional constant $\const_\chi$ independent of $q$, $k$, $i$, and $\vec{l}$ such that for all $(x,t)\in\supp\psi_{i,q}\tilde\chi_{i,k,q}$, the support of $\zeta_{q,\diamond,i,k,\xi,\vec{l}}\, (\cdot,t)$ satisfies
    \begin{equation}\label{eq:checkerboard:support}
        \textnormal{diam} ( \supp ( \zeta_{q,\diamond,i,k,\xi,\vec{l}}\,(\cdot,t) ) ) \lesssim \Gamma_q^{-8}\lambda_{q}^{-1} \, . 
    \end{equation}
\end{enumerate}
\end{lemma}
\begin{proof}[Proof of Lemma~\ref{lem:checkerboard:estimates}]
The proof of \eqref{item:check:1} is immediate from \eqref{eq:checkerboard:definition}.  The first equality in \eqref{eq:checkerboard:partition} follows from \eqref{item:check:1} and the definition of the Reynolds cutoffs in Step 3 above. The second and third equalities follow from \eqref{item:check:1} and the definition of the current cutoffs in Step 4 above. To verify \eqref{item:check:3}, the only nontrivial calculations are those including the differential operator $\xi^\ell A_{\ell}^j\partial_j$. Using the Leibniz rule, the contraction
$$ \xi^\ell A_\ell^j \partial_j \zeta_{q,\diamond,i,k,\xi,\vec{l}} = \xi^\ell A_\ell^j (\partial_m \mathcal{X}_{q,\xi,\vec{l},\diamond})(\Phi_{i,k,q}) \partial_j \Phi_{i,k,q}^m= \xi^m (\partial_m \mathcal{X}_{q,\xi,\vec{l},\diamond})(\Phi_{i,k,q}) \, , $$
the diameter of the cutoffs defined in Steps 1 and 2 above, and \eqref{eq:Lagrangian:Jacobian:5}--\eqref{eq:Lagrangian:Jacobian:6} gives the desired estimate. The proof of \eqref{eq:checkerboard:support} follows from the construction of $\mathcal{X}_{q,\xi,\vec{l},\diamond}$ and the Lipschitz bound obeyed by $\hat u_q$ on the support of $\psi_{i,q}$; see for example \eqref{eq:diameter:inequality}.
\end{proof}

We may similarly obtain estimates on the flowed cutoff functions $\etab_\xi^I$ which come from Definition~\ref{def:etab}. The proof is quite similar to the one above, and we omit the details.\index{$\zetab_{\xi}^I$}
\begin{lemma}[\bf Strongly\index{$I$} anisotropic checkerboard cutoff function]\label{lem:finer:checkerboard:estimates}
The cutoff functions $\etab_\xi^I\circ\Phiik$ satisfy the following properties:
\begin{enumerate}[(1)]
    \item\label{item:check:check:1} The material derivative $\Dtq (\etab_\xi^I\circ \Phiik)$ vanishes.
    \item\label{item:check:check:2} For all fixed values of $q,i,k,\xi$, each $t\in\mathbb{R}$, and all $x=(x_1,x_2,x_3)\in\mathbb{T}^3$,
    \begin{equation}\label{eq:checkerboard:partition:check}
    \sum_{I} (\etab_\xi^I \circ \Phiik)^6(x,t) = 1 \, .
    \end{equation}
    \item\label{item:check:check:3} 
    Let $A=(\nabla\Phi_{(i,k)})^{-1}$.  Then we have the spatial derivative estimate
    \begin{align}\label{eq:checkerboard:derivatives:check}
        \bigl\| D^{N_1} \Dtq^M  ({\xi^\ell A_\ell^j \partial_j} )^{N_2} \etab_\xi^I \circ \Phiik \bigr\|_{L^\infty\left(\supp \psi_{i,q}\tilde\chi_{i,k,q} \right)} &\lesssim \lambda_{q+\lfloor \sfrac \bn 2 \rfloor}^{N_1}  \MM{M,\Nindt,\Gamma_{q}^{i}\tau_q^{-1},\Tau_q^{-1}\Gamma_{q}^{-1}} \, .
    \end{align}
    for all $N_1+N_2+M\leq \sfrac{3\Nfin}{2}+1$.
    \item\label{item:check:check:4} There exists an implicit dimensional constant $\const_\chi$ independent of $q$, $k$, $i$, and $\xi$ such that for all $(x,t)\in\supp\psi_{i,q}\tilde\chi_{i,k,q}$, the support of $\zetab_{\xi}^I \circ \Phiik(\cdot,t)$ satisfies
    \begin{equation}\label{eq:checkerboard:support:check}
        \textnormal{diam} ( \supp ( \zetab_{\xi}^I \circ \Phiik(\cdot,t) ) ) \lesssim \Gamma_q^{-8}\lambda_{q}^{-1} \, . 
    \end{equation}
\end{enumerate}
\end{lemma}
We also need the following lemma that bounds the cardinality of these anisotropic cut-offs.
\begin{lemma}\label{lem.cardinality}
For fixed $q,i,k,\xi$, we have that
\begin{align}
    \# \left\{(\vecl,I) : \supp \left(\zeta_{q,i,k,\xi,\vecl} \, \etab^I_{\xi} \circ \Phi_{(i,k)}\right) \neq \emptyset \right\} \lec \Ga_q^8 \la_{q} \la_{q+\half}^2 \,.
\end{align}
\end{lemma}
\begin{proof}
Note first that for a fixed $I$, there are at most $4$ values of $l_0^\perp$ such that $\supp (\mathcal{X}_{q,\xi,l_0^\perp} \etab^I_{\xi}) \neq \emptyset$. Also note that for a fixed $l_0^\perp$, we have $\# \{\vecl : l^\perp = l_0^\perp \} \lec \la_q \Ga_q^8$. Putting these together along with the bound on the number of $I$ given by Remark~\ref{rem:strong:cardinality}, we get that
$$ \# \{ (\vecl,I) : \supp(\mathcal{X}_{q,\xi,\vecl,\diamond} \etab^I_{\xi}) \neq \emptyset \} \lec \Ga_q^8 \la_{q} \la_{q+\half}^2 \,. $$
Now the desired conclusion follows as all these cut-offs are flowed by the same $\Phi_{(i,k)}$.
\end{proof}

\subsection{Definition of the cumulative cutoff function}\label{sec:cutoff:total:definitions}
Finally, combining the cutoff functions defined in subsection~\ref{sec:cutoff:inductive}, Definition~\ref{def:psi:i:q:def}, Definition~\ref{def:pressure:cutoff}, \eqref{eq:chi:cut:def}, and the previous subsection, we define the cumulative cutoff functions by\index{$\eta_{i,j,k,\xi,\vecl,\diamond}$}\index{cumulative cutoff function}
\begin{align}
 \eta_{i,j,k,\xi,\vec{l},\diamond}\,(x,t) &= \psi_{i,q}^\diamond(x,t) \omega_{j,q}^\diamond(x,t) \chi_{i,k,q}^\diamond(t)\zeta_{q,\diamond,i,k,\xi,\vec{l}}\,(x,t) \label{def:cumulative:current}
  \, , 
\end{align}
where the $\diamond$ in the superscript of the first three functions is equal to $2$ if $\diamond=\varphi$ (so that they are cubic-summable to $1$) and $3$ if $\diamond=R$ (so that they are square-summable to $1$). 
We conclude this section with estimates on the $L^p$ norms of the cumulative cutoff functions. 

\begin{lemma}[\bf Cumulative support bounds for cutoff functions]\label{lemma:cumulative:cutoff:Lp}
For $r_1,r_2\in [1,\infty]$ with $\frac{1}{r_1}+\frac{1}{r_2}=1$ and any $0\leq i \leq \imax$, $0\leq j,\leq \jmax$, $\xi\in\Xi,\Xi'$, and $\diamond=\varphi,R$, we have that for each $t$,
\begin{align}
    \sum_{\vecl} \left| \supp_x \left( \eta_{i,j,k,\xi,\vecl,\diamond}(t,x) \right) \right| &\lesssim \Gamma_{q}^{\frac{-3i + \CLebesgue}{r_1} + \frac{-3j}{r_2}+3} \, . \label{eq:supp:cumul:varphi}
\end{align}
We furthermore have that
\begin{align}
    \sum_{i,j,k,\xi,\vecl,I,\diamond} \mathbf{1}_{\supp \eta_{i,j,k,\xi,\vecl,\diamond} \rhob_\pxi^\diamond \zetab_\xi^I} \approx \sum_{i,j,k,\xi,\vecl,\diamond} \mathbf{1}_{\supp \eta_{i,j,k,\xi,\vecl,\diamond} \rhob_\pxi^\diamond} \lesssim 1 \, . \label{eq:desert:cowboy:sum}
\end{align}
\end{lemma}
\begin{proof}[Proof of Lemma~\ref{lemma:cumulative:cutoff:Lp}]
We shall prove the first bound for $\diamond=\varphi$. Then from \eqref{def:cumulative:current}, the only differences between $\diamond=R$ and $\diamond=\varphi$ are the powers to which various cutoff functions are raised, and so we shall omit the proof for $\diamond=R$. To prove the bound for $\diamond=\varphi$, we have that
\begin{align*}
    \sum_{\vecl} \left| \supp \eta_{i,j,k,\xi,\vecl,\varphi} \right| 
    &\lesssim 
    \left\| (\psi_{i-1,q}^6 + \psi_{i,q}^6 + \psi_{i+1,q}^6)^{\sfrac{1}{6}} (\omega_{j-1,q}^6 + \omega_{j,q}^6 + \omega_{j+1,q}^6)^{\sfrac{1}{6}}  \right\|_{L^1} \notag\\
     &\lesssim 
    \left\| (\psi_{i-1,q}^6 + \psi_{i,q}^6 + \psi_{i+1,q}^6)^{\sfrac{1}{6}} \right\|_{L^{r_1}} \left\| (\omega_{j-1,q}^6 + \omega_{j,q}^6 + \omega_{j+1,q}^6)^{\sfrac{1}{6}}  \right\|_{L^{r_2}} \notag\\
    &\lesssim  \Gamma_{q}^{\frac{-3(i-1)+\CLebesgue}{r_1}} \Gamma_{q}^{\frac{-3(j-1)}{r_2}} \,.
\end{align*}
To achieve the final inequality, we have used interpolation, \eqref{eq:psi:i:q:support:old} at level $q$, and \eqref{eq:omega:support}. Using that $\frac{1}{r_1}+\frac{1}{r_2}=1$ gives the desired estimate.  Finally, to prove \eqref{eq:desert:cowboy:sum}, we appeal to \eqref{eq:inductive:partition} at level $q$, \eqref{eq:chi:cut:partition:unity} and \eqref{e:chi:overlap}, \eqref{eq:omega:cut:partition:unity}, item~\eqref{i:bundling:2} from Proposition~\ref{prop:bundling}, Definition~\ref{def:etab}, and Lemma~\ref{lem:checkerboard:estimates}.
\end{proof}

\subsection{Cutoff aggregation lemmas}

\begin{corollary}[\bf Aggregated $L^p$ estimates]\label{rem:summing:partition}
Let $\theta\in{(}0,3]$, and $\theta_1,\theta_2\geq 0$ with $\theta_1+\theta_2=\theta$. Let $H=H_{i,j,k,\xi,\vecl,\diamond}$ or $H=H_{i,j,k,\xi,\vecl,I,\diamond}$ be a function with\index{aggregation lemmas}
\begin{align}
    \supp H_{i,j,k,\xi,\vecl,\diamond} \subseteq \supp \eta_{i,j,k,\xi,\vecl,\diamond} \qquad \textnormal{or} \qquad \supp H_{i,j,k,\xi,\vecl,I,\diamond} \subseteq \supp  \eta_{i,j,k,\xi,\vecl,\diamond} \zetab_{\xi}^{I,\diamond}{\circ\Phiik} \, . \label{eq:agg:assump:1}
\end{align}
Let $p\in[1,{\infty)}$ and let $\theta_1,\theta_2\in[0,3]$ be such that $\theta_1+\theta_2=\sfrac 3p$. Assume that there exists $\const_H,N_*,M_*,N_x,M_t$ and $\lambda,\Lambda,\tau,\Tau$ such that
\begin{subequations}
\begin{align}
  \left\| D^N \Dtq^M H_{i,j,k,\xi,\vecl,\diamond} \right\|_{L^p} &\lesssim \sup_{t\in\R} \left( \left| \supp_x \left( \eta_{i,j,k,\xi,\vecl,\diamond} (t,x) \right) \right|^{\sfrac 1p} \right) \notag\\
  &\qquad \qquad \times \const_H \Gamma_q^{\theta_1 i + \theta_2 j}  \MM{N,N_x,\lambda,\Lambda} \MM{M,M_t,\tau^{-1}\Gamma_q^i,\Tau^{-1}} \label{eq:agg:assump:2} \\
  \left\| D^N \Dtq^M H_{i,j,k,\xi,\vecl,I,\diamond} \right\|_{L^p} &\lesssim \sup_{t\in\R} \left( \left| \supp_x \left( \eta_{i,j,k,\xi,\vecl,\diamond} \zetab_\xi^{I,\diamond}\circ\Phiik  (t,x)\right) \right|^{\sfrac 1p} \right) \notag\\
  &\qquad \qquad \times \const_H \Gamma_q^{\theta_1 i + \theta_2 j} \MM{N,N_x,\lambda,\Lambda} \MM{M,M_t,\tau^{-1}\Gamma_q^i,\Tau^{-1}} \,  \label{eq:agg:assump:3}
\end{align}
\end{subequations}
for $N\leq N_*,M\leq M_*$. Then in the same range of $N$ and $M$,
\begin{subequations}
\begin{align}\label{eq:agg:conc:1}
 \left\| \psi_{i,q} \sum_{i',j,k,\xi,\vecl,\diamond} D^N \Dtq^M H_{i',j,k,\xi,\vecl,\diamond} \right\|_{L^p} &\lesssim \Gamma_q^{{3+\theta_1 \CLebesgue}} \const_H \MM{N,N_x,\lambda,\Lambda} \MM{M,M_t,\tau^{-1}\Gamma_q^{i+1},\Tau^{-1}} \\
 \label{eq:agg:conc:2}
 \left\| \psi_{i,q} \sum_{i',j,k,\xi,\vecl,I,\diamond} D^N \Dtq^M H_{i',j,k,\xi,\vecl,I,\diamond} \right\|_{L^p} &\lesssim \Gamma_q^{{3+\theta_1 \CLebesgue}} \const_H \MM{N,N_x,\lambda,\Lambda} \MM{M,M_t,\tau^{-1}\Gamma_q^{i+1},\Tau^{-1}} \, .
\end{align}
\end{subequations}
\end{corollary}
\begin{proof}
We prove only \eqref{eq:agg:conc:2}, as \eqref{eq:agg:conc:1} is slightly easier and follows the same method. Using \eqref{eq:agg:assump:1}, \eqref{eq:inductive:partition} at level $q$, \eqref{eq:agg:assump:3}, Lemma~\ref{lemma:cumulative:cutoff:Lp} with $r_1=\frac{3}{p\theta_1},r_2=\frac{3}{p\theta_2}$, $\theta_1+\theta_2=\sfrac 3p$, we may write that
\begin{align}
    &\left\| \psi_{i,q} \sum_{i',j,k,\xi,\vecl,I,\diamond} D^N \Dtq^M H_{i',j,k,\xi,\vecl,I,\diamond} \right\|_p^p \notag\\
    &\quad \leq \sup_{t\in\R} \int_{\T^3} \psi_{i,q} \left| \sum_{\substack{i-1\leq i' \leq i+1 \\ j,k,\xi,\vecl,I,\diamond }} D^N \Dtq^M H_{i',j,k,\xi,\vecl,I,\diamond} \right|^p(t,x) \, dx \notag\\
    &\quad \leq \sup_{t\in\R} \sum_{\substack{i-1\leq i' \leq i+1 \\ j,k,\xi,\vecl,I,\diamond }} \left| \supp_x \left( \eta_{i,j,k,\xi,\vecl,\diamond} \zetab_\xi^{I,\diamond} \circ\Phiik(t,x) \right) \right| \const_H^p \Gamma_q^{p\theta_1 i + p\theta_2 j} \notag\\
    &\qquad \qquad \times \left(\MM{N,N_x,\lambda,\Lambda} \MM{M,M_t,\tau^{-1}\Gamma_q^i,\Tau^{-1}}\right)^p \notag\\
    &\quad \lesssim \sup_{t\in\R} \sum_{\substack{i-1\leq i' \leq i+1 \\ j,k,\xi,\vecl,\diamond }} \left| \supp_x \left( \eta_{i,j,k,\xi,\vecl,\diamond} (t,x) \right) \right| \const_H^p \Gamma_q^{p\theta_1 i + p\theta_2 j} \left(\MM{N,N_x,\lambda,\Lambda} \MM{N,N_t,\tau^{-1}\Gamma_q^i,\Tau^{-1}}\right)^p \notag\\
    &\quad\leq \const_H^p \Gamma_q^{{p \theta_1}\CLebesgue+3p} \left(\MM{N,N_x,\lambda,\Lambda} \MM{M,M_t,\tau^{-1}\Gamma_q^i,\Tau^{-1}}\right)^p \, , \notag 
\end{align}
concluding the proof.
\end{proof}

\begin{remark}[\bf Aggregated  $L^1$ estimates with $\Ga_q^i$]\label{rem:agg:current}
Assume that \eqref{eq:agg:assump:1}--\eqref{eq:agg:assump:3} hold for $p=\sfrac32$, but with $\const_H = \Ga_q^i \td\const_H$. Then we can obtain the $L^1$ estimates
\begin{subequations}
\begin{align}\label{eq:agg:conc:1:rem}
 \left\| \psi_{i,q} \sum_{i',j,k,\xi,\vecl,\diamond} D^N \Dtq^M H_{i',j,k,\xi,\vecl,\diamond} \right\|_{1} &\lesssim \td\const_H \Gamma_q^{2\CLebesgue+3} \MM{N,N_x,\lambda,\Lambda} \MM{M,M_t,\tau^{-1}\Gamma_q^i,\Tau^{-1}} \\
 \label{eq:agg:conc:2:rem}
 \left\| \psi_{i,q} \sum_{i',j,k,\xi,\vecl,I,\diamond} D^N \Dtq^M H_{i',j,k,\xi,\vecl,I,\diamond} \right\|_{1} &\lesssim \td\const_H \Gamma_q^{2\CLebesgue+3} \MM{N,N_x,\lambda,\Lambda} \MM{M,M_t,\tau^{-1}\Gamma_q^i,\Tau^{-1}}\, .
\end{align}
\end{subequations}
Indeed, considering \eqref{eq:agg:conc:2:rem}, we have
\begin{align*}
    &\left\| \psi_{i,q} \sum_{i',j,k,\xi,\vecl,I,\diamond} D^N \Dtq^M H_{i',j,k,\xi,\vecl,I,\diamond} \right\|_1\\
    &\quad \leq \sup_{t\in\R}\sum_{\substack{i-1\leq i' \leq i+1 \\ j,k,\xi,\vecl,I,\diamond }} \int_{\T^3} \psi_{i,q} 
    \textbf{1}_{\supp_x \left( \eta_{i,j,k,\xi,\vecl,\diamond} \zetab_\xi^{I,\diamond}\circ\Phiik \right)}
    \left|  D^N \Dtq^M H_{i',j,k,\xi,\vecl,I,\diamond} \right|(t,x) \, dx \\
    &\quad \leq \sup_{t\in\R}
    \left[\sum_{\substack{i-1\leq i' \leq i+1 \\ j,k,\xi,\vecl,I,\diamond }} 
    \Ga_q^{3i}\norm{\psi_{i,q}\textbf{1}_{\supp_x \left( \eta_{i,j,k,\xi,\vecl,\diamond} \zetab_\xi^{I,\diamond} \circ\Phiik\right)}}_3^3\right]^{\sfrac13}
    \left[
    \sum_{\substack{i-1\leq i' \leq i+1 \\ j,k,\xi,\vecl,I,\diamond }}
    \Ga_q^{-\sfrac32 i}\norm{  D^N \Dtq^M H_{i',j,k,\xi,\vecl,I,\diamond}}_{\sfrac32}^{\sfrac32} \right]^{\sfrac23}\\
     &\quad \lesssim \sup_{t\in\R}
     \left[
     \sum_{ j,k,\xi,\vecl,\diamond } 
     \left| \supp_x \left( \eta_{i,j,k,\xi,\vecl,\diamond} (t,x) \right) \right| \Ga_q^{3i}
     \right]^{\sfrac13}
     \left[
     \sum_{ j,k,\xi,\vecl,\diamond } \left| \supp_x \left( \eta_{i,j,k,\xi,\vecl,\diamond} (t,x) \right) \right|  \Gamma_q^{\sfrac32(\theta_1 i + \theta_2 j)}\right]^{\sfrac23}\\
     &\hspace{5cm}\cdot
\td\const_H\MM{N,N_x,\lambda,\Lambda} \MM{M,M_t,\tau^{-1}\Gamma_q^i,\Tau^{-1}}  \\
    &\quad\leq \td\const_H \Gamma_q^{2\CLebesgue+3} \MM{N,N_x,\lambda,\Lambda} \MM{M,M_t,\tau^{-1}\Gamma_q^i,\Tau^{-1}}
    \, . 
\end{align*}
In the last inequality, we used Lemma~\ref{lemma:cumulative:cutoff:Lp} 
with $r_1=1, r_2=\infty$ and with $r_1=\frac{3}{p\theta_1},r_2=\frac{3}{p\theta_2}$, and $\theta_1+\theta_2=\sfrac 3p$.
\end{remark}

We now state two similar corollaries which allow us to aggregate pointwise estimates.

\begin{corollary}[\bf Aggregated pointwise estimates]\label{lem:agg.pt}
Let $H=H_{i,j,k,\xi,\vecl,\diamond}$ or $H=H_{i,j,k,\xi,\vecl,I,\diamond}$ be a function with 
\begin{align}
    \supp H_{i,j,k,\xi,\vecl,\diamond} \subseteq \supp \eta_{i,j,k,\xi,\vecl,\diamond} \qquad \textnormal{or} \qquad \supp H_{i,j,k,\xi,\vecl,I,\diamond} \subseteq \supp  \eta_{i,j,k,\xi,\vecl,\diamond} \zetab_{\xi}^{I,\diamond}\circ\Phiik \,  \label{eq:aggpt:assump:1}
\end{align}
and let $\varpi = \varpi_{i,j,k,\xi,\vecl,\diamond}$ or $\varpi = \theta_{i,j,k,\xi,\vecl,I,\diamond}$ be a {non-negative} function such that
\begin{align}
    \supp \varpi_{i,j,k,\xi,\vecl,\diamond} \subseteq \supp \eta_{i,j,k,\xi,\vecl,\diamond} \qquad \textnormal{or} \qquad \supp \varpi_{i,j,k,\xi,\vecl,I,\diamond} \subseteq \supp  \eta_{i,j,k,\xi,\vecl,\diamond} \zetab_{\xi}^{I,\diamond} \circ\Phiik\,  \label{eq:aggpt:assump:2}
\end{align}
Let $p\in({0},\infty)$ and assume that there exists $\lambda,\Lambda,\tau$ such that
\begin{subequations}
\begin{align}
  |D^N \Dtq H_{i,j,k,\xi,\vecl,\diamond}| &\lesssim \varpi_{i,j,k,\xi,\vecl,\diamond}^p \MM{N,N_x,\lambda,\Lambda} \MM{N,N_t,\tau^{-1}\Gamma_q^i,\Tau^{-1}} \, \label{eq:aggpt:assump:3} \\
  |D^N \Dtq H_{i,j,k,\xi,\vecl,I,\diamond}| &\lesssim \varpi_{i,j,k,\xi,\vecl,I,\diamond}^p \MM{N,N_x,\lambda,\Lambda} \MM{N,N_t,\tau^{-1}\Gamma_q^i,\Tau^{-1}} \, \label{eq:aggpt:assump:4} 
\end{align}
\end{subequations}
for $N\leq N_*,M\leq M_*$. Then in the same range of $N$ and $M$,
\begin{subequations}
\begin{align}\label{eq:aggpt:conc:1}
 \left|\psi_{i,q} \sum_{i',j,k,\xi,\vecl,\diamond} D^N \Dtq^M H_{i',j,k,\xi,\vecl,\diamond} \right| &\lesssim \left(\sum_{i,j,k,\xi,\vecl,\diamond}\varpi_{i,j,k,\xi,\vecl,\diamond}\right)^p \MM{N,N_x,\lambda,\Lambda} \MM{M,M_t,\tau^{-1}\Gamma_q^{i+1},\Tau^{-1}} \\
 \label{eq:aggpt:conc:2}
 \left| \psi_{i,q} \sum_{i',j,k,\xi,\vecl,I,\diamond} D^N \Dtq^M H_{i',j,k,\xi,\vecl,I,\diamond} \right| &\lesssim \left(\sum_{i,j,k,\xi,\vecl,I,\diamond}\varpi_{i,j,k,\xi,\vecl,I,\diamond}\right)^p \MM{N,N_x,\lambda,\Lambda} \MM{M,M_t,\tau^{-1}\Gamma_q^{i+1},\Tau^{-1}} \, .
\end{align}
\end{subequations}
\end{corollary}

\begin{corollary*}[\bf Aggregated pointwise estimates with $\Ga_q^i$]\label{lem:agg.Dtq}
Let $H=H_{i,j,k,\xi,\vecl,I,\diamond}$ be a function with 
\begin{align}
     \supp H_{i,j,k,\xi,\vecl,\diamond} \subseteq \supp \eta_{i,j,k,\xi,\vecl,\diamond} \qquad \textnormal{or} \qquad \supp H_{i,j,k,\xi,\vecl,I,\diamond} \subseteq \supp  \eta_{i,j,k,\xi,\vecl,\diamond} \zetab_{\xi}^{I,\diamond}\circ\Phiik \,  \label{eq:aggDtq:assump:1}
\end{align}
and let $\varpi$ be a {non-negative} function and assume that there exists $\lambda,\Lambda,\tau,\Tau$ such that for $H = H_{i,j,k,\xi,\vecl,\diamond}$ or $H_{i,j,k,\xi,\vecl,I,\diamond}$
\begin{subequations}
\begin{align}
  \left|D^N \Dtq^M H\right| &\lesssim \tau_q^{-1} \Ga_q^i \psi_{i,q} \varpi \MM{N,N_x,\lambda,\Lambda} \MM{M,M_t,\tau^{-1}\Gamma_q^i,\Tau^{-1}} \, \label{eq:aggDtq:assump:4} 
\end{align}
\end{subequations}
for $N\leq N_*,M\leq M_*$. Then in the same range of $N$ and $M$,
\begin{subequations}
\begin{align}\label{eq:aggDtq:conc:1.0}
 \left|\psi_{i,q} \sum_{i',j,k,\xi,\vecl,\diamond} D^N \Dtq^M H_{i',j,k,\xi,\vecl,\diamond} \right| &\lesssim \Ga_q r_{q}^{-1} \la_q \left(\pi_q^q\right)^{\sfrac12} \varpi \MM{N,N_x,\lambda,\Lambda} \MM{M,M_t,\tau^{-1}\Gamma_q^{i+1},\Tau^{-1}} \\
\label{eq:aggDtq:conc:1}
 \left| \psi_{i,q} \sum_{i',j,k,\xi,\vecl,I,\diamond} D^N \Dtq^M H_{i',j,k,\xi,\vecl,I,\diamond} \right| &\lesssim \Ga_q r_{q}^{-1} \la_q \left(\pi_q^q\right)^{\sfrac12} \varpi \MM{N,N_x,\lambda,\Lambda} \MM{M,M_t,\tau^{-1}\Gamma_q^{i+1},\Tau^{-1}} \, .
\end{align}
\end{subequations}
\end{corollary*}
\begin{proof}[Proofs of Corollaries~\ref{lem:agg.pt} and \ref{lem:agg.Dtq}]
We will give the full details for estimate \eqref{eq:aggDtq:conc:1} from Corollary~\ref{lem:agg.Dtq}, since the proofs of all the other estimates are slightly easier and follow the same method.  We first note that summing the estimate in \eqref{eq:aggDtq:assump:4} over $j,k,\xi,\vecl,I,\diamond$ and using \eqref{eq:omega:cut:partition:unity}, \eqref{e:chi:overlap}, \eqref{eq:checkie:squared}, \eqref{eq:checkie:sixthed}, and \eqref{eq:sa:summability}, we find that
\begin{align*}
    \left|\sum_{j,k,\xi,\vecl,I,\diamond} D^N \Dtq^M H_{i,j,k,\xi,\vecl,I,\diamond}\right| &\lesssim \psi_{i\pm,q} \tau_q^{-1} \Ga_q^i \varpi  \MM{N,N_x,\lambda,\Lambda} \MM{M,M_t,\tau^{-1}\Gamma_q^i,\Tau^{-1}} 
\end{align*}
since $\supp H_{i,j,k,\xi,\vecl,I,\diamond} \subseteq \supp  \eta_{i,j,k,\xi,\vecl,\diamond} \zetab_{\xi}^{I,\diamond}\circ\Phiik \subseteq \supp \psi_{i,q}$ and $\psi_{i\pm,q}=(\psi_{i-1,q}^6 + \psi_{i,q}^6 + \psi_{i+1,q}^6)^{\sfrac 16}$. Now summing on $i$ and using \eqref{eq:inductive:partition} and Remark~\ref{rem.summing.psi}, we find that
\begin{align*}
    \left|\psi_{i,q} \sum_{i',j,k,\xi,\vecl,I,\diamond} D^N \Dtq^M H_{i',j,k,\xi,\vecl,I,\diamond} \right|
    &\lec  \left(\sum_i \Ga_q^{i} \tau_q^{-1} \psi_{i\pm,q}\right) \varpi \MM{N,N_x,\lambda,\Lambda} \MM{M,M_t,\tau^{-1}\Gamma_q^i,\Tau^{-1}}\\
    &\lec  \Ga_q r_{q}^{-1} (\pi_q^q)^{\sfrac12}\la_q \varpi \MM{N,N_x,\lambda,\Lambda} \MM{M,M_t,\tau^{-1}\Gamma_q^{i+1},\Tau^{-1}} \, .
\end{align*}
\end{proof}

\section{Velocity increment}\label{s.corr-ec-tor}
In this section, we define and estimate the velocity increment.  The first subsection contains the definition of $w_{q+1}$, save for the choice of placements of the bundles (see Remark~\ref{rem:sat:sat:sat}), which is addressed in the second subsection.  The final subsection then estimates both the pre-mollified velocity increment $w_{q+1}$ and the mollified velocity increment $\hat w_\qbn$.

\subsection{Definition of the corrector}\label{ss:corr-ec-tor}
In this subsection, we define the premollified velocity increment $w_{q+1}$, \emph{except} for the choice of placement, which we treat in \cite[section~4]{GKN23}; see also the discussion following Lemma~\ref{lem:dodging}. None of the discussion or properties in this subsection depend on the choice of placement.

\opsubsubsection{Definition of the current corrector}\label{ss:current corrector:definition}

For any fixed values of $i$, $k$, we recall the constant $c_3$ from \eqref{eq:summy:summ:0} and define
\begin{equation}\label{eq:phiqnpj}
\ph_{q,i,k} = {\sfrac{-1}{c_3}}\nabla\Phi_{(i,k)}\ph_\ell
\,.
\end{equation}
Let $\xi\in\Xi'$, cf. Proposition~\ref{prop:geo2}. For all $\xi\in\Xi'$, we define the coefficient function $a_{\xi,i,j,k,\vecl,\varphi}$ by\index{$a_{\xi,i,j,k,\vecl,\varphi}$}\index{$a_{\pxi,\diamond}$}
\begin{equation}
a_{\xi,i,j,k,\vecl,\ph} = a_{(\xi),\ph} = 
\delta_{q+\bn}^{\sfrac 12}r_q^{-\sfrac13} \Gamma^{j-1}_{q} \psi_{i,q}^\ph \omega_{j,q}^{\varphi} \chi_{i,k,q}^\ph \zeta_{q,\varphi,i,k,\xi,\vecl}\,
|\na \Phi_{(i,k)}^{-1} \xi|^{-\sfrac23}
\td\gamma_{\xi}\left(\frac{\ph_{q,i,k}}{\delta_{q+\bn}^{\sfrac32}r_q^{-1}\Gamma^{3j-3}_{q}}\right) \, ,
\label{eq:a:xi:phi:def}
\end{equation}
where $\tilde \gamma_\xi$ is defined in Proposition~\ref{prop:geo2}, $\zeta_{q,\varphi,i,k,\xi,\vecl}$ is defined in Definition~\ref{def:checkerboard}, and
\begin{align}
    \psi_{i,q}^\varphi &:= \psi_{i,q}^2 \, , \qquad 
    \omega_{j,q}^\varphi := \omega_{j,q}^2  \, , \qquad
    {\chi_{i,k,q}^\varphi := \chi_{i,k,q}^2} \, .
    \label{eq:current:coeff:defs}
\end{align}
From Corollary~\ref{cor:D:Dt:Rn:sharp:new} and estimate \eqref{eq:Lagrangian:Jacobian:1} from Corollary~\ref{cor:deformation}, we have that $|\ph_\ell| \lesssim \Gamma_{q}^{3j-7} \delta_{q+\bn}^{\sfrac32} r_q^{-1}$, and so $\varphi_{q,i,k}$ is well-defined on the support of $\psi_{i,q}^\ph\omega_{j,q}^\ph$ once $\lambda_0$ is sufficiently large.

The coefficient function $a_{(\xi),\varphi}$ is then multiplied by an intermittent pipe bundle $ \nabla \Phi_{(i,k)}^{-1}  \BB_{(\xi),\ph} \circ \Phi_{(i,k)}$, where we have used Proposition~\ref{prop:pipe.flow.current} (with $\lambda=\lambda_{q+\bn}$ and $r=r_{q}$), Definition~\ref{defn:pipe.bundle}, and the shorthand notation
\begin{align} 
  \mathbb{B}_{(\xi),\varphi} = \rhob_{(\xi)}^\varphi \sum_I \zetab_{\xi}^{I,\varphi} \WW_{(\xi),\varphi}^{I}
\label{eq:W:xi:q+1:phi:def}
\end{align}
to refer to the pipe bundle associated with the region {$\Omega_0= \supp\zeta_{q,\varphi,i,k,\xi,\vecl}\cap \{t=k\tau_q\Gamma_q^{-i}\}$ and the index $j$. The choice of placement of this pipe bundle\index{intermittent Mikado bundle}\index{$\BB_{\pxi,\diamond}$}\index{$w_{q+1,\varphi}$} will be detailed in subsection~\ref{ss:stress:oscillation:2}.} We will use $\UU_{(\xi),\ph}^I$ to denote the potential satisfying $\curl \UU_{(\xi),\ph}^I=\WW_{(\xi),\ph}^I$. Applying the algebraic identity \eqref{eq:pipes:flowed:1} from Proposition~\ref{prop:pipeconstruction}, we define the principal part of the current corrector by
\begin{equation}\label{wqplusoneonephip}
    w_{q+1,\ph}^{(p)} = \sum_{i,j,k,\xi,\vecl,I} \underbrace{a_{(\xi),\ph} \left(\chib_{(\xi)}^{\ph} \etab^{I,\ph}_\xi\right) \circ \Phiik  \curl \left( \nabla\Phi_{(i,k)}^T \mathbb{U}_{(\xi),\ph}^I \circ \Phi_{(i,k)} \right)}_{=: w_{(\xi),\varphi}^{(p),I}} \, .
\end{equation}
The notation $w_{(\xi),\ph}^{(p),I}$ refers to \emph{fixed} values of the indices $i,j,k,\xi,\vecl,I$. We add the divergence corrector
\begin{equation}\label{wqplusoneonecphi}
    w_{q+1,\ph}^{(c)} = \sum_{i,j,k,\xi,\vecl,I} \underbrace{\nabla \left( a_{(\xi),\ph} \left(\chib_{(\xi)}^\ph \etab^{I,\ph}_\xi \right) \circ \Phiik \right)  \times 
    \left( \nabla\Phi_{(i,k)}^T \mathbb{U}_{(\xi),\ph}^I \circ \Phi_{(i,k)} \right)}_{ =: w_{(\xi),\ph}^{(c),I}}  \, ,
\end{equation}
so that the mean-zero, divergence-free total current corrector is given by
\begin{equation}\label{wqplusoneonephi}
    w_{q+1,\ph} = w_{q+1,\ph}^{(p)} + w_{q+1,\ph}^{(c)} = \sum_{i,j,k,\xi,\vecl,I} \underbrace{\curl \left( a_{(\xi),\ph} \left(\chib_{(\xi)}^\ph \etab^{I,\ph}_\xi\right) \circ \Phiik \nabla\Phi_{(i,k)}^T \mathbb{U}_{(\xi),\ph}^I \circ \Phi_{(i,k)} \right)}_{ =: w_{(\xi),\varphi}^I} \, .
\end{equation}

\subsubsection{Definition of the Euler-Reynolds corrector}\label{ss:stress:definition}

For any fixed values of $i$, $k$, we recall \eqref{eq:checkie:sixthed} and define
\begin{align}
&R_{q,i,k} = -\nabla\Phi_{(i,k)}
\Biggl{(}{R}_\ell {-\pi_\ell}\Id  \notag\\
& + \sum_{\substack{\xi',i',j' \\ k',l'}}  \frac{\delta_{q+\bn} \Gamma_q^{2j'-2} {C\Gamma_q^{-2}} }{\left| \nabla\Phi^{-1}_{(i',k')} \xi'\right|^{\sfrac 43}} \psi_{i',q}^4 \omega_{j',q}^4 \chi_{i',k',q}^4 \mathcal{X}_{q,\xi',l'}^4 \circ \Phi_{i',k',q} \tilde \gamma_{\xi'}^2 \nabla\Phi_{(i',k')}^{-1}\xi'\otimes \xi'\left(\nabla\Phi_{(i',k')}^{-T}\right)
\Biggr{)}\nabla\Phi_{(i,k)}^T \, , \label{eq:rqnpj}
\end{align}
where the constant $C=c_0c_1c_2$ is geometric and bounded independently of $q$; see \eqref{eq:BB.ph.decomp}. For all $\xi\in\Xi_R$, we define the coefficient function\index{$w_{q+1,R}$} $a_{\xi,i,j,k,\vecl,R}$ by\index{$a_{\xi,i,j,k,\vecl,R}$}
\begin{align}
a_{\xi,i,j,k,\vecl,R}
&=a_{(\xi),R}=\delta_{q+\bn}^{\sfrac 12}\Gamma^{j-1}_{q} \psi_{i,q}^R \omega_{j,q}^R  \chi_{i,k,q}^R \zeta_{q,R,i,k,\xi,\vecl}\, \gamma_{\xi,\Gamma_q^9}\left(\frac{R_{q,i,k}}{\delta_{q+\bn}\Gamma_q^{2j-2}}\right)
\label{eq:a:xi:def}
\end{align}
where $\gamma_{\xi,\Gamma_q^9}$ is defined in Proposition~\ref{p:split} with the parameter choice $K=\Gamma_q^9$, and 
\begin{align}
    \psi_{i,q}^R &:= \psi_{i,q}^3 \, , \qquad \omega_{j,q}^R := \omega_{j,q}^3 \, , \qquad \chi_{i,k,q}^R := \chi_{i,k,q}^3 \, . \label{eq:Reynolds:coeff:defs}
\end{align}
In order to show that \eqref{eq:a:xi:def} is well-defined, we first recall \eqref{pt.est.pi.lem.lower} from Lemma~\ref{lem:D:Dt:Rn:sharp}, which gives that $\pi_\ell |_{\supp \omega_{j,q}} \geq \sfrac 14 \Gamma_q^{2j}\delta_{q+\bn}$. Using this in combination with Corollary~\ref{cor:D:Dt:Rn:sharp:new}, we find that for all $j$,
\begin{equation}\label{pi:top:bottom}
    \Gamma_q \leq \frac{\pi_\ell|_{\supp \omega_{j,q}}}{\delta_{q+\bn}\Gamma_q^{2j-2}} \leq \Gamma_q^9 \, .
\end{equation}
Furthermore, from \eqref{eq:rqnpj}, \eqref{eq:omega:cut:partition:unity}, and Corollary~\ref{cor:deformation}, we have that the second term in \eqref{eq:rqnpj} is pointwise bounded by $2C\delta_{q+\bn}\Gamma_q^{2j-2}$, or upon division by $\delta_{q+\bn}\Gamma_q^{2j-2}$ is bounded above by $2C$.  Finally, from \eqref{pt.est.R.lem}, we have that $\nabla\Phiik R_\ell \nabla\Phiik^T$ is pointwise bounded by $\delta_{q+\bn}\Gamma_q^{2j-3}$, or upon division by $\delta_{q+\bn}\Gamma_q^{2j-2}$ is pointwise bounded by $\Gamma_q^{-1}$.  Combining the above arguments, we find that 
\begin{align}
    \left| \frac{R_{q,i,k}}{\delta_{q+\bn}\Gamma_q^{2j-2}} - \frac{\pi_\ell}{\delta_{q+\bn}\Gamma_q^{2j-2}} \Id \right| &\leq \Gamma_q \, ,  \notag 
\end{align}
and so Proposition~\ref{p:split} may be applied with $K=\Gamma_q^{9}$ since $\frac{R_{q,i,k}}{\delta_{q+\bn}\Gamma_q^{2j-2}}$ belongs to the ball of radius $\Gamma_q$ around $\frac{\pi_\ell \Id}{\delta_{q+\bn}\Gamma_q^{2j-2}}$, which itself is a multiply of the identity bounded between $1$ and $\Gamma_q^9$ from \eqref{pi:top:bottom}.

The coefficient function $a_{(\xi),R}$ is then multiplied by an intermittent pipe bundle $\nabla \Phi_{(i,k)}^{-1}  \BB_{(\xi),R} \circ \Phi_{(i,k)}$, where we have used Proposition~\ref{prop:pipeconstruction} (with $\lambda=\lambda_{q+\bn}$ and $r=r_{q}$), Definition~\ref{defn:pipe.bundle}, and the shorthand notation
\begin{align} 
\mathbb{B}_{(\xi),R} = \rhob_{(\xi)}^R \sum_I \zetab_{\xi}^{I,R} \WW_{(\xi),R}^{I}
\label{eq:W:xi:q+1:nn:def}
\end{align}
to refer to the pipe bundle associated with the region $\Omega_0= \supp\zeta_{q,R,i,k,\xi,\vecl}\cap \{t=k\tau_q\Gamma_q^{-i}\}$ and the index $j$. We will use $\UU_{(\xi),R}^I$ to denote the potential satisfying $\curl \UU_{(\xi),R}^I=\WW_{(\xi),R}^I$. Applying \eqref{eq:pipes:flowed:1} from Proposition~\ref{prop:pipeconstruction}, we define the principal part of the Reynolds corrector by
\begin{equation}\label{wqplusoneonep}
    w_{q+1,R}^{(p)} = \sum_{i,j,k,\xi,\vecl,I} \underbrace{a_{(\xi),R} \left(\chib_{(\xi)}^R \etab_{\xi}^{I,R} \right) \circ \Phiik  \curl \left(\nabla\Phi_{(i,k)}^T \mathbb{U}_{(\xi),R}^I \circ \Phi_{(i,k)} \right)}_{=: w_{(\xi),R}^{(p),I}}
    \, .
\end{equation}
The notation $w_{(\xi),R}^{(p),I}$ refers to \emph{fixed} values of $i,j,k,\xi,\vecl,I$.  We add the divergence corrector
\begin{equation}\label{wqplusoneonec}
    w_{q+1,R}^{(c)} =\sum_{i,j,k,\xi,\vecl,I} \underbrace{\nabla  \left( a_{(\xi),R} \left(\chib_{(\xi)}^R \etab_{\xi}^{I,R} \right) \circ \Phiik \right) \times \left(\nabla\Phi_{(i,k)}^T \mathbb{U}_{(\xi),R}^I \circ \Phi_{(i,k)} \right)}_{=: w_{(\xi),R}^{(c),I}}
    \, ,
\end{equation}
so that the mean-zero, divergence-free total Reynolds corrector is given by
\begin{equation}\label{wqplusoneone}
    w_{q+1,R} = \sum_{i,j,k,\xi,\vecl,I} \underbrace{\curl \left( a_{(\xi),R} \left(\chib_{(\xi)}^R \etab_{\xi}^{I,R} \right) \circ \Phiik  \nabla\Phi_{(i,k)}^T \mathbb{U}_{(\xi),R}^I \circ \Phi_{(i,k)} \right)}_{=: w_{(\xi),R}^I} \, .
\end{equation}

\subsubsection{Definition of the complete corrector}
We shall sometimes want to aggregate pieces of the Reynolds and current velocity correctors as
\begin{align}\label{defn:w}
    w_{q+1} = w_{q+1,R} + w_{q+1,\varphi} \, , \qquad  w_{q+1}^{(p)} := w_{q+1,R}^{(p)} + w_{q+1,\varphi}^{(p)} \, , \qquad w_{q+1}^{(c)} := w_{q+1,R}^{(c)} + w_{q+1,\varphi}^{(c)} \, .
\end{align}

\subsection{Dodging for new velocity increment}
 \label{ss:stress:oscillation:2}
In this section, we define a mollified velocity increment $\hat w_{q+\bn}$. We then introduce Lemma \ref{lem:dodging}, which is in fact a stronger statement than Hypothesis \ref{hyp:dodging1}.

\begin{definition}[\bf Definition of $\hat w_\qbn$ and $u_{q+1}$]\label{def:wqbn}
Let $\mathcal{\tilde P}_{q+\bn,x,t}$\index{$\mathcal{\tilde P}_{q+\bn,x,t}$} denote a space-time mollifier which is a product of compactly supported kernels at spatial scale $\lambda_{q+\bn}^{-1}\Gamma_{q+\bn-1}^{-\sfrac 12}$ and temporal scale $\Tau_{q+1}^{-1}$.  We again assume that both kernels have vanishing moments up to $10\Nfin$ and are $C^{10\Nfin}$ differentiable and define\index{$\hat w_\qbn$}\index{$u_{q+1}$}
\begin{equation}\label{def.w.mollified}
    \hat w_{q+\bn} := \mathcal{\tilde P}_{q+\bn,x,t} w_{q+1} \, , \qquad u_{q+1} = u_q + \hat w_{\qbn} \, .
\end{equation}
\end{definition}
\noindent 
We also recall from \eqref{eq:space:time:balls} the notations $B(\Omega,\lambda^{-1})$ and $B(\Omega,\lambda^{-1},\tau)$ for space and space-time balls, respectively, around a space-time set $\Omega$.  Using these notations, we may write that
\begin{equation}
    \supp \hat w_\qbn \subseteq B\left( \supp w_{q+1}, \sfrac 12 \lambda_\qbn^{-1} , \sfrac 12 \Tau_q \right) \, . \label{eq:dodging:useful:support}
\end{equation}
Now recalling the formula in \eqref{eq:WW:explicit} for an intermittent Mikado flow, \eqref{eq:W:xi:q+1:phi:def}, and \eqref{eq:W:xi:q+1:nn:def}, we set
\begin{align}\label{eq:pipez:thursday}
    \varrho_{(\xi),\diamond}^{I} := \xi \cdot \WW_{(\xi),\diamond}^I \, .
\end{align}
Next, in slight conflict with \eqref{eq:space:time:balls}, we shall also use the notation
\begin{align}\label{eq:ballz:useful}
    B\left(\supp \varrho_{\pxi,\diamond}^I,\lambda^{-1}\right) := \left\{ x\in\T^3 \, : \, \exists y \in \supp \varrho_{\pxi,\diamond}^I \, , |x-y| \leq \lambda^{-1} \right\}
\end{align}
throughout this section, despite the fact that $\supp\varrho_{\pxi,\diamond}^I$ is not a set in space-time, but merely a set in space. We shall also use the same notation but with $\varrho_{\pxi,\diamond}^I$ replaced by $\rhob_\xi^\diamond$. Finally, for any smooth set $\Omega\subseteq \mathbb{T}^3$ and any flow map $\Phi$ defined in Definition~\ref{def:transport:maps}, we use the notation 
\begin{equation}\label{eq:flowing:sets}
\Omega \circ \Phi := \left\{(y,t): t\in \R, \Phi(y,t)\in \Omega\right\} = \supp \left(\mathbf{1}_{\Omega}\circ \Phi\right) \, .
\end{equation}
In other words, for any smooth set $\Omega\subseteq\T^3$, $\Omega\circ\Phi$ is a space-time set whose characteristic function is annihilated by $\Dtq$.

We can now introduce the workhouse which will help us verify Hypotheses~\ref{hyp:dodging1} and \ref{hyp:dodging2}.  The full proof is contained in~\cite[section~4]{GKN23}, although we outline the main idea following the statement.

\begin{lemma}[\bf Dodging and preventing self-intersections for $w_{q+1}$ and $\hat w_\qbn$]\label{lem:dodging} 
We construct $w_{q+1}$ so that the following hold.\index{dodging}\index{effective dodging}
\begin{enumerate}[(i)]
    \item\label{item:dodging:more:oldies} Let {$q+1\leq q' \leq q+ \sfrac \bn 2$} and fix indices $\diamond,i,j,k,\xi,\vecl$, which we abbreviate by $(\pxi,\diamond)$, for a coefficient function $a_{\pxi,\diamond}$ (cf.~\eqref{eq:a:xi:phi:def}, \eqref{eq:a:xi:def}).  Then
    \begin{equation}
        B\left( \supp \hat w_{q'}, \frac 12 {\lambda_{q+1}^{-1}\Ga_q^2}, {2 \Tau_q} \right) \cap \supp \left( \tilde \chi_{i,k,q} \zeta_{q,\diamond,i,k,\xi,\vecl} \, \rhob_{\pxi}^{\diamond}\circ \Phi_{(i,k)} \right) = \emptyset \, . \label{eq:oooooldies}
    \end{equation}
    \item\label{item:dodging:1} Let $q'$ satisfy $q+1\leq q' \leq q+\bn-1$, fix indices $(\pxi,\diamond,I)$, and assume that $\Phiik$ is the identity at time $t_{\pxi}$, cf. Definition~\ref{def:transport:maps}. Then we have that
\begin{align}
    B \left(\supp \hat w_{q'}, \frac 14 \lambda_{q'}^{-1} \Gamma_{q'}^2, 2\Tau_{q} \right)  \cap  \supp &\left( \tilde \chi_{i,k,q} \zeta_{q,\diamond,i,k,\xi,\vecl} \left(\rhob_{(\xi)}^\diamond \zetab_{\xi}^{I,\diamond} \right)\circ \Phiik \right) \notag \\
    &\cap
    B\left( \supp \varrho^I_{(\xi),\diamond} , \frac 12 {\lambda_{q'}^{-1} \Gamma_{q'}^2}\right)\circ \Phiik
    = \emptyset \, . \label{eq:dodging:oldies:prep}
    \end{align}
As a consequence we have
\begin{equation}\label{eq:dodging:oldies}
        B\left( \supp \hat w_{q'}, \frac 14 {\lambda_{q'}^{-1} \Gamma_{q'}^2}, 2\Tau_q \right) \cap  \supp w_{q+1} = \emptyset \, .
    \end{equation}
    \item\label{item:dodging:2} Consider the set of indices $\{(\pxi,\diamond,I)\}$, whose elements we use to index the correctors constructed in \eqref{wqplusoneonephi} and \eqref{wqplusoneone}, and let $\ttl, \ov \ttl \in \{p,c\}$ denote either principal or divergence corrector parts. Then if $(\ov\diamond,(\ov \xi), \ov I) \neq (\diamond,(\xi),I)$, we have that for any $\ttl, \ov \ttl$,
    \begin{equation}\label{eq:dodging:newbies}
        \supp w_{\pxi,\diamond}^{(\ttl),I} \cap \supp w_{(\ov \xi),\ov \diamond}^{(\ov \ttl), \ov I} = \emptyset \, .
    \end{equation}
    \item\label{item:dodging:zero}  $\hat w_\qbn$ satisfies Hypothesis~\ref{hyp:dodging2} with $q$ replaced by $q+1$.
\end{enumerate}
\end{lemma}
\begin{remark}[\bf Verifying Hypothesis~\ref{hyp:dodging1}]\label{rem:checking:hyp:dodging:1}
We claim that \eqref{eq:dodging:oldies} and \eqref{eq:dodging:useful:support} imply that Hypothesis~\ref{hyp:dodging1} holds with $q+1$ replacing all instances of $q$.  To check this, we must show that \eqref{eq:ind:dodging} holds for $q',q''\leq \qbn$ and $0<|q'-q''|\leq \bn-1$.  By induction on $q$ and the symmetry of $q''$ and $q'$, the only case we must check is the case that $q+\bn=q''$ and $0<\qbn-q'\leq\bn-1$. But it is a simple exercise in set theory to check that for $q+1\leq q'\leq \qbn-1$, \eqref{eq:dodging:oldies} is equivalent to $\supp \hat w_{q'} \cap B(\supp w_{q+1},\sfrac 14 \lambda_{q'}^{-1}\Ga_{q'}^2,2\Tau_q) = \emptyset$. Then using \eqref{eq:dodging:useful:support} and the inequalities $\lambda_{q'}^{-1}\Ga_{q'}^2 \geq \la_\qbn^{-1}$, $b<2 \implies \Ga_{q'+1}\ll \Ga_{q'}^2$ implies that \eqref{eq:ind:dodging} holds.
\end{remark}
\begin{proof}[Idea behind the proof of Lemma~\ref{lem:dodging}]
We shall give the idea behind the proof of Hypothesis~\ref{hyp:dodging1}, as the precise statements written above are technical variants on this idea and can be found in \cite[section~4]{GKN23}.  Consider the support of a single mildly anisotropic cutoff $\zeta_{q,\diamond,i,k,\xi,\vecl}$\index{$\zeta_{q,\diamond,i,k,\xi,\vecl}$} from Definition~\ref{def:checkerboard} of dimensions $(\lambda_{q+1}\Gamma_q^{-5})^{-1} \times (\lambda_{q+1}\Gamma_q^{-5})^{-1}\times (\lambda_q\Gamma_q^8)^{-1}$. The prism contains pipes from $\hat w_{q+1},\dots \hat w_{q+\half}$, and we want to place a new set of bundling pipes $\rhob_{\pxi}^\diamond$ from Proposition~\ref{prop:bundling} of thickness $\lambda_{q+1}^{-1}\Gamma_q$ and spacing $\lambda_{q+1}^{-1}\Gamma_q^4$ disjoint from these pipes.  To this end, we divide the face $[0,\lambda_{q+1}^{-1}\Gamma_q^5]^2$ of the prism perpendicular to $\vec e_3$ into the grid of squares of sidelength $\lambda_{q+1}^{-1}\Gamma_q$ (the thickness of the support of $\rhob_{\pxi}^\diamond$). Since the support of $\rhob_{\pxi}^\diamond$ will be placed $\mathbb{T}^2/(\la_{q+1}^{-1}\Gamma_q^4)$-periodically, 
\begin{align*}
        \text{the possible number of placements of the support}
        = \left(\frac{\la_{q+1}^{-1} \Gamma_q^4}{\la_{q+1}^{-1}\Gamma_q} \right)^2 = \Gamma_q^6 \, .
\end{align*}
The pipes that we want to dodge have spacing/thickness between $\lambda_{q-\half}^{-1}$/$\lambda_q^{-1}$ (correponding to $\hat w_{q}$) and $\lambda_{q}^{-1}$/$\lambda_{q+\half}^{-1}$ (corresponding to $\hat w_{q+\half}$); note that each of these has spacing greater than $(\lambda_q\Ga_q^8)^{-1}$, which is the longest side length of the prism. Then from Hypothesis~\ref{hyp:dodging2}, at most a constant number of such pipes can intersect the prism. Upon projecting these pipes onto the face $[0,\lambda_{q+1}^{-1}\Gamma_q^5]^2$ perpendicular to $\vec e_3$, each pipe projection will be contained in a $\lambda_{q+1}^{-1}$-neighborhood of a line of length $\lambda_{q+1}^{-1}\Gamma_q^5$. Counting the number of grid squares of size $\lambda_{q+1}^{-1}\Gamma_q$ taken by these projections, we obtain
\begin{align*}
        \sim \frac{\lambda_{q+1}^{-1}\Gamma_q^5}{\lambda_{q+1}^{-1}} \lesssim \Gamma_q^5 \, ,
\end{align*}
which is less than the possible number of placements. Therefore we can place the support of the bundling pipe $\rhob_\pxi^\diamond$ so that it is disjoint from $\hat w_{q+1},\dots \hat w_{q+\half}$ on the support of $\zeta_{q,\diamond,i,k,\xi,\vecl}$.

To enact the dodging with pipes from $\hat w_{q+\half+1},\dots ,\hat w_\qbn$ of thickness/spacing $\la_{q+\half+1}^{-1}$/$\la_{q+1}^{-1}$, \dots, \\ $\la_{\qbn}^{-1}$/$\la_{q+\half}^{-1}$, we follow the exact same method, only replacing the mildly anistropic cutoff $\zeta_{q,\diamond,i,k,\xi,\vecl}$ with the highly anistropic cutoff $\zetab_\xi^\diamond$ from Lemma~\ref{lem:finer:checkerboard:estimates}, and the mildly intermittent bundling pipe $\rhob_\pxi^\diamond$ with the highly intermittent pipes $\WW_{\pxi,\diamond}^I$ from Propositions~\ref{prop:pipeconstruction} and \ref{prop:pipe.flow.current}. We leave further details to the reader.
\end{proof}

\subsection{Estimates for \texorpdfstring{$w_{q+1}$}{wqn} and \texorpdfstring{$\hat w_{q+\bn}$}{hwqn}}\label{ss:stress:w:estimates}

\begin{lemma}[\bf Coefficient function estimates]
\label{lem:a_master_est_p}
For $N,N',N'',M$ with $N'',N'\in\{0,1\}$ and $N,M \leq {\sfrac{\Nfin} {3}}$, we have the following estimates.
\begin{subequations}\label{e:a_master_est_p}
\begin{align}
&\left\|D^{N-N''} D_{t,q}^M (\xi^\ell A_\ell^h \partial_h)^{N'} D^{N''} a_{\xi,i,j,k,\vec{l},\varphi}\right\|_{ r} \notag\\
&\qquad\lessg \left| \supp \eta_{i,j,k,\xi,\vecl,\varphi} \right|^{\sfrac 1r} \delta_{q+\bn}^{\sfrac 12} r_q^{-\sfrac 13} \Gamma^{j- {1}}_{q}  \left({\Gamma_{q}^{-5}\lambda_{q+1}}\right)^N  \left({\Gamma_q^{  5}\Lambda_{q}}\right)^{N'} \MM{M, \NindSmall, \tau_{q}^{-1}\Gamma_{q}^{i+4}, \Tau_{q}^{-1}}\label{e:a_master_est_p_phi} \, , \\
&\left\| D^{N-N''} D_{t,q}^M (\xi^\ell A_\ell^h \partial_h)^{N'} D^{N''} \left( a_{\xi,i,j,k,\vecl,\varphi} \left( \rhob_{(\xi)}^\varphi \zetab_\xi^{I,\varphi} \right) \circ \Phiik \right) \right\|_{ r}\notag\\
&\qquad\lessg \left| \supp \left( \eta_{i,j,k,\xi,\vecl,\varphi} \zetab_\xi^{I,\varphi} \right) \right|^{\sfrac 1r} \delta_{q+\bn}^{\sfrac 12} r_q^{-\sfrac 13} \Gamma_q^{j {+1}} \left({\lambda_{q+\lfloor \sfrac \bn 2 \rfloor}}\right)^N \left({\Gamma_q^{  5}\Lambda_{q}}\right)^{N'} \MM{M, \NindSmall, \tau_{q}^{-1}\Gamma_{q}^{i+4}, \Tau_{q}^{-1}}\label{e:a_master_est_p_phi:zeta}
\, ,\\
&\left\| D^{N-N''} D_{t,q}^M (\xi^\ell A_\ell^h \partial_h)^{N'} D^{N''} a_{\xi,i,j,k,\vecl,R}\right\|_{ r} \notag\\
&\qquad\lesssim \left| \supp \eta_{i,j,k,\xi,\vecl,R} \right|^{\sfrac 1r} \delta_{q+\bn}^{\sfrac 12} \Gamma_q^{j+ {4}} \left({\Gamma_{q}^{-5}\lambda_{q+1}}\right)^N \left({\Gamma_q^{ {13}}\Lambda_{q}}\right)^{N'} \MM{M, \NindSmall, \tau_{q}^{-1}\Gamma_{q}^{i+ {13}}, \Tau_{q}^{-1} {\Ga_q^8}}\label{e:a_master_est_p_R} \, , \\
&\left\| D^{N-N''} D_{t,q}^M (\xi^\ell A_\ell^h \partial_h)^{N'} D^{N''} \left( a_{\xi,i,j,k,\vecl,R} \left( \rhob_{(\xi)}^R \zetab_\xi^{I,R} \right) \circ \Phiik \right) \right\|_{ r}\notag\\
&\qquad\lessg \left| \supp \left( \eta_{i,j,k,\xi,\vecl,R} \zetab_\xi^{I,R} \right) \right|^{\sfrac 1r} \delta_{q+\bn}^{\sfrac 12} \Gamma_q^{j+ {7}} \left({\lambda_{q+\lfloor \sfrac \bn 2 \rfloor}}\right)^N \left({\Gamma_q^{ {13}}\Lambda_{q}}\right)^{N'} \MM{M, \NindSmall, \tau_{q}^{-1}\Gamma_{q}^{i+ {13}}, \Tau_{q}^{-1} {\Ga_q^8}}\label{e:a_master_est_p_R:zeta} \, .
\end{align}
\end{subequations}
In the case that $r=\infty$, the above estimates give that
\begin{subequations}
\begin{align}
\left\| D^{N-N''} D_{t,q}^M (\xi^\ell A_\ell^h \partial_h)^{N'} D^{N''} a_{\xi,i,j,k,\vec{l},R}\right\|_{ \infty} 
&\lessg {\Gamma_q^{\frac{\badshaq}{2}+ {7}}} \left(\Gamma_{q}^{-5}\lambda_{q+1}\right)^N \notag\\
&\qquad \qquad \times \left({\Gamma_q^{ {13}}\Lambda_{q}} \right)^{N'} \MM{M,\Nindt, \tau_{q}^{-1}\Gamma_{q}^{i+ {13}},\Tau_q^{-1} {\Ga_q^8}} \label{e:a_master_est_p_uniform_R} \, . \\
\left\| D^{N-N''} D_{t,q}^M (\xi^\ell A_\ell^h \partial_h)^{N'} D^{N''} a_{\xi,i,j,k,\vec{l},\varphi}\right\|_{ \infty} 
&\lessg {\Gamma_q^{\frac{\badshaq}{2}+ {2}}} r_q^{-\sfrac 13} \left(\Gamma_{q}^{-5}\lambda_{q+1}\right)^N \notag\\
&\qquad \qquad \times \left({\Gamma_q^8\Lambda_{q}} \right)^{N'} \MM{M,\Nindt, \tau_{q}^{-1}\Gamma_{q}^{i+4},\Tau_q^{-1}} \label{e:a_master_est_p_uniform_phi} \, ,
\end{align}
\end{subequations}
with analogous estimates (incorporating a loss of $\Ga_q^3$ for $\diamond=R$ and $\Ga_q^2$ for $\diamond=\varphi$)\index{$\diamond$} holding for the product $a_{(\xi),\diamond}\zetab_{\xi}^{I,\diamond}\rhob_{(\xi)}^\diamond$. Finally, we have the pointwise estimates
\begin{subequations}\label{e:a_master_est_p_pointwise}
\begin{align}
    \left| D^{N-N''} D_{t,q}^M (\xi^\ell A_\ell^h \partial_h)^{N'} D^{N''} a_{\xi,i,j,k,\vecl,R}\right| &\lesssim \Gamma_q^{12} \pi_\ell^{\sfrac 12} \left({\Gamma_{q}^{-5}\lambda_{q+1}}\right)^N \left({\Gamma_q^{ {13}}\Lambda_{q}}\right)^{N'} \MM{M, \NindSmall, \tau_{q}^{-1}\Gamma_{q}^{i+ {13}}, \Tau_{q}^{-1} {\Ga_q^8}} \label{e:a_master_est_p_R_pointwise} \\     \left| D^{N-N''} D_{t,q}^M (\xi^\ell A_\ell^h \partial_h)^{N'} D^{N''} a_{\xi,i,j,k,\vecl,\varphi}\right| &\lesssim \Gamma_q^{ {12}} \pi_\ell^{\sfrac 12} r_q^{-\sfrac{1}{3}} \left({\Gamma_{q}^{-5}\lambda_{q+1}}\right)^N \left({\Gamma_q^{ {5}}\Lambda_{q}}\right)^{N'} \MM{M, \NindSmall, \tau_{q}^{-1}\Gamma_{q}^{i+4}, \Tau_{q}^{-1}} \, . \label{e:a_master_est_p_phi_pointwise}
\end{align}
\end{subequations}
\end{lemma}
\begin{proof}
[Proof of Lemma~\ref{lem:a_master_est_p}]
We first prove \eqref{e:a_master_est_p_phi} and \eqref{e:a_master_est_p_phi:zeta}, since a portion of $a_{\pxi,\varphi}$ appears in the definition of the Reynolds corrector in \eqref{eq:rqnpj}. We further simplify by computing \eqref{e:a_master_est_p_phi} for the case $r=\infty$ first. Recalling estimate \eqref{eq:D:Dt:ph:sharp:new}, we have that for all $N,M\leq \sfrac{\Nfin}{2}$,
\begin{align}
\bigl\| D^N D_{t,q}^M \varphi_\ell \bigr\|_{L^{\infty}(\supp \psi_{i,q}\omega_{j,q})}
&\lessg \delta_{q+\bn}^{\sfrac 32} r_q^{-1} \Gamma^{3j-7}_{q}
 \left(\Gamma_{q}\Lambda_q \right)^N
\MM{M, \NindSmall, \tau_{q}^{-1}\Gamma_{q}^{i}, \Tau_{q}^{-1}} \, . \notag  
\end{align}
Thus from definition \eqref{eq:phiqnpj}, the Leibniz rule, and Corollary~\ref{cor:deformation}, and the fact that $\supp \eta_{i,j,k,\xi,\vecl,\varphi}$ is contained in $\supp \psi_{i,q}\omega_{j,q}\chi_{i,k,q}$ we have that for $N,M\leq \sfrac{\Nfin}{2}$,
\begin{align}
&\norm{D^N D_{t,q}^M \vp_{q,i,k} }_{L^{\infty}(\supp \eta_{i,j,k,\xi,\vecl,\vp})} \lessg \delta_{q+\bn}^{\sfrac 32} r_q^{-1} \Gamma^{3j-7}_{q} \left(\Gamma_{q}\Lambda_q \right)^N \MM{M, \NindSmall, \tau_{q}^{-1}\Gamma_{q}^{i}, \Tau_{q}^{-1}} \label{eq:davidc:1:(}
\,.\end{align}
The above estimates allow us to apply \cite[Lemma~A.5]{BMNV21} with $N'=M'=\sfrac{\Nfin}{2}$, $\psi = \td\gamma_{\xi,}$, $\Gamma_\psi=1$, $v = \hat u_q$, $D_t = D_{t,q}$, $h(x,t) = \vp_{q,i,k}(x,t)$, $C_h = \delta_{q+\bn}^{\sfrac 32}r_q^{-1}\Gamma_{q}^{3j-6} = \Gamma^2$, $\lambda=\tilde\lambda = \Lambda_q\Gamma_{q}$, $\mu = \tau_{q}^{-1} \Gamma_{q}^{i}$, $\tilde \mu = \Tau_{q}^{-1}$, and $N_t=\Nindt$. We obtain that for all $N,M\leq \sfrac{3\Nfin}{4}$,
\begin{align}
\norm{ D^ND_{t,q}^M \td\gamma_{\xi}\left(\frac{\vp_{q,i,k}}{\delta_{q+\bn}^{\sfrac 32}r_q^{-1}\Gamma_q^{3j-3}}\right)}_{L^{\infty}(\supp \eta_{i,j,k,\xi,\vecl,\vp})} &\lesssim \left(\Gamma_{q}\Lambda_q \right)^N \MM{M, \NindSmall, \tau_{q}^{-1}\Gamma_{q}^{i}, \Tau_{q}^{-1}} \,. \label{eq:gamma:phi:ell}
\end{align}
Finally, from Corollary~\ref{cor:deformation} and an application of the mixed derivative Fa'a di Bruno formula from \cite[Lemma~A.5]{BMNV21} with $\psi(\cdot):B_{\sfrac 12}(\xi)\rightarrow \R$ defined by $\psi(\cdot)=|\cdot|^{-\sfrac 43}$, $\Gamma_\psi=1$, $v=\hat u_q$, $\Gamma=1$, $\lambda=\tilde\lambda=\Lambda_q$, $\mu = \tau_q^{-1}\Gamma_q^i$, $\tilde \mu = \Gamma_q^{-1}\Tau_q^{-1}$, $N_x=0$, $N_t=\Nindt$, $h=\nabla\Phiik^{-1}\xi$, and $\const_h=1$, we have that for all $N+M\leq \sfrac{3\Nfin}{2}$,
\begin{align*}
\left\| D^N D_{t,q}^M \left( \left| \nabla \Phiik^{-1} \xi \right|^{-\sfrac 43} \right) \right\|_{L^\infty(\supp(\psi_{i,q}^{ {\vp}}\chi^{ {\vp}}_{i,k,q}))} &\les {\Lambda}_q^{N} \MM{M,\NindSmall,\Gamma_{q}^{i} \tau_q^{-1},\Tau_q^{-1}\Gamma_{q}^{-1}}.
\end{align*}
From the above three bounds, definition \eqref{eq:a:xi:phi:def}, the Leibniz rule, estimate \eqref{eq:nasty:Dt:psi:i:q:orangutan} at level $q$, \eqref{eq:chi:cut:dt}, \eqref{eq:D:Dt:omega:sharp}, and \eqref{eq:checkerboard:derivatives}, we obtain that for $N'=0,1$ and $N,M\leq \sfrac{\Nfin}{2}$,
\begin{align}
\bigl\| D^N D_{t,q}^M ( \xi^\ell A_\ell^j \partial_j)^{N'} a_{\xi,i,j,k,\vec{l},\vp}\bigr\|_{\infty}
& \lessg \delta_{q+\bn}^{\sfrac 12} \Gamma_{q}^{j- {1}} r_q^{-\sfrac 13} (\Gamma_q^{-5}\lambda_{q+1})^N ({\Gamma_{q}^5\Lambda_q})^{N'} \MM{M,\Nindt, \tau_{q}^{-1}\Gamma_{q}^{i+4},\Tau_q^{-1}} \, . \label{eq:helping:me:pointwise}
\end{align}
Using \eqref{ineq:jmax:use}, we obtain \eqref{e:a_master_est_p_uniform_phi}. When $r\neq \infty$, we use $\left\| f \right\|_{L^r}\leq \left\| f \right\|_{L^\infty} | \{ \supp f \} |^{\sfrac 1r} $ and the demonstrated bound for $r=\infty$ to obtain \eqref{e:a_master_est_p_phi} for the full range of $r$ and for $N''=0$.  The estimate in \eqref{e:a_master_est_p_phi:zeta} for $N''=0$ follows in the same way using \eqref{e:fat:pipe:estimates:1} for $p=\infty$ and \eqref{eq:checkerboard:derivatives:check}. Similar estimates for $N''=1$ in both cases are nearly identical, and we omit the details

We now compute \eqref{e:a_master_est_p_R} for the case $r=\infty$, from which the remaining bounds in \eqref{e:a_master_est_p_R:zeta} and \eqref{e:a_master_est_p_uniform_R} will follow as before. Recalling estimates \eqref{eq:D:Dt:pi:sharp:new} and \eqref{eq:D:Dt:R:sharp:new}, we have that for all $N,M\leq \sfrac{\Nfin}{2}$,
\begin{align}
&\bigl\| D^N D_{t,q}^M R_\ell \bigr\|_{L^{\infty}(\supp \eta_{i,j,k,\xi,\vecl,R})} + \bigl\| D^N D_{t,q}^M \pi_\ell \bigr\|_{L^{\infty}(\supp \eta_{i,j,k,\xi,\vecl,R})} \notag\\
&\qquad \lessg \delta_{q+\bn}\Gamma^{2j+6}_{q}
 \left(\Gamma_{q}\Lambda_q \right)^N
\MM{M, \NindSmall, \tau_{q}^{-1}\Gamma_{q}^{i}, \Tau_{q}^{-1}} \, . \notag  
\end{align}
From \eqref{eq:nasty:Dt:psi:i:q:orangutan} and \eqref{eq:inductive:partition} at level $q$, \eqref{eq:omega:cut:partition:unity}, \eqref{eq:D:Dt:omega:sharp}, \eqref{eq:chi:cut:dt}, \eqref{eq:checkie:sixthed}, \eqref{cor:deformation}, and \eqref{eq:gamma:phi:ell}, we find that
\begin{align}
&\left\| D^N D_{t,q}^M \sum_{i',j',k',\xi',\vecl'} \frac{\delta_{q+\bn}\Gamma_q^{2j'-4}C}{\left| \nabla \Phi_{i',k'} \xi' \right|^{\sfrac 43}} \psi_{i',q}^4 \omega_{j',q}^4 \chi_{i',k',q}^4 \mathcal{X}_{q,\xi',l'}^4 \circ \Phi_{i',k',q} \tilde \gamma_\xi^2  \nabla\Phi_{(i',k')}^{-1} \xi' \otimes \xi' \nabla \Phi_{(i',k')}^{-T} \right\|_{L^\infty(\supp\eta_{i,j,k,\xi,\vecl,R})} \notag\\
&\qquad \lessg \delta_{q+\bn}\Gamma^{2j-4}_{q}
 \left(\Gamma_{q}^{ {5}}\Lambda_q \right)^N
\MM{M, \NindSmall, \tau_{q}^{-1}\Gamma_{q}^{i+ {5}}, \Tau_{q}^{-1}} \, . \notag  
\end{align}
Thus from the Leibniz rule and definition \eqref{eq:rqnpj}, we find that for $N,M\leq \sfrac \Nfin 2$,
\begin{align}
&\norm{D^N D_{t,q}^M R_{q,i,k} }_{L^{\infty}(\supp \eta_{i,j,k,\xi,\vecl,R})} \lessg \delta_{q+\bn} \Gamma^{2j+6}_{q} \left(\Gamma_{q}^{ {5}}\Lambda_q \right)^N \MM{M, \NindSmall, \tau_{q}^{-1}\Gamma_{q}^{i+ {5}}, \Tau_{q}^{-1}} \label{eq:davidc:1}
\,;\end{align}
the loss of $\Ga_q$ in the sharp material derivative cost comes from the fact that the sum includes $\psi_{i',q}$ and is estimated on the supported of $\psi_{i,q}$. The above estimates allow us to apply \cite[Lemma~A.5]{BMNV21} with $N'=M'=\sfrac \Nfin 2$, $\psi = \Ga_q^{-5}\gamma_{\xi,\Ga_q^9}$ as in \eqref{def:gamma:xi:R},\footnote{Since $\gamma_{\xi,\Ga_q^9}$ and all its derivatives are bounded by $\Ga_q^5$ from \eqref{eq:gamma:xi:derivative:bounds}, we first rescale by $\Ga_q^{-5}$ on the outside and then apply the Faa di Bruno lemma, which requires $\psi$ to be bounded in between $0$ and $1$.  Rescaling back then produces the desired bound.} $\Gamma_\psi=1$, $v = \hat u_q$, $D_t = D_{t,q}$, $h(x,t) = R_{q,i,k}(x,t)$, $C_h = \delta_{q+\bn}\Gamma_{q}^{2j {+6}}$, $\Gamma^2=\delta_\qbn \Ga_q^{2j-2}$, $\lambda=\tilde\lambda = \Lambda_q\Gamma_{q}^5$, $\mu = \tau_{q}^{-1} \Gamma_{q}^{i+5}$, $\tilde \mu = \Tau_{q}^{-1}$, and $N_t=\Nindt$. We obtain that for all $N,M\leq \sfrac{\Nfin}{2}$,
\begin{align*}
\norm{ D^ND_{t,q}^M \gamma_{\xi,\Gamma_q^9}\left(\frac{R_{q,i,k}}{\delta_{q+\bn}\Gamma_q^{2j-2}}\right)}_{L^{\infty}(\supp \eta_{i,j,k,\xi,\vecl,R})} &\lesssim \Gamma_q^{ {5}} \left(\Gamma_{q}^{ {13}}\Lambda_q \right)^N \MM{M, \NindSmall, \tau_{q}^{-1}\Gamma_{q}^{i+ {13}}, \Tau_{q}^{-1} {\Ga_q^8}} \,.
\end{align*}
From the above bound, definition \eqref{eq:a:xi:def}, the Leibniz rule, estimate \eqref{eq:nasty:Dt:psi:i:q:orangutan} at level $q$, \eqref{eq:Lagrangian:Jacobian:6}, \eqref{eq:chi:cut:dt}, \eqref{eq:D:Dt:omega:sharp}, and \eqref{eq:checkerboard:derivatives}, we obtain that for $N'=0,1$ and $N,M\leq \sfrac{\Nfin}{2}$,
\begin{align*}
\bigl\| D^ND_{t,q}^M ( \xi^\ell A_\ell^j \partial_j)^{N'} a_{\xi,i,j,k,\vec{l},R}\bigr\|_{L^\infty}
& \lessg \delta_{q+\bn}^{\sfrac 12} \Gamma_{q}^{j+ {4}} (\Gamma_q^{-5}\lambda_{q+1})^N ({\Gamma_{q}^{ {13}}\Lambda_q})^{N'} \MM{M,\Nindt, \tau_{q}^{-1}\Gamma_{q}^{i+ {13}},\Tau_q^{-1} {\Ga_q^8}} \, .
\end{align*}
Using \eqref{ineq:jmax:use}, we obtain \eqref{e:a_master_est_p_uniform_R} for $N''=0$. When $r\neq \infty$, we use $\left\| f \right\|_{L^r}\leq \left\| f \right\|_{L^\infty} | \{ \supp f \} |^{\sfrac 1r} $ and the demonstrated bound for $r=\infty$ to obtain \eqref{e:a_master_est_p_R} for the full range of $r$ and $N''=0$.  The estimate in \eqref{e:a_master_est_p_R:zeta} follows in the same way using \eqref{e:fat:pipe:estimates:1} for $p=\infty$ and \eqref{eq:checkerboard:derivatives:check} and the fact that $\zetab_\xi^{I,R}\leq 1$. Estimates for $N''=1$ are again nearly identical, and we omit further details.

Finally, we prove the pointwise estimates.  Recalling that the left-hand side of \eqref{eq:helping:me:pointwise} is supported inside the support of $\omega_{j,q}$ and using \eqref{eq:omega:cut:partition:unity} and \eqref{pt.est.pi.cutoff} proves the claim for $\diamond=\varphi$. Arguing analogously for $\diamond=R$ concludes the proof.
\end{proof}

\begin{corollary}[\bf Full velocity increment estimates]
\label{cor:corrections:Lp}
For $N,M\leq {\sfrac \Nfin 4}$, we have the estimates
\begin{subequations}\label{eq:w:oxi:ps}
\begin{align}
\norm{D^ND_{t,q}^M w_{(\xi),\diamond}^{(p),I}}_{L^r} 
&\lessg \left| \supp \left(\eta_{i,j,k,\xi,\vecl,\diamond}\zetab_\xi^{I,\diamond} \right) \right|^{\sfrac 1r} \delta_{q+\bn}^{\sfrac 12}\Gamma_q^{j+ {7}} r_q^{\frac2r-1} \lambda_{q+\bn}^N  \MM{M, \NindSmall, \tau_{q}^{-1}\Gamma_{q}^{i+ {13}}, \Tau_{q}^{-1} {\Ga_q^8}}\label{eq:w:oxi:est:master}
\\
\norm{D^ND_{t,q}^M w_{(\xi),\diamond}^{(p),I}}_{L^\infty}
&\lessg \Gamma_q^{\frac \badshaq 2 +  {10}}
r_{q}^{-1} \lambda_{q+\bn}^N \MM{M, \NindSmall, \tau_{q}^{-1}\Gamma_{q}^{i+ {13}}, \Tau_{q}^{-1} {\Ga_q^8}} \, .
\label{eq:w:oxi:unif:master}
\end{align}
\end{subequations}
Also, for $N,M\leq {\sfrac \Nfin 4}$, we have that
\begin{subequations}\label{eq:w:oxi:cs}
\begin{align}
\norm{D^ND_{t,q}^M w_{(\xi),\diamond}^{(c),I}}_{L^r}
&\lessg r_q \left| \supp \left(\eta_{i,j,k,\xi,\vecl,\diamond}\zetab_\xi^{I,\diamond} \right) \right|^{\sfrac 1r} \delta_{q+\bn}^{\sfrac 12} \Gamma_q^{j+ {7}} r_q^{\frac2r-1} \lambda_{q+\bn}^N  \MM{M, \NindSmall, \tau_{q}^{-1}\Gamma_{q}^{i+ {13}}, \Tau_{q}^{-1} {\Ga_q^8}} \label{eq:w:oxi:est:master:c}
\\
\norm{D^ND_{t,q}^M w_{(\xi),\diamond}^{(c),I}}_{L^\infty}
&\lessg \Gamma_q^{\frac \badshaq 2 +  {10}} \lambda_{q+\bn}^N \MM{M, \NindSmall, \tau_{q}^{-1}\Gamma_{q}^{i+ {13}}, \Tau_{q}^{-1} {\Ga_q^8}} \, .
\label{eq:w:oxi:unif:master:c}
\end{align}
\end{subequations}
\end{corollary}
\begin{proof}[Proof of Corollary~\ref{cor:corrections:Lp}]
Recalling the definition of $w_{(\xi),\diamond}^{(p),I}$ from \eqref{wqplusoneonephip} and \eqref{wqplusoneonep}, we shall prove \eqref{eq:w:oxi:est:master} by applying Lemma~\ref{l:slow_fast} with
\begin{align}
    N_*=M_*= \sfrac \Nfin 4 \, , \qquad f=a_{(\xi),\diamond}\left(\rhob_{(\xi)}^{\diamond}\zetab_{\xi}^{I, \diamond}\right) \circ \Phi_{(i,k)} \nabla\Phi_{(i,k)}^{-1} \, , \qquad \Phi=\Phi_{(i,k)} \, , \notag \\
    \lambda=\lambda_{q+\lfloor \sfrac \bn 2 \rfloor} \, , \qquad \tau^{-1}=\tau_q^{-1}\Gamma_q^{i+13} \, , \qquad \Tau = \Tau_q\Ga_q^{-8} \, , \qquad  \const_{f,R}=\left|\supp \eta_{(\xi),R}\zetab_{\xi}^{I,R}\right|^{\sfrac1r}\de_{q+\bn}^{\sfrac12}\Gamma_q^{j+ {7}} \notag \\
    \const_{f,\varphi}=\left|\supp \eta_{(\xi),\varphi}\zetab_{\xi}^{I,\varphi}\right|^{\sfrac1r} \de_{q+1}^{\sfrac12} r_q^{-\sfrac 13} \Gamma_q^{j+ {7}} \, , \qquad v = \hat u_q \, , \qquad \varphi = \WW_{(\xi),\diamond}^I \, , \qquad \mu = \lambda_{q+\lfloor \sfrac \bn 2 \rfloor} \Gamma_q \, , \notag \\
    \qquad \Upsilon = \Lambda = \lambda_{q+\bn} \, ,   \qquad
    \const_{\varrho,R} = r_q^{\frac 2r -1} \, , \qquad     \const_{\varrho,\varphi} = r_q^{\frac 2r - \frac 23} \, , \qquad N_t = \Nindt \, . \notag 
\end{align}
From \eqref{e:a_master_est_p}, Corollary~\ref{cor:deformation}, and \eqref{eq:checkerboard:derivatives:check}, we have that for $N,M \leq \sfrac \Nfin 4$,
\begin{align}
&\norm{D^ND_{t,q}^M \left( a_{(\xi),\diamond}\left(\rhob_{(\xi)}^{\diamond}\zetab_{\xi}^{I,\diamond}\right)\circ \Phiik \right) }_r \notag\\
&\qquad \qquad \lec \left|\supp \eta_{(\xi),\diamond}\zetab_{\xi}^{I,\diamond}\right|^{\sfrac1r}\de_{q+1}^{\sfrac12}\Gamma_q^{j+ {7}} \lambda_{q+\lfloor \sfrac \bn 2 \rfloor}^N \MM{M, \NindSmall, \tau_{q}^{-1}\Gamma_{q}^{i+ {13}}, \Tau_{q}^{-1} {\Ga_q^8}}\label{e:a_master_est_p_repeat}\\
&\left\| D^N D_{t,q}^M (D \Phi_{(i,k)})^{-1} \right\|_{L^\infty(\supp(\psi_{i,q}\tilde\chi_{i,k,q}))} 
\leq \Lambda_q^N \MM{M,\NindSmall,\Gamma_{q}^i\tau_q^{-1},\Tau_q^{-1}\Ga_q^{-1}},
\label{eq:D:N+1:Phi}\\
&{\norm{D^N\Phi_{(i,k)} }_{L^\infty(\supp(\psi_{i,q}\tilde\chi_{i,k,q} ))} } + \norm{D^N\Phi^{-1}_{(i,k)} }_{L^\infty(\supp(\psi_{i,q}\tilde\chi_{i,k,q} ))} 
\lesssim \Gamma_{q}^{-1} \Lambda_q^{N-1} \, , \label{eq:Lagrangian:Jacobian:2:repeat}
\end{align}
showing that \eqref{eq:slow_fast_0}, \eqref{eq:slow_fast_1}, and \eqref{eq:slow_fast_2} are satisfied. From Proposition~\ref{prop:pipeconstruction} and \ref{prop:pipe.flow.current}, we have that from $\WW_{(\xi),\diamond}^I$ is periodic to scale $\lambda_{q+\lfloor \sfrac \bn 2 \rfloor}\Gamma_q$, in addition to the estimates \eqref{e:pipe:estimates:2} and \eqref{e:pipe:estimates:2:current}, and so \eqref{eq:slow_fast_4} is satisfied for $\diamond=R,\varphi$.  Next, from \eqref{condi.Ndec0} and \eqref{condi.Nfin0}, the assumptions \eqref{eq:slow_fast_3} and \eqref{eq:slow_fast_3_a} are satisfied. We may thus apply Lemma \ref{l:slow_fast} to obtain that for $N,M \leq \sfrac \Nfin 4$, \eqref{eq:w:oxi:est:master} is satisfied.  Applying \eqref{ineq:jmax:use} then gives \eqref{eq:w:oxi:unif:master}. 

The argument for the corrector is similar, save for the fact that $D_{t,q}$ will land on $\nabla a_{(\xi)}$, and so we require an extra commutator estimate from Lemma~\ref{lem:cooper:2}, specifically Remark~\ref{rem:cooper:2:sum}. We omit the details of this commutator estimates and refer the reader to \cite[Corollary~8.2]{BMNV21}. However, we note that the gain in amplitude comes from the quotient of a spatial derivative cost of $\lambda_{q+\lfloor \sfrac \bn 2 \rfloor}$ on the low-frequency function, and a gain of $\lambda_{q+\bn}$ from \eqref{e:pipe:estimates:2} or \eqref{e:pipe:estimates:2:current}.  Using the definition of $r_q$ gives a net gain of $r_q\Gamma_q^{-1}$, concluding the proof.
\end{proof}

Now we estimate the mollified velocity increment given in~Definition~\ref{def:wqbn}.

\begin{lemma}[\bf Estimates on $\hat w_{q+\bn}$]\label{lem:mollifying:w}
We have that $\hat w_{q+\bn}$ satisfies the following properties.
\begin{enumerate}[(i)]
    \item\label{item:moll:vel:1} For all $N+M\leq 2\Nfin$, we have that
\begin{subequations}\label{eq:vellie:upgraded:statement}
\begin{align}
&\norm{ D^N D_{t,q+\bn-1}^M  \hat w_{q+\bn} }_{L^3(\supp \psi_{i,q+\bn-1})}  \notag\\
&\qquad\qquad \les \Gamma_q^{20} \delta_{q+\bn}^{\sfrac 12} r_q^{-\sfrac 13} \left(\lambda_{q+\bn}\Gamma_{q+\bn-1}\right)^N \MM{M, \NindRt, \Gamma_{q+\bn-1}^{i-1} \tau_{q+\bn-1}^{-1}, \Tau_{q+\bn-1}^{-1} \Ga_{q+\bn-1} }
\label{eq:vellie:inductive:dtq-1:upgraded:statement} \\
&\norm{ D^N D_{t,q+\bn-1}^M  \hat w_{q+\bn} }_{L^\infty(\supp \psi_{i,q+\bn-1})}  \notag\\
&\qquad\qquad \les \Gamma_q^{\sfrac{\badshaq}{2}+16} r_q^{-1} \left(\lambda_{q+\bn}\Gamma_{q+\bn-1}\right)^N \MM{M, \NindRt, \Gamma_{q+\bn-1}^{i-1} \tau_{q+\bn-1}^{-1}, \Tau_{q+\bn-1}^{-1} \Ga_{q+\bn-1} } \, . \label{eq:vellie:inductive:dtq-1:uniform:upgraded:statement}
\end{align}
\end{subequations}
    \item\label{item:moll:vel:2} For all $N+M\leq \sfrac{\Nfin}{4}$, we have that
    \begin{align}
    \norm{D^N D_{t,q+\bn-1}^M \left(w_{q+1}- \hat w_{q+\bn}\right)}_\infty &\lec \delta_{q+3\bn}^3 \Tau_\qbn^{25\Nindt} \left(\lambda_{q+\bn}\Gamma_{q+\bn-1}\right)^N  \notag\\
    &\qquad \qquad \times \MM{M, \NindRt, \tau_{q+\bn-1}^{-1}, \Tau_{q+\bn-1}^{-1} \Ga_{q+\bn-1} }  \, . \label{eq:diff:moll:vellie:statement}
\end{align}
\end{enumerate}
\end{lemma}
\begin{proof}[Proof of Lemma~\ref{lem:mollifying:w}]
We prove items~\eqref{item:moll:vel:1}--\eqref{item:moll:vel:2} in steps. First, we apply Corollary~\eqref{rem:summing:partition} with $\theta=1$, $\theta_1=0$, $\theta_2=1$, $H_{i,j,k,\xi,\vecl,I,\diamond}=w_{(\xi),\diamond}^{(\bullet),I}$ with $\bullet=p,c$, $p=3$, $\const_H=\delta_{q+\bn}^{\sfrac 12}\Gamma_q^{12}r_q^{-\sfrac 13}$, $N_*=M_*=\sfrac{\Nfin}{4}$, $M_t=\Nindt$, $N_x=\infty$, $\lambda=\Lambda=\lambda_{q+\bn}$, $\tau^{-1}=\tau_q^{-1}\Gamma_q^4$, $\Tau=\Tau_q$.  From the definition of $w_{\pxi,\diamond}^{(\bullet),I}$ and Corollary~\ref{cor:corrections:Lp}, we have that\eqref{eq:agg:assump:1}--\eqref{eq:agg:assump:3} are satisfied, and so from \eqref{eq:agg:conc:2}, we conclude that for $N,M\leq\sfrac{\Nfin}{4}$
\begin{align}
    \left\| \psi_{i,q} D^N \Dtq^M w_{q+1} \right\|_{3} \lesssim \Ga_q^{20} \delta_{q+\bn}^{\sfrac 12}r_q^{-\sfrac 13} \lambda_{q+\bn}^N \MM{N,\Nindt,\tau_q^{-1}\Ga_q^{i+14},\Tau_q^{-1}\Ga_q^8} \, . \label{eq:moll:vel:onesie}
\end{align}
In the case $p=\infty$, we may aggregate estimates from Corollary~\ref{cor:corrections:Lp} using the fact that only a finite, $q$-independent number of terms $w_{\pxi,\diamond}^{(\bullet),I}$ are non-zero at any fixed point in space-time to give the bound
\begin{align}
    \left\| \psi_{i,q} D^N \Dtq^M w_{q+1} \right\|_{\infty} \lesssim \Ga_q^{\frac{\badshaq}{2}+16} r_q^{-1} \lambda_{q+\bn}^N \MM{N,\Nindt,\tau_q^{-1}\Ga_q^{i+14},\Tau_q^{-1}} \, . \label{eq:moll:vel:twosie}
\end{align}
Next, from \eqref{eq:dodging:oldies}, which asserts that $\supp w_{q+1} \cap \supp \hat w_{q'}=\emptyset$ for $q+1\leq q'\leq q+\bn-1$, and from \eqref{eq:inductive:timescales} applied with $q'=q+\bn-1$ and $q''=q$, we may upgrade \eqref{eq:moll:vel:onesie}--\eqref{eq:moll:vel:twosie} to 
\begin{subequations}\label{eq:moll:vel:threesie}
\begin{align}
     \left\| D^N \Dtqnm^M w_{q+1} \right\|_{L^3(\supp \psi_{i,q+\bn-1})} &\lesssim \Ga_q^{20} \delta_{q+\bn}^{\sfrac 12}r_q^{-\sfrac 13} \lambda_{q+\bn}^N \MM{N,\Nindt,\tau_{q+\bn-1}^{-1}\Ga_{q+\bn-1}^{i-2},\Tau_q^{-1}} \\
     \left\| D^N \Dtqnm^M w_{q+1} \right\|_{L^\infty(\supp \psi_{i,q+\bn-1})} &\lesssim \Ga_q^{\frac{\badshaq}{2}+16} r_q^{-1} \lambda_{q+\bn}^N \MM{N,\Nindt,\tau_{q+\bn-1}^{-1}\Ga_{q+\bn-1}^{i-2},\Tau_q^{-1}} \, .
\end{align}
\end{subequations}
We now apply Proposition~\ref{lem:mollification:general} with the choices
\begin{align*}
    &p=3, \infty \, , \quad N_{\rm g}, N_{\rm c} \textnormal{ as in \eqref{i:par:12}} \, , \quad M_t = \Nindt\,  , \quad N_* = \sfrac{\Nfin}{4} \, , \\
    &N_\gamma = 2\Nfin \, , \quad  \Omega = \supp \psi_{i,q+\bn-1}\, , \quad v = \hat u_{q+\bn-1}\, , \quad i=i \, , \\
    &\lambda = \lambda_{q+\bn} \, , \quad \Lambda = \lambda_{q+\bn}\Ga_{q+\bn-1} \,, \quad \Ga = \Ga_{q+\bn-1}, \quad \tau = \tau_{q+\bn-1}\Gamma_{q+\bn-1}^{-2} \, , \quad \Tau = \Tau_{q+\bn-1}\, ,\\
    &f = w_{q+1} \, , \quad \const_{f,3} = \Ga_q^{20} \de_{q+\bn}^{\sfrac 12}r_q^{-\sfrac 13} \, , \quad \const_{f,\infty} = \tilde \const_f =  \Ga_q^{\sfrac{\badshaq}{2} + 16} r_q^{-1} \, , \quad \const_v = \Lambda_{q+\bn-1}^{\sfrac 12} \,  .
\end{align*}
From \eqref{i:par:12} and \eqref{v:global:par:ineq}, we have that \eqref{eq:moll:assumps:1} is satisfied. {From \eqref{eq:bobby:old}, we have that \eqref{moll.assum.v.est} is satisfied.} From \eqref{eq:moll:vel:threesie}, we have that \eqref{eq:moll:f:1} is satisfied. In order to verify \eqref{eq:moll:f:2}, we apply Remark~\ref{rem:upgrade.material.derivative.end} with the following choices. We set $p=\infty$, $N_x=N_t=\infty$, $N_*=\sfrac{\Nfin}{4}$, $\Omega=\T^3\times\R$, $v=w=\hat u_{q+\bn-1}$, $\const_w = \Gamma_{q+\bn-1}^{\imax+2}\delta_{q+\bn-1}^{\sfrac 12}\lambda_{q+\bn-1}^2$, $\lambda_w=\tilde\lambda_w=\Lambda_{q+\bn-1}$, $\mu_w=\tilde \mu_w = \Gamma_{q+\bn-1}^{-1}\Tau_{q+\bn-1}^{-1}$ in \eqref{eq:cooper:w}, while in \eqref{eq:cooper:2:v} and \eqref{eq:cooper:2:f} we set $v=\hat u_{q+\bn-1}$, $\const_v=\const_w$, $\lambda_v=\tilde\lambda_v=\Lambda_{q+\bn-1}$, $\mu_v=\tilde\mu_v = \Gamma_{q+\bn-1}^{-1}\Tau_{q+\bn-1}^{-1}$, $f=w_{q+1}$, $\const_f = \Gamma_q^{\sfrac{\badshaq}{2}+16}r_q^{-1}$, $\lambda_f=\tilde\lambda_f=\lambda_{q+\bn}$, $\mu_f=\tilde \mu_f= \Tau_{q}^{-1}$.  Then \eqref{eq:cooper:2:v} and \eqref{eq:cooper:2:f} are satisfied from \eqref{eq:nasty:D:vq:old} at level $q+\bn-1$, \eqref{eq:moll:vel:threesie}, \eqref{eq:imax:old}, and \eqref{v:global:par:ineq}. Next, \eqref{eq:cooper:w} is satisfied from \eqref{eq:bob:Dq':old} at level $q+\bn-1$. Thus from \eqref{eq:cooper:f:mat} and \eqref{v:global:par:ineq}, we obtain that 
\begin{align}\label{eq:w:f'ed:up}
\left\| D^N \partial_t^M w_{q+1} \right\|_\infty \lesssim \Ga_q^{\sfrac{\badshaq}{2} + 16} r_q^{-1} \lambda_{q+\bn}^N \Tau_{q+\bn-1}^{-M}
\end{align}
for $N+M \leq \sfrac{\Nfin}{4}$, thus verifying the final assumption \eqref{eq:moll:f:2} from Lemma~\ref{lem:mollification:general}.

We first apply \eqref{eq:moll:conc:1} to conclude that \eqref{eq:vellie:upgraded:statement} holds. Finally, we have from \eqref{eq:moll:conc:2} and \eqref{eq:darnit:2} that the difference $w_{q+1}- \hat w_{q+\bn}$ satisfies \eqref{eq:diff:moll:vellie:statement}.
\end{proof}

\opsection{Abstract construction of intermittent pressure}\label{opsection:pressure}

As in all convex integration schemes for the Euler equations, part of the goal of the pressure $\pi_\ell$ in our setting is to ensure that $R_\ell-\pi_\ell\Id$ is negative definite.  Then the low-frequency portion of $w_{q+1}\otimes w_{q+1}$, which is positive-definite, cancels $R_\ell- \pi_\ell\Id$ via Proposition~\ref{p:split}; see~\eqref{eq:cancellation:plus:pressure:nn}.  The simplest way to define $\pi_\ell$ for this purpose is to set $\pi_\ell \approx |R_\ell|$.  However, in order to ensure additionally that $\pi_\ell$ dominates the Reynolds stress and the gradient of velocity via estimates such as \eqref{eq:wednezday} (see also~\eqref{eq:pressure:inductive}--\eqref{eq:psi:q:q'}), one must include in the definition of $\pi_\ell$ derivative estimates on stresses and velocities, similar to the procedure described in Remark~\ref{rem:toozday}.  This is part of the goal of Lemmas~\ref{lem:pr.st} and~\ref{lem:pr.cu} and \texttt{Step 1} from Proposition~\ref{lem:pr.vel}.  The first of these two lemmas carries out this task for stress errors, while the latter does the same for current errors.  For example, Lemma~\ref{lem:pr.st} defines a positive scalar function $\sigma_S^+$ which dominates a stress error $S$ (for example part of $R_q^q$) via an estimate such as \eqref{est.S.by.pr}.  We also have that $\sigma_S^+$ dominates itself via an estimate such as \eqref{est.pr.S}.  

One should view $\sigma_S^+$ as essentially identical to $\delta_\qbn \Ga_q^{2j}$ from \eqref{heatsie:cheatsie}. However, due to the fact that $\sigma_S^+$ is positive, and no effort has been made yet to keep track of its active frequencies, one will never be able to effectively invert the divergence on any term containing $\sigma_S^+$.  For the method of proof described in Remark~\ref{rem:toozday}, or the iterations in \cite{NV22, BMNV21}, this was not an issue.  However, the relaxed local energy inequality~\ref{ineq:relaxed.LEI} throws a rather large wrench into this method.  Namely, the addition of $w_{q+1}$ into this equation will produce an error term of the form $(\pa_t + \hat u_q \cdot \nabla)|w_{q+1}|^2$, which can only be handled by inverting the divergence to create a new current error term.  This is the role of $\kappa_q^q$ in \eqref{ineq:relaxed.LEI}, which is essentially equal to $-\mathbb{P}_{\leq \la_q}(|w_{q+1}|^2)$.  Indeed then
$$  (\pa_t + \hat u_q\cdot\nabla)\left( \kappa_q^q + |w_{q+1}|^2 \right) \approx (\pa_t + \hat u_q\cdot\nabla)\left( \mathbb{P}_{> \la_q} \left( |w_{q+1}|^2 \right) \right) \, , $$
and so we can effectively invert the divergence on this term. But the appearance of the term $(\pa_t + \hat u_q\cdot \nabla)\kappa_q^q$ in \eqref{ineq:relaxed.LEI} means that one must have created current errors at earlier stages of the iteration by adding in $\div^{-1}\mathbb{P}_{\neq 0}\left((\pa_t + \hat u_q \cdot \nabla) \kappa_q^q \right)$.  Commuting for the moment the projection operator past the material derivative, this means that one must be able to estimate $\div^{-1}\left((\pa_t + \hat u_q \cdot \nabla) \mathbb{P}_{\neq 0} \kappa_q^q \right)$, which we refer to as a ``pressure current error.''\index{pressure current error} This will only be possible if we have accurate information on the frequency support of $\kappa_q^q$, i.e. accurate frequency support information on the scalar function $\sigma_S^+$ which is approximately equal to $-R_\ell + \pi_\ell \Id$.  Therefore, rather than simply adding $\sigma_S^+$ to dominate $-R_\ell + \pi_\ell \Id$, we must add $\sigma_S = \sigma_S^+ - \sigma_S^-$, where $\sigma_S$ is essentially mean-zero and $\sigma_S^-$ is low-frequency; see~\eqref{heatsie:stress}.  We then record an estimate of the form~\eqref{heatsie:low}, which asserts that $\sigma_S^-$ can be dominated by old intermittent pressure.  This is the second main goal of Lemmas~\ref{lem:pr.st} and \ref{lem:pr.cu}; to show that the low-frequency portion of the pressure increment can be absorbed by old intermittent pressure. 

Now that the pressure increment\index{pressure increment} $\sigma_S = \sigma_S^+ - \sigma_S^-$ defined in Lemma~\ref{lem:pr.st} is effectively mean-zero, we can apply a material derivative and invert the divergence.  This is the content of Proposition~\ref{lem.pr.invdiv2}, which contains several steps.  The first step is to use the inverse divergence from Proposition~\ref{prop:intermittent:inverse:div} to produce an error term $S$.  The second step is to apply Lemma~\ref{lem:pr.st} to produce a mean-zero pressure increment $\sigma_S$.  The final step is to apply a material derivative to $\sigma_S$ and invert the divergence.  Since this procedure has to be carried out for essentially every stress error term, one is forced to write a rather abstract, intricate result like Proposition~\ref{lem.pr.invdiv2} which can be applied over and over again.  Proposition~\ref{lem.pr.invdiv2.c} carries out a similar procedure, except for the current error.  Proposition~\ref{lem:pr.vel} creates the pressure increment for the velocity field, and since one need only apply this result one time at each step $q\mapsto q+1$, Proposition~\ref{lem:pr.vel} is analogous to the combination of Lemma~\ref{lem:pr.st} and Proposition~\ref{lem.pr.invdiv2} for the stress.  It would be reasonable for the reader to read only the proofs of Lemma~\ref{lem:pr.st} and Proposition~\ref{lem.pr.invdiv2}, as the remainder of the section is identical in character to these results.

\begin{lemma*}[\bf Pressure increment for stress error]\label{lem:pr.st} 
Let $v$ be an incompressible vector field on\index{intermittent pressure} $\R\times\T^3$. Denote its material derivative by $D_{t}=\partial_t + v\cdot\nabla$. We use large positive integers $N_\dagger \geq M_\dagger \gg M_t$ for counting derivatives and specify additional constraints that they must satisfy in assumptions \eqref{sample:item:assump:0}--\eqref{sample:item:assump:3}.

Suppose a stress error $S = H \, \rho\circ\Phi$ and a non-negative, continuous function $\pi$ are given such that the following hold.
\begin{enumerate}[(i)]
    \item\label{sample:item:assump:0}  There exist constants $\const_{G,p}$ and $\const_{\rho,p}$\footnote{In practice, $\const_{\rho,p}=\const_{\ast,p}\zeta^{-2}\xi \Lambda^\alpha$ from \eqref{eq:inverse:div:sub:1}. We shall also assume that these constants are ordered in the obvious way, i.e. $\const_{\bullet,\sfrac 32}\leq \const_{\bullet,\infty}$.} for $p=\sfrac 32$ and $p=\infty$ and frequency paramaters $\lambda,\Lambda,\nu,\nu'$ such that
    \begin{subequations}
    \begin{align}
    \norm{D^N D_{t}^M H}_p 
    &\les \const_{G,p} \la ^N\MM{M,M_{t},\nu,\nu'} \label{est.G} \\ 
    \left|D^N D_{t}^M H\right| 
    &\les \pi \la ^N\MM{M,M_{t},\nu,\nu'}\label{est.G.pt}\\
    \norm{ D^N  \rho}_{p} 
    &\lec \const_{\rho,p} \Lambda^N  \label{est.rho} \\
     \norm{S}_p &\les \const_{G,p} \const_{\rho,p} =: \delta_{S,p} \, . \label{est.S}
    \end{align}
    \end{subequations}
for all $N \leq N_\dagger$, $M\leq M_\dagger$.
\item\label{sample:item:assump:2} There exist a frequency parameter $\mu$, a parameter $\Gamma$ for measuring small losses in derivative costs,\footnote{In practice, $\Gamma=\Gamma_{q'}$ for some $q'$, which then makes $\Gamma$ a small power of $\lambda$ or $\Lambda$.} and a positive integer $\Ndec$ such that $\rho$ is $(\sfrac{\T}{\mu})^3$-periodic and $\la\ll \mu\leq \Lambda$, whereby we mean that
\begin{equation}\label{eq:sample:1:decoup}
     (\Lambda \Gamma)^{4}  \leq  \left( \frac{\mu}{4 \pi \sqrt{3} (\lambda\Gamma)}\right)^{\Ndec} \, .
\end{equation}
\item\label{sample:item:assump:1} Let
$\Phi$ be a volume preserving diffeomorphism of $\T^3$ such that $D_{t} \Phi =0$ and $\Phi$ is the identity at a time slice which intersects the support of $H$, and
\begin{subequations}
\begin{align}
\norm{D^{N+1}   \Phi}_{L^\infty\left(\supp H\right)} 
&+ \norm{D^{N+1}   \Phi^{-1}}_{L^\infty\left(\supp H\right)} 
\les \lambda^{N}
\label{eq:DDpsi:sample}\\
\norm{D^ND_{t}^M D  v}_{L^\infty\left(\supp H\right)}
&\les  \nu \lambda^N \MM{M,M_{t},\nu,\nu'}
 \label{eq:DDv:sample}
\end{align}
\end{subequations}
for all $N \leq N_\dagger$, $M\leq M_\dagger$.  
\item\label{sample:item:assump:3} There exist positive integers $\NcutLarge,\NcutSmall$ and a small parameter $\delta_{\rm tiny}\leq 1$ such that\footnote{The choice of $\NcutSmall$ is such that $\Gamma^{-\NcutSmall}$ can absorb a Sobolev loss from $H$ or $\rho$, or help absorb small remainder terms into the miniscule constant $\delta_{\rm tiny}$.}
\begin{subequations}\label{eq:sample:1:Ncut}
\begin{align}
\NcutSmall&\leq \NcutLarge \, , \label{eq:sample:1:Ncut:1} \\
 \left(\const_{G,\infty}+1\right)\left(\const_{\rho,\infty}+1\right)\Ga^{-\NcutSmall} &\leq \delta_{\rm tiny} \, , \const_{G,\sfrac32} \, ,
  \const_{\rho,\sfrac32} \, ,  \label{eq:sample:1:Ncut:2} \\
 2 \Ndec + 4 \leq N_\dagger - \NcutLarge, &\quad {\NcutSmall \leq}\, {M_t} \, . \label{eq:sample:1:Ncut:3}
\end{align}
\end{subequations}
\end{enumerate}
Then one can construct a pressure increment $\si_S=\si_{S}^+ - \si_S^-$ associated to the stress error $S$, where
\begin{subequations}\label{heatsie:stress}
\begin{align}
    \si_{S}
    &:=  \pr(H)
\left(\pr(\rho)\circ\Phi-\langle \pr(\rho)\rangle\right) \, , \label{d:st:pr:1} \\
\si_S^+ &:= \pr(H)\pr(\rho)\circ\Phi \, , \label{d:st:pr:2}
\end{align} 
\end{subequations}
and    
\begin{subequations}
\begin{align}
\pr(H) &:=\left( (\const_{G,\infty}\Ga^{-\NcutSmall})^2 + \sum_{N=0}^{\NcutLarge}\sum_{M=0}^{\NcutSmall} (\la \Ga)^{-2N}  (\nu\Ga) ^{-2M} |D^N D_{t}^M H|^2 \right)^\frac12-\const_{G,\infty}\Ga^{-\NcutSmall} \, , \label{d:st:pr:3} \\
\pr(\rho) &:= \left( (\const_{\rho,\infty}\Ga^{-\NcutSmall})^2 + \sum_{N=0}^{\NcutLarge} (\Lambda\Gamma )^{-2N}   |D^N \rho|^2 \right)^\frac12 
- \const_{\rho,\infty}\Ga^{-\NcutSmall} \, , \label{d:st:pr:4}
\end{align}
\end{subequations}
and which has the properties listed below.
\begin{enumerate}[(i)] 
\item\label{sample:item:2} $\si_S^+$ dominates derivatives of $S$ with suitable weights, so that for all $N \leq N_\dagger$ and $M\leq M_\dagger$,
\begin{align}\label{est.S.by.pr}
    \left|D^N D_{t}^M S\right|
    \lec (\si_S^+  + \delta_{\rm tiny}) (\Lambda\Ga)^N\MM{M,M_{t},\nu\Ga,\nu'\Ga} \, .
\end{align}
\item\label{sample:item:3} $\si_S^+$ dominates derivatives of itself with suitable weights, so that for all $N \leq N_\dagger-\NcutLarge$, $M\leq M_\dagger-\NcutSmall$,
\begin{align}\label{est.pr.S}
    \left|D^N D_{t}^M \si_S^+\right|
    \lec (\si_S^{+}+ \delta_{\rm tiny}) (\Lambda\Ga)^N\MM{M,M_{t}-\NcutSmall,\nu\Ga,\nu'\Ga} \, .
\end{align}
\item\label{sample:item:1} $\si_S^+$ and $\si_S^-$ have the same size as $S$, so that
    \begin{align}
        \norm{\si_S^+}_{p}\lec \de_{S,p}, \quad
        &\norm{\si_S^-}_{p}\lec 
        \de_{S,p} \, . \label{eq:sample:1:conc:1}
    \end{align}
Furthermore $\pr(H)$ and $\pr(\rho)$ have the same size as $H$ and $\rho$, so that for $N\leq N_\dagger -\NcutLarge$, $M\leq M_\dagger-\NcutSmall$, and $p=\sfrac 32,\infty$
\begin{align}\label{eq:sample:1:conc:2}
    \norm{D^N D_t^M \pr(H)}_{p}\lec \const_{G,p} (\lambda\Gamma)^N \MM{M,M_t-\NcutSmall,\nu\Gamma,\nu'\Gamma} \, , \qquad &\norm{D^N \pr(\rho)}_p \lec \const_{\rho,p} (\Lambda\Ga)^N \, . 
\end{align}
We note also that $\pr(\rho)$ is $\left(\sfrac{\T}{\mu}\right)^3$-periodic.
\item\label{sample:item:4} $\pi$ dominates $\si_S^-$ and $\pr(H)$ and their derivatives with suitable weights, so that for all $N \leq N_\dagger-\NcutLarge$ and $M\leq M_\dagger-\NcutSmall$,
\begin{subequations}\label{heatsie:low}
\begin{align}\label{est.pr.S-} 
    \left|D^N D_{t}^M \si_S^-\right|
    &\lec \pi \norm{\pr(\rho)}_1  (\la\Gamma)^N\MM{M,M_{t}-\NcutSmall,\nu\Gamma,\nu'\Gamma}\, ,\\
    \left|D^N D_{t}^M \pr(H)\right|
    &\lec \pi (\la\Gamma)^N\MM{M,M_{t}-\NcutSmall,\nu\Gamma,\nu'\Gamma}\, . \label{est.pr.H}
\end{align}
\end{subequations}
\item\label{sample:item:6} $\si_S^+$ and $\si_S^-$ are supported on $\supp(S)$ and $\supp(H)$, respectively.
\end{enumerate}
\end{lemma*}
\begin{proof}[Proof of Lemma~\ref{lem:pr.st}]
We break the proof into steps in which we prove each of the items \eqref{sample:item:2}--\eqref{sample:item:6}.

\noindent\textbf{Proof of \eqref{sample:item:2}:}  We first use \eqref{eq:DDpsi:sample} and $D_t\Phi=0$ from \eqref{sample:item:assump:1} and Lemma \ref{l:useful_ests} to deduce that for $N\leq N_\dagger$ and $M\leq M_\dagger$,
\begin{align}
|D^ND_t^M S |
&= |D^N( (D_t^M H) (\rho)\circ\Phi)|
\leq \sum_{N_1+N_2=N} |D^{N_1}(D_t^M H)| |D^{N_2}(\rho\circ\Phi))|\nonumber\\    
&\lec  \sum_{N_1+N_2=N} |D^{N_1}(D_t^M H)| \sum_{n_2=1}^{N_2} (\lambda\Gamma)^{N_2-n_2}\abs{(D^{n_2} \rho)\circ \Phi }\, . \label{exp.der.S}
\end{align}
Estimate \eqref{est.S.by.pr} will then follow from \eqref{exp.der.S} and the following claims;
\begin{subequations}
\begin{align}
\pr(H) &\lec \const_{G,\infty} \label{est.prH.inf.0} \\
\pr(\rho) &\lec \const_{\rho,\infty} \label{est.prrho.inf.0}\\
|D^{N_1}D_t^M H|
&\lec (\pr(H) + \const_{G,\infty}\Ga^{-\NcutSmall}) (\lambda\Ga)^{N_1} \MM{M, M_t, \nu\Ga, \nu'\Ga} \label{est.G.by.pr}\\
\la^{N_2-n_2}|D^{n_2} \rho|
&\lec (\pr(\rho) + \const_{\rho,\infty}\Ga^{-\NcutSmall}) (\Lambda\Ga)^{N_2} \label{est.rho.by.pr}
\end{align}
\end{subequations}
for any integers $0\leq N_1, n_2\leq N_\dagger$, $M\leq M_\dagger$.   Indeed, the above claims, \eqref{eq:sample:1:Ncut:1}--\eqref{eq:sample:1:Ncut:2}, and \eqref{exp.der.S} give that for $N\leq N_\dagger$ and $M\leq M_\dagger$,
\begin{align}
    |D^N D_t^M S| &\lesssim (\pr(H) + \const_{G,\infty}\Ga^{-\NcutSmall})(\pr(\rho)\circ \Phi + \const_{\rho,\infty}\Ga^{-\NcutSmall}) (\Lambda\Gamma)^{N} \MM{M,M_t,\nu\Gamma,\nu'\Gamma} \notag\\
    &\lesssim \left( \pr(H) \pr(\rho)\circ \Phi + \Gamma^{-\NcutSmall}\left(\const_{G,\infty}\pr(\rho)\circ\Phi + \const_{\rho,\infty} \pr(H) + \const_{G,\infty} \const_{\rho,\infty} \Gamma^{-\NcutSmall} \right) \right) \notag\\
    &\qquad \qquad \qquad \times (\Lambda\Gamma)^{N} \MM{M,M_t,\nu\Gamma,\nu'\Gamma} \notag\\
    &\lesssim (\sigma_s^+ + \delta_{\rm tiny} ) (\Lambda\Gamma)^{N} \MM{M,M_t,\nu\Gamma,\nu'\Gamma} \, . \notag 
\end{align}
The proofs of the claims are then given as follows. The first is immediate from the definition of $\pr(H)$ and the computation
\begin{align*}
    \pr(H) &\lesssim \const_{G,\infty} \\
    \impliedby \qquad \left(\pr(H) + \const_{G,\infty}\Ga^{-\NcutSmall} \right)^2 &\lesssim \const_{G,\infty}^2  \\
    \impliedby \qquad  (\lambda\Gamma)^{-2N}(\nu\Gamma)^{-2M}|D^ND_t^M H|^2 &\lesssim \const_{G,\infty}^2 \, ,
\end{align*}
which holds for $N\leq \NcutLarge$ and $M\leq \NcutSmall$ from \eqref{est.G}. A similar computation holds for $\pr(\rho)$. For the next two claims, if $M\leq \NcutSmall$ and $N_1, N_2\leq \NcutLarge$, an argument quite similar to the above computation shows that
\begin{subequations}
\begin{align}
    |D^{N_1}(D_t^M H)|
    &\lec (\pr(H) +\const_{G,\infty}\Ga^{-\NcutSmall})(\la \Ga)^{N_1}(\nu  \Ga)^M \, , \label{est.good.for.comm} \\
    \lambda^{N_2-n_2}\abs{(D^{n_2} \rho)\circ \Phi } &\lec(\Lambda \Ga)^{N_2} \left(\pr(\rho)\circ\Phi +\const_{\rho,\infty}\Ga^{-\NcutSmall}\right)\, .
    \label{est.G.low.der}
\end{align}
\end{subequations}
If however $M> \NcutSmall$, $N_1> \NcutLarge$, or $N_2> \NcutLarge$, we use \eqref{eq:sample:1:Ncut:1}--\eqref{eq:sample:1:Ncut:2} and \eqref{est.G} in the first two cases and \eqref{est.rho} in the third case to obtain, respectively, that
\begin{subequations}
\begin{align}
\norm{D^{N_1}(D_t^M H)}_{L^\infty}
&\lec \const_{G,\infty} \la^{N_1}
\MM{M,M_{t},\nu,\nu'}
\lec  \Ga^{-\NcutSmall}\const_{G,\infty} \la ^{N_1}\MM{M,M_{t},\nu\Ga,\nu'\Ga} \label{est.G.high.mt.der}\\
\norm{D^{N_1}(D_t^M H)}_{L^\infty}
&
\lec \Ga^{-\NcutSmall}\const_{G,\infty}(\la\Ga)^{N_1}  \MM{M,M_{t},\nu,\nu'} \label{est.G.high.der}\\
\la^{N_2-n_2}\norm{D^{n_2}\rho}_{L^\infty}
&\lec \Ga^{-\NcutSmall}\const_{\rho,\infty}(\Lambda \Ga)^{N_2}  \, , \label{est.rho.high.der}
\end{align}
\end{subequations}
concluding the proof of the claims and thus \eqref{est.S.by.pr}.

\noindent\textbf{Proof of \eqref{sample:item:3}: } We first show by induction that for integers $K\geq 0$ and $N,M$ such that $N+M=K,N\leq N_\dagger-\NcutLarge$, and $M\leq M_\dagger-\NcutSmall$,
\begin{align}\label{ind.G}
    |D^N D_t^M \pr(H)|\lec \left( \pr(H) + \const_{G,\infty}\Ga^{-\NcutSmall} \right) (\la\Ga)^N \MM{M, M_t-\NcutSmall, \nu\Ga, \nu'\Ga} \, . 
\end{align}
When $K=0$ the claim is immediate. Now, suppose by induction that \eqref{ind.G} holds true for any $K\leq K_0$, $K_0\in \mathbb{N}\cup\{0\}$.
To obtain \eqref{ind.G} for $K_0+1$, we first note that for $N'',M''$ such that $0<N''+M''$, $|D^{N''} D_t^{M''} \pr(H)| = |D^{N''} D_t^{M''} (\pr(H)+\const_{G,\infty}\Ga^{-\NcutSmall})|$. We then obtain the inequality
\begin{equation}
    \begin{split}\label{leib.prH}
\left|D^N D_{t}^M \pr(H) \right| &= \left|D^N D_{t}^M \left( \pr(H)+\const_{G,\infty}\Ga^{-\NcutSmall} \right) \right| \\
&\lec \frac{1}{\left|\pr(H)+\const_{G,\infty}\Ga^{-\NcutSmall}\right|} \bigg{[} \left|D^N D_{t}^M\left((\pr(H)+\const_{G,\infty}\Ga^{-\NcutSmall})^2 \right)\right| \\
&\qquad\qquad\qquad\qquad\qquad + \sum_{\substack{0 \leq N' \leq N\\0\leq M' \leq M\\0<N'+M' \leq K_0}}  \left|D^{N'} D_{t}^{M'} \pr(H)\right| \left| D^{N-N'} D_{t}^{M - M'} \pr(H)\right| \bigg{]} \, ,
    \end{split}
\end{equation}
which follows from Lemma \ref{lem:leib} with $p=2$ and the positivity of $\left|\pr(H)+\const_{G,\infty}\Ga^{-\NcutSmall}\right|$. Using the inductive assumption \eqref{ind.G}, which is valid since $0<N'+M'\leq K_0$, and \eqref{eq:sample:1:Ncut:2}, the second term can be controlled by
\begin{align}
& \frac{1}{\left|\pr(H)+\const_{G,\infty}\Ga^{-\NcutSmall}\right|} \left( \pr(H) + \const_{G,\infty}\Ga^{-\NcutSmall} \right)^2 (\lambda \Gamma)^N \MM{M,M_t-\NcutSmall,\Gamma \nu, \Gamma \nu'} \notag\\
&\qquad \les \left( \pr(H) + \const_{G,\infty}\Ga^{-\NcutSmall} \right) (\lambda \Gamma)^N \MM{M,M_t-\NcutSmall,\Gamma \nu, \Gamma \nu'} \, .
\label{est.prG.sec}
\end{align}
As for the first term, we have that
\begin{align}
&\frac{\left|D^N D_t^M \left((\pr(H)+\const_{G,\infty}\Ga^{-\NcutSmall})^2 \right) \right|}{\left|\pr(H)+\const_{G,\infty}\Ga^{-\NcutSmall}\right|} \notag\\
& \leq  \frac{1}{\left|\pr(H)+\const_{G,\infty}\Ga^{-\NcutSmall}\right|} \sum_{n=0}^{\NcutLarge}\sum_{m=0}^{\NcutSmall} (\la \Ga)^{-2n}  (\nu\Ga) ^{-2m} \left|D^N D_t^M \left|D^n D_{t}^m H\right|^2\right|\nonumber \\
& = \frac{1}{\left|\pr(H)+\const_{G,\infty}\Ga^{-\NcutSmall}\right|} \sum_{n=0}^{\NcutLarge}\sum_{m=0}^{\NcutSmall}
\sum_{\substack{0\leq N'\leq N\\0\leq M'\leq M}} 
(\la \Ga)^{-2n}  (\nu\Ga) ^{-2m}
\left|D^{N'} D_t^{M'} D^n D_{t}^m H\right| 
\left|D^{N-N'} D_t^{M-M'} D^n D_{t}^m H\right|  \, .
\label{exp.prG2}
\end{align}
To bound the quantity above, we first claim that for multi-indices $\alpha,\beta \in \mathbb{N}^k$ with $k\geq 2$, $|\alpha|\leq N_\dagger$, and $|\beta|\leq M_\dagger$, 
\begin{equation}
    \left| \prod_{i=1}^k D^{\alpha_i} D_t^{\beta_i} H \right| (x) \les \left(\pr(H)(x) + \const_{G,\infty}\Ga^{-\NcutSmall} \right) (\lambda \Gamma)^{|\alpha|} \MM{|\beta|,M_t,\nu\Gamma,\nu'\Gamma} \, . \label{eq:sample:mess:1}
\end{equation}
To prove this claim, let $\Omega(x){\subseteq\supp(H)}$ be a closed set containing $x$. Then applying Lemma~\ref{lem:cooper:2} with $p=\infty$, $N_t=M_t$, $N_*=N_\dagger$, $M_*=M_\dagger$, $\Omega=\Omega(x)$, $\const_v=\nu\lambda^{-1}$, $\lambda_v=\tilde\lambda_v=\lambda$, $\mu_v=\nu$, $\tilde \mu_v=\nu'$, $f=H$, $\const_f = \sup_{\Omega(x)}(\pr(H)+\const_{G,\infty}\Ga^{-\NcutSmall})$, $\lambda_f = \tilde\lambda_f = \lambda\Gamma$, $\mu_f =  \nu\Gamma$, and $\tilde \mu_f = \nu'\Gamma$, we have that \eqref{eq:cooper:2:v} is satisfied from \eqref{eq:DDv:sample}, and \eqref{eq:cooper:2:f} is satisfied by \eqref{est.G.by.pr} and the assumption on $|\alpha|,|\beta|$.  Then \eqref{eq:cooper:2:f:2} gives that 
\begin{equation}
    \left| \prod_{i=1}^k D^{\alpha_i} D_t^{\beta_i} H \right| (x) \les \left( \sup_{\Omega(x)}\pr(H) +\const_{G,\infty}\Ga^{-\NcutSmall} \right) (\lambda \Gamma)^{|\alpha|} \MM{|\beta|,M_t,\nu\Gamma,\nu'\Gamma} \, .
\end{equation}
Since $\Omega(x)$ is arbitrary and $\pr(H)$ is continuous, we have proven \eqref{eq:sample:mess:1}. Plugging the bound in \eqref{eq:sample:mess:1} into \eqref{exp.prG2}, we find that 
\begin{align}
\frac{\left|D^N D_t^M \left((\pr(H)+\const_{G,\infty}\Ga^{-\NcutSmall})^2 \right) \right|}{\left|\pr(H)+\const_{G,\infty}\Ga^{-\NcutSmall}\right|} &\lesssim \frac{1}{\left|\pr(H)+\const_{G,\infty}\Ga^{-\NcutSmall}\right|} \left( \pr(H)(x) + \const_{G,\infty}\Ga^{-\NcutSmall}\right)^2 \notag\\
&\qquad \qquad \times (\lambda\Gamma)^N \MM{M,M_t-\NcutSmall,\nu\Gamma,\nu'\Gamma} \, , \notag 
\end{align}
which matches the desired bound in \eqref{ind.G}.  This concludes the proof of \eqref{ind.G}.

Arguing in a similar way (in fact the proof is simpler since only spatial derivatives are required), we also have that for each integer $0\leq  N\leq N_\dagger-\NcutLarge$, 
\begin{subequations}
\begin{align}
    \left|D^N \pr(\rho)\right| &\lec \left(\pr(\rho) + \const_{\rho,\infty}\Ga^{-\NcutSmall} \right) \ (\Lambda\Ga)^N  \,,\label{ind.rho.0} \\
    \left|D^N (\pr(\rho)\circ \Phi)\right| &\lec \left(\pr(\rho)\circ \Phi + \const_{\rho,\infty}\Ga^{-\NcutSmall} \right) \ (\Lambda\Ga)^N  \, .\label{ind.rho}
\end{align}
\end{subequations}
Combining \eqref{ind.G}, \eqref{ind.rho}, and the choice of $\delta_{\rm tiny}$ from \eqref{eq:sample:1:Ncut:2}, we obtain the desired estimate \eqref{est.pr.S}.  

\noindent\textbf{Proof of \eqref{sample:item:1}: } 
Observe that by the construction of $\pr(H)$, \eqref{est.G}, and a computation similar to that used to produce \eqref{est.prH.inf.0}, we have $\norm{\pr(H)+\const_{G,\infty}\Ga^{-\NcutSmall}}_{p} \lec \const_{G,p}$ for $p=\sfrac 32,\infty$, and so $\norm{\pr(H)}_{p} \lec \const_{G,p}$. It follows from \eqref{ind.G} and \eqref{eq:sample:1:Ncut:2} that
\begin{align}\label{ind.G.p}
    \norm{D^N D_t^M \pr(H)}_{p} \lec \const_{G,p} (\la\Gamma)^N \MM{M,M_t-\NcutSmall,\nu\Gamma,\nu'\Gamma}
\end{align}
for $N\leq N_\dagger-\NcutLarge$ and $M\leq M_\dagger-\NcutSmall$. Similarly, by the construction of $\pr(\rho)$, \eqref{est.rho} and \eqref{ind.rho.0}, we have that $\norm{\pr(\rho)}_{p} \lec \const_{\rho, p}$, and so
\begin{align}\label{ind.rho.p}
    \norm{D^N \pr(\rho)}_{p}  \lec \const_{\rho,p} (\Lambda\Gamma)^N
\end{align}
for $N\leq N_\dagger-\NcutLarge$. Thus \eqref{eq:sample:1:conc:2} is verified. Also, by the construction of $\pr(\rho)$, its periodicity easily follows from \eqref{sample:item:assump:2}.
Next, we can immediately deduce from the definition of $\sigma_S^-$ the easier bound
\begin{align*}
\norm{\si^-_S}_{p} \lec 
\norm{\pr(H)}_{p}
\norm{\pr(\rho)}_1
\lec \const_{G,p} \const_{\rho,p} = \de_{S,p} \, .   
\end{align*}
In the case of $\si_S^+$ and $p=\sfrac 32$, we additionally apply Lemma \ref{l:slow_fast} by setting
\begin{align}
    &N_*=N_\dagger-\NcutLarge, \quad M_*= M_\dagger-\NcutSmall, \quad f=\pr(H), \quad \Phi=\Phi \, , \notag \\
    &\lambda=\la \Gamma, \quad \tau^{-1}=\nu \Gamma, \quad \Tau^{-1} = \nu'\Gamma,  \notag \\
    &\const_{f}=\const_{G,\sfrac 32}, \quad v = v, \quad \varrho = \pr(\rho), \quad \mu = \mu, \notag \\
    &\Upsilon = \Lambda = \Lambda\Gamma,   \quad
    \const_{\varrho} = \const_{\rho,\sfrac 32}, \quad  N_t = M_t-\NcutSmall \, . \notag 
\end{align}
Then \eqref{eq:slow_fast_0} is verified from \eqref{ind.G.p}, \eqref{eq:slow_fast_1}--\eqref{eq:slow_fast_2} follow from \eqref{eq:DDpsi:sample},  \eqref{eq:slow_fast_4} follows from \eqref{ind.rho.p} 
and the periodicity of $\pr(\rho)$, \eqref{eq:slow_fast_3} follows from \eqref{eq:sample:1:decoup}, and \eqref{eq:slow_fast_3_a} follows from \eqref{eq:sample:1:Ncut:3}. We then obtain from \eqref{eq:slow_fast_5} that
\begin{align*}
\norm{\si^+_S}_{\sfrac32} \lec \const_{G,\sfrac 32} \const_{\rho,\sfrac 32} = \de_{S,\sfrac32} \, .
\end{align*}
Finally, the estimate for $\norm{\si^+_S}_{\infty}$ is trivial, so that \eqref{eq:sample:1:conc:1} holds and \eqref{sample:item:1} is totally verified.

\noindent\textbf{Proof of \eqref{sample:item:4}:} We first prove \eqref{est.pr.H} by induction; namely,
for each integer $K=N+M\geq 0$, $N\leq N_\dagger-\NcutLarge$, $M\leq M_\dagger-\NcutSmall$,
\begin{align}\label{ind.G.pt}
    |D^ND_t^M \pr(H)|\lec \pi (\la\Ga)^N \MM{M, M_t-\NcutSmall, \nu\Ga, \nu\Ga'}  \, .
\end{align}
The proof uses an argument quite similar to the proof of \eqref{ind.G}.  The base case follows from writing that
\begin{align*}
    \pr(H) &\les \pi \\
    \iff \pr(H) + \const_{G,\infty}\Ga^{-\NcutSmall} &\leq C\pi + \const_{G,\infty}\Ga^{-\NcutSmall} \\
    \impliedby \left( \pr(H) + \const_{G,\infty}\Ga^{-\NcutSmall} \right)^2 &\leq C^2\pi^2 + \const_{G,\infty}^2\Ga^{-2\NcutSmall} \, ,
\end{align*}
for some absolute constant $C=C(\NcutSmall, \NcutLarge)$ which can be seen to hold from the definition of $\pr(H)$ and \eqref{est.G.pt}. For the inductive step, we argue starting from \eqref{leib.prH}, although with slightly different steps to follow. Using the inductive assumption from \eqref{ind.G.pt} to control one term and the bound \eqref{ind.G} to control the other term, and \eqref{eq:sample:1:Ncut:2}, we have that the second term from \eqref{leib.prH} may be bounded by
\begin{align}
& \frac{1}{\left|\pr(H)+\const_{G,\infty}\Ga^{-\NcutSmall} \right|}  \pi \left( \pr(H) +  \const_{G,\infty}\Ga^{-\NcutSmall} \right) (\lambda \Gamma)^N \MM{M,M_t-\NcutSmall,\Gamma \nu, \Gamma \nu'} \notag\\
&\qquad \les \pi (\lambda \Gamma)^N \MM{M,M_t-\NcutSmall,\Gamma \nu, \Gamma \nu'} \, .
\label{est.prG.sec.redux}
\end{align}
Thus it remains to control the first term from \eqref{leib.prH}.  Towards this end, we claim that for multi-indices $\alpha,\beta\in\mathbb{N}^k$ with $k\geq 2$, $|\alpha|\leq N_\dagger$, and $|\beta|\leq M_\dagger$, 
\begin{equation}
    \left| \prod_{i=1}^k D^{\alpha_i} D_t^{\beta_i} H \right| (x) \les \pi(x) (\lambda \Gamma)^{|\alpha|} \MM{|\beta|,M_t,\nu\Gamma,\nu'\Gamma} \, . \label{eq:sample:mess:1.redux}
\end{equation}
We apply Lemma~\ref{lem:cooper:2} with precisely the same choices as in the proof of \eqref{eq:sample:mess:1}, save for the choice of $\const_f = \sup_{\Omega(x)}\pi$.  Then \eqref{eq:cooper:2:v} is satisfied from \eqref{eq:DDv:sample}, and \eqref{eq:cooper:2:f} is satisfied by \eqref{est.G.pt}. Then applying \eqref{eq:cooper:2:f:2}, shrinking $\Omega(x)$ to a point, and using the continuity of $\pi$ provides \eqref{eq:sample:mess:1.redux}. Plugging this bound into \eqref{exp.prG2} and using \eqref{eq:sample:mess:1} and \eqref{eq:sample:1:Ncut:2}, we find that for $N\leq N_\dagger-\NcutLarge$ and $M\leq M_\dagger-\NcutSmall$,
\begin{align}
&\frac{\left|D^N D_t^M \left((\pr(H)+\const_{G,\infty}\Ga^{-\NcutSmall})^2 \right) \right|}{\left|\pr(H)+\const_{G,\infty}\Ga^{-\NcutSmall}\right|} \notag\\
&\qquad \les \frac{1}{\left|\pr(H)+\const_{G,\infty}\Ga^{-\NcutSmall}\right|} \pi \left( \pr(H) + \const_{G,\infty}\Gamma^{-\NcutSmall} \right) (\lambda\Gamma)^{N} \MM{M,M_t-\NcutSmall,\nu\Gamma,\nu'\Gamma} \notag\\
&\qquad \les \pi (\lambda\Gamma)^{N} \MM{M,M_t-\NcutSmall,\nu\Gamma,\nu'\Gamma} \, , \notag 
\end{align}
which combined with \eqref{est.prG.sec.redux} concludes the proof of \eqref{est.pr.H}.  To prove \eqref{est.pr.S-}, we use \eqref{est.pr.H} and the definition of $\sigma_S^-$.

\noindent\textbf{Proof of \eqref{sample:item:6}:} By the definition of $\pr(H)$ and $\pr(\rho)$, it is easy to see that $\supp(\pr(H))\subseteq \supp(H)$ and $\supp(\pr(\rho))\subseteq \supp(\rho)$, and so \eqref{sample:item:6} is verified.
\end{proof}

\begin{lemma*}[\bf Pressure increment for current error]\label{lem:pr.cu} 
Let $v$ be an incompressible vector field on $\R\times\T^3$. Denote its material derivative by $D_{t}=\partial_t + v\cdot\nabla$. We use large positive integers $N_*\geq M_*\gg M_t$ for counting derivatives and specify additional constraints that they must satisfy in assumptions \eqref{sample:item:assump:0:c}--\eqref{sample:item:assump:3:c}.

Suppose a current error $\phi= H \, \rho\circ\Phi$ and a non-negative, continuous function $\pi$ are given such that the following hold.
\begin{enumerate}[(i)]
    \item\label{sample:item:assump:0:c} There exist constants $\const_{G,p}$ and $\const_{\rho,p}$ for $p=1,\infty$, frequency parameters $\lambda,\Lambda,\nu,\nu'$, and intermittency parameters $0<r_G,r_{\phi}\leq 1$ such that
    \begin{subequations}
    \begin{align}
    \norm{D^N D_{t}^M H}_p 
    &\les \const_{G,p} \la ^N\MM{M,M_{t},\nu,\nu'} \label{est.G.c} \\ 
    \left|D^N D_{t}^M H\right| 
    &\les \pi^{\sfrac 32} r_G^{-1} \la ^N\MM{M,M_{t},\nu,\nu'}\label{est.G.c.pt}\\
    \norm{ D^N  \rho}_{p} 
    &\lec \const_{\rho,p} \Lambda^N  \label{est.rho.c} \\
    \left\| \phi \right\|_p &\lesssim \const_{G,p} \const_{\rho,p} =: \delta_{\phi,p}^{\sfrac 32} r_{\phi}^{-1} \label{est.phi.c}
    \end{align}
    \end{subequations}
for all $N \leq N_*$, $M\leq M_*$.
\item\label{sample:item:assump:2:c} There exist a frequency parameter $\mu$, a parameter $\Gamma$ for measuring small losses in derivative costs, and a positive integer $\Ndec$ such that $\rho$ is $(\sfrac{\T}{\mu})^3$-periodic and $\la\ll \mu\leq\Lambda$, whereby we mean that
\begin{equation}\label{eq:sample:1:decoup:c}
     (\Lambda \Gamma)^{4}  \leq  \left( \frac{\mu}{4 \pi \sqrt{3} (\lambda\Gamma)}\right)^{\Ndec} \, .
\end{equation}
\item\label{sample:item:assump:1:c} Let $\Phi$ be a volume preserving diffeomorphism of $\T^3$ such that $D_{t} \Phi =0$ and $\Phi$ is the identity at a time slice which intersects the support of $H$, and 
\begin{subequations}
\begin{align}
 \norm{D^{N+1} \Phi}_{L^\infty\left(\supp H\right)} 
 &+ \norm{D^{N+1}   \Phi^{-1}}_{L^\infty\left(\supp H\right)} 
 \les \lambda^{N}
 \label{eq:DDpsi:sample:c}\\
 \norm{D^ND_{t}^M D  v}_{L^\infty\left(\supp H\right)}
&\les  \nu\lambda^{N}\MM{M,M_{t},\nu,\nu'}
\label{eq:DDv:sample:c}
\end{align}
\end{subequations}
for all $N \leq N_*$, $M\leq M_*$.  
\item\label{sample:item:assump:3:c} There exist positive integers $\NcutLarge,\NcutSmall$ and a small parameter $\delta_{\rm tiny}\leq 1$ such that
\begin{subequations}\label{eq:sample:1:Ncut:c}
\begin{align}
    \NcutLarge &\geq \NcutSmall \label{eq:sample:1:Ncut:1:c} \\
    \left( \const_{G,\infty}+1 \right)\left( \const_{\rho,\infty} +1 \right) \Gamma^{-\NcutSmall} &\leq \delta_{\rm tiny}^{\sfrac 32} \, , \const_{G,1} \, , \const_{\rho,1} \, , \label{eq:sample:1:Ncut:2:c} \\
    2\Ndec + 4 \leq N_*- \NcutLarge - 4,  &\quad \NcutSmall \leq \, M_t \, . \label{eq:sample:1:Ncut:3:c}
\end{align}
\end{subequations}
\end{enumerate}

Then one can construct a pressure increment $\si_\phi$ associated to the current error $\phi$, where
\begin{subequations}
\begin{align}
    \si_{\phi} &= r_\phi^{\sfrac23} \pr(H) \left(\pr(\rho)\circ\Phi-\langle 
    \pr(\rho)\rangle\right) \, , \label{eq:cu:d:1} \\
    \sigma_\phi^+ &:= r_\phi^{\sfrac 23} \pr(H) \pr(\rho)\circ \Phi \, , \label{eq:cu:d:2}
\end{align}
\end{subequations}
and
\begin{subequations}
\begin{align}
\pr(H) &:= \left( \left(\const_{G,\infty}\Ga^{-\NcutSmall}\right)^2 + \sum_{N=0}^{\NcutLarge}\sum_{M=0}^{\NcutSmall} (\la \Ga)^{-2N}  (\nu\Ga) ^{-2M} |D^N D_{t}^M H|^2 \right)^\frac13 - \left(\const_{G,\infty}\Ga^{-\NcutSmall}\right)^{\sfrac23} \, , \label{eq:cu:d:3} \\
\pr(\rho) &:= \left( \left(\const_{\rho,\infty}\Ga^{-\NcutSmall}\right)^2 + \sum_{N=0}^{\NcutLarge} (\Lambda\Gamma )^{-2N}   |D^N \rho|^2 \right)^\frac13 -  \left(\const_{\rho,\infty}\Ga^{-\NcutSmall}\right)^{\sfrac 23} \, , \label{eq:cu:d:4}
\end{align}
\end{subequations}
and which has the properties listed below.
\begin{enumerate}[(i)]
\item\label{sample:item:2:c} $\si_\phi^+$ dominates derivatives of $\phi$ with suitable weights, so that for all $N\leq N_*$ and $M\leq M_*$,
\begin{align}\label{est.phi.by.pr}
    \left|D^N D_{t}^M \phi\right|
    \lec \left( (\si_\phi^+)^{\sfrac 32} r_{\phi}^{-1}  + \de_{\rm tiny} \right)  (\Lambda\Ga)^N\MM{M,M_{t},\nu\Ga,\nu'\Ga} \, .
\end{align}
\item\label{sample:item:3:c} $\si_\phi^+$ dominates derivatives of itself with suitable weights, so that for all $N \leq N_*-\NcutLarge$, $M\leq M_*-\NcutSmall$,
\begin{align}\label{est.pr.phi}
    \left| D^N D_{t}^M \si_\phi^+\right|
    \lec \left(\si_\phi^{+}+ \de_{\rm tiny}\right) (\Lambda\Ga)^N\MM{M,M_{t},\nu,\nu'} \, .
\end{align}
\item\label{sample:item:1:c} $\si_\phi^+$ and $\si_\phi^-$ have size comparable to $\phi$, so that
\begin{subequations}
\begin{align}
        \norm{\si_\phi^+}_{\sfrac32}\lec \delta_{\phi,1} \,, \qquad
        \norm{\si_\phi^-}_{\sfrac32}\lec \delta_{\phi,1} \, , \label{eq:sample:1:conc:1:c} \\
        \norm{\si_\phi^+}_{\infty}\lec \de_{\phi,\infty}\, , \qquad
        \norm{\si_\phi^-}_{\infty}\lec \de_{\phi,\infty} \, . \label{eq:sample:1:conc:1:c:infty}
\end{align}
\end{subequations}
Furthermore, $\pr(H)$ and $\pr(\rho)$ have size comparable to $H$ and $\rho$, respectively, so that for all $N\leq N_*-\NcutLarge$ and $M\leq M_*-\NcutSmall$,
\begin{subequations}
\begin{align}
    \norm{D^N D_t^M \pr(H)}_{\sfrac32}&\lec \const_{G,1}^{\sfrac23} (\lambda\Gamma)^N \MM{M,M_t-\NcutSmall,\nu\Gamma,\nu'\Gamma} \, , \qquad &\norm{D^N \pr(\rho)}_{\sfrac32}\lec \const_{\rho,1}^{\sfrac23} (\Lambda\Gamma)^N \, , \label{eq:sample:1:conc:2:c} \\
    \norm{D^N D_t^M \pr(H)}_{\infty}&\lec \const_{G,\infty}^{\sfrac23} (\lambda\Gamma)^N \MM{M,M_t-\NcutSmall,\nu\Gamma,\nu'\Gamma} \, , \qquad &\norm{D^N \pr(\rho)}_{\infty}\lec \const_{\rho,\infty}^{\sfrac23} (\Lambda\Gamma)^N \, , \label{eq:sample:1:conc:2:c:infty}
\end{align}
\end{subequations}
We note also that $\pr(\rho)$ is $\left( \sfrac{\T}{\mu}\right)^3$-periodic.
\item\label{sample:item:4:c} $\pi$ dominates $\si_\phi^-$ and $\pr(H)$ and their derivatives with suitable weights, so that for all $N\leq N_*-\NcutLarge$ and $M\leq M_* - \NcutSmall$,
\begin{subequations}
\begin{align}\label{est.pr.phi-}
    \left| D^N D_{t}^M \si_\phi^-\right|
    &\lec \left( \frac{r_\phi}{r_G} \right)^{\sfrac 23} \pi \norm{\pr(\rho)}_1  (\la\Gamma)^N\MM{M,M_{t}-\NcutSmall,\nu\Gamma,\nu'\Gamma}\, ,\\
    \left|D^N D_{t}^M \pr(H)\right|
    &\lec r_G^{-\sfrac 23} \pi (\la\Gamma)^N\MM{M,M_{t}-\NcutSmall,\nu\Gamma,\nu'\Gamma}\, . \label{est.pr.H.c}
\end{align}
\end{subequations}
\item\label{sample:item:6:c} $\si_\phi^+$ and $\si_\phi^-$ are supported on $\supp(\phi)$ and $\supp(H)$, respectivly. 
\end{enumerate}
\end{lemma*}

\begin{proof}[Proof of Lemma~\ref{lem:pr.cu}]
We break the proof into steps in which we prove each of the items \eqref{sample:item:2:c}--\eqref{sample:item:6:c}. The proof follows quite closely the proof of Lemma~\ref{lem:pr.st}, save for various rescalings related to the different scalings for current errors versus stress errors.

\noindent\textbf{Proof of \eqref{sample:item:2:c}:}  We first use \eqref{eq:DDpsi:sample:c} and $D_t\Phi=0$ from \eqref{sample:item:assump:1:c} and Lemma \ref{l:useful_ests} to deduce that for $N\leq N_*$ and $M\leq M_*$,
\begin{align}
|D^ND_t^M \phi |
&= |D^N( (D_t^M H) (\rho)\circ\Phi)|
\leq \sum_{N_1+N_2=N} |D^{N_1}(D_t^M H)| |D^{N_2}(\rho\circ\Phi))|\nonumber\\    
&\lec  \sum_{N_1+N_2=N} |D^{N_1}(D_t^M H)| \sum_{n_2=1}^{N_2} (\lambda\Gamma)^{N_2-n_2}\abs{(D^{n_2} \rho)\circ \Phi }\, . \label{exp.der.c}
\end{align}
Estimate \eqref{est.phi.by.pr} will then follow from \eqref{exp.der.c} and the following claims;
\begin{subequations}
\begin{align}
\pr(H) &\lec \const_{G,\infty}^{\sfrac 23} \label{eq:something:I:need} \\
\pr(\rho) &\lec \const_{\rho,\infty}^{\sfrac 23}\label{eq:something:else:I:need} \\
|D^{N_1}D_t^M H|
&\lec \left(\pr^{\sfrac 32}(H) + \const_{G,\infty}\Ga^{-\NcutSmall}\right) (\lambda\Ga)^{N_1} \MM{M, M_t, \nu\Ga, \nu'\Ga} \label{est.G.by.pr.c}\\
\la^{N_2-n_2}|D^{n_2} \rho|
&\lec \left(\pr^{\sfrac 32}(\rho) + \const_{\rho,\infty}\Ga^{-\NcutSmall}\right) (\Lambda\Ga)^{N_2} \label{est.rho.by.pr.c}
\end{align}
\end{subequations}
for any integers $0\leq N_1, n_2\leq N_*$, $M\leq M_*$.  Indeed, the above claims, \eqref{eq:sample:1:Ncut:1:c}--\eqref{eq:sample:1:Ncut:2:c}, and \eqref{exp.der.c} give that for $N\leq N_*$ and $M\leq M_*$,
\begin{align}
    \left| D^N D_t^M \phi \right| &\lesssim \left(\pr^{\sfrac 32}(H) + \const_{G,\infty}\Ga^{-\NcutSmall}\right) \left(\pr^{\sfrac 32}(\rho)\circ \Phi + \const_{\rho,\infty}\Ga^{-\NcutSmall}\right) (\Lambda\Gamma)^{N} \MM{M,M_t,\nu\Gamma,\nu'\Gamma} \notag\\
    &\lesssim \left( \left(\pr(H) \pr(\rho)\circ \Phi\right)^{\sfrac 32} + \Gamma^{-\NcutSmall}\left(\const_{G,\infty}\pr^{\sfrac 32}(\rho)\circ \Phi + \const_{\rho,\infty} \pr^{\sfrac 32}(H) + \const_{G,\infty} \const_{\rho,\infty} \Gamma^{-\NcutSmall} \right) \right) \notag\\
    &\qquad \qquad \qquad \times (\Lambda\Gamma)^{N} \MM{M,M_t,\nu\Gamma,\nu'\Gamma} \notag\\
    &\lesssim \left( (\sigma_s^+)^{\sfrac 32} r_\phi^{-1} + \delta_{\rm tiny} \right) (\Lambda\Gamma)^{N} \MM{M,M_t,\nu\Gamma,\nu'\Gamma} \, . \notag
\end{align}
The proofs of the claims are then given as follows. The first is immediate from the definition of $\pr(H)$ and the computation
\begin{align*}
    \pr(H) &\lesssim \const_{G,\infty}^{\sfrac 23} \\
    \impliedby \qquad \left(\pr(H) + \left(\const_{G,\infty}\Gamma^{-\NcutSmall}\right)^{\sfrac 23} \right)^3 &\lesssim \const_{G,\infty}^2 \\
    \impliedby \qquad  (\lambda\Gamma)^{-2N}(\nu\Gamma)^{-2M}|D^ND_t^M H|^2 &\lesssim \const_{G,\infty}^2 \, ,
\end{align*}
which holds for $N\leq \NcutLarge$ and $M\leq \NcutSmall$ from \eqref{est.G.c}. A similar computation holds for $\pr(\rho)$. Next, if $M\leq \NcutSmall$ and $N_1, N_2\leq \NcutLarge$, a computation similar to the one above shows that
\begin{subequations}
\begin{align}
    |D^{N_1}(D_t^M H)|
    &\lec \left(\pr(H) + \left(\const_{G,\infty}\Ga^{-\NcutSmall}\right)^{\sfrac 23}\right)^{\sfrac 32} (\la \Ga)^{N_1}(\nu  \Ga)^M \, , \label{est.good.for.comm.c} \\
    \lambda^{N_2-n_2}\abs{(D^{n_2} \rho)\circ \Phi } &\lec(\Lambda \Ga)^{N_2} \left(\pr(\rho)\circ\Phi + \left(\const_{\rho,\infty}\Ga^{-\NcutSmall}\right)^{\sfrac 23}\right)^{\sfrac 32}\, .
    \label{est.G.low.der.c}
\end{align}
\end{subequations}
If however $M> \NcutSmall$, $N_1> \NcutLarge$, or $N_2> \NcutLarge$, we use \eqref{eq:sample:1:Ncut:1:c}--\eqref{eq:sample:1:Ncut:2:c} and \eqref{est.G.c} in the first two cases and \eqref{est.rho.c} in the third case to obtain, respectively, that
\begin{subequations}
\begin{align}
\norm{D^{N_1}(D_t^M H)}_{L^\infty}
&\lec \const_{G,\infty} \la^{N_1}
\MM{M,M_{t},\nu,\nu'}
\lec  \Ga^{-\NcutSmall}\const_{G,\infty} \la ^{N_1}\MM{M,M_{t},\nu\Ga,\nu'\Ga} \label{est.G.high.mt.der.c}\\
\norm{D^{N_1}(D_t^M H)}_{L^\infty}
&
\lec \Ga^{-\NcutSmall}\const_{G,\infty}(\la\Ga)^{N_1}  \MM{M,M_{t},\nu,\nu'} \label{est.G.high.der.c}\\
\la^{N_2-n_2}\norm{D^{n_2}\rho}_{L^\infty}
&\lec \Ga^{-\NcutSmall}\const_{\rho,\infty}(\Lambda \Ga)^{N_2}  \, , \label{est.rho.high.der.c}
\end{align}
\end{subequations}
concluding the proof of the claims and thus of \eqref{est.phi.by.pr}.

\noindent\textbf{Proof of \eqref{sample:item:3:c}: } We first show by induction that for integers $K\geq 0$ and $N,M$ such that $N+M=K,N\leq N_*-\NcutLarge$, and $M\leq M_*-\NcutSmall$,
\begin{align}\label{ind.G.c}
    |D^N D_t^M \pr(H)|\lec \left( \pr(H) + (\const_{G,\infty}\Gamma^{-\NcutSmall})^{\sfrac 23} \right) (\la\Ga)^N \MM{M, M_t-\NcutSmall, \nu\Ga, \nu'\Ga} \, . 
\end{align}
When $K=0$ the claim is immediate. Now, suppose by induction that \eqref{ind.G.c} holds true for any $K\leq K_0$, $K_0\in \mathbb{N}\cup\{0\}$.
To obtain \eqref{ind.G.c} for $K_0+1$, we first note that for $N'',M''$ such that $0<N''+M''$, $|D^{N''} D_t^{M''} \pr(H)| = |D^{N''} D_t^{M''} (\pr(H)++ (\const_{G,\infty}\Gamma^{-\NcutSmall})^{\sfrac 23})|$. We then obtain the inequality
\begin{align}
\left|D^N D_{t}^M \pr(H) \right| &= \left|D^N D_{t}^M \left( \pr(H)+ (\const_{G,\infty}\Gamma^{-\NcutSmall})^{\sfrac 23} \right) \right| \notag \\
&\lec \frac{1}{\left|\pr(H)+ (\const_{G,\infty}\Gamma^{-\NcutSmall})^{\sfrac 23}\right|^2} \bigg{[} \left|D^N D_{t}^M \left((\pr(H)+ (\const_{G,\infty}\Gamma^{-\NcutSmall})^{\sfrac 23})^3 \right)\right| \notag \\
&\qquad \qquad \qquad \qquad \qquad + \sum_{\left\{\substack{\alpha , \beta  \, : \, \sum_{i=1}^3 \alpha_i = N \, , \\ \sum_{i=1}^3 \beta_i = M \, , \\  \alpha_i+\beta_i < N+M \, \forall \, i }\right\}} \prod_{i=1}^3 \left|D^{\alpha_i} D_{t}^{\beta_i}\left( \pr(H)+ (\const_{G,\infty}\Gamma^{-\NcutSmall})^{\sfrac 23} \right)\right| \bigg{]} \, , \label{eq:liebniz:current}
\end{align}
which follows from Lemma \ref{lem:leib} with $p=3$ and the positivity of $\left|\pr(H)+ (\const_{G,\infty}\Gamma^{-\NcutSmall})^{\sfrac 23}\right|$. Using the inductive assumption \eqref{ind.G.c}, which is valid since $0<N'+M'\leq K_0$, and \eqref{eq:sample:1:Ncut:2:c}, the second term can be controlled by
\begin{align}
& \frac{1}{\left|\pr(H)+ (\const_{G,\infty}\Gamma^{-\NcutSmall})^{\sfrac 23} \right|^2} \left( \pr(H) + (\const_{G,\infty}\Gamma^{-\NcutSmall})^{\sfrac 23} \right)^3 (\lambda \Gamma)^N \MM{M,M_t-\NcutSmall,\Gamma \nu, \Gamma \nu'} \notag\\
&\qquad \les \left( \pr(H) + (\const_{G,\infty}\Gamma^{-\NcutSmall})^{\sfrac 23} \right) (\lambda \Gamma)^N \MM{M,M_t-\NcutSmall,\Gamma \nu, \Gamma \nu'} \, .
\label{est.prG.sec.c}
\end{align}
As for the first term, we have that
\begin{align}
&\frac{\left|D^N D_t^M \left((\pr(H)+ (\const_{G,\infty}\Gamma^{-\NcutSmall})^{\sfrac 23})^3 \right) \right|}{\left|\pr(H)+ (\const_{G,\infty}\Gamma^{-\NcutSmall})^{\sfrac 23}\right|^2} \notag\\
&\qquad \leq  \frac{1}{\left|\pr(H)+ (\const_{G,\infty}\Gamma^{-\NcutSmall})^{\sfrac 23}\right|^2} \sum_{n=0}^{\NcutLarge}\sum_{m=0}^{\NcutSmall} (\la \Ga)^{-2n}  (\nu\Ga) ^{-2m} \left|D^N D_t^M \left|D^n D_{t}^m H\right|^2\right|\nonumber \\
&\qquad = \frac{1}{\left|\pr(H)+ (\const_{G,\infty}\Gamma^{-\NcutSmall})^{\sfrac 23}\right|^2} \sum_{n=0}^{\NcutLarge}\sum_{m=0}^{\NcutSmall}
\sum_{\substack{0\leq N'\leq N\\0\leq M'\leq M}} 
(\la \Ga)^{-2n}  (\nu\Ga)^{-2m}
\left|D^{N'} D_t^{M'} D^n D_{t}^m H\right| \notag\\
&\qquad \qquad \qquad \qquad \qquad \qquad \qquad \qquad \qquad \qquad \qquad \qquad  \times 
\left|D^{N-N'} D_t^{M-M'} D^n D_{t}^m H\right|  \, .
\label{exp.prG2.c}
\end{align}
To bound the quantity above, we first claim that for multi-indices $\alpha,\beta \in \mathbb{N}^k$ with $k\geq 2$, $|\alpha|\leq N_*$, and $|\beta|\leq M_*$, 
\begin{equation}
    \left| \prod_{i=1}^k D^{\alpha_i} D_t^{\beta_i} H \right| (x) \les \left(\pr(H)^{\sfrac 32}(x) + \const_{G,\infty}\Gamma^{-\NcutSmall} \right) (\lambda \Gamma)^{|\alpha|} \MM{|\beta|,M_t,\nu\Gamma,\nu'\Gamma} \, . \label{eq:sample:mess:1.c}
\end{equation}
To prove this claim, let $\Omega(x)\subseteq \supp H$ be a closed set containing $x$. Then applying Lemma~\ref{lem:cooper:2} with $p=\infty$, $N_t=M_t$, $N_*=\NcutLarge$, $M_*=\NcutSmall$, $\Omega=\Omega(x)$, $\const_v=\nu\lambda^{-1}$, $\lambda_v=\tilde\lambda_v=\lambda$, $\mu_v=\nu$, $\tilde \mu_v=\nu'$, $f=H$, $\const_f = \sup_{\Omega(x)}\left(\pr^{\sfrac 32}(H)+\const_{G,\infty}\Gamma^{-\NcutSmall}\right)$, $\lambda_f = \tilde\lambda_f = \lambda\Gamma$, $\mu_f =  \nu\Gamma$, and $\tilde \mu_f = \nu'\Gamma$, we have that \eqref{eq:cooper:2:v} is satisfied from \eqref{eq:DDv:sample:c}, and \eqref{eq:cooper:2:f} is satisfied by \eqref{est.G.by.pr.c}.  Then \eqref{eq:cooper:2:f:2} gives that 
\begin{equation}
    \left| \prod_{i=1}^k D^{\alpha_i} D_t^{\beta_i} H \right| (x) \les \left(\sup_{\Omega(x)}\pr(H)^{\sfrac 32} + \const_{G,\infty}\Gamma^{-\NcutSmall} \right) (\lambda \Gamma)^{|\alpha|} \MM{|\beta|,M_t,\nu\Gamma,\nu'\Gamma} \, . \label{eq:damn:commutators}
\end{equation}
Since $\Omega(x)$ is arbitrary and $\pr(H)$ is continuous, we have proven \eqref{eq:sample:mess:1.c}. 
Plugging this bound into \eqref{exp.prG2.c}, we find that 
\begin{align}
\frac{\left|D^N D_t^M \left((\pr(H)+ (\const_{G,\infty}\Gamma^{-\NcutSmall})^{\sfrac 23})^3 \right) \right|}{\left|\pr(H)+ (\const_{G,\infty}\Gamma^{-\NcutSmall})^{\sfrac 23}\right|^2} &\lesssim \frac{1}{\left|\pr(H)+ (\const_{G,\infty}\Gamma^{-\NcutSmall})^{\sfrac 23}\right|^2} \left( \pr^{\sfrac 32}(H) + \const_{G,\infty}\Gamma^{-\NcutSmall}\right)^2 \notag\\
&\qquad \qquad \times (\lambda\Gamma)^N \MM{M,M_t-\NcutSmall,\nu\Gamma,\nu'\Gamma} \, , \notag 
\end{align}
which implies the desired bound in \eqref{ind.G.c} concluding its proof.

Arguing in a similar way (in fact the proof is simpler since only spatial derivatives are required), we also have that for each integer $0\leq  N\leq N_*-\NcutLarge$,
\begin{subequations}
\begin{align}\label{ind.rho.c}
    \left|D^N (\pr(\rho)\circ \Phi)\right| \lec \left(\pr(\rho)\circ \Phi + (\const_{\rho,\infty}\Gamma^{-\NcutSmall})^{\sfrac 23} \right) \ (\Lambda\Ga)^N  \, , \\
    \label{ind.rho.c.I.need}
    \left|D^N \pr(\rho)\right| \lec \left(\pr(\rho) + (\const_{\rho,\infty}\Gamma^{-\NcutSmall})^{\sfrac 23} \right) \ (\Lambda\Ga)^N  \, .
\end{align}
\end{subequations}
Combining \eqref{ind.G.c}, \eqref{ind.rho.c}, and the choice of $\delta_{\rm tiny}$ from \eqref{eq:sample:1:Ncut:2:c}, we obtain the desired estimate \eqref{est.pr.phi}.  

\noindent\textbf{Proof of \eqref{sample:item:1:c}: }
Observe that by the construction of $\pr(H)$, \eqref{est.G.c}, and a computation similar to that used to produce \eqref{eq:something:I:need}, we have $\norm{\pr(H)+ (\const_{G,\infty}\Gamma^{-\NcutSmall})^{\sfrac 23}}_{\sfrac 32} \lec \const_{G,1}^{\sfrac 23}$, and so $\norm{\pr(H)}_{\sfrac 32} \lec \const_{G,1}^{\sfrac 23}$, with analogous bounds holding for $\rho$. It follows from \eqref{ind.G.c} and \eqref{eq:sample:1:Ncut:2:c} that
\begin{align}\label{ind.G.p.c}
    \norm{D^N D_t^M \pr(H)}_{\sfrac 32} \lec \const_{G,1}^{\sfrac 23} (\la\Gamma)^N \MM{M,M_t-\NcutSmall,\nu\Gamma,\nu'\Gamma}
\end{align}
for $N\leq N_*-\NcutLarge$ and $M\leq M_*-\NcutSmall$. If the left-hand side is measured instead in $L^\infty$, we may appeal to \eqref{eq:something:I:need} to deduce that \eqref{ind.G.p.c} holds with $\const_{G,\infty}$ in place of $\const_{G,1}$. Arguing similarly for $\pr(\rho)$ but appealing to \eqref{ind.rho.c} and \eqref{eq:something:else:I:need}, we have that \eqref{eq:sample:1:conc:2:c}--\eqref{eq:sample:1:conc:2:c:infty} are verified. Also, by the construction of $\pr(\rho)$, its periodicity easily follows from \eqref{sample:item:assump:2:c}. Next, we can immediately deduce from the definition of $\sigma_S^-$ and for $p=\sfrac 32,\infty$ the easier bound
\begin{align*}
\norm{\si^-_S}_{p} \lec r_\phi^{\sfrac 23}
\norm{\pr(H)}_{p}
\norm{\pr(\rho)}_1 \, ,
\end{align*}
which matches the desired bounds in \eqref{eq:sample:1:conc:1:c}--\eqref{eq:sample:1:conc:1:c:infty} for $\sigma_\phi^-$ after using the aforementioned bounds for $\pr(H),\pr(\rho)$ and recalling the definition of $\delta_{\phi,\cdot}$ from \eqref{est.phi.c}.
In the case of $\si_\phi^+$ and $p=\sfrac 32$, we additionally apply Lemma \ref{l:slow_fast} by setting
\begin{align}
    &N_*=N_*-\NcutLarge, \quad M_*= M_*-\NcutSmall, \quad f=\pr(H), \quad \Phi=\Phi \, , \notag \\
    &\lambda=\la \Gamma, \quad \tau^{-1}=\nu \Gamma, \quad \Tau^{-1} = \nu'\Gamma,  \notag \\
    &\const_{f}=\const_{G,1}^{\sfrac 23}, \quad v = v, \quad \varrho = \pr(\rho), \quad \mu = \mu, \notag \\
    &\Upsilon = \Lambda = \Lambda\Gamma,   \quad
    \const_{\varrho} = \const_{\rho,1}^{\sfrac 23}, \quad  N_t = M_t-\NcutSmall \, . \notag 
\end{align}
Then \eqref{eq:slow_fast_0} is verified from \eqref{ind.G.p.c}, \eqref{eq:slow_fast_1}--\eqref{eq:slow_fast_2} follow from \eqref{eq:DDpsi:sample:c}, \eqref{eq:slow_fast_4} follows from \eqref{ind.rho.c.I.need} and the periodicity of $\pr(\rho)$, \eqref{eq:slow_fast_3} follows from \eqref{eq:sample:1:decoup:c}, and \eqref{eq:slow_fast_3_a} follows from \eqref{eq:sample:1:Ncut:3:c}. We then obtain from \eqref{eq:slow_fast_5} that
\begin{align*}
\norm{\si^+_S}_{\sfrac32} \lec r_\phi^{\sfrac 23} \const_{G,1}^{\sfrac 23} \const_{\rho,1}^{\sfrac 23} = \de_{\phi,1} \, .
\end{align*}
Finally, the estimate for $\norm{\si^+_S}_{\infty}$ is trivial, so that \eqref{eq:sample:1:conc:1:c}--\eqref{eq:sample:1:conc:1:c:infty} holds for $\sigma_\phi^+$, and \eqref{sample:item:1:c} is totally verified.

\noindent\textbf{Proof of \eqref{sample:item:4:c}:} We first prove \eqref{est.pr.H.c} by induction; namely,
for each integer $K=N+M\geq 0$, $N\leq N_*-\NcutLarge$, $M\leq M_*-\NcutSmall$,
\begin{align}\label{ind.G.pt.c}
    |D^ND_t^M \pr(H)|\lec r_G^{-\sfrac 23}\pi  (\la\Ga)^N \MM{M, M_t-\NcutSmall, \nu\Ga, \nu\Ga'}  \, .
\end{align}
The proof uses an argument quite similar to the proof of \eqref{ind.G.c}.  The base case follows from writing that
\begin{align*}
    \pr(H) &\les \pi \, r_G^{-\sfrac 23} \\
    \iff \pr(H) + \left(\const_{G,\infty}\Gamma^{-\NcutSmall}\right)^{\sfrac 23} &\les \pi \, r_G^{-\sfrac 23} + \left(\const_{G,\infty}\Gamma^{-\NcutSmall}\right)^{\sfrac 23} \\
    \impliedby \left( \pr(H) + \left(\const_{G,\infty}\Gamma^{-\NcutSmall}\right)^{\sfrac 23} \right)^3 &\les \pi^3 r_G^{-2} + \left(\const_{G,\infty}\Gamma^{-\NcutSmall}\right)^{2} \, ,
\end{align*}
which can be seen to hold from the definition of $\pr(H)$ and \eqref{est.G.c.pt}. For the inductive step, we argue starting from \eqref{eq:liebniz:current}, although with slightly different steps to follow. Using the inductive assumption from \eqref{ind.G.pt.c} to control the term from the trilinear product in the second term with the \emph{highest} number of derivatives,\footnote{In fact any term which has been differentiated at all will suffice, so that we may replace $\pr(H)+\const_{G,1}^{\sfrac 23}$ with simply $\pr(H)$.} the bound \eqref{ind.G.c} to control the other two terms from the trilinear product, and \eqref{eq:sample:1:Ncut:2:c}, we have that the second term from \eqref{eq:liebniz:current} may be bounded by
\begin{align}
& \frac{1}{\left|\pr(H)+\left(\const_{G,\infty}\Gamma^{-\NcutSmall}\right)^{\sfrac 23} \right|^2} r_G^{-\sfrac 23} \pi \left( \pr(H) +  \left(\const_{G,\infty}\Gamma^{-\NcutSmall}\right)^{\sfrac 23} \right)^2 (\lambda \Gamma)^N \MM{M,M_t-\NcutSmall,\Gamma \nu, \Gamma \nu'} \notag\\
&\qquad \les r_G^{-\sfrac 23} \pi (\lambda \Gamma)^N \MM{M,M_t-\NcutSmall,\Gamma \nu, \Gamma \nu'} \, .
\label{est.prG.sec.c.redux}
\end{align}
Thus it remains to control the first term from \eqref{eq:liebniz:current}.  Towards this end, we claim that for multi-indices $\alpha,\beta\in\mathbb{N}^k$ with $k\geq 2$, $|\alpha|\leq N_*$, and $|\beta|\leq M_*$, 
\begin{equation}
    \left| \prod_{i=1}^k D^{\alpha_i} D_t^{\beta_i} H \right| (x) \les \pi^{\sfrac 32}(x) r_G^{-1} (\lambda \Gamma)^{|\alpha|} \MM{|\beta|,M_t,\nu\Gamma,\nu'\Gamma} \, . \label{eq:sample:mess:1.c.redux}
\end{equation}
As in the proof of \eqref{eq:sample:mess:1.c}, we apply Lemma~\ref{lem:cooper:2} with precisely the same choices as led to the bound in \eqref{eq:damn:commutators}, save for the choice of $\const_f = \sup_{\Omega(x)}\pi^{\sfrac 32} r_G^{-1}$.  Then \eqref{eq:cooper:2:v} is satisfied from \eqref{eq:DDv:sample:c}, and \eqref{eq:cooper:2:f} is satisfied by \eqref{est.G.c.pt}. Then applying \eqref{eq:cooper:2:f:2}, shrinking $\Omega(x)$ to a point, and using the continuity of $\pi$ provides \eqref{eq:sample:mess:1.c.redux}. Then plugging this bound into \eqref{exp.prG2.c} and using \eqref{eq:sample:mess:1.c} and \eqref{eq:sample:1:Ncut:2:c}, we find that for $N\leq N_*-\NcutLarge$ and $M\leq M_*-\NcutSmall$,
\begin{align}
&\frac{\left|D^N D_t^M \left(\pr(H)+(\const_{G,\infty}\Ga^{-\NcutSmall})^{\sfrac 23} \right)^3 \right|}{\left|\pr(H)+\left(\const_{G,\infty}\Gamma^{-\NcutSmall}\right)^{\sfrac 23}\right|^2} \notag\\
&\qquad \les \frac{1}{\left|\pr(H)+\left(\const_{G,\infty}\Gamma^{-\NcutSmall}\right)^{\sfrac 23}\right|^2} \pi r_G^{-\sfrac 23} \left( \pr^{\sfrac 32}(H) + \const_{G,\infty}\Gamma^{-\NcutSmall} \right)^{\sfrac 43} (\lambda\Gamma)^{N} \MM{M,M_t-\NcutSmall,\nu\Gamma,\nu'\Gamma} \notag\\
&\qquad \les \pi r_G^{-\sfrac 23} \frac{\pr^2(H) + \left({\const_{G,\infty}\Gamma^{-\NcutSmall}}\right)^{\sfrac 43}}{\left|\pr(H)+\left(\const_{G,\infty}\Gamma^{-\NcutSmall}\right)^{\sfrac 23}\right|^2} (\lambda\Gamma)^{N} \MM{M,M_t-\NcutSmall,\nu\Gamma,\nu'\Gamma} \notag\\
&\qquad \les \pi r_G^{-\sfrac 23}  (\lambda\Gamma)^{N} \MM{M,M_t-\NcutSmall,\nu\Gamma,\nu'\Gamma} \, , \notag 
\end{align}
which combined with \eqref{est.prG.sec.c.redux} concludes the proof of \eqref{est.pr.H.c}.  To prove \eqref{est.pr.phi-}, we use \eqref{est.pr.H.c} and the definition of $\sigma_\phi^-$.

\noindent\textbf{Proof of \eqref{sample:item:6:c}:} By the definition of $\pr(H)$ and $\pr(\rho)$, it is easy to see that $\supp(\pr(H))\subseteq \supp(H)$ and $\supp(\pr(\rho))\subseteq \supp(\rho)$, and so \eqref{sample:item:6:c} is verified.
\end{proof}

\begin{proposition*}[\bf Pressure increment and upgrade error from velocity increment potential]\label{lem:pr.vel} We begin with \emph{assumptions} which allow for the construction of a pressure increment and an upgrade current error.  Then we delineate a number of properties satisfied by the \emph{pressure increment}, before applying the material derivative and inverse divergence to produce a \emph{current error} satisfying additional properties.

\noindent\textbf{Part 1: Assumptions}

\noindent Let $v$ be an incompressible vector field on $\R\times\T^3$. Denote its material derivative by $D_{t}=\partial_t + v\cdot\nabla$. 
We use large positive integers $N_{**}$, $\dpot$, $K_\circ$, $N_*\geq M_*\gg M_t$, and $1\leq M_\circ\leq N_\circ\leq \sfrac{1}{2}(M_*-\NcutSmall-1-N_{**})$ and specify additional constraints that they must satisfy below. Suppose a velocity increment potential $\hat\upsilon= G (\rho\circ\Phi)$ and a non-negative continuous function $\pi$ are given such that the following hold.
\begin{enumerate}[(i)]
    \item There exist constants $\const_{G,p}$ and $\const_{\rho,p}$ for $p =3, \infty$, frequency parameters $\la, \La, \nu, \nu'$, and intermittency parameters $r_G, r_{\hat \upsilon}\leq 1$ such that
    \begin{subequations}
    \begin{align}
    \norm{D^N D_{t}^M G}_p 
    &\les \const_{G,p} \la ^N\MM{M,M_{t},\nu,\nu'} \label{est.G.sample3} \\ 
    \left|D^N D_{t}^M G\right| 
    &\les \pi^\frac12r_{G}^{-\frac13} \la ^N\MM{M,M_{t},\nu,\nu'}\label{est.G.pt.sample3}\\
    \norm{ D^N  \rho}_{p} 
    &\les \const_{\rho,p} \Lambda^N  \label{est.rho.sample3}\\
    \norm{\hat\upsilon}_p &\les
    \const_{G,p}\const_{\rho,p} 
    =:\de_{{\hat\upsilon},p}^\frac12r_{{\hat\upsilon}}^{-\frac13} \label{eq:desert:def:svi}
    \end{align}
    \end{subequations}
   
for all $N \leq N_*$, $M\leq M_*$. 

\item There exist frequency parameters $\mu$ and $\lambda'$, a parameter $\Gamma=\Lambda^\alpha$ for $0<\alpha\ll 1$ for measuring small losses in derivative costs, and a positive integer $\Ndec$ such that $\rho$ is $(\sfrac{\T}{\mu})^3$-periodic and $\la,\la'\ll \mu\leq \Lambda$, whereby we mean that
\begin{subequations}
\begin{align}
\label{eq:sample:5:decoup:vel}
\max(\lambda,\lambda')\Ga \mu^{-1} \leq 1 \, , \qquad     (\Lambda \Gamma)^{4}  &\leq  \left( \frac{\mu}{4 \pi \sqrt{3} \max(\lambda',\lambda)\Gamma}\right)^{\Ndec} \, .
\end{align}
\end{subequations}
\item 
Let $\Phi$ be a volume preserving diffeomorphism of $\T^3$ such that $D_{t} \Phi =0$ and $\Phi$ is the identity at a time slice which intersects the support of $G$, and
\begin{subequations}
\begin{align}
\norm{D^{N+1}   \Phi}_{L^\infty\left(\supp G\right)} 
+ \norm{D^{N+1}   \Phi^{-1}}_{L^\infty\left(\supp G\right)} 
&\les \lambda'^{N}
\label{eq:DDpsi:sample3}\\
\norm{D^ND_{t}^M D  v}_{L^\infty\left(\supp G\right)}
&\les  \nu \la'^{N}\MM{M,M_{t},\nu,\nu'}
 \label{eq:DDv:sample3}
\end{align}
for all $N \leq N_*$, $M\leq M_*$.  Furthermore, assume that we have the lossy estimate
\begin{align}
    \norm{D^N \partial_t^M v}_{L^\infty}\les  \const_v \lambda'^N (\nu')^M  \, , \qquad \const_{v}\la'\lec \nu'
\label{eq:inverse:div:v:global.sample3}
\end{align}
for all $M\leq M_\circ$ and $N+M \leq N_\circ+M_\circ$.
\end{subequations}
\item There exist positive integers $\NcutLarge$, $\NcutSmall$ and a small parameter $\delta_{\rm tiny} \leq 1$ such that 
\begin{subequations}
\begin{align}
\NcutSmall&\leq \NcutLarge \, , \label{par.con.sample3.Ncut} \\
    (\const_{G,\infty}^2+1)(\const_{\rho,\infty}^2+1)\Ga^{-2\NcutSmall  }&\leq \delta_{\rm tiny} \, ,  \, \const_{G,3}^2 \, , \,  \const_{\rho,3}^2
    \,, \label{par.con.sample3}\\
2\Ndec + 4 \leq N_* - \NcutLarge - N_{**} \, , &\qquad \NcutSmall \leq M_t-1  \, . \label{par.con.sample3.dec}
\end{align}
\end{subequations}
\item Let an increasing sequence of frequencies $\{\mu_0, \cdots, \mu_{\bar m}\}$, $\mu < \mu_0<\cdots<\mu_{\bar m-1}<\Lambda\Ga <\mu_{\bar m}$ be given satisfying
\begin{align}
\max(\lambda,\lambda')\Ga \mu_{m-1}^{-2} \mu_m \leq 1
 \label{eq:sample:prop:parameters:0:vel}
\end{align}
for all $1\leq m <\bar m$.
\item Assume that $\dpot$ and $N_{**}$ are sufficiently large so that
\begin{subequations}
\begin{align}
        \nu\Ga \const_{G,p}^2 \const_{\rho,p}^2 (\max(\lambda,\lambda')\Ga)^{\lfloor \sfrac \dpot 2 \rfloor} \mu^{-\lfloor \sfrac \dpot 2 \rfloor}  &(\Lambda\Ga)^{5+K_\circ} \left(1 + \frac{\max\{ \nu'\Ga, \const_v \Lambda\Ga \}}{\nu\Ga
}\right)^{M_\circ} \leq 1\, ,
\label{eq:sample:riots:4:vel:1} 
\\
        \nu\Ga \const_{G,p}^2 \const_{\rho,p}^2 (\max(\lambda,\lambda')\Ga)^{\lfloor \sfrac \dpot 2 \rfloor} (\mu_m \mu_{m-1}^{-2})^{\lfloor \sfrac \dpot 2 \rfloor}  &(\Lambda\Ga)^{5+K_\circ} \left(1 + \frac{\max\{ \nu'\Ga, \const_v \Lambda\Ga \}}{\nu\Ga
}\right)^{M_\circ} \leq 1\, ,
\label{eq:sample:riots:4:vel:2}\\
\nu\Ga \const_{G,\infty}^2 \const_{\rho,3}^2 ((\La\Ga)\mu_{\bar m}^{-1})^{N_{**}}   &(\Lambda\Ga)^{5+K_\circ} \left(1 + \frac{\max\{ \nu'\Ga, \const_v \Lambda\Ga \}}{\nu\Ga
}\right)^{M_\circ} \leq 1\, ,
\label{eq:sample:riots:4:vel:3}
\end{align}
\end{subequations}
for $1\leq m\leq \bar m$.
\end{enumerate}
\smallskip

\noindent\textbf{Part 2: Pressure increment}

\noindent There exists a pressure increment $\si_{\hat\upsilon}=\si_{\hat\upsilon}^+ - \si_{\hat\upsilon}^-$ associated to the velocity increment potential $\hat\upsilon$ which is defined by
\begin{subequations}\label{eq:desert:some:stuff}
\begin{align}
    \si_{\hat\upsilon}
    &:= r_{\hat\upsilon}^{2} \pr(G)
\left(\pr(\rho)\circ\Phi-\langle 
    \pr(\rho)\rangle\right) 
    =:
    \si_{\hat\upsilon}^+ -\si_{\hat\upsilon}^-\, , \\
\pr(G) &:= \sum_{N=0}^{\NcutLarge}\sum_{M=0}^{\NcutSmall} (\la \Ga)^{-2N}  (\nu\Ga) ^{-2M} |D^N D_{t}^M G|^2 \, , \\
\pr(\rho) &:=  \sum_{N=0}^{\NcutLarge} (\Lambda\Gamma )^{-2N}   |D^N \rho|^2  \, , 
\end{align}
may be decomposed as
\begin{align}\label{eq:desert:decomp:sv}
\sigma_{\hat \upsilon} = \sigma_{\hat \upsilon}^* + \sum_{m=0}^{\bm} \sigma_{\hat \upsilon}^m \, ,    
\end{align}
\end{subequations}
and satisfies the properties listed below.
\begin{enumerate}[(i)] 
\item\label{sample3:item:2} $(\si_{\hat\upsilon}^+)^{\sfrac12}$ dominates derivatives of $\hat\upsilon$ with suitable weights, so that
\begin{align}\label{est.w.by.pr}
    \left|D^N D_{t}^M {\hat\upsilon}\right|
    \lec (\si_{\hat\upsilon}^+  + \de_{\rm tiny})^{\sfrac12}r_{\hat\upsilon}^{-1} (\Lambda\Ga)^N\MM{M,M_{t},\nu\Ga,\nu'\Ga}.
\end{align}
for all $N \leq N_*$, $M\leq M_*$.
\item\label{sample3:item:3} $\si_{\hat\upsilon}^+$ dominates derivatives of itself with suitable weights, so that
\begin{align}\label{est.pr.w}
    |D^N D_{t}^M \si_{\hat\upsilon}^+|
    \lec (\si_{\hat\upsilon}^{+}+ \de_{\rm tiny}) (\Lambda\Ga)^N\MM{M,M_{t}-\NcutSmall,\nu\Ga,\nu'\Ga}
\end{align}
for all $N \leq N_*-\NcutLarge$, $M\leq M_*-\NcutSmall$.
\item\label{sample3:item:1} Let $(p,p')=(3,\sfrac32)$ or $(\infty, \infty)$. Then $\si_{\hat\upsilon}^+$ and $\si_{\hat\upsilon}^-$ satisfy
\begin{align*}
        \norm{\si_{\hat\upsilon}^+}_{p'}\lec \de_{{\hat\upsilon},p} r_{\hat\upsilon}^{\sfrac43}, \quad
        &\norm{\si_{\hat\upsilon}^-}_{p'}\lec \de_{{\hat\upsilon},p}r_{\hat\upsilon}^{\sfrac43} \, .
\end{align*}
We note also that $\pr(\rho)$ is $\left(\sfrac{\T}{\mu}\right)^3$-periodic. Furthermore, $\pr(G)$ and $\pr(\rho)$ have the same size as $G$ and $\rho$, so that for $N\leq N_* -\NcutLarge$ and $M\leq M_*-\NcutSmall$,
\begin{align}\label{eq:sample:1:conc:2:v}
    \norm{D^N D_t^M \pr(G)}_{p'}\lec \const_{G,p}^2 (\lambda\Gamma)^N \MM{M,M_t-\NcutSmall,\nu\Gamma,\nu'\Gamma} \, , \qquad &\norm{D^N \pr(\rho)}_{p'} \lec \const_{\rho,p}^2 (\Lambda\Ga)^N \, . 
\end{align}
\item\label{sample3:item:4} $\pi$ dominates $\si_{\hat\upsilon}^-$ and $\pr(G)$ and its derivatives with suitable weights, so that
\begin{subequations}
\begin{align} 
       \left|D^N D_{t}^M \pr(G)\right|
    &\lec  \pi  r_G^{-\sfrac{2}{3}}(\la\Gamma)^N\MM{M,M_{t}-\NcutSmall,\nu\Gamma,\nu'\Gamma}\, , \label{est.pr.w-.0} \\
    | D^N D_{t}^M \si_{\hat\upsilon}^-|
    &\lec  \pi r_G^{-\sfrac{2}{3}} \norm{\pr(\rho)}_1 r_{\hat\upsilon}^2 (\la\Gamma)^N\MM{M,M_{t}-\NcutSmall,\nu\Gamma,\nu'\Gamma} \label{est.pr.w-}
\end{align}
\end{subequations}
for all $N \leq N_*-\NcutLarge$, $M\leq M_*-\NcutSmall$. 
\item We have the support properties
\begin{align}\label{supp:pr:vecl}
    \supp(\si_{\hat\upsilon}^+)\subset\supp(\hat\upsilon) \, , \quad
    \supp(\si_{\hat\upsilon}^-) \subseteq \supp(G) \, .
\end{align}
\end{enumerate}
\smallskip

\noindent\textbf{Part 3: Current error}

\noindent There exists an upgrade current error $\phi_{\hat\upsilon}$ which satisfies the following properties.
\begin{enumerate}[(i)]
\item \label{sample3:item:decomp}  We have the decomposition and equalities
\begin{subequations}
\begin{align}
\phi_{\hat\upsilon} &= \underbrace{\phi_{\hat\upsilon}^*}_{\textnormal{nonlocal}} + \underbrace{\sum_{m=0}^{\bar m} \phi_{\hat\upsilon}^m}_{\textnormal{local}}\\ 
\div \left( \phi_{\hat\upsilon}^m(t,x) + \divR(D_t \sigma^m_{\hat \upsilon})(t,x) \right) &= D_t \sigma_{\hat \upsilon}^m(t,x) - \int_{\T^3} D_t \sigma_{\hat \upsilon}^m(t,x') \, dx' \, , \\
\div \left( \phi_{\hat\upsilon}^*(t,x) - \sum_{m=0}^{\bm} \divR(D_t \sigma^m_{\hat \upsilon})(t,x) \right) &= D_t \sigma_{\hat \upsilon}^*(t,x) - \int_{\T^3} D_t \sigma_{\hat \upsilon}^*(t,x') \, dx' \, .
\end{align}
\end{subequations}
\item \label{sample3:item:7} Let $(p,p')=(3,\sfrac32)$ or $(\infty, \infty)$. The current error $\phi_{\hat\upsilon}^{m}$ satisfies
\begin{subequations}
\begin{align}
\norm{D^ND_t^M\phi_{\hat\upsilon}^{0}}_{p'}
&\lec \nu\Ga^2 \const_{G,p}^2\const_{\rho,3}^2 {\left(\frac{\mu_0}{\mu}\right)^{\frac 43 - \frac{2}{p'}}} r_{\hat\upsilon}^2\mu^{-1} \mu_0^N \MM{M,M_{t}-\NcutSmall-1,\nu\Gamma,\nu'\Gamma}\, ,
\label{pr:current:vel:loc:p:0}\\
\left|D^ND_t^M\phi_{\hat\upsilon}^{0}\right|
&\lec \nu \Ga^2 \pi r_G^{-\sfrac23}  \const_{\rho,3}^2 \left(\frac{\mu_0}{\mu}\right)^{\sfrac43} r_{\hat\upsilon}^{2} \mu^{-1} 
\mu_0^N \MM{M,M_{t}-\NcutSmall-1,\nu\Gamma,\nu'\Gamma}
\label{pr:current:vel:loc:pt:0}
\, , \\
\norm{D^ND_t^M\phi_{\hat\upsilon}^{m}}_{p'}
&\lec \nu\Ga^2 \const_{G,p}^2\const_{\rho,3}^2 {\left(\frac{\min(\mu_m,\Lambda\Gamma)}{\mu}\right)^{\frac 43 - \frac{2}{p'}}} r_{\hat\upsilon}^2(\mu_{m-1}^{-2}\mu_m) \notag\\
&\qquad \qquad \times \min(\mu_m, \La\Ga)^N \MM{M,M_{t}-\NcutSmall-1,\nu\Gamma,\nu'\Gamma}\, ,
\label{pr:current:vel:loc:p}
\\
\left|D^ND_t^M\phi_{\hat\upsilon}^{m}\right|
&\lec \nu \Ga^2 \pi r_G^{-\sfrac 23} \const_{\rho,3}^2 \left(\frac{\min(\mu_m,\La\Ga)}{\mu}\right)^{\sfrac43} r_{\hat\upsilon}^{2}  \mu_{m-1}^{-2}\mu_m \notag\\
&\qquad \qquad \times
(\min(\mu_m, \La\Ga))^N \MM{M,M_{t}-\NcutSmall-1,\nu\Gamma,\nu'\Gamma}
\label{pr:current:vel:loc:pt}
\, , 
\end{align}
\end{subequations}
for any $1\leq  m \leq \bar m$, $N\leq N_*-\halfd - \NcutLarge-N_{**}$, and $M\leq M_*-\NcutSmall-1-N_{**}$.  Furthermore, we have that $\phi_{\hat\upsilon}^{ *}$ satisfies
\begin{align}
    \norm{D^ND_t^M\phi_{\hat\upsilon}^{*}}_{\infty}
    &\lec {\mu_0^{-K_\circ}} (\La\Ga)^N (\nu\Ga)^M \label{pr:current:vel:nonloc:infty1}
\end{align}
for all $N\leq N_\circ$ and $M\leq M_\circ$.
\item\label{sample3:item:6}
We have the support properties\footnote{For any $\Omega\in\T^3$, we use $\Omega\circ \Phiik$ to refer to the space-time set $\Phiik^{-1}(t,\cdot)\Omega$ whose characteristic function is annihilated by $D_t$.}
\begin{align}
    \supp(\phi_{\hat\upsilon}^{0})
    \subseteq \supp(G) \, , \quad
    \supp(\phi_{\hat\upsilon}^{m}) \subseteq  \supp G \cap B\left( \supp \rho, 2\mu_{m-1}^{-1} \right) \circ \Phi \label{supp:pr:vel:cur}
\end{align}
for all $0< m \leq \bar m $. 
\item For all $M\leq M_* -\NcutSmall - 1$, we have that the mean $\langle D_t \si_{\hat \upsilon} \rangle$ satisfies
\begin{equation}\label{est:mean.Dtsivel}
    \left|\frac{d^M}{dt^M}\langle D_t \si_{\hat \upsilon} \rangle\right|\lec (\La\Ga)^{-K_\circ}
    \MM{M,M_t - \NcutSmall, -1, \nu\Ga, \nu'\Ga}\, . 
\end{equation}
\end{enumerate}
\end{proposition*}

\begin{proof}\, 
\noindent\texttt{Step 1:} Constructing $\si_{\hat\upsilon}$ and verifying the properties in Part 2.

\noindent For the moment we ignore the decomposition in \eqref{eq:desert:decomp:sv} and handle the rest of the conclusions in Part 2. Towards a proof of \eqref{sample3:item:2}, we first have that $\pr(G)\lec \const_{G,\infty}^2$ and $\pr(\rho)\lec \const_{\rho,\infty}^2$.  The proof of these is similar to \eqref{est.prH.inf.0} and \eqref{est.prrho.inf.0}, and we omit the details. Also, using a method of proof similar to that used to obtain \eqref{est.G.by.pr} and \eqref{est.rho.by.pr}, we can show that
\begin{subequations}
\begin{align}
|D^{N_1}D_t^M G|
&\lec (\pr(G) + \const_{G,\infty}^2\Ga^{-2\NcutSmall})^{\sfrac12} (\lambda\Ga)^{N_1} \MM{M, M_t, \nu\Ga, \nu'\Ga} \label{eq:home:G:est} \\
\la^{N_2-n_2}|D^{n_2} \rho|
&\lec (\pr(\rho) + \const_{\rho,\infty}^2\Ga^{-2\NcutSmall})^{\sfrac12} (\Lambda\Ga)^{N_2} \label{eq:home:rho:est}
\end{align}
\end{subequations}
for any integers $0\leq N_1, N_2\leq N_*$, $0\leq n_2\leq N_2$ and $M\leq M_*$.
Then, \eqref{sample3:item:2} follows as in the proof of \eqref{est.S.by.pr}. 

Next, to prove \eqref{sample3:item:3}, we again claim that for $N\leq N_*-\NcutLarge$ and $M\leq M_*-\NcutSmall$,
\begin{subequations}
\begin{align}
     |D^N D_t^M \pr(G)|
     &\lec \left( \pr(G) + \const_{G,\infty}^2\Ga^{-2\NcutSmall} \right) (\la\Ga)^N \MM{M, M_t-\NcutSmall, \nu\Ga, \nu'\Ga} \label{ind.G.w}\\
     \left|D^N \pr(\rho)\right| 
       &\lec \left(\pr(\rho) + \const_{\rho,\infty}^2\Ga^{-2\NcutSmall} \right) \ (\Lambda\Ga)^N  \, \label{ind.rho.w}\\
       \left|D^N (\pr(\rho)\circ \Phi)\right| 
       &\lec \left(\pr(\rho)\circ \Phi + \const_{\rho,\infty}^2\Ga^{-2\NcutSmall} \right) \ (\Lambda\Ga)^N  \, . \label{ind.rho.w1}
\end{align}
\end{subequations}
The proof of the claims is similar to, and in fact easier, than the proofs of the analogous estimates in \eqref{ind.G} and \eqref{ind.rho}. Indeed, instead of \eqref{leib.prH}, we simply have from the Leibniz rule that
\begin{align*}
    \left|D^N D_{t}^M \pr(G) \right|
     &\leq   \sum_{n=0}^{\NcutLarge}\sum_{m=0}^{\NcutSmall} (\la \Ga)^{-2n}  (\nu\Ga) ^{-2m} \left|D^N D_t^M \left|D^n D_{t}^m G\right|^2\right|\nonumber \\
&=  \sum_{n=0}^{\NcutLarge}\sum_{m=0}^{\NcutSmall}
\sum_{\substack{0\leq N'\leq N\\0\leq M'\leq M}} 
(\la \Ga)^{-2n}  (\nu\Ga) ^{-2m}
\left|D^{N'} D_t^{M'} D^n D_{t}^m G\right| 
\left|D^{N-N'} D_t^{M-M'} D^n D_{t}^m G\right|  \, ,
\end{align*}
at which point we apply \eqref{eq:home:G:est}.  A similar argument produces the other two bounds listed above. Then \eqref{ind.G.w}--\eqref{ind.rho.w1} imply \eqref{sample3:item:3} as in the proof of Proposition \ref{lem:pr.st}. 

Regarding \eqref{sample3:item:1}, as before, the estimate for $G$ in \eqref{eq:sample:1:conc:2:v} follows from \eqref{est.G.sample3}, \eqref{ind.G.w}, and \eqref{par.con.sample3}.  The estimate for $\pr(\rho)$ follows similarly from \eqref{est.rho.sample3}, \eqref{ind.rho.w}, and \eqref{par.con.sample3}. Therefore, \eqref{eq:sample:1:conc:2:v} is verified, and as a consequence $\norm{\sigma_{\hat\upsilon}^-}_{p'} \lec \delta_{\hat \upsilon, p}r_{\hat\upsilon}^{\sfrac43}$ follows after using \eqref{eq:desert:def:svi}. The periodicity of $\pr(\rho)$ is immediate from the definition and the periodicity assumption on $\rho$. To obtain $\norm{\sigma_{\hat\upsilon}^+}_{3} \lec \delta_{\hat \upsilon, \sfrac32}r_{\hat\upsilon}^{\sfrac43}$, we use Lemma~\ref{l:slow_fast} as in the proof of \eqref{eq:sample:1:conc:1:c}, for example. The assumptions in the lemma can be verified using \eqref{eq:sample:1:conc:2:v}, \eqref{eq:DDpsi:sample3}, \eqref{eq:sample:5:decoup:vel}, and
\eqref{par.con.sample3.dec} and the recently observed periodicity. Therefore, the desired estimate for $\sigma_{\hat\upsilon}^+$ in $L^{\sfrac32}$ follows from \eqref{eq:slow_fast_5}. The $L^\infty$ estimate follows trivially from \eqref{eq:sample:1:conc:2:v}.  

Next, we consider \eqref{sample3:item:4}.
Similar to the proof of \eqref{ind.G.w}, one can obtain
\begin{align}\label{ind.G.pt.w}
    |D^ND_t^M \pr(G)|\lec \pi r_{G}^{-\sfrac 23}   (\la\Ga)^N \MM{M, M_t-\NcutSmall, \nu\Ga, \nu\Ga'}  \, 
\end{align}
for any integer $N\leq N_*- \NcutLarge$ and $M\leq M_*-\NcutSmall$. Then we have \eqref{est.pr.w-.0}, and hence \eqref{est.pr.w-} holds. Finally, \eqref{supp:pr:vecl} is immediate from the definitions in \eqref{eq:desert:some:stuff}, concluding the proof of all claims in Part 2 except \eqref{eq:desert:decomp:sv}.

\smallskip

\noindent\texttt{Step 2:} Constructing the current errors $\phi_{\hat\upsilon}^m$ and verifying the properties in Part 3.

\noindent We first define $\si_{\hat\upsilon}^m$ in order to verify \eqref{eq:desert:decomp:sv}. Using the synthetic Littlewood-Paley decomposition from \eqref{eq:decomp:showing} and Definition~\ref{def:synth:LP}, we write
\begin{equation}
     \mathbb{P}_{\neq 0}\pr(\rho) = \tilde{\mathbb{P}}_{\mu_0}\mathbb{P}_{\neq 0} (\pr(\rho)) + \left( \sum_{m=1}^{\bar m}  \tilde{\mathbb{P}}_{(\mu_{m-1},\mu_m]} (\pr(\rho)) \right) + \underbrace{\left( \Id -  \tilde{\mathbb{P}}_{\mu_{\bar m}} \right)}_{=:\mathbb{P}^*}(\pr(\rho)) \, . \label{eq:decomp:showing.vel} 
\end{equation}
For convenience, we use the abbreviations $\mathbb{P}_{0}$ for $\tilde{\mathbb{P}}_{\mu_0}\mathbb{P}_{\neq 0}$ and $\mathbb{P}_{m}$ for $\tilde{\mathbb{P}}_{(\mu_{m-1},\mu_m]}$ for $1\leq m\leq \bar m$.
Define $\si_{\hat\upsilon}^m$, $\si_{\hat \upsilon}^*$, $\phi_{\hat\upsilon}^m$, and $\phi_{\hat\upsilon}^*$ by
\begin{align*}
    \si_{\hat\upsilon}
    &= \si_{\hat\upsilon}^* + \sum_{m=0}^{\bar m} \si_{\hat\upsilon}^m
    := r_{\hat\upsilon}^2\pr(G)( \mathbb{P}^*\pr(\rho))\circ\Phi) + r_{\hat\upsilon}^2
    \sum_{m=0}^{\bar m} \pr(G)   (\tilde{\mathbb{P}}_{m} (\pr(\rho)) \circ\Phi) \, , \\
    \phi_{\hat\upsilon}^m&:=\divH(D_t \si_{\hat\upsilon}^m),\quad 
    \phi_{\hat\upsilon}^*:= (\divH+\divR)\si_{\hat\upsilon}^* + \sum_{m=0}^{\bar m}\divR(D_t \si_{\hat\upsilon}^m) \, .
\end{align*}
Assuming that everything above is well-defined, we have verified \eqref{sample3:item:decomp}. We aim to apply Proposition \ref{prop:intermittent:inverse:div} with Remarks \ref{rem:scalar:inverse:div} and \ref{rem:pointwise:inverse:div} in separate cases according to which projector is being applied above.  In order to apply the inverse divergence, we may however first treat the low-frequency assumptions from Part 1, which are the same in all cases (irrespective of which projector is being applied). We therefore set
\begin{align*}
    &\overline N_* = N_*  - \NcutLarge - N_{**} \, , \quad \overline M_*=M_*-\NcutSmall-1 - N_{**} \, , \quad
    \overline M_t = M_t-\NcutSmall-1
    \\
    &\overline G = D_t \pr(G), \quad \overline \const_{G,\sfrac32} = {\nu\Ga \const_{G,3}^2}, \quad \overline\const_{G,{\infty}} = {\nu\Ga \const_{G,\infty}^2},\quad \overline\mu=\mu \, , \quad \ov \la' = \la' \, , \\
    &\overline\Phi=\Phi,\quad \overline\la=\max (\la,\la')\Ga,\quad \overline\nu=\nu\Ga,\quad \overline\nu'=\nu'\Ga \, , \quad \overline\pi = \nu{\Ga} \pi r_G^{-\sfrac23}\, , \quad \ov v = v \, ,
\end{align*}
where we have used the convention set out in Remark~\ref{rem:bar:variables} to rewrite the symbols from Lemma~\ref{lem:pr.st} with bars above on the left-hand side of the equalities below, while the right-hand side are parameters given in the assumptions of this Lemma. Then we have that \eqref{eq:inv:div:NM} is verified from the assumption $N_*\geq M_*$ and \eqref{par.con.sample3.Ncut}, \eqref{eq:inverse:div:DN:G} follows from conclusion \eqref{eq:sample:1:conc:2:v}, and \eqref{eq:inv:div:extra:pointwise} follows from conclusion~\eqref{est.pr.w-.0}. Next, we see that \eqref{eq:DDpsi2}, \eqref{eq:DDpsi}, \eqref{eq:DDv}, and \eqref{eq:inverse:div:v:global} hold from \eqref{eq:DDpsi:sample3}--\eqref{eq:inverse:div:v:global.sample3}. At this point we split into cases based on which projector is applied and address parts 2-4 of Proposition~\ref{prop:intermittent:inverse:div} in order to conclude the proof of this Lemma.
\smallskip

\noindent\texttt{Step 2a: Lowest shell.} For the case $m=0$, we appeal to Lemma~\ref{lem:special:cases} with $q=\sfrac 32$, $\lambda=\Lambda\Gamma$, $\rho = \mathbb{P}_{\neq 0} \pr(\rho)$, and $\alpha$ such that $\lambda^\alpha$ in \eqref{eq:lowest:shell:inverse} is equal to $\Gamma$.  Specifically, to verify the assumptions in Part 2 of Proposition~\ref{prop:intermittent:inverse:div}, we set for $p'=\sfrac 32,\infty$
\begin{align*}
    &\ov{\varrho} = \mathbb{P}_0 \pr(\rho) \, , \quad \ov{\vartheta}\textnormal{ as defined in \eqref{eq:lowest:shell:inverse}} \, , \quad \ov{\const}_{*,p'} = \const_{\rho,3}^{2} \left( \frac{\mu_0}{\mu} \right)^{\frac 43 - \frac{2}{p'}} \, , \\
    &\ov \mu = \mu \, , \quad \ov \Upsilon = \ov \Upsilon' = \mu \, , \quad \ov \Lambda = \mu_0 \, , \quad \ov \dpot = \dpot \, .
\end{align*}
Then \eqref{eq:moll:1:as} is satisfied with $\const_{p,\sfrac 32}=\const_{\rho,3}^2$ and $\lambda=\Lambda\Gamma$ from standard Littlewood-Paley theory, \eqref{eq:sample:1:conc:2:v}, and the choices from Step 1 which led to that conclusion, and so from \eqref{eq:lowest:shell:inverse} we have that \eqref{eq:DN:Mikado:density} is satisfied. From \eqref{eq:sample:5:decoup:vel}, \eqref{par.con.sample3.dec}, and the choice of $\ov N_*$ above, we have that \eqref{eq:inverse:div:parameters:0}--\eqref{eq:inverse:div:parameters:1} are satisfied. Continuing onto the nonlocal assumptions from Proposition~\ref{prop:intermittent:inverse:div}, we have that \eqref{eq:inv:div:wut}--\eqref{eq:inverse:div:v:global:parameters} are satisfied from \eqref{eq:inverse:div:v:global.sample3} and the assumptions from Part 1 on $M_\circ$ and $N_\circ$.  We have that \eqref{eq:riots:4} is satisfied from \eqref{eq:sample:riots:4:vel:1}. We then appeal to the conclusions \eqref{eq:inverse:div}--\eqref{eq:inverse:div:error:1} and \eqref{eq:inverse:div:error:stress}--\eqref{eq:inverse:div:error:stress:bound} to conclude as follows. From \eqref{eq:inverse:div:stress:1}, we obtain \eqref{pr:current:vel:loc:p:0}. The pointwise bound in \eqref{pr:current:vel:loc:pt:0} holds due to \eqref{eq:inv:div:extra:conc}, \eqref{eq:inverse:div:sub:1}, and \eqref{eq:divH:formula}. Next, we obtain \eqref{pr:current:vel:nonloc:infty1} for the portion of $\phi_{\hat\upsilon}^*$ coming from this case $m=0$ from \eqref{eq:inverse:div:error:stress:bound}. Finally, we obtain \eqref{supp:pr:vel:cur} from \eqref{eq:inverse:div:linear}, concluding the proof of the desired conclusions for $m=0$
.
\smallskip

\noindent\texttt{Step 2b: Intermediate shells.} For the cases $1\leq m\leq \bar m$, we appeal to Lemma~\ref{lem:LP.est} with $q=\sfrac 32$ and $\rho = \mathbb{P}_{\neq 0} \pr(\rho)$. Specifically, to verify the assumptions in Part 2 of Proposition~\ref{prop:intermittent:inverse:div}, we set for $p'=\sfrac 32, \infty$
\begin{align*}
    &\overline \varrho = \mathbb{P}_{m} \pr(\rho), \quad 
    \overline\const_{*,\sfrac32} = \const_{\rho,3}^2, \quad
    \overline\const_{*,{\infty}} =\min( (\sfrac{\mu_m}{\mu})^{\sfrac43}\const_{\rho,3}^2,\const_{\rho,\infty}^2)\, , \quad \overline\Upsilon=\mu_{m-1} \, , \\
    &\overline\Upsilon'=\overline\Lambda=\min(\mu_m, \La\Ga) \, , \quad \vartheta \textnormal{ as defined in Lemma~\ref{lem:LP.est}} \, , \qquad \alpha \textnormal{ as in the previous substep} \, .
\end{align*}
Then \eqref{eq:gen:id:assump} is satisfied with $\const_{p,\sfrac 32}=\const_{\rho,3}^2$ as in the last substep, and so from \eqref{eq:LP:div:estimates} we have that \eqref{eq:DN:Mikado:density} is satisfied. From \eqref{eq:sample:5:decoup:vel}, \eqref{par.con.sample3.dec}, \eqref{eq:sample:prop:parameters:0:vel}, and the choice of $\ov N_*$ above, we have that \eqref{eq:inverse:div:parameters:0}--\eqref{eq:inverse:div:parameters:1} are satisfied. Continuing onto the nonlocal assumptions from Proposition~\ref{prop:intermittent:inverse:div}, we have that \eqref{eq:inv:div:wut}--\eqref{eq:inverse:div:v:global:parameters} are satisfied as in the last substep.  We have that \eqref{eq:riots:4} is satisfied from \eqref{eq:sample:riots:4:vel:2}. We then appeal to the conclusions \eqref{eq:inverse:div}--\eqref{eq:inverse:div:error:1} and \eqref{eq:inverse:div:error:stress}--\eqref{eq:inverse:div:error:stress:bound} to conclude as follows. From \eqref{eq:inverse:div:stress:1}, we obtain \eqref{pr:current:vel:loc:p}. The pointwise bound in \eqref{pr:current:vel:loc:pt} holds due to \eqref{eq:inv:div:extra:conc}, \eqref{eq:inverse:div:sub:1}, and \eqref{eq:divH:formula}. Next, we obtain \eqref{pr:current:vel:nonloc:infty1} for the portion of $\phi_{\hat\upsilon}^*$ coming from this case $1\leq m \leq \bm$ from \eqref{eq:inverse:div:error:stress:bound}. Finally, we obtain \eqref{supp:pr:vel:cur} from \eqref{eq:inverse:div:linear} and \eqref{eq:LP:div:support}, concluding the proof of the desired conclusions for $1\leq m \leq \bm$.
\smallskip

\noindent\texttt{Step 2c: Highest shell.} For the case $m=\bm$, we appeal to Lemma~\ref{lem:special:cases} with $q=\sfrac 32$, $\lambda=\Lambda\Gamma$, $\rho = \mathbb{P}_{\neq 0} \pr(\rho)$, and $\alpha$ such that $\lambda^\alpha$ in \eqref{eq:lowest:shell:inverse} is equal to $\Gamma$.  Specifically, to verify the assumptions in Part 2 of Proposition~\ref{prop:intermittent:inverse:div}, we set for $p'=\infty$
\begin{align*}
    &\ov{\varrho} = \mathbb{P}^* \mathbb{P}_0 \pr(\rho) \, , \quad \ov{\vartheta}\textnormal{ as defined in \eqref{eq:remainder:inverse}} \, , \quad \ov{\const}_{*,p'} = \const_{\rho,3}^{2} (\Lambda\Ga)^3 \left( \frac{\Lambda\Gamma}{\mu_{\bm}} \right)^{N_**} \, , \\
    &\ov \mu = \ov \Upsilon = \ov \Upsilon' = \mu \, , \quad \ov \Lambda = \Lambda\Gamma \, , \quad \ov \dpot = 0 \, .
\end{align*}
Then \eqref{eq:moll:1:as} is satisfied as in the previous substeps, and so from \eqref{eq:remainder:inverse} we have that \eqref{eq:DN:Mikado:density} is satisfied. We have that \eqref{eq:inverse:div:parameters:0}--\eqref{eq:inverse:div:parameters:1} are satisfied as in the first substep.  The nonlocal assumptions are satisfied as in the previous substeps, except that we now have \eqref{eq:riots:4} from \eqref{eq:sample:riots:4:vel:3}.  The only conclusion we require at this point is to produce a bound matching \eqref{pr:current:vel:nonloc:infty1}, which follows from \eqref{eq:inverse:div:error:stress:bound}.

\noindent\texttt{Step 3: Verification of \eqref{est:mean.Dtsivel}.}  Since the vector field $v$ is incompressible, $\frac{d^M}{dt^M} \langle D_t \si_{\hat \upsilon} \rangle = \langle D_t^{M+1} \si_{\hat \upsilon} \rangle$. Since $\pr(\rho)$ is periodic in $(\sfrac{\T}{\mu})^2$, we have that for $M+1 \leq M_* - \NcutSmall - 1$,
\begin{align*}
\bigg|\int_{\T^3}&D_t^{M+1}\pr\left(G\right) \left(\mathbb{P}_{\neq 0} \pr(\rho)\right)\circ \Phi \, dx \bigg|\\     &=\left|\int_{\T^3}D_t^{M+1}\pr\left(G\right)\circ \Phi^{-1} \Delta^{\floor{\frac{\dpot}4}}\Delta^{-\floor{\frac{\dpot}4}}\left(\mathbb{P}_{\neq 0} \pr(\rho)\right) \,  dx\right|\\
     &=\left|\int_{\T^3}\Delta^{\floor{\frac{\dpot}4}}\left(D_t^{M+1}\pr\left(G\right)\circ \Phi^{-1}\right) \Delta^{-\floor{\frac{\dpot}4}}\left(\mathbb{P}_{\neq 0} \pr(\rho)\right) dx\right|\\
     &\lec \norm{\Delta^{\floor{\frac{\dpot}4}}\left(D_t^{M+1}\pr\left(G\right)\circ \Phi^{-1}\right)}_{\sfrac32}
     \norm{\Delta^{-\floor{\frac{\dpot}4}}\left(\mathbb{P}_{\neq 0} \pr(\rho)\right)}_{1}\\
     &\lec \const_{G,\sfrac32} (\max(\la, \la')\Ga)^{\halfd}
     \mu^{-\halfd}  \const_{*,\sfrac32}\Upsilon^{-2}\Upsilon' 
     \MM{M+1,M_t - \NcutSmall, \nu\Ga, \nu'\Ga}
     \\
     &\leq (\La\Ga)^{-K_\circ} \MM{M,M_t - \NcutSmall-1, \nu\Ga, \nu'\Ga} \, . 
\end{align*}
Here, we have used Lemma~\ref{lem:decoup}, \eqref{eq:sample:1:conc:2:v}, \eqref{eq:DDpsi:sample3}, \eqref{eq:sample:riots:4:vel:2}, and standard Littlewood-Paley theory. 
\end{proof}

\begin{proposition*}[\bf Pressure increment and upgrade error for stress error]\label{lem.pr.invdiv2}
We begin with \emph{preliminary assumptions}, which include all of the assumptions and conclusions from the inverse divergence in Proposition~\ref{prop:intermittent:inverse:div} and the pointwise bounds in Remark~\ref{rem:pointwise:inverse:div}. We then include \emph{additional assumptions}, which allow for the application of Lemma~\ref{lem:pr.st} to the stress error and Proposition~\ref{prop:intermittent:inverse:div} to the material derivative of the output. We thus obtain a \emph{pressure increment} which satisfies a number of properties.  Finally, the material derivative of this pressure increment produces a \emph{current error} which itself satisfies a number of properties.
\smallskip

\noindent\textbf{Part 1: Preliminary assumptions}
\begin{enumerate}[(i)]
    \item\label{i:st:sample:1} There exists a vector field $G$, constants $\const_{G,p}$ for $p=\sfrac 32,\infty$, and parameters $M_t,\lambda,\nu,\nu',N_*,M_*$ such that \eqref{eq:inv:div:NM} and \eqref{eq:inverse:div:DN:G} are satisfied. There exists a smooth, non-negative scalar function $\pi$ such that \eqref{eq:inv:div:extra:pointwise} holds.
    \item\label{i:st:sample:2} There exists an incompressible vector field $v$, associated material derivative $D_t=\partial_t + v\cdot \nabla$, a volume preserving diffeomorphism $\Phi$, inverse flow $\Phi^{-1}$, and parameter $\lambda'$ such that \eqref{eq:DDpsi2}--\eqref{eq:DDv} are satisfied.
    \item\label{i:st:sample:3} There exists a zero mean scalar function $\varrho$, a mean-zero tensor potential $\vartheta$, constants $\const_{*,p}$ for $p=\sfrac 32,\infty$, and parameters $\mu,\Upsilon,\Upsilon',\Lambda,\Ndec,\dpot$ such that \eqref{item:inverse:i}--\eqref{item:inverse:iii} and \eqref{eq:DN:Mikado:density}--\eqref{eq:inverse:div:parameters:1} are satisfied.
    \item\label{i:st:sample:5} The symmetric stress $S=\divH(G\varrho\circ\Phi)$ and nonlocal error $E$ satisfy the conclusions in \eqref{eq:inverse:div}, \eqref{item:div:local:i}--\eqref{item:div:nonlocal}, as well as the conclusion \eqref{eq:inv:div:extra:conc} from Remark~\ref{rem:pointwise:inverse:div}.
    \item\label{i:st:sample:6} There exist integers $N_\circ,M_\circ,K_\circ$ such that \eqref{eq:inv:div:wut}--\eqref{eq:riots:4} 
    are satisfied, and as a consequence conclusions \eqref{eq:inverse:div:error:stress}--\eqref{eq:inverse:div:error:stress:bound} hold.
\end{enumerate}
\smallskip

\noindent\textbf{Part 2: Additional assumptions}
\begin{enumerate}[(i)]
\item There exists a large positive integer $N_{**}$ and integers positive $\NcutLarge, \NcutSmall$ such that we have the additional inequalities 
\begin{subequations}\label{helping:matt:shorten}
\begin{align}
N_*-2\dpot - \NcutLarge - N_{**} - 3 \geq M_* \, , \label{i:st:sample:wut} \\
M_*-\NcutSmall - 1 \geq 2 N_\circ \, , \label{i:st:sample:wut:wut} \\
N_{**} \geq 2\dpot + 3
\label{i:st:sample:wut:wut:wut}
\end{align}
\end{subequations}
\item\label{i:st:sample:7} There exist parameters $\Gamma=\Lambda^\alpha$ for $0<\alpha\ll 1$ and $\de_{\rm tiny}$ satisfying
\begin{subequations}
\begin{align}
\NcutSmall&\leq \NcutLarge \, , \label{eq:sample:prop:Ncut:1} \\
\left(\const_{G,\infty} + 1 \right) \left( {\const_{*,\infty}\Upsilon'\Upsilon^{-2}} + 1 \right) \Gamma^{-\NcutSmall} &\leq \de_{\rm tiny}\, ,  \const_{G,\sfrac 32} \, , \const_{*,\sfrac 32} {\Upsilon'\Upsilon^{-2}} \, , \label{eq:sample:prop:Ncut:2} \\
 2 \Ndec + 4 \leq N_* - N_{**} - \NcutLarge - 3 \dpot - 3  \, , &\qquad \NcutSmall \leq M_t - 1  \, , \label{eq:sample:prop:Ncut:3} \\
\label{eq:sample:prop:decoup}
(\Lambda \Gamma)^{4}  &\leq  \left( \frac{\mu}{2 \pi \sqrt{3} \Gamma \max(\lambda,\lambda')}\right)^{\Ndec} \, .
\end{align}
\end{subequations}
\item\label{i:st:sample:8} There exists a parameter $\bar{m}$ and an increasing sequence of frequencies $\{\mu_0, \cdots, \mu_{\bar m}\}$ satisfying
\begin{subequations}
\begin{align}
\mu < \mu_0 &< \cdots< \mu_{\bar m -1} \leq \Lambda <\Lambda\Gamma < \mu_{\bar m} \, , \label{eq:sample:prop:par:00} \\
\max(\lambda,&\lambda') \Ga \left(  \mu_{m-1}^{-2} \mu_m + {\mu}^{-1} \right) \leq 1 \, , 
\label{eq:sample:prop:parameters:0} \\
\const_{G,\sfrac 32} \const_{*,\sfrac 32} & \nu \Gamma (\max(\lambda,\lambda')\Ga)^{\lfloor \sfrac \dpot { 4} \rfloor} \left(\max\left({\mu}^{-1},\mu_m \mu_{m-1}^{-2}\right) \right)^{\lfloor \sfrac \dpot { 4} \rfloor} \notag\\
&\qquad \times (\mu_{\bar m})^{5+K_\circ} \left(1 + \frac{\max\{\nu', \const_v \mu_{\bar m} \}}{\nu
}\right)^{M_\circ} \leq 1 \,, \label{eq:sample:riots:4}\\
\const_{G,\sfrac 32} \nu \Gamma \const_{*,\sfrac 32} &\left( \frac{\Lambda\Gamma}{\mu_{\bar m}} \right)^{N_{**}} (\mu_{\bar m})^{8+K_\circ} \left(1 + \frac{\max\{\nu', \const_v \mu_{\bar{m}} \}}{\nu
}\right)^{M_\circ} \leq 1 \, , \label{eq:sample:riot:4:4}
\end{align}
\end{subequations}
for all ${1\leq m \leq \bar m}$.
\end{enumerate}
\smallskip

\noindent\textbf{Part 3: Pressure increment}
\begin{enumerate}[(i)]
\item There exists a pressure increment $\si_S$, where we have a decomposition
\begin{equation}\label{d:press:stress:sample}
\si_S=\si_S^+ - \si_S^- = \sigma_S^* + \sum_{m = 0}^{\bar{m}} \si_S^{m} \, .
\end{equation}
\item\label{sample2.item2} $\si_S^{+}$ dominates derivatives of $S$ with suitable weights, so that
\begin{align}\label{est.S.by.pr.final2}
    \left|D^N D_{t}^M S\right|
    \lec (\si_S^{+}  + \de_{\rm tiny}) \left(\Lambda\Ga\right)^N\MM{M,M_{t},\nu\Gamma,\nu'\Gamma} \, .
\end{align}
for all $N \leq N_* - \floor{\sfrac{\dpot}{2}}$, $M\leq M_*$.
\item\label{sample2.item4} $\si_S^+$ dominates derivatives of itself with suitable weights, so that
\begin{align}
    \label{est.S.prbypr.pt}
    \left|D^N D_{t}^M \si_S^{+}\right|
    &\lec (\si_S^{+}+ \de_{\rm tiny}) \left(\Lambda\Ga\right)^N\MM{M,M_{t}-\NcutSmall,\nu\Ga,\nu'\Ga}
\end{align}
for all $N \leq N_* - \floor{\sfrac{\dpot}{2}} - \NcutLarge $, $M\leq M_* - \NcutSmall$. 
\item\label{sample2.item1} $\si_S^+$ and $\si_S^-$ have the same size as $S$, so that for $p=\sfrac 32,\infty$,
    \begin{align}
        \label{est.S.pr.p}
        \norm{\si_S^{+}}_{p},\, \norm{\si_S^{-}}_{p}&\lec  \const_{G,p} \const_{*,p}  \Upsilon' \Upsilon^{-2} \, .
    \end{align}
\item\label{sample2.item5} $\pi$ dominates $\si_S^-$ and its derivatives with suitable weights, so that
\begin{align}
    \label{est.S.prminus.pt}
    \left|D^N D_{t}^M \si_S^{-}\right|
    &\lec {\const_{*,\sfrac 32} \Upsilon^{-2} \Upsilon'} \pi
    (\max (\la,\la')\Ga)^N\MM{M,M_{t}-\NcutSmall,\nu\Gamma,\nu'\Gamma}
\end{align}
for all $N \leq N_* - \floor{\sfrac{\dpot}{2}} - \NcutLarge $, $M\leq M_* - \NcutSmall$. 
\item\label{sample2.item7} 
We have the support properties 
\begin{align}
    \supp(\si_S^+) &\subseteq \supp(S) \, , \quad
    \supp(\si_S^-) \subseteq \supp(G) \label{est.S.pr.p.support:1} \, .
\end{align}
\end{enumerate}
\smallskip

\noindent\textbf{Part 4: Current error}
\begin{enumerate}[(i)]
\item There exists a current error $\phi$, where we have the decomposition and equalities
\begin{subequations}\label{d:cur:error:stress:sample}
\begin{align}
\phi &= \phi_S^* + \sum_{m=0}^{\bar m} \phi_S^m \\
\div \phi^m_S(t,x) &= D_t \sigma^{m}_S(t,x) - \int_{\T^3} D_t \sigma_S^m(t,x') \,dx' \, ,  \\
\div \phi^*_S(t,x) &= D_t \sigma^{*}_S(t,x) - \int_{\T^3} D_t \sigma_S^*(t,x') \,dx' \, .
\end{align}
\end{subequations}
\item\label{sample2.item3} $\phi^m_{S}$ can be written as $\phi^m_S = \phi^{m,l}_S + \phi^{m,*}_S$, and for $1\leq m \leq \bar m$, these satisfy
\begin{subequations}
\begin{align}
\label{est.S.phi.1}
\norm{D^N D_{t}^M \phi_{S}^m}_{\sfrac 32}
&\lec {\nu\Ga^2 \const_{G,\sfrac 32} \const_{*,\sfrac 32}  \Upsilon' \Upsilon^{-2} \mu_{m-1}^{-2}\mu_m}  \left(\min(\mu_m,\Lambda\Gamma)\right)^N\MM{M,M_{t}-\NcutSmall-1,\nu\Ga,\nu'\Ga}\, ,\\
\norm{D^N D_{t}^M \phi_{S}^m}_{\infty}
&\lec {\nu\Ga^2 \const_{G,\infty} \const_{*,\sfrac 32} \Upsilon' \Upsilon^{-2} \left(\frac{\min(\mu_m,\Lambda\Gamma)}{ \mu}\right)^{\sfrac 43} \mu_{m-1}^{-2}\mu_m} \notag\\
&\qquad \qquad \qquad \times \left(\min(\mu_m,\Lambda\Gamma)\right)^N\MM{M,M_{t}-\NcutSmall-1,\nu\Ga,\nu'\Ga}\, , \label{est.S.phi.infty} \\
\left|D^N D_{t}^M \phi_{S}^{m,l}\right|
&\lec {\nu\Ga^2 \pi \const_{*,\sfrac 32} \Upsilon' \Upsilon^{-2} \left(\frac{\min(\mu_m,\Lambda\Gamma)}{ \mu}\right)^{\sfrac 43} \mu_{m-1}^{-2}\mu_m}  \notag\\
&\qquad \qquad \qquad \times \left(\min(\mu_m,\Lambda\Gamma)\right)^N\MM{M,M_{t}-\NcutSmall-1,\nu\Ga,\nu'\Ga}\, , \label{est.S.by.pr.final3}
\end{align}
\end{subequations}
for all $N \leq N_* - 2\dpot-{\NcutLarge}$, $M\leq M_* - {\NcutSmall}-1$. For $m=0$ and the same range of $N$ and $M$, $\phi_S^{m}$ and $\phi_S^{m,l}$ satisfy identical bounds but with $\mu_{m-1}^2\mu_m$ replaced with $\Gamma\mu^{-1}$ and $\min(\mu_m,\Lambda\Gamma)$ replaced with $\mu_0$ in all three bounds. Furthermore, the nonlocal portions satisfy the improved estimate
\begin{align}
\label{est.S.by.pr.final4}
    \norm{D^N D_{t}^M \phi_{S}^{m,*}}_\infty
    &\lec  \left(\min(\mu_m,\Lambda\Gamma)\right)^{N-K_\circ} (\max(\lambda,\lambda')\Ga)^{\lfloor \sfrac \dpot { 4} \rfloor} \left(\max\left({\mu}^{-1},\mu_m \mu_{m-1}^{-2}\right) \right)^{\lfloor \sfrac \dpot { 4} \rfloor} (\nu\Gamma)^M
\end{align}
for all $N \leq N_\circ, M\leq M_\circ$, and the remainder term $\phi_S^*$ satisfies the improved estimate
\begin{align}
\label{est.S.by.pr.final.star}
    \norm{D^N D_{t}^M \phi_{S}^{*}}_\infty
    &\lec (\Lambda\Gamma)^{-K_\circ}  (\max(\lambda,\lambda')\Ga)^{\lfloor \sfrac \dpot { 4} \rfloor} \left(\max\left({\mu}^{-1},\mu_m \mu_{m-1}^{-2}\right) \right)^{\lfloor \sfrac \dpot { 4} \rfloor}\left(\Lambda\Ga\right)^N (\nu\Gamma)^M
\end{align}
in the same range of $N$ and $M$.
\item We have the support properties\footnote{For any $\Omega\in\T^3$, we use $\Omega\circ \Phiik$ to refer to the space-time set $\Phiik^{-1}(t,\cdot)\Omega$ whose characteristic function is annihilated by $D_t$.}
\begin{align}
    \supp(\phi^{m,l}_S)
    &\subseteq \supp G \cap B\left( \supp \vartheta, 2\mu_{m-1}^{-1} \right) \circ \Phi \,  \textnormal{ for } \, 1\leq m \leq \bar{m} \, , \qquad \supp\left(\phi_S^{0,l} \right) \subseteq \supp G \,  .  \label{est.S.pr.p.support:2}
\end{align}
\item For all $M\leq M_* -\NcutSmall - 1$, we have that the mean $\langle D_t \si_S \rangle$ satisfies
\begin{equation}\label{est:mean.DtsiS}
    \left|\frac{d^M}{dt^M}\langle D_t \si_S \rangle\right|\lec (\La\Ga)^{-K_\circ}(\max(\lambda,\lambda')\Ga)^{\lfloor \sfrac \dpot { 4} \rfloor} \mu^{-\lfloor \sfrac \dpot { 4} \rfloor}
    \MM{M,M_t - \NcutSmall, -1, \nu\Ga, \nu'\Ga}\, . 
\end{equation}
\end{enumerate}
\end{proposition*}

\begin{proof}
\textbf{Step 1: Defining and estimating $\si_S$ to verify \eqref{est.S.by.pr.final2}--\eqref{est.S.pr.p.support:1}.} \, From~\eqref{eq:divH:formula} of Proposition \ref{prop:intermittent:inverse:div}, we have that $S$ can be written as
$$ S = \sum_{j=0}^{\mathcal{C}_\mathcal{H}} H^{\alpha{(j)}} \rho^{\beta{(j)}} \circ \Phi \, , $$
where $H^{\alpha{(j)}}$ and $\rho^{\beta{(j)}}$ satisfy the bounds in~\eqref{eq:inverse:div:sub:1}, \eqref{eq:inverse:div:sub:2}. In addition, we have the pointwise bounds on $H^{\alpha{(j)}}$ in terms of $\pi$ given by~\eqref{eq:inv:div:extra:conc} in Remark~\ref{rem:pointwise:inverse:div}. For each $0\leq j \leq \mathcal{C}_\mathcal{H}$, we shall apply Lemma~\ref{lem:pr.st} with the following choices, where we have used the convention set out in Remark~\ref{rem:bar:variables} to rewrite the symbols from Lemma~\ref{lem:pr.st} with bars above on the left-hand side of the equalities below, while the right-hand side are parameters given in the assumptions of this Proposition:
\begin{align*}
    &\overline v=v \, , \quad \overline N_\dagger = N_*-\lfloor\sfrac \dpot 2 \rfloor \, , \quad \ov{M}_\dagger = M_* \, , \quad \ov{M}_t = M_t \, ,  \\
    &\ov{H} = H^{\alpha{(j)}}, \quad \ov{\const}_{G,\sfrac 32} = {\const}_{G,\sfrac 32} \, , \quad \ov{\const}_{G,\infty} = \const_{G,\infty} \, , \\
    &\ov{\rho} = \rho^{\beta{(j)}}, \quad \ov{\const}_{\rho,\sfrac 32} = \const_{*,\sfrac 32} \Upsilon' \Upsilon^{-2} \, , \quad \ov{\const}_{\rho,\infty} = \const_{*\infty} \Upsilon' \Upsilon^{-2} , \\
    &\ov{\la} = \max (\la, \la') \, , \quad \ov{\Lambda}=\Lambda, \quad \ov{\Ga}=\Ga \, , \quad  \ov{\Phi}=\Phi \, ,\\
    &\ov{\pi} = \pi \, , \quad \ov{\nu} = \nu \, , \quad {\ov{\nu}'=\nu'} \, , \quad \ov{\mu}= \mu \, , \quad \ov{\Ndec} = \Ndec \, ,
\end{align*}
and $\NcutLarge$, $\NcutSmall$, and $\delta_{\rm tiny}$ as in preliminary assumption~\eqref{i:st:sample:7}. From \eqref{eq:inverse:div:sub:main}, \eqref{eq:inv:div:extra:conc}, and \eqref{eq:inverse:div:stress:1}, we have that \eqref{est.G}-\eqref{est.S} are satisfied. Assumption \eqref{eq:sample:1:decoup} is satisfied from \eqref{eq:sample:prop:decoup}. All the assumptions in~\eqref{sample:item:assump:1} are satisfied from preliminary assumption~\eqref{i:st:sample:2} from this proposition. Finally, all assumptions in~\eqref{sample:item:assump:3} are satisfied from the additional assumption~\eqref{i:st:sample:7} from this Proposition.

We may then apply \eqref{d:st:pr:1}--\eqref{d:st:pr:4} from Lemma~\ref{lem:pr.st} to obtain for $0\leq j \leq \const_H$ the pressure increments $\si_S^j = \si_S^{+,j} - \si_S^{-,j}$, and we then collect terms to define
\begin{align*}
    \si_S^+ := \sum_{j=0}^{\mathcal{C}_\mathcal{H}} \si_S^{+,j}, \qquad \si_S^- := \sum_{j=0}^{\mathcal{C}_\mathcal{H}} \si_S^{-,j}, \qquad \si_S := \si_S^{+} - \si_S^{-} \, .
\end{align*}
From conclusions~\eqref{sample:item:2}--\eqref{sample:item:6} of Lemma~\ref{lem:pr.st}, we have that \eqref{est.S.by.pr.final2}--\eqref{est.S.pr.p.support:1} are satisfied.
\smallskip

\noindent\textbf{Step 2: Decomposing $\sigma_S$ to verify \eqref{d:press:stress:sample}, and defining and estimating $\phi_S^m$ to verify \eqref{d:cur:error:stress:sample}--\eqref{est.S.pr.p.support:2}.} \, From \eqref{d:st:pr:1}--\eqref{d:st:pr:2}, we have that
\begin{equation}\label{rep.siS}
\si_S = \sum_{j=0}^{\mathcal{C}_{\mathcal{H}}} \pr\left(H^{\alpha{(j)}}\right) \left(\mathbb{P}_{\neq 0} \pr(\rho^{\beta{(j)}})\right)\circ \Phi \, .
\end{equation}
Note further that $\pr(\rho^{\beta{(j)}})$ is $(\sfrac{\T}{\mu})^3$-periodic and has derivative cost $\Lambda\Ga$ from \eqref{eq:sample:1:conc:2}, conclusion~\eqref{sample:item:1} from Lemma~\ref{lem:pr.st}. So we use the sequence of frequencies $\mu_0,\dots,\mu_{\bar{m}}$ to apply the synthetic Littlewood-Paley decomposition (\`{a} la \eqref{eq:decomp:showing}) to $\pr(\rho^{\beta{(j)}})$ and write
\begin{align}
     \pr(\rho^{\beta{(j)}}) &= \tilde{\mathbb{P}}_{\mu_0} (\pr(\rho^{\beta{(j)}})) + \left( \sum_{m=1}^{\bar m}  \tilde{\mathbb{P}}_{(\mu_{m-1},\mu_m]} (\pr(\rho^{\beta{(j)}})) \right) + \left( \Id -  \tilde{\mathbb{P}}_{\mu_{\bar m}} \right)\pr(\rho^{\beta{(j)}}) \, . \label{eq:decomp:showing:redux} 
\end{align}
From now on, we shall abbreviate notation by writing $\mathbb{P}_{0}$ for $\tilde{\mathbb{P}}_{\mu_0}$, $\mathbb{P}_{m}$ for $\tilde{\mathbb{P}}_{(\mu_{m-1},\mu_m]}$ for $1\leq m\leq \bar m$, and $\mathbb{P}^*$ for $\Id - \tilde{\mathbb{P}}_{\mu_{\bar m}}$, so that we may use \eqref{eq:decomp:showing:redux} to write
\begin{align}
\sigma_S &= \sigma_S^* + \sum_{m=0}^{\bar m} \sigma_S^m := \sum_{j=0}^{\const_H}\pr\left(H^{\alpha(j)} \right) \mathbb{P}^* \left(\pr\left(\rho^{\beta{(j)}}\right)\right) \circ \Phi + \sum_{m=0}^{\bar m} \sum_{j=0}^{\const_H} \pr\left( H^{\alpha(j)} \right) \mathbb{P}_m \left(\pr\left(\rho^{\beta{(j)}}\right)\right)\circ \Phi \, . \label{eq:decomp:showing:redux:2}
\end{align}
We aim to apply Proposition~\ref{prop:intermittent:inverse:div} with Remarks~\ref{rem:scalar:inverse:div}, \ref{rem:pointwise:inverse:div} to the material derivative of each of the terms in \eqref{eq:decomp:showing:redux:2}, which would produce
\begin{align}
    \phi := {\phi_S^*} + \sum_{m=0}^{\bar m} \phi_S^m &=: \sum_{j=0}^{\mathcal{C}_{\mathcal{H}}} \underbrace{(\divH+\divR)\left( D_t \pr(H^{\alpha{(j)}}) \left( \mathbb{P}^*\mathbb{P}_{\neq 0} \pr(\rho^{\beta{(j)}})\right)\circ \Phi \right)}_{=:\phi^{*,j}}\notag\\
    &\qquad \qquad + \sum_{m=0}^{\bar m}\sum_{j=0}^{\mathcal{C}_{\mathcal{H}}} \underbrace{(\mathcal{H}+\divR) \left( D_t \pr(H^{\alpha{(j)}}) \left(\mathbb{P}_{m}\mathbb{P}_{\neq 0} \pr(\rho^{\beta{(j)}})\right)\circ \Phi \right)}_{=: \phi^{m,j} } \notag\\
    &= (\divH + \divR)(D_t \sigma_S^*) + \sum_{m=0}^{\bar m} (\divH+\divR)(D_t\sigma_S^m) \, . \notag 
\end{align}
Assuming that we succeed in doing so, we have at least verified \eqref{d:press:stress:sample} and \eqref{d:cur:error:stress:sample}.  Now in order to apply the inverse divergence with the pointwise bounds from Remark~\ref{rem:pointwise:inverse:div}, we first treat the low-frequency assumptions from Part 1, which are the same in all cases (irrespective of the projector on $\pr(\rho^{\beta(j)})$). Specifically, we shall use the convention from Remark~\ref{rem:bar:variables} and in all cases set
\begin{align*}
    &\ov{p}=\sfrac 32, \infty \, , \quad \ov v = v \, , \quad  \ov{N}_* = N_* - \dpot - \lfloor \sfrac{\dpot}{2} \rfloor - \NcutLarge \, , \quad \ov{M}_*=M_*-\NcutSmall -1 \, , \quad \ov M_t = M_t - \NcutSmall - 1 \, , \\
    &\ov{G} = D_t \pr(H^{\alpha{(j)}}), \quad \ov{\const}_{G,{p}} = {\nu\Ga \const_{G,p}} \, , \quad \ov \mu=\mu \, , \quad \ov\la=\max (\la,\la')\Ga\, , \quad \ov \Phi = \Phi \, , \quad \ov \lambda' = \lambda' \, , \\
    &\ov \nu=\nu\Ga\, , \quad \ov \nu' = \nu'\Gamma \, , \quad \ov\Phi=\Phi \, , \quad \ov\pi = \nu{\Ga} \pi \, , \quad \ov \Ndec = \Ndec \, , \quad \ov \dpot = \dpot \, .
\end{align*}
Then \eqref{eq:inv:div:NM} is satisfied from the additional assumption~\eqref{i:st:sample:wut}, and \eqref{eq:inverse:div:DN:G} is satisfied from the conclusion \eqref{eq:sample:1:conc:2} and the parameter choices from Step 1 which led to that conclusion. The estimates in \eqref{eq:DDpsi2}, \eqref{eq:DDpsi} and \eqref{eq:DDv} hold from assumption~\eqref{i:st:sample:2} from this Proposition. The pointwise bound in \eqref{eq:inv:div:extra:pointwise} holds with $\ov M_t = M_t - \NcutSmall-1$ and $\ov \pi = \nu \Gamma \pi$ due to \eqref{est.pr.H}, which was verified in Step 1. At this point we split into cases based on which projector is applied to $\mathbb{P}_{\neq 0} \pr(\rho^{\beta(j)})$ in \eqref{eq:decomp:showing:redux:2} and address parts 2-4 of Proposition~\ref{prop:intermittent:inverse:div}.
\smallskip

\noindent\textbf{Step 2a: Lowest shell.} \,  For the case $m=0$, we appeal to Lemma~\ref{lem:special:cases} with $q=\sfrac 32$, $\lambda=\Lambda\Gamma$, $\rho = \mathbb{P}_{\neq 0} \pr(\rho^{\beta(j)})$, and $\alpha$ such that $\lambda^\alpha$ in \eqref{eq:lowest:shell:inverse} is equal to $\Gamma$.  Specifically, to verify the assumptions in Part 2 of Proposition~\ref{prop:intermittent:inverse:div}, we set for $p=\sfrac 32,\infty$
\begin{align*}
    &\ov{\varrho} = \mathbb{P}_0 \mathbb{P}_{\neq 0} \pr(\rho^{\beta(j)}) \, , \quad \ov{\vartheta}\textnormal{ as defined in \eqref{eq:lowest:shell:inverse}} \, , \quad \ov{\const}_{*,p} = \Gamma \const_{*,\sfrac 32} \Upsilon^{-2}\Upsilon' \left( \frac{\mu_0}{\mu} \right)^{\frac 43 - \frac 2p} \, , \\
    &\ov \mu = \mu \, , \quad \ov \Upsilon = \ov \Upsilon'= \mu \, , \quad \ov \Lambda = \mu_0 \, , \quad \ov \dpot = \dpot \, .
\end{align*}
Then \eqref{eq:moll:1:as} is satisfied with $\const_{p,\sfrac 32}=\const_{*,\sfrac 32}\Upsilon^{-2}\Upsilon'$ and $\lambda=\Lambda\Gamma$ from standard Littlewood-Paley theory, \eqref{eq:sample:1:conc:2}, and the choices from Step 1 which led to that conclusion, and so from \eqref{eq:lowest:shell:inverse} we have that \eqref{eq:DN:Mikado:density} is satisfied. From \eqref{eq:sample:prop:decoup}, \eqref{eq:sample:prop:par:00}, \eqref{eq:sample:prop:parameters:0}, the choice of $\ov N_*$ above, \eqref{eq:sample:1:conc:2}, and \eqref{eq:sample:prop:Ncut:3}, we have that \eqref{eq:inverse:div:parameters:0}--\eqref{eq:inverse:div:parameters:1} are satisfied. Continuing onto the nonlocal assumptions from Proposition~\ref{prop:intermittent:inverse:div}, we have that \eqref{eq:inv:div:wut}--\eqref{eq:inverse:div:v:global:parameters} are satisfied from preliminary assumption~\eqref{i:st:sample:6} and \eqref{i:st:sample:wut:wut}.  We have that \eqref{eq:riots:4} is satisfied from \eqref{eq:sample:riots:4}. We then appeal to the conclusions \eqref{eq:inverse:div}--\eqref{eq:inverse:div:error:1} and \eqref{eq:inverse:div:error:stress}--\eqref{eq:inverse:div:error:stress:bound} to conclude as follows. First, we set
$$ \phi_S^{0,l} = \divH (D_t \sigma^0_S) \, , \qquad \phi_S^{0,*} = \divR (D_t \sigma_S^0) \, . $$
From \eqref{eq:inverse:div:stress:1}, we obtain both \eqref{est.S.phi.1} and \eqref{est.S.phi.infty}, but with the appropriate modifications for $m=0$ as indicated. The pointwise bound in \eqref{est.S.by.pr.final3} holds due to \eqref{eq:inv:div:extra:conc}, \eqref{eq:inverse:div:sub:1}, and \eqref{eq:divH:formula}. Next, we obtain \eqref{est.S.by.pr.final4} for $m=0$ from \eqref{eq:inverse:div:error:stress:bound}. Finally, we obtain \eqref{est.S.pr.p.support:2} from \eqref{eq:inverse:div:linear}, concluding the proof of the desired conclusions for $m=0$
.
\smallskip

\noindent\textbf{Step 2b: Intermediate shells.} \, For the cases $1\leq m\leq \bar m$, we appeal to Lemma~\ref{lem:LP.est} with $q=\sfrac 32$ and $\rho = \mathbb{P}_{\neq 0} \pr(\rho^{\beta(j)})$. Specifically, to verify the assumptions in Part 2 of Proposition~\ref{prop:intermittent:inverse:div}, we set for $p=\sfrac 32, \infty$
\begin{align*}
&\ov{\varrho} = \mathbb{P}_m \mathbb{P}_{\neq 0} \pr(\rho^{\beta(j)})  \, , \quad \ov\vartheta = \mu_{m-1}^{-\dpot}\Theta_\rho^{\mu_{m-1},\mu_m} \,  \textnormal{ as defined in Lemma~\ref{lem:LP.est}} \, , \\
&\ov{\const}_{*,p} = \const_{*,\sfrac 32} \Upsilon^{-2}\Upsilon' \left( \frac{\min(\mu_m,\Lambda\Gamma)}{\mu} \right)^{\frac 43 - \frac 2p} \, , \quad  \ov \Upsilon = \mu_{m-1} \, , \quad \ov \Upsilon' = \ov \Lambda = \min(\mu_{m},\Gamma\Lambda) \, , \\
&\ov \dpot = \dpot \, , \quad \ov \mu = \mu \, , \qquad \alpha \textnormal{ as in the previous substep} . 
\end{align*}
Then \eqref{eq:gen:id:assump} is satisfied exactly as in the previous substep, and so from \eqref{eq:LP:equality}--\eqref{eq:LP:div:estimates} we have that \eqref{eq:DN:Mikado:density} is satisfied. As before, we use \eqref{eq:sample:prop:decoup}, \eqref{eq:sample:prop:par:00}, \eqref{eq:sample:prop:parameters:0}, the choice of $\ov N_*$ above, \eqref{eq:sample:1:conc:2}, and \eqref{eq:sample:prop:Ncut:3} to see that \eqref{eq:inverse:div:parameters:0}--\eqref{eq:inverse:div:parameters:1} are satisfied. Continuing onto the nonlocal assumptions from Proposition~\ref{prop:intermittent:inverse:div}, we have that \eqref{eq:inv:div:wut}--\eqref{eq:inverse:div:v:global:parameters} are satisfied as in the previous substep, and \eqref{eq:riots:4} is satisfied from \eqref{eq:sample:riots:4}. We then appeal to the conclusions \eqref{eq:inverse:div}--\eqref{eq:inverse:div:error:1} and \eqref{eq:inverse:div:error:stress}--\eqref{eq:inverse:div:error:stress:bound} to conclude as follows. First, we set
$$ \phi_S^{m,l} = \divH (D_t \sigma^m_S) \, , \qquad \phi_S^{m,*} = \divR (D_t \sigma_S^m) \, . $$
From \eqref{eq:inverse:div:stress:1}, we obtain both \eqref{est.S.phi.1} and \eqref{est.S.phi.infty}. The pointwise bound in \eqref{est.S.by.pr.final3} holds due to \eqref{eq:inv:div:extra:conc}, \eqref{eq:inverse:div:sub:1}, and \eqref{eq:divH:formula}. Next, we obtain \eqref{est.S.by.pr.final4} from \eqref{eq:inverse:div:error:stress:bound}. Finally, we obtain \eqref{est.S.pr.p.support:2} from \eqref{eq:inverse:div:linear} and \eqref{eq:LP:div:support}, concluding the proof for $1\leq m \leq \bar m$.
\smallskip

\noindent\textbf{Step 2c: Highest shell.} \, For the case with the highest shell, corresponding to the projector $\mathbb{P}^*$ from \eqref{eq:decomp:showing:redux:2}, we appeal to Lemma~\ref{lem:special:cases} with $q=\sfrac 32$, $\lambda=\Lambda\Gamma$, $\rho = \mathbb{P}_{\neq 0} \pr(\rho^{\beta(j)})$. Specifically, to verify the assumptions in Part 2 of Proposition~\ref{prop:intermittent:inverse:div}, we set for $p=\sfrac 32,\infty$
\begin{align*}
&\ov{\varrho} = \mathbb{P}^* \mathbb{P}_{\neq 0} \pr(\rho^{\beta(j)})  \, , \quad \ov\vartheta = \vartheta \,  \textnormal{ as defined in \eqref{eq:remainder:inverse}} \, , \\
&\ov{\const}_{*,p} = \left( \frac{\Lambda\Gamma}{\mu_{\bar m}} \right)^{N_{**}} \const_{*,\sfrac 32} \Upsilon^{-2}\Upsilon' (\Lambda\Gamma)^3 \, , \quad  \ov \Upsilon =  \ov \Upsilon' = \mu \, , \quad \ov \Lambda = \Gamma\Lambda \, , \\
&\ov \dpot = {0} \, , \quad \ov N_* = N_* - \NcutLarge - N_{**} - 3 \, . 
\end{align*}
We note that we have altered the definition of $N_*$ compared to the previous two substeps for convenience. But from \eqref{i:st:sample:wut:wut:wut}, we have in fact made it \emph{smaller}, so that the low-frequency assumptions from the inverse divergence are still satisfied. Then \eqref{eq:moll:1:as} is satisfied exactly as in the first substep, and so from \eqref{eq:remainder:inverse} we have that \eqref{eq:DN:Mikado:density} is satisfied. We use \eqref{eq:sample:prop:decoup}, \eqref{eq:sample:prop:par:00}, \eqref{eq:sample:prop:parameters:0}, the altered choice of $\ov N_*$ above, \eqref{eq:sample:1:conc:2}, and \eqref{eq:sample:prop:Ncut:3} to see that \eqref{eq:inverse:div:parameters:0}--\eqref{eq:inverse:div:parameters:1} are satisfied. Continuing onto the nonlocal assumptions from Proposition~\ref{prop:intermittent:inverse:div}, we have that \eqref{eq:inv:div:wut}--\eqref{eq:inverse:div:v:global:parameters} are satisfied as in the previous substep, and \eqref{eq:riots:4} is satisfied from \eqref{eq:sample:riot:4:4}. We then appeal to the conclusions \eqref{eq:inverse:div}--\eqref{eq:inverse:div:error:1} and \eqref{eq:inverse:div:error:stress}--\eqref{eq:inverse:div:error:stress:bound} to conclude as follows. First, we set
$$ \phi_S^* = (\divH+\divR) (D_t \sigma_S^*) \, . $$
We may ignore \eqref{eq:inverse:div:stress:1} since $\dpot=0$. Then the only conclusion we require is \eqref{est.S.by.pr.final.star}, which follows from \eqref{eq:inverse:div:error:stress:bound}.

\medskip

\noindent\textbf{Step 3: Verification of \eqref{est:mean.DtsiS}.}  Since the vector field $v$ is incompressible, $\frac{d^M}{dt^M} \langle D_t \si_S \rangle = \langle D_t^{M+1} \si_S \rangle$.
From \eqref{rep.siS}, we have
\begin{align*}
    D_t^{M+1}\si_S = \sum_{j=0}^{\mathcal{C}_{\mathcal{H}}} D_t^{M+1}\pr\left(H^{\alpha{(j)}}\right) \left(\mathbb{P}_{\neq 0} \pr(\rho^{\beta{(j)}})\right)\circ \Phi \, .
\end{align*}
Since $\pr(\rho^{\beta{(j)}})$ is periodic in $(\sfrac{\T}{\mu})^2$, we have that for $M+1 \leq M_* - \NcutSmall - 1$
\begin{align*}
\bigg|\int_{\T^3}&D_t^{M+1}\pr\left(H^{\alpha{(j)}}\right) \left(\mathbb{P}_{\neq 0} \pr(\rho^{\beta{(j)}})\right)\circ \Phi dx\bigg|\\     &=\left|\int_{\T^3}D_t^{M+1}\pr\left(H^{\alpha{(j)}}\right)\circ \Phi^{-1} \Delta^{\floor{\frac{\dpot}4}}\Delta^{-\floor{\frac{\dpot}4}}\left(\mathbb{P}_{\neq 0} \pr(\rho^{\beta{(j)}})\right) \,  dx\right|\\
     &=\left|\int_{\T^3}\Delta^{\floor{\frac{\dpot}4}}\left(D_t^{M+1}\pr\left(H^{\alpha{(j)}}\right)\circ \Phi^{-1}\right) \Delta^{-\floor{\frac{\dpot}4}}\left(\mathbb{P}_{\neq 0} \pr(\rho^{\beta{(j)}})\right) dx\right|\\
     &\lec \norm{\Delta^{\floor{\frac{\dpot}4}}\left(D_t^{M+1}\pr\left(H^{\alpha{(j)}}\right)\circ \Phi^{-1}\right)}_{\sfrac32}
     \norm{\Delta^{-\floor{\frac{\dpot}4}}\left(\mathbb{P}_{\neq 0} \pr(\rho^{\beta{(j)}})\right)}_{1}\\
     &\lec \const_{G,\sfrac32} (\max(\la, \la')\Ga)^{\halfd}
     \mu^{-\halfd}  \const_{*,\sfrac32}\Upsilon^{-2}\Upsilon' 
     \MM{M+1,M_t - \NcutSmall, \nu\Ga, \nu'\Ga}
     \\
     &\leq (\La\Ga)^{-K_\circ}(\max(\lambda,\lambda')\Ga)^{\lfloor \sfrac \dpot { 4} \rfloor} {\mu}^{-\lfloor \sfrac \dpot { 4} \rfloor} \MM{M,M_t - \NcutSmall-1, \nu\Ga, \nu'\Ga}\, . 
\end{align*}
Here, we have used Lemma~\ref{lem:decoup}, \eqref{eq:sample:1:conc:2}, \eqref{eq:DDpsi:sample}, \eqref{eq:sample:riots:4}, and standard Littlewood-Paley theory. 
\end{proof}

\begin{proposition*}[\bf Pressure increment and upgrade error from current error]\label{lem.pr.invdiv2.c}
We begin with \emph{preliminary assumptions}, which include all of the assumptions and conclusions from the inverse divergence in Proposition~\ref{prop:intermittent:inverse:div} and the pointwise bounds in Remark~\ref{rem:pointwise:inverse:div}. We then include \emph{additional assumptions}, which allow for the application of Lemma~\ref{lem:pr.cu} to the current error and Proposition~\ref{prop:intermittent:inverse:div} to the material derivative of the output. We thus obtain a \emph{pressure increment} which satisfies a number of properties.  Finally, the material derivative of this pressure increment produces a \emph{current error} which itself satisfies a number of properties.
\smallskip

\noindent\textbf{Part 1: Preliminary assumptions}
\begin{enumerate}[(i)]
    \item\label{i:st:sample:1:c} There exists a scalar field $G$, constants $\const_{G,p}$ for $p=1,\infty$, and parameters $M_t,\lambda,\nu,\nu',N_*,M_*$ such that \eqref{eq:inv:div:NM} and \eqref{eq:inverse:div:DN:G} are satisfied. There exists a smooth, non-negative scalar function $\pi$ and a parameter $r_G$ such that \eqref{est.G.c.pt} holds with $H$ replaced by $G$.
    \item\label{i:st:sample:2:c} There exists an incompressible vector field $v$, associated material derivative $D_t=\partial_t + v\cdot \nabla$, a volume preserving diffeomorphism $\Phi$, inverse flow $\Phi^{-1}$, and parameter $\lambda'$ such that \eqref{eq:DDpsi2}--\eqref{eq:DDv} are satisfied.
    \item\label{i:st:sample:3:c} There exists a zero mean scalar function $\varrho$, a mean-zero tensor potential $\vartheta$, constants $\const_{*,p}$ for $p=1,\infty$, and parameters $\mu,{\Upsilon,\Upsilon'},\Lambda,\Ndec,\dpot$ such that \eqref{item:inverse:i}--\eqref{item:inverse:iii} and \eqref{eq:DN:Mikado:density}--\eqref{eq:inverse:div:parameters:1} are satisfied.
    \item\label{i:st:sample:5:c} The current error $\varphi=\divH(G\varrho\circ\Phi)$ and nonlocal error $E$ satisfy the conclusions in \eqref{eq:inverse:div}, \eqref{item:div:local:i}--\eqref{item:div:nonlocal}, as well as the conclusion \eqref{eq:inv:div:extra:conc} from Remark~\ref{rem:pointwise:inverse:div} with $\pi$ replaced by $\pi^{\sfrac 32} r_G^{-1}$.
    \item\label{i:st:sample:6:c} There exist integers $N_\circ,M_\circ,K_\circ$ such that \eqref{eq:inv:div:wut}--\eqref{eq:riots:4} are satisfied, and as a consequence conclusions \eqref{eq:inverse:div:error:stress}--\eqref{eq:inverse:div:error:stress:bound} hold.
\end{enumerate}
\smallskip

\noindent\textbf{Part 2: Additional assumptions}
\begin{enumerate}[(i)]
\item There exists a large positive integer $N_{**}$ and positive integers $\NcutLarge, \NcutSmall$ such that we have the additional inequalities 
\begin{subequations}
\begin{align}
N_*-2\dpot - \NcutLarge - N_{**} - 3 \geq M_* \, , \label{i:st:sample:wut:c} \\
M_*-\NcutSmall - 1 \geq 2 N_\circ \, , \label{i:st:sample:wut:wut:c} \\
N_{**} \geq {2\dpot} + 3
\label{i:st:sample:wut:wut:wut:c}
\end{align}
\end{subequations}
\item\label{i:st:sample:7:c} There exist parameters $\Gamma=\Lambda^\alpha$ for $0<\alpha\ll 1$, $\de_{\rm tiny}, r_\phi$, and $\delta_{\phi,p}$ for $p=1,\infty$ satisfying
\begin{subequations}
\begin{align}
0< r_\phi \leq 1 \, , &\qquad \delta_{\phi,p}^{\sfrac 32} = \const_{G,p} \const_{\ast,p}{\Upsilon'\Upsilon^{-2}} r_\phi \, , \label{eq:sample:prop:de:phi:c} \\
\NcutSmall&\leq \NcutLarge \, , \label{eq:sample:prop:Ncut:1:c} \\
\left(\const_{G,\infty} + 1 \right) \left( {\const_{*,\infty}\Upsilon'\Upsilon^{-2}} + 1 \right) \Gamma^{-\NcutSmall} &\leq {\delta_{\rm tiny}}^{\sfrac 32} \, , \const_{G,1} \, , {\const_{*,1}\Upsilon'\Upsilon^{-2}} \, , \label{eq:sample:prop:Ncut:2:c} \\
\qquad 2 \Ndec + 4 \leq N_* - N_{**} - \NcutLarge - 3 \dpot - 3  \, , &\qquad \NcutSmall \leq M_t - 1 \label{eq:sample:prop:Ncut:3:c} \\
\label{eq:sample:prop:decoup:c}
(\Lambda \Gamma)^{4}  &\leq  \left( \frac{\mu}{2 \pi \sqrt{3} \Gamma \max(\lambda,\lambda')}\right)^{\Ndec} \, .
\end{align}
\end{subequations}
\item\label{i:st:sample:8:c} There exists a parameter $\bar{m}$ and an increasing sequence of frequencies $\{\mu_0, \cdots, \mu_{\bar m}\}$ satisfying
\begin{subequations}
\begin{align}
\mu < \mu_0 &< \cdots< \mu_{\bar m -1} \leq \Lambda <\Lambda\Gamma < \mu_{\bar m} \, , \label{eq:sample:prop:par:00:c} \\
\max(\lambda,&\lambda') \Ga \left(  \mu_{m-1}^{-2} \mu_m + {\mu}^{-1} \right) \leq 1 \, , 
\label{eq:sample:prop:parameters:0:c} \\
\left(\const_{G,1} \const_{*,1} r_\phi\right)^{\sfrac 23}& \nu \Gamma (\max(\lambda,\lambda')\Ga)^{\lfloor \sfrac \dpot {{4}} \rfloor} \left(\max\left({\mu}^{-1},\mu_m \mu_{m-1}^{-2}\right) \right)^{\lfloor \sfrac \dpot {{4}} \rfloor} \notag\\
&\qquad \times (\mu_{\bar m})^{5+K_\circ} \left(1 + \frac{\max\{\nu', \const_v \mu_{\bar m} \}}{\nu
}\right)^{M_\circ} \leq 1 \,, \label{eq:sample:riots:4:c}\\
\left(\const_{G,1}\const_{*,1}r_\phi\right)^{\sfrac 23}  \nu \Gamma &\left( \frac{\Lambda\Gamma}{\mu_{\bar m}} \right)^{N_{**}} (\mu_{\bar m})^{8+K_\circ} \left(1 + \frac{\max\{\nu', \const_v \mu_{\bar{m}} \}}{\nu
}\right)^{M_\circ} \leq 1 \, , \label{eq:sample:riot:4:4:c}
\end{align}
\end{subequations}
for all ${1\leq m \leq \bar m}$.
\end{enumerate}
\smallskip

\noindent\textbf{Part 3: Pressure increment}
\begin{enumerate}[(i)]
\item There exists a pressure increment $\si_\varphi$, where we have a decomposition
\begin{equation}\label{d:press:stress:sample:c}
\si_\varphi=\si_\varphi^+ - \si_\vp^- = \sigma_\vp^* + \sum_{m = 0}^{\bar{m}} \si_\vp^{m} \, .
\end{equation}
\item\label{sample2.item2:c} $\si_\vp^{+}$ dominates derivatives of $\varphi$ with suitable weights, so that
\begin{align}\label{est.S.by.pr.final2:c}
    \left|D^N D_{t}^M \varphi\right|
    \lec \left((\si_\varphi^{+})^{\sfrac 32} r_\phi^{-1}  + \de_{\rm tiny}\right) \left(\Lambda\Ga\right)^N\MM{M,M_{t},\nu\Gamma,\nu'\Gamma} \, .
\end{align}
for all $N \leq N_* - \floor{\sfrac{\dpot}{2}}$, $M\leq M_*$.
\item\label{sample2.item4:c} $\si_\vp^+$ dominates derivatives of itself with suitable weights, so that
\begin{align}
    \label{est.S.prbypr.pt:c}
    \left|D^N D_{t}^M \si_\vp^{+}\right|
    &\lec (\si_\vp^{+}+ \de_{\rm tiny}) \left(\Lambda\Ga\right)^N\MM{M,M_{t}-\NcutSmall,\nu\Ga,\nu'\Ga}
\end{align}
for all $N \leq N_* - \floor{\sfrac{\dpot}{2}} - \NcutLarge $, $M\leq M_* - \NcutSmall$. 
\item\label{sample2.item1:c} $\si_\vp^+$ and $\si_\vp^-$ have size comparable to $\varphi$, so that
    \begin{align}
        \label{est.S.pr.p:c}
        \norm{\si_\vp^{+}}_{\sfrac 32},\, \norm{\si_\vp^{-}}_{\sfrac 32}&\lec \delta_{\phi,1} \, , \qquad \norm{\si_\vp^{+}}_{\infty},\, \norm{\si_\vp^{-}}_{\infty}\lec \delta_{\phi,\infty} \, .
    \end{align}
\item\label{sample2.item5:c} $\pi$ dominates $\si_\varphi^-$ and its derivatives with suitable weights, so that
\begin{align}
    \label{est.S.prminus.pt:c}
    \left|D^N D_{t}^M \si_\varphi^{-}\right|
    &\lec \left( \frac{r_\phi}{r_G} \right)^{\sfrac 23} \left({\const_{*,1} {\Upsilon^{-2}} \Upsilon'}\right)^{\sfrac 23} \pi
    (\max (\la,\la')\Ga)^N\MM{M,M_{t}-\NcutSmall,\nu\Gamma,\nu'\Gamma}
\end{align}
for all $N \leq N_* - \floor{\sfrac{\dpot}{2}} - \NcutLarge $, $M\leq M_* - \NcutSmall$. 
\item\label{sample2.item7:c} 
We have the support properties 
\begin{align}
    \supp(\si_\vp^+) &\subseteq \supp(\varphi) \, , \quad
    \supp(\si_\vp^-) \subseteq \supp(G) \label{est.S.pr.p.support:1:c} \, .
\end{align}
\end{enumerate}
\smallskip

\noindent\textbf{Part 4: Current error}
\begin{enumerate}[(i)]
\item There exists a current error $\phi_{\ph}$, where we have the decomposition and equalities
\begin{subequations}\label{d:cur:error:stress:sample:c}
\begin{align}
\phi_{\ph} &= \phi_\vp^* + \sum_{m=0}^{\bar m} \phi_\vp^m  
\label{S:pr:current:dec}\\
\div \phi^m_\vp(t,x) &= D_t \sigma^{m}_\vp(t,x) - \int_{\T^3} D_t \sigma_\vp^m(t,x') \,dx' \, ,  \\
\div \phi^*_\vp(t,x) &= D_t \sigma^{*}_\vp(t,x) - \int_{\T^3} D_t \sigma_\vp^*(t,x') \,dx' \, , 
\end{align}
\end{subequations}
\item\label{sample2.item3:c} $\phi^m_{\vp}$ can be written as $\phi^m_\vp = \phi^{m,l}_\vp + \phi^{m,*}_\vp$ and {for $1\leq m \leq \bar m$} these satisfy
\begin{subequations}
\begin{align}
\label{est.S.phi.1:c}
\norm{D^N D_{t}^M \phi_{\vp}^m}_{\sfrac 32}
&\lec {\nu\Ga^2 \left(\const_{G,1} \const_{*,1}  {\Upsilon'\Upsilon^{-2}} r_\phi \right)^{\sfrac 23} \mu_{m-1}^{-2}\mu_m}  \left(\min(\mu_m,\Lambda\Gamma)\right)^N\MM{M,M_{t}-\NcutSmall-1,\nu\Ga,\nu'\Ga}\, ,\\
\norm{D^N D_{t}^M \phi_{\vp}^m}_{\infty}
&\lec {\nu\Ga^2 \left(\const_{G,\infty} \const_{*,1} {\Upsilon'\Upsilon^{-2}} r_\phi\right)^{\sfrac 23} \left(\frac{\min(\mu_m,\Lambda\Gamma)}{ \mu}\right)^{\sfrac 43} \mu_{m-1}^{-2}\mu_m} \notag\\
&\qquad \qquad \qquad \times \left(\min(\mu_m,\Lambda\Gamma)\right)^N\MM{M,M_{t}-\NcutSmall-1,\nu\Ga,\nu'\Ga}\, , \label{est.S.phi.infty:c} \\
\left|D^N D_{t}^M \phi_{\vp}^{m,l}\right|
&\lec {\nu\Ga^2 \pi \left( \frac{r_\phi}{r_G}\right)^{\sfrac 23} \left(\const_{*,1} {\Upsilon'\Upsilon^{-2}} \right)^{\sfrac 23} \left(\frac{\min(\mu_m,\Lambda\Gamma)}{ \mu}\right)^{\sfrac 43} \mu_{m-1}^{-2}\mu_m}  \notag\\
&\qquad \qquad \qquad \times \left(\min(\mu_m,\Lambda\Gamma)\right)^N\MM{M,M_{t}-\NcutSmall-1,\nu\Ga,\nu'\Ga}\, , \label{est.S.by.pr.final3:c}
\end{align}
\end{subequations}
for all $N \leq N_* - 2\dpot-{\NcutLarge}$, $M\leq M_* - {\NcutSmall}-1$. For {$m=0$} and the same range of $N$ and $M$, $\phi_\varphi^{m}$ and $\phi_\varphi^{m,l}$ satisfy identical bounds but with $\mu_{m-1}^2\mu_m$ replaced with ${\Gamma\mu^{-1}}$ and $\min(\mu_m,\Lambda\Gamma)$ replaced with {$\mu_0$} in all three bounds. Furthermore, the nonlocal portions satisfy the improved estimate
\begin{align}
\label{est.S.by.pr.final4:c}
    \norm{D^N D_{t}^M \phi_{\vp}^{m,*}}_\infty
    &\lec  \left(\min(\mu_m,\Lambda\Gamma)\right)^{N-K_\circ} 
    {(\max(\lambda,\lambda')\Ga)^{\lfloor \sfrac \dpot {{4}} \rfloor} \left(\max\left({\mu}^{-1},\mu_m \mu_{m-1}^{-2}\right) \right)^{\lfloor \sfrac \dpot {{4}} \rfloor}}
    (\nu\Gamma)^M \, ,
\end{align}
for all $N \leq N_\circ, M\leq M_\circ$, and the remainder term $\phi_\vp^*$ satisfies the improved estimate
\begin{align}
\label{est.S.by.pr.final.star:c}
    \norm{D^N D_{t}^M \phi_{\vp}^{*}}_\infty
    &\lec (\Lambda\Gamma)^{-K_\circ}  
    {(\max(\lambda,\lambda')\Ga)^{\lfloor \sfrac \dpot {{4}} \rfloor} \left(\max\left({\mu}^{-1},\mu_m \mu_{m-1}^{-2}\right) \right)^{\lfloor \sfrac \dpot {{4}} \rfloor}}
    \left(\Lambda\Ga\right)^N (\nu\Gamma)^M
\end{align}
in the same range of $N$ and $M$.
\item\label{sample2.item4:c:redux} We have the support properties
\begin{align}
    \supp(\phi^{m,l}_\varphi)
    &\subseteq \supp G \cap B\left( \supp \vartheta, 2\mu_{m-1}^{-1} \right) \circ \Phi \,  \textnormal{ for } \, 1\leq m \leq \bar{m} \, , \qquad \supp\left(\phi_\varphi^{0,l} \right) \subseteq \supp G \,  . \label{est.S.pr.p.support:2:c}
\end{align}

\item\label{sample2.item5:c:redux} 
For all $M\leq M_*-\NcutSmall-1$, we have that the mean $\langle D_t \si_S \rangle$ satisfies
\begin{equation}\label{est:mean.Dtsiph}
    \left|\frac{d^M}{dt^M}\langle D_t \si_\ph \rangle\right|\les
    (\La\Ga)^{-K_\circ}
    {(\max(\lambda,\lambda')\Ga)^{\lfloor \sfrac \dpot {{4}} \rfloor} \mu^{-\lfloor \sfrac \dpot {{4}} \rfloor}}
    \MM{M,M_t - \NcutSmall-1, \nu\Ga, \nu'\Ga}
\end{equation}
\end{enumerate}
\end{proposition*}
\begin{proof}
\texttt{Step 1: Defining and estimating $\si_\vp$ to verify \eqref{est.S.by.pr.final2:c}--\eqref{est.S.pr.p.support:1:c}.} \, From~\eqref{eq:divH:formula} of Proposition \ref{prop:intermittent:inverse:div}, we have that $\vp$ can be written as
$$ \vp = \sum_{j=0}^{\mathcal{C}_\mathcal{H}} H^{\alpha{(j)}} \rho^{\beta{(j)}} \circ \Phi \, , $$
where $H^{\alpha{(j)}}$ and $\rho^{\beta{(j)}}$ satisfy the bounds in~\eqref{eq:inverse:div:sub:1}, \eqref{eq:inverse:div:sub:2}. In addition, we have the pointwise bounds on $H^{\alpha{(j)}}$ in terms of $\pi^{\sfrac 32}r_G^{-1}$ given by~\eqref{eq:inv:div:extra:conc} in Remark~\ref{rem:pointwise:inverse:div}, but with the modifications listed in preliminary assumption~\eqref{i:st:sample:1:c}. For each $0\leq j \leq \mathcal{C}_\mathcal{H}$, we shall apply Lemma~\ref{lem:pr.cu} with the following choices, where we have used the convention set out in Remark~\ref{rem:bar:variables} to rewrite the symbols from Lemma~\ref{lem:pr.cu} with bars above on the left-hand side of the equalities below, while the right-hand side are parameters given in the assumptions of this Proposition:
\begin{align*}
    &\overline v=v \, , \quad \overline N_* = N_*-\lfloor\sfrac \dpot 2 \rfloor \, , \quad \ov{M}_* = M_* \, , \quad \ov{M}_t = M_t \, ,  \\
    &\ov{H} = H^{\alpha{(j)}}, \quad \ov{\const}_{G,1} = {\const}_{G,1} \, , \quad \ov{\const}_{G,\infty} = \const_{G,\infty} \, , \\
    &\ov{\rho} = \rho^{\beta{(j)}}, \quad \ov{\const}_{\rho,1} = \const_{*,1}\Upsilon^{-2}\Upsilon' \, , \quad \ov{\const}_{\rho,\infty} = \const_{*,\infty} \Upsilon^{-2}\Upsilon' \, , \quad \ov r_G = r_G \, , \quad \ov r_\phi = r_\phi \\
    &\ov{\la} = \max (\la, \la') \, , \quad \ov{\Lambda}=\Lambda, \quad \ov{\Ga}=\Ga \, , \quad  \ov{\Phi}=\Phi \, ,\\
    &\ov{\pi} = \pi \, , \quad \ov{\nu} = \nu \, , \quad {\ov{\nu}'=\nu'} \, , \quad \ov{\mu}= \mu \, , \quad \ov{\Ndec} = \Ndec \, ,
\end{align*}
and $\NcutLarge$, $\NcutSmall$, and $\delta_{\rm tiny}$ as in preliminary assumption~\eqref{i:st:sample:7:c}. From \eqref{eq:inverse:div:sub:main}, the modified version of \eqref{eq:inv:div:extra:conc}, which is listed in preliminary assumption~\eqref{i:st:sample:1:c}, \eqref{eq:inverse:div:stress:1}, and \eqref{eq:sample:prop:de:phi:c}, we have that \eqref{est.G.c}--\eqref{est.phi.c} are satisfied. Assumption \eqref{eq:sample:1:decoup:c} is satisfied from \eqref{eq:sample:prop:decoup:c}. All the assumptions in~\eqref{sample:item:assump:1:c} are satisfied from preliminary assumption~\eqref{i:st:sample:2:c} from this proposition. Finally, all assumptions in~\eqref{sample:item:assump:3:c} are satisfied from the additional assumption~\eqref{i:st:sample:7:c} from this Proposition.

We may then apply \eqref{eq:cu:d:1}--\eqref{eq:cu:d:4} from Lemma~\ref{lem:pr.cu} to obtain for $0\leq j \leq \const_H$ the pressure increments $\si_\vp^j = \si_\vp^{+,j} - \si_\vp^{-,j}$, and we then collect terms to define
\begin{align*}
    \si_\vp^+ := \sum_{j=0}^{\mathcal{C}_\mathcal{H}} \si_\vp^{+,j}, \qquad \si_\vp^- := \sum_{j=0}^{\mathcal{C}_\mathcal{H}} \si_\vp^{-,j}, \qquad \si_\vp := \si_\vp^{+} - \si_\vp^{-} \, .
\end{align*}
From conclusions~\eqref{sample:item:2:c}--\eqref{sample:item:6:c} of Lemma~\ref{lem:pr.cu}, we have that \eqref{est.S.by.pr.final2:c}--\eqref{est.S.pr.p.support:1:c} are satisfied.
\smallskip

\noindent\texttt{Step 2: Decomposing $\sigma_\vp$ to verify \eqref{d:press:stress:sample:c}, and defining and estimating $\phi_\vp^m$ to verify \eqref{d:cur:error:stress:sample:c}--\eqref{est.S.pr.p.support:2:c}} \, From \eqref{eq:cu:d:1}--\eqref{eq:cu:d:4}, we have that
\begin{equation}\label{rep.Dtsiph} 
\si_\vp = r_\phi^{\sfrac 23} \sum_{j=0}^{\mathcal{C}_{\mathcal{H}}} \pr\left(H^{\alpha{(j)}}\right) \left(\mathbb{P}_{\neq 0} \pr(\rho^{\beta{(j)}})\right)\circ \Phi \, .
\end{equation}
Note further that $\pr(\rho^{\beta{(j)}})$ is $(\sfrac{\T}{\mu})^3$-periodic and has derivative cost $\Lambda\Ga$ from \eqref{eq:sample:1:conc:2:c}, conclusion~\eqref{sample:item:1:c} from Lemma~\ref{lem:pr.cu}. So we decompose as in \eqref{eq:decomp:showing:redux} to write
\begin{align}
     \pr(\rho^{\beta{(j)}}) &= \tilde{\mathbb{P}}_{\mu_0} (\pr(\rho^{\beta{(j)}})) + \left( \sum_{m=1}^{\bar m}  \tilde{\mathbb{P}}_{(\mu_{m-1},\mu_m]} (\pr(\rho^{\beta{(j)}})) \right) + \left( \Id -  \tilde{\mathbb{P}}_{\mu_{\bar m}} \right)\pr(\rho^{\beta{(j)}}) \, . \label{eq:decomp:showing:redux:c} 
\end{align}
Using the same abbreviations used in \eqref{eq:decomp:showing:redux:2}, from \eqref{eq:decomp:showing:redux:c} we may write
\begin{align}
\sigma_\vp &= \sigma_\vp^* + \sum_{m=0}^{\bar m} \sigma_\vp^m := r_\phi^{\sfrac 23} \sum_{j=0}^{\const_H}\pr\left(H^{\alpha(j)} \right) \mathbb{P}^* \left(\pr\left(\rho^{\beta{(j)}}\right)\right) \circ \Phi + r_\phi^{\sfrac 23} \sum_{m=0}^{\bar m} \sum_{j=0}^{\const_H} \pr\left( H^{\alpha(j)} \right) \mathbb{P}_m \left(\pr\left(\rho^{\beta{(j)}}\right)\right)\circ \Phi \, . \label{eq:decomp:showing:redux:2:c}
\end{align}
We aim to apply Proposition~\ref{prop:intermittent:inverse:div} with Remarks~\ref{rem:scalar:inverse:div}, \ref{rem:pointwise:inverse:div} to the material derivative of each of the terms in \eqref{eq:decomp:showing:redux:2:c}, which would produce
\begin{align}
    \phi_\ph := \phi_\vp^* + \sum_{m=0}^{\bar m} \phi_\vp^m &=: r_\phi^{\sfrac 23} \sum_{j=0}^{\mathcal{C}_{\mathcal{H}}} \underbrace{(\divH+\divR)\left( D_t \pr(H^{\alpha{(j)}}) \left( \mathbb{P}^*\mathbb{P}_{\neq 0} \pr(\rho^{\beta{(j)}})\right)\circ \Phi \right)}_{=:\phi^{*,j}}\notag\\
    &\qquad \qquad + r_\phi^{\sfrac 23} \sum_{m=0}^{\bar m}\sum_{j=0}^{\mathcal{C}_{\mathcal{H}}} \underbrace{(\mathcal{H}+\divR) \left( D_t \pr(H^{\alpha{(j)}}) \left(\mathbb{P}_{m}\mathbb{P}_{\neq 0} \pr(\rho^{\beta{(j)}})\right)\circ \Phi \right)}_{=: \phi^{m,j} } \notag\\
    &= (\divH + \divR)(D_t \sigma_\vp^*) + \sum_{m=0}^{\bar m} (\divH+\divR)(D_t\sigma_\vp^m) \, . \notag 
\end{align}
Assuming that we succeed in doing so, we have at least verified \eqref{d:press:stress:sample:c} and \eqref{d:cur:error:stress:sample:c}.  Now in order to apply the inverse divergence with the pointwise bounds from Remark~\ref{rem:pointwise:inverse:div}, we again first treat the low-frequency assumptions from Part 1, which are the same in all cases (irrespective of the projector on $\pr(\rho^{\beta(j)})$). Specifically, we shall use the convention from Remark~\ref{rem:bar:variables} and in all cases set
\begin{align*}
    &\ov{p}=\sfrac 32, \infty \, , \quad \ov v = v \, , \quad  \ov{N}_* = N_* - \dpot - \lfloor \sfrac{\dpot}{2} \rfloor - \NcutLarge \, , \quad \ov{M}_*=M_*-\NcutSmall -1 \, , \quad \ov M_t = M_t - \NcutSmall - 1 \, , \\
    &\ov{G} = r_\phi^{\sfrac 23} D_t \pr(H^{\alpha{(j)}}), \quad \ov{\const}_{G,\sfrac 32} = r_\phi^{\sfrac 23} {\nu\Ga \const_{G,1}^{\sfrac 23}} \, , \quad \ov \mu=\mu \, , \quad \ov\la=\max (\la,\la')\Ga\, , \quad \ov \Phi = \Phi \, , \quad \ov \lambda' = \lambda' \, , \\
    &\ov \nu=\nu\Ga\, , \quad \ov \nu' = \nu'\Gamma \, , \quad \ov\Phi=\Phi \, , \quad \ov\pi = \nu{\Ga} \pi r_G^{-\sfrac 23} \, , \quad \ov \Ndec = \Ndec \, , \quad \ov \dpot = \dpot \, , \quad \ov \const_{G,\infty} = r_\phi^{\sfrac 23} \nu \Gamma \const_{G,\infty}^{\sfrac 23} \,  .
\end{align*}
Then \eqref{eq:inv:div:NM} is satisfied from the additional assumption~\eqref{i:st:sample:wut:c}, and \eqref{eq:inverse:div:DN:G} is satisfied from the conclusion \eqref{eq:sample:1:conc:2:c} and the parameter choices from Step 1 which led to that conclusion. The estimates in \eqref{eq:DDpsi2}, \eqref{eq:DDpsi} and \eqref{eq:DDv} hold from assumption~\eqref{i:st:sample:2:c} from this Proposition. The pointwise bound in \eqref{eq:inv:div:extra:pointwise} holds with $\ov M_t = M_t - \NcutSmall-1$ and $\ov \pi = \nu\Gamma\pi r_G^{-\sfrac 23}$ due to \eqref{est.pr.H.c}, which was verified in Step 1. At this point we split into cases based on which projector is applied to $\mathbb{P}_{\neq 0} \pr(\rho^{\beta(j)})$ in \eqref{eq:decomp:showing:redux:2:c} and address parts 2-4 of Proposition~\ref{prop:intermittent:inverse:div}.
\smallskip

\noindent\texttt{Step 2a: Lowest shell.} \,  For the case $m=0$, we appeal to Lemma~\ref{lem:special:cases} with $q=\sfrac 32$, $\lambda=\Lambda\Gamma$, $\rho = \mathbb{P}_{\neq 0} \pr(\rho^{\beta(j)})$, and $\alpha$ such that $\lambda^\alpha$ in \eqref{eq:lowest:shell:inverse} is equal to $\Gamma$.  Specifically, to verify the assumptions in Part 2 of Proposition~\ref{prop:intermittent:inverse:div}, we set for $p=\sfrac 32,\infty$
\begin{align*}
    &\ov{\varrho} = \mathbb{P}_0 \mathbb{P}_{\neq 0} \pr(\rho^{\beta(j)}) \, , \quad \ov{\vartheta}\textnormal{ as defined in \eqref{eq:lowest:shell:inverse}} \, , \quad \ov{\const}_{*,p} = \Gamma \left(\const_{*,1} \Upsilon^{-2}\Upsilon'\right)^{\sfrac 23} \left( \frac{\mu_0}{\mu} \right)^{\frac 43 - \frac 2p} \, , \\
    &\ov \mu = \mu \, , \quad \ov \Upsilon = \ov \Upsilon' = \mu \, , \quad \ov \Lambda = \mu_0 \, , \quad \ov \dpot = \dpot \, .
\end{align*}
Then \eqref{eq:moll:1:as} is satisfied with $\const_{p,\sfrac 32}=\left(\const_{*,1}\Upsilon^{-2}\Upsilon'\right)^{\sfrac 23}$ and $\lambda=\Lambda\Gamma$ from standard Littlewood-Paley theory, \eqref{eq:sample:1:conc:2:c}, and the choices from Step 1 which led to that conclusion, and so from \eqref{eq:lowest:shell:inverse} we have that \eqref{eq:DN:Mikado:density} is satisfied. From \eqref{eq:sample:prop:decoup:c}, \eqref{eq:sample:prop:par:00:c}, \eqref{eq:sample:prop:parameters:0:c}, the choice of $\ov N_*$ above, \eqref{eq:sample:1:conc:2:c} and \eqref{eq:sample:1:conc:2:c:infty}, and \eqref{eq:sample:prop:Ncut:3:c}, we have that \eqref{eq:inverse:div:parameters:0}--\eqref{eq:inverse:div:parameters:1} are satisfied. Continuing onto the nonlocal assumptions from Proposition~\ref{prop:intermittent:inverse:div}, we have that \eqref{eq:inv:div:wut}--\eqref{eq:inverse:div:v:global:parameters} are satisfied from preliminary assumption~\eqref{i:st:sample:6:c} and \eqref{i:st:sample:wut:wut:c}.  We have that \eqref{eq:riots:4} is satisfied from \eqref{eq:sample:riots:4:c}. We then appeal to the conclusions \eqref{eq:inverse:div}--\eqref{eq:inverse:div:error:1} and \eqref{eq:inverse:div:error:stress}--\eqref{eq:inverse:div:error:stress:bound} to conclude as follows. First, we set
$$ \phi_\vp^{0,l} = \divH (D_t \sigma^0_\vp) \, , \qquad \phi_\vp^{0,*} = \divR (D_t \sigma_\vp^0) \, . $$
From \eqref{eq:inverse:div:stress:1}, we obtain both \eqref{est.S.phi.1:c} and \eqref{est.S.phi.infty:c}, but with the appropriate modifications for $m=0$ as indicated. The pointwise bound in \eqref{est.S.by.pr.final3:c} holds due to \eqref{eq:inv:div:extra:conc}, \eqref{eq:inverse:div:sub:1}, and \eqref{eq:divH:formula}. Next, we obtain \eqref{est.S.by.pr.final4:c} for $m=0$ from \eqref{eq:inverse:div:error:stress:bound}. Finally, we obtain \eqref{est.S.pr.p.support:2:c} from \eqref{eq:inverse:div:linear}, concluding the proof of the desired conclusions for $m=0$
.
\smallskip

\noindent\texttt{Step 2b: Intermediate shells.} \, For the cases $1\leq m\leq \bar m$, we appeal to Lemma~\ref{lem:LP.est} with $q=\sfrac 32$ and $\rho = \mathbb{P}_{\neq 0} \pr(\rho^{\beta(j)})$. Specifically, to verify the assumptions in Part 2 of Proposition~\ref{prop:intermittent:inverse:div}, we set for $p=\sfrac 32,\infty$
\begin{align*}
&\ov{\varrho} = \mathbb{P}_m \mathbb{P}_{\neq 0} \pr(\rho^{\beta(j)})  \, , \quad \ov\vartheta = \mu_{m-1}^{-\dpot}\Theta_\rho^{\mu_{m-1},\mu_m} \,  \textnormal{ as defined in Lemma~\ref{lem:LP.est}} \, , \\
&\ov{\const}_{*,p} = \left(\const_{*,1} \Upsilon^{-2}\Upsilon'\right)^{\sfrac 23} \left( \frac{\min(\mu_m,\Lambda\Gamma)}{\mu} \right)^{\frac 43 - \frac 2p} \, , \quad  \ov \Upsilon = \mu_{m-1} \, , \quad \ov \Upsilon' = \ov \Lambda = \min(\mu_{m},\Gamma\Lambda) \, , \\
&\ov \dpot = \dpot \, , \qquad \ov \mu = \mu \, , \qquad \alpha \textnormal{ as in the previous substep} \,  . 
\end{align*}
Then \eqref{eq:gen:id:assump} is satisfied exactly as in the previous substep, and so from \eqref{eq:LP:equality}--\eqref{eq:LP:div:estimates} we have that \eqref{eq:DN:Mikado:density} is satisfied. As before, we use \eqref{eq:sample:prop:decoup:c}, \eqref{eq:sample:prop:par:00:c}, \eqref{eq:sample:prop:parameters:0:c}, the choice of $\ov N_*$ above, \eqref{eq:sample:1:conc:2:c} and \eqref{eq:sample:1:conc:2:c:infty}, and \eqref{eq:sample:prop:Ncut:3:c} to see that \eqref{eq:inverse:div:parameters:0}--\eqref{eq:inverse:div:parameters:1} are satisfied. Continuing onto the nonlocal assumptions from Proposition~\ref{prop:intermittent:inverse:div}, we have that \eqref{eq:inv:div:wut}--\eqref{eq:inverse:div:v:global:parameters} are satisfied as in the previous substep, and \eqref{eq:riots:4} is satisfied from \eqref{eq:sample:riots:4:c}. We then appeal to the conclusions \eqref{eq:inverse:div}--\eqref{eq:inverse:div:error:1} and \eqref{eq:inverse:div:error:stress}--\eqref{eq:inverse:div:error:stress:bound} to conclude as follows. First, we set
$$ \phi_\vp^{m,l} = \divH (D_t \sigma^m_\vp) \, , \qquad \phi_\vp^{m,*} = \divR (D_t \sigma_\vp^m) \, . $$
From \eqref{eq:inverse:div:stress:1}, we obtain both \eqref{est.S.phi.1:c} and \eqref{est.S.phi.infty:c}. The pointwise bound in \eqref{est.S.by.pr.final3:c} holds due to \eqref{eq:inv:div:extra:conc}, \eqref{eq:inverse:div:sub:1}, and \eqref{eq:divH:formula}. Next, we obtain \eqref{est.S.by.pr.final4:c} from \eqref{eq:inverse:div:error:stress:bound}. Finally, we obtain \eqref{est.S.pr.p.support:2:c} from \eqref{eq:inverse:div:linear} and \eqref{eq:LP:div:support}, concluding the proof for $1\leq m \leq \bar m$.
\smallskip

\noindent\texttt{Step 2c: Highest shell.} \, For the case with the highest shell, corresponding to the projector $\mathbb{P}^*$ from \eqref{eq:decomp:showing:redux:2:c}, we appeal to Lemma~\ref{lem:special:cases} with $q=\sfrac 32$, $\lambda=\Lambda\Gamma$, $\rho = \mathbb{P}_{\neq 0} \pr(\rho^{\beta(j)})$. Specifically, to verify the assumptions in Part 2 of Proposition~\ref{prop:intermittent:inverse:div}, we set for $p=\sfrac 32,\infty$
\begin{align*}
&\ov{\varrho} = \mathbb{P}^* \mathbb{P}_{\neq 0} \pr(\rho^{\beta(j)})  \, , \quad \ov\vartheta = \vartheta \,  \textnormal{ as defined in \eqref{eq:remainder:inverse}} \, , \\
&\ov{\const}_{*,p} = \left( \frac{\Lambda\Gamma}{\mu_{\bar m}} \right)^{N_{**}} \left(\const_{*,1} \Upsilon^{-2}\Upsilon'\right)^{\sfrac 23} (\lambda\Gamma)^3 \, , \quad  \ov \Upsilon =  \ov \Upsilon' = \mu \, , \quad \ov \Lambda = \Gamma\Lambda \, , \\
&\ov \dpot = {0} \, ,\quad \ov N_* = N_* - \NcutLarge - N_{**} - 3 \, . 
\end{align*}
We note that we have altered the definition of $N_*$ compared to the previous two substeps for convenience. But from \eqref{i:st:sample:wut:wut:wut:c}, we have in fact made it \emph{smaller}, so that the low-frequency assumptions from the inverse divergence are still satisfied. Then \eqref{eq:moll:1:as} is satisfied exactly as in the first substep, and so from \eqref{eq:remainder:inverse} we have that \eqref{eq:DN:Mikado:density} is satisfied. We use \eqref{eq:sample:prop:decoup:c}, \eqref{eq:sample:prop:par:00:c}, \eqref{eq:sample:prop:parameters:0:c}, the altered choice of $\ov N_*$ above, \eqref{eq:sample:1:conc:2:c} and \eqref{eq:sample:1:conc:2:c:infty}, and \eqref{eq:sample:prop:Ncut:3:c} to see that \eqref{eq:inverse:div:parameters:0}--\eqref{eq:inverse:div:parameters:1} are satisfied. Continuing onto the nonlocal assumptions from Proposition~\ref{prop:intermittent:inverse:div}, we have that \eqref{eq:inv:div:wut}--\eqref{eq:inverse:div:v:global:parameters} are satisfied as in the previous substep, and \eqref{eq:riots:4} is satisfied from \eqref{eq:sample:riot:4:4:c}. We then appeal to the conclusions \eqref{eq:inverse:div}--\eqref{eq:inverse:div:error:1} and \eqref{eq:inverse:div:error:stress}--\eqref{eq:inverse:div:error:stress:bound} to conclude as follows. First, we set
$$ \phi_\vp^* = (\divH+\divR) (D_t \sigma_\vp^*) \, . $$
We may ignore \eqref{eq:inverse:div:stress:1} since $\dpot=0$. Then the only conclusion we require is \eqref{est.S.by.pr.final.star:c}, which follows from \eqref{eq:inverse:div:error:stress:bound}.

\medskip

\noindent\texttt{Step 3: Verification of \eqref{est:mean.Dtsiph}}. 
The proof is similar to \eqref{est:mean.DtsiS}. Indeed, we have
\begin{align*}
r_\phi^{\sfrac23}\bigg|\int_{\T^3}&D_t^{M+1}\pr\left(H^{\alpha{(j)}}\right) \left(\mathbb{P}_{\neq 0} \pr(\rho^{\beta{(j)}})\right)\circ \Phi dx\bigg| \\   
&\lec r_\phi^{\sfrac23}\norm{\Delta^{\floor{\frac{\dpot}4}}\left(D_t^{M+1}\pr\left(H^{\alpha{(j)}}\right)\circ \Phi^{-1}\right)}_{\sfrac32}
     \norm{\Delta^{-\floor{\frac{\dpot}4}}\left(\mathbb{P}_{\neq 0} \pr(\rho^{\beta{(j)}})\right)}_{\sfrac32}\\
     &\lec r_\phi^{\sfrac23} \const_{G,1}^{\sfrac23} (\max(\la, \la')\Ga)^{\halfd}
     \mu^{-\halfd}  (\const_{*,1}\Upsilon^{-2}\Upsilon')^{\sfrac23} 
     \MM{M+1,M_t - \NcutSmall, \nu\Ga, \nu'\Ga}
     \\
     &\lec (\La\Ga)^{-K_\circ}(\Upsilon^{-2}\Upsilon' )^{\sfrac23}
     {(\max(\la, \la')\Ga)^{\lfloor\sfrac{\dpot}4\rfloor}
     \mu^{-\lfloor\sfrac{\dpot}4\rfloor} }
     \MM{M,M_t - \NcutSmall-1, \nu\Ga, \nu'\Ga}\, . 
\end{align*}
using Lemma \ref{lem:decoup}, \eqref{eq:sample:1:conc:2:c}, \eqref{eq:DDpsi:sample}, \eqref{eq:sample:riots:4:c} with standard Littlewood-Paley theory. 
Then, recalling $\frac{d^M}{dt^M} \langle D_t \si_\ph \rangle = \langle D_t^{M+1} \si_\ph \rangle$ and using the representation \eqref{rep.Dtsiph} of $D_t\si_\ph$, we obtain \eqref{est:mean.Dtsiph}. 
\end{proof}

\section{Error estimates}\label{sec:ER:main}
In this section, we will define and estimate a number of error terms, as well as the pressure increments and pressure current errors. Such estimates will require repeated application of the inverse divergence operator from Proposition~\ref{prop:intermittent:inverse:div}, and the pressure creation and pressure current error estimates from section~\ref{opsection:pressure}. First, in subsection~\ref{sec:error.ER}, we add $\hat w_\qbn$ to the Euler-Reynolds system and identify the remaining error terms.  These include the oscillation stress error, the transport and Nash stress errors, the divergence corrector errors, and the mollification error.  We estimate these error terms and define and estimate the related pressure increments and current errors in subsections~\ref{ss:ssO}, \ref{ss:ER:TN}, \ref{sss:dce}, and \ref{ss:essemm}, respectively. The reader who is only interested in the proof of Theorem~\ref{thm:main:ER} following the strategy outlined in Remark~\ref{rem:toozday} can ignore all the results from these sections labeled with an asterisk.  The reader who is interested in the proof of Theorem~\ref{thm:main:ER} following the strategy of Proposition~\ref{prop:main:ER}, i.e. a strategy which includes the construction of $\pi_q^q$, should read the asterisked lemmas with the subtitle ``pressure increment'' but can skip the lemmas with the subtitle ``pressure current,'' as these estimate the current errors generated by new pressure increments. Then in subsection~\ref{dodge4verification}, we upgrade material derivatives and check Hypothesis~\ref{hyp:dodging4}, while in subsection~\ref{sec:new.pressure.stress}, we collect all the pressure increments and pressure current errors created so far in this section. Then in subsection~\ref{op:tnce}, we estimate a number of error terms, known as the transport-Nash current errors, which are related to the Reynolds stress errors and which will appear in the relaxed local energy inequality; we refer to~\cite[subsection~5.1]{GKN23} for a full derivation.  Since many of these error terms require precise knowledge of the structure of the Reynolds stress, we include the estimates in this section. Finally, subsection~\ref{op:moll:ss} contains estimates for mollification errors which appear in the relaxed local energy inequality.  

\subsection{Defining new Euler-Reynolds error terms}\label{sec:error.ER}


We define $S_{q+1}$ by adding $\hat w_{q+\bn}$ to the Euler-Reynolds system for $(u_q, p_q, R_q, -\pi_q)$ in \eqref{eqn:ER:LEI:new} (recall also \eqref{eqn:ER}) and collecting various error terms, which we shall show are well-defined in the remainder of this section.
\begin{align}
\div(S_{q+1}) 
&= \pa_t \hat w_{q+\bn} + (u_{q}\cdot \na )\hat w_{q+\bn} + (\hat w_{q+\bn}\cdot \na) u_q + \div(\hat w_{q+\bn}\otimes \hat w_{q+\bn} + R_\ell-\pi_\ell \Id)  \notag\\
&\qquad + \div \left( R_q^q - R_\ell+\left( \pi_\ell - \pi_q^q \right) \Id \right)  \notag \\
&= \underbrace{(\pa_t + \hat u_q \cdot \na) w_{q+1} + w_{q+1}\cdot\nabla \hat u_q}_{=:\, \div S_{TN}}
+\underbrace{\div \left(\wp_{q+1}\otimes \wp_{q+1} + R_\ell-\pi_\ell\Id\right)}_{=:\, \div S_O} \notag \\
&\qquad + \underbrace{\div \left(\wp_{q+1}\otimes_s \wc_{q+1}+\wc_{q+1}\otimes \wc_{q+1}\right)
}_{=:\, \div S_C} 
+ \underbrace{\div \left( R_q^q - R_\ell + \left( \pi_\ell - \pi_q^q \right) \Id \right) }_{=:\div S_{M1}}
\label{ER:new:error}\\
&\qquad + 
\underbrace{(\pa_t + \hat u_{q}\cdot \na )(\hat w_{q+\bn}-w_{q+1}) 
+ ((\hat w_{q+\bn}-w_{q+1})\cdot \na) \hat u_q + \div(\hat w_{q+\bn}\otimes \hat w_{q+\bn} - w_{q+1}\otimes w_{q+1} )}_{=:\div S_{M2}}
\, . \notag
\end{align}
In the second equality,\index{$S_{q+1}$}\index{$\ov R_{q+1}$} we used \eqref{eq:dodging:oldies} to exchange $u_q$ and $\hat u_q$. (Recall also \eqref{eq:hat:no:hat}.) We note that the symmetric stresses $S_O$ and $S_C$ are not simply the quantities inside parentheses and take some care to construct; see subsections~\ref{ss:ssO}, \ref{sss:dce}.  
Also, we note that $\pa_t w_{q+1} + (\hat u_q\cdot\na)w_{q+1} + w_{q+1}\cdot \na \hat u_q$ has mean-zero, so that it can be written in divergence form $\div S_{TN}$; see subsection~\ref{ss:ER:TN}. This is because the second and third terms can be written in divergence form, and $w_{q+1}$ is given by the curl of a vector-valued function (see \eqref{wqplusoneonephi} and \eqref{wqplusoneone}.) The same reasoning works for the terms in $\div S_{M2}$. 

With the above definitions, we set\index{$\ov R_{q+1}$}\index{$S_{q+1}$}
\begin{align}\label{eq:ovR:def}
    \ov R_{q+1} := R_q - R_q^q + S_{q+1} \, .
\end{align}
We have notated the error with an overline as $\ov R_{q+1}$ in order to be consistent with the notation from \cite{GKN23}, where the stress error $\ov R_{q+1}$ will be adjusted slightly in \cite[Section~7]{GKN23} in order to produce the final Reynolds stress $R_{q+1}$ needed to complete the proof of Theorem~\ref{thm:main}. We can now see that $(u_{q+1}, p_q, \ov R_{q+1}, -(\pi_q-\pi_q^q))$ solves the Euler-Reynolds system (recall from \eqref{def.w.mollified} that $u_{q+1} = u_q + \hat w_{q+\bn}$)
\begin{align}
    \partial_t u_{q+1} 
    + \div \left( u_{q+1} \otimes u_{q+1} \right) + \nabla p_q    
    = \div (  -(\pi_q-\pi_q^q)\Id + \ov R_{q+1}) \, , \qquad \div \, u_{q+1} = 0 \, . \label{ER:new:equation}
\end{align}
We will show in the remainder of this section that the new stress error $S_{q+1}$ can be decomposed into components ${S}_{q+1}^k$ as
\begin{align*}
    S_{q+1} &=\sum_{k=q+1}^{q+\bn} {S}_{q+1}^k \, . 
\end{align*}

\subsection{Oscillation stress error \texorpdfstring{$S_O$}{essoh}}\label{ss:ssO}
In order to define and analyze $S_O$, cf. \eqref{ER:new:error}, we first consider
\begin{align}\label{eqn:wpwp}
    \div\left(\wp_{q+1}\otimes \wp_{q+1}\right)^\bullet &= \sum_{\xi,i,j,k,\vecl,\diamond} \partial_\alpha \left(  a_{(\xi),\diamond} (\nabla\Phi_{(i,k)}^{-1})_{\theta}^\alpha \BB_{(\xi),\diamond}^\theta(\Phi_{(i,k)}) \,  a_{(\xi),\diamond} (\nabla\Phi_{(i,k)}^{-1})_\gamma^\bullet \BB_{(\xi),\diamond}^\gamma(\Phi_{(i,k)}) \right) \, ,
\end{align}
where $\bullet$ denotes the unspecified components of a vector field and we have used \eqref{eq:dodging:newbies} from Lemma~\ref{lem:dodging} to eliminate all cross terms. Recalling from \eqref{eq:W:xi:q+1:phi:def} and \eqref{eq:W:xi:q+1:nn:def} that
$\BB_{(\xi), \diamond}=\chib_{(\xi)}^\diamond \sum_I \etab_\xi^{I,\diamond} \WW_{(\xi),\diamond}^I$, that the $\WW_{(\xi),\diamond}^I$'s are identical up to a shift, and the notational convention for $\rhob_{(\xi)}^\diamond$ from Remark~\ref{rem:notational:conventions}, we decompose
\begin{align*}
(\BB\otimes\BB)_{(\xi),\diamond}
&= \left(\chib_{(\xi)}^{\diamond}\right)^2 \sum_{I} (\etab_{\xi}^{I,\diamond})^{2} \mathbb{P}_{\neq 0}(\WW_{(\xi),\diamond}^I\otimes  \WW_{(\xi), \diamond}^I)
+  \left(\chib_{(\xi)}^{\diamond}\right)^2 \mathbb{P}_{\neq 0} \left(\sum_I (\etab_{\xi}^{I,\diamond})^{2} \right)
\left\langle \WW_{\pxi,\diamond}^I\otimes  \WW_{\pxi,\diamond}^I \right\rangle
\notag\\
&\quad + \mathbb{P}_{\neq 0} \left(\chib_{(\xi)}^{\diamond}\right)^2 \left\langle\sum_I (\etab_{\xi}^{I,\diamond})^{2} \right\rangle
\left\langle \WW_{\pxi,\diamond}^I\otimes  \WW_{\pxi, \diamond}^I \right\rangle 
+ 
\left\langle \left(\chib_{(\xi)}^{\diamond}\right)^2 \right\rangle \left\langle\sum_I (\etab_{\xi}^{I,\diamond})^{2}\right\rangle
\left\langle \WW_{\pxi,\diamond}^I\otimes  \WW_{\pxi,\diamond}^I \right\rangle \, .
\end{align*}
In particular, using \eqref{i:bundling:3} and the definitions of $\rhob_{(\xi)}^\diamond$ and $\ov \rhob_{(\xi)}:=\ov \rhob_{\xi,k}$ from Proposition~\ref{prop:bundling}, \eqref{item:pipe:4} from Proposition~\ref{prop:pipeconstruction}, \eqref{item:pipe:5:current} from Proposition~\ref{prop:pipe.flow.current}, Definition~\ref{defn:pipe.bundle}, and
\eqref{eq:sa:summability}, we obtain that
\begin{subequations}
\begin{align}
(\BB\otimes\BB)_{(\xi),R}
&= \left(\ov\chib_{(\xi)}\right)^6 \sum_{I}(\etab_{\xi}^I)^{6}\mathbb{P}_{\neq 0}(\WW_{(\xi),R}^I\otimes  \WW_{(\xi),R}^I)
+ \left(\left(\ov\chib_{(\xi)}\right)^6 - 1\right)\xi\otimes \xi 
+ \xi\otimes \xi \, ,
\label{eq:BB.decomp}\\
(\BB\otimes\BB)_{(\xi),\ph}
&
=\left(\ov\chib_{(\xi)}\right)^4 \sum_{I}(\etab_{\xi}^I)^{4} \mathbb{P}_{\neq 0}(\WW_{(\xi),\varphi}^I\otimes  \WW_{(\xi),\varphi}^I)
+ c_0 \left(\ov\chib_{(\xi)}\right)^4 r_q^\frac23 \xi\otimes \xi \,  \mathbb{P}_{\neq 0} \left(\sum_I (\etab_{\xi}^{I})^4 \right)
\notag\\
&\qquad + c_0 c_1 \mathbb{P}_{\neq 0} \left(\left(\ov\chib_{(\xi)}\right)^4 \right) r_q^\frac23
 \xi\otimes \xi 
+ c_0 c_1 c_2 r_q^\frac23 \Gamma_q^{-2} \xi\otimes \xi \, ,
\label{eq:BB.ph.decomp}
\end{align}
\end{subequations}
for dimensional constants $c_0$, $c_1$, and $c_2$ which are bounded independently of $q$ and depend only on the dimensional constants in \eqref{e:fat:pipe:estimates:1} and \eqref{e:pipe:estimates:1:current} and the mean of $\sum_I (\zetab_\xi^I)^4$. Since each vector field used to define the simple symmetric tensors in \eqref{eq:BB.decomp} and \eqref{eq:BB.ph.decomp} does not vary in the $\xi$-direction (see, \eqref{eq:derivative:along:pipe}, \eqref{i:bundling:1}, and Definition~\ref{def:etab}), each simple symmetric tensor satisfies $\xi\cdot\nabla (\BB\otimes \BB)_{\pxi,\diamond}=0$. Then using that each vector field in \eqref{eq:BB.decomp} and \eqref{eq:BB.ph.decomp} has been composed with $\Phiik$ and the identity $\partial_\alpha \left( (\nabla\Phiik^{-1})_\theta^\alpha (\BB\otimes \BB)_{\pxi,\diamond}\circ \Phiik \xi^\theta\right) = \xi^\theta (\partial_\theta (\BB\otimes \BB)_{\pxi,\diamond})\circ \Phiik =0$, we have that \eqref{eqn:wpwp} can be expanded as
\begin{subequations}
\label{eqn:osc.expand}
\begin{align}
 \div\left(\wp_{q+1}\otimes \wp_{q+1}\right)^\bullet
&=
\sum_{\xi,i,j,k,\vecl} \partial_\alpha \left(  a_{(\xi),R}^2 
(\nabla\Phi_{(i,k)}^{-1})_{\theta}^\alpha 
(\nabla\Phi_{(i,k)}^{-1})_\gamma^\bullet
(\xi^\theta \xi^\gamma) \right)
\label{eqn:osc.exp1}\\
&\quad+\sum_{\xi,i,j,k,\vecl} \partial_\alpha \left(  a_{(\xi),\ph}^2 
(\nabla\Phi_{(i,k)}^{-1})_{\theta}^\alpha 
(\nabla\Phi_{(i,k)}^{-1})_\gamma^\bullet
{c_0c_1c_2} \Gamma_q^{-2} r_q^\frac23
(\xi^\theta \xi^\gamma)
\right)
\label{eqn:osc.exp1.ph}
\\
&\quad+
\sum_{\xi,i,j,k,\vecl} B_{(\xi),R}^{\bullet}
\left(\mathbb{P}_{\neq 0} \ov\chib_{(\xi)}^6 \right) \circ \Phi_{(i,k)}
\label{eqn:osc.exp2}\\
&\quad+
\sum_{\xi,i,j,k,\vecl} B_{(\xi),\ph}^{\bullet}
\left( \mathbb{P}_{\neq 0} \ov\chib_{\xi}^4 \right) \circ  (\Phi_{(i,k)}) {c_0 c_1} r_q^\frac23
\label{eqn:osc.exp2.ph}\\
&\quad + {c_0}
\sum_{\xi,i,j,k,\vecl} B_{(\xi),\ph}^{\bullet}
 r_q^\frac23 \left(\ov\chib_{(\xi)}^4\mathbb{P}_{\neq 0} 
\sum_I (\etab_\xi^I)^4 \right) {\circ \Phiik}
\label{eqn:osc.exp.2.5}\\
&\quad+
\sum_{\xi,i,j,k,\vecl,\diamond} B_{(\xi),\diamond}^{\bullet}
\left( \left(\chib_{(\xi)}^\diamond\right)^2 \sum_{I} (\etab_{\xi}^{I,\diamond})^{2} \mathbb{P}_{\neq 0}(\varrho_{(\xi),\diamond}^I)^2\right) \circ \Phi_{(i,k)}
\label{eqn:osc.exp3}
\end{align}
\end{subequations}
where for convenience we set
\begin{equation}\label{eq:rhoxidiamond:def}
B_{(\xi),\diamond}^{\bullet} := \xi^\theta \xi^\gamma\partial_\alpha \left(  a_{(\xi),\diamond}^2 
(\nabla\Phi_{(i,k)}^{-1})_{\theta}^\alpha 
(\nabla\Phi_{(i,k)}^{-1})_\gamma^\bullet\right)\, , \qquad \varrho_{(\xi),\diamond}^{I}:=\xi \cdot \WW_{(\xi),\diamond}^I \, .
\end{equation}

The first and second terms above in \eqref{eqn:osc.exp1} and \eqref{eqn:osc.exp1.ph} cancel out $-R_\ell {+\pi_\ell\Id}$ from \eqref{ER:new:error} as follows:
\begin{align}
    &\sum_{\xi,i,j,k,\vecl} a_{(\xi),R}^2 \nabla\Phi_{(i,k)}^{-1} \left(\xi\otimes \xi\right) \nabla\Phi_{(i,k)}^{-\top}\notag\\
     &\quad \underset{\eqref{eq:a:xi:def}}{=} \sum_{\xi,i,j,k,\vecl} \delta_{q+\bn}\Gamma_{q}^{2j-2} \psi_{i,q}^6 \omega_{j,q}^6 \chi_{i,k,q}^6 \zeta_{q,R,i,k,\xi,\vecl}^2 \gamma_{\xi,\Ga_q^9}^2\left(\frac{R_{q,i,k}}{\delta_{q+\bn}\Gamma_{q}^{2j-2}}\right)\nabla\Phi_{(i,k)}^{-1}\left(\xi\otimes\xi\right)\nabla\Phi_{(i,k)}^{-\top} \notag\\
    &\quad \underset{\eqref{e:split}, \eqref{eq:rqnpj}, \eqref{eq:summy:summ:1}}{=}  -\sum_{i,j,k} \psi_{i,q}^6 \omega_{j,q}^6 \chi_{i,k,q}^6
   \bigg{(} {R}_\ell -\pi_\ell\Id \notag\\
   &\qquad
    + \sum_{\substack{\xi',i',j' \\ k',l'}}  \frac{\delta_{q+\bn} \Gamma_q^{2j'-2} {C\Gamma_q^{-2}} }{\left| \nabla\Phi^{-1}_{(i',k')} \xi'\right|^{\sfrac 43}} \psi_{i',q}^4 \omega_{j',q}^4 \chi_{i',k',q}^4 \mathcal{X}_{q,\xi',l'}^4 \circ \Phi_{i',k',q} \tilde \gamma_{\xi'}^2 \nabla\Phi_{(i',k')}^{-1}\xi'\otimes \xi'\left(\nabla\Phi_{(i',k')}^{-T}\right)
    \bigg{)} \notag\\
    &\quad\underset{\eqref{eq:inductive:partition}, \eqref{eq:omega:cut:partition:unity}, \eqref{eq:chi:cut:partition:unity}}{=} \pi_\ell\Id -{R}_\ell \notag\\
    &\qquad - \sum_{\substack{\xi',i',j' \\ k',l'}}  \frac{\delta_{q+\bn} \Gamma_q^{2j'-2} {C\Gamma_q^{-2}} }{\left| \nabla\Phi^{-1}_{(i',k')} \xi'\right|^{\sfrac 43}} \psi_{i',q}^4 \omega_{j',q}^4 \chi_{i',k',q}^4 \mathcal{X}_{q,\xi',l'}^4 \circ \Phi_{i',k',q} \tilde \gamma_{\xi'}^2 \nabla\Phi_{(i',k')}^{-1}\xi'\otimes \xi'\left(\nabla\Phi_{(i',k')}^{-T}\right) \notag\\
    &\quad\underset{\eqref{eq:summy:summ:2}}{=} \pi_\ell\Id -{R}_\ell 
    - \sum_{\substack{\xi',i',j' \\ k',\vecl'}}  \frac{\delta_{q+\bn} \Gamma_q^{2j'-2} {C\Gamma_q^{-2}} }{\left| \nabla\Phi^{-1}_{(i',k')} \xi'\right|^{\sfrac 43}} \psi_{i',q}^4 \omega_{j',q}^4 \chi_{i',k',q}^4 \zeta_{q,\varphi,i',k',\xi',\vecl'}^2 \tilde \gamma_{\xi'}^2 \nabla\Phi_{(i',k')}^{-1}\xi'\otimes \xi'\left(\nabla\Phi_{(i',k')}^{-T}\right) \notag\\
    &\quad\underset{\eqref{eq:a:xi:phi:def}}{=} \pi_\ell\Id -{R}_\ell 
    - \underbrace{\sum_{\substack{\xi',i',j' \\ k',\vecl'}}  a_{\pxi,\varphi}^2 c_0 c_1 c_2 \Ga_q^{-2} r_q^{\sfrac 23} \nabla\Phi_{(i',k')}^{-1}\xi'\otimes \xi'\left(\nabla\Phi_{(i',k')}^{-T}\right)}_{=\eqref{eqn:osc.exp1.ph}}
    \,.\label{eq:cancellation:plus:pressure:nn}
\end{align}
The inverse divergence of the remaining terms \eqref{eqn:osc.exp2}-\eqref{eqn:osc.exp3} will therefore form the oscillation stress errors. 

\begin{lemma}[\bf Applying inverse divergence]\label{lem:oscillation:general:estimate}
There exist symmetric stresses $S_O^m$ for $m=1,\dots,q+\bn$ such that the following hold.\index{$S_O^m$}
\begin{enumerate}[(i)]
    \item\label{item:oscee:1} $\displaystyle \div \left( w_{q+1}^{(p)} \otimes w_{q+1}^{(p)} + R_\ell - \pi_\ell \Id \right) = \sum_{m=q+1}^{q+\bn} \div S_O^m$, where $S_O^m$ can be split into local and non-local errors as $S_O^m = S_O^{m,l}+S_O^{m,*}$.
    \item For $m=q+1, \dots, q+\bn$ and $N,M\leq \sfrac{\Nfin}{10}$, the local parts ${S}^{m,l}_{O}$ satisfy
\begin{subequations}
\begin{align}
    \left\| \psi_{i,q} D^N D_{t, q}^M {S}^{m,l}_{O} \right\|_{\sfrac32}
    &\lesssim \Gamma_{m}^{-9}  \delta_{m+\bn} \lambda_{m}^N \MM{M, \Nindt, \tau_{q}^{-1}\Gamma_{q}^{i+14}, \Tau_{q}^{-1}\Ga_q^9}\label{eq:Onpnp:estimate:2} \\
\left\|\psi_{i,q} D^N D_{t, q}^M {S}^{m,l}_{O} \right\|_{\infty} 
    &\lesssim \Gamma_{m}^{\badshaq-9} \lambda_{m}^N
    \MM{M, \Nindt, \tau_{q}^{-1}\Gamma_{q}^{i+14}, \Tau_{q}^{-1}\Ga_q^9} \, . \label{eq:Onpnp:estimate:2:new}
\end{align}   \end{subequations}
When $m=q+2,\dots,q+\bn$ and $q+1\leq q' \leq m-1$, the local parts satisfy
\begin{align}
B\left( \supp \hat w_{q'}, \lambda_{q'}^{-1} \Gamma_{q'+1} \right) \cap \supp S_O^{m,l} = \emptyset \, . \label{osc:support:first}
\end{align}
\item For $m=q+1, \dots, q+\bn$ and $N,M\leq 2\Nind$, the non-local parts ${S}^{m,*}_{O}$ satisfy
\begin{align}
\left\| D^N D_{t, q}^M {S}^{m,*}_{O} \right\|_{L^\infty}
    &\leq { \Tau_{q+\bn}^{4\Nindt}}\delta_{q+3\bn}\la_{m}^{N}\tau_{q}^{-M} \, .
    \label{eq:Onpnp:estimate:1}
\end{align}
\end{enumerate}
\end{lemma}
\begin{remark*}[\bf Abstract formulation of the oscillation stress error]\label{rem.ct.osc}
For the purposes of analyzing the transport and Nash current errors subsection~\ref{op:tnce} and streamlining the creation of pressure increments, it will be useful to abstract the properties of these error terms. First, there exists a $q$-independent constant $\const_{\mathcal{H}}$ such that
\begin{subequations}\label{osc.loc.pot}
\begin{align}
    S^{m,l}_O &= \sum_{i,j,k,\xi,\vecl,\diamond} \sum_{j'=0}^{\const_{\mathcal{H}}} H_{i,j,k,\xi,\vecl,\diamond}^{\alpha{(j')}} \rho_{i,j,k,\xi,\vecl,\diamond}^{\beta{(j')}} \circ \Phi_{(i,k)} \qquad \textnormal{ if } \qquad m=q+1, \, q+\half \, , \label{osc.loc.pot.1} \\
    S^{m,l}_O &= \sum_{i,j,k,\xi,\vecl,I,\diamond} \sum_{j'=0}^{\const_{\mathcal{H}}} H_{i,j,k,\xi,\vecl,I,\diamond}^{\alpha{(j')}} \rho_{i,j,k,\xi,\vecl,I,\diamond}^{\beta{(j')}} \circ \Phi_{(i,k)} \qquad \textnormal{ if } q+\half+1\leq m \leq \qbn \, . \label{osc.loc.pot.2}
\end{align}
\end{subequations}
For the remaining values of $m$, $S_{O}^{m,l}$ is zero. These equalities will be proven in the course of proving Lemma~\ref{lem:oscillation:general:estimate}, \ref{lem:osc:no:pressure}, and \ref{lem:oscillation:pressure}. The pointwise estimate \eqref{eq:H:eckel:osc} will be proved in Lemma~\ref{lem:osc:no:pressure} and \ref{lem:oscillation:pressure}, and the rest of the claims in this remark will be proved in Lemma~\ref{lem:oscillation:general:estimate}. Note that the proof of \eqref{rho:eckel:estimate} will also require Remark \ref{rem:potential.rho}.

Next, the functions $H$ and $\rho$ (with subscripts and superscripts suppressed for convenience) defined above satisfy the following.
\begin{enumerate}[(i)]
    \item For all $N, M \leq \sfrac{\Nfin}{10}$, 
    \begin{align}\label{eq:H:eckel:osc}
        \left| D^N \Dtq^M H \right| &\lec \pi_\ell \Ga_q^{100} \Lambda_q \bar\lambda^N \MM{M, \Nindt, \tau_{q}^{-1}\Gamma_{q}^{i+13}, \Tau_{q}^{-1}\Ga_q^8} \, ,
    \end{align}
    where $\bar \la = \la_{q+1}\Ga_q^{-5}$ for $m=q+1, q+\half$ while $\bar \la = \la_{q+\half}$ for $m \geq q+\half+1$. 
    \item We have that
    \begin{subequations}\label{H:eckel:support:osc}
    \begin{align}
        \supp H &\subseteq \supp \eta_{i,j,k,\xi,\vecl,\diamond}  \quad \textnormal{ if } \quad m=q+1, q+\half \label{H:eckel:support:osc:1} \\
        \supp H &\subseteq \supp \eta_{i,j,k,\xi,\vecl,\diamond} \zetab_\xi^{I,\diamond} \quad \textnormal{ if } \quad q+\half+1\leq m \leq \qbn \label{H:eckel:support:osc:2}
    \end{align}
    \end{subequations}
    \item For $\dpot$ as in \eqref{i:par:10}, there exist a tensor potential $\vartheta$ (we suppress the indices at the moment for convenience) such that $\rho=\partial_{i_1 \dots i_{\dpot}}\vartheta^{(i_1,\dots,i_\dpot)}$.  Furthermore, $\vartheta$ is $(\T / \la_{q+1}\Ga_q^{-4} )^3$-periodic in the case $m=q+1$, $(\T / \la_{q+\half} )^3$-periodic in the case $m=q+\half$, and $(\T / \la_{q+\half}\Ga_q )^3$-periodic in the remaining cases. Finally, $\vartheta$ satisfies the estimates
    \begin{subequations}\label{rho:eckel:estimate}
    \begin{align}
        \norm{D^N \partial_{i_1}\dots \partial_{i_k} \vartheta^{(i_1,\dots, i_\dpot)}}_{L^p} &\les \Ga_q^{12} (\la_{q+1}\Ga_q^{-4})^{k-\dpot-1} \MM{N, \dpot - k , \la_{q+1}\Ga_q^{-4}, \la_{q+1}\Ga_q^{-1}} \quad \textnormal{if} \quad m=q+1 \label{rho:eckel:estimate:1} \\
        \norm{D^N \partial_{i_1}\dots \partial_{i_k} \vartheta^{(i_1,\dots, i_\dpot)}}_{L^p} &\les \Ga_q^5 \la_{q+\half}^{k-\dpot-1} \la_{q+\half}^N \quad \textnormal{if} \quad m=q+\half \label{rho:eckel:estimate:2} \\
        \norm{D^N \partial_{i_1}\dots \partial_{i_k} \vartheta^{(i_1,\dots, i_\dpot)}}_{L^p} &\les {\left(\frac{\la_{q+\half+1}}{\la_{q+\bn}r_q}\right)^{2-\sfrac2p}} \Ga_q^2 \la_{q+\half}^{-1} (\la_{q+\half}\Ga_q)^{k-\dpot} \notag\\
        &\qquad \times \MM{N, \dpot - k ,\la_{q+\half}\Ga_q, \la_{q+\half+1}} \quad \textnormal{if} \quad  m=q+\half+1 \label{rho:eckel:estimate:3}\\
        \norm{D^N \partial_{i_1}\dots \partial_{i_k} \vartheta^{(i_1,\dots, i_\dpot)}}_{L^p} &\les {\left(\frac{\min(\la_m,\la_\qbn)}{\la_{q+\bn}r_q}\right)^{2-\sfrac2p}} \Ga_q^2 \la_{m-1}^{-2}\la_m \la_{m-1}^{k-\dpot} \la_m^{N} \notag\\
        &\qquad \qquad \textnormal{if} \quad  q+\half+2 \leq m\leq \qbn \label{rho:eckel:estimate:4}
    \end{align}
    \end{subequations}
      for $p=\sfrac 32,\infty$, all $N\leq \sfrac{\Nfin}{5}$, and $0\leq k \leq \dpot$. 
    \item In the cases $m=q+1,q+\half,q+\half+1$, we claim no special support properties for the potential $\vartheta$. In the cases $q+\half+2\leq m \leq \qbn$, we have that
    \begin{align}\label{rho:eckel:support}
        \supp \left( H \rho\circ \Phi \right) \cap B\left( \supp \hat w_{q'}, \la_{q'}^{-1} \Ga_{q'+1} \right) = \emptyset
    \end{align}
    for all $q+1 \leq q' \leq m-1$ (where $m$ refers to the index in $S_O^{m,l}$ from \eqref{osc.loc.pot.1}). 
\end{enumerate}
\end{remark*}
\begin{proof}[Proof of Lemma~\ref{lem:oscillation:general:estimate}]
To define $S_O$, we recall the synthetic Littlewood-Paley decomposition (cf. Section \ref{sec:LP}). Indeed, since $\varrho_{\pxi,\diamond}^I$ depends only on the variables in the plane $\xi^\perp$ from \eqref{eq:derivative:along:pipe} and is periodized to scale $\left(\lambda_{q+\bn}r_{q}\right)^{-1} {=(\la_{q+\half}\Ga_q)^{-1}}$, we can decompose $\mathbb{P}_{\neq 0}$ in front of $(\varrho_{\pxi,\diamond}^{I})^2$ in \eqref{eqn:osc.exp3} into\index{$\tP_{\la_m}$}\index{$\tP_{(m-1,m]}$}
\begin{align}
\mathbb{P}_{\neq 0}
   = \tP_{\la_{q+\half{+1}}}^\xi\mathbb{P}_{\neq 0} + {\sum_{m={q+\half+{2}}}^{q+\bn+1} } \tP_{(\lambda_{m-1}, \lambda_m]}^\xi + (\Id-\tP_{\lambda_{q+\bn+1}}^\xi)  \notag  \\
    =: \tP_{{q+\half{+1}}}^\xi + {\sum_{m={q+\half+{2}}}^{q+\bn+1} } \tP_{({m-1}, m]}^\xi + (\Id-\tP_{{q+\bn+1}}^\xi) \, .
\end{align}
Assuming we can apply the inverse divergence from Proposition \ref{prop:intermittent:inverse:div}, we define
\begin{subequations}
\begin{align}
     {S}^{q+1}_{O}
    &:=(\divH + \divR)\left[
    \sum_{\xi,i,j,k,\vecl} 
    B_{(\xi),R}
    \left( \mathbb{P}_{\neq 0} \ov\chib_{\xi}^6 \right) \circ
    \Phi_{(i,k)}
     +\sum_{\xi,i,j,k,\vecl} B_{(\xi),\ph}{c_0c_1} r_q^\frac23 \left( \mathbb{P}_{\neq 0}\ov\chib_{\xi}^4\right)\circ \Phi_{(i,k)} 
    \right] \label{osc.err.low} \\
     {S}^{q+\half}_{O}
    &:=(\divH + \divR)
    \left[
    \sum_{\xi,i,j,k,\vecl} B_{(\xi),\ph} {c_0} r_q^\frac23 \left(\ov\chib_{\xi}^4  \mathbb{P}_{\neq 0} \left(\sum_I (\etab_\xi^I)^4 \right) \right) \circ \Phi_{(i,k)}
    \right]
    \label{osc.err.med1-}
    \\
     {S}^{q+\half+1}_{O}
    &:= (\divH + \divR)\left[
    \sum_{\xi,i,j,k,\vecl, I,\diamond} 
    B_{(\xi),\diamond}
     \left( \left( \chib_{\pxi}^{\diamond} \right)^2 \left(\etab_{\xi}^{I,\diamond}\right)^{2} \tP_{q+\half+1}^\xi \mathbb{P}_{\neq 0} (\varrho_{\pxi,\diamond}^I)^2\right)
    \circ \Phi_{(i,k)} \right]\label{osc.err.med0}\\
     {S}^{m}_{O}
    &:= (\divH + \divR)\left[
    \sum_{\xi,i,j,k,\vecl, I} 
    B_{(\xi),\diamond}
     \left( \left( \chib_{\pxi}^{\diamond} \right)^2 \left(\etab_{\xi}^{I,\diamond}\right)^{2} \tP_{(m-1, m]}^\xi(\varrho_{\pxi, \diamond}^I)^2\right)
    \circ \Phi_{(i,k)} \right]\label{osc.err.med}\\
     {S}^{q+\bn}_{O}
    &:=\sum_{m=q+\bn}^{q+\bn +1} (\divH + \divR)\left[
    \sum_{\xi,i,j,k,\vecl, I,\diamond} 
    B_{(\xi),\diamond} \left( \left( \chib_{\pxi}^{\diamond} \right)^2 \left(\etab_{\xi}^{I,\diamond}\right)^{2} \tP_{(m-1, m]}^\xi(\varrho_{\pxi,\diamond}^I)^2\right) \circ \Phi_{(i,k)} \right]\label{osc.err.high}\\
    &\quad + (\divH + \divR)\left[
    \sum_{\xi,i,j,k,\vecl, I,\diamond} 
    B_{(\xi),\diamond}
     \left( \left( \chib_{\pxi}^{\diamond} \right)^2 \left(\etab_{\xi}^{I,\diamond}\right)^{2}
(\Id-\tP_{q+\bn+1}^\xi) (\varrho_{\pxi, \diamond}^I)^2\right) \circ \Phi_{(i,k)}     \right] \label{osc.err.high2}
\end{align}
\end{subequations}
for $m=q+ \half+2, \cdots, q+\bn -1$. 
For $q+1 \leq m <q+\bn$, we decompose $  {S}^{m}_{O}$ into  the local part $ {S}^{m, l}_{O}$ which involves the operator $\divH$ and the nonlocal part $ {S}^{m,*}_{O}$ containing the remaining terms. In the case of $m=q+\bn$, we set 
\begin{align}
     {S}^{q+\bn,l}_{O}
    &:=\sum_{m=q+\bn}^{q+\bn +1} \divH\left[
    \sum_{\xi,i,j,k,\vecl, I,\diamond} 
    B_{(\xi),\diamond} \left( \left( \chib_{\pxi}^{\diamond} \right)^2 \left(\etab_{\xi}^{I,\diamond}\right)^{2} \tP_{(m-1, m]}^\xi(\varrho_{\pxi,\diamond}^I)^2\right) \circ \Phi_{(i,k)} \right] 
\end{align}
and absorb the $\divR$ terms in \eqref{osc.err.high} and all the terms in \eqref{osc.err.high2} into $ {S}_{O}^{\qbn,*}$. For the undefined $  {S}^{m}_{O}$ corresponding to $m=q+2,\cdots, q+\half-1$, we set them as identically zero. 

The desired estimates will follow from applying Proposition \ref{prop:intermittent:inverse:div}. While many of the parameter choices will vary depending on the case, we fix the following choices throughout the proof:
\begin{subequations}\label{eq:osc:general:choices}
\begin{align}
    &p=\sfrac 32, \infty \, , \quad v= \hat u_q \, , \quad D_t = D_{t,q} \, , \quad {N_* = \sfrac{\Nfin}{4} \, , \quad M_* = \sfrac{\Nfin}{5}} \ , \\
    &\la' = \Lambda_q \, , \quad M_t = \Nindt \, , \quad \nu' = \Tau_q^{-1} \Ga_q^{{8}} \, , \quad \Ndec \textnormal{ as in \eqref{i:par:9}} \, , \\
    &\ {M_\circ = N_\circ = 2\Nind \, , \quad K_\circ \textnormal{ as in \eqref{i:par:9.5}} \, , \quad\const_v = \La_q^{\sfrac12}\, .} 
\end{align}
\end{subequations}

\noindent\texttt{Case 1:} Estimates for \eqref{osc.err.low}.
Fix values of $i,j,k,\xi,\vecl$ and consider the term which includes $B_{(\xi),R}$, where we have abbreviated $B_{(\xi),R}^{\bullet} =B_{(\xi, i, j, k, \vecl),R}^{\bullet}$.
We apply Proposition \ref{prop:intermittent:inverse:div} with the low-frequency choices
\begin{align*}
    &G^\bullet = B_{(\xi),R}^{\bullet}\, , \quad \const_{G,\sfrac 32} =\left| \supp(\eta_{i,j,k,\xi,\vecl,R}^2) \right|^{\sfrac 23} \delta_{q+\bn} \Gamma_{q}^{2j+21} \La_q \, , \quad \const_{G,\infty}=\Gamma_{q}^{\badshaq+30}\La_q \, , \\
    &\la =\la_{q+1}\Gamma_q^{-5} \, , \quad \nu = \tau_q^{-1} \Gamma_q^{i+{13}} \, , \quad \Phi= \Phi_{(i,k)} \, ,
\end{align*}
and the choices from \eqref{eq:osc:general:choices}. We have that \eqref{eq:inv:div:NM} is satisfied by definition.  Next, to check \eqref{eq:inverse:div:DN:G}, we observe that in $B_{\pxi,R}^{\bullet}$, the differential operator on $a^2_{\xi}$ is $\xi^\theta (\nabla\Phi^{-1}_{(i,k)})_\theta^\alpha \partial_\theta$. Therefore $G$ satisfies \eqref{eq:inverse:div:DN:G} for $p=\sfrac 32$ from \eqref{e:a_master_est_p_R} and for $p=\infty$ from the same inequality and \eqref{ineq:jmax:use}. By Corollary  \ref{cor:deformation}, $\Phi_{(i,k)}$ satisfies \eqref{eq:DDpsi2} and \eqref{eq:DDpsi} for $\la' = \Lambda_q$, and by \eqref{eq:nasty:D:vq:old} at level $q$, we have that \eqref{eq:DDv} is satisfied.

To check the high-frequency assumptions, we set
\begin{subequations}
\begin{align}
    &\varrho = \left( \mathbb{P}_{\neq 0} \ov \rhob_\xi^6 \right) \, , \quad \dpot \textnormal{ as in \eqref{i:par:10}} \, , \quad \vartheta=\delta_{i_1 i_2} \delta_{i_3 i_4} \dots \delta_{i_{\dpot-1} i_\dpot} \Delta^{-\sfrac \dpot 2} \varrho \, , \label{im:defining:a:potential} \\
    &\mu = \Upsilon=\Upsilon'= \lambda_{q+1}\Gamma_q^{-4} \, , \quad \ov\Lambda = \lambda_{q+1}\Gamma_q^{-1} \, , \quad \const_{*,p} = \Gamma_q^6 {\la_{q+1}^\alpha} \, ,
\end{align}
\end{subequations}
where $\alpha$ is chosen as in \eqref{eq:choice:of:alpha}. Then from Proposition~\ref{prop:bundling} and standard Littlewood-Paley theory, we have that \eqref{eq:DN:Mikado:density} is satisfied. Next, we have that \eqref{eq:inverse:div:parameters:0} is satisfied by definition and from \eqref{condi.Nfin0}. In addition, we have that \eqref{eq:inverse:div:parameters:1} is satisfied from \eqref{condi.Ndec0}.  In order to check the nonlocal assumptions in Part 4, we first appeal to \eqref{condi.Nfin0}, which gives \eqref{eq:inv:div:wut}.  We have that \eqref{eq:inverse:div:v:global} is satisfied from \eqref{eq:bobby:old}, and \eqref{eq:inverse:div:v:global:parameters} is satisfied from \eqref{v:global:par:ineq} and \eqref{eq:imax:old}. Finally, we have that \eqref{eq:riots:4} is satisfied from \eqref{ineq:dpot:1}.

We therefore may appeal to the local conclusions \eqref{item:div:local:0}--\eqref{item:div:nonlocal} and the nonlocal outputs from \eqref{eq:inverse:div:error:stress}--\eqref{eq:inverse:div:error:stress:bound}, from which we have the following. First, we note that from \eqref{item:div:local:ii}, we have that \eqref{osc.loc.pot.1} is satisfied. Next, abbreviating $G \varrho \circ \Phi$ as $T_{i,j,k,\xi,\vecl,R}$, we have from \eqref{eq:inverse:div} and \eqref{eq:inverse:div:stress:1} that for $N\leq \frac{\Nfin}{4}-\dpot$ and $M\leq \frac{\Nfin}{5}$, 
\begin{align*}
    \left\| D^N D_{t, q}^M \mathcal{H} T_{i,j,k,\xi,\vecl,R} \right\|_{\sfrac32} 
    &\lesssim \left| \supp(\eta_{i,j,k,\xi,\vecl,R}^2) \right|^{\sfrac 23} \delta_{q+\bn} \Gamma_{q}^{2j-2} \La_q \Gamma_q^{42} \\& \qquad\qquad \times \lambda_{q+1}^{-1} \lambda_{q+1}^{\alpha+N} \MM{M, \Nindt, \tau_{q}^{-1}\Gamma_{q}^{i+13}, \Tau_{q}^{-1}\Ga_q^8} \\
\left\| D^N D_{t,q}^M \mathcal{H} T_{i,j,k,\xi,\vecl,R} \right\|_{\infty} 
    &\lesssim \Gamma_q^{\badshaq + 48} \Lambda_q \lambda_{q+1}^{-1} \lambda_{q+1}^{\alpha+N} \MM{M, \Nindt, \tau_{q}^{-1}\Gamma_{q}^{i+13}, \Tau_{q}^{-1}\Ga_q^8} \, , \notag \\
    &\lesssim \Gamma_{q+1}^{\badshaq-9} \lambda_{q+1}^N \MM{M, \Nindt, \tau_{q}^{-1}\Gamma_{q}^{i+13}, \Tau_{q}^{-1}\Ga_q^8} \, ,
\end{align*}
where we have used \eqref{eq:prepping:badshaq} to achieve the last inequality.  Notice that from \eqref{item:div:local:i}, the support of $\div \mathcal{H} T_{i,j,k,\xi,\vecl,R}$ is contained in the support of $T_{i,j,k,\xi,\vecl,R}$, which itself is contained in the support of $\eta_{i,j,k,\xi,\vecl,R}$. From this observation, we have that \eqref{H:eckel:support:osc:1} is satisfied. Finally, we have that \eqref{rho:eckel:estimate:1} holds after defining a potential $\vartheta$ as in \eqref{im:defining:a:potential} and appealing to standard Littlewood-Paley estimates and \eqref{eq:inverse:div:sub:1}.

Now we may apply the aggregation Corollaries~\ref{rem:summing:partition} and \ref{lem:agg.pt} with $H = \mathcal{H}T_{i,j,k,\xi,\vecl,R}$ and $\theta = \theta_2 = 2$, $p=\sfrac 32$ in the first case, or $\varpi=\Ga_{q+1}^{\badshaq-9}$ in the second case, to estimate
\begin{align}\notag 
     S^{q+1,l}_{O,R} := \sum_{i,j,k,\xi,\vecl} \mathcal{H} T_{i,j,k,\xi,\vecl,R} \, .
\end{align}
From \eqref{eq:agg:conc:1} and \eqref{eq:agg:conc:2} in the case $p=\sfrac 32$, and \eqref{eq:aggpt:conc:1} in the case $p=\infty$, we thus have that for $N,M$ in the same range as above,
\begin{align*}
    \left\| \psi_{i,q} D^N D_{t, q}^M  {S}^{q+1,l}_{O,R} \right\|_{\sfrac32} 
    &\lesssim \delta_{q+\bn} \La_q \Gamma_q^{50} \lambda_{q+1}^{-1} \lambda_{q+1}^{\alpha+N} \MM{M, \Nindt, \tau_{q}^{-1}\Gamma_{q}^{i+14}, \Tau_{q}^{-1}\Ga_q^8} \\
    \left\| \psi_{i,q} D^N D_{t, q}^M  {S}^{q+1,l}_{O,R} \right\|_{\infty} 
    &\lesssim \Gamma_{q+1}^{\badshaq-9} \lambda_{q+1}^N \MM{M, \Nindt, \tau_{q}^{-1}\Gamma_{q}^{i+14}, \Tau_{q}^{-1}\Ga_q^8} \, ,
\end{align*}
and so \eqref{eq:Onpnp:estimate:2} and \eqref{eq:Onpnp:estimate:2:new} follow for this term from \eqref{la.beats.de} and \eqref{condi.Nfin0}. 

For the nonlocal term, we first note that the left-hand side of the equality in \eqref{item:oscee:1} has zero mean, and so we may ignore the means of individual terms that get plugged into the inverse divergence since their sum will vanish. Then from \eqref{eq:inverse:div:error:stress}, \eqref{eq:inverse:div:error:stress:bound}, Remark~\ref{rem:lossy:choices}, and Lemma \ref{lem.cardinality}, we have that for $N,M\leq 2\Nind$,
\begin{align}
    \left\| D^N \Dtq^M \sum_{i,j,k,\xi,\vecl} \divR T_{i,j,k,\xi,\vecl,R} \right\|_\infty \leq 
    \la_{q+\bn}^{-5} \delta_{q+3\bn}^{\sfrac 32} \Tau_{q+\bn}^{4\Nindt}
    {\lambda_{q+1}^N} \tau_q^{-M} \, , \notag
\end{align}
matching the desired estimate in \eqref{eq:Onpnp:estimate:1}.

Finally, we must estimate the terms which include $B_{\pxi,\varphi}$ from \eqref{osc.err.low}.  However, we note that from Lemma~\ref{lem:a_master_est_p} $a_{\pxi,\varphi}^2$, differs in size relative to $a_{\pxi,R}^2$ by a factor of $r_q^{-\sfrac 23}$, which is exactly balanced out by the factor of $r_q^{\sfrac 23}$ in \eqref{osc.err.low}; the other differences in size actually make the estimates for $a_{\pxi,\varphi}^2$ stronger than for $a_{\pxi,R}^2$.  We therefore may argue exactly as above (in fact the estimates are slightly better since $\ov \rhob_\xi^4 < \ov \rhob_\xi^6$ and the power on $\Ga_q$ is smaller), and we omit further details.
\medskip

\noindent\texttt{Case 2:} Estimates for  \eqref{osc.err.med1-}.  As before, we fix $i,j,k,\xi,\vecl$. We apply Proposition \ref{prop:intermittent:inverse:div} with the low-frequency choices
\begin{subequations}
\begin{align}
     &G^\bullet = B_{\pxi,\varphi}^{\bullet} c_0 r_q^\frac23 \ov\chib_{\xi}^{4} (\Phi_{(i,k)}) \, , \quad \const_{G,\sfrac 32} = \left| \supp \eta_{i,j,k,\xi,\vecl,\vp}^2 \right|^{\sfrac 23} \delta_{q+\bn}\Gamma_{q}^{2j+25} \La_q \, , \quad \const_{G,\infty}= \Gamma_q^{\badshaq+35} \Lambda_q \, , \\
     &\lambda = \lambda_{q+1} \Gamma_q^{-1} \, , \quad \nu = \tau_q^{-1}\Gamma_q^{i+13} \, , \quad \Phi = \Phiik \, ,
\end{align}
\end{subequations}
as well as the choices from \eqref{eq:osc:general:choices}. The estimates in \eqref{eq:inverse:div:DN:G} and the assumption in \eqref{eq:inv:div:NM} hold due to Proposition~\ref{prop:bundling} and the estimates for $B_{\pxi,\varphi}r_q^{\sfrac 23}$ from \texttt{Case 1}.  \eqref{eq:DDpsi2}, \eqref{eq:DDpsi}, and \eqref{eq:DDv} are satisfied as in the previous substep.

To check the high-frequency assumptions, we set
\begin{subequations}
\begin{align}
    &\varrho = \mathbb{P}_{\neq 0} \left( \sum_I (\etab_{\xi}^I)^4 \right) \, , \quad \dpot \textnormal{ as in \eqref{i:par:11}} \, , \quad \vartheta = \delta_{i_1i_2} \delta_{i_3i_4} \dots \delta_{i_{\dpot-1}i_{\dpot}} \Delta^{-\sfrac \dpot 2} \varrho \, , \\
    &\mu = \Upsilon = \Upsilon' = \Lambda = \lambda_{q+\half} \, , \quad \const_{*,\sfrac 32} = \const_{*,\infty} = {\la_{q+\half}^\alpha} \, ,
    \end{align}
\end{subequations}
where $\alpha$ is chosen as in \eqref{eq:choice:of:alpha}. Then from Definition~\ref{def:etab}, standard Littlewood-Paley theory, and the same inequalities involving $\Ndec$ as in \texttt{Case 1}, we have that \eqref{eq:DN:Mikado:density} is satisfied, as well as the other high-frequency assumptions in \eqref{item:inverse:i}--\eqref{item:inverse:iv}.  The nonlocal assumptions are identical to those of \texttt{Case 1}, and are satisfied trivially.

We therefore may appeal to the local conclusions \eqref{item:div:local:0}--\eqref{item:div:nonlocal} and \eqref{eq:inverse:div:error:stress}--\eqref{eq:inverse:div:error:stress:bound}, from which we have the following.  First, we note that from \eqref{item:div:local:ii}, we have that \eqref{osc.loc.pot.1} is satisfied. Next, abbreviating $G \varrho \circ \Phi$ as $T_{i,j,k,\xi,\vecl,\varphi}$, we have from \eqref{eq:inverse:div} and \eqref{eq:inverse:div:stress:1} that for $N\leq \frac{\Nfin}{4}-\dpot$ and $M\leq \frac{\Nfin}{5}$, 
\begin{align*}
    \left\| D^N D_{t, q}^M \mathcal{H} T_{i,j,k,\xi,\vecl,\varphi} \right\|_{\sfrac32} 
    &\lesssim \left| \supp(\eta_{i,j,k,\xi,\vecl,\vp}^2) \right|^{\sfrac 23} \delta_{q+\bn} \Gamma_{q}^{2j-2} \La_q \Gamma_q^{50} \\& \qquad\qquad \times \lambda_{q+\half}^{-1} \lambda_{q+\half}^{N+\alpha} \MM{M, \Nindt, \tau_{q}^{-1}\Gamma_{q}^{i+13}, \Tau_{q}^{-1}\Ga_q^8} \\
\left\| D^N D_{t,q}^M \mathcal{H} T_{i,j,k,\xi,\vecl,\vp} \right\|_{\infty} 
    &\lesssim \Gamma_q^{\badshaq + 60} \Lambda_q \lambda_{q+\half}^{-1} \lambda_{q+\half}^{N+\alpha} \MM{M, \Nindt, \tau_{q}^{-1}\Gamma_{q}^{i+13}, \Tau_{q}^{-1}\Ga_q^8} \, , \notag \\
    &\lesssim \Gamma_{q+\half}^{\badshaq-9} \lambda_{q+\half}^N \MM{M, \Nindt, \tau_{q}^{-1}\Gamma_{q}^{i+13}, \Tau_{q}^{-1}\Ga_q^8} \, ,
\end{align*}
where we have used \eqref{eq:prepping:badshaq} to achieve the last inequality.  Notice that from \eqref{item:div:local:i}, the support of $\div \mathcal{H} T_{i,j,k,\xi,\vecl,\vp}$ is contained in the support of $T_{i,j,k,\xi,\vecl,\vp}$, which itself is contained in the support of $\eta_{i,j,k,\xi,\vecl,\vp}$. From this observation, we have that \eqref{H:eckel:support:osc:1} is satisfied.  Finally, we have that \eqref{rho:eckel:estimate:2} is satisfied from \eqref{eq:inverse:div:sub:1} after arguing in a manner similar to that in \texttt{Case 1}.

Now we may apply the aggregation Corollaries~\ref{rem:summing:partition} and \ref{lem:agg.pt} as in \texttt{Case 1} to estimate
\begin{align}\notag 
     S^{q+\half,l}_{O} := \sum_{i,j,k,\xi,\vecl} \mathcal{H} T_{i,j,k,\xi,\vecl,\vp} \, .
\end{align}
We find that for $N,M$ in the same range as above,
\begin{align*}
    \left\| \psi_{i,q} D^N D_{t, q}^M  {S}^{q+\half,l}_{O} \right\|_{\sfrac32} 
    &\lesssim \delta_{q+\bn} \La_q \Gamma_q^{60} \lambda_{q+\half}^{-1} \lambda_{q+\half}^{N+\alpha} \MM{M, \Nindt, \tau_{q}^{-1}\Gamma_{q}^{i+14}, \Tau_{q}^{-1}\Ga_q^8} \\
    \left\| \psi_{i,q} D^N D_{t, q}^M  {S}^{q+\half,l}_{O} \right\|_{\infty}
    &\lesssim \Gamma_{q+1}^{\badshaq-9} \lambda_{q+\half}^N \MM{M, \Nindt, \tau_{q}^{-1}\Gamma_{q}^{i+14}, \Tau_{q}^{-1}\Ga_q^8} \, ,
\end{align*}
and so \eqref{eq:Onpnp:estimate:2} and \eqref{eq:Onpnp:estimate:2:new} follow for this term from \eqref{la.beats.de} and \eqref{condi.Nfin0}. Finally, we must verify \eqref{osc:support:first} for $ S_O^{q+\half,l}$. This however follows from \eqref{item:div:local:ii}, which asserts that the support of $ S_O^{q+\half,l}$ is contained in the support of $\cup_{\pxi} a_{\pxi,\vp} \rhob_{\pxi}^{\vp} \circ \Phiik$, and \eqref{item:dodging:more:oldies} of Lemma~\ref{lem:dodging}. Finally, the nonlocal conclusions for $ S_O^{q+\half,l}$ follow in much the same way as in \texttt{Case 1}, and we omit further details.

\medskip

\noindent\texttt{Case 3:} Estimates for \eqref{osc.err.med0}, \eqref{osc.err.med}, and \eqref{osc.err.high} and $\diamond=R$. Fix $i,j,k,\xi,\vecl,I$ and set
\begin{align}
&G^{\bullet}= B_{\xi, i,j,k, \vecl, R}^{\bullet} \left((\chib_{(\xi)}^{R})^2(\etab_{\xi}^{I,R})^{2}\right) \circ \Phi_{(i,k)} \, , \quad
\Phi = \Phi_{(i,k)} \, , \quad {\nu = \tau_q^{-1}\Gamma_q^{i+13}} \, ,  \notag \\
&\const_{G,\sfrac 32} =\left| \supp(\eta_{i,j,k,\xi,\vecl,R}^2(\zetab_\xi^{I,R})^2) \right|^{\sfrac 23} \delta_{q+\bn} \Gamma_{q}^{2j+38}  {\La_q} + \la_\qbn^{-10} \, , \quad \const_{G,\infty}=
{\Ga_{q}^{\badshaq + 40}\La_q} \, , \quad \la = \la_{q+\half} \, , \label{shortening:some:stuff}
\end{align}
as well as the choices from \eqref{eq:osc:general:choices}. We then have that \eqref{eq:inv:div:NM} is satisfied as in the last step. Next, we have that \eqref{eq:inverse:div:DN:G} is satisfied by combining the corresponding bounds for $G^\bullet$ from the last step with the bounds for $\zetab_\xi^{I,R}$ from Definition~\ref{def:etab}.\footnote{We have added the extra $\la_\qbn^{-10}$ in the $\const_{G,\sfrac 32}$ bound in order to facilitate the creation of a pressure increment later.} The bounds in \eqref{eq:DDpsi2}--\eqref{eq:DDv} hold as before without any modifications. Finally, we have that the nonlocal assumptions in \eqref{eq:inv:div:wut}--\eqref{eq:riots:4} are satisfied for the same reasons as the previous cases.  At this point, we split the argument into subcases based on the differing synthetic Littlewood-Paley projectors in \eqref{osc.err.med}--\eqref{osc.err.high2}.
\smallskip

\noindent\texttt{Case 3a: } Estimates for \eqref{osc.err.med0} and $\diamond=R$. In order to set up the high-frequency assumptions for this case, we set
\begin{align*}
&\mu = \la_{q+\half}\Ga_q = \la_\qbn r_q \, , \quad \varrho = \tP_{q+\half+1}^\xi\mathbb{P}_{\neq0}(\varrho_{\pxi, R}^I)^2 \, , \quad \vartheta \textnormal{ as in Lemma~\ref{lem:special:cases}} \, , \quad \dpot \textnormal{ as in item~\eqref{i:par:10}} \\
&\const_{*,\sfrac 32} =  \la_{q+\half+1}^\alpha \, , \quad \const_{*, \infty}= {\left(\frac{\la_{q+\half+1}}{\la_{q+\bn}r_q}\right)^{2}}\la_{q+\half+1}^\alpha
\, ,\quad \Upsilon = \Upsilon' = \mu,\quad \Lambda= \la_{q+\half+1} \, ,
\end{align*}
where $\alpha$ is chosen as in \eqref{eq:choice:of:alpha}. We then have that \eqref{eq:DN:Mikado:density} is satisfied by appealing to estimate~\eqref{eq:lowest:shell:inverse} from Lemma~\ref{lem:special:cases} with $q=1$ and $p=\sfrac 32$, where we note that the assumption in \eqref{eq:moll:1:as} is satisfied with $\const_{\rho,q}=1$ and $\lambda=\lambda_\qbn$ from Proposition~\ref{prop:pipeconstruction}. We have in addition that \eqref{eq:inverse:div:parameters:0} and \eqref{eq:inverse:div:parameters:1} are satisfied by definition and by appealing to the same parameter inequalities as the previous steps. Finally, we have that the nonlocal assumption in \eqref{eq:riots:4} is satisfied from \eqref{ineq:dpot:1}.

We therefore may appeal to the local conclusions \eqref{item:div:local:0}--\eqref{item:div:nonlocal} of Proposition \ref{prop:intermittent:inverse:div} and \eqref{eq:inverse:div:error:stress}--\eqref{eq:inverse:div:error:stress:bound}, from which we have the following.  First, we note that from item~\eqref{item:div:local:iii}, \eqref{osc.loc.pot.2} is satisfied. Next, abbreviating $G \varrho \circ \Phi$ as $T_{i,j,k,\xi,\vecl,I,R}$, we have from \eqref{eq:inverse:div} and \eqref{eq:inverse:div:stress:1} that for $N\leq \frac{\Nfin}{4}-\dpot$ and $M\leq \frac{\Nfin}{5}$,
\begin{align*}
    \left\| D^N D_{t, q}^M \mathcal{H} T_{i,j,k,\xi,\vecl,I,R} \right\|_{\sfrac32} 
    &\lesssim \left(\left| \supp\left(\eta_{i,j,k,\xi,\vecl,R}^2 (\zetab_\xi^{I,R})^2 \right) \right|^{\sfrac 23} \delta_{q+\bn} \Gamma_{q}^{2j+39} \La_q + \la_\qbn^{-10} \right) \\
    & \qquad \times \left(\frac{\la_{q+\half+1}}{\la_{q+\half}}\right)^{\sfrac23} \la_{q+\half}^{-1} \la_{q+\half+1}^{N+\alpha} \MM{M, \Nindt, \tau_{q}^{-1}\Gamma_{q}^{i+13}, \Tau_{q}^{-1}\Ga_q^8} \\
    \left\| D^N D_{t,q}^M \mathcal{H} T_{i,j,k,\xi,\vecl,I,R} \right\|_{\infty} 
    &\lesssim \Gamma_q^{\badshaq + 40} \left(\frac{\la_{q+\half+1}}{\la_{q+\half}}\right)^{2} \Lambda_q \la_{q+\half}^{-1} \la_{q+\half+1}^{N+\alpha}  \MM{M, \Nindt, \tau_{q}^{-1}\Gamma_{q}^{i+13}, \Tau_{q}^{-1}\Ga_q^8} \, , \notag \\
    &\leq \Gamma_{q+\half}^{\badshaq-9} \la_{q+\half+1}^{N}  \MM{M, \Nindt, \tau_{q}^{-1}\Gamma_{q}^{i+13}, \Tau_{q}^{-1}\Ga_q^8} \, .
\end{align*}
We have used \eqref{eq:par:div:2} to simplify the second inequality. Notice that from \eqref{item:div:local:i}, the support of $\div \mathcal{H} T_{i,j,k,\xi,\vecl,I,R}$ is contained in the support of $T_{i,j,k,\xi,\vecl,I,R}$, which itself is contained in the support of $\eta_{i,j,k,\xi,\vecl,R}\zetab_\xi^{I,R}$.  From this observation, we have that \eqref{H:eckel:support:osc:2} is satisfied.  Finally, we have that \eqref{rho:eckel:estimate:3} is satisfied from \eqref{eq:inverse:div:sub:1} and Lemma~\ref{lem:special:cases} applied with $q=p=\sfrac 32,\infty$.

Now we may again apply the aggregation Corollaries~\ref{rem:summing:partition} and \ref{lem:agg.pt} to estimate
\begin{align}\notag 
     S^{q+\half+1,l}_{O,R} := \sum_{i,j,k,\xi,\vecl,I} \mathcal{H} T_{i,j,k,\xi,\vecl,I,R} \, .
\end{align}
From \eqref{eq:agg:conc:2} and \eqref{eq:aggpt:conc:2}, we then have that for $N,M$ in the same range as above,
\begin{align*}
    \left\| \psi_{i,q} D^N D_{t, q}^M  {S}^{q+\half+1,l}_{O,R} \right\|_{\sfrac32}
    &\lesssim  
    \delta_{q+\bn} \La_q \Gamma_q^{50} \left(\frac{\la_{q+\half+1}}{\la_{q+\half}}\right)^{\sfrac23} (\la_{q+\bn}r_q)^{-1} \notag\\
    &\qquad \qquad \times \la_{q+\half+1}^{N} \MM{M, \Nindt, \tau_{q}^{-1}\Gamma_{q}^{i+ {14}}, \Tau_{q}^{-1}\Ga_q^9} \, , \\
    &\leq \Ga_{q+\half+1}^{-10} \de_{q+\half+1+\bn}  \la_{q+\half+1}^{N} \MM{M, \Nindt, \tau_{q}^{-1}\Gamma_{q}^{i+ {14}}, \Tau_{q}^{-1}\Ga_q^9} \,, \\
     \left\| \psi_{i,q} D^N D_{t, q}^M  {S}^{q+\half+1,l}_{O,R} \right\|_{\infty}
    &\lesssim  
{\Gamma_{q+\half+1}^{\badshaq-9} \la_{q+\half+1}^N \MM{M, \Nindt, \tau_{q}^{-1}\Gamma_{q}^{i+14}, \Tau_{q}^{-1}\Ga_q^9}}
 \,, 
\end{align*}
where we have used \eqref{ineq:osc:general} to simplify the first inequality. Finally, the nonlocal conclusions follow in much the same way as in the previous cases, and so we omit further details.

\smallskip

\noindent\texttt{Case 3b: } Estimates for \eqref{osc.err.med} and \eqref{osc.err.high} and $\diamond=R$.  In order to set up the high-frequency assumptions for this case, we consider for the moment the cases when $m>q+\half+2$ and set
\begin{align}
&\mu = \la_{q+\half}\Ga_q = \la_\qbn r_q \, , \quad \varrho = \tP_{(m-1,m]}^\xi (\varrho_{\pxi, R}^I)^2 \, , \quad \vartheta \textnormal{ as in Lemma~\ref{lem:LP.est}} \, , \quad \dpot \textnormal{ as in item~\eqref{i:par:10}} \notag \\
&\const_{*,\sfrac 32} = {\left(\frac{\min(\la_m,\la_\qbn)}{\la_{q+\bn}r_q}\right)^{\sfrac23}} \, , \quad \const_{*, \infty}= {\left(\frac{\min(\la_m,\la_\qbn)}{\la_{q+\bn}r_q}\right)^{2}}\la_{q+\half+1}^\alpha \, , \notag \\
&\Upsilon = \la_{m-1} \, , \quad \Upsilon' = \la_m \, , \quad \Lambda = \min(\la_m,\la_\qbn) \, . \label{more:shortening}
\end{align}
We then have that \eqref{eq:DN:Mikado:density} is satisfied by appealing to \eqref{eq:LP:div:estimates} with $q=1$ and $p=\sfrac 32,\infty$; we note that \eqref{eq:gen:id:assump} is satisfied for $q=1$ and $\const_{\rho,q}=1$ and $\la=\la_\qbn$ as in the last step. Next, we have that \eqref{eq:inverse:div:parameters:0}--\eqref{eq:inverse:div:parameters:1} are satisfied by definition and immediate computation and the same inequalities as in the previous steps.  Finally, we have that the nonlocal assumption in \eqref{eq:riots:4} is satisfied from \eqref{ineq:dpot:1}.

In the case of $m= q+\half +2$, we have to take an extra step to minimize the gap between $\Upsilon$ and $\Upsilon'$ in order to ensure that the second inequality in \eqref{eq:inverse:div:parameters:0} is satisfied.  Towards this end, we decompose the synthetic Littlewood-Paley operator further as 
\begin{equation}\label{eq:damn:projection}
\tP_{(q+\half+1, q+\half+2]}^\xi
:=\tP_{(q+\half+1, q+\half+{\sfrac32}]}^\xi
+\tP_{(q+\half+{\sfrac32}, q+\half+2]}^\xi \, ,
\end{equation}
where the $q+\half+\sfrac 32$ portion of the projector correponds to the frequency which is the geometric means of $\la_{q+\half+1}$ and $\la_{q+\half+2}$.  This extra division helps us minimize the gap between $\Upsilon$ and $\Upsilon'$.  Then we can set
\begin{align*}
&\mu = \la_{q+\half}\Ga_q = \la_\qbn r_q \, , \quad \varrho = \tP_{\bullet}^\xi (\varrho_{\pxi, R}^I)^2 \, , \quad \vartheta \textnormal{ as in Lemma~\ref{lem:LP.est}} \, , \quad \dpot \textnormal{ as in item~\eqref{i:par:10}} \\
&\const_{*,\sfrac 32} = {\left(\frac{\la_{q+\half+2}}{\la_{q+\bn}r_q}\right)^{\sfrac23}} \, , \quad \const_{*, \infty}= {\left(\frac{\la_{q+\half+2}}{\la_{q+\bn}r_q}\right)^{2}}\la_{q+\half+1}^\alpha \, ,  \\
&\Upsilon = \la_{q+\half+1} \, , \quad \Upsilon' = \la_{q+\half+\sfrac 32} \textnormal{  if $\bullet$ corresponds to the first projector} \, , \\
&\Upsilon = \la_{q+\half+\sfrac 32} \, , \quad \Upsilon' = \la_{q+\half+2} \textnormal{  if $\bullet$ corresponds to the second projector} \, .
\end{align*}
We then have that \eqref{eq:DN:Mikado:density} is satisfied by appealing to \eqref{eq:LP:div:estimates} with $q=1$ and $p=\sfrac 32,\infty$ as before. Next, we have that \eqref{eq:inverse:div:parameters:0}--\eqref{eq:inverse:div:parameters:1} are satisfied by definition and immediate computation (here we crucially use the extra subdivision to ensure that the second inequality in \eqref{eq:inverse:div:parameters:0} holds) and the same inequalities as in the previous steps.  Finally, we again have that the nonlocal assumption in \eqref{eq:riots:4} is satisfied from \eqref{ineq:dpot:1}.

We therefore may appeal to the local conclusions \eqref{item:div:local:0}--\eqref{item:div:nonlocal} of Proposition \ref{prop:intermittent:inverse:div} and \eqref{eq:inverse:div:error:stress}--\eqref{eq:inverse:div:error:stress:bound}, from which we have the following.  First, we note that from item~\eqref{item:div:local:iii}, \eqref{osc.loc.pot.2} is satisfied. Next, abbreviating $G \varrho \circ \Phi$ as $T_{i,j,k,\xi,\vecl,I,R}$, we have from \eqref{eq:inverse:div} and \eqref{eq:inverse:div:stress:1} that for $N\leq \frac{\Nfin}{4}-\dpot$ and $M\leq \frac{\Nfin}{5}$,
\begin{align*}
    \left\| D^N D_{t, q}^M \mathcal{H} T_{i,j,k,\xi,\vecl,I,R} \right\|_{\sfrac32} 
    &\lesssim \left( \left| \supp(\eta_{i,j,k,\xi,\vecl,R}^2) \right|^{\sfrac 23} \delta_{q+\bn} \Gamma_{q}^{2j+39} \La_q + \la_\qbn^{-10} \right) \left(\frac{ \min(\la_{m}, \la_{q+\bn})}{\la_{q+\bn}r_q}\right)^{2/3} \lambda_{m-1}^{-2} \\& \qquad\qquad \times ( \min(\la_{m}, \la_{q+\bn}))^{N+1} \MM{M, \Nindt, \tau_{q}^{-1}\Gamma_{q}^{i+13}, \Tau_{q}^{-1}\Ga_q^8} \, , \\
\left\| D^N D_{t,q}^M \mathcal{H} T_{i,j,k,\xi,\vecl,I,R} \right\|_{\infty} 
    &\lesssim \Gamma_q^{\badshaq + 40} \left(\frac{ \min(\la_{m}, \la_{q+\bn})}{\la_{q+\bn}r_q}\right)^2 \Lambda_q \lambda_{m-1}^{- {2}} \notag\\
    &\qquad \qquad \times ( {\min(\la_{m}, \la_{q+\bn})})^{N+1} \MM{M, \Nindt, \tau_{q}^{-1}\Gamma_{q}^{i+13}, \Tau_{q}^{-1}\Ga_q^8} \, , \notag \\
    &\lesssim \Gamma_{q+\half}^{\badshaq-9} ( \min(\la_{m}, \la_{q+\bn}))^N \MM{M, \Nindt, \tau_{q}^{-1}\Gamma_{q}^{i+13}, \Tau_{q}^{-1}\Ga_q^8} \, ,
\end{align*}
where we have used \eqref{eq:par:div:2} to achieve the last inequality. Notice that from \eqref{item:div:local:i}, the support of $\div \mathcal{H} T_{i,j,k,\xi,\vecl,I,R}$ is contained in the support of $T_{i,j,k,\xi,\vecl,I,R}$, which itself is contained in the support of $\eta_{i,j,k,\xi,\vecl,R}\zetab_\xi^{I,R}$.  From this observation, we have that \eqref{H:eckel:support:osc:2} is satisfied.  Furthermore, we have that \eqref{rho:eckel:estimate:4} is satisfied from \eqref{eq:inverse:div:sub:1} and Lemma~\ref{lem:special:cases} applied with $q=p=\sfrac 32,\infty$. Finally, we have that \eqref{rho:eckel:support} is satisfied due to item~\eqref{item:div:local:i} and \eqref{eq:LP:div:support}. We note also that \eqref{osc:support:first} follows from \eqref{rho:eckel:support} and \eqref{eq:dodging:oldies}.

Now we may again apply the aggregation Corollaries~\ref{rem:summing:partition} and \ref{lem:agg.pt} to estimate
\begin{align}\notag 
     S^{m,l}_{O,R} := \sum_{i,j,k,\xi,\vecl,I} \mathcal{H} T_{i,j,k,\xi,\vecl,I,R} \, .
\end{align}
From \eqref{eq:agg:conc:2} and \eqref{eq:aggpt:conc:2}, we then have that for $N,M$ in the same range as above,
\begin{align*}
\left\| \psi_{i,q} D^N D_{t, q}^M  {S}^{m,l}_{O,R} \right\|_{\sfrac32}
&\lesssim \delta_{q+\bn} \La_q \Gamma_q^{50} \left(\frac{ {\min(\la_{m}, \la_{q+\bn})}}{\la_{q+\bn}r_q}\right)^{\sfrac23} \lambda_{m-1}^{-2} \notag\\
&\qquad \qquad \times {\min(\la_{m}, \la_{q+\bn})}^{N+1} \MM{M, \Nindt, \tau_{q}^{-1}\Gamma_{q}^{i+ {14}}, \Tau_{q}^{-1}\Ga_q^9} \\
&\les \Ga_m^{-10} \de_{m+\bn}  {(\min(\la_{m}, \la_{q+\bn}))}^{N} \MM{M, \Nindt, \tau_{q}^{-1}\Gamma_{q}^{i+ {14}}, \Tau_{q}^{-1}\Ga_q^9} \,, \\
\left\|\psi_{i,q} D^N D_{t, q}^M  {S}^{m,l}_{O,R} \right\|_{\infty} 
&\lesssim   {\Gamma_{m}^{\badshaq-9} (\min(\la_{m}, \la_{q+\bn}))^N \MM{M, \Nindt, \tau_{q}^{-1}\Gamma_{q}^{i+14}, \Tau_{q}^{-1}\Ga_q^9}}
\end{align*}
where we have used \eqref{ineq:osc:general} to simplify the first inequality. Finally, the nonlocal conclusions follow in much the same way as in the previous cases, and so we omit further details.
\medskip

\noindent\texttt{Case 4: } Estimates for \eqref{osc.err.med0}, \eqref{osc.err.med}, and \eqref{osc.err.high} and $\diamond=\ph$.
Estimates for these follow from similar arguments as in the cases when $\diamond=R$. Indeed, the only significant differences are that the estimates for $a_{\pxi,\varphi}^2$ than those of $a_{\pxi,R}^2$ are \emph{worse} by a factor of $r_q^{-\sfrac 23}$ from Lemma~\ref{lem:a_master_est_p}, while the estimates for $\varrho$ encoded in the constants $\const_{*,\sfrac 32}$ and $\const_{*,\infty}$ are \emph{better} by a factor of $r_q^{\sfrac 23}$ from Proposition~\ref{prop:pipe.flow.current}. Therefore, to compensate such loss or gain, we define $G^\bullet =B_{\xi, i,j,k, \vecl, \ph}^{\bullet} \left((\chib_{(\xi)}^{\ph})^2(\etab_{\xi}^{I,\ph})^{2}\right) \circ \Phi_{(i,k)}r_q^{\sfrac23}$ with the extra factor $r_q^{\sfrac23}$ and define $\varrho$ analogous to the case $\diamond=R$ but with the extra factor $r_q^{-\sfrac23}$. Then, the same choice of parameters and functions as in the case of $\diamond=R$ will lead to the desired estimates. We omit further details.
\medskip

\noindent\texttt{Case 5:} Estimates for \eqref{osc.err.high2}. Here we apply Proposition~\ref{prop:intermittent:inverse:div} with $p=\infty$ and the following choices.  The low-frequency assumptions in Part 1 are exactly the same as the $L^\infty$ low-frequency assumptions in \texttt{Case 3} and \texttt{Case 4}.  For the high-frequency assumptions, we recall the choice of $N_{**}$ from \eqref{i:par:10} and set
\begin{align}
    &\varrho_R = (\Id - \tilde{\mathbb{P}}_{q+\bn+1}^\xi ) \mathbb{P}_{\neq 0} \left( \varrho_{\pxi,R}^I \right)^2  , \quad
    \varrho_\ph = (\Id - \tilde{\mathbb{P}}_{q+\bn+1}^\xi ) \mathbb{P}_{\neq 0} \left( \varrho_{\pxi,\ph}^I \right)^2 r_q^{-\sfrac23}  , \quad
    \vartheta_\diamond^{i_1i_2\dots i_{\dpot-1}i_\dpot} = \delta^{i_1i_2\dots i_{\dpot-1}i_\dpot} \Delta^{-\sfrac \dpot 2}\varrho_\diamond  ,  \notag\\
    &\Lambda=\lambda_{q+\bn}\, , \quad \mu = \Upsilon=\Upsilon' = \lambda_{q+\half}\Ga_q \, , \, , \quad \const_{*,\infty} = \left( \frac{\lambda_\qbn}{\lambda_{q+\bn+1}} \right)^{N_{**}} \lambda_\qbn^3 \, , \quad \Ndec \textnormal{ as in \eqref{i:par:9}} \, , \quad \dpot=0 \, . \notag
\end{align}
Then we have that item~\eqref{item:inverse:i} is satisfied by definition, item~\eqref{item:inverse:ii} is satisfied as in the previous steps, \eqref{eq:DN:Mikado:density} is satisfied using Propositions~\ref{prop:pipeconstruction} and \ref{prop:pipe.flow.current} and \eqref{eq:remainder:inverse} from Lemma~\ref{lem:special:cases}, \eqref{eq:inverse:div:parameters:0} is satisfied by definition and as in the previous steps, and \eqref{eq:inverse:div:parameters:1} is satisfied by \eqref{condi.Ndec0}.  For the nonlocal assumptions, we choose $M_\circ,N_\circ=2\Nind$ so that \eqref{eq:inv:div:wut}--\eqref{eq:inverse:div:v:global:parameters} are satisfied as in \texttt{Case 1}, and \eqref{eq:riots:4} is satisfied from \eqref{ineq:Nstarz:1}. We have thus satisfied all the requisite assumptions, and we therefore obtain nonlocal bounds very similar to those from the previous steps, which are consistent with \eqref{eq:Onpnp:estimate:1} at level $q+\bn$. We omit further details.
\end{proof}

\begin{lemma*}[\bf Low shells have no pressure increment]\label{lem:osc:no:pressure}
The errors $ S^{q+1}_O$ and $ S^{q+\half}_O$ require no pressure increment as they are already dominated by intermittent pressure from the previous step. More precisely, we have that for {$N,M\leq \sfrac{\Nfin}{10}$,}
\begin{subequations}
\begin{align}
    \label{eq:lowshell:nopr:1}
    \left|\psi_{i,q} D^N \Dtq^M S^{q+1,l}_O \right| &\leq \Ga^{-100}_{q+1} \pi_{q}^{q+1}
    \la_{q+1}^N \MM{M, \Nindt, \tau_{q}^{-1}\Gamma_{q}^{i+14}, \Tau_{q}^{-1}\Ga_q^8} \, , \\
    \label{eq:lowshell:nopr:2}
    \left|\psi_{i,q} D^N \Dtq^M S^{q+\half,l}_O \right| &\leq \Ga^{-100}_{q+\half} \pi_{q}^{q+\half} \la_{q+\half}^N \MM{M, \Nindt, \tau_{q}^{-1}\Gamma_{q}^{i+14}, \Tau_{q}^{-1}\Ga_q^8} \, .
\end{align}
\end{subequations}
\end{lemma*}
\begin{proof}
We first note that the application of Proposition \ref{prop:intermittent:inverse:div} in \texttt{Case 1} of the proof of Lemma~\ref{lem:oscillation:general:estimate} can be supplemented with Remark~\ref{rem:pointwise:inverse:div}.  Specifically, we may set
\begin{align}\label{eq:pi:osc:lowest}
    \ov \pi = \pi_\ell \Gamma_q^{40} \Lambda_q \, ,
\end{align}
so that \eqref{eq:inv:div:extra:pointwise} follows from the definition of $B_{\pxi,R}$ in \eqref{eq:rhoxidiamond:def} and \eqref{e:a_master_est_p_R_pointwise}.  Then from \eqref{eq:divH:formula}, \eqref{eq:inverse:div:sub:1}, and \eqref{eq:inv:div:extra:conc}, we have that
\begin{align*}
    \left|D^N \Dtq^M \mathcal{H}T_{i,j,k,\xi,\vecl,R} \right| &\lec \pi_\ell \Gamma_q^{50} \Lambda_q \lambda_{q+1}^{-1} \la_{q+1}^N \MM{M, \Nindt, \tau_{q}^{-1}\Gamma_{q}^{i+13}, \Tau_{q}^{-1}\Ga_q^8} \, .
\end{align*}
We pause also to note that \eqref{eq:H:eckel:osc} in this case follows from \eqref{eq:divH:formula} and \eqref{eq:inv:div:extra:conc}.  Now applying the aggregation Corollary~\ref{lem:agg.pt} with $H= \mathcal{H}T_{i,j,k,\xi,\vecl,R}$, $\varpi = \pi_\ell \Gamma_q^{50}\Lambda_q$, and $p=1$ along with \eqref{eq:ind.pr.anticipated}, \eqref{ind:pi:lower}, and \eqref{la.beats.de} gives \eqref{eq:lowshell:nopr:1}.

The proof of \eqref{eq:lowshell:nopr:2} follows similarly from supplementing \texttt{Case 2} of the proof of Lemma~\ref{lem:oscillation:general:estimate} with pointwise assumptions. We omit further details.
\end{proof}

\begin{lemma*}[\bf Pressure increment]\label{lem:oscillation:pressure}
For every $q+ \half +1 \leq m \leq q+\bn$, there exists a function $\si_{S^{m}_O} = \si_{S^{m}_O}^+ - \si_{S^{m}_O}^-$ such that the following hold.\index{$\sigma_{S^m_O}$}
\begin{enumerate}[(i)]
    \item We have that
\begin{subequations}
\begin{align}
    \label{eq:o.p.1}
    \left|\psi_{i,q} D^N \Dtq^M  S^{ m,l}_{O}\right| &< \left(\si_{S^{ m}_O}^+ + \de_{q+3\bn}\right) \left(\lambda_{m}\Gamma_q\right)^N \MM{M,\Nindt,\tau_q^{-1}\Gamma_{q}^{i+16},\Tau_q^{-1} {\Ga_q^9}}\\
    \label{eq:o.p.2}
    \left|\psi_{i,q} D^N \Dtq^M \si_{S^{ m}_O}^+\right| &< \left(\si_{S^{ m}_O}^+ +\de_{q+3\bn}\right) \left(\lambda_{m}\Gamma_q\right)^N \MM{M,\Nindt,\tau_q^{-1}\Gamma_{q}^{i+16},\Tau_q^{-1}\Ga_q^9}\\
    \label{eq:o.p.3}
    \norm{\psi_{i,q} D^N \Dtq^M \si_{S^{ m}_O}^+}_{\sfrac32} &\leq \Ga_{m}^{-9} \de_{ m+\bn} \left(\lambda_{m}\Gamma_q \right)^N \MM{M,\Nindt,\tau_q^{-1}\Gamma_{q}^{i+16},\Tau_q^{-1}\Ga_q^9}\\
    \label{eq:o.p.3.inf}
     {
    \norm{D^N D_{t}^M \si_{S^{ m}_O}^+}_{\infty} }
    &\leq \Gamma_{q+1}^{\badshaq-9} (\la_{m} \Ga_q)^N \MM{M,\Nindt,\tau_q^{-1}\Ga_q^{i+16},\Tau_q^{-1}\Gamma_q^9} \\
    \label{eq:o.p.4}
    \left|\psi_{i,q} D^N \Dtq^M \si_{S^{ m}_O}^-\right| &\lec \Ga_{q+\half}^{-100} \pi_q^{q+\half}  \left(\lambda_{q+\half}\Gamma_q\right)^N \MM{M,\Nindt,\tau_q^{-1}\Gamma_{q}^{i+16},\Tau_q^{-1}\Ga_q^9}
\end{align}
\end{subequations}
for all $N,M < \sfrac{\Nfin}{100}$.
\item For $ m\geq q+\half +2$, we  have that
\begin{equation}\label{eq:os.p.6}
\begin{split}
B\left( \supp \hat w_{q'}, \lambda_{q'}^{-1} \Gamma_{q'+1} \right) \cap \supp (\si_{S^{m}_O}^+) &= \emptyset \qquad \forall q+1\leq q' \leq  m-1\\
B\left( \supp \hat w_{q'}, \lambda_{q'}^{-1} \Gamma_{q'+1} \right) \cap \supp (\si_{S^{m}_O}^-) &= \emptyset \qquad \forall q+1\leq q' \leq  q+\half\, .
\end{split}
\end{equation}
\item \label{i:pc:5:ER} Define 
\begin{equation}\label{def:bmu:ER:osc}
    \bmu_{\sigma_{S_O^m}}(t) = \int_0^t \left \langle \Dtq \sigma_{S_O^m}  \right \rangle (s) \, ds \, .
\end{equation}
Then we have that
    \begin{align}\label{eq:thurzday:night}
      \left|\frac{d^{M+1}}{dt^{M+1}} \bmu_{\sigma_{S_O^m}} \right| 
      \leq (\max(1, T))^{-1}\delta_{q+3\bn}^2 \MM{M,\Nindt,\tau_q^{-1},\Tau_{q+1}^{-1}}
    \end{align}
    for $0\leq M\leq 2\Nind$. 
\end{enumerate}
\end{lemma*}
\begin{proof}[Proof of Lemma~\ref{lem:oscillation:pressure}]
We follow the case numbering from Lemma~\ref{lem:oscillation:general:estimate}.  Since we have shown in Lemma~\ref{lem:osc:no:pressure} that the low shells have no pressure increment, we only need to analyze \texttt{Cases 3} and \texttt{4}.  Since the only difference between \texttt{Case 3} and \texttt{Case 4} is the rebalancing of $r_q^{\sfrac 23}$, we shall only hint at the proofs in \texttt{Case 4} and focus on the case $\diamond=R$. We divide into subcases \texttt{3a} and \texttt{3b} and apply Proposition \ref{lem.pr.invdiv2}.
\smallskip

\noindent\texttt{Case 3a:} pressure increment for \eqref{osc.err.med0} and $\diamond=R$.
Recall that Part 1 of Proposition \ref{lem.pr.invdiv2} requires preliminary assumptions which are the same as those from the inverse divergence, along with pointwise bounds corresponding to Remark~\ref{rem:pointwise:inverse:div}. Since we have already chosen parameters corresponding to the inverse divergence, we simply set $\ov \pi = \pi_\ell \Gamma_q^{50} \Lambda_q$, which verifies \eqref{eq:H:eckel:osc} in this case.  Then the assumption in \eqref{eq:inv:div:extra:pointwise} follows from the pointwise estimates for $B_{\pxi,R}$ used in Lemma~\ref{lem:osc:no:pressure} along with Proposition~\ref{prop:bundling}, Lemma~\ref{lem:finer:checkerboard:estimates}, and Corollary~\ref{cor:deformation} to estimate $\left((\chib_{(\xi)}^{R})^2(\etab_{\xi}^{I,R})^{2}\right) \circ \Phi_{(i,k)}$.

In order to check the additional assumptions from Part 2, we set
\begin{align}
    &N_{**} \textnormal{ as in \eqref{i:par:10}} \,, \quad \NcutLarge, \NcutSmall \textnormal{ as in \eqref{i:par:6}} \, , \quad \Ga=\Ga_q^{\sfrac 12} \, , \quad \delta_{\rm tiny} = \delta^2_{q+3\bn} \, , \label{eq:desert:choices:1:ER}\\
    &\bm = 1 \, , \quad \mu_0 = \lambda_{q+\half+1}\Ga_q^{-1} \, , \quad \mu_{\bm} = \mu_1= \lambda_{q+\half+1}\Ga_q^{2} \, . \notag
\end{align}
Then \eqref{i:st:sample:wut}--\eqref{i:st:sample:wut:wut} hold from \eqref{condi.Nfin0}, \eqref{i:st:sample:wut:wut:wut} holds from \eqref{ineq:Nstarstar:dpot}, \eqref{eq:sample:prop:Ncut:1} holds from \eqref{condi.Ncut0.1}, \eqref{eq:sample:prop:Ncut:2} holds from \eqref{condi.Ncut0.2}, \eqref{eq:sample:prop:Ncut:3} holds from \eqref{condi.Nfin0}, \eqref{eq:sample:prop:decoup} holds from \eqref{condi.Ndec0}, \eqref{eq:sample:prop:par:00} holds by definition, \eqref{eq:sample:prop:parameters:0} holds by definition and immediate computation, \eqref{eq:sample:riots:4} holds due to \eqref{ineq:dpot:1}, and \eqref{eq:sample:riot:4:4} holds due to \eqref{ineq:Nstarz:1}.

At this point, we appeal to the conclusions from Part 3 to construct a pressure increment and delineate its properties.  First, from \eqref{d:press:stress:sample}--\eqref{est.S.by.pr.final2} and \eqref{condi.Nfin0}, we have that there exists a pressure increment $\sigma_{ \divH T_{i,j,k,\xi,\vecl,I,R}^{q+\half+1}} =\sigma_{ \divH T_{i,j,k,\xi,\vecl,I,R}^{q+\half+1}}^+-\sigma_{ \divH T_{i,j,k,\xi,\vecl,I,
R}^{q+\half+1}}^-$ such that for $N,M\leq \sfrac{\Nfin}{7}$, 
\begin{align}\label{eq:desert:cowboy:1:ER}
    \left| D^N \Dtq^M \divH T_{i,j,k,\xi,\vecl,I,R}^{q+\half+1}  \right| \lesssim \left( \sigma_{ \divH T_{i,j,k,\xi,\vecl,I,R}^{q+\half+1}}^+ + \delta_{q+3\bn}^2 \right) (\lambda_{q+\half+1}\Ga_q)^N \MM{M,\Nindt,\tau_q^{-1}\Ga_q^{i+14},\Tau_q^{-1}\Ga_q^{9}} \, .
\end{align}
From \eqref{eq:inverse:div:linear} and \eqref{est.S.pr.p.support:1}, we have that
\begin{align}
    \supp \left( \sigma_{ \divH T_{i,j,k,\xi,\vecl,I,R}^{q+\half+1}}^+ \right) \subseteq \supp \left( \divH T_{i,j,k,\xi,\vecl,I,R}^{q+\half+1}  \right) \subseteq \supp \left( a_{\pxi,R} \left(\rhob_\pxi^R \zetab_\xi^{I} \right)\circ\Phiik \right) \, . \label{eq:ocdc:support:ER} 
\end{align}

Now define
\begin{align}\label{eq:desert:pressure:def:ER}
    \sigma_{ S_{O,R}^{q+\half+1}}^{\pm} = \sum_{i,j,k,\xi,\vecl,I} \sigma_{ \divH T_{i,j,k,\xi,\vecl,I,R}^{q+\half+1}}^\pm \, . 
\end{align}
Then \eqref{eq:oooooldies} gives that \eqref{eq:os.p.6} is satisfied for $m=q+\half+1$. 
From \eqref{eq:desert:cowboy:1:ER}, \eqref{eq:desert:cowboy:sum}, \eqref{eq:inductive:partition}, and Corollary~\ref{lem:agg.pt} with 
\begin{align}
&H=\divH T_{i,j,k,\xi,\vecl,I,R}^{q+\half+1} \, , \qquad \varpi=\left[\sigma_{ \divH T_{i,j,k,\xi,\vecl,I,R}^{q+\half+1}}^+ + \delta_{q+3\bn}^2\right] \mathbf{1}_{\supp a_{\pxi,R}\rhob_\pxi^R\zetab_\xi^I}  \,  , \qquad p=1 \, , \notag 
\end{align}
we have that for $N,M\leq \sfrac{\Nfin}{7}$,
\begin{align}\label{eq:desert:cowboy:2:ER}
    \left| \psi_{i,q} D^N \Dtq^M \sum_{i',j,k,\xi,\vecl,I} \divH T_{i',j,k,\xi,\vecl,I,R}^{q+\half+1}  \right| &\lesssim \left(  \sigma_{ S_{O,R}^{q+\half+1}}^+  + \delta_{q+3\bn}^2 \right) \notag\\
    &\qquad \times (\lambda_{q+\half+1}\Ga_q)^N \MM{M,\Nindt,\tau_q^{-1}\Ga_q^{i+15},\Tau_q^{-1}\Ga_q^{9}} \, .
\end{align}
We therefore have that \eqref{eq:o.p.1} is satisfied for $m=q+\half+1$. From \eqref{est.S.prbypr.pt}, \eqref{condi.Nfin0}, and \eqref{condi.Nindt}, we have that for $N,M\leq \sfrac{\Nfin}{7}$,
\begin{align}
     \left| D^N \Dtq^M \sigma_{ \divH T_{i,j,k,\xi,\vecl,I,R}^{q+\half+1}}^+ \right| \lesssim \left( \sigma_{ \divH T_{i,j,k,\xi,\vecl,I,R}^{q+\half+1}}^+ + \delta_{q+3\bn}^2 \right) (\lambda_{q+\half+1}\Ga_q)^N \MM{M,\Nindt,\tau_q^{-1}\Ga_q^{i+15},\Tau_q^{-1}\Ga_q^{9}} \, . \label{eq:desert:cowboy:4:ER}
\end{align}
From \eqref{eq:desert:cowboy:4:ER}, \eqref{eq:desert:cowboy:sum}, \eqref{eq:inductive:partition}, and Corollary~\ref{lem:agg.pt} with 
\begin{align}
&H=\sigma^+_{\divH T_{i,j,k,\xi,\vecl,I,R}^{q+\half+1}} \, , \qquad \varpi= \left[ H + \delta_{q+3\bn}^2 \right] \mathbf{1}_{\supp a_{\pxi,R}\rhob_\pxi^R\zetab_\xi^I}  \,  , \qquad p=1 \, , \notag 
\end{align}
we have that \eqref{eq:o.p.2} is satisfied for $m=q+\half+1$.

Next, from \eqref{est.S.pr.p}, we have that
\begin{align}
    \left\| \sigma^{\pm}_{\divH T_{i,j,k,\xi,\vecl,I,R}^{q+\half+1}} \right\|_{\sfrac 32} &\les \left( \left| \supp(\eta_{i,j,k,\xi,\vecl,R}^2(\zetab_\xi^{I,R})^2) \right|^{\sfrac 23} \delta_{q+\bn} \Gamma_{q}^{2j+38}  {\La_q} + \la_\qbn^{-10} \right) {\left(\frac{\la_{q+\half+1}}{\la_{q+\bn}r_q}\right)^{\sfrac23}}\la_{q+\half+1}^\alpha \la_{q+\half}^{-1} \notag \, .
\end{align}
Now from \eqref{eq:desert:pressure:def:ER}, \eqref{ineq:osc:general}, and Corollary~\ref{rem:summing:partition} with $\theta=2$, $\theta_1=0$, $\theta_2=2$, $H=\sigma^{\pm}_{\divH T_{i,j,k,\xi,\vecl,I,\diamond}^{q+\half+1}}$, and $p=\sfrac 32$, we have that
\begin{align}
    \left\| \psi_{i,q} \sigma_{ S_{O,R}^{q+\half+1}}^{\pm} \right\|_{\sfrac 32} &\les \delta_{\qbn+\half+1} \Ga_{q+\half+1}^{-10} \, . \notag 
\end{align}
Combined with \eqref{eq:o.p.2}, this verifies \eqref{eq:o.p.3} at level $q+\half+1$.  Arguing now for $p=\infty$ from \eqref{est.S.pr.p}, we have that
\begin{align}
    \left\| \sigma^{\pm}_{\divH T_{i,j,k,\xi,\vecl,I,R}^{q+\half+1}} \right\|_{\infty} &\les \Ga_q^{\badshaq+40}\La_q {\left(\frac{\la_{q+\half+1}}{\la_{q+\bn}r_q}\right)^{2}}\la_{q+\half+1}^\alpha \la_{q+\half}^{-1} \notag \, .
\end{align}
Now from \eqref{eq:desert:pressure:def:ER}, \eqref{eq:par:div:2}, and Corollary~\ref{lem:agg.pt} with $H=\sigma^{\pm}_{\divH T_{i,j,k,\xi,\vecl,I,R}^{q+\half+1}}$, $\varpi=\mathbf{1}_{\supp a_{\pxi,R}\rhob_\pxi^R\zetab_\xi^I}$ and $p=1$, we have that
\begin{align}
    \left\| \psi_{i,q} \sigma_{ S_{O,R}^{q+\half+1}}^{\pm} \right\|_\infty &\les \Ga_q^{\badshaq+40}\La_q {\left(\frac{\la_{q+\half+1}}{\la_{q+\bn}r_q}\right)^{2}}\la_{q+\half+1}^\alpha \la_{q+\half}^{-1} \leq \Ga_{q+\half+1}^{\badshaq-100} \, . \notag 
\end{align}
Combined again with \eqref{eq:o.p.2}, this verifies \eqref{eq:o.p.3.inf} at level $q+\half+1$. 

Finally, from \eqref{est.S.prminus.pt}, \eqref{condi.Nindt}, \eqref{condi.Nfin0}, \eqref{o:5}, \eqref{ind:pi:lower}, and \eqref{eq:ind.pr.anticipated}, we have that for $N,M\leq \sfrac{\Nfin}{7}$, 
\begin{align}
    \left| D^N \Dtq^M \sigma^{-}_{\divH T_{i,j,k,\xi,\vecl,I,R}^{q+\half+1}} \right| &\les {\left(\frac{\la_{q+\half+1}}{\la_{q+\bn}r_q}\right)^{\sfrac23}}\la_{q+\half+1}^\alpha \la_{q+\half}^{-1} \pi_\ell \Ga_q^{50} \La_q \notag\\
    &\qquad \qquad \times (\la_{q+\half}\Ga_q)^N \MM{M,\Nindt,\tau_q^{-1}\Ga_q^{i+15},\Tau_q^{-1}\Ga_q^9} \notag\\
    &\leq \Ga_q^{-100} \pi_q^{q+\half} (\la_{q+\half}\Ga_q)^N \MM{M,\Nindt,\tau_q^{-1}\Ga_q^{i+15},\Tau_q^{-1}\Ga_q^9} \, . \notag 
\end{align}
Applying \eqref{eq:desert:pressure:def:ER}, Corollary~\ref{lem:agg.pt} with $H=\sigma^{-}_{\divH T_{i,j,k,\xi,\vecl,I,R}^{q+\half+1}}$, $\varpi=\Ga_q^{-100} \pi_q^{q+\half} \mathbf{1}_{\supp a_{\pxi,R}\rhob_\pxi^R \zetab_\xi^I}$ and $p=1$, and \eqref{ind:pi:lower}, we have that \eqref{eq:o.p.4} is verified at level $m=q+\half+1$.  The estimate for $\bmu_{\sigma_{S_O^m}}$ in item~\eqref{i:pc:5:ER} in these cases follows from \eqref{est:mean.DtsiS}, \eqref{i:par:9.5}, and a large choice of $a_*$ in item~\eqref{i:choice:of:a} to ensure that we can gain the advantageous prefactor of  $\max(1,T)^{-1}$.
\smallskip

\noindent\texttt{Case 3b: }pressure increment for \eqref{osc.err.med} and \eqref{osc.err.high} and $\diamond=R$.
We set $\ov \pi = \pi_\ell \Gamma_q^{50} \Lambda_q$ as in the previous case since the low-frequency portion of the error term is identical. Since all the preliminary assumptions in Part 1 are now satisfied, we need to check the additional assumptions from Part 2. In order to do so, we set
\begin{align}
    &N_{**} \textnormal{ as in \eqref{i:par:10}} \,, \quad \NcutLarge, \NcutSmall \textnormal{ as in \eqref{i:par:6}} \, , \quad \Ga=\Ga_q^{\sfrac 12} \, , \quad \delta_{\rm tiny} = \delta^2_{q+3\bn} \, , 
   \quad \mu =\la_{q+\half}\Ga_q \, , \notag\\
   &\mu_0 = \lambda_{q+\half+1} \, , \quad \mu_1 = \lambda_{q+\half+\sfrac 32}\Ga_q^2 \, , \notag\\
    &\mu_{m'} = \lambda_{q+\half+{m'}}\Ga_q^2 \quad \textnormal{if} \quad 2\leq m' \leq \half \, , \notag\\
    &\bm = 1 \quad \textnormal{for the first projector in \eqref{eq:damn:projection} if} \quad  m=q+\half+2 \, , \notag \\
    &\bm = 2 \quad \textnormal{for the second projector in \eqref{eq:damn:projection} if} \quad  m=q+\half+2 \, , \notag \\
    &\bm = m-q-\half \quad \textnormal{if} \quad  m>q+\half+2 \, . \label{eq:desert:choices:1:ER:ER}
\end{align}
Then \eqref{i:st:sample:wut}--\eqref{eq:sample:prop:Ncut:1} hold as in the previous case, \eqref{eq:sample:prop:Ncut:2} holds from \eqref{condi.Ncut0.2}, \eqref{eq:sample:prop:Ncut:3}--\eqref{eq:sample:prop:decoup} hold as in the previous case, \eqref{eq:sample:prop:par:00} holds by definition, \eqref{eq:sample:prop:parameters:0} holds by definition and immediate computation, \eqref{eq:sample:riots:4} holds due to \eqref{ineq:dpot:1}, and \eqref{eq:sample:riot:4:4} holds due to \eqref{ineq:Nstarz:1}.

At this point, we appeal to the conclusions from Part 3 to construct a pressure increment and delineate its properties.  First, from \eqref{d:press:stress:sample}--\eqref{est.S.by.pr.final2} and \eqref{condi.Nfin0}, we have that for $q+\half+2\leq m \leq q+\bn+1$, there exists a pressure increment $\sigma_{ \divH T_{i,j,k,\xi,\vecl,I,R}^{m}}=\sigma_{ \divH T_{i,j,k,\xi,\vecl,I,R}^{m}}^+-\sigma_{ \divH T_{i,j,k,\xi,\vecl,I,
R}^{m}}^-$ such that for $N,M\leq \sfrac{\Nfin}{7}$, 
\begin{align}\label{eq:desert:cowboy:1:ER:ER}
    \left| D^N \Dtq^M \divH T_{i,j,k,\xi,\vecl,I,R}^{m}  \right| \lesssim \left( \sigma_{ \divH t_{i,j,k,\xi,\vecl,I,R}^{m}}^+ + \delta_{q+3\bn}^2 \right) (\min(\lambda_{m},\la_\qbn)\Ga_q)^N \MM{M,\Nindt,\tau_q^{-1}\Ga_q^{i+14},\Tau_q^{-1}\Ga_q^{9}} \, .
\end{align}
From \eqref{eq:inverse:div:linear}, \eqref{est.S.pr.p.support:1}, and \eqref{eq:LP:div:support}, we have that
\begin{align}
    \supp \left( \sigma_{ \divH T_{i,j,k,\xi,\vecl,I,R}^{m}}^+ \right) \subseteq \supp \left( \divH T_{i,j,k,\xi,\vecl,I,R}^{m}  \right) \subseteq \supp \left( a_{\pxi,R} \left(\rhob_\pxi^R \zetab_\xi^I \right)\circ\Phiik \right) \cap B\left( \supp \varrho_{\pxi,R}^I, \lambda_{m-1}^{-1} \right) \, . \label{eq:ocdc:support:ER:ER} 
\end{align}

Now define
\begin{subequations}
\label{eq:desert:pressure:def:ER:ER}
\begin{align}
    \sigma_{S_{O,R}^{m}}^{\pm} &= \sum_{i,j,k,\xi,\vecl,I} \sigma_{ \divH T_{i,j,k,\xi,\vecl,I,R}^{m}}^\pm \quad \textnormal{ if } \quad m\neq q+\bn \, ,\\
    \sigma_{S_{O,R}^{m}}^{\pm} &= \sum_{\td m=\qbn}^{\qbn+1} \sum_{i,j,k,\xi,\vecl,I} \sigma_{ \divH T_{i,j,k,\xi,\vecl,I,R}^{m'}}^\pm \quad \textnormal{ if } \quad m=\qbn \, .
\end{align}
\end{subequations}
Then \eqref{eq:oooooldies} and \eqref{eq:dodging:oldies} give that \eqref{eq:os.p.6} is satisfied for $q+\half+2\leq m \leq \qbn$. 
From \eqref{eq:desert:cowboy:1:ER:ER}, \eqref{eq:desert:cowboy:sum}, \eqref{eq:inductive:partition}, and Corollary~\ref{lem:agg.pt} with 
\begin{align}
&H=\divH T_{i,j,k,\xi,\vecl,I,R}^{m} \, , \qquad \varpi=\left[\sigma_{ \divH T_{i,j,k,\xi,\vecl,I,R}^{m}}^+ + \delta_{q+3\bn}^2\right] \mathbf{1}_{\supp a_{\pxi,R}\rhob_\pxi^R\zetab_\xi^I}  \,  , \qquad p=1 \, , \notag 
\end{align}
we have that for $N,M\leq \sfrac{\Nfin}{7}$,
\begin{align}
    \left| \psi_{i,q} D^N \Dtq^M \sum_{i',j,k,\xi,\vecl,I} \divH T_{i',j,k,\xi,\vecl,I,R}^{m}  \right| &\lesssim \left(  \sigma_{S_{O,R}^{m}}^{+}  + \delta_{q+3\bn}^2 \right) \notag\\
    &\qquad \qquad \times (\min(\lambda_{m},\la_\qbn)\Ga_q)^N \MM{M,\Nindt,\tau_q^{-1}\Ga_q^{i+15},\Tau_q^{-1}\Ga_q^{9}} \, . \label{eq:desert:cowboy:2:ER:ER}
\end{align}
We therefore have that \eqref{eq:o.p.1} is satisfied for $q+\half+2\leq m \leq q+\bn$. From \eqref{est.S.prbypr.pt}, \eqref{condi.Nfin0}, and \eqref{condi.Nindt}, we have that for $N,M\leq \sfrac{\Nfin}{7}$,
\begin{align}
     \left| D^N \Dtq^M \sigma_{ \divH T_{i,j,k,\xi,\vecl,I,R}^{m}}^+ \right| \lesssim \left( \sigma_{ \divH T_{i,j,k,\xi,\vecl,I,R}^{m}}^+ + \delta_{q+3\bn}^2 \right) (\min(\lambda_{m},\la_\qbn)\Ga_q)^N \MM{M,\Nindt,\tau_q^{-1}\Ga_q^{i+15},\Tau_q^{-1}\Ga_q^{9}} \, . \label{eq:desert:cowboy:4:ER:ER}
\end{align}
From \eqref{eq:desert:cowboy:4:ER:ER}, \eqref{eq:desert:cowboy:sum}, \eqref{eq:inductive:partition}, and Corollary~\ref{lem:agg.pt} with 
\begin{align}
&H=\sigma^+_{\divH T_{i,j,k,\xi,\vecl,I,R}^{m}} \, , \qquad \varpi= \left[ H + \delta_{q+3\bn}^2 \right] \mathbf{1}_{\supp a_{\pxi,R}\rhob_\pxi^R\zetab_\xi^I}  \,  , \qquad p=1 \, , \notag 
\end{align}
we have that \eqref{eq:o.p.2} is satisfied for $q+\half+2\leq m \leq \qbn$.

Next, from \eqref{est.S.pr.p}, we have that
\begin{align}
    \left\| \sigma^{\pm}_{\divH T_{i,j,k,\xi,\vecl,I,R}^{m}} \right\|_{\sfrac 32} &\les \left( \left| \supp(\eta_{i,j,k,\xi,\vecl,R}^2(\zetab_\xi^{I,R})^2) \right|^{\sfrac 23} \delta_{q+\bn} \Gamma_{q}^{2j+38}  {\La_q} + \la_\qbn^{-10} \right) {\left(\frac{\min(\la_{m},\la_\qbn)}{\la_{q+\bn}r_q}\right)^{\sfrac23}}  \la_{m-1}^{-2}\la_m \notag \, .
\end{align}
Now from \eqref{eq:desert:pressure:def:ER:ER}, \eqref{ineq:osc:general}, and Corollary~\ref{rem:summing:partition} with $\theta=2$, $\theta_1=0$, $\theta_2=2$, $H=\sigma^{\pm}_{\divH T_{i,j,k,\xi,\vecl,I,\diamond}^{m}}$, and $p=\sfrac 32$, we have that
\begin{align}
    \left\| \psi_{i,q} \sigma_{S_{O,R}^{m}}^{\pm} \right\|_{\sfrac 32} &\les \delta_{m+\bn} \Ga_{m}^{-10} \, . \notag 
\end{align}
Combined with \eqref{eq:o.p.2}, this verifies \eqref{eq:o.p.3} at level $m$.  Arguing now for $p=\infty$ from \eqref{est.S.pr.p}, we have that
\begin{align}
    \left\| \sigma^{\pm}_{\divH T_{i,j,k,\xi,\vecl,I,R}^{m}} \right\|_{\infty} &\les \Ga_q^{\badshaq+40}\La_q {\left(\frac{\min(\la_{m},\la_\qbn)}{\la_{q+\bn}r_q}\right)^{2}}\la_{m}^\alpha \la_{m-1}^{-2}\la_m \notag \, .
\end{align}
Now from \eqref{eq:desert:pressure:def:ER:ER}, \eqref{eq:par:div:2}, and Corollary~\ref{lem:agg.pt} with $H=\sigma^{\pm}_{\divH T_{i,j,k,\xi,\vecl,I,R}^{m}}$, $\varpi=\mathbf{1}_{\supp a_{\pxi,R}\rhob_\pxi^R\zetab_\xi^I}$ and $p=1$, we have that
\begin{align}
    \left\| \psi_{i,q} \sigma_{S_{O,R}^{m}}^{\pm} \right\|_\infty &\les \Ga_q^{\badshaq+40}\La_q {\left(\frac{\min(\la_{m},\la_\qbn)}{\la_{q+\bn}r_q}\right)^{2}}\la_{m}^\alpha \la_{q+\half}^{-1} \leq \Ga_{m}^{\badshaq-100} \, . \notag 
\end{align}
Combined again with \eqref{eq:o.p.2}, this verifies \eqref{eq:o.p.3.inf} at level $m$. 

Finally, from \eqref{est.S.prminus.pt}, \eqref{condi.Nindt}, \eqref{condi.Nfin0}, \eqref{o:5}, \eqref{ind:pi:lower}, and \eqref{eq:ind.pr.anticipated}, we have that for $N,M\leq \sfrac{\Nfin}{7}$,
\begin{align}
    \left| D^N \Dtq^M \sigma^{-}_{\divH T_{i,j,k,\xi,\vecl,I,R}^{m}} \right| &\les {\left(\frac{\min(\la_{m},\la_\qbn)}{\la_{q+\bn}r_q}\right)^{\sfrac23}} \la_{m-1}^{-2} \la_m \pi_\ell \Ga_q^{50} \La_q \notag\\
    &\qquad \qquad \times \min(\la_{m},\la_\qbn)\Ga_q)^N \MM{M,\Nindt,\tau_q^{-1}\Ga_q^{i+15},\Tau_q^{-1}\Ga_q^9} \notag\\
    &\leq \Ga_q^{-100} \pi_q^{q+\half} (\la_{q+\half}\Ga_q)^N \MM{M,\Nindt,\tau_q^{-1}\Ga_q^{i+15},\Tau_q^{-1}\Ga_q^9} \, . \notag 
\end{align}
Applying \eqref{eq:desert:pressure:def:ER:ER}, Corollary~\ref{lem:agg.pt} with $H=\sigma^{-}_{\divH T_{i,j,k,\xi,\vecl,I,R}^{m}}$, $\varpi=\Ga_q^{-100} \pi_q^{q+\half} \mathbf{1}_{\supp a_{\pxi,R}\rhob_\pxi^R \zetab_\xi^I}$ and $p=1$, and \eqref{ind:pi:lower}, we have that \eqref{eq:o.p.4} is verified at levels $q+\half+2\leq m \leq \qbn$. The bounds in item~\eqref{i:pc:5:ER} follow much as in the previous case, and we omit further details.
\smallskip

\noindent\texttt{Case 4: } pressure increment for $\diamond=\vp$. As we noted in the beginning of the proof, the only differences between $\diamond=\vp$ and $\diamond=R$ arise from the redistribution of $r_q^{\sfrac 23}$.  We may therefore define $\sigma_{S_{O,\vp}^m}$ for $q+\half+1\leq m \leq \qbn$ and set
\begin{align}\notag
    \sigma_{S_{O}^{m}}^\pm = \sigma_{S_{O,R}^{m}}^\pm + \sigma_{S_{O,\vp}^{m}}^\pm \, ,
\end{align}
from which \eqref{eq:o.p.1}--\eqref{eq:thurzday:night} follow.
\end{proof}

\begin{lemma*}[\bf Pressure current]\label{lem:oscillation:pressure:current}
For every $m\in\{q+\half+1,\dots,q+\bn\}$, there exists a current error $\phi_{{S_O^{m}}}$ associated to the pressure increment $\si_{S_O^{m}}$ defined by Lemma~\ref{lem:oscillation:pressure} which satisfies the following properties.
\begin{enumerate}[(i)]
    \item\label{i:pc:2:ER} We have the decompositions and equalities
    \begin{subequations}
    \begin{align}\label{eq:desert:decomp:ER}
        \phi_{{S_O^{m}}} &= \phi_{{S_O^{m}}}^* + \sum_{m'=q+\half+1}^{{m}} \phi_{{S_O^{m}}}^{m'} \, , \quad 
         \phi_{S_O^{m}}^{m'} = \phi_{S_O^{m}}^{m',l} + \phi_{S_O^{m}}^{m',*} \, \\
         \div \phi_{{S_O^{m}}}
          &= D_{t,q}\si_{{S_O^{m}}} - \langle D_{t,q}\si_{{S_O^{m}}}\rangle\, . 
    \end{align}
    
    \end{subequations}
    \item\label{i:pc:3:ER} For $q+\half+1 \leq m' \leq m$ and $N,M\leq  2\Nind$,
    \begin{subequations}
    \begin{align}
        &\left|\psi_{i,q} D^N \Dtq^M \phi_{S_O^{m}}^{m',l} \right| < \Ga_{m'}^{-100} \left(\pi_q^{m'}\right)^{\sfrac 32} r_{m'}^{-1} (\la_{m'} \Ga_{m'}^2)^M \MM{M,\Nindt,\tau_q^{-1}\Ga_q^{i+17},\Tau_q^{-1}\Ga_q^9} \label{o:p:c:pt} \\
        &\left\| D^N \Dtq^M \phi_{S_O^{m}}^{m',*} \right\|_\infty + \left\| D^N\Dtq^M \phi_{S_O^{m}}^{*} \right\|_\infty < \Tau_\qbn^{2\Nindt} \delta_{q+3\bn}^{\sfrac 32} (\la_{m}\Ga_{m}^2)^N \tau_q^{-M} \label{e:p:c:nonlocal} \, .
    \end{align}
    \end{subequations}
    \item\label{i:pc:4:ER} For all $q+\half+1\leq m' \leq m$ and all $q+1\leq q' \leq m'-1$, 
    \begin{align}
        B\left( \supp \hat w_{q'}, \sfrac 12 \lambda_{q'}^{-1} \Ga_{q'+1} \right) \cap \supp \left( \phi^{m',l}_{S_O^{m}} \right) = \emptyset \label{o:p:c:supp} \, .
    \end{align}
\end{enumerate}
\end{lemma*}
\begin{proof}
We utilize the case numbering from Lemma~\ref{lem:oscillation:pressure}.  Note that the only cases which require a pressure increments were \texttt{Cases 3a} and \texttt{3b}, which correspond to the analysis of \eqref{osc.err.med0}--\eqref{osc.err.high} and $\diamond=R$, and \texttt{Case 4}, which corresponds to the same terms but with $\diamond=\vp$. We combine the analysis for $\diamond=R$ and $\diamond=\vp$ into a single argument, since as explained in the previous lemmas, the estimates are essentially the same.
\smallskip

\noindent\texttt{Case 3a/4a: }pressure current error from \eqref{osc.err.med0} and $\diamond=R,\vp$.  In this case, we recall from \eqref{eq:desert:choices:1:ER} that we have chosen $\bm=1$ in item~\eqref{i:st:sample:8}, $\mu_0=\lambda_{q+\half+1}\Ga_q^{-1}$, and $\mu_{\bm}=\mu_1=\lambda_{q+\half+1}\Ga_q^2$. We therefore have from \eqref{d:press:stress:sample} that
\begin{align}
    \sigma_{ \divH T_{i,j,k,\xi,\vecl,I,\diamond}^{q+\half+1}} = \sigma_{ \divH T_{i,j,k,\xi,\vecl,I,\diamond}^{q+\half+1}}^+ - \sigma_{ \divH T_{i,j,k,\xi,\vecl,I,\diamond}^{q+\half+1}}^- = \sigma_{ \divH T_{i,j,k,\xi,\vecl,I,\diamond}^{q+\half+1}}^* + \sigma_{ \divH T_{i,j,k,\xi,\vecl,I,\diamond}^{q+\half+1}}^{0} + \sigma_{ \divH T_{i,j,k,\xi,\vecl,I,\diamond}^{q+\half+1}}^{1} \, . \notag
\end{align}
We then define
\begin{align}
    \sigma_{S_O^{q+\half+1}}^* := \sum_{i,j,k,\xi,\vecl,I,\diamond} \sigma_{ \divH T_{i,j,k,\xi,\vecl,I,\diamond}^{q+\half+1}}^* \, , \qquad   \sigma_{S_O^{q+\half+1}}^{q+\half+1} := \sum_{\substack{i,j,k,\xi,\vecl,I,\diamond \\ \bullet=0,1}} \sigma_{ \divH T_{i,j,k,\xi,\vecl,I,\diamond}^{q+\half+1}}^{\bullet} \, , \notag 
\end{align}
so that then using \eqref{d:cur:error:stress:sample}, we may define the current errors
\begin{subequations}
\begin{align}
      \phi_{S_O^{q+\half+1}}^* := \sum_{i,j,k,\xi,\vecl,I,\diamond} \phi_{S_{i,j,k,\xi,\vecl,I,\diamond}^{q+\half+1}}^* &:= \sum_{i,j,k,\xi,\vecl,I,\diamond} (\divH+\divR)\left(\Dtq\sigma_{ \divH T_{i,j,k,\xi,\vecl,I,\diamond}^{q+\half+1}}^*\right) \, , \notag\\   \phi_{S_O^{q+\half+1}}^{q+\half+1} :=\sum_{\substack{i,j,k,\xi,\vecl,I,\diamond \\ \bullet=0,1}} \phi_{S_{i,j,k,\xi,\vecl,I,\diamond}^{q+\half+1}}^{\bullet} &:=  \sum_{\substack{i,j,k,\xi,\vecl,I,\diamond \\ \bullet=0,1}} ( \divH + \divR ) \left( \Dtq \sigma_{ \divH T_{i,j,k,\xi,\vecl,I,\diamond}^{q+\half+1}}^{\bullet} \right) \notag \\
      &= \underbrace{\phi_{S_O^{m}}^{q+\half+1,l}}_{\textnormal{all the $\divH$ terms}} + \underbrace{\phi_{S_O^{m}}^{q+\half+1,*}}_{\textnormal{all the $\divR$ terms}}\, ,  \notag
\end{align}
\end{subequations}
which satisfy
    \begin{align*}
        \div \phi_{{S_O^{q+\half+1}}}^* &= \Dtq \sigma_{S_O^{q+\half+1}}^* - \int_{\T^3} \Dtq \sigma_{S_O^{q+\half+1}}^*(t,x') \,dx' \, , \\
         \div \phi_{{S_O^{q+\half+1}}}^{q+\half+1} &= \Dtq \sigma_{S_O^{{q+\half+1}}}^{q+\half+1} - \int_{\T^3} \Dtq \sigma_{S_O^{q+\half+1}}^{q+\half+1}(t,x') \,dx' \, .
    \end{align*}
We decompose the current error further into $\phi_{S_O^{q+\half+1}}^{q+\half+1} = \phi_{S_O^{q+\half+1}}^{q+\half+1,l} + \phi_{S_O^{q+\half+1}}^{q+\half+1,*}$ using item~\ref{sample2.item3}.

In order to check \eqref{o:p:c:pt}, we recall the parameter choices from \texttt{Case 3a} of Lemma~\ref{lem:oscillation:general:estimate} and the choice of $\ov \pi = \pi_\ell \Ga_q^{50}\La_q$ from Lemma~\ref{lem:oscillation:pressure} apply Part 4 of Proposition~\ref{lem.pr.invdiv2}, specifically \eqref{est.S.by.pr.final3}.  We then have from \eqref{condi.Nfin0} that for each $i,j,k,\xi,\vecl,I,\diamond,\bullet$ and $M,N\leq 2\Nind$ (after appending a superscript $l$ to refer to the local portion),
\begin{align}
    \left| D^N \Dtq^M \phi_{S_{i,j,k,\xi,\vecl,I,\diamond}^{q+\half+1}}^{\bullet,l} \right| &\leq {\tau_q^{-1} \Ga_q^{i+70}} {\pi_\ell} \La_q {\left(\frac{\la_{q+\half+1}}{\la_{q+\bn}r_q}\right)^{2}} \la_{q+\half}^{-1} \notag\\
    &\qquad \qquad \times (\la_{q+\half+1}\Ga_q)^N \MM{M,\Nindt-\NcutSmall-1,\tau_{q}^{-1}\Ga_q^{i+14},\Tau_q^{-1}\Ga_q^9} \, . \label{eq:desert:parsing:ER}
\end{align}
Next, from \eqref{est.S.pr.p.support:2}, we have that
\begin{align}
    \supp \left( \phi_{S_{i,j,k,\xi,\vecl,I,\diamond}^{q+\half+1}}^{\bullet,l} \right) &\subseteq B\left( \divH T^{q+\half+1}_{i,j,k,\xi,\vecl,I,\diamond}, 2\lambda_{q+\half+1}\Ga_q^{-1} \right) \notag\\
    &\subseteq B\left( \supp \left( a_{\pxi,\diamond}(\varrho_\pxi^\diamond\zetab_\xi^I)\circ\Phiik \right), 2\lambda_{q+\half+1}\Ga_q^{-1} \right) \,. \notag
\end{align}
Then applying \eqref{eq:oooooldies}, we have that \eqref{o:p:c:supp} is verified for $m=m'=q+\half+1$. Returning to the proof of \eqref{o:p:c:pt}, we can now apply Corollary~\ref{lem:agg.Dtq} with
\begin{align*}
    H= \phi_{S_{i,j,k,\xi,\vecl,I,\diamond}^{q+\half+1}}^{\bullet,l} \, , \qquad \varpi = \Ga_q^{70} {\pi_\ell} \La_q {\left(\frac{\la_{q+\half+1}}{\la_{q+\bn}r_q}\right)^{2}} \la_{q+\half}^{-1} \, .
\end{align*}
From \eqref{eq:aggDtq:conc:1}, \eqref{condi.Nindt}, \eqref{ind:pi:lower}, \eqref{eq:ind.pr.anticipated}, \eqref{ineq:r's:eat:Gammas}, and \eqref{ineq:in:the:morning:ER}, we have that
\begin{align}
    &\left| \psi_{i,q} \sum_{i',j,k,\xi,\vecl,I,R,\bullet} \divH \left( \Dtq \sigma_{ \divH t_{i,j,k,\xi,\vecl,I,R}^{q+\half+1}}^{\bullet} \right) \right|\notag\\
    &\quad \underset{\eqref{eq:aggDtq:conc:1}}{\lesssim} \underbrace{r_q^{-1} \la_q \left( \pi_q^q \right)^{\sfrac 12}}_{\textnormal{cost of $\Dtq$}} \underbrace{\pi_\ell}_{\substack{\textnormal{dominates} \\ \textnormal{low-freq. coeff's}}} \underbrace{\La_q \la_{q+\half}^{-1} }_{\textnormal{freq. gain}} \underbrace{\Ga_q^{76}}_{\textnormal{lower order}} \underbrace{\left( \frac{\la_{q+\half+1}\Ga_q}{\la_{q+\half}} \right)^{2}}_{\textnormal{intermittency losses}} \underbrace{\la_{q+\half}^{-1}}_{\textnormal{inv. div. gain}} \notag\\
    &\quad \qquad \times (\la_{q+\half+1}\Ga_q)^N \MM{M,\Nindt-\NcutSmall-1,\tau_{q}^{-1}\Ga_q^{i+15},\Tau_q^{-1}\Ga_q^9} \notag\\
    &\quad \underset{\eqref{ind:pi:lower},\eqref{eq:ind.pr.anticipated}}{\les} r_q^{-1} \Ga_q^{100} \left(\pi_q^{q+\half+1} \frac{\delta_{q+\bn}}{\delta_{q+\half+1+\bn}}\right)^{\sfrac 32} \La_q^2 \left( \frac{\la_{q+\half+1}\Ga_q}{\la_{q+\half}} \right)^{2} \la_{q+\half}^{-2} \notag\\
    &\qquad \qquad \times (\la_{q+\half+1}\Ga_q)^N \MM{M,\Nindt-\NcutSmall-1,\tau_{q}^{-1}\Ga_q^{i+15},\Tau_q^{-1}\Ga_q^9} \notag\\
    &\quad \underset{\eqref{condi.Nindt},\eqref{ineq:in:the:morning:ER},\eqref{ineq:r's:eat:Gammas}}{\leq} \Ga_q^{-150} r_{q+\half+1}^{-1} \left( \pi_q^{q+\half+1} \right)^{\sfrac 32} (\la_{q+\half+1}\Ga_q)^N \MM{M,\Nindt,\tau_{q}^{-1}\Ga_q^{i+16},\Tau_q^{-1}\Ga_q^9} \label{eq:desert:parsing:2:ER} 
\end{align}
for $N,M\leq 2\Nind$ from \eqref{condi.Nfin0}, which verifies \eqref{o:p:c:pt} at level $q+\half+1$. In order to achieve \eqref{e:p:c:nonlocal}, we appeal to \eqref{est.S.by.pr.final4}--\eqref{est.S.by.pr.final.star}, the choice of $K_\circ$ in item~\eqref{i:par:9.5}, \eqref{condi.Nfin0}, and an aggregation quite similar to previous nonlocal aggregations.
\smallskip

\noindent\texttt{Case 3b/4b: } pressure current error from~\eqref{osc.err.med} and~\eqref{osc.err.high} and $\diamond=R,\vp$.  In this case we consider the higher shells from the oscillation error.  The general principle is that the estimate will only be sharp in the $m=m'=\qbn$ double endpoint case, for which the intermittency loss is most severe.  We now explain why this is the case by parsing estimates \eqref{eq:desert:parsing:ER} and \eqref{eq:desert:parsing:2:ER}. We incur a material derivative cost of $\tau_q^{-1}\Ga_q^{i+70}$, which is converted into $r_q^{-1}\lambda_q (\pi_q^q)^{\sfrac 12}$ using \eqref{eq:psi:q:q'} and the rough definition of $\tau_q^{-1}=\delta_q^{\sfrac 12}\lambda_q r_q^{-\sfrac 13}$, or equivalently Corollary~\ref{lem:agg.Dtq}. The $L^{\sfrac 32}$ size of the high-frequency coefficients from the oscillation error is $(\la_m \la_{q+\half}^{-1})^{\sfrac 23}$; this encodes the intermittency loss from $L^1$ to $L^{\sfrac 32}$ of a squared, $\leq \lambda_m$ frequency projected, $L^2$ normalized pipe flow with minimum frequency $\la_{q+\half}$ -- see also the choices of $\const_{*,\sfrac 32}$ from Lemma~\ref{lem:oscillation:general:estimate}.  This accounts for $\sfrac 23$ of the squared power in the intermittency losses. The low-frequency coefficient function from a quadratic oscillation error incurs a derivative cost of $\Lambda_q$ (which we have grouped with ``frequency gain") and is dominated by $\pi_\ell$. The negative power in the frequency gain will be $\lambda_m$ and is determined by which shell (indexed by $m$) of the oscillation error is being considered.  The lower order terms may be ignored.  Next, we have an $L^{\sfrac 32}\rightarrow L^\infty$ intermittency loss of $(\la_{m'}\la_{q+\half}^{-1})^{\sfrac 43}$, which accounts for $\sfrac 43$ of the power in the intermittency losses and is used to pointwise dominate the high-frequency portion (at frequency $\la_{m'}$ due to the frequency projector) of the pressure increment using the $L^{\sfrac 32}$ norm.  By simply pointwise dominating the high-frequency portion of the pressure increment, using this to compute the $L^1$ norm of the resulting current error, and showing that the result is dominated by existing pressure, we prevent a loop of new current error and new pressure creation. Finally, we have an inverse divergence gain depending on which synthetic Littlewood-Paley shell of the pressure increment we are considering. The net effect is that the $\La_q$ from ``frequency gain" and the $\la_{m'}^{-1}$ from ``inv. div. gain" upgrade the $\pi_\ell^{\sfrac 32}$ to $(\pi_q^{m'})^{\sfrac 32}$, and the remaining $\la_q\la_{m}^{-1}$ from the $\Dtq$ cost and the frequency gain is strong enough to absorb the intermittency loss since $m'\leq m$, with a perfect balance in the case
$$  m=m'=\qbn \qquad \implies \qquad \left(\frac{\la_{\qbn}}{\la_{q+\half}}\right)^2  \la_q \la_{\qbn}^{-1} \approx 1 \, . $$

In order to fill in the details, we now recall the choices of $\bm$ and $\mu_{m'}$ from \eqref{eq:desert:choices:1:ER:ER}.  For the sake of brevity we ignore the slight variation in the case of the first projector for $m=q+\half+2$ and focus on the second projector for $m=q+\half+2$ and the other cases $q+\half+2<m\leq \qbn+1$. We have from \eqref{d:press:stress:sample} that
\begin{align}
    \sigma_{ \divH T_{i,j,k,\xi,\vecl,I,\diamond}^{m}} = \sigma_{ \divH T_{i,j,k,\xi,\vecl,I,\diamond}^{m}}^+ - \sigma_{ \divH T_{i,j,k,\xi,\vecl,I,\diamond}^{m}}^- = \sigma_{ \divH T_{i,j,k,\xi,\vecl,I,\diamond}^{m}}^* + \sum_{\iota=0}^{m-q-\half} \sigma_{ \divH T_{i,j,k,\xi,\vecl,I,\diamond}^{m}}^{\iota} \, . \notag
\end{align}
We then define the frequency-projected pressure increments by
\begin{align}
   \sigma_{S_O^{m}}^* &= \sum_{i,j,k,\xi,\vecl,I,\diamond} \sigma_{ \divH T_{i,j,k,\xi,\vecl,I,\diamond}^{m}}^* \, , \qquad   \sigma_{S_O^{m}}^{q+\half+1} = \sum_{{i,j,k,\xi,\vecl,I,\diamond}} \sigma_{ \divH T_{i,j,k,\xi,\vecl,I,\diamond}^{m}}^{0} \, , \notag \\
   \sigma_{S_O^{m}}^{q+\half+2} &= \sum_{\substack{i,j,k,\xi,\vecl,I,\diamond \\ \iota=1,2}} \sigma_{ \divH T_{i,j,k,\xi,\vecl,I,\diamond}^{m}}^{\iota} \, , \notag\\
   \sigma_{S_O^{m}}^{q+\half+m'} &= \sum_{\substack{i,j,k,\xi,\vecl,I,\diamond \\ \iota=m'}} \sigma_{ \divH T_{i,j,k,\xi,\vecl,I,\diamond}^{m}}^{\iota} \quad \textnormal{if $q+\half+m' = q+\half+\iota \leq m \leq \qbn-1$} \, , \\
   \sigma_{S_O^m}^\qbn &= \sum_{\substack{i,j,k,\xi,\vecl,I,\diamond \\ \iota=}} \quad \textnormal{if $\iota$ $m=\qbn, \qbn+1$} \, . \notag 
\end{align}
Using \eqref{d:cur:error:stress:sample}, we may define the current errors
\begin{align}
   \phi_{S_O^{m}}^* &= \sum_{i,j,k,\xi,\vecl,I,R} (\divH + \divR ) \left( \Dtq \sigma_{ \divH T_{i,j,k,\xi,\vecl,I,\diamond}^{m}}^* \right) \, , \quad   \phi_{S_O^{m}}^{q+\half+1} = \sum_{{i,j,k,\xi,\vecl,I,\diamond}} (\divH + \divR) \left( \Dtq \sigma_{ \divH T_{i,j,k,\xi,\vecl,I,\diamond}^{m}}^{0} \right) \, , \notag \\
   \phi_{S_O^{m}}^{q+\half+2} &= \sum_{\substack{i,j,k,\xi,\vecl,I,\diamond \\ \iota=1,2}} (\divH + \divR) \left( \Dtq \sigma_{ \divH T_{i,j,k,\xi,\vecl,I,\diamond}^{m}}^{\iota} \right) \, , \notag\\
   \phi_{S_O^{m}}^{q+\half+m'} &= \sum_{\substack{i,j,k,\xi,\vecl,I,\diamond \\ \iota=m'}} (\divH + \divR) \left( \Dtq \sigma_{ \divH T_{i,j,k,\xi,\vecl,I,\diamond}^{m}}^{\iota} \right) \quad \textnormal{if $q+\half+m' = q+\half+\iota < m$} \, , \notag \\
    \phi_{S_O^{m}}^{\qbn} &= \sum_{\substack{i,j,k,\xi,\vecl,I,\diamond \\ \iota=m-q-\half,m-q-\half+1}} (\divH + \divR) \left( \Dtq \sigma_{ \divH T_{i,j,k,\xi,\vecl,I,\diamond}^{m}}^{\iota} \right) \, . \notag 
\end{align}
As in the previous case, we may append superscripts of $l$ and $*$ for $q+\half+1\leq m \leq q+\bn$ corresponding to the $\divH$ and $\divR$ portions, respectively. We have thus verified item~\eqref{i:pc:2:ER} immediately from these definitions and from \eqref{d:cur:error:stress:sample} and item~\eqref{sample2.item3}. In order to check \eqref{o:p:c:pt}, we define the temporary notation $m'(\iota)$ to make a correspondence between the value of $\iota$ above and the superscript on the left-hand side, which determines which bin the current errors go into.  Specifically, we set $m'(0)=1$, $m'(1)=m'(2)=2$, $m'(\iota)=\iota$ if $q+\half+\iota < m$, and $m'(m-q-\half)=m'(m-q-\half+1)=m-q-\half$. Then from Part 4 of Proposition~\ref{lem.pr.invdiv2}, specifically \eqref{est.S.by.pr.final3}, and \eqref{condi.Nfin0}, we have that for each $i,j,k,\xi,\vecl,I,\diamond,\iota$ and $M,N\leq 2\Nind$,
\begin{align}
    &\left| D^N \Dtq^M (\divH + \divR) \left( \Dtq \sigma^\iota_{\divH T^{m}_{i,j,k,\xi,\vecl,I,\diamond}} \right) \right| \notag\\
    & \leq \tau_q^{-1} \Ga_q^{i+70} \pi_\ell \La_q \left( \frac{\min(\la_m,\la_\qbn)}{\la_{q+\half}} \right)^{\sfrac 23} \la_{m-1}^{-2} \la_{m} \left( \frac{\min(\la_{q+\half+m'(\iota)},\la_\qbn)\Ga_q}{\la_{q+\half}} \right)^{\sfrac 43}  \notag\\
    &\quad \times \la_{q+\half+m'(\iota)-1}^{-2} \la_{q+\half+m'(\iota)} \left(\min(\la_{q+\half+m'(\iota)},\la_{m})\Ga_q\right)^N \MM{M,\Nindt-\NcutSmall-1,\tau_{q}^{-1}\Ga_q^{i+14},\Tau_q^{-1}\Ga_q^9} \, . \notag 
\end{align}
Next, from \eqref{est.S.pr.p.support:2} and the fact that $q+\half+m'(\iota)\leq m$, we have that
\begin{align}
    \supp \left( \divH \left( \Dtq \sigma^\iota_{\divH T^{m}_{i,j,k,\xi,\vecl,I,\diamond}} \right) \right) &\subseteq B\left( \divH T^{m}_{i,j,k,\xi,\vecl,I,\diamond}, 2\lambda_{q+\half+m'(\iota)-1}\Ga_q^{-2} \right) \notag\\
    &\subseteq B\left( \supp \left( a_{\pxi,\diamond}(\varrho_\pxi^\diamond\zetab_\xi^I)\circ\Phiik \rho_{\pxi,\diamond}^I  \right), \lambda_{m-1}^{-1}+2\lambda_{q+\half+m'(\iota)-1}\Ga_q^{-2} \right) \notag\\
     &\subseteq B\left( \supp \left( a_{\pxi,\diamond}(\varrho_\pxi^\diamond\zetab_\xi^I)\circ\Phiik \rho_{\pxi,\diamond}^I  \right), 2\lambda_{q+\half+m'(\iota)-1} \right)\,. \notag
\end{align}
Then applying \eqref{eq:oooooldies}, we have that \eqref{o:p:c:supp} is verified for $m'=q+\half+m'(\iota)$. Returning to the proof of \eqref{o:p:c:pt}, we can now apply Corollary~\ref{lem:agg.Dtq} with
\begin{align*}
    H&= \divH \left( \Dtq \sigma^\iota_{\divH T^{m}_{i,j,k,\xi,\vecl,I,\diamond}} \right) \, , \notag\\
    \varpi &= \Ga_q^{70} \pi_\ell \La_q \left( \frac{\min(\la_m,\la_\qbn)}{\la_{q+\half}} \right)^{\sfrac 23} \la_{m-1}^{-2} \la_{m} \left( \frac{\min(\la_{q+\half+m'(\iota)},\la_\qbn)\Ga_q}{\la_{q+\half}} \right)^{\sfrac 43}  \la_{q+\half+m'(\iota)-1}^{-2} \la_{q+\half+m'(\iota)} \, .
\end{align*}
From \eqref{ineq:r's:eat:Gammas}, \eqref{eq:aggDtq:conc:1}, \eqref{ind:pi:lower}, \eqref{eq:ind.pr.anticipated}, and \eqref{ineq:in:the:afternoon:ER}, we have that
\begin{align}
    &\left| \psi_{i,q} \sum_{i',j,k,\xi,\vecl,I,\diamond} \divH \left( \Dtq \sigma_{ \divH T^{m}_{i,j,k,\xi,\vecl,I,\diamond}}^{\iota} \right) \right|\notag\\
    &\quad \underset{\eqref{eq:aggDtq:conc:1}}{\lesssim} \Ga_q^{76} r_q^{-1} \la_q \left( \pi_q^q \right)^{\sfrac 12} \pi_\ell \La_q \left( \frac{\min(\la_m,\la_\qbn)}{\la_{q+\half}} \right)^{\sfrac 23} \la_{m-1}^{-2} \la_{m} \left( \frac{\min(\la_{q+\half+m'(\iota)},\la_\qbn)\Ga_q}{\la_{q+\half}} \right)^{\sfrac 43}  \notag\\
    &\qquad \times \la_{q+\half+m'(\iota)-1}^{-2} \la_{q+\half+m'(\iota)} \left(\min(\la_{q+\half+m'(\iota)},\la_{m})\Ga_q\right)^N \MM{M,\Nindt,\tau_{q}^{-1}\Ga_q^{i+16},\Tau_q^{-1}\Ga_q^9} \notag\\
    &\quad \underset{\eqref{ind:pi:lower},\eqref{eq:ind.pr.anticipated}}{\les} \Ga_q^{76} r_q^{-1} \la_q \left(\pi_q^{q+\half+m'(\iota)} \frac{\delta_{q+\bn}}{\delta_{q+\half+m'(\iota)+\bn}} \right)^{\sfrac 32} \La_q \left( \frac{\min(\la_m,\la_\qbn)}{\la_{q+\half}} \right)^{\sfrac 23} \notag\\
    &\qquad \times \la_{m-1}^{-2} \la_{m} \left( \frac{\min(\la_{q+\half+m'(\iota)},\la_\qbn)\Ga_q}{\la_{q+\half}} \right)^{\sfrac 43}  \la_{q+\half+m'(\iota)-1}^{-2} \la_{q+\half+m'(\iota)} \notag\\
    &\qquad \qquad \times \left(\min(\la_{q+\half+m'(\iota)},\la_{m'})\Ga_q\right)^N \MM{M,\Nindt,\tau_{q}^{-1}\Ga_q^{i+16},\Tau_q^{-1}\Ga_q^9} \notag\\
    &\quad \underset{\eqref{ineq:in:the:afternoon:ER},\eqref{ineq:r's:eat:Gammas}}{\leq} \Ga_{m'}^{-150} r_{m'}^{-1} \left( \pi_q^{q+\half+m'(\iota)} \right)^{\sfrac 32} \left(\min(\la_{q+\half+m'(\iota)},\la_{m})\Ga_q\right)^N \MM{M,\Nindt,\tau_{q}^{-1}\Ga_q^{i+16},\Tau_q^{-1}\Ga_q^9} \, , \notag 
\end{align}
for  $N,M\leq 2\Nind$ from \eqref{condi.Nfin0}, which verifies \eqref{o:p:c:pt} at level $m'$. In order to achieve \eqref{e:p:c:nonlocal}, we appeal to \eqref{est.S.by.pr.final4}--\eqref{est.S.by.pr.final.star}, the choice of $K_\circ$ in item~\ref{i:par:9.5}, and \eqref{condi.Nfin0}.
\end{proof}

\subsection{Transport and Nash stress errors \texorpdfstring{$S_{TN}$}{esstee}}\label{ss:ER:TN}

\begin{lemma}[\bf Applying inverse divergence]\label{l:transport:error}
There exist symmetric stresses $S_{TN}=S_{TN}^l+S_{TN}^*$ which satisfy the following.\index{$S_{TN}$}
\begin{enumerate}[(i)]
    \item For all $N,M\leq \sfrac{\Nfin}{10}$, the local part $S_{TN}^l$ satisfies
\begin{subequations}
\begin{align}
\left\| \psi_{i,q} D^N \Dtq^M S_{TN}^l \right\|_{{\sfrac 32}} &\les \Gamma_{q+\bn}^{-100} \delta_{q+2\bn} \lambda_{q+\bn}^N \MM{M,\Nindt, \tau_q^{-1}\Gamma_{q}^{i + 15},\Tau_q^{-1}\Ga_q^9} \label{eq:trans:L1:est} \\
\left\| \psi_{i,q} D^N \Dtq^M S_{TN}^l \right\|_{\infty}  &\les \Gamma_{q+\bn}^{\badshaq-100} \lambda_{q+\bn}^N \MM{M,\Nindt, \tau_q^{-1}\Gamma_{q}^{i+15},\Tau_q^{-1}\Ga_q^9} \label{eq:trans:Loo:est} \, .
\end{align}
\end{subequations}
Furthermore, we have that
\begin{subequations}
\begin{align}
   B \left( \supp \hat w_{q'}, \lambda_{q'}^{-1} \Gamma_{q'+1} \right) \cap \supp S_{TN}^l = \emptyset \label{eq:trans:supp}
\end{align}
\end{subequations}
for all $q+1\leq q' \leq q+\bn-1$.
\item For $N,M\leq 2\Nind$ the nonlocal part satisfies
\begin{equation}\label{eq:trans:nonlocal}
    \left\| D^N \Dtq^M S_{TN}^{*} \right\|_{\infty} \leq  \Tau_\qbn^{4\Nindt} \delta_{q+3\bn}^2 \lambda_{q+\bn}^N \tau_q^{-M} \, .
\end{equation}
\end{enumerate}
\end{lemma}
\begin{remark*}[\bf Abstract formulation of the transport and Nash stress errors]\label{rem.ct.tn}
For the purposes of analyzing the transport and Nash current errors in subsection~\ref{op:tnce} and streamlining the creation of pressure increments, it will again be useful to abstract the properties of these error terms. We will prove every one of the following claims in the course of of proving Lemma~\ref{l:transport:error}. First, there exists a $q$-independent constant $\const_{\mathcal{H}}$ such that 
\begin{align}\label{trans.loc.pot}
    S^{l}_{TN} = \sum_{i,j,k,\xi,\vecl,I,\diamond} \sum_{j'=0}^{\const_{\mathcal{H}}} H_{i,j,k,\xi,\vecl,I,\diamond}^{\alpha{(j')}} \rho_{i,j,k,\xi,\vecl,I,\diamond}^{\beta{(j')}} \circ \Phi_{(i,k)} \, .
\end{align}
Next, the functions $H$ and $\rho$ (with subscripts and superscripts suppressed for convenience) defined above satisfy the following.
\begin{enumerate}[(i)]
    \item $H$ satisfies
    \begin{align}
        \left| D^N \Dtq^M H \right| &\lec \pi_\ell \La_q \la_{q+\half}^N \MM{M, \Nindt, \tau_{q}^{-1}\Gamma_{q}^{i+14}, \Tau_{q}^{-1}\Ga_q^8} \label{eq:H:eckel:TN}
    \end{align}
    for all $N, M \leq \sfrac{\Nfin}{10}$. 
    \item We have that
    \begin{align}\label{eq:H:eckel:trans}
        \supp H \subseteq \supp\left(\eta_{i,j,k,\xi,\vecl,\diamond} \etab^{I,\diamond}_{\xi}\right) \, .
    \end{align}
    \item For $\dpot$ as in \eqref{i:par:10}, there exist a tensor potential $\vartheta$ (we suppress the indices at the moment for convenience) such that $\rho=\partial_{i_1 \dots i_{\dpot}}\vartheta^{(i_1,\dots,i_\dpot)}$. Furthermore, $\vartheta$ is $(\T / \la_{q+\half}\Ga_q )^3$-periodic and satisfies the estimates
    \begin{align}\label{rho:eckel:estimate:TN}
        \norm{D^N \partial_{i_1}\dots \partial_{i_k} \vartheta^{(i_1,\dots, i_\dpot)}}_{L^p} \les r_q^{\sfrac 2p -2} \la_{q+\bn}^{-1+N+k-\dpot} \, .
    \end{align}
   for $p=\sfrac 32,\infty$, all $N\leq \sfrac{\Nfin}{5}$, and $0\leq k\leq \dpot$.
   \item We have that
   \begin{align}\label{rho:eckel:support:TN}
    \supp \left( H \rho\circ \Phi \right) \cap B\left( \supp \hat w_{q'}, \la_{q'}^{-1} \Ga_{q'+1} \right) = \emptyset
    \end{align}
    for all $q+1 \leq q' \leq \qbn-1$. 
\end{enumerate}
\end{remark*}

\begin{proof}[Proof of Lemma~\ref{l:transport:error}]
We start by considering either a Reynolds or current corrector defined in subsection~\ref{ss:corr-ec-tor} and expanding
\begin{align}
    \Dtq w_{q+1,\diamond} &=\Dtq \biggl(\sum_{i,j,k,\xi,\vecl,I} \curl \left( a_{(\xi),\diamond} (\rhob_{(\xi)}^\diamond \zetab^{I,\diamond}_\xi)\circ \Phiik \nabla\Phi_{(i,k)}^{T}\UU_{(\xi),\diamond}^{I} \circ \Phi_{(i,k)} \right) \biggr) \nonumber\\
    &= \sum_{i,j,k,\xi,\vecl,I} \Dtq \left(a_{(\xi),\diamond} \nabla\Phi_{(i,k)}^{-1}\right) (\rhob_{(\xi)}^\diamond \etab^{I,\diamond}_\xi)\circ \Phiik \WW_{(\xi),\diamond}^{I} \circ \Phi_{(i,k)} \notag \\ 
    &\qquad + \sum_{i,j,k,\xi,\vecl,I} \Dtq \nabla \left( (\rhob_{(\xi)}^\diamond \zetab^{I,\diamond}_\xi)\circ \Phiik a_{(\xi),\diamond} \right) \times \left( \nabla\Phiik \UU_{(\xi),\diamond}^{I}\circ \Phiik \right) \notag\\
    &\qquad + \sum_{i,j,k,\xi,\vecl,I} \nabla \left((\rhob_{(\xi)}^\diamond \zetab^{I,\diamond}_\xi)\circ \Phiik a_{(\xi),\diamond} \right) \times \left( \Dtq \nabla\Phiik \UU_{(\xi),\diamond}^{I} \circ \Phiik \right) \label{eq:transport:estimate:1}
\end{align}
and 
\begin{align}
    w_{q+1,\diamond} \cdot \nabla \hat u_q &= \sum_{i,j,k,\xi,\vecl,I} \curl \left( a_{(\xi),\diamond} (\rhob_{(\xi)}^\diamond \zetab^{I,\diamond}_\xi)\circ \Phiik \nabla\Phi_{(i,k)}^{T} \UU_{(\xi),\diamond}^{I} \circ \Phi_{(i,k)} \right) \cdot \nabla \hat u_q \nonumber\\
    &= \sum_{i,j,k,\xi,\vecl,I} \left( a_{(\xi),\diamond} \nabla\Phi_{(i,k)}^{-1}  (\rhob_{(\xi)}^\diamond \zetab^{I,\diamond}_\xi)\circ \Phiik \WW_{(\xi),\diamond}^{I} \circ \Phi_{(i,k)} \right) \cdot \nabla \hat u_q \notag \\ 
    &\qquad + \sum_{i,j,k,\xi,\vecl,I} \left( \nabla \left( (\rhob_{(\xi)}^\diamond \zetab^{I,\diamond}_\xi)\circ \Phiik a_{(\xi),\diamond} \right) \times \left( \nabla\Phiik \UU_{(\xi),\diamond}^{I} \circ \Phiik \right) \right) \cdot \nabla \hat u_q  \, . \label{eq:Nash:estimate:1}
\end{align}
We shall only consider the worst terms, which are the ones containing $\WW_{(\xi),\diamond}^{I}$.  Since $\Dtq w_{q+1,\diamond}$ and $w_{q+1,\diamond} \cdot \nabla \hat u_q$ are mean-zero (see the argument below the display in \eqref{ER:new:error}), we can apply $\divH$ and $\divR$ from Proposition~\ref{prop:intermittent:inverse:div} to each term in \eqref{eq:transport:estimate:1} while ignoring the last term in \eqref{eq:inverse:div:error:stress}. 

We now fix values of $i$, $j$, $k$, $\xi$, $\vecl$, $I$, and $\diamond$ so that we are simply considering 
\begin{align}\label{eq:trans:simplified}
    T_{i,j,k,\xi,\vecl,I,\diamond} &:= \Dtq \left(a_{(\xi),\diamond }\nabla\Phi_{(i,k)}^{-1}\right) (\chib^\diamond_{(\xi)}\etab^{I,\diamond}_\xi)\circ \Phiik \WW_{(\xi),\diamond}^I \circ \Phi_{(i,k)}\\
    &\qquad + \na \hat u_q \cdot \left(a_{(\xi),\diamond }\nabla\Phi_{(i,k)}^{-1}\right) (\chib^\diamond_{(\xi)}\etab^{I,\diamond}_\xi)\circ \Phiik \WW_{(\xi),\diamond}^I \circ \Phi_{(i,k)}\, .\notag
\end{align}
We apply Proposition~\ref{prop:intermittent:inverse:div} along with Remark~\ref{rem:pointwise:inverse:div} with the following choices.  Let $p\in\{\sfrac 32,\infty\}$.  We set $v=\hat u_q$, and $D_t=\Dtq=\partial_t+\hat u_q\cdot\nabla$.  In order to verify the low-frequency assumptions from Part 1 of Proposition~\ref{prop:intermittent:inverse:div} and Remark~\ref{rem:pointwise:inverse:div}, we set
\begin{align}
    &G_{i,j,k,\xi,\vecl,I,R} = r_q \left[\Dtq \left(a_{(\xi),R}\nabla\Phi_{(i,k)}^{-1}\right) ( \chib_{(\xi)}^{R} \zetab_\xi^{I,R})\circ\Phiik \xi  + \na \hat u_q \cdot \left(a_{(\xi),R}\nabla\Phi_{(i,k)}^{-1}\right) ( \chib_{(\xi)}^{R} \zetab_\xi^{I,R})\circ\Phiik \xi \right] \, , \notag \\
    &G_{i,j,k,\xi,\vecl,I,\vp} = r_q^{\sfrac 43} \left[\Dtq \left(a_{(\xi),\vp}\nabla\Phi_{(i,k)}^{-1}\right) ( \chib_{(\xi)}^{\vp} \zetab_\xi^{I,\vp})\circ\Phiik \xi  + \na \hat u_q \cdot \left(a_{(\xi),\vp}\nabla\Phi_{(i,k)}^{-1}\right) ( \chib_{(\xi)}^{\vp} \zetab_\xi^{I,\vp})\circ\Phiik \xi \right] \, , \notag  \\
    &N_*={\Nfin/4} \, , \quad M_*={\Nfin/5} \, , \quad \const_{G,\sfrac 32} = r_q \left| \supp(\eta_{i,j,k,\xi,\vecl,\diamond}\zetab_\xi^{I,\diamond}) \right|^{\sfrac 23} \delta_{q+\bn}^{\sfrac{1}{2}} \Gamma_{q}^{i+j+20} \tau_q^{-1} + r_q \la_\qbn^{-10} \, , \notag \\
    &\const_{G,\infty} = \La_q \Gamma_q^{2+\badshaq} \, , \quad \lambda=\lambda_{q+\lfloor \sfrac \bn 2 \rfloor} \, , \quad \nu=\tau_q^{-1}\Gamma_{q}^{i+13} \, , \quad M_t=\Nindt\, , \quad \nu'=\Tau_{q}^{-1}\Ga_q^8 \, ,\notag \\
    &v =\hat u_q \, , \quad \Phi = \Phiik \, , \quad D_t = \Dtq \, , \quad \lambda' = \La_q \, , \quad \pi = \pi_\ell \La_q \, .\label{eq:eckel:TN:choices}
\end{align}
Then we have that \eqref{eq:inv:div:NM} is satisfied by definition, and \eqref{eq:DDpsi2}--\eqref{eq:DDv} are satisfied as in the proof of Lemma~\ref{lem:oscillation:general:estimate}. In order to check \eqref{eq:inverse:div:DN:G}, we appeal to Lemma~\ref{lem:a_master_est_p}, estimate \eqref{eq:Lagrangian:Jacobian:6} for $(\nabla\Phiik)^{-1}$, estimate~\eqref{eq:checkerboard:derivatives:check} from Lemma~\ref{lem:finer:checkerboard:estimates} to estimate $\zetab_\xi^{I,\diamond}\circ\Phiik$, Proposition~\ref{prop:bundling}, and \eqref{eq:nasty:D:vq:old}.  Specifically, we have that for all $N,M \leq 9\Nind$,
\begin{align}
\Biggl\| D^N \Dtq^M G_{i,j,k,\xi,\vecl,I,\diamond} \Biggr\|_{\sfrac 32} &\lessg \const_{G,\sfrac 32} \lambda_{q+\lfloor\sfrac \bn 2 \rfloor}^N \MM{M,\Nindt-1,\tau_q^{-1}\Gamma_{q}^{i+13},\Tau_q^{-1}\Ga_q^8} \notag\\
&\lesssim \const_{G,\sfrac 32}\lambda_{q+\lfloor\sfrac \bn 2 \rfloor}^N \MM{M,\Nindt,\tau_q^{-1}\Gamma_{q}^{i+14},\Tau_q^{-1}\Ga_q^8}
\,
,\label{eq:david:transport:0}\\
\Biggl| D^N \Dtq^M G_{i,j,k,\xi,\vecl,I,\diamond} \Biggr| &\lessg r_q \Ga_q^{50} \pi_\ell^{\sfrac 12} \tau_q^{-1}\Ga_q^i \lambda_{q+\lfloor\sfrac \bn 2 \rfloor}^N \MM{M,\Nindt-1,\tau_q^{-1}\Gamma_{q}^{i+13},\Tau_q^{-1}\Ga_q^8} \notag\\
&\lessg r_q r_{q-\bn}^{-1} \Ga_q^{100} \pi_\ell \La_q \lambda_{q+\lfloor\sfrac \bn 2 \rfloor}^N \MM{M,\Nindt-1,\tau_q^{-1}\Gamma_{q}^{i+13},\Tau_q^{-1}\Ga_q^8} \notag\\
&\lesssim \pi_\ell \La_q \lambda_{q+\lfloor\sfrac \bn 2 \rfloor}^N \MM{M,\Nindt,\tau_q^{-1}\Gamma_{q}^{i+14},\Tau_q^{-1}\Ga_q^8} 
\, , \label{eckel:TN:1}
\end{align}
where we have used \eqref{condi.Nindt} to upgrade the sharp derivatives to $\Nindt$ in both inequalities, \eqref{eq:psi:q:q'}, \eqref{ineq:tau:q}, and \eqref{ind:pi:lower} to convert $\tau_q^{-1}\Ga_q^i$ into $\pi_\ell^{\sfrac 12}\Ga_q^{50}\La_q r_{q-\bn}^{-1}$ in the pointwise bounds, and \eqref{ineq:r's:eat:Gammas} to absorb the $\Ga_q^{100}$. In order to obtain an $L^\infty$ bound, we can appeal to \eqref{eckel:TN:1} and \eqref{eq:pressure:inductive:dtq-1:uniform:upgraded}. Thus we have that \eqref{eq:inverse:div:DN:G} and \eqref{eq:inv:div:extra:pointwise} are satisfied in all cases.

In order to verify the high-frequency assumptions from Part 2 of Proposition~\ref{prop:intermittent:inverse:div}, we set 
\begin{align}
    &r_q \varrho_R = \varrho_{\pxi,R}^I \, , \quad r_q \vartheta_R \textnormal{ as defined in item~\eqref{item:pipe:1} from Proposition~\ref{prop:pipeconstruction}} \notag \\
    &r_q^{\sfrac 43} \varrho_\vp = \varrho_{\pxi,\vp}^I \, , \quad r_q^{\sfrac 43} \vartheta_\vp \textnormal{ defined similarly but adjusted to fit Proposition~\ref{prop:pipe.flow.current}} \notag  \\
    &\Ndec\textnormal{ as in \eqref{i:par:9}}  \, , \quad \dpot\textnormal{ as in \eqref{i:par:10}} \, , \quad \const_{*,\sfrac 32} = r_q^{-\sfrac 23} \, , \quad \const_{*,\infty} = r_q^{-2} \, ,  \notag \\
    &\mu = \la_\qbn r_q = \la_{q+\half} \Ga_q \, , \quad \Upsilon= \Upsilon' = \La = \la_{\qbn} \, . \label{more:more:shortening}
\end{align}
Then we have that \eqref{item:inverse:i} is satisfied from \eqref{eq:WW:explicit}, \eqref{item:inverse:ii} is satisfied by the construction of $w_{q+1}$ in subsection~\ref{ss:corr-ec-tor}, and \eqref{eq:DN:Mikado:density} is satisfied from Proposition~\ref{prop:pipeconstruction} or the corresponding estimates in Proposition~\ref{prop:pipe.flow.current}. Finally, we have that \eqref{eq:inverse:div:parameters:0} follows by definition and from \eqref{condi.Nfin0}, while \eqref{eq:inverse:div:parameters:1} is satisfied from \eqref{condi.Ndec0}.

We therefore may appeal to the local conclusions \eqref{item:div:local:0}--\eqref{item:div:nonlocal} and \eqref{eq:inverse:div:error:stress}--\eqref{eq:inverse:div:error:stress:bound}, from which we have the following. First, we note that from \eqref{item:div:local:ii}, we have that \eqref{trans.loc.pot} is satisfied. Next, we have from \eqref{eq:inverse:div}, \eqref{eq:inverse:div:stress:1}, and \eqref{eq:inv:div:extra:conc} that for $N\leq \frac{\Nfin}{4}-\dpot$ and $M\leq \frac{\Nfin}{5}$,
\begin{align}
\left\| D^N \Dtq^M \left( \divH \left(  T_{i,j,k,\xi,\vecl,I,\diamond} \right) \right) \right\|_{\sfrac 32} &\lessg \left( \left| \supp(\eta_{i,j,k,\xi,\vecl,I,\diamond}\zetab_\xi^{I,\diamond})\right|^{\sfrac 23} \delta_{q+\bn}^{\sfrac 12} r_q^{\sfrac 13} \Gamma_q^{i+j+25} \tau_q^{-1} + \la_\qbn^{-10} \right) \notag\\
&\qquad \times  \la_\qbn^{-1+N} \MM{M,\Nindt,\tau_q^{-1}\Gamma_{q}^{i+14},\Tau_q^{-1}\Ga_q^8} \, , \label{eq:david:transport:2} \\
\left| D^N \Dtq^M \left( \divH \left(  T_{i,j,k,\xi,\vecl,I,\diamond} \right) \right) \right| &\lessg \pi_\ell \La_q r_q^{-2} \la_\qbn^{-1+N} \MM{M,\Nindt,\tau_q^{-1}\Gamma_{q}^{i+14},\Tau_q^{-1}\Ga_q^8} \, . \label{eq:david:transport:2:2}
\end{align}
Notice that from \eqref{item:div:local:i}, the support of $\div \mathcal{H} T_{i,j,k,\xi,\vecl,I,R}$ is contained in the support of $T_{i,j,k,\xi,\vecl,I,R}$, which itself is contained in the support of $\eta_{i,j,k,\xi,\vecl,R}\zetab_\xi^{I,R}$.  From this observation, we have that \eqref{eq:H:eckel:trans} is satisfied.  Furthermore, we have that \eqref{rho:eckel:estimate:TN} is satisfied from \eqref{eq:inverse:div:sub:1} and the estimates from Proposition~\ref{prop:pipeconstruction} and \ref{prop:pipe.flow.current}.  Next, we have that \eqref{eq:H:eckel:TN} is satisfied from \eqref{eq:inv:div:extra:conc}. Finally, we have that \eqref{rho:eckel:support:TN} holds due to item~\eqref{item:div:local:i} and item~\eqref{item:pipe:3.5} from Proposition~\ref{prop:pipeconstruction}. We note also that \eqref{eq:trans:supp} follows from \eqref{eq:H:eckel:trans}, \eqref{rho:eckel:support:TN}, and \eqref{eq:dodging:oldies}. 

In order to aggregate $L^{\sfrac 32}$ estimates, we appeal to Corollary~\ref{rem:summing:partition} with $\theta_1=\theta_2=1$, $H=\divH\left( T_{i,j,k,\xi,\vecl,I,\diamond} \right)$, \eqref{eq:inductive:partition} at level $q$, and \eqref{ineq:transport:basic} to write that
\begin{align}
    &\left\|\psi_{i,q} \sum_{i',j,k,\xi,\vecl,I,\diamond}  D^N \Dtq^M \left( \divH \left(  T_{i',j,k,\xi,\vecl,I,\diamond} \right) \right) \right\|_{\sfrac 32} \notag\\
    &\qquad\lesssim \Gamma_q^{50+ \CLebesgue} \delta_{q+\bn}^{\sfrac 12} r_q^{\sfrac 13} \lambda_{q+\bn}^{-1+N} \tau_q^{-1}  \MM{M,\Nindt,\tau_q^{-1}\Gamma_{q}^{i+15},\Tau_q^{-1}\Ga_q^8} \notag\\
    &\qquad\lesssim \Ga_{q+\bn}^{-25} \delta_{q+2\bn} \lambda_{q+\bn}^{N} \MM{M,\Nindt,\tau_q^{-1}\Gamma_{q}^{i+15},\Tau_q^{-1}\Ga_q^8} \, .
\end{align}
In order to aggregate pointwise estimates, we appeal to Corollary~\ref{lem:agg.pt} with the same choice of $H$ and $\varphi = \pi_\ell\La_q r_q^{-2} \mathbf{1}_{\supp (\eta_{i,j,k,\xi,\vecl,R}\zetab_\xi^{I,R})}$.  Then from \eqref{eq:aggpt:conc:2}, \eqref{eq:desert:cowboy:sum}, \eqref{eq:pressure:inductive:dtq:uniform:upgraded}, and \eqref{eq:par:div:2}, we have that
\begin{align*}
    \left|\psi_{i,q} \sum_{i',j,k,\xi,\vecl,I,\diamond}  D^N \Dtq^M \left( \divH \left(  T_{i',j,k,\xi,\vecl,I,\diamond} \right) \right) \right| &\lesssim \pi_\ell \La_q r_q^{-2} \lambda_{q+\bn}^{-1+N} \MM{M,\Nindt,\tau_q^{-1}\Gamma_{q}^{i+15},\Tau_q^{-1}\Ga_q^8} \\
    &\leq \Ga_\qbn^{\badshaq-200} \lambda_{q+\bn}^N \MM{M,\Nindt,\tau_q^{-1}\Gamma_{q}^{i+15},\Tau_q^{-1}\Ga_q^8} \, .
\end{align*}

To conclude the proof for the leading order term from $\Dtq w_{q+1}$, we must still estimate the nonlocal $\divR$ portion of the inverse divergence. In order to check the nonlocal assummptions, we again set
\begin{align*}
    M_\circ = N_\circ = 2\Nind \, , \quad K_\circ \textnormal{ as in \eqref{i:par:9.5}} \, .
\end{align*}
Then from \eqref{ineq:dpot:1} and Remark~\ref{rem:lossy:choices}, we have that \eqref{eq:inv:div:wut}--\eqref{eq:riots:4} are satisfied. We note that $\Dtq w_{q+1} + w_{q+1}\cdot \na \hat u_q$ has zero mean, and so we may ignore the means of individual terms that get plugged into the inverse divergence since their sum will vanish. Then from \eqref{eq:inverse:div:error:stress}, \eqref{eq:inverse:div:error:stress:bound}, and Remark~\ref{rem:lossy:choices}, we have that for $N,M\leq 2\Nind$,
\begin{align}
    \left\| D^N \Dtq^M \sum_{i,j,k,\xi,\vecl} \divR T_{i,j,k,\xi,\vecl,\diamond} \right\|_\infty \leq \delta_{q+3\bn}^2 \Tau_\qbn^{2\Nindt} {\lambda_{q+\bn}^N} \tau_q^{-M} \, , \notag
\end{align}
matching the desired estimate in \eqref{eq:trans:nonlocal}.
\end{proof}

At this point, we can construct the pressure increment and associated current error coming from the Nash and transport errors.  Since the proofs of both lemmas are completely analogous to the proofs of the corresponding lemmas for the \emph{highest frequency shell} from \eqref{osc.err.high} of the oscillation error, we omit the majority of the details and merely note the minor differences required in a combined proof.

\begin{lemma*}[\bf Pressure increment]\label{lem:transport:pressure}
There exists a function $\si_{S_{TN}} = \si_{S_{TN}}^+ - \si_{S_{TN}}^-$ such that the following hold.\index{$\sigma_{S_{TN}}$}
\begin{enumerate}[(i)]
    \item We have that
\begin{subequations}
\begin{align}
    \label{eq:t.p.1}
    \left|\psi_{i,q} D^N \Dtq^M S_{TN}\right| &< \left(\si_{S_{TN}}^+ + \de_{q+3\bn}\right) \left(\lambda_{q+\bn}\Gamma_q \right)^N \MM{M,\Nindt,\tau_q^{-1}\Gamma_{q}^{i+17},\Tau_q^{-1}\Ga_q^9}\\
    \label{eq:t.p.2}
    \left|\psi_{i,q} D^N \Dtq^M \si_{S_{TN}}^+\right| &< \left(\si_{S_{TN}}^+ +\de_{q+3\bn}\right) \left(\lambda_{q+\bn}\Gamma_q \right)^N \MM{M,\Nindt,\tau_q^{-1}\Gamma_{q}^{i+17},\Tau_q^{-1}\Ga_q^9}\\
    \label{eq:t.p.3}
    \norm{\psi_{i,q} D^N \Dtq^M \si_{S_{TN}}^+}_{\sfrac32} &\leq \Ga_{q+\bn}^{-9}\de_{q+2\bn} \left(\lambda_{q+\bn}\Gamma_q \right)^N \MM{M,\Nindt,\tau_q^{-1}\Gamma_{q}^{i+17},\Tau_q^{-1}\Ga_q^9}\\
    \label{eq:t.p.3.inf}
    \norm{\psi_{i,q} D^N \Dtq^M \si_{S_{TN}}^+}_{\infty} &\leq \Ga_{q+\bn}^{\badshaq-9} \left(\lambda_{q+\bn}\Gamma_q \right)^N \MM{M,\Nindt,\tau_q^{-1}\Gamma_{q}^{i+17},\Tau_q^{-1}\Ga_q^9}\\
    \label{eq:t.p.4}
    \left|\psi_{i,q} D^N \Dtq^M \si_{S_{TN}}^-\right| &\leq \Ga_{q+\half}^{-100} \pi_q^{q+\half}  \left(\lambda_{q+\half}\Gamma_q \right)^N \MM{M,\Nindt,\tau_q^{-1}\Gamma_{q}^{i+17},\Tau_q^{-1}\Ga_q^9}
\end{align}
\end{subequations}
for all $N,M < \sfrac{\Nfin}{100}$.
\item For all $q+1\leq q'\leq q+\half$ and $q+1\leq q'' \leq \qbn-1$, we have that
\begin{align}\label{eq:t.p.6}
     B\left( \supp \hat w_{q'}, \la_{q'}^{-1} \Ga_{q'+1} \right) \cap \supp \sigma_{S_{TN}}^- =  B\left( \supp \hat w_{q''}, \la_{q''}^{-1} \Ga_{q''+1} \right) \cap \supp \sigma_{S_{TN}}^+ = \emptyset \, .
\end{align}
\item\label{i:pc:5:ER:TN} Define 
\begin{equation}\label{def:bmu:ER:TN}
    \bmu_{\sigma_{S_{TN}}}(t) = \int_0^t \left \langle \Dtq \sigma_{S_{TN}}  \right \rangle (s) \, ds \, .
\end{equation}
Then we have that for $0\leq M\leq 2\Nind$,
    \begin{align}\label{th:billys:2}
      \left|\frac{d^{M+1}}{dt^{M+1}} \bmu_{\sigma_{S_{TN}}} \right| 
      \leq (\max(1, T))^{-1}\delta_{q+3\bn} \MM{M,\Nindt,\tau_q^{-1},\Tau_{q+1}^{-1}} \, .
    \end{align}
\end{enumerate}
\end{lemma*}

\begin{lemma*}[\bf Pressure current]\label{lem:transport:pressure:current}
There exists a current error $\phi_{S_{TN}}$ associated to the pressure increment $\sigma_{S_{TN}}$ defined by Lemma~\ref{lem:transport:pressure} which satisfies the following properties.
\begin{enumerate}[(i)]
    \item We have the decomposition and equalities
       \begin{subequations}
    \begin{align}\label{eq:desert:decomp:ER:TN}
        \phi_{{S_{TN}}} &= \phi_{{S_{TN}}}^* + \sum_{m'=q+\half+1}^{{\qbn}} \phi_{{S_{TN}}}^{m'} \, , \qquad
         \phi_{S_{TN}}^{m'} = \phi_{S_{TN}}^{m',l} + \phi_{S_{TN}}^{m',*} \, \\
         \div \phi_{{S_{TN}}}
          &= D_{t,q}\si_{{S_{TN}}} - \langle D_{t,q}\si_{{S_{TN}}}\rangle\, . 
    \end{align}
    \end{subequations}
    \item For all $N,M \leq 2\Nind$,
\begin{align}
    &\left|\psi_{i,q} D^N \Dtq^M \phi_{S_{TN}}^{k',l}\right| < \Ga_{ k'}^{-100} r_{ k'}^{-1} \left(\pi_q^{ k'}\right)^{\sfrac32}\left(\lambda_{ k'}\Gamma_{m'}^2\right)^N \MM{M,\Nindt,\tau_q^{-1}\Gamma_{q}^{i+18},\Tau_q^{-1}\Ga_q^9} \, , \label{t:p:c:pt} \\
    &\norm{D^N \Dtq^M \phi_{S_{TN}}^{k',*}}_{L^\infty}
    < \Tau_{q+\bn}^{2\Nindt}
    \delta_{q+3\bn}^{\sfrac 32} (\la_\qbn\Ga_q^2)^N \tau_q^{-M} \, . \label{t:p:c:nonlocal}
    \end{align}
\item For all $m',q'$ with $q+1 \leq q'\leq m'-1$ and $q+\half+1\leq m' \leq \qbn$, we have that
\begin{align}\label{t:p:c:supp}
    B\left(\supp \hat w_{q'}, \sfrac 12 \la_{q'}^{-1} \Ga_{q'+1} \right) \cap \supp \phi_{S_{TN}}^{k',l} = \emptyset \, .
\end{align}
\end{enumerate}
\end{lemma*}

\begin{proof}[Proofs of Lemmas~\ref{lem:transport:pressure} and \ref{lem:transport:pressure:current}]
As in Lemmas~\ref{lem:oscillation:pressure} and \ref{lem:oscillation:pressure:current}  in the case $m=\qbn$, the proofs of Lemmas~\ref{lem:transport:pressure} and \ref{lem:transport:pressure:current} use Proposition~\ref{lem.pr.invdiv2} to estimate a single error term indexed by $i,j,k,\xi,\vecl,I,\diamond$, and then aggregate estimates according to Corollaries~\ref{rem:summing:partition}--\ref{lem:agg.Dtq}. We now identify the minor differences between the applications of these various tools to the transport/Nash error and the oscillation error.

We first check the preliminary assumptions from Part 1 of Proposition~\ref{lem.pr.invdiv2}. Let us first compare the low-frequency parameter choices for the transport error in \eqref{eq:eckel:TN:choices} to the low-frequency parameter choices for the error terms in \eqref{osc.err.high}, which was analyzed in \texttt{Case 3b} from Lemma~\ref{lem:oscillation:general:estimate}.  First, we have that the vector field $G$ in \eqref{eq:eckel:TN:choices} is different than the vector field in \eqref{shortening:some:stuff}, but it retains the exact same support properties due to the presence of $\rhob_\xi^\diamond \zetab_\xi^\diamond$ in both.  Next, we claim that $\const_{G,p}$ is effectively \emph{smaller} in \eqref{eq:eckel:TN:choices} than in \eqref{shortening:some:stuff}. In the case $p=\infty$, this is immediate, so we focus on the case $\sfrac 32$. We use \eqref{ineq:tau:q}, \eqref{ineq:r's:eat:Gammas}, and \eqref{eq:par:div:1} to write that
$$  \tau_q^{-1} r_q \leq \Ga_q^{50} \la_q \delta_q^{\sfrac 12} r_{q-\bn}^{-\sfrac 13} r_q \leq \delta_{q+\bn}^{\sfrac 12} \La_q \Ga_q^{-50} \, . $$
The difference between $\Ga_q^{i+j}$ in \eqref{eq:eckel:TN:choices} and $\Ga_q^{2j}$ in \eqref{shortening:some:stuff} only matters in the application of Corollaries~\ref{rem:summing:partition}--\ref{lem:agg.Dtq}. Indeed, trading a $j$ for an $i$ simply necessitates a difference choice of $\theta_1$ and $\theta_2$, and the only difference in the output is the factor of  $\Ga_q^{\theta_1 \CLebesgue}$ which must be absorbed in the latter case. The reader is invited to check inequalities \eqref{ineq:in:the:afternoon:ER}, \eqref{ineq:in:the:morning:ER}, \eqref{eq:par:div:2}, \eqref{ineq:osc:general}, and \eqref{o:5}, each of which has a $\Ga_q^{5\CLebesgue}$ on the left-hand side that can therefore absorb this extra insignificant factor. Next, we have that the choices of $M_t,M_*,N_*,\lambda,\nu,\nu'$ are the same, and the choice of $\varpi=\pi_\ell\Ga_q^{50}\La_q$ from the beginning of Lemma~\ref{lem:oscillation:pressure} is larger than the choice of $\varpi$ from \eqref{eq:eckel:TN:choices} for the transport error.  Finally, the vector field $v$ and associated material derivative $D_t$ from item~\eqref{i:st:sample:2} are identical in both cases.

Next, we compare the high-frequency parameter choices from item~\eqref{i:st:sample:3} in the case of the oscillation error in \eqref{more:shortening} to the choices for the transport error in \eqref{more:more:shortening}. The potential $\vartheta$ in \eqref{more:more:shortening} is supported in a $\la_\qbn^{-1}$ neighborhood of $\varrho_{\pxi,\diamond}^I$, while for the oscillation error, the support is larger due to the presence of the synthetic Littlewood-Paley projector $\td{\mathbb{P}}_{(\la_{\qbn-1},\qbn]}$ applied to $(\varrho_{\pxi,\diamond}^I)^2$.  Thus the potential for transport error has more advantageous support properties than that of the oscillation error. Next, the choices of $\mu$ and $\Lambda$ are identical, while the choices of $\Upsilon$ and $\Upsilon'$ are \emph{more} advantageous for the transport error than they are for the oscillation error in the case $m=\qbn$.  Indeed, this is because the inverse divergence gain in the transport error is a full $\la_\qbn$ from \eqref{eq:WW:explicit}, while the highest shell of the oscillation error only gains $\la_{\qbn-1}$ due to the presence of the synthetic Littlewood-Paley projector. Next, the choices of $\const_{*,p}$ are identical due to our choice of rescaling in the transport error, and the choices of $\Ndec$ and $\dpot$ are identical as well.  Therefore, all assumptions from item~\eqref{i:st:sample:3} are \emph{stronger} for the transport error than the oscillation error.  Finally, we note that the nonlocal assumptions in item~\eqref{i:st:sample:6} are not changed in any significant way, and so we may treat the nonlocal transport error terms in the same way as the nonlocal oscillation error terms.

Moving to the additional assumptions from Part 2 of Proposition~\ref{lem.pr.invdiv2}, we have that all inequalities in \eqref{helping:matt:shorten}, \eqref{eq:sample:prop:Ncut:1}, \eqref{eq:sample:prop:Ncut:3}, \eqref{eq:sample:prop:decoup} are identical. The inequality in \eqref{eq:sample:prop:Ncut:2} follows in the same was as in the oscillation error; indeed, all nonlocal error bounds can be treated in the same way via a large choice of $\dpot$ or $N_{**}$. The inequalities in item~\eqref{i:st:sample:8} are the same for the transport error as for the highest shell of the oscillation error, since these inequalities relate to the synthetic Littlewood-Paley projection of a function which oscillates at frequency $\approx\La = \la_\qbn$.  

Now that we have highlighted the unimportant differences in the set-up, we merely note that the sharp material derivative cost in Lemmas~\ref{l:transport:error}--\ref{lem:transport:pressure:current} is worse by a factor of $\Ga_q$ than the corresponding estimates in Lemmas~\ref{lem:oscillation:general:estimate}--\ref{lem:oscillation:pressure:current}.  This is due to the fact that the transport error loses a material derivative.  This concludes the proofs of Lemmas~\ref{lem:transport:pressure} and \ref{lem:transport:pressure:current}.
\end{proof}

\subsection{Divergence corrector error \texorpdfstring{$S_C$}{essee}}\label{sss:dce}

We will write the divergence corrector error as 
\begin{equation}
S_C = S_{C1} + S_{C2} \, , \qquad  \textnormal{for} \qquad \div S_{C1} = \div \left( w_{q+1}^{(p)}\otimes_s w_{q+1}^{(c)}\right) \, , \qquad S_{C2} = w_{q+1}^{(c)} \otimes w_{q+1}^{(c)} \, , \label{eq:div:cor:expand}
\end{equation}
and estimate them in the following lemma.\index{$S_C$}

\begin{lemma}[\bf Basic estimates and applying inverse divergence]\label{l:divergence:corrector:error}
There exist symmetric stresses $S_C^{m}$ for $m\in\{q+\lfloor\sfrac{\bn}{2}\rfloor+1, \dots , q+\bn\}$ such that the following hold.
\begin{enumerate}[(i)]
    \item $\div \left( w_{q+1}^{(p)}\otimes_s w_{q+1}^{(c)} +  w_{q+1}^{(c)} \otimes w_{q+1}^{(c)} \right) = \sum_{m=q+\lfloor\sfrac{\bn}{2}\rfloor+1}^{q+\bn} \div S_C^{m}$, where $S_C^m$ can be split into local and non-local errors as $S_C^m = S_C^{m,l}+S_C^{m,*}$. 
\item For the same range of $m$ and for all $N,M\leq\sfrac{\Nfin}{10}$, the local parts $S_C^{m,l}$ satisfy
\begin{subequations}
\begin{align}
&\left\| \psi_{i,q} D^N \Dtq^M S_C^{m,l}
\right\|_{{\sfrac 32}} \les \Gamma_m^{-9} \delta_{m+\bn}\lambda_{m}^N \MM{M,\Nindt, \tau_q^{-1}\Gamma_{q}^{i + 15},\Tau_q^{-1}\Ga_q^8}
\label{eq:div:corrector:L1}
\\
&\left\| \psi_{i,q} D^N \Dtq^M S_C^{m,l} \right\|_{\infty} \les \Gamma_m^{-9}  \lambda_{m}^N \MM{M,\Nindt, \tau_q^{-1}\Gamma_{q}^{i + 15},\Tau_q^{-1}\Ga_q^8} \, . \label{eq:div:corrector:Linfty}
\end{align}
\end{subequations}
\item For $q+\half +1 \leq m \leq q+\bn$ and $q+1 \leq q'\leq m-1$, the local parts satisfy
\begin{equation}\label{eq:dc:ER:supp}
    \supp S_C^{m,l} \cap 
    B \left( \supp \hat w_{q'}, \lambda_{q'}^{-1}\Gamma_{q'+1} \right)  = \emptyset \, .
\end{equation}
\item For the same range of $m$ and $N,M\leq 2\Nind$, the nonlocal parts $S_C^{m,\ast}$ satisfy
\begin{align}\label{eq:divER:nonlocal}
    \left\| D^N \Dtq^M S_C^{m,\ast} \right\|_{\infty} \leq \Tau_\qbn^{4\Nindt}\delta_{q+3\bn} \lambda_m^N \tau_q^{-M} \, .
\end{align}
\end{enumerate}
\end{lemma}

\begin{remark*}[\bf Abstract formulation of the divergence corrector errors]\label{rem.ct.divcorr}
For the purposes of analyzing the transport and Nash current errors in subsection~\ref{op:tnce} and streamlining the creation of pressure increments, it is useful again to abstract the properties of these error terms. As we shall see in the course of the proof in Lemma~\ref{l:divergence:corrector:error}, however, these error terms may be decomposed and analyzed in \emph{exactly} the same way as the oscillation errors. This is not surprising, since both error terms are quadratic in $w_{q+1}$, and morally speaking, one expects the estimates for terms involving divergence correctors to be slightly better.  Therefore we refer the reader to Remark~\ref{rem.ct.osc} rather than reproduce it in entirety here.
\end{remark*}

\begin{proof}[Proof of Lemma~\ref{l:divergence:corrector:error}]

The analysis in the proof generally follows that of the divergence corrector errors in \cite{NV22}, and we shall occasionally refer to algebraic identities from those arguments. The main difference is that we have to incorporate the synthetic Littlewood-Paley projector in certain terms before applying  the inverse divergence operator in order to upgrade the material derivatives later. However, synthetic Littlewood-Paley projectors have already been applied to terms which are quadratic in high frequency objects in Lemma~\ref{lem:oscillation:general:estimate}, and so we may pirate a significant portion of the analysis from there as well.
\smallskip

\noindent\texttt{Step 1}. We first consider $\div \bigl( w_{q+1}^{(p)} \otimes_s w_{q+1}^{(c)}\bigr)$. We write that
\begin{align}
\div \bigl( w_{q+1}^{(p)} \otimes_s w_{q+1}^{(c)} \bigr)^{\bullet} &= \sum_{\diamond,i,j,k,\xi,\vecl,I}
\partial_m \bigg{(} a_{(\xi),\diamond} \left(\rhob_{(\xi)}^\diamond \zetab_{\xi}^{I,\diamond} \varrho_{(\xi),\diamond}^{I} \right)\circ \Phiik \xi^\ell \bigl( A_{\ell}^m  \epsilon_{\bullet p r} + A_{\ell}^\bullet  \epsilon_{m p r} \bigr) \notag\\
&\qquad \qquad \qquad \times \partial_p \left( a_{(\xi),\diamond} \left( \rhob_{(\xi)}^\diamond \zetab_\xi^{I,\diamond} \right)\circ \Phiik\right)  \partial_r \Phi_{(i,k)}^s (\mathbb{U}_{(\xi),\diamond}^{I})^s \circ \Phi_{(i,k)} 
\bigg{)} \, ,
\label{eq:nets:suck:1}
\end{align} 
where we have used Lemma~\ref{lem:dodging}, the definition of $\WW_{(\xi),\diamond}^I$ in \eqref{eq:WW:explicit} (and the corresponding version for $L^3$ normalized pipes), $\epsilon_{i_1i_2i_3}$ is the Levi-Civita alternating tensor, we implicitly contract the repeated indices $\ell,m,p,r,s$, and the $\bullet$ refers to the indices of the vectors on either side of the above display. Using that $\{\xi,\xi',\xi''\}$ is an orthonormal basis associated with the direction vector $\xi$ with $\xi \times \xi' = \xi''$ and decomposing as in \cite[(7.50)]{NV22}, we have that
\begin{align}
    \partial_p &\left( a_{(\xi),\diamond} \left( \rhob_{(\xi)}^\diamond \zetab_\xi^{I,\diamond} \right)\circ \Phiik\right)  =   \underbrace{\partial_p \Phiik^n \xi^n  \xi^{\ell} A_\ell^j \partial_j \left( a_{(\xi),\diamond} \left( \rhob_{(\xi)}^\diamond \zetab_\xi^{I,\diamond} \right)\circ \Phiik\right)  }_{=:  a_{(\xi),\diamond}^{p, \rm good}} \qquad      \label{eq:nets:suck:2}\\
    & +     \underbrace{\partial_p \Phiik^n (\xi^\prime)^n (\xi^\prime)^{\ell} A_\ell^j \partial_j \left( a_{(\xi),\diamond} \left( \rhob_{(\xi)}^\diamond \zetab_\xi^{I,\diamond} \right)\circ \Phiik\right) + \partial_p \Phiik^n (\xi^{\prime \prime})^n (\xi^{\prime \prime})^\ell A_\ell^j \partial_j \left( a_{(\xi),\diamond} \left( \rhob_{(\xi)}^\diamond \zetab_\xi^{I,\diamond} \right)\circ \Phiik\right)  }_{=: a_{(\xi),\diamond}^{p, \rm bad} }
     \,, \notag
\end{align}
where we have also set $A = A_{(i,k)} =  (\nabla \Phiik)^{-1}$. Indeed, the good differential operator appearing in $a_{(\xi),\diamond}^{p,\rm good}$ only costs $\Lambda_q\Gamma_q^{13}$ (see Lemma~\ref{lem:a_master_est_p}), so that we will leave $a_{(\xi),\diamond}^{p, \rm good}$ inside the divergence and dump the symmetric stress inside of the divergence into $S_C^{q+\bn}$. On the other hand, $a_{(\xi),\diamond}^{p, \rm bad}$ contains an expensive derivative at $\lambda_{q+\lfloor \sfrac \bn 2 \rfloor}$, but $\xi^\ell A_\ell^m \partial_m$ only costs $\Lambda_q \Gamma_q^{13}$, which will be crucially used below.

Splitting the terms involved with $a_{(\xi),\diamond}^{p, \rm bad}$ from \eqref{eq:nets:suck:1} as in \cite[(7.52)]{NV22}, we further analyze
\begin{align}
&\sum_{\diamond,i,j,k,\xi,\vec{l},I}
\partial_m \left( a_{(\xi),\diamond} \left(\rhob_{(\xi)}^\diamond \zetab_{\xi}^{I,\diamond} \varrho_{(\xi)}^{I,\diamond} \right)\circ \Phiik \xi^\ell \bigl( A_{\ell}^m  \epsilon_{\bullet p r} + A_{\ell}^\bullet  \epsilon_{m p r} \bigr) a_{(\xi),\diamond}^{p, \rm bad} \partial_r \Phi_{(i,k)}^s (\mathbb{U}_{(\xi),\diamond}^{I})^s \circ \Phi_{(i,k)}
\right) 
= \mathbf{V}_1^\bullet + \mathbf{V}_2^\bullet
\label{eq:nets:suck:3}
\end{align}
where $\mathbf{V}_1$ contains $A_{\ell}^m  \epsilon_{\bullet p r}$, and $\mathbf{V}_2$ contains $A_{\ell}^\bullet \epsilon_{m p r}$.  To analyze $\mathbf{V}_1$, we use that $\partial_m$ and $\xi^\ell A_\ell^m$ commute, so that
\begin{align*}
   \xi^\ell A_{\ell}^m \partial_m \left((\varrho_{(\xi),\diamond}^{I} (\UU_{(\xi),\diamond}^{I})^s) \circ \Phiik \right)  &= 0 
   \, .
\end{align*}
Furthermore, the differential operator $\xi^\ell A_\ell^m \partial_m$ landing anywhere else costs only $\Lambda_{q}\Gamma_q^{13}$ from \eqref{e:a_master_est_p}. Then we have in total that
\begin{align}
\mathbf{V}_1^\bullet
&= \sum_{\diamond,i,j,k,\xi,\vec{l},I} \partial_m \left( a_{(\xi),\diamond} \left(\rhob_{(\xi)}^\diamond \zetab_{\xi}^{I,\diamond}  \right)\circ \Phiik \xi^\ell A_{\ell}^m  \epsilon_{\bullet p r} a_{(\xi),\diamond}^{p, \rm bad} \partial_r \Phi_{(i,k)}^s 
\right) \left( \varrho_{(\xi),\diamond}^{I} (\mathbb{U}_{(\xi),\diamond}^{I})^s \right) \circ \Phi_{(i,k)}\label{eq:nets:suck:4}\\
&=:\sum_{\diamond,i,j,k,\xi,\vec{l},I}
(C_{(\xi), \diamond}^{1, I})^{\bullet s}\left( \varrho_{(\xi),\diamond}^{I} (\mathbb{U}_{(\xi),\diamond}^{I})^s \right) \circ \Phi_{(i,k)} \nonumber
\end{align}
is a product of a high-frequency, mean-zero potential which has gained one factor of $\lambda_{q+\bn}$, and a low-frequency object which has lost one costly derivative at frequency $\lambda_{q+\lfloor \sfrac \bn 2 \rfloor}$, and one cheap derivative at frequency $\Lambda_{q}\Gamma_q^{13}$. To analyze $\mathbf{V}_2$, we follow  \cite[7.56]{NV22} to get
\begin{align}
\mathbf{V}_2^\bullet   &= \sum_{\diamond,i,j,k,\xi,\vecl,I} \partial_m \left( a_{(\xi),\diamond} \left(\rhob_{(\xi)}^\diamond \zetab_{\xi}^{I,\diamond} \varrho_{(\xi)}^{I,\diamond} \right)\circ \Phiik \xi^\ell A_{\ell}^\bullet \epsilon_{m p r} a_{(\xi),\diamond}^{p, \rm bad} \partial_r \Phi_{(i,k)}^s (\mathbb{U}_{(\xi)}^{I,\diamond})^s \circ \Phi_{(i,k)}
\right)  \notag\\ 
&= \sum_{\diamond,i,j,k,\xi,\vecl,I}  \Bigl(   \partial_m \bigl(   \xi^\ell   A_{\ell}^\bullet  \epsilon_{m p r} \partial_r \Phi_{(i,k)}^s    \bigr) a_{(\xi),\diamond} (\rhob_{(\xi)}^\diamond \zetab_\xi^{I,\diamond})\circ\Phiik a_{(\xi),\diamond}^{p, \rm bad} +  a_{(\xi),\diamond}^{m, \rm good}  \xi^\ell   A_{\ell}^\bullet  \epsilon_{m p r} a_{(\xi),\diamond}^{p, \rm bad} \partial_r \Phi_{(i,k)}^s \notag\\ 
&\qquad \qquad \qquad \qquad \qquad  -  a_{(\xi),\diamond} (\rhob_{(\xi)}^\diamond \zetab_\xi^{I,\diamond})\circ\Phiik \xi^\ell   A_{\ell}^\bullet  \epsilon_{m p r}   \partial_m (a_{(\xi),\diamond}^{p, \rm good})   \partial_r \Phi_{(i,k)}^s  \Bigr)  \left(\varrho_{(\xi)}^{I,\diamond} (\mathbb{U}_{(\xi)}^{I,\diamond})^s  \right)\circ \Phi_{(i,k)}  \notag\\  
&\quad+  \sum_{\diamond,i,j,k,\xi,\vecl,I} a_{(\xi),\diamond} \xi^\ell A_{\ell}^\bullet  \epsilon_{m p r}    a_{(\xi),\diamond}^{p, \rm bad}   \partial_r \Phi_{(i,k)}^s  \partial_m \left(\varrho_{(\xi)}^{I,\diamond} (\mathbb{U}_{(\xi)}^{I,\diamond})^s  \right)\circ \Phi_{(i,k)}   \,.  \label{eq:nets:suck:5}\\
&=: \sum_{\diamond,i,j,k,\xi,\vecl,I}
(C_{(\xi), \diamond}^{2, I})^{\bullet s}\left( \varrho_{(\xi),\diamond}^{I} (\mathbb{U}_{(\xi),\diamond}^{I})^s \right) \circ \Phi_{(i,k)} \nonumber 
\end{align}
In the second equality above we have used the identities $\epsilon_{m p r}\partial_m  (a_{(\xi),\diamond}^{p,\rm bad}) = - \epsilon_{mpr} \partial_m (a_{(\xi),\diamond}^{p, \rm good})$, which follows from \eqref{eq:nets:suck:2}, and $\epsilon_{m p r} a_{(\xi),\diamond}^{m, \rm bad} a_{(\xi),\diamond}^{p, \rm bad} = 0$. Furthermore, we recall from \cite[pgs. 42-43]{NV22} that the last term on the right-hand side of the second equality vanishes. As before, the slow function $C_{(\xi), \diamond}^{2, I}$ contains two spatial derivatives, one cheap and one expensive.
\smallskip

\noindent\texttt{Step 2}. We now define the stress error $S_C^m$ from the divergence corrector. From \eqref{item:pipe:means} of Proposition \ref{prop:pipeconstruction} and \eqref{item:pipe:means:current} of Proposition \ref{prop:pipe.flow.current}, we know that $\varrho_{(\xi),\diamond}^{I} (\mathbb{U}_{(\xi),\diamond}^{I})^s$ has zero mean. As in the oscillation stress error, we decompose $\varrho_{(\xi),\diamond}^{I} (\mathbb{U}_{(\xi),\diamond}^{I})^s$, applying the synthetic Littlewood-Paley decomposition suggested in \eqref{eq:decomp:showing}, and set for $q+\half +1<m <q+\bn$, 
\begin{subequations}
\begin{align}
     S_C^{q+\half+1}
    &:= (\divH + \divR)
    \left[
    \sum_{\diamond,i,j,k,\xi,\vecl,I}
(C_{(\xi), \diamond}^{1, I}+C_{(\xi), \diamond}^{2, I})^{\bullet s}\tP_{q+\half+1} \left( \varrho_{(\xi),\diamond}^{I} (\mathbb{U}_{(\xi),\diamond}^{I})^s \right) \circ \Phi_{(i,k)}\right] \label{div.err.low1}
\end{align}
\begin{align}
S_C^{m}
    &:= (\divH + \divR)
    \left[
    \sum_{\diamond,i,j,k,\xi,\vecl,I}
(C_{(\xi), \diamond}^{1, I}+C_{(\xi), \diamond}^{2, I})^{\bullet s}\tP_{(m-1,m]} \left( \varrho_{(\xi),\diamond}^{I} (\mathbb{U}_{(\xi),\diamond}^{I})^s \right) \circ \Phi_{(i,k)}\right] \label{div.err.mid1}
\end{align}
\begin{align}
&S_C^{q+\bn}:= w_{q+1}^{(c)}\otimes w_{q+1}^{(c)} \label{div.err.high3}\\
&+\sum_{\diamond,i,j,k,\xi,\vecl,I}
  a_{(\xi),\diamond} \left( \rhob_{(\xi)}^\diamond \zetab_{\xi}^{I,\diamond} \varrho_{(\xi),\diamond}^{I} \right)\circ \Phiik \xi^\ell \bigl( A_{\ell}^m  \epsilon_{\bullet p r} + A_{\ell}^\bullet  \epsilon_{m p r} \bigr) a_{(\xi),\diamond}^{p, \rm good} \partial_r \Phi_{(i,k)}^s (\mathbb{U}_{(\xi),\diamond}^{I})^s \circ \Phi_{(i,k)} \label{eq:div:corrector:form:2}\\
& +\sum_{m=q+\bn}^{q+\bn+1}(\divH + \divR)
    \left[
    \sum_{\diamond,i,j,k,\xi,\vecl,I}
(C_{(\xi), \diamond}^{1, I}+C_{(\xi), \diamond}^{2, I})^{\bullet s}(\tP_{(m-1,m]} + \Id - \tP_{q+\bn+1}) \left( \varrho_{(\xi),\diamond}^{I} (\mathbb{U}_{(\xi),\diamond}^{I})^s \right) \circ \Phi_{(i,k)}\right]\, .  \label{div.err.high1}
\end{align}
\end{subequations}
Here, the terms involved with the operators $\divR$ or $\Id -\tP_{q+\bn+1}$ will go into the nonlocal part and all the remaining terms will be included in the local parts. 

The conclusions of Lemma \ref{l:divergence:corrector:error} 
for the terms \eqref{div.err.low1}, \eqref{div.err.mid1}, and the terms involving $\tP_{(m-1,m]}$ in \eqref{div.err.high1}
follow similarly to \texttt{Case 3} from the proof of Lemma \ref{lem:oscillation:general:estimate}. Indeed, we fix indices $i,j,k,\xi,\vecl,I, s$, $\diamond=R$, and apply Proposition \ref{prop:intermittent:inverse:div} to
\begin{align*}
G^\bullet_R 
= \la_{q+\bn}^{-1}(C_{(\xi), R}^{1, I}+C_{(\xi), R}^{2, I})^{\bullet s}, \quad
\varrho_R
=\begin{cases}
\la_{q+\bn}\tP_{q+\half+1} \left( \varrho_{(\xi),R}^{I} (\mathbb{U}_{(\xi),R}^{I})^s\right) &\text{for \eqref{div.err.low1}} \\
\la_{q+\bn}\tP_{(m-1,m]} \left( \varrho_{(\xi),R}^{I} (\mathbb{U}_{(\xi),R}^{I})^s\right) &\text{for \eqref{div.err.mid1}, \eqref{div.err.high1}}, \\
\end{cases}
\end{align*}
with the same choice of the rest of parameters as in \texttt{Case 3}. In the case of $\diamond=\ph$, as in \texttt{Case 3}, $G_\ph$ and $\varrho_\ph$ will have extra $r_q^{\sfrac23}$ and $r_q^{-\sfrac23}$, respectively, with the replacement of $R$ with $\ph$ in $C_{(\xi), R}^{1, I}$, $C_{(\xi), R}^{1, I}$, and $\varrho_{(\xi),R}^{I} (\mathbb{U}_{(\xi),R}^{I})^s$. The assumptions in \eqref{eq:inverse:div:DN:G} and \eqref{eq:DN:Mikado:density} of Proposition \ref{prop:intermittent:inverse:div} can be verified using Lemma~\ref{lem:a_master_est_p}, Lemma~\ref{lem:special:cases}, Lemma~\ref{lem:LP.est},
item~\eqref{item:pipe:5} from Proposition~\ref{prop:pipeconstruction} and item~\eqref{item:pipe:5:current} from Proposition~\ref{prop:pipe.flow.current}.\footnote{Note that we have traded $\la_\qbn$ between $G_R^\bullet$ and $\rho_R$ so that the parameter choices are the same as the oscillation error.  We also note that thanks to the extra gain $\sfrac{\la_{q+\half}}{\la_{q+\bn}}$ in the estimate of $G_R$ and $G_\ph$ compared with \texttt{Case 3}, all the error terms are actually
small enough in amplitude to absorbed into the highest shell. The only reason to use the synthetic Littlewood-Paley decomposition here is to ensure that we can upgrade material derivatives via dodging later.} The rest of the assumptions follow exactly as in \texttt{Case 3} from the proof of Lemma~\ref{lem:oscillation:general:estimate}. We note now that the support of the low-frequency function $G$ is the same as in the oscillation error due to the presence of $\rhob_\pxi^\diamond\zetab_\xi^\diamond$ and their derivatives.  In addition, the support of the high-frequency potentials is the same as in the oscillation error since $\UU_{\pxi,\diamond}^I$ and $\varrho_{\pxi,\diamond}^I$ are both supported in a $2\la_\qbn^{-1}$ neighborhood of the pipe potential from \eqref{eq:WW:explicit} and item~\eqref{item:pipe:3.5}. Finally, to deal with the remaining term in \eqref{div.err.high1}, we may use the same type of arguments as in \texttt{Case 4} in the proof of Lemma \ref{lem:oscillation:general:estimate}. For the sake of both the readers and authors, we omit these details.
  
Lastly, we consider \eqref{div.err.high3} and \eqref{eq:div:corrector:form:2}, which are absorbed into $S_C^{q+\bn,l}$. From Lemma~\ref{lem:dodging}, we have that
\begin{align}
w_{q+1}^{(c)} \otimes w_{q+1}^{(c)} = \sum_{\diamond,i,j,k,\xi,\vecl,I} &  \left( \nabla \left( a_{(\xi),\diamond} (\rhob_{(\xi)}^\diamond \zetab^{I,\diamond}_\xi)\circ \Phiik \right) \times \left( \nabla\Phi_{(i,k)}^{T}\UU_{(\xi),\diamond}^{I} \circ \Phi_{(i,k)} \right) \right)     \notag\\
& \qquad \otimes \left( \nabla \left( a_{(\xi),\diamond} (\rhob_{(\xi)}^\diamond \zetab^{I,\diamond}_\xi)\circ \Phiik \right) \times \left( \nabla\Phi_{(i,k)}^{T} \UU_{(\xi),\diamond}^{I} \circ \Phi_{(i,k)} \right) \right) \, . \label{eq:div:corrector:form:1}
\end{align}
It follows immediately from estimate~\eqref{eq:w:oxi:cs} with $r=3,\infty$, \eqref{eq:inductive:partition} at level $q$, and Lemma~\ref{lemma:cumulative:cutoff:Lp} with $r_1=\infty,r_2=1$ that for $N,M\leq \sfrac{\Nfin}{10}$,
\begin{align}
\norm{\psi_{i,q} D^N \Dtq^M \left( w_{q+1}^{(c)} \otimes w_{q+1}^{(c)} \right)  }_{\infty} & \lesssim \Gamma_q^{\badshaq+9} \lambda_{q+\bn}^{N}\MM{M,\NindSmall, \tau_{q}^{-1}\Gamma_{q}^{i+15}, \Tau_{q}^{-1}\Ga_q^8} \notag \\
\norm{\psi_{i,q} D^N \Dtq^M \left( w_{q+1}^{(c)} \otimes w_{q+1}^{(c)} \right)  }_{\sfrac 32}^{\sfrac 32}
&\les r_{q}^{2} \sum_{\diamond,i,j,k,\xi,\vecl,I} \left|  \supp\left( \eta_{i,j,k,\xi,\vecl,\diamond} \zetab^{I,\diamond}_\xi \right) \right| \delta_{q+\bn}^{\sfrac 32} \Gamma_{q}^{3j+21} 
\lambda_{q+\bn}^{\sfrac{3N}2} \notag\\
&\qquad \qquad \times (\MM{M,\NindSmall, \tau_{q}^{-1}\Gamma_{q}^{i+13}, \Tau_{q}^{-1}\Ga_q^8})^{\sfrac32} \notag \\
&\lesssim r_q^2 \delta_{q+\bn}^{\sfrac 32} \Gamma_q^{30} \lambda_{q+\bn}^{\sfrac{3N}2} (\MM{M,\NindSmall, \tau_{q}^{-1}\Gamma_{q}^{i+15}, \Tau_{q}^{-1}\Ga_q^8})^{\sfrac32} \,. \notag 
\end{align}
The estimate for the $L^\infty$ norm matches \eqref{eq:div:corrector:Linfty} for $m=q+\bn$ after using \eqref{eq:par:div:2}.  For the $L^{\sfrac 32}$ estimate, taking cube roots and using the parameter inequality \eqref{eq:par:div:1} matches \eqref{eq:div:corrector:L1} for $m=q+\bn$. Finally, we have that the support of this error term is contained in $w_{q+1}$; then \eqref{eq:dc:ER:supp} is immediate from Lemma~\ref{lem:dodging}. On the other hand, one can observe that \eqref{eq:div:corrector:form:2} enjoys the exact same properties as $w_{q+1}^{(c)}\otimes w_{q+1}^{(c)}$, and hence we get the desired conclusion in a similar way.
\end{proof}

\begin{lemma*}[\bf Pressure increment]\label{lem:corrector:pressure}
For every $q+\half+1\leq m \leq \qbn$, there exists a function $\si_{S^{m}_C} = \si_{S^{m}_C}^+ - \si_{S^{m}_C}^-$ such that the following hold.\index{$\sigma_{S^m_C}$}
\begin{enumerate}[(i)]
    \item We have that for all $N,M < \Nfin/100$ and $q+\half+1\leq m \leq \qbn-1$,
\begin{subequations}
\begin{align}
    \label{eq:c.p.1}
    \left|\psi_{i,q} D^N \Dtq^M S^{m,l}_{C}\right| &< \left(\si_{S^{{m}}_C}^+  + \de_{q+3\bn}\right) \left(\lambda_{m}\Gamma_q\right)^N \MM{M,\Nindt,\tau_q^{-1}\Gamma_{q}^{i+16},\Tau_q^{-1}\Ga_q^9} \, , \\
    \label{eq:c.p.1.5}
    \left|\psi_{i,q} D^N \Dtq^M S^{q+\bn,l}_{C}\right| &< \left(\si_{S^{q+\bn}_C}^+ + \si_{\upsilon}^+ + \de_{q+3\bn}\right) \left(\lambda_\qbn\Gamma_\qbn\right)^N \MM{M,\Nindt,\tau_q^{-1}\Gamma_{q}^{i+16},\Tau_q^{-1}\Ga_q^9} \, ,
\end{align}
\end{subequations}
where $\si_{\upsilon}^+$ is defined as in \eqref{defn:si.upsilon}. Furthermore, for any integer $q+\half <m\leq q+\bn$ and for all $N,M < \Nfin/100$,
\begin{subequations}
\begin{align}    
    \label{eq:c.p.2}
    \left|\psi_{i,q} D^N \Dtq^M \si_{S^{m}_C}^+\right| &< \left(\si_{S^{m}_C}^+ +\de_{q+3\bn}\right)  \left(\lambda_{m}\Gamma_q\right)^N \MM{M,\Nindt,\tau_q^{-1}\Gamma_{q}^{i+16},\Tau_q^{-1}\Ga_q^9}\\
    \label{eq:c.p.3}
    \norm{\psi_{i,q} D^N \Dtq^M \si_{S^{m}_C}^+}_{\sfrac32} &\leq \Ga_m^{-9} \de_{m+\bn} \left(\lambda_{m}\Gamma_q\right)^N \MM{M,\Nindt,\tau_q^{-1}\Gamma_{q}^{i+16},\Tau_q^{-1}\Ga_q^9}\\
    \label{eq:c.p.3.5}
    \norm{D^N \Dtq^M \si_{S^{m}_C}^+}_{\infty} &\leq \Ga_m^{\badshaq-9}  \left(\lambda_{m}\Gamma_q\right)^N \MM{M,\Nindt,\tau_q^{-1}\Gamma_{q}^{i+16},\Tau_q^{-1}\Ga_q^9}\\
    \label{eq:c.p.4}
    \left|\psi_{i,q} D^N \Dtq^M \si_{S^{m}_C}^-\right| &<
    \Ga_{q+\half}^{-100}\pi_q^{q+\half}
  \left(\lambda_{q+\half}\Gamma_q\right)^N \MM{M,\Nindt,\tau_q^{-1}\Gamma_{q}^{i+16},\Tau_q^{-1}\Ga_q^9} \, .
\end{align}
\end{subequations}
\item For $q+\half+1\leq m \leq \qbn$, we have that
\begin{align}\label{eq:c.p.6}
   \begin{split}
B\left( \supp \hat w_{q'}, \lambda_{q'}^{-1} \Gamma_{q'+1} \right) \cap \supp (\si_{S^{m}_C}^+) &= \emptyset \qquad \forall q+1\leq q' \leq  m-1 \\
B\left( \supp \hat w_{q'}, \lambda_{q'}^{-1} \Gamma_{q'} \right) \cap \supp (\si_{S^{m}_C}^-) &= \emptyset \qquad \forall q+1\leq q' \leq  q+\half \, .
\end{split}
\end{align}
\item\label{item:thursday:night} Define 
\begin{equation}\label{def:bmu:ER:C}
    \bmu_{\sigma_{S_C^m}}(t) = \int_0^t \left \langle \Dtq \sigma_{S_C^m}  \right \rangle (s) \, ds \, .
\end{equation}
Then we have that for $0\leq M\leq 2\Nind$,
    \begin{align}\label{th:billys:1}
      \left|\frac{d^{M+1}}{dt^{M+1}} \bmu_{\sigma_{S_C^m}} \right| 
      \leq (\max(1, T))^{-1}\delta_{q+3\bn} \MM{M,\Nindt,\tau_q^{-1},\Tau_{q+1}^{-1}} \, .
    \end{align}
\end{enumerate}
\end{lemma*}

\medskip

\begin{lemma*}[\bf Pressure current]\label{lem:corrector:pressure:current}
For every $q+\half < m \leq q+\bn$, there exists a current error $\phi_{S_C^{m}}$ associated to the pressure increment $\sigma_{S_C^{m}}$ defined by Lemma~\ref{lem:corrector:pressure} which satisfies the following properties.
\begin{enumerate}[(i)]
    \item We have the decompositions and equalities
\begin{subequations}
    \begin{align}\label{eq:desert:decomp:ER:C}
        \phi_{{S_C^{m}}} &= \phi_{{S_C^{m}}}^* + \sum_{k=q+\half+1}^{{m}} \phi_{{S_C^{m}}}^{k} \, , \quad 
         \phi_{S_C^{m}}^{k} = \phi_{S_C^{m}}^{k,l} + \phi_{S_C^{m}}^{k,*} \, \\
         \div \phi_{{S_C^{m}}}
          &= D_{t,q}\si_{{S_C^{m}}} - \langle D_{t,q}\si_{{S_C^{m}}}\rangle\, . 
    \end{align}
\end{subequations}
\item For $q+\half+1\leq k \leq m$ and $N,M\leq 2\Nind$,
\begin{subequations}
\begin{align}
    &\left|\psi_{i,q} D^N \Dtq^M \phi_{S_C^{m}}^{k,l}\right| < \Ga_{ k}^{-100} r_{ k}^{-1} \left(\pi_q^{ k}\right)^{\sfrac32}\left(\lambda_{ k}\Gamma_q\right)^N \MM{M,\Nindt,\tau_q^{-1}\Gamma_{q}^{i+16},\Tau_q^{-1}\Ga_q^9}\label{c:p:c:pt} \\
    &\norm{D^N \Dtq^M \phi_{S_C^{m}}^{k,*}}_{L^\infty}
    \leq 
    \delta_{q+3\bn}^{\sfrac 32} \Tau_{q+\bn}^{2\Nindt} {\la_m^N} \tau_q^{-M} \, . \label{e:c:p:c:nonlocal}
\end{align}
\end{subequations}
\item For all $q+\half+1\leq k \leq m$ and all $q+1\leq q' \leq k-1$, 
\begin{align}
        B\left( \supp \hat w_{q'}, \sfrac 12 \lambda_{q'}^{-1} \Ga_{q'+1} \right) \cap \supp \left( \phi^{k,l}_{S_C^{m}} \right) = \emptyset \label{c:p:c:supp} \, .
\end{align}
\end{enumerate}
\end{lemma*}

\begin{proof}[Proofs of Lemmas \ref{lem:corrector:pressure}-\ref{lem:corrector:pressure:current}]
\texttt{Case 0: }pressure for \eqref{div.err.low1}, \eqref{div.err.mid1}, and \eqref{div.err.high1}. The pressure increment and the current error associated to each piece in the local part of \eqref{div.err.low1}, \eqref{div.err.mid1}, and \eqref{div.err.high1} can be constructed in the same way as in Lemma \ref{lem:oscillation:pressure}-\ref{lem:oscillation:pressure:current}. Indeed, the proof relies on Proposition \ref{lem.pr.invdiv2}, and  $(G_R,\varrho_R)$, $(G_\ph, \varrho_\ph)$ given in the proof of Lemma \ref{l:divergence:corrector:error} have the exact same properties required in the proposition as the one given in \texttt{Case 3} of the proof of Lemma \ref{lem:oscillation:general:estimate}. In particular, the preliminary assumptions \eqref{i:st:sample:5} holds with $\bar \pi$ given as in \eqref{eq:pi:osc:lowest} due to \eqref{e:a_master_est_p_pointwise}. Therefore, we get the same conclusions by repeating the same arguments. In particular, all conclusions from Lemma \ref{lem:corrector:pressure}--\ref{lem:corrector:pressure:current} are obtained in the cases $m<q+\bn$. Furthermore, when $m=q+\bn$, we denote the pressure increment and the current error associated to \eqref{div.err.high1} by $\si_{\eqref{div.err.high1}}= \si_{\eqref{div.err.high1}}^{+}- \si_{\eqref{div.err.high1}}^{-}$ and $\phi_{\eqref{div.err.high1}}^{k}
=\phi_{\eqref{div.err.high1}}^{k,l}
+\phi_{\eqref{div.err.high1}}^{k,*}
$, respectively. Since these error terms are defined using the same parameter choices as the oscillation error, we obtain estimates consistent with \eqref{eq:c.p.2}--\eqref{c:p:c:supp} for these error terms.  We note also that we obtain a version of \eqref{eq:c.p.1.5} which does not require the introduction of $\sigma_\upsilon^+$ on the right-hand side; later error terms will require $\sigma_\upsilon^+$.
\medskip

\noindent\texttt{Case 1: }\eqref{div.err.high3} needs no new pressure increment. From
\eqref{est.vel.inc.c.by.pr}, we have that
\begin{align*}
\left|\psi_{i,q} D^N \Dtq^M \eqref{div.err.high3}\right| &\lec \Ga_q^{-2}( \si_{\upsilon}^+ + \de_{q+3\bn}) \left(\lambda_\qbn\Gamma_\qbn\right)^N \MM{M,\Nindt,\tau_q^{-1}\Gamma_{q}^{i+16},\Tau_q^{-1}\Ga_q^9}\,  
\end{align*}
for $N,M\leq \sfrac{\Nfin}{100}$. This estimate is consistent with \eqref{eq:c.p.1.5}, and since no pressure increment is created here, we need not check any of the conclusion in \eqref{eq:c.p.2}--\eqref{eq:c.p.6}.

\noindent\texttt{Case 2: }pressure for \eqref{eq:div:corrector:form:2}. The general idea for this error term is that since it is given as a product of two slightly altered velocity increments, we can apply Proposition~\ref{lem:pr.vel} (which was used to construct pressure increments for velocity increments already in subsection~\ref{sec:vel.inc}) to construct pressure increments $\si_{\eqref{eq:div:corrector:form:2}}^\pm$ and current errors $\phi_{\eqref{eq:div:corrector:form:2}}^{k}$. So we fix the indices $i,j,k,\xi,\vecl,I,\diamond$ and apply Proposition~\ref{lem:pr.vel} to the functions $\hat\upsilon_{b,\diamond} = \hat\upsilon_{b, i,j,k,\xi,\vecl,I,\diamond}$ defined by $\hat \upsilon_{b,\diamond} = G_{b\diamond}\rho_{b\diamond}\circ\Phiik$, $b=1,2$, where
\begin{align*}
    \hat\upsilon_{1,\diamond}&:= r_q^{\sfrac13}{\la_q^{\sfrac13}\la_{q+\bn}^{-\sfrac13}}
    a_{(\xi),\diamond} \left( \rhob_{(\xi)}^\diamond \zetab_{\xi}^{I,\diamond} \varrho_{(\xi),\diamond}^{I} \right)\circ \Phiik\\
    \hat\upsilon_{2,\diamond}&:= r_q^{-\sfrac13}{\la_q^{-\sfrac13}\la_{q+\bn}^{\sfrac13}}\xi^\ell \bigl( A_{\ell}^m  \epsilon_{\bullet p r} + A_{\ell}^\bullet  \epsilon_{m p r} \bigr) a_{(\xi),\diamond}^{p, \rm good} \partial_r \Phi_{(i,k)}^s (\mathbb{U}_{(\xi),\diamond}^{I})^s \circ \Phi_{(i,k)} \\
    G_{1R} &:= {\la_q^{\sfrac13}\la_{q+\bn}^{-\sfrac13}} a_{(\xi),R} \left( \rhob_{(\xi)}^R \zetab_{\xi}^{I,R}  \right)\circ \Phiik, \quad
    \rho_{1R} 
    := r_q^{\sfrac13}\varrho_{(\xi),R}^{I} \\
    G_{1\ph} &:= r_q^{\sfrac13}{\la_q^{\sfrac13}\la_{q+\bn}^{-\sfrac13}} a_{(\xi),\ph} \left( \rhob_{(\xi)}^\ph \zetab_{\xi}^{I,\ph}  \right)\circ \Phiik, \quad
    \rho_{1\ph} 
    := \varrho_{(\xi),\ph}^{I} \\
    G_{2R} &:= r_q^{-\sfrac23}{\la_q^{-\sfrac13}\la_{q+\bn}^{\sfrac13}}\la_{q+\bn}^{-1}\xi^\ell \bigl( A_{\ell}^m  \epsilon_{\bullet p r} + A_{\ell}^\bullet  \epsilon_{m p r} \bigr) a_{(\xi),R}^{p, \rm good} \partial_r \Phi_{(i,k)}^s, \quad
    \rho_{2R} :=    
    r_q^{\sfrac13}\la_{q+\bn}(\mathbb{U}_{(\xi),R}^{I})^s \\
    G_{2\ph} &:= r_q^{-\sfrac13}{\la_q^{-\sfrac13}\la_{q+\bn}^{\sfrac13}}\la_{q+\bn}^{-1}\xi^\ell \bigl( A_{\ell}^m  \epsilon_{\bullet p r} + A_{\ell}^\bullet  \epsilon_{m p r} \bigr) a_{(\xi),\ph}^{p, \rm good} \partial_r \Phi_{(i,k)}^s, \quad
    \rho_{2\ph} :=    
    \la_{q+\bn}(\mathbb{U}_{(\xi),\ph}^{I})^s  \, .
\end{align*}
We then set the following choices for the application of Proposition~\ref{lem:pr.vel}:
\begin{align*}
N_* &= M_* = \sfrac{\Nfin}{10}, \quad M_t = \Nindt, \quad M_\circ = N_\circ = 2\Nind , \quad K_\circ \textnormal{ as in \eqref{i:par:9.5}} \, ,
\\    
\Phi &= \Phiik, \quad v = \hat u_q, \quad D_t = D_{t,q}, \quad \la' = \La_q,\quad
    \nu' = \Tau_q^{-1}\Ga_8, \quad 
      \quad\const_v = \La_q^{\sfrac12} \, , \quad \Ga = \Ga_q^{\sfrac {1}{10}}, \\
    \const_{G,3}
    &=\left| \supp \left( \eta_{i,j,k,\xi,\vecl,\diamond} \zetab_\xi^{I,\diamond} \right) \right|^{\sfrac 13} (\delta_{q+2\bn}\Ga_{q+\bn}^{-20})^{\sfrac 12} \Gamma_q^{j} + \la_{q+2\bn}^{-10}, \quad
    \const_{G,\infty}
    =\Gamma_{q+\bn}^{\frac{\badshaq}{2}-20}r_q^{\sfrac23}, 
    \quad \pi= \pi_\ell \Ga_q^{30} r_q^{-\sfrac23}\La_q^{\sfrac23}\la_{q+\bn}^{-\sfrac23},
    \\
    \const_{\rho,3}&:= 1, \quad \const_{\rho,\infty}= r_q^{-\sfrac23}, \quad
    \la= \la_{q+\half}, \quad\La = \la_{q+\bn},\quad \nu = \tau_q^{-1}\Ga_q^{i+13},  \quad r_G = r_{\hat\upsilon}={1}, \quad \mu = \la_{q+\bn}r_q\\
    \de_{\rm tiny} &=  {\de_{q+3\bn} },\quad
    \bar m = m+1-(q+\half), \quad \mu_0 = \la_{q+\half+1}, \quad \mu_1 = \la_{q+\half+\sfrac32},
    \quad \mu_k = \la_{q+\half+k}\,,\\
    \NcutLarge, &\NcutSmall \textnormal{ as in \eqref{i:par:6}} \, , \quad \Ndec \textnormal{ as in \eqref{i:par:9}} \, , \quad  \dpot, N_{**}\textnormal{ as in \eqref{i:par:10}} \, .
\end{align*}
First, the verification of the assumptions from part 1 of Proposition~\ref{lem:pr.vel} can be done in a similar manner as in the proofs of Lemmas~\ref{lem:pr.inc.vel.inc.pot} and \ref{lem:pr.current.vel.inc}.  We omit further details, but note that in this case, the intermittency parameters are chosen as $1$ and $G$ has extra factor $\la_q^{\sfrac13}\la_{q+\bn}^{-\sfrac 13}$ instead. From the definitions, the support properties of the low frequency functions $G_{b\diamond}$ and the high frequency functions $\rho_{b\diamond}$ are essentially the same as those of the corresponding functions in Lemmas~\ref{lem:pr.inc.vel.inc.pot} and \ref{lem:pr.current.vel.inc}.

As a consequence of \eqref{est.w.by.pr}, we have pressure increments associated to $\hat\upsilon_{b,\diamond}$, $b=1,2$, which satisfies
\begin{align*}
    \left|D^N \Dtq^M \hat \upsilon_{b,\diamond}\right| 
    \lec (\si^+_{\hat\upsilon_{b,\diamond}} + \de_{q+3\bn})^{\sfrac12} (\la_{q+\bn}\Ga_q)^N \MM{M, \Nindt, \tau_q^{-1}\Ga_q^{i+15}, \Tau_q^{-1}\Ga_q^9}
\end{align*}
for any $N, M\leq \sfrac{\Nfin}{10}$. This implies that
\begin{align*}
\left|D^N \Dtq^M (\hat \upsilon_{1,\diamond} \hat \upsilon_{2,\diamond})\right|
\lec (\si^+_{\hat\upsilon_{1,\diamond}} + \si^+_{\hat\upsilon_{2,\diamond}} + \de_{q+3\bn}) (\la_{q+\bn}\Ga_q)^N \MM{M, \Nindt, \tau_q^{-1}\Ga_q^{i+15}, \Tau_q^{-1}\Ga_q^9}
\end{align*}
for any $N, M\leq \sfrac{\Nfin}{10}$. Then appealing to the same conclusions used in \eqref{est.pr.vel.inc:piece}--\eqref{est.pr.vel.inc-:piece}, we have that
\begin{align*}
 \left|D^N D_{t,q}^M \si_{\hat\upsilon_b}^{+}\right|
    &\lec (\si_{\hat\upsilon_b}^{+}+ \de_{q+3\bn})  ( \la_{q+\bn}\Ga_q)^{N} \MM{M,\Nindt,\tau_q^{-1}\Ga_q^{i+15},\Tau_q^{-1}\Gamma_q^9}\\
    \norm{D^N D_{t,q}^M \si_{\hat\upsilon_b}^{+}}_{\sfrac32}
    &\lesssim 
    \left[\left| \supp \left( \eta_{i,j,k,\xi,\vecl,\diamond} \zetab_\xi^{I,\diamond} \right) \right|^{\sfrac 23} \delta_{q+2\bn}\Gamma_{q+\bn}^{-20} \Gamma_q^{2j} + \de_{q+3\bn}\right] \nonumber\\
    &\qquad \times 
    ( \la_{q+\bn}\Ga_q)^{N} \MM{M,\Nindt,\tau_q^{-1}\Ga_q^{i+14},\Tau_q^{-1}\Gamma_q^9} \\
\norm{D^N D_{t,q}^M \si_{\hat\upsilon_b}^{+}}_{\infty} 
    &\lesssim \Gamma_{q+\bn}^{\badshaq-40} ( \la_{q+\bn}\Ga_q)^{N} \MM{M,\Nindt,\tau_q^{-1}\Ga_q^{i+14},\Tau_q^{-1}\Gamma_q^9} \,\\
    \left|D^N D_{t,q}^M \si_{\hat\upsilon_b}^{-}\right|
    &\lec  
    \pi_\ell \Gamma_q^{41} \la_q^{\sfrac13}\la_{q+\bn}^{-\sfrac13} (\la_{q+\half}\Ga_q)^N\MM{M,\Nindt,\tau_q^{-1}\Ga_q^{i+14},\Tau_q^{-1}\Gamma_q^9}
    \end{align*}
for all $N,M \leq \sfrac{\Nfin}{100}$. We reintroduce the indices $i,j,k,\xi,\vecl,I$ and define the pressure increment associated to \eqref{eq:div:corrector:form:2} by
\begin{align*}
\si_{\eqref{eq:div:corrector:form:2}}^{\pm} : = 
\sum_{i,j,k,\xi,\vecl,I,b,\diamond} \si^{\pm}_{\hat \upsilon_{b,i,j,k,\xi,\vecl,I,\diamond}} .
\end{align*}
The estimates \eqref{eq:c.p.1} and \eqref{eq:c.p.2} associated to \eqref{eq:div:corrector:form:2} follow using an aggregation procedure identical to that used in the proofs of Lemmas~\ref{lem:pr.inc.vel.inc.pot} and \ref{lem:pr.current.vel.inc}, and so we omit further details.

Lastly, we define $\phi_{\eqref{eq:div:corrector:form:2}}^{k,l}$ and $\phi_{\eqref{eq:div:corrector:form:2}}^{k,*}$ as in the proofs of Lemmas~\ref{lem:pr.inc.vel.inc.pot} and \ref{lem:pr.current.vel.inc} and obtain \eqref{c:p:c:pt}, \eqref{e:c:p:c:nonlocal}, and \eqref{c:p:c:supp} as in the cited Lemmas.  Setting 
\begin{align*}
    \si_{S^{q+\bn}_C}^{\pm}
    :=\si_{\eqref{div.err.high1}}^{\pm}+    \si_{\eqref{eq:div:corrector:form:2}}^{\pm}, \quad
    \phi_{S^{q+\bn,l}_C}^{k}
    :=\phi_{\eqref{div.err.high1}}^{k}+    \phi_{\eqref{eq:div:corrector:form:2}}^{k}
\end{align*}
and collecting the properties of these objects obtained above, we conclude \eqref{eq:c.p.2}--\eqref{c:p:c:supp} and \eqref{eq:c.p.1.5}.
\end{proof}

\subsection{Mollification error \texorpdfstring{$S_M$}{essemm}}\label{ss:essemm}

Recalling from subsection \ref{sec:error.ER} that $\div S_{M2}$ has mean-zero, we use Proposition~\ref{prop:intermittent:inverse:div}, Remark~\ref{rem:inverse.div.spcial} to first define the mollification error $S_M = S_{M1} + S_{M2}$ by
\begin{align}
    S_{M1}&:=R_q^q - R_\ell + \left( \pi_\ell - \pi_q^q \right) \Id =:S_{M}^{q+1, *}
   \label{ER:new:error:moll1} \\
    S_{M2}&:= \divR  \left[(\pa_t + \hat u_{q}\cdot \na )(\hat w_{q+\bn}-w_{q+1})
    + (\hat w_{q+\bn}-w_{q+1})\otimes \hat u_q \right] + \hat w_{q+\bn}\otimes \hat w_{q+\bn} - w_{q+1}\otimes w_{q+1}
    =:S_{M}^{q+\bn, *}
    \, . \notag
\end{align}
For the undefined mollification stress errors $S_{M}^{k,l}$, $S_{M}^{k,*}$, we set them as zero.

\begin{lemma}[\bf Basic estimates and applying inverse divergence]\label{lem:moll:ER} The mollification error 
$S_{M}^{q+1,*}$ and
$S_{M}^{q+\bn,*}$ satisfy\index{$S_M$}
\begin{subequations}
    \begin{align}
        \norm{D^N D_{t,q}^M S_{M}^{q+1, *}}_\infty
        &\leq \Ga_{q+1}^9 \delta_{q+3\bn} \Tau_{q+1}^{2\Nindt} \left(\lambda_{q+1}\Gamma_{q+1}\right)^N \MM{M, \NindRt, \tau_{q}^{-1}, \Tau_{q}^{-1} } \, . 
    \label{est:stress.mollification1}\\
        \norm{ D^N D_{t,q+\bn-1}^M S_{M}^{q+\bn, *}}_\infty
        &\leq \Ga_\qbn^9 \delta_{q+3\bn}\Tau_\qbn^{2\Nindt} \left(\lambda_{q+\bn}\Gamma_{q+\bn}\right)^N \MM{M, \NindRt, \tau_{q+\bn-1}^{-1}, \Tau_{q+\bn-1}^{-1} } \, . 
    \label{est:stress.mollification}
    \end{align}
\end{subequations}
for all $N+ M\leq 2\Nind$.
\end{lemma}
\begin{proof}[Proof of Lemma~\ref{lem:moll:ER}]
From \eqref{eq:Rcomm:bounds}, we have
\begin{align*}
     \left\| D^N \Dtq^M S_{M} \right\|_{\infty} 
    \les \Gamma_{q+1} \Tau_{q+1}^{2\Nindt} \delta_{q+3\bn}^2 \lambda_{q+1}^N \MM{M,\Nindt, \tau_{q}^{-1},\Gamma_{q}^{-1}\Tau_{q}^{-1}}
\end{align*}
for all $N+M\leq 2\Nind$, which immediately leads to \eqref{est:stress.mollification1}.
    
To deal with $S_{M2}$, we recall from \eqref{eq:diff:moll:vellie:statement} that
    \begin{align*}
    \norm{D^N D_{t,q+\bn-1}^M \left(w_{q+1}- \hat w_{q+\bn}\right)}_\infty \lec \delta_{q+3\bn}^3 \Tau_\qbn^{25\NindRt} \left(\lambda_{q+\bn}\Gamma_{q+\bn-1}\right)^N \MM{M, \NindRt, \tau_{q+\bn-1}^{-1}, \Tau_{q+\bn-1}^{-1} }  \, .
\end{align*}
for all $N+M\leq \sfrac{\Nfin}{4}$. Using Lemma \ref{lem:dodging}, we note that $D_{t,q-\bn-1}w_{q+1} = D_{t,q}w_{q+1}$ and $D_{t,q-\bn-1}\hat w_{q+\bn} = D_{t,q}\hat w_{q+\bn}$. Then, writing $\hat w_{q+\bn}\otimes \hat w_{q+\bn} - w_{q+1}\otimes w_{q+1}
=(\hat w_{q+\bn}-w_{q+1})\otimes \hat w_{q+\bn} + w_{q+1}\otimes (\hat w_{q+\bn}-w_{q+1}) $ and using \eqref{eq:vellie:upgraded:statement} and \eqref{eq:moll:vel:threesie},
we have
\begin{align}
    &\norm{\psi_{i, q+\bn-1}
    D^N D_{t,q+\bn-1}^M
    [\hat w_{q+\bn}\otimes \hat w_{q+\bn} - w_{q+1}\otimes w_{q+1}]
    }_\infty\notag\\
    &\qquad\leq \delta_{q+3\bn} \Tau_\qbn^{2\NindRt} \left(\lambda_{q+\bn}\Gamma_{q+\bn}\right)^N \MM{M, \NindRt, \tau_{q+\bn-1}^{-1}, \Tau_{q+\bn-1}^{-1} }\, ,
    \label{est:SM21.proof}
\end{align}
for all $N+M\leq 2\Nind$.

As for the remaining term, we first upgrade the material derivative in the estimate for $\hat u_q$.
Applying Lemma \ref{lem:upgrading.material.derivative} to $F^l = 0$, $F^* = \hat u_q$, $k=q+\bn$, $N_\star = \sfrac{3\Nfin}{4}$ with \eqref{eq:bob:Dq':old}, we get
\begin{align*}
    \norm{ D^N D_{t,q+\bn-1}^M \hat u_q}_\infty
    \lec \Tau_{q}^{-1} \lambda_{q+\bn}^N \Tau_{q+\bn-1}^{-M}
\end{align*}
Here, we used \eqref{v:global:par:ineq}. Then, we use Remark \ref{rem:inverse.div.spcial} with \eqref{v:global:par:ineq}, setting 
\begin{align*}
    G &=D_{t,q+\bn-1} (\hat w_{q+\bn} - w_{q+1})
    \quad (\text{or } G= (\hat w_{q+\bn} - w_{q+1}) \otimes \hat u_q)
    , \quad v = \hat u_{q+\bn-1}
\\
    \const_{G, \infty}&= \de_{q+\bn}^3 \Tau_{q+\bn}^{20\Nindt}, \quad
    \la=\la' = \la_{q+\bn}\Ga_{q+\bn-1}, \quad 
    M_t = \Nindt,\quad
    \nu =\nu' = \Tau_{q+\bn}^{-1}, \quad \const_v = \La_{q+\bn-1}^{\sfrac12} \\
    &\hspace{2cm}N_* = \sfrac{\Nfin}{9}, \quad
M_* = \sfrac{\Nfin}{10}, \quad
N_\circ = M_\circ = 2\Nind\,  .
\end{align*}
As a result, with a suitable choice of positive integer $K_\circ$ to have
\begin{align*}
    \de_{q+\bn}^3 \Tau_{q+\bn}^{20\Nindt }\la_{q+\bn}^{5} 2^{2\Nind}
    \leq \la_{q+\bn}^{-K_\circ}
    \leq \de_{q+3\bn} \Tau_\qbn^{10\NindRt}\, ,
\end{align*}
we get
\begin{align}\label{est:SM22.proof}
    \norm{D^N D_{t,q+\bn-1}^M\divR ( D_{t,q} (\hat w_{q+\bn} - w_{q+1})) }_\infty
    &=
    \norm{D^N D_{t,q+\bn-1}^M\divR ( D_{t,q+\bn-1} (\hat w_{q+\bn} - w_{q+1})) }_\infty\\
    &\lec \de_{q+3\bn} \Tau_\qbn^{10\Nindt} (\la_{q+\bn}\Ga_\qbn)^N \Tau_{q+\bn}^{-M}\\
    &\leq \de_{q+3\bn}\Tau_\qbn^{2\Nindt}(\la_{q+\bn}\Ga_\qbn)^N \MM{M, \Nindt,\tau_{q+\bn}^{-1} ,\Tau_{q+\bn}^{-1}}\, ,  
\end{align}
for all $N+M\leq 2\Nind$. This completes the proof of \eqref{est:stress.mollification}.
\end{proof}

\subsection{Upgrading material derivatives and
Hypothesis \ref{hyp:dodging4}}\label{dodge4verification}

\begin{definition}[\bf Definition of $\ov R_{q+1}$ and $S_{q+1}^m$]\label{def:of:new:stresses} Recalling Lemma~\ref{lem:oscillation:general:estimate}, Lemma~\ref{l:transport:error}, Lemma~\ref{l:divergence:corrector:error}, and Lemma~\ref{lem:moll:ER}, we define $S_{q+1}^m := S^{m,l}_{q+1}+ S^{m,*}_{q+1}$ for all $q+1 \leq m \leq q+\bn$ by\index{$S_{q+1}$}
\begin{subequations}\label{defn:newstress}
\begin{align}
     S^{m,l}_{q+1} &:= S^{m,l}_{O} + S^{m,l}_{TN} + S^{m,l}_{C} + S^{m,l}_M \, ,\\
    S^{m,*}_{q+1} &:= S^{m,*}_{O} + S^{m,*}_{TN} + S^{m,*}_{C} + S^{m,*}_M\, .
\end{align}
\end{subequations}
Here, any undefined terms are taken to be $0$. We then define the primitive stress error $\ov R_{q+1}$ at $q+1$ step by
\begin{align}\label{defn:primitive.stress}
\ov R_{q+1} := \sum_{m=q+1}^{q+\bn} \ov R_{q+1}^{m} \, , \qquad
\ov R_{q+1}^m = R_q^m + S_{q+1}^m \, .
\end{align}
The local part $R_{q+1}^{m,l}$ and the non-local part $\ov R_{q+1}^{m,*}$ are defined by
\begin{align}\label{defn:local.stress}
R_{q+1}^{m,l} := R_q^{m,l} + S_{q+1}^{m,l}\, , \qquad \ov R_{q+1}^{m,*} := R_q^{m,*} + S_{q+1}^{m,*}\, .
\end{align}
\end{definition}
We note that by the above definition, we have that
\begin{align}
    \ov R_{q+1}^m = R_{q+1}^{m,l} + \ov R_{q+1}^{m,*} \,.
\end{align}
We sometimes also use the notation $\ov R_{q+1}^{m,l}$ to denote $R_{q+1}^{m,l}$, since it will be shown later that the local portion of $\ov R_{q+1}^{m,l}$ remains unchanged throughout the rest of the analysis.

\begin{lemma}[\bf Upgrading material derivatives and verifying Hypothesis~\ref{hyp:dodging4}]\label{l:divergence:stress:upgrading}  The new stress errors $S_{q+1}^m=S_{q+1}^{m,l}+S_{q+1}^{m,*}$ satisfy the following.
\begin{enumerate}[(i)]
    \item\label{upgrade:item:1} $R_{q+1}^{m,l}$ satisfies Hypothesis~\ref{hyp:dodging4} with $q$ replaced by $q+1$.
    \item\label{upgrade:item:2} For $q+2\leq m\leq q+\half$, the symmetric stresses $S_{q+1}^{m,l}$ obey the estimates
\begin{align}\label{eq:lo:upgrade:1}
    \left|\psi_{i,m-1} D^N D_{t,m-1}^M S_{q+1}^{m,l} \right| 
    &\lec \Ga^{-50}_{m} \pi_{q}^{m}
    \La_{m}^N \MM{M,\Nindt, \Ga_{m-1}^{i-5}\tau_{m-1}^{-1} , \Tau_{q}^{-1}\Ga_q^9}
\end{align}
for $N,M \leq \sfrac{\Nfin}{10}$.  For the same range of $N,M$, the symmetric stress $S_{q+1}^{q+1,l}$ obeys the estimates
\begin{align}\label{eq:lo:upgrade:2}
    \left|\psi_{i,q} D^N D_{t,q}^M S_{q+1}^{q+1,l} \right| 
    &\lec \Ga^{-50}_{{q+1}} \pi_{q}^{q+1}
    \La_{{q+1}}^N \MM{M,\Nindt, \Ga_{q}^{i+19}\tau_{q}^{-1}, \Tau_{q}^{-1}\Ga_q^9} \, .
\end{align}
\item\label{upgrade:item:3} For $q+\half+1\leq m \leq \qbn$ and $N,M\leq \sfrac{\Nfin}{100}$, the symmetric stresses $S_{q+1}^{m,l}$ obey the estimates
\begin{subequations}
\begin{align}
\left| \psi_{i,m-1} D^N D_{t,m-1}^M S^{m,l}_{q+1} \right| &\lesssim \left( \sigma_{S_O^m}^+ + \sigma_{S_C^{m,l}}^+ + \mathbf{1}_{\{m=\qbn\}} \left( \sigma_{S_{TN}}^+ + \sigma_\upsilon^+ \right) + \delta_{q+3\bn} \right) \notag\\
&\qquad \times (\lambda_{m}\Ga_m)^N \MM{M,\Nindt,\Gamma_{m-1}^{i-5} \tau_{m-1}^{-1}, \Tau_{q}^{-1} \Ga_q^9 }   \, .  \label{eq:stress:Linfty:upgraded}
\end{align}
\end{subequations}
\item\label{upgrade:item:4} For all $q+1\leq m \leq \qbn$ and $N+M\leq 2\Nind$, the symmetric stresses $S_{q+1}^{m,*}$
\begin{align}
\left\| D^N D_{t, m-1}^M S^{m,*}_{q+1} \right\|_{L^\infty}
    &\leq \Ga_{q+1}^2 \Tau_{q+1}^{2\Nindt} \delta_{q+3\bn}^2 \la_m^N \MM{M,\Nindt,\tau_{m-1}^{-1},\Tau_{m-1}^{-1}}  \, .
    \label{eq:nlstress:upgraded}
\end{align}
\end{enumerate}
\end{lemma}
\begin{proof}[Proof of Lemma~\ref{l:divergence:stress:upgrading}]
In order to prove the claim in item~\eqref{upgrade:item:1}, note that for the portion of $R_{q+1}^{m,l}$ coming from $R_q^{m,l}$ (c.f. \eqref{defn:primitive.stress}), the claim follows by the inductive hypothesis itself.  For the portion coming from $S_{q+1}^{m,l}$, we may appeal to \eqref{defn:newstress} and \eqref{osc:support:first}, \eqref{eq:trans:supp}, and \eqref{eq:dc:ER:supp}. 

Next, we may prove \eqref{eq:lo:upgrade:2} directly from \eqref{eq:lowshell:nopr:1}, since from Lemma~\ref{l:transport:error} and Lemma~\ref{l:divergence:corrector:error}, the transport, Nash, and divergence corrector errors do not contribute to $S_{q+1}^{q+1,l}$.  In order to prove \eqref{eq:lo:upgrade:1}, we note that from Lemma~\ref{l:transport:error} and Lemma~\ref{l:divergence:corrector:error}, the transport, Nash, and divergence corrector errors do not contribute to $S_{q+1}^{m,l}$ for $q+2\leq m \leq q+\half$. Then from Lemmas~\ref{lem:oscillation:general:estimate} and \ref{lem:osc:no:pressure}, we need only consider the case $m=q+\half$, for which we have that for $N,M\leq \sfrac{\Nfin}{10}$,
\begin{align}
    \left| \psi_{i,m-1} D^N D_{t,m-1}^M {S}^{m,l}_{q+1} \right| 
    &\underset{\eqref{eq:inductive:partition}}{=}
    \left| \psi_{i,m-1} \sum_{i'}\psi_{i',q}^6 D^N D_{t,m-1}^M {S}^{m,l}_{q+1} \right| \notag\\
    &\underset{\eqref{osc:support:first}}{\lec} \sum_{i':\psi_{i',q}\psi_{i,m-1}\neq 0 } \left| \psi_{i',q} D^N D_{t,q}^M {S}^{m,l}_{q+1} \right| \notag\\
    &\underset{\eqref{eq:lowshell:nopr:2},\eqref{eq:inductive:timescales}}{\lec} \Gamma_{m}^{-100}  \pi_q^m \lambda_{m}^N \MM{M, \Nindt, \tau_{m-1}^{-1}\Gamma_{m-1}^{i-5}, \Tau_{q}^{-1}\Ga_q^9} \, . \label{eq:Onpnp:estimate:2:0}
\end{align}

In order to prove \eqref{eq:stress:Linfty:upgraded}, we utilize a very similar argument to the one used to produce \eqref{eq:Onpnp:estimate:2:0}.  The only difference is that instead of appealing to \eqref{eq:lowshell:nopr:2}, we appeal to \eqref{eq:o.p.1}, \eqref{eq:t.p.1}, \eqref{eq:c.p.1}, and \eqref{eq:c.p.1.5}.  We omit further details.

Lastly, the proof of \eqref{eq:nlstress:upgraded} is very similar to \eqref{eq:Rcomm:bounds}, and so we omit further details. 
\end{proof}

\opsubsection{Total pressure increment and current from stress errors}\label{sec:new.pressure.stress}

We collect the pressure increments generated by new stress errors and new velocity increment potentials. Recall that Lemmas~\ref{lem:oscillation:pressure}, \ref{lem:transport:pressure}, and \ref{lem:corrector:pressure} defined pressure increments ($\sigma_{S_O^m}$, $\sigma_{S_{TN}}$, and $\sigma_{S_C^m}$, respectively) associated to various stress errors. Fixing $m$ such that $q+\half+1 \leq m \leq q+\bn$, we define
\begin{align}
    \si_{S^m} &:= \si_{S^{m}_O} + \si_{S^{m}_C} + \mathbf{1}_{\{m=\qbn\}} \sigma_{S_{TN}} \,.  \label{defn:si.m.stress}
\end{align}
Recalling that every pressure increment referenced above has a decomposition $\si_\bullet = \si_\bullet^+ - \si_\bullet^-$, we define $\si_{m,q+1}^+$ and $\si_{m,q+1}^-$ in the obvious way. 

Next, associated to each pressure increment $\sigma_\bullet$ listed above is a function of time $\bmu_{\sigma_\bullet}$ which satisfies $\bmu_{\sigma_\bullet}' = \langle \Dtq \sigma_\bullet \rangle$ (see Lemmas~\ref{lem:oscillation:pressure}, \ref{lem:transport:pressure}, \ref{lem:corrector:pressure}), and so we define\index{$\bmu_{\sigma_{S^m}}$}
\begin{align}
    \bmu_{\sigma_{S^m}} &:= \bmu_{\sigma_{S^{m}_O}} + \bmu_{\sigma_{S^{m}_C}} + \mathbf{1}_{\{m=\qbn\}} \bmu_{\sigma_{S_{TN}}} \, . \label{defn:bmu.m.stress}
\end{align}
Furthermore, recall that Lemmas~\ref{lem:oscillation:pressure:current}, \ref{lem:transport:pressure:current}, and \ref{lem:corrector:pressure:current} defined current errors associated to various stress error pressure increments. Then fixing $m,m'$ such that $q+\half+1 \leq m' \leq m \leq \qbn$, we define
\begin{subequations}\label{defn:phi.m.m'.stress}
\begin{align}
    \phi^{m',l}_{S^m} &:= \phi^{m',l}_{{S^{m}_O}} + \phi^{m',l}_{S^{m}_C} + \mathbf{1}_{\{m=\qbn\}} \phi_{S_{TN}}^{m',l}\\
    \phi^{m',*}_{S^m} &:= \phi^{m',*}_{{S^{m}_O}} + \phi^{m',*}_{S^{m}_C} +\mathbf{1}_{\{m=\qbn\}} \phi_{S_{TN}}^{m',*} \notag \\
   &\qquad \qquad + \mathbf{1}_{\{m'=m\}} \bigg{(} \phi^{*}_{{S^{m}_O}} + \phi^{*}_{S^{m}_C} +\mathbf{1}_{\{m=\qbn\}} \phi_{S_{TN}}^{*} \bigg{)} \, .
\end{align}
\end{subequations}
Now we set
\begin{equation}\label{def:phi:m:qplus.stress}
    \phi_{S^m} := \sum_{m'=q+\half{+1}}^m \phi^{m',l}_{S^m} + \phi^{m',*}_{S^m} \, ,
\end{equation}
so that the aforementioned lemmas give the equality 
\begin{align}\label{exp.Dtq.si.stress}
\div \phi_{S^m} = \Dtq \si_{S^m} - \bmu'_{S^m}
=\Dtq \si_{S^m}-\langle \Dtq \si_{S^m} \rangle \, .
\end{align}
By appealing to the lemmas mentioned above, we have that the $\si_{S^m}$'s satisfy the properties listed in the following lemma.

\begin{lemma}[\bf Collected properties of stress error terms and pressure increments]\label{lem:prop.si.pre.stress}
For each $q+\half+1 \leq m \leq q+\bn$, $\si_{S^m}$ satisfies the following properties.
\begin{enumerate}[(i)]
\item For any $0\leq k \leq \dpot$, we have that
    \begin{subequations}
        \begin{align}
        \left|\psi_{i,q}D^N D_{t,q}^M S^{m,l}_{q+1}\right|
    &\lec \left(\si_{S^m}^+ +  {\de_{q+3\bn}}\right)  (\la_{m}\Ga_m)^N  \label{est.S.m.pt.sim.stress}
    \MM{M, \Nindt, \Ga_q^{i+18} \tau_q^{-1}, \Tau_q^{-1}\Ga_q^9} 
    \end{align}
    \end{subequations}
where the bound holds for $N+M\leq 2\Nind$.
\item For $N, M \leq {\sfrac{\Nfin}{200}}$, we have that
\begin{subequations}
    \begin{align}
        \norm{\psi_{i,q}D^N D_{t,q}^M \si_{S^m}^+}_{\sfrac32}
    &\lec \Ga_m^{-9} \de_{m+\bn}  (\la_{m}\Ga_m)^N 
    \MM{M, \Nindt, \Ga_q^{i+18} \tau_q^{-1}, \Tau_q^{-1}\Ga_q^9} \label{est:si.m+.32.Dtq.stress}\\
    \norm{\psi_{i,q}D^N D_{t,q}^M \si_{S^m}^+}_{\infty}
    &\lec {\Ga_m^{\badshaq-9}}   (\la_{m}\Ga_m)^N 
    \MM{M, \Nindt, \Ga_q^{i+18} \tau_q^{-1}, \Tau_q^{-1}\Ga_q^9} \label{est:si.m+.infty.Dtq.stress}\\
    \left|\psi_{i,q}D^N D_{t,q}^M \si_{S^m}^+\right|
    &\lec \left(\si_{S^m}^+ +  {\de_{q+3\bn}}\right)  (\la_{m}\Ga_m)^N 
    \MM{M, \Nindt, \Ga_q^{i+18} \tau_q^{-1}, \Tau_q^{-1}\Ga_q^9} \label{est:si.m+.pt.stress}\\
    \left|\psi_{i,q}D^N D_{t,q}^M \si_{S^m}^-\right|
    &\lec {\Ga_{q+\half}^{-100}} \pi_q^{q+\half}  (\la_{q+\half}\Ga_{q+\half})^N 
    \MM{M, \Nindt, \Ga_q^{i+18} \tau_q^{-1}, \Tau_q^{-1}\Ga_q^9} \, .
\end{align}
\end{subequations}
\item\label{item:sat:one} $\si_{S^m}$ and $\si_{S^m}^+$ have the support properties
\begin{subequations}
\begin{align}
    B(\supp \hat w_{q'}, \la_{q'}^{-1}{\Ga_{q'+1}})
    \cap \si_{S^m} = \emptyset \qquad \forall q+1\leq q'\leq q+\half  \, , \label{supp.si.m.stress}\\
    B(\supp \hat w_{q'}, \la_{q'}^{-1}{\Ga_{q'+1}})
    \cap \si_{S^m}^+ = \emptyset \qquad \forall q+1\leq q'\leq m-1\, . \label{supp.si.m+.stress}
\end{align}
\end{subequations}
\item\label{upgrade:item:5} The function of time $\bmu_{\sigma_{S^m}}$ defined in \eqref{defn:bmu.m.stress} satisfies 
\begin{align}\label{est:bmu:stress}
      \left|\frac{d^{M+1}}{dt^{M+1}} \bmu_{\sigma_{S^m}} \right| 
      \lec (\max(1, T))^{-1}\delta_{q+3\bn} \MM{M,\Nindt,\tau_q^{-1},\Tau_{q+1}^{-1}} \, ,
    \end{align}
    for $0\leq M\leq 2\Nind$.
\end{enumerate}
\end{lemma}

\begin{lemma*}[\bf Total pressure current from stress errors]\label{lem:stress:pressure:current.stress}
For every $m\in\{q+\half+1,\dots,q+\bn\}$, the current error $\phi_{{S^{m}}}$ defined in \eqref{def:phi:m:qplus.stress} satisfies the following properties.
\begin{enumerate}[(i)]
    \item\label{i:pc:2:ER.stress} We have the decompositions and equalities
    \begin{subequations}
    \begin{align}\label{eq:desert:decomp:ER.stress}
        \phi_{{S^{m}}} &= \phi_{{S^{m}}}^* + \sum_{m'=q+\half+1}^{{m}} \phi_{{S^{m}}}^{m'} \, , \quad 
         \phi_{S^{m}}^{m'} = \phi_{S^{m}}^{m',l} + \phi_{S^{m}}^{m',*} \, \\
         \div \phi_{{S^{m}}}
          &= D_{t,q}\si_{{S^{m}}} - \langle D_{t,q}\si_{{S^{m}}}\rangle\, . 
    \end{align}
    
    \end{subequations}
    \item\label{i:pc:3:ER.stress} For $q+\half+1 \leq m' \leq m$ and $N,M\leq  2\Nind$,
    \begin{subequations}
    \begin{align}
        &\left|\psi_{i,q} D^N \Dtq^M \phi_{S^{m}}^{m',l} \right| < \Ga_{m'}^{-100} \left(\pi_q^{m'}\right)^{\sfrac 32} r_{m'}^{-1} (\la_{m'} \Ga_{m'}^2)^M \MM{M,\Nindt,\tau_q^{-1}\Ga_q^{i+17},\Tau_q^{-1}\Ga_q^9} \label{s:p:c:pt.stress} \\
        &\left\| D^N \Dtq^M \phi_{S^{m}}^{m',*} \right\|_\infty + \left\| D^N\Dtq^M \phi_{S^{m}}^{*} \right\|_\infty < \Tau_\qbn^{2\Nindt} \delta_{q+3\bn}^{\sfrac 32} (\la_{m}\Ga_{m}^2)^N \tau_q^{-M} \label{e:p:c:nonlocal.stress} \, .
    \end{align}
    \end{subequations}
    \item\label{i:pc:4:ER.stress} For all $q+\half+1\leq m' \leq m$ and all $q+1\leq q' \leq m'-1$, 
    \begin{align}
        B\left( \supp \hat w_{q'}, \sfrac 12 \lambda_{q'}^{-1} \Ga_{q'+1} \right) \cap \supp \left( \phi^{m',l}_{S^{m}} \right) = \emptyset \label{s:p:c:supp.stress} \, .
    \end{align}
\end{enumerate}
\end{lemma*}

\opsubsection{Transport/Nash current error}\label{op:tnce}

Recall the definitions of the stress error terms $\ov R_{q+1}$ and $S_{q+1}$ from \eqref{ER:new:error} and \eqref{eq:ovR:def}.  Since $\div \left( (R_q-\pi_q \Id) \hat u_q \right)$ appears in the relaxation \eqref{ineq:relaxed.LEI} of the local energy inequality, the new Reynolds stresses $\ov R_{q+1}$ and $S_{q+1}$ will create current error terms.  For this reason, we must estimate the Nash current error, which is given by
\begin{equation}\label{eq:nash:cur:def}
    \div \ov\phi_{N} +\bmu'_N := \na \hat u_q : (w_{q+1}\otimes w_{q+1} + R_q {-\pi_q^q\Id} - \ov R_{q+1} ) \, .
\end{equation}
The function of time $\bmu'_N$ is defined by
\begin{equation}\label{bmun:def}
     \bmu_N(t) := \int_0^t \left\langle \na \hat u_q : (w_{q+1}\otimes w_{q+1} + R_q -  R_{q+1} ) \right\rangle(s) \, ds
\end{equation}
and ensures that the error can be put in divergence form.  In addition, we must estimate a similar error term called the transport current error, which is given by
\begin{equation}\label{eq:transport:cur:def}
    \div \ov\phi_{T} +\bmu'_T:= (\pa_t + \hat u_q \cdot \na ) \left(\frac 12 |w_{q+1}|^2+ {\ka_q^q}  - \frac12\tr(S_{q+1})\right) \, .
\end{equation}
As before, we set 
\begin{equation}\label{bmut:def}
        \bmu_T(t) := \int_0^t \left\langle (\pa_t + \hat u_q \cdot \na ) \left(\frac 12 |w_{q+1}|^2+ \ka_\ell  - \frac12\tr(S_{q+1})\right) \right\rangle(s) \, ds
\end{equation}
to ensure that the error can be put in divergence form.  For a detailed derivation of how these error terms arise by adding $\hat w_\qbn$ to the relaxed local energy inequality, we refer to \cite[subsection~5.1]{GKN23}.

We now carefully decompose these error terms.  
Recall that from \eqref{eqn:osc.expand} and  \eqref{eq:cancellation:plus:pressure:nn}, we have
\begin{subequations}\label{eqn:w2}
\begin{align}
 \left(\wp_{q+1}\otimes \wp_{q+1}\right)^{\alpha,\bullet}
 -\pi_\ell\Id + R_\ell
&=
\sum_{\xi,i,j,k,\vecl} A_{(\xi),R}^{\alpha,\bullet}
(\mathbb{P}_{\neq 0}\chib_{\xi}^6 )(\Phi_{(i,k)})
\label{eqn:w2.exp2}\\
&\quad+
\sum_{\xi,i,j,k,\vecl} A_{(\xi),\ph}^{\alpha,\bullet}
\mathbb{P}_{\neq 0}\chib_{\xi}^4(\Phi_{(i,k)}) c_0 {c_1} r_q^\frac23
\label{eqn:w2.exp2.ph}\\
&\quad + {c_0}
\sum_{\xi,i,j,k,\vecl} A_{(\xi),\ph}^{\alpha,\bullet} 
 r_q^\frac23 \left( \chib_{\xi}^4\mathbb{P}_{\neq 0} 
\sum_I (\etab_\xi^I)^4\right) \circ \Phi_{(i,k)} 
\label{eqn:w2.exp.2.5}\\
&\quad+
\sum_{\xi,i,j,k,\vecl,\diamond} A_{(\xi),\diamond}^{\alpha,\bullet}
\left(\chib_{\xi}^{2\diamond}\sum_{I}(\etab_{\xi}^I)^{2\diamond}
\mathbb{P}_{\neq 0}(\varrho_{\xi,\diamond}^I)^2\right)(\Phi_{(i,k)})
\label{eqn:w2.exp3}
\end{align}
\end{subequations}
where $A_{(\xi),\diamond}^{\alpha,\bullet} := \xi^\theta \xi^\gamma \left(  a_{(\xi),\diamond}^2 
(\nabla\Phi_{(i,k)}^{-1})_{\theta}^\alpha 
(\nabla\Phi_{(i,k)}^{-1})_\gamma^\bullet\right)$.  To shorten notation, we define the operator 
\begin{align}\label{def:LTN}
    L_{TN} := (\pa_t + \hat u_q \cdot \na ) \frac{1}{2} \tr + (\na \hat u_q) : \, .
\end{align}
Using \eqref{eq:ovR:def}, we then write 
\begin{align}
     (\pa_t + &\hat u_q \cdot \na ) \left(\frac 12 |w_{q+1}|^2 + \ka_q^q  - \frac {\tr (S_{q+1})}2 \right) + (\na \hat u_q) : \left(w_{q+1}\otimes w_{q+1} + R_q - \pi_q^q \Id - \ov R_{q+1}\right) \notag \\
     &= L_{TN} \left( w_{q+1}^{(p)} \otimes w_{q+1}^{(p)} + {R_\ell - \pi_\ell \Id} \right)  
     \label{eq:ct:split:one}
     \\
     &\quad + L_{TN}  \left( w_{q+1}^{(p)} \otimes_s w_{q+1}^{(c)} \right) 
     \label{eq:ct:split:two}
     \\
     &\quad + L_{TN} \left( R_q^q - \pi_q^q \Id - {R_\ell + \pi_\ell \Id} - S_{q+1} + w_{q+1}^{(c)} \otimes w_{q+1}^{(c)} \right) \, . \label{eq:ct:split:three}
\end{align}
From \eqref{eqn:w2}, we have that \eqref{eq:ct:split:one} is actually equal to 
\begin{equation}
 \eqref{eq:ct:split:one}= 
L_{TN} {\bigl(  \eqref{eqn:w2.exp2} + \eqref{eqn:w2.exp2.ph} + \eqref{eqn:w2.exp.2.5} + \eqref{eqn:w2.exp3}  \bigr)}
    \, . \label{eq:ct:split:one:alt}
\end{equation}
Since $\Dtq$ can never land on the high-frequency object in these terms, we will estimate them directly using the inverse divergence. 
We will estimate \eqref{eq:ct:split:two} directly using the inverse divergence, and the fact that the high-frequency part of a product of principal and corrector parts has zero mean from Proposition~\ref{prop:pipeconstruction}, item~\ref{item:pipe:means} and Proposition~\ref{prop:pipe.flow.current}, item~\ref{item:pipe:means:current}. The last term, on the other hand, can be written as
\begin{align}
    \eqref{eq:ct:split:three} = - L_{TN} \left(S_O + S_{TN} + S_{C1} + S_{M2} \right) \label{eq:ct:split:four}
\end{align}
using \eqref{ER:new:error:moll1} and \eqref{eq:div:cor:expand}. We now split the analysis of these error terms into several lemmas.

\opsubsubsection{Transport/Nash current error from principal part of the velocity increment}

\begin{lemma*}[\bf Current error and pressure increment from \eqref{eq:ct:split:one}]\label{lem:ct:general:estimate}
There exists a vector field $\ov\phi_{TNW}$ and a function $\bmu_{TNW}$ of time such that
\begin{align}
    L_{TN} \left( w_{q+1}^{(p)} \otimes w_{q+1}^{(p)} + R_\ell - \pi_\ell \Id \right) &=L_{TN} \bigl(  \eqref{eqn:w2.exp2} + \eqref{eqn:w2.exp2.ph} + \eqref{eqn:w2.exp.2.5} + \eqref{eqn:w2.exp3}  \bigr) \notag\\
    &= \div \ov \phi_{TNW} + \bmu_{TNW}' \, , \notag  \\
    \ov \phi_{TNW} &= \sum_{m=q+1}^{q+\bn}  \ov \phi_{TNW}^m \, , \notag
\end{align}
where $\ov \phi_{TNW}^m = \ov \phi_{TNW}^{m,l} + \ov \phi_{TNW}^{m,*} $ for $m\in\{q+1,\dots,q+\bn\}$ satisfy the following.
\begin{enumerate}[(i)]
    \item The errors $\ov \phi^{q+1}_{TNW}$ and $\ov \phi^{q+\floor{\bn/2}}_{TNW}$ require no pressure increment. More precisely, we have that for {$N,M\leq \sfrac{\Nfin}{100}$,}
\begin{subequations}
\begin{align}
    \label{eq:ct:lowshell:nopr:1}
    \left|\psi_{i,q} D^N \Dtq^M \ov \phi^{q+1,l}_{TNW} \right| &< \Ga_{q+1}^{-100} \left(\pi_q^{q+1}\right)^{\sfrac32} r_{q+1}^{-1} \la_{q+1}^N \MM{M,\Nindt, \tau_q^{-1} \Ga_q^{i+15}, \Tau_q^{-1}\Ga_q^8} \, , \\
    \label{eq:ct:lowshell:nopr:2}
    \left|\psi_{i,q} D^N \Dtq^M \ov \phi^{q+\floor{\bn/2},l}_{TNW} \right| &< \Ga_{q+\half}^{-100} \left(\pi_q^{q+\half}\right)^{\sfrac32} r_{q+\half}^{-1} \la_{q+\floor{\bn/2}}^N \MM{M,\Nindt, \tau_q^{-1} \Ga_q^{i+15}, \Tau_q^{-1}\Ga_q^8} \, .
\end{align}
\end{subequations}
\item For $q+\half+1 \leq m \leq \bn$, there exists functions $\si_{\ov \phi^m_{TNW}} = \si_{\ov \phi^m_{TNW}}^+ - \si_{\ov \phi^m_{TNW}}^-$ such that
\begin{subequations}
\begin{align}
    \label{eq:ct.p.1}
    \left|\psi_{i,q} D^N \Dtq^M \ov \phi^{m,l}_{{TNW}}\right| &\les \left( (\si_{\ov \phi^{m}_{TNW}}^+)^{\sfrac 32} r_m^{-1} + {\de_{q+3\bn}^2}\right)  \left(\lambda_{m}\Gamma_q\right)^N \MM{M,\Nindt,\tau_q^{-1}\Gamma_{q}^{i+16},\Tau_q^{-1}{\Ga_q^9}}\\
    \label{eq:ct.p.2}
    \left|\psi_{i,q} D^N \Dtq^M \si_{\ov \phi^m_{TNW}}^+ \right| &\les \left(\si_{\ov\phi^m_{TNW}}^+ +\de_{q+3\bn}\right) \left(\lambda_{m}\Gamma_q\right)^N \MM{M,\Nindt,\tau_q^{-1}\Gamma_{q}^{i+{17}},\Tau_q^{-1}{\Ga_q^9}}\\
    \label{eq:ct.p.3}
    \norm{\psi_{i,q} D^N \Dtq^M \si_{\ov\phi^m_{TNW}}^+}_{\sfrac32} &\les \de_{m+\bn} \Gamma_{m}^{-9} \left(\lambda_{m}\Gamma_q\right)^N \MM{M,\Nindt,\tau_q^{-1}\Gamma_{q}^{i+{17}},\Tau_q^{-1}{\Ga_q^9}}\\
    \label{eq:ct.p.3.1}
    \norm{\psi_{i,q} D^N \Dtq^M \si_{\ov\phi^m_{TNW}}^+}_{\infty} &\les \Gamma_{q+1}^{\badshaq-9} \left(\lambda_{m}\Gamma_q\right)^N \MM{M,\Nindt,\tau_q^{-1}\Gamma_{q}^{i+{17}},\Tau_q^{-1}{\Ga_q^9}}\\
    \label{eq:ct.p.4}
    \left|\psi_{i,q} D^N \Dtq^M \si_{\ov \phi^m_{TNW}}^- \right| &\les \left(\frac{\la_q}{\la_{q+\floor{\bn/2}}}\right)^{\sfrac23} \pi_q^q  \left(\lambda_{q+\floor{\bn/2}}\Gamma_q\right)^N \MM{M,\Nindt,\tau_q^{-1}\Gamma_{q}^{i+{17}},\Tau_q^{-1}{\Ga_q^9}}
\end{align}
\end{subequations}
for all $N,M \leq \sfrac{\Nfin}{100}$. 
Furthermore, we have that for $q+1 \leq m' \leq m-1$ and $q+1\leq q'' \leq q+\half$,
\begin{align}\label{eq:ct.p.6}
    \supp \si_{\ov\phi^m_{TNW}}^- \cap B\left( \supp \hat w_{q''}, \la_{q''}^{-1}\Ga_{q''+1} \right) =
    \supp \si_{\ov \phi^m_{TNW}}^+ \cap B\left(\supp \hat w_{m'}, \la_{m'}^{-1}\Ga_{m'+1} \right) = \emptyset \, .
\end{align}
\item When $m=q+2,\dots,q+\bn$ and $q+1\leq q' \leq m-1$, the local parts satisfy
\begin{align}
B\left( \supp \hat w_{q'}, \lambda_{q'}^{-1} \Gamma_{q'+1} \right) \cap \supp \ov \phi_{TNW}^{m,l} = \emptyset \, . \label{ct:support:first}
\end{align}
\item For $m=q+1, \dots, q+\bn$ and $N,M\leq 2\Nind$, the non-local parts ${\ov \phi}^{m,*}_{O}$ satisfy
\begin{align}
\left\| D^N D_{t, q}^M {\ov \phi}^{m,*}_{{TNW}} \right\|_{L^\infty}
    &\leq { \Tau_{q+\bn}^{2\Nindt}}\delta_{q+3\bn}^{\sfrac32} \la_{m}^{N}\tau_{q}^{-M} \, .
    \label{eq:ctnl:estimate:1}
\end{align}

\item For $M \leq 2\Nind$, the time function $\bmu_{TNW}$ satisfies
\begin{align}\label{sat:morn:10:57}
\bmu_{TNW}(t) = \int_0^t \langle\textnormal{(\ref{eq:ct:split:one})} (s) \rangle\, ds \, , \quad \left| \frac{d^{M+1}}{dt^{M+1}} \bmu_{TNW} \right| \leq \left( \max(1,T) \right)^{-1} \delta_{q+3\bn}^2 \tau_q^{-M} \, .
\end{align}
\end{enumerate}
\end{lemma*}

\begin{proof}
The analysis of this error is similar to that of the oscillation stress error dealt with in subsection~\ref{ss:ssO}, Lemmas~\ref{lem:oscillation:general:estimate}--\ref{lem:oscillation:pressure:current}. We will invert the divergence on this error term using Proposition~\ref{prop:intermittent:inverse:div} and apply Proposition~\ref{lem.pr.invdiv2.c} to construct the pressure increment. Let us define
\begin{subequations}
\begin{align}
    \ov \phi^{q+1}_{{TNW}} &:= (\divH + \divR) \left[\sum_{\xi,i,j,k,\vecl} L_{TN} \left ( A_{(\xi),R}^{\alpha,\bullet}\right)
(\mathbb{P}_{\neq 0}\chib_{\xi}^6 )(\Phi_{(i,k)}) \right] \notag\\
&\qquad + (\divH + \divR) \left[\sum_{\xi,i,j,k,\vecl} L_{TN} \left ( A_{(\xi),\ph}^{\alpha,\bullet}\right)
\mathbb{P}_{\neq 0}\chib_{\xi}^4(\Phi_{(i,k)}) c_0 {c_1} r_q^\frac23 \right] \label{ct.err.low} \\
\label{ct.err.med1-}
\ov\phi^{q+\floor{n/2}}_{{TNW}} &:= (\divH + \divR) \left[ \sum_{\xi,i,j,k,\vecl} L_{TN} \left ( A_{(\xi),\ph}^{\alpha,\bullet}\right) 
 {c_0}r_q^\frac23 \left(\chib_{\xi}^4\mathbb{P}_{\neq 0}
\sum_I (\etab_\xi^I)^4 \right)\circ \Phiik\right]\\
\label{ct.err.med0}
\ov\phi^{q+\floor{n/2}+1}_{{TNW}} &:= (\divH + \divR) \left[\sum_{\xi,i,j,k,\vecl,I,\diamond} L_{TN} \left (A_{(\xi),\diamond}^{\alpha,\bullet}\right)
    \left(\chib_{\xi}^{2\diamond}(\etab_{\xi}^I)^{2\diamond}
    \tP_{q+\bn+1}^\xi \mathbb{P}_{\neq 0}  (\varrho_{\xi,\diamond}^I)^2\right)(\Phi_{(i,k)}) \right]\\
    \label{ct.err.med}
    \ov\phi^{m}_{{TNW}} &:= (\divH + \divR) \left[\sum_{\xi,i,j,k,\vecl,I,\diamond} L_{TN} \left (A_{(\xi),\diamond}^{\alpha,\bullet}\right)
    \left(\chib_{\xi}^{2\diamond}(\etab_{\xi}^I)^{2\diamond}
    \tP_{(m-1,m]}^\xi (\varrho_{\xi,\diamond}^I)^2\right)(\Phi_{(i,k)}) \right]\\
    \label{ct.err.high}
    \ov\phi^{q+\bn}_{{TNW}} &:= \sum_{m=q+\bn}^{q+\bn+1}(\divH + \divR) \left[\sum_{\xi,i,j,k,\vecl,I,\diamond} L_{TN} \left (A_{(\xi),\diamond}^{\alpha,\bullet}\right)
    \left(\chib_{\xi}^{2\diamond}(\etab_{\xi}^I)^{2\diamond}
    \tP_{(m-1,m]}^\xi(\varrho_{\xi,\diamond}^I)^2\right)(\Phi_{(i,k)}) \right]\\
    &\qquad + (\divH + \divR) \left[\sum_{\xi,i,j,k,\vecl,I,\diamond} L_{TN} \left (A_{(\xi),\diamond}^{\alpha,\bullet}\right)
    \left(\chib_{\xi}^{2\diamond}(\etab_{\xi}^I)^{2\diamond}\left(\Id-
    \tP_{q+\bn+1}^\xi\right) (\varrho_{\xi,\diamond}^I)^2\right)(\Phi_{(i,k)}) \right] \label{ct.err.high2}
\end{align}
\end{subequations}
for $m=q+\half+2, \cdots, q+\bn -1$. We decompose $\ov\phi^{m}_{{TNW}}$ into the nonlocal part $\ov\phi^{m,*}_{{TNW}}$ which involves the operator $\divR$ or $\Id-\tP_{q+\bn+1}^\xi$ and the local part $\ov\phi^{m, l}_{{TNW}}$ containing the remaining terms. For the undefined $ \ov\phi^{m}_{{TNW}}$ corresponding to $m=q+2,\cdots, q+\half-1$, we set them as identically zero. 

The construction of the pressure increment and the desired estimates will follow from applying Propositions~\ref{prop:intermittent:inverse:div} and \ref{lem.pr.invdiv2.c}. While many of the parameter choices will vary depending on the case, we fix the following choices throughout the proof:
\begin{subequations}\label{eq:ct:general:choices}
\begin{align}
    &v= \hat u_q \, , \quad D_t = D_{t,q} \, , \quad {N_* = \sfrac{\Nfin}{4} \, , \quad M_* = \sfrac{\Nfin}{5}} \ , \\
    &\la' = \Lambda_q \, , \quad M_t = \Nindt \, , \quad \nu' = \Tau_q^{-1}\Ga_q^8 \, , \quad \Ndec \textnormal{ as in \eqref{i:par:9}} \, .
\end{align}
\end{subequations}

\noindent\texttt{Case 1: Estimates for \eqref{ct.err.low}.}
Fix values of $i,j,k,\xi,\vecl$ and consider the term which includes $L_{TN} A_{(\xi),R}$.
We apply Proposition \ref{prop:intermittent:inverse:div} with the low-frequency choices
\begin{align*}
    &G = L_{TN} A_{(\xi),R}\, , \quad \const_{G,3/2} =\left| \supp(\eta_{i,j,k,\xi,\vecl,R}^2) \right| \tau_q^{-1}\Ga_q^{i+13} \delta_{q+\bn} \Gamma_{q}^{2j+8} \, , \quad \const_{G,\infty}=\Gamma_{q}^{\badshaq+14}\tau_q^{-1}\Ga_q^{\imax+13} \, , \\
    &\pi = {\Ga_q^{50} \tau_q^{-1}\Ga_q^i \psi_{i,q} \pi_\ell} \, , \quad \la =\la_{q+1}\Gamma_q^{-5} \, , \quad \nu = \tau_q^{-1} \Gamma_q^{i+14} \, , \quad \Phi= \Phi_{(i,k)} \, ,
\end{align*}
and the choices from \eqref{eq:ct:general:choices}. By Corollary \ref{cor:deformation}, $\Phi_{(i,k)}$ satisfies \eqref{eq:DDpsi2} and \eqref{eq:DDpsi}, and by \eqref{eq:nasty:D:vq:old} at level $q$ and \eqref{ineq:tau:q}, we have that \eqref{eq:DDv} is satisfied. To check \eqref{eq:inverse:div:DN:G}, we observe that $L_{TN}$ involves a material derivative and a multiplication by $\na \hat u_q$. Therefore, by \eqref{eq:nasty:D:vq:old}, $G$ satisfies \eqref{eq:inverse:div:DN:G} for $p=\sfrac 32$ from \eqref{e:a_master_est_p_R} and for $p=\infty$ from the same inequality and \eqref{ineq:jmax:use}. Also, \eqref{eq:inv:div:extra:pointwise} is satisfied by \eqref{e:a_master_est_p_pointwise}. To check the high-frequency assumptions, we set (exactly as in the analogous case for the oscillation stress error - see Lemmas~\ref{lem:oscillation:general:estimate}--\ref{lem:oscillation:pressure:current})
\begin{subequations}
\begin{align*}
    &\varrho = \left( \mathbb{P}_{\neq 0} \ov \rhob_\xi^6 \right) \, , \quad \dpot \textnormal{ as in \eqref{i:par:10}} \, , \quad \vartheta=\delta_{i_1 i_2} \delta_{i_3 i_4} \dots \delta_{i_{\dpot-1} i_\dpot} \Delta^{-\sfrac \dpot 2} \varrho \, , \\
    &\mu = \Upsilon=\Upsilon'= \lambda_{q+1}\Gamma_q^{-4} \, , \quad \ov\Lambda = \lambda_{q+1}\Gamma_q^{-1} \, , \quad {\const_{*,1} = \Ga_q^6 \la_{q+1}^{\alpha}} \, .
\end{align*}
\end{subequations}
Since the choice of parameters is exactly the same as in the oscillation stress error, we see that the other high frequency assumptions are satisfied. In order to check the nonlocal assumptions, we set
\begin{align}
    M_\circ = N_\circ = 2\Nind \, , \quad K_\circ \textnormal{ as in \eqref{i:par:9.5}} \, , \quad \const_v = \Lambda_q^{\sfrac12} \, .
\end{align}
Then from \eqref{ineq:dpot:1} and Remark~\ref{rem:lossy:choices}, we have that \eqref{eq:inv:div:wut}--\eqref{eq:riots:4} are satisfied.

{We can therefore apply Remark~\ref{rem:pointwise:inverse:div}. Note that \eqref{eq:inv:div:extra:pointwise} follows from the definition of $L_{TN} A_{\pxi,R}$ in \eqref{eq:rhoxidiamond:def} and \eqref{e:a_master_est_p_R_pointwise}.  Then, abbreviating $G \varrho \circ \Phi$ as $t_{i,j,k,\xi,\vecl,R}$, from \eqref{eq:divH:formula}, \eqref{eq:inverse:div:sub:1}, and \eqref{eq:inv:div:extra:conc}, we have that for all $N \leq \frac{\Nfin}4 - \dpot$ and $M \leq \frac{\Nfin}5$
\begin{align*}
    \left|D^N \Dtq^M \mathcal{H}t_{i,j,k,\xi,\vecl,R} \right| &\lec \tau_q^{-1} \Ga_q^i \psi_{i,q} \pi_\ell \Gamma_q^{60} \lambda_{q+1}^{-1} \la_{q+1}^{N+\alpha} \MM{M,\Nindt,\tau_q^{-1} \Ga_q^{i+14}, \Tau_q^{-1}\Ga_q^8} \, .
\end{align*}}
Notice that from \eqref{item:div:local:i}, we have
\begin{align}\label{supp.osc.current.piece}
    \supp(\div \mathcal{H} t_{i,j,k,\xi,\vecl,R})
    \subseteq \supp t_{i,j,k,\xi,\vecl,R}
    \subseteq \supp \eta_{i,j,k,\xi,\vecl,R} \, .
\end{align}

As for the terms which include $A^{\alpha,\bullet}_{\pxi,\varphi}$ from \eqref{ct.err.low}, we note that from Lemma~\ref{lem:a_master_est_p} $a_{\pxi,\varphi}^2$ differs in size relative to $a_{\pxi,R}^2$ by a factor of $r_q^{-\sfrac 23}$, which is exactly balanced out by the factor of $r_q^{\sfrac 23}$ in \eqref{ct.err.low}.  We therefore may argue exactly as above (in fact the estimates are slightly better since $\ov \rhob_\xi^4 < \ov \rhob_\xi^6$), and we omit further details. In this case, we use the abbreviation $t_{i,j,k,\xi,\vecl,\ph}$ instead of $t_{i,j,k,\xi,\vecl,R}$, which will satisfy an analogous support property to \eqref{supp.osc.current.piece}. 
 
We now set 
\begin{align}\notag 
    \ov\phi^{q+1,l}_{{TNW}} = \sum_{i,j,k,\xi,\vecl,\diamond} \mathcal{H} t_{i,j,k,\xi,\vecl,\diamond} \, .
\end{align}
Using \eqref{supp.osc.current.piece} and applying the aggregation Corollary~\ref{lem:agg.Dtq} with $H = \mathcal{H}t_{i,j,k,\xi,\vecl,\diamond}$ and 
$$ \varpi = \pi_\ell \Ga^{60} \la_{q+1}^{-1+\alpha}, \quad \la=\Lambda=\la_{q+1},\quad \tau = \tau_q \Ga_q^{-{14}}, \quad \Tau = \Tau_q \Ga_q^{-8} $$
to get an estimate from \eqref{eq:aggDtq:conc:1.0},
\begin{align*}
    \left| \psi_{i,q} D^N D_{t, q}^M \ov\phi^{q+1,l}_{{TNW}} \right| 
    &\lesssim r_q^{-1} \la_q (\pi_q^q)^{\sfrac12} \pi_\ell \Gamma_q^{61} \lambda_{q+1}^{-1} \la_{q+1}^{N+\alpha} \MM{M,\Nindt,\tau_q^{-1} \Ga_q^{i+{15}}, \Tau_q^{-1}\Ga_q^8} \, .
\end{align*}
for $N,M$ in the same range as above. Then, \eqref{eq:ct:lowshell:nopr:1} follows from this term using \eqref{ind:pi:lower}, \eqref{eq:ind.pr.anticipated.2} and \eqref{condi.Nfin0}. 

For the non-local term, from \eqref{eq:inverse:div:error:stress:bound}, and Remark~\ref{rem:lossy:choices}, we have that for $N,M\leq 2\Nind$,
\begin{align}
    \left\| D^N \Dtq^M \sum_{i,j,k,\xi,\vecl} \divR t_{i,j,k,\xi,\vecl,R} \right\|_\infty \leq \delta_{q+3\bn}^{\sfrac 32} \Tau_q^{2\Nindt} {\lambda_{q+1}^N} \tau_q^{-M} \, , \notag
\end{align}
matching the desired estimate in \eqref{eq:ctnl:estimate:1}. The estimate in \eqref{sat:morn:10:57} follows using Remark~\ref{rem:est.mean} and a large choice of $a_*$, and we omit further details. The version of these estimates in the later cases will again be similar, and so we do not address them again.
\smallskip

\noindent\texttt{Case 2: Estimates for  \eqref{ct.err.med1-}.} As before, we fix $i,j,k,\xi,\vecl$. We apply Proposition \ref{prop:intermittent:inverse:div} with Remark~\ref{rem:pointwise:inverse:div} with the low-frequency choices
\begin{subequations}
\begin{align}
     &G = L_{TN} A_{\pxi,\varphi} c_0 r_q^{\sfrac23} \ov\chib_{\xi}^{4} (\Phi_{(i,k)}) \, , \quad \const_{G,3/2} = \left| \supp \eta_{i,j,k,\xi,\vecl,\vp}^2 \right| \tau_q^{-1} \Ga_q^i \delta_{q+\bn}\Gamma_{q}^{20} \, , \quad \const_{G,\infty}= \Gamma_q^{\badshaq+20} \tau_q^{-1} \Ga_q^{\imax} \, , \\
     &\pi = \Ga_q^{50} \tau_q^{-1}\Ga_q^i \psi_{i,q} \pi_\ell \, , \quad \lambda = \lambda_{q+1} \Gamma_q^{-1} \, , \quad \nu = \tau_q^{-1}\Gamma_q^{i+13} \, , \quad \Phi = \Phiik \, ,
\end{align}
\end{subequations}
as well as the choices from \eqref{eq:ct:general:choices}. As in the previous substep, \eqref{eq:DDpsi2}, \eqref{eq:DDpsi}, and \eqref{eq:DDv} are satisfied. The estimates in \eqref{eq:inverse:div:DN:G} hold due to Proposition~\ref{prop:bundling} and the estimates for $L_{TN} A_{\pxi,\varphi}$ from \texttt{Case 1}.

To check the high-frequency assumptions, we set the parameters and functions exactly as in \texttt{Case 2} in the proof of Lemma \ref{lem:oscillation:general:estimate}. Since we work with $p=1$ instead of $p=\frac32$, the only difference is that $\const_{*,1} := \const_{*,\infty} = \la_{q+\half}^\al$ instead of $\const_{*,\sfrac32}$. Then, as before, high-frequency assumptions in \eqref{item:inverse:i}--\eqref{item:inverse:iv} can be verified. The nonlocal assumptions are identical to those of \texttt{Case 1}, and are satisfied trivially. The non-local parameters are set to be the same as in the previous case.

We therefore may appeal to the local conclusions \eqref{item:div:local:0}--\eqref{item:div:nonlocal} and \eqref{eq:inverse:div:error:stress}--\eqref{eq:inverse:div:error:stress:bound}, from which we have the following.  First, abbreviating $G \varrho \circ \Phi$ as $t_{i,j,k,\xi,\vecl,\varphi}$, we have from \eqref{eq:inverse:div} and \eqref{eq:inverse:div:stress:1} that for $N\leq \frac{\Nfin}{4}-\dpot$ and $M\leq \frac{\Nfin}{5}$, 
\begin{align*}
    \left| D^N D_{t, q}^M \mathcal{H} t_{i,j,k,\xi,\vecl,\varphi} \right|
    &\lesssim \tau_q^{-1}\Ga_q^i \psi_{i,q} \pi_\ell \Ga_q^{50} \la_{q+\half}^{-1} \la_{q+\half}^{N+\alpha} \MM{M, \Nindt, \tau_{q}^{-1}\Gamma_{q}^{i+14}, \Tau_{q}^{-1}\Ga_q^8} \, ,
\end{align*}
 Notice that from \eqref{item:div:local:i}, the support of $\div \mathcal{H}t_{i,j,k,\xi,\vecl,\vp}$ is contained in  $\supp t_{i,j,k,\xi,\vecl,\vp} \subset \supp \left(\eta_{i,j,k,\xi,\vecl,\vp}  \right) $.  Thus as before we may apply the aggregation Corollary~\ref{lem:agg.Dtq} with $H = \mathcal{H}t_{i,j,k,\xi,\vecl,R}$ and 
$$ \varpi = \pi_\ell \Ga^{50} \la_{q+\half}^{-1}, \quad \la=\Lambda=\la_{q+\half},\quad \tau = \tau_q \Ga_q^{-14}, \quad \Tau = \Tau_q \Ga_q^{-8} $$
to estimate
\begin{align}\notag 
    \ov\phi^{q+\half,l}_{{TNW}} = \sum_{i,j,k,\xi,\vecl} \mathcal{H} t_{i,j,k,\xi,\vecl,\ph} \, .
\end{align}
From \eqref{eq:aggDtq:conc:1.0}, we thus have that for $N,M$ in the same range as above,
\begin{align*}
    \left| \psi_{i,q} D^N D_{t, q}^M \ov\phi^{q+\half,l}_{{TNW}} \right| 
    &\lesssim r_q^{-1} \la_q (\pi_q^q)^{\sfrac12} \pi_\ell \Gamma_q^{50} \lambda_{q+\half}^{-1} \la_{q+\half}^{N+\alpha} \MM{M,\Nindt,\tau_q \Ga_q^{i+15}, \Tau_q^{-1}\Ga_q^8} \, .
\end{align*}
and so we can conclude \eqref{eq:ct:lowshell:nopr:2} as before.  we must verify \eqref{ct:support:first} for $\ov\phi^{q+\half,l}_{{TNW}}$. This however follows from \eqref{item:div:local:ii}, which asserts that the support of $\ov\phi_{TNW}^{q+\half,l}$ is contained in $\cup_{\pxi} \supp(a_{\pxi,\vp} \rhob_{\pxi}^{\vp} \circ \Phiik)$, and \eqref{item:dodging:more:oldies} of Lemma~\ref{lem:dodging}. The non-local conclusions also follow in much the same way as in \texttt{Case 1}, and we omit further details.
\smallskip

\noindent\texttt{Case 3: Estimates of the local portions of  \eqref{ct.err.med0}, \eqref{ct.err.med}, and \eqref{ct.err.high}.} Fix $\xi$, $i$, $j$, $k$, $\vecl$, $I$, and $\diamond$. In order to check the low-frequency, preliminary assumptions in Part 1 of Proposition~\ref{lem.pr.invdiv2.c}, we set
\begin{align}
    &p=1,\infty \, , \quad G_{R} =  L_{TN} \left (A_{(\xi),\diamond}^{\alpha,\bullet}\right)
    \left(\chib_{\xi}^{2\diamond}(\etab_{\xi}^I)^{2\diamond}
    \right)(\Phi_{(i,k)}) \, , \quad G_{\vp} = L_{TN} \left (A_{(\xi),\diamond}^{\alpha,\bullet}\right)
    \left(\chib_{\xi}^{2\diamond}(\etab_{\xi}^I)^{2\diamond}
    \right)(\Phi_{(i,k)}) r_q^{\sfrac23}\, , \notag\\
    &\const_{G_\diamond,1} = \delta_\qbn \tau_q^{-1} \Ga_q^{i+2j+20} \left| \supp \left( \eta_{i,j,k,\xi,\vecl,\diamond} \zetab_\xi^{I,\diamond} \right) \right| +\la_{q+\bn}^{-10} \, , \quad \const_{G_\diamond,\infty} = \delta_\qbn  \tau_q^{-1} \Ga_q^{i+2j+20} \, ,\notag\\
    &\la = \la_{q+\half}\, , \quad \nu=\tau_q^{-1}\Ga_q^{i+14}\, ,  \quad  \Phi=\Phiik \, , \quad \pi=\Ga_q^{50} \pi_\ell \la_q^{\sfrac 23} \, , \quad r_G=r_q \, .  \label{eq:ctn:desert:choices}
\end{align}
Then we have that \eqref{eq:inv:div:NM} is satisfied by definition, \eqref{eq:inverse:div:DN:G} is satisfied by \eqref{e:a_master_est_p_phi:zeta}, \eqref{e:a_master_est_p_R:zeta},  Corollary~\ref{cor:deformation}, \eqref{e:fat:pipe:estimates:1}, and Definition~\ref{def:etab}, \eqref{eq:DDpsi2}--\eqref{eq:DDv} hold from Corollary~\ref{cor:deformation} and \eqref{eq:nasty:D:vq:old} at level $q$, and \eqref{est.G.c.pt} holds from \eqref{e:a_master_est_p_pointwise}, Remark~\ref{rem.summing.psi}, and \eqref{ind:pi:lower}.

In order to check the high-frequency, preliminary assumptions in Part 1 of Proposition~\ref{lem.pr.invdiv2.c}, we choose parameters and functions exactly same as in \texttt{Case 3} and \texttt{Case 4} of Lemma \ref{lem:oscillation:general:estimate}. The only difference is that we use $C_{*,1}$ instead of $C_{*,\sfrac32}$. Indeed, we choose $C_{*,1} =\la_{q+\half+1}^\al$ in both cases \texttt{Case 3a} and \texttt{Case 3b}. Then, it is enough to check \eqref{eq:DN:Mikado:density}, which holds true due to Propositions~\ref{prop:pipeconstruction} and \ref{prop:pipe.flow.current} and estimate \eqref{eq:lowest:shell:inverse} from Lemma~\ref{lem:special:cases} or \ref{lem:LP.est} with $q=1$. In order to check the additional assumptions in Part 2 of Proposition~\ref{lem.pr.invdiv2.c}, we again choose the same parameters and functions as in as in \texttt{Case 3} and \texttt{Case 4} of Lemma \ref{lem:oscillation:pressure}, and set the extra parameters as $\de_{\phi, p}$ and $r_\phi$ are 
\begin{align*}
    \delta_{\phi,p}^{\sfrac 32} = \const_{G_\diamond,p}\const_{*,p}\Upsilon'\Upsilon^{-2}r_{\min(m,\qbn)}\, , \quad r_\phi = r_{\min(m,\qbn)}\, . 
\end{align*}
Compared to Proposition \ref{lem.pr.invdiv2}, we only need to check \eqref{eq:sample:prop:Ncut:2:c}, \eqref{eq:sample:riots:4:c}, and \eqref{eq:sample:riot:4:4:c}, which can be verified by \eqref{condi.Ncut0.2}, \eqref{ineq:dpot:1}, and \eqref{ineq:Nstarz:1}.

Using the abbreviation $t_{i,j,k,\xi,\vecl,I,\diamond}^{m}$ for $G\varrho\circ \Phi$ at the level of $q+\half+2\leq m \leq q+\bn+1$, as a consequence of \eqref{d:press:stress:sample:c}--\eqref{est.S.prbypr.pt:c}, \eqref{condi.Nfin0}, \eqref{condi.Nindt}, \eqref{ind:pi:lower}, and \eqref{ineq:r's:eat:Gammas}, there exists a pressure increment $\sigma_{ \divH t_{i,j,k,\xi,\vecl,I,\diamond}^{m}}^+$ such that for $N,M\leq \sfrac{\Nfin}{7}$, 
\begin{align}
    \left| D^N \Dtq^M \divH t_{i,j,k,\xi,\vecl,I,\diamond}^{m}  \right| &\lesssim \left( \left(\sigma_{ \divH t_{i,j,k,\xi,\vecl,I,\diamond}^{m}}^+\right)^\frac32 r_{\min(m, q+\bn)}^{-1} + \delta_{q+3\bn}^2 \right) \notag\\
    &\qquad \qquad \times (\min(\lambda_m,\lambda_{q+\bn})\Ga_q)^N \MM{M,\Nindt,\tau_q^{-1}\Ga_q^{i+{15}},\Tau_q^{-1}\Ga_q^{9}} \, .  \label{eq:desert:cowboy:1:redux:T}\\
     \left| D^N \Dtq^M \sigma_{ \divH t_{i,j,k,\xi,\vecl,I,\diamond}^{m}}^+ \right| &\lesssim \left( \sigma_{ \divH t_{i,j,k,\xi,\vecl,I,\diamond}^{m}}^+ + \delta_{q+3\bn}^2 \right) (\min(\lambda_{m},\la_\qbn)\Ga_q)^N \MM{M,\Nindt,\tau_q^{-1}\Ga_q^{i+{16}},\Tau_q^{-1}\Ga_q^{9}} \, . \label{eq:desert:cowboy:4:redux:TN}\\
     \left| D^N \Dtq^M \sigma^{-}_{\divH t_{i,j,k,\xi,\vecl,I,\diamond}^{m}} \right| 
     &\les \left(\frac{r_{\min(m,\qbn)}}{r_q} \right)^{\sfrac 23} \Ga_q^{28} \pi_\ell \La_q^{\sfrac 23} \left(\lambda_{m-1}^2\lambda_m^{-1}\right)^{-\sfrac 23} \notag \\
     &\qquad \qquad \times (\la_{q+\half}\Ga_q)^N \MM{M,\Nindt,\tau_q^{-1}\Ga_q^{i+{16}},\Tau_q^{-1}\Ga_q^9} \notag\\
    &\leq \Ga_q^{-10}\left( \frac{\la_q}{\la_{q+\half}} \right)^{\sfrac 23} \pi_q^q (\la_{q+\half}\Ga_q)^N \MM{M,\Nindt,\tau_q^{-1}\Ga_q^{i+{16}},\Tau_q^{-1}\Ga_q^9} \, . \notag 
\end{align}
From \eqref{eq:inverse:div:linear}, \eqref{est.S.pr.p.support:1:c}, and \eqref{eq:LP:div:support}, we have that 
\begin{align}
    \supp \left( \sigma_{ \divH t_{i,j,k,\xi,\vecl,I,\diamond}^{m}}^+ \right) &\subseteq \supp \left( \divH t_{i,j,k,\xi,\vecl,I,\diamond}^{m} \right) \notag\\
    &\subseteq \supp \left( a_{\pxi,\diamond} \left(\rhob_\pxi^\diamond \zetab_\xi^I \right)\circ\Phiik \right) \cap B\left( \supp \varrho_{\pxi,\diamond}^I, \lambda_{m-1}^{-1} \right){\circ\Phiik} \, \label{eq:desert:support:T}, \\
    \supp\left( \sigma_{ \divH t_{i,j,k,\xi,\vecl,I,\diamond}^{m}}^- \right) 
    &\subseteq \supp \left( a_{\pxi,\diamond} \left(\rhob_\pxi^\diamond \zetab_\xi^I \right)\circ\Phiik \right) \label{eq:desert:support:T2}
\end{align}
Then, we can obtain the desired estimates for 
\begin{align*}
    \ov \phi_{TNW}^{m.l} = \sum_{i,j,k,\xi, \vecl,I, \diamond}\divH t_{i,j,k,\xi,\vecl,I,\diamond}^{m}, \qquad 
    \ov \phi_{TNW}^{q+\bn.l} = \sum_{m=q+\bn}^{q+\bn+1} \sum_{i,j,k,\xi, \vecl,I, \diamond}\divH t_{i,j,k,\xi,\vecl,I,\diamond}^{m}, \,,\\
    \sigma_{\ov\phi_{TNW}^{m}}^{\pm} = \sum_{i,j,k,\xi,\vecl,I,\diamond} \sigma_{ \divH t_{i,j,k,\xi,\vecl,I,\diamond}^{m}}^{\pm} \, \qquad
    \sigma_{\ov\phi_{TNW}^{q+\bn}}^{\pm} = \sum_{m=q+\bn}^{q+\bn+1} \sum_{i,j,k,\xi,\vecl,I,\diamond} \sigma_{ \divH t_{i,j,k,\xi,\vecl,I,\diamond}^{m}}^{\pm} 
\end{align*}
for $q+\half+1 \leq m <q+\bn$ by applying Corollary~\ref{lem:agg.pt} with $p=1$ and
\begin{align*}
    H =\divH t_{i,j,k,\xi,\vecl,I,\diamond}^{m} \, , \qquad \varpi = \divH t_{i,j,k,\xi,\vecl,I,\diamond}^{m}\mathbf{1}_{\supp a_{\pxi,\diamond}(\rhob_\pxi^\diamond\zetab_\xi^I)\circ \Phiik}, \qquad&\text{for }\eqref{eq:ct.p.1}\\
    H= \sigma^+_{\divH t_{i,j,k,\xi,\vecl,I,\diamond}^{m}} \, , \qquad \varpi= \left[ H + \delta_{q+3\bn}^2 \right] \mathbf{1}_{\supp a_{\pxi,\diamond}(\rhob_\pxi^\diamond\zetab_\xi^I)\circ \Phiik},  \qquad &\text{for }\eqref{eq:ct.p.2} \\
    H=\sigma^{-}_{\divH t_{i,j,k,\xi,\vecl,I,\diamond}^{m}}\, , \qquad \varpi=\left( \frac{\la_q}{\la_{q+\half}} \right)^{\sfrac 23} \pi_\ell \mathbf{1}_{\supp a_{\pxi,\diamond}(\rhob_\pxi^\diamond \zetab_\xi^I)\circ\Phiik}, \qquad&\text{for }\eqref{eq:ct.p.4}\, . 
\end{align*}
Also, \eqref{eq:oooooldies}--\eqref{eq:dodging:oldies},  \eqref{eq:desert:support:T}, and \eqref{eq:desert:support:T2} give that \eqref{ct:support:first} and \eqref{eq:ct.p.6} are satisfied for $q+\half+1 \leq m'\leq q+\bn$.

Next, from \eqref{est.S.pr.p:c}, we have that
\begin{align}
    \left\| \sigma^{\pm}_{\divH t_{i,j,k,\xi,\vecl,I,\diamond}^{m}} \right\|_{\sfrac 32} &\les \delta_\qbn^{\sfrac23} \tau_q^{-\sfrac 23} \Ga_q^{\sfrac23 (i+2j+24)} \left| \supp \left( \eta_{i,j,k,\xi,\vecl,\diamond} \zetab_\xi^{I,\diamond} \right) \right|^{\sfrac 23} \left(\lambda_{m-1}^2\lambda_m^{-1}\right)^{-\sfrac 23} r_{\min(m,\qbn)}^{\sfrac 23} \notag \,,\\
    \left\| \sigma^{\pm}_{\divH t_{i,j,k,\xi,\vecl,I,\diamond}^{m}} \right\|_{\infty} &\les \delta_\qbn^{\sfrac23} \tau_q^{-\sfrac 23} \Ga_q^{\sfrac23 (i+2j+24)} \left( \frac{\min(\la_m,\la_\qbn)}{\la_{q+\half}\Ga_q} \right)^{\sfrac 43} \left(\lambda_{m-1}^2\lambda_m^{-1}\right)^{-\sfrac 23} r_{\min(m,\qbn)}^{\sfrac 23}  \notag \\
    &\les  \Ga_{q-\bn}^{\frac23\left(\frac{\badshaq}2 +18\right)}r_{q-\bn}^{-\sfrac23}r_{\min(m,\qbn)}^{\sfrac 23}
    \Ga_q^{\frac23(40+\badshaq)} \left( \frac{\min(\la_{m},\la_\qbn)}{\la_{q+\half}\Ga_q} \right)^{\sfrac 43} \Lambda_q^{\sfrac 23}  \left(\lambda_{m-1}^2\lambda_m^{-1}\right)^{-\sfrac 23}  \\ &\leq \Ga_{q+\half+1}^{\badshaq-11}\notag \, . 
\end{align}
The last two inequalities follow from \eqref{eq:imax:old}, \eqref{ineq:jmax:use} and \eqref{eq:par:div:2}. 
Then, we apply Corollary~\ref{rem:summing:partition} to $\theta=2$, $\theta_1= \sfrac23 $, $\theta_2= \sfrac43$, $H=\sigma^{\pm}_{\divH t_{i,j,k,\xi,\vecl,I,\diamond}^{m}}$, and $p=\sfrac 32$, which gives 
\begin{align}
    \left\| \psi_{i,q} \sigma^{\pm}_{\ov \phi_{TNW}^{m}} \right\|_{\sfrac 32} &\les \delta_\qbn^{\sfrac23} \tau_q^{-\sfrac 23} \Ga_q^{20+ \sfrac23 \CLebesgue} \left(\lambda_{m-1}^2\lambda_m^{-1}\right)^{-\sfrac 23} r_{\min(m,\qbn)}^{\sfrac 23} \leq\delta_{m+\bn} \Ga_{m}^{-10} \, . \notag 
\end{align}
from \eqref{eq:desert:ineq}. Combined with \eqref{eq:ct.p.2}, this verifies \eqref{eq:ct.p.3} for $q+\half+2\leq m' \leq q+\bn$. On the other hand, from Corollary~\ref{lem:agg.pt} with $H=\sigma^{\pm}_{\divH t_{i,j,k,\xi,\vecl,I,\diamond}^{m}}$, $\varpi=\Ga_{q+\half+1}^{\badshaq-11}\mathbf{1}_{\supp a_{\pxi,\diamond}\rhob_\pxi^\diamond\zetab_\xi^I}$ and $p=1$, we have that
\begin{align}
    \left\| \psi_{i,q} \sigma^{\pm}_{\ov \phi_O^{m}} \right\|_\infty 
    &\leq \Ga_{q+\half+1}^{\badshaq-10} \, . \notag 
\end{align}
Combined again with \eqref{eq:ct.p.2}, this verifies \eqref{eq:ct.p.3.1} at level $q+\half+1 \leq m' \leq q+\bn$. Lastly, we have that \eqref{eq:ctnl:estimate:1} at level $m'$ for $q+\half+1 \leq m' < q+\bn$ and for the nonlocal part of \eqref{ct.err.high} are satisfied by an argument essentially identical to that of the previous case.
\smallskip

\noindent\texttt{Case 4: Estimate of \eqref{ct.err.high2}.} Here we apply Proposition~\ref{prop:intermittent:inverse:div} with $p=\infty$ and the following choices.  The low-frequency assumptions in Part 1 are exactly the same as the $L^\infty$ low-frequency assumptions in the previous two steps.  For the high-frequency assumptions, we recall the choice of $N_{**}$ from \eqref{i:par:10} and set
\begin{align}
    &\varrho_R = (\Id - \tilde{\mathbb{P}}_{q+\bn+1}^\xi ) \mathbb{P}_{\neq 0} \left( \varrho_{\pxi,R}^I \right)^2 \, , \quad
    \varrho_\ph = (\Id - \tilde{\mathbb{P}}_{q+\bn+1}^\xi ) \mathbb{P}_{\neq 0} \left( \varrho_{\pxi,\ph}^I \right)^2r_q^{-\sfrac23} \, , \\
    &\vartheta^{i_1i_2\dots i_{\dpot-1}i_\dpot}_\diamond = \delta^{i_1i_2\dots i_{\dpot-1}i_\dpot} \Delta^{-\sfrac \dpot 2}\varrho_\diamond  \, , \quad \Lambda=\lambda_{q+\bn}  \, , \quad \dpot=0 \, , \\
    &\mu = \Upsilon=\Upsilon' = \lambda_{q+\half}\Ga_q \, , \quad \const_{*,\infty} = \left( \frac{\lambda_\qbn}{\lambda_{q+\bn+1}} \right)^{N_{**}} \lambda_\qbn^3 \, , \quad \Ndec \textnormal{ as in \eqref{i:par:9}}\, . \notag
\end{align}
Then we have that item~\eqref{item:inverse:i} is satisfied by definition, item~\eqref{item:inverse:ii} is satisfied as in the previous steps, \eqref{eq:DN:Mikado:density} is satisfied using Propositions~\ref{prop:pipeconstruction} and \ref{prop:pipe.flow.current} and \eqref{eq:remainder:inverse} from Lemma~\ref{lem:special:cases}, \eqref{eq:inverse:div:parameters:0} is satisfied by definition and as in the previous steps, and \eqref{eq:inverse:div:parameters:1} is satisfied by \eqref{condi.Ndec0}.  For the non-local assumptions, we choose $M_\circ,N_\circ=2\Nind$ so that \eqref{eq:inv:div:wut}--\eqref{eq:inverse:div:v:global:parameters} are satisfied as in Case 1, and \eqref{eq:riots:4} is satisfied from \eqref{ineq:Nstarz:1}. We have thus satisfied all the requisite assumptions, and we therefore obtain non-local bounds very similar to those from the previous steps, which are consistent with \eqref{eq:ctnl:estimate:1} at level $q+\bn$. We omit further details.
\end{proof}

\opsubsubsection{Transport/Nash current error from the divergence corrector part of the velocity increment}

\begin{lemma*}[\bf Current error and pressure increment from divergence correctors]\label{lem:ctdiv:general:estimate}
There exist vector fields $\ov\phi_{TNC}$ and a function $\bmu_{TNC}$ of time such that 
\begin{equation}
L_{TN} \left( w_{q+1}^{(p)} \otimes_s w_{q+1}^{(c)} \right) =  \div \left( \ov\phi_{TNC} \right) + \bmu_{TNC}' \, ,
 \qquad
  \ov\phi_{TNC}
  = \sum_{m=q+\half+1}^{\qbn} \div \ov\phi_{TNC}^m \,, 
\end{equation}
where $\ov\phi_{TNC}^m = \ov\phi_{TNC}^{m,l} + \ov\phi_{TNC}^{m,*}$ for $q+\half+1\leq m \leq \qbn$ satisfy the following.
\begin{enumerate}[(i)]
    \item For $q+\half+1\leq m \leq \qbn$, there exist functions $\si_{\ov\phi^m_{TNC}} = \si_{\ov\phi^m_{TNC}}^+ - \si_{\ov\phi^m_{TNC}}^-$ such that 
\begin{subequations}
\begin{align}
    \label{eq:ct.pc.1}
    \left|\psi_{i,q} D^N \Dtq^M \ov \phi^m_{TNC}\right| &\les \left( (\si_{\ov\phi^{m}_{TNC}}^+)^{\sfrac 32} r_m^{-1} + {\de_{q+3\bn}^2}\right) \left(\lambda_{m}\Gamma_{q}\right)^N \MM{M,\Nindt,\tau_q^{-1}\Gamma_{q}^{i+16},\Tau_q^{-1}{\Ga_q^9}}\\
    \label{eq:ct.pc.2}
    \left|\psi_{i,q} D^N \Dtq^M \si_{\ov\phi^m_{TNC}}^+\right| &\les \left(\si_{\ov\phi^m_{TNC}}^+ +\de_{q+3\bn}\right) \left(\lambda_{m}\Gamma_{q}\right)^N \MM{M,\Nindt,\tau_q^{-1}\Gamma_{q}^{i+17},\Tau_q^{-1}\Ga_q^9}\\
    \label{eq:ct.pc.3}
    \norm{\psi_{i,q} D^N \Dtq^M \si_{\ov\phi^m_{TNC}}^+}_{\sfrac32} &\lec \de_{m+\bn} \Gamma_{m}^{-9} \left(\lambda_{m}\Gamma_q\right)^N \MM{M,\Nindt,\tau_q^{-1}\Gamma_{q}^{i+17},\Tau_q^{-1}\Ga_q^9}\\
    \label{eq:ct.pc.3.1}
    \norm{\psi_{i,q} D^N \Dtq^M \si_{\ov\phi^m_{TNC}}^+}_{\infty} &\lec \Gamma_{m}^{\badshaq-9} \left(\lambda_{m}\Gamma_q\right)^N \MM{M,\Nindt,\tau_q^{-1}\Gamma_{q}^{i+17},\Tau_q^{-1}\Ga_q^9}\\
    \label{eq:ct.pc.4}
    \left|\psi_{i,q} D^N \Dtq^M \si_{\ov\phi^m_{TNC}}^-\right| &\lec \left(\frac{\la_q}{\la_{q+\floor{\bn/2}}}\right)^\frac23 \pi_q^q  \left(\lambda_{q+\floor{\bn/2}}\Gamma_q\right)^N \MM{M,\Nindt,\tau_q^{-1}\Gamma_{q}^{i+17},\Tau_q^{-1}\Ga_q^9}
\end{align}
\end{subequations}
for all $N,M \leq\sfrac{\Nfin}{100}$. Furthermore, we have that for $q+1\leq m' \leq m-1$ and $q+1\leq q'' \leq q+\half$,
\begin{align}\label{eq:ct.pc.6}
      \supp \si_{\ov\phi^m_{TNC}}^- \cap B\left( \supp \hat w_{q''}, \la_{q''}^{-1}\Ga_{q''+1} \right) =
    \supp \si_{\ov \phi^m_{TNC}}^+ \cap B\left(\supp \hat w_{m'}, \la_{m'}^{-1}\Ga_{m'+1} \right) = \emptyset \, . 
\end{align}
\item  For all $q+\half+1\leq m \leq \qbn$ and $q+1\leq q' \leq m-1$, the local parts satisfy
\begin{align}
B\left( \supp \hat w_{q'}, \lambda_{q'}^{-1} \Gamma_{q'+1} \right) \cap \supp \ov\phi_{TNC}^{m,l} = \emptyset \, . \label{ct.pc:support:first}
\end{align}
\item For all $q+\half+1\leq m \leq \qbn$ and $N,M\leq 2\Nind$, the non-local parts ${\ov\phi}^{m,*}_{TNC}$ satisfy
\begin{align}
\left\| D^N D_{t, q}^M {\ov\phi}^{m,*}_{TNC} \right\|_{L^\infty}
    &\leq { \Tau_{q+\bn}^{2\Nindt}}\delta_{q+3\bn}^{\sfrac32} \la_{m}^{N}\tau_{q}^{-M} \, .
    \label{eq:ct.pc:estimate:1}
\end{align}
\item For $M \leq 2\Nind$, the time function $\bmu_{TNC}$ satisfies
\begin{align}\label{sat:morn:11:14}
\bmu_{TNC}(t) = \int_0^t \langle\textnormal{(\ref{eq:ct:split:two})} (s)\rangle \, ds \, , \quad \left| \frac{d^{M+1}}{dt^{M+1}} \bmu_{TNC} \right| \leq \left( \max(1,T) \right)^{-1} \delta_{q+3\bn}^2 \MM{M,\Nindt, \tau_q^{-1}, \Tau_{q+1}^{-1}} \, .
\end{align}
\end{enumerate}
\end{lemma*}
\begin{proof}
The proof is similar to Step 2 of the proof of Lemma~\ref{l:divergence:corrector:error}. In fact, it is much simpler since the $\Dtq$ in $L_{TN}$ is always a ``good" derivative. We provide a few details below.

First note that
\begin{align*}
    L_{TN} \left(w_{q+1}^{(p)} \otimes_s w_{q+1}^{(c)}\right)  &= \sum_{\diamond,i,j,k,\xi,\vecl,I}
L_{TN} \bigg{[} a_{(\xi),\diamond} \left(\rhob_{(\xi)}^\diamond \zetab_{\xi}^{I,\diamond} \varrho_{(\xi),\diamond}^{I} \right)\circ \Phiik \xi^\ell \bigl( A_{\ell}^m  \epsilon_{\bullet p r} + A_{\ell}^\bullet  \epsilon_{m p r} \bigr) \notag\\
&\qquad \qquad \qquad \times \partial_p \left( a_{(\xi),\diamond} \left( \rhob_{(\xi)}^\diamond \zetab_\xi^{I,\diamond} \right)\circ \Phiik\right)  \partial_r \Phi_{(i,k)}^s (\mathbb{U}_{(\xi),\diamond}^{I})^s \circ \Phi_{(i,k)} 
\bigg{]}\\
&= \sum_{\diamond,i,j,k,\xi,\vecl,I} L_{TN} \bigg{[} a_{(\xi),\diamond} \left(\rhob_{(\xi)}^\diamond \zetab_{\xi}^{I,\diamond}  \right)\circ \Phiik \xi^\ell \bigl( A_{\ell}^m  \epsilon_{\bullet p r} + A_{\ell}^\bullet  \epsilon_{m p r} \bigr) \notag\\
&\qquad \qquad \qquad \times \partial_p \left( a_{(\xi),\diamond} \left( \rhob_{(\xi)}^\diamond \zetab_\xi^{I,\diamond} \right)\circ \Phiik\right)  \partial_r \Phi_{(i,k)}^s \bigg{]} (\varrho_{(\xi),\diamond}^{I} \mathbb{U}_{(\xi),\diamond}^{I})^s \circ \Phi_{(i,k)} \\
&=: \sum_{\diamond,i,j,k,\xi,\vecl,I} G_{\diamond,i,j,k,\xi,\vecl,I} (\varrho_{(\xi),\diamond}^{I} \mathbb{U}_{(\xi),\diamond}^{I})^s \circ \Phi_{(i,k)}
\end{align*}
We note that $(\varrho_{(\xi),\diamond}^{I} \mathbb{U}_{(\xi),\diamond}^{I})^s$ has mean $0$ (by property \eqref{item:pipe:means} of Proposition~\ref{prop:pipeconstruction} and \eqref{item:pipe:means:current} of Proposition~\ref{prop:pipe.flow.current}) and is $\frac{\T^d}{\la_{q+\half}\Ga_q}$ -periodic. So just as in the divergence corrector stress error, we apply the synthetic Littlewood-Paley decomposition suggested in \eqref{eq:decomp:showing} and define the current errors as follows:
\begin{align*}
    \ov\phi_{TNC}^{q+\half+1} &:= \sum_{\diamond,i,j,k,\xi,\vecl,I} \left(\divH + \divR \right) \left(G_{\diamond,i,j,k,\xi,\vecl,I} \tilde{\mathbb{P}}_{\lambda_{q+\half +1}} (\varrho_{(\xi),\diamond}^{I} \mathbb{U}_{(\xi),\diamond}^{I})^s \circ \Phi_{(i,k)}\right) \, ,\\
    \ov\phi_{TNC}^{m} &:= \sum_{\diamond,i,j,k,\xi,\vecl,I} \left(\divH + \divR \right) \left(G_{\diamond,i,j,k,\xi,\vecl,I} \tilde{\mathbb{P}}_{(\lambda_{m-1},\lambda_{m}]} (\varrho_{(\xi),\diamond}^{I} \mathbb{U}_{(\xi),\diamond}^{I})^s \circ \Phi_{(i,k)}\right)  \, ,\\
    \ov\phi_{TNC}^{q+\bn} &:= \sum_{m=q+\bn}^{q+\bn+1}\sum_{\diamond,i,j,k,\xi,\vecl,I} \left(\divH + \divR \right) \left(G_{\diamond,i,j,k,\xi,\vecl,I} \tilde{\mathbb{P}}_{(\lambda_{m-1},\lambda_{m}]} (\varrho_{(\xi),\diamond}^{I} \mathbb{U}_{(\xi),\diamond}^{I})^s \circ \Phi_{(i,k)}\right) \, , \\
    & \qquad + \sum_{\diamond,i,j,k,\xi,\vecl,I} \left(\divH + \divR \right) \left(G_{\diamond,i,j,k,\xi,\vecl,I} \left( \Id - \tilde{\mathbb{P}}_{\lambda_{q+\bn+1}} \right) (\varrho_{(\xi),\diamond}^{I} \mathbb{U}_{(\xi),\diamond}^{I})^s \circ \Phi_{(i,k)}\right) \, .
\end{align*}
We shall apply the inverse divergence operator to each term in the sum separately with the following choices. In all cases, we set 
\begin{align*}
G_R=\la_{q+\bn}^{-1}G_{R,i,j,k,\xi,\vecl,I}\, , \quad 
G_\ph=\la_{q+\bn}^{-1}r_q^{\sfrac23}G_{\ph,i,j,k,\xi,\vecl,I}\, . 
\end{align*}
We choose the high-frequency potentials as in \texttt{Step 2} of the proof of Lemma~\ref{l:divergence:corrector:error}, and choose the rest of parameters and functions required in Proposition~\ref{lem.pr.invdiv2.c} the same as in \texttt{Case 3} of the proof of Lemma~\ref{lem:ct:general:estimate}. 
In fact, the size of $G_{\diamond,1}$ and $G_{\diamond,\infty}$ is smaller than the one in \texttt{Case 3}.
By the same argument as in \texttt{Case 3}, we then get the same conclusion as in Lemma~\ref{lem:ct:general:estimate} for $\ov\phi_{TNC}^m$. We omit further details. 
\end{proof}

\opsubsubsection{Transport/Nash current error from oscillation, transport, Nash, divergence corrector, and mollification stress errors}

\begin{lemma*}[\bf Current error and pressure increment from \eqref{eq:ct:split:three}]\label{lem:ctosc:general:estimate}  
There exist vector field $\ov \phi_{TNS}$ and a function $\bmu_{TNS}$ of time such that
\begin{align*}
    \eqref{eq:ct:split:three} =- L_{TN} \left(S_O + S_{TN} + S_{C1} + S_{M2} \right)
    = \div \ov \phi_{TNS} + \bmu_{TNS}' \, , \qquad
    \ov \phi_{TNS} = \sum_{m=q+1}^{q+\bn} \ov \phi_{TNS}^m \, , 
\end{align*}
where $\ov\phi_{TNS}^m = \ov\phi_{TNS}^{m, l} + \ov\phi_{TNS}^{m, *} + \ov\phi_{TNS}^{ *}$ satisfies the following properties.
\begin{enumerate}[(i)]
    \item For $m=q+1,\, q+\half$, the local part $\ov\phi_{TNS}^m$ satisfies
\begin{align}
    \left| \psi_{i,q} D^N D_{t,q}^M \ov\phi_{TNS}^{m,l} \right| 
    &\lesssim \Ga_q^{-12} (\pi_q^m)^{\sfrac32} r_q^{-1} \la_m^N \MM{M,\Nindt, \tau_q^{-1}\Ga_q^{{i+14}}, \Tau_q^{-1}{\Ga_q^8}} \label{eq:ctn.osc.low} 
\end{align}
for $M, N \leq \sfrac{\Nfin}{100}$. 
\item For $m=q+\half+1,\dots,q+\bn$, there exist functions $\si_{\ov\phi_{TNS}^m} = \si_{\ov\phi_{TNS}^m}^+ - \si_{\ov\phi_{TNS}^m}^-$ such that
\begin{subequations}\label{eq:ctn.oc.general}
\begin{align}
\left| \psi_{i,q} D^N D_{t,q}^M \ov\phi_{TNS}^{m,l} \right|
&\lesssim  \left( (\si_{\ov\phi_{TNS}^m}^+)^{\sfrac 32} r_m^{-1} + {\de_{q+3\bn}^2} \right) \left(\lambda_{m}\Gamma_q\right)^N \MM{M,\Nindt,\tau_q^{-1}\Gamma_{q}^{i+{17}},\Tau_q^{-1}{\Ga_q^9}}\\
\left|\psi_{i,q} D^N \Dtq^M \si_{\ov\phi_{TNS}^m}^+ \right| &< \left( \si_{\ov\phi_{TNS}^m}^+ +\de_{q+3\bn}^2\right) \left(\lambda_{m}\Gamma_q\right)^N \MM{M,\Nindt,\tau_q^{-1}\Gamma_{q}^{i+{18}},\Tau_q^{-1}{\Ga_q^9}}\\
\norm{\psi_{i,q} D^N \Dtq^M \si_{\ov\phi_{TNS}^m}^+ }_{\sfrac32} &< \de_{m+\bn} \Gamma_{m}^{-9} \left(\lambda_{m}\Gamma_q\right)^N \MM{M,\Nindt,\tau_q^{-1}\Gamma_{q}^{i+{18}},\Tau_q^{-1}{\Ga_q^9}}\\
\norm{\psi_{i,q} D^N \Dtq^M \si_{\ov\phi_{TNS}^m}^+ }_{\infty} &< \Gamma_{q+\half+1}^{\badshaq-9} \left(\lambda_{m}\Gamma_q\right)^N \MM{M,\Nindt,\tau_q^{-1}\Gamma_{q}^{i+{18}},\Tau_q^{-1}{\Ga_q^9}}\\
\left|\psi_{i,q} D^N \Dtq^M \si_{\ov\phi_{TNS}^m}^- \right| &< \left(\frac{\la_q}{\la_{q+\floor{\bn/2}}}\right)^\frac23 \pi_q^q  \left(\lambda_{q+\floor{\bn/2}}\Gamma_q\right)^N \MM{M,\Nindt,\tau_q^{-1}\Gamma_{q}^{i+{18}},\Tau_q^{-1}{\Ga_q^9}} 
\end{align}
\end{subequations}
for $M, N \leq \sfrac{\Nfin}{200}$.
\item For $q+\half+1 \leq m \leq q+\bn$, $q+1\leq m' \leq m-1$, $q+1\leq q'' \leq q+\half$, $q+1\leq k \leq \qbn$, and $q+1\leq k' \leq k-1$, we have that
    \begin{subequations}\label{ctn:support:general}
    \begin{align}
      \supp \si_{\ov\phi_{TNS}^m}^- \cap B\left( \supp \hat w_{q''}, \la_{q''}^{-1}\Ga_{q''+1} \right) &=
    \supp \si_{\ov\phi_{TNS}^m}^+ \cap B\left(\supp \hat w_{m'}, \la_{m'}^{-1}\Ga_{m'+1} \right) = \emptyset \, . \label{ctn:support} \\
    B\left( \supp \hat w_{k'}, \lambda_{k'}^{-1} \Gamma_{k'+1} \right) \cap \supp \ov\phi_{TNS}^{k,l} &= \emptyset \, . \label{ct.p:support:first}
    \end{align}
    \end{subequations}
\item For $m=q+1,\dots,q+\bn$, the non-local parts satisfy
\begin{subequations}\label{eq:ctn.nonlocal}
    \begin{align}
\left\| D^N D_{t,q}^M \ov\phi_{TNS}^{m,*} \right\|_{\infty} 
    &\leq \Tau_{q+\bn}^{2\Nindt}\delta_{q+3\bn}^{\sfrac32} \la_{m}^{N}\tau_{q}^{-M} \, , \\
\left\| D^N D_{t,q+\bn-1}^M \ov\phi_{TNS}^* \right\|_\infty &\leq \de_{q+3\bn}^{\frac32}(\la_{q+\bn}\Gamma_{q+\bn-1})^N\MM{M,\Nindt, \tau_{q+\bn-1}^{-1}, \Tau_{q+\bn-1}^{-1}\Gamma_{q+\bn-1}}
\end{align}
\end{subequations}
for all $N,M \leq \sfrac{\Nind}4$.

\item For $M \leq 2\Nind$, the time function $\bmu_{TNS}$ satisfies
\begin{align}\label{sat:morn:11:24}
\bmu_{TNS}(t) = \int_0^t \langle\textnormal{(\ref{eq:ct:split:three})} (s)\rangle \, ds \, , \quad \left| \frac{d^{M+1}}{dt^{M+1}} \bmu_{TNS} \right| \leq \left( \max(1,T) \right)^{-1} \delta_{q+3\bn}^2 \MM{M,\Nindt, \tau_q^{-1}, \Tau_{q+1}^{-1}} \, .
\end{align}
\end{enumerate}
\end{lemma*}
\begin{proof}
Recall from \eqref{eq:ct:split:four} that \eqref{eq:ct:split:three} consists of $-L_{TN}(S_\triangle)$ where $\triangle$ represents $O$, $TN$, $C1$, or $M2$. We first consider the terms involving the local part of $S_\triangle$, and then deal with the terms with the non-local parts. 

\noindent\texttt{Case 1.} Current error from the terms $-L_{TN}(S_\triangle^{m,l})$ with $m=q+1$ or $m=q+\half$. In this case, we first note that $S_\triangle^{m,l}$ is non-trivial only when $\triangle=O$. Recall the expression of $S^{m,l}_O$ from \eqref{osc.loc.pot.1} of Remark~\ref{rem.ct.osc}, which gives
\begin{align*}
    L_{TN} S^{m,l}_O = \sum_{i,j,k,\xi,\vecl,\diamond} \sum_{j'=0}^{\const_{\mathcal{H}}} (L_{TN} H_{i,j,k,\xi,\vecl,\diamond}^{\alpha_{(j')}}) \rho_{i,j,k,\xi,\vecl,\diamond}^{\beta_{(j')}} \circ \Phi_{(i,k)} 
   \, . 
\end{align*}
In order to get the associated current error, we fix indices $j',\diamond,i,j,k,\xi,\vecl$ and apply the inverse divergence Proposition~\ref{prop:intermittent:inverse:div} and Remark~\ref{rem:pointwise:inverse:div} with the following choice of parameters and functions. Set 
\begin{align*}
    G= -(\la_{q+1}\Ga_q^{-4})^{-1} L_{TN} H_{i,j,k,\xi,\vecl,\diamond}^{\alpha_{(j')}}, \quad
    \varrho &= \la_{q+1}\Ga_q^{-4}\rho_{i,j,k,\xi,\vecl,\diamond}^{\beta_{(j')}}, \quad m=q+1\\
    G= -\la_{q+\half}^{-1} L_{TN} H_{i,j,k,\xi,\vecl,\diamond}^{\alpha_{(j')}}, \quad
    \varrho &= \la_{q+\half}\rho_{i,j,k,\xi,\vecl,\diamond}^{\beta_{(j')}}, \qquad m=q+\half\, . 
\end{align*}
We choose the rest of parameters and functions the same as in \texttt{Case 1} and \texttt{Case 2} in the proof of Lemma~\ref{lem:ct:general:estimate},
except for $N_* = \sfrac{\Nfin}{50}$ and $M_*= \sfrac{\Nfin}{100}$. 
\footnote{In fact, the actual size of $G$ is smaller than the one in \texttt{Case 1} and \texttt{Case 2}.} 
With this change, \eqref{eq:inv:div:NM} and \eqref{eq:inverse:div:parameters:0} still hold from \eqref{condi.Nfin0}. The rest of assumptions are satisfied as in \texttt{Case 1,2}. As a result, in the case of $m=q+1$ or $m=q+\half$, we obtain the associated current error $\ov \phi_{TNS}^m = \ov \phi_{TNS}^{m,l} + \ov \phi_{TNS}^{m, *}$ which satisfy 
\begin{align}
    \div \ov \phi_{TNS}^{m}  = -L_{TN} S_O^{m, l}
    +\langle L_{TN} S_O^{m, l} \rangle
\end{align}
and the same properties as $\ov\phi_{{TNW}}^m$ have, except that the range of $N$ and $M$ in the estimates are restricted to $N, M \leq \sfrac{\Nfin}{100}$. In particular, \eqref{eq:ctn.osc.low}, \eqref{ct.p:support:first} for $k=q+1, q+\half$, and \eqref{eq:ctn.nonlocal} with $m=q+1, q+\half$ hold. Finally, \eqref{sat:morn:11:24} holds due to similar arguments as in previous lemmas, and we omit further details throughout this proof.
\smallskip

\noindent\texttt{Case 2.} Current error and pressure increment from the terms $-L_{TN}(S_\triangle^{m,l})$ with $q+\half+1\leq m \leq  q+\bn$. Since $S_{M2}$ only have the non-local parts, we consider only when $\triangle = O, TN, C1$. Recall from Remarks~\ref{rem.ct.osc}, \ref{rem.ct.tn} and \ref{rem.ct.divcorr} that for $\triangle=O,TN,C1$, we have
\begin{align}\label{rep:LTNS}
    L_{TN} S_{\triangle}^{m,l} = \sum_{i,j,k,\xi,\vecl,I,\diamond} \sum_{j'=0}^{\const_{\mathcal{H}}} (L_{TN} H_{\triangle,i,j,k,\xi,\vecl,I,\diamond}^{\alpha_{(j')}})\rho_{\triangle,i,j,k,\xi,\vecl,I,\diamond}^{\beta_{(j')}} \circ \Phi_{(i,k)}
    \, . 
\end{align}
With this representation \eqref{rep:LTNS}, we fix indices $\triangle,j',\diamond,i,j,k,\xi,\vecl$ and apply Proposition~\ref{lem.pr.invdiv2.c} to construct desired current errors and pressure increments. 

\texttt{Case 2-1.} Consider $\triangle = O, C1$. Observe that $H_{\triangle,i,j,k,\xi,\vecl,I,\diamond}^{\alpha_{(j')}}$ and $\rho_{\triangle,i,j,k,\xi,\vecl,I,\diamond}^{\beta_{(j')}}$, $\triangle = O, C1$, have the same properties in Remark~\ref{rem.ct.osc}, \ref{rem.ct.divcorr}. Set the parameters and functions in the proposition the same as in \texttt{Case 3} in the proof of Lemma~\ref{lem:ct:general:estimate}, except for $N_* = \sfrac{\Nfin}{50}$, $M_*= \sfrac{\Nfin}{100}$,  
\begin{align*}
    G = -(\la_{q+\half}\Ga_q)^{-1}L_{TN} H_{\triangle,i,j,k,\xi,\vecl,I,\diamond}^{\alpha_{(j')}}\, , \quad
    \varrho &= \la_{q+\half}\Ga_q \rho_{\triangle,i,j,k,\xi,\vecl,I,\diamond}^{\beta_{(j')}}, \qquad \text{when }m=q+\half+1\\
    G = -(\la_{m-1}^{2}\la_m^{-1})^{-1}L_{TN} H_{\triangle,i,j,k,\xi,\vecl,I,\diamond}^{\alpha_{(j')}}\, , \quad
    \varrho &= \la_{m-1}^{2}\la_m^{-1}  \rho_{\triangle,i,j,k,\xi,\vecl,I,\diamond}^{\beta_{(j')}}, \qquad \text{otherwise }\, . 
\end{align*}
Then, \eqref{eq:inv:div:NM}, \eqref{eq:inverse:div:parameters:0}, \eqref{i:st:sample:wut:c}, \eqref{i:st:sample:wut:wut:c}, \eqref{eq:sample:prop:Ncut:3:c} still hold from \eqref{condi.Nfin0} and \eqref{condi.Nfin0}. The rest of assumptions are all satisfied as we see in \texttt{Case 3}. Therefore, as before, in each case of $m$, we obtain the associated current error $\ov\phi_{TN\triangle}^m=\ov\phi_{TN\triangle}^{m,l} + \ov\phi_{TN\triangle}^{m,*}$ and pressure increment $\si_{\ov\phi_{TN\triangle}^m} = \si_{\ov\phi_{TN\triangle}^m}^+ - \si_{\ov\phi_{TN\triangle}^m}^-$,
which satisfy
\begin{align}
    -L_{TN} S_{\triangle}^{m} + \left  \langle L_{TN} S_{\triangle}^{m} \right \rangle &= \div \ov\phi_{TN\triangle}^m, 
\end{align}
and share the same properties as $\ov\phi_{TNW}^m$ and $\si_{\ov\phi_{TNW}^m}$ have in the restricted range of $N, M$. In particular, \eqref{eq:ctn.oc.general}, \eqref{ctn:support:general}, and \eqref{eq:ctn.nonlocal} holds with the replacement of $\ov\phi_{TNS}^{m,l}$ and $\si_{\ov\phi_{TNS}^m}^{\pm}$ with $\ov\phi_{TN\triangle}^{m,l}$ and $\si_{\ov\phi_{TN\triangle}^m}^{\pm}$. 

\texttt{Case 2-2.} Consider $\triangle=TN$. Comparing the properties of $H_{\triangle,i,j,k,\xi,\vecl,I,\diamond}^{\alpha_{(j')}}$ and $\rho_{\triangle,i,j,k,\xi,\vecl,I,\diamond}^{\beta_{(j')}}$ in Remark~\ref{rem.ct.osc} when $m=q+\bn$ with those in Remark~\ref{rem.ct.divcorr}, one can see that
\begin{align*}
   G = -\la_{q+\bn}^{-1} L_{TN} H_{\triangle,i,j,k,\xi,\vecl,I,\diamond}^{\alpha_{(j')}}\, , \quad 
   \varrho = \la_{q+\bn} \rho_{\triangle,i,j,k,\xi,\vecl,I,\diamond}^{\beta_{(j')}}
\end{align*}
satisfies the same estimates as $G$ and $\varrho$ defined in \texttt{Case 2-1} when $m=q+\bn$, except that $G$ when $\triangle=TN$ has more expensive sharp material derivative cost by $\Ga_q$. Thereefore, repeating the same argument, we can obtain the associated current error $\ov\phi_{TN\triangle}^{q+\bn}=\ov\phi_{TN\triangle}^{q+\bn,l} + \ov\phi_{TN\triangle}^{q+\bn,*}$ and pressure increment $\si_{\ov\phi_{TN\triangle}^{q+\bn}} = \si_{\ov\phi_{TN\triangle}^{q+\bn}}^+ - \si_{\ov\phi_{TN\triangle}^{q+\bn}}^-$,
which satisfy
\begin{align}
    -L_{TN} S_{\triangle}^{{q+\bn}} + \left  \langle L_{TN} S_{\triangle}^{{q+\bn}} \right \rangle &= \div \ov\phi_{TN\triangle}^{q+\bn}, 
\end{align}
and share the same properties as $\ov\phi_{TNW}^{q+\bn}$ and $\si_{\ov\phi_{TNW}^{q+\bn}}$ have in the restricted range of $N, M$ expect that the sharp material derivative have extra $\Ga_q$ cost. In particular, \eqref{eq:ctn.oc.general}, \eqref{ctn:support:general}, and \eqref{eq:ctn.nonlocal} holds with the replacement of $\ov\phi_{TNS}^{{q+\bn},l}$ and $\si_{\ov\phi_{TNS}^{q+\bn}}^{\pm}$ with $\ov\phi_{TN\triangle}^{{q+\bn},l}$ and $\si_{\ov\phi_{TN\triangle}^{q+\bn}}^{\pm}$. 

Lastly, we define 
\begin{align*}
    \ov\phi_{TNS}^{m} :=\ov\phi_{TNO}^{m} + \ov\phi_{TNC1}^{m} + \ov\phi_{TNTN}^{m}, \quad
    \si_{\ov\phi_{TNS}^{m}} := 
    \si_{\ov\phi_{TNO}^{m}} + \si_{\ov\phi_{TNC1}^{m}}
    +\si_{\ov\phi_{TNT}^{m}}
\end{align*}
and the local and nonlocal parts of $\ov\phi_{TNS}^{m}$ and the superscript $\pm$ part of $\si_{\ov\phi_{TNS}^{m}}$ analogously. Here, we set undefined current errors
$\ov\phi_{TN\triangle}^{m}$ and pressure increments $\si_{\ov\phi_{TN\triangle}^{m}}=0$ as zero. Then, combining the analysis in \texttt{Case 2-1, 2-2}, \eqref{eq:ctn.oc.general}, \eqref{ctn:support:general}, and \eqref{eq:ctn.nonlocal} for $\ov\phi_{TNS}^{m,*}$ can be verified. 
\smallskip

\noindent\texttt{Case 3.} Current error from the terms $-L_{TN}(S_\triangle^{m, *})$ with $q+1\leq m \leq q+\bn$. Lastly, we construct $\ov\phi_{TNS}^*$ satisfying
\begin{align*}
    \div \ov\phi_{TNS}^* 
    =-\sum_{m=q+1}^{q+\bn}\mathbb{P}_{\neq 0} L_{TN} \left( S_{O}^{m,*} + S_{TN}^{m,*} + S_{C1}^{m,*} + S_{M2}^{m,*} \right)
\end{align*}
and \eqref{eq:ctn.nonlocal}. 
The terms on the right-hand side are not be intermittent, so there is no pressure increment generated from them. We fix $\triangle$ and $m$, and 
apply Remark~\ref{rem:inverse.div.spcial} of Proposition~\ref{prop:intermittent:inverse:div}. We first consider when $\triangle\neq M2$. Set $N_*=M_*=\Nind-1$, $M_\circ = N_\circ = \sfrac{\Nind}4$,
\begin{align*}
    &G = -L_{TN}  S_{\triangle}^{m,*}\, ,\quad \const_{G,\infty} = \tau_q^{-1}  \Tau_{q+\bn}^{4\Nindt} \de_{q+3\bn}, \quad \la = \la_{q+\bn},\quad \nu= \nu' = \Tau_q^{-1}\, ,\\
   & v=\hat u_q\,, \quad D_t=\Dtq\,, \quad \lambda'=\la_q\Ga_q
    \,, \quad\const_v = \Lambda_q^{\sfrac12}\, , 
\end{align*}
and choose a natural number $K_\circ$ such that
\begin{align*}
    \Tau_{q+\bn}^{2\Nindt}\de_{q+\bn}^{\sfrac32}
    \leq \la_{q+\bn}^{-K_\circ}
    \leq\Tau_{q+\bn}^{2\Nindt+1}\de_{q+\bn}^{\sfrac32}
\end{align*}
Then, all the assumptions are satisfied by \eqref{eq:Onpnp:estimate:1}, \eqref{eq:trans:nonlocal}, \eqref{eq:divER:nonlocal}, \eqref{eq:nasty:D:vq:old}, Corollary \ref{cor:deformation}. In particular, \eqref{eq:riots:4} can be verified by the choice of sufficiently large $a$. 
As a result of Remark~\ref{rem:inverse.div.spcial}, summing over $m$, we have $\ov\phi_{TN\triangle}^*$ which satisfies
\begin{align*}
\div \ov\phi_{TN\triangle}^*
= - \sum_{m=q+1}^{q+\bn}\mathbb{P}_{\neq 0} L_{TN} S^{m,*}_{\triangle}, \qquad
\norm{D^N D_{t,q}^M \ov\phi_{TN\triangle}^*}_\infty 
\leq  \Tau_{q+\bn}^{2\Nindt+1} \de_{q+3\bn}^{\sfrac32} \la_{q+\bn}^{N} \Tau_q^{-M}
\end{align*}
for $N,M\leq \sfrac{\Nind}4$. Lastly, we apply Lemma~\ref{lem:upgrading.material.derivative} to $\ov\phi_{TN\triangle}^*$, we have
\begin{align*}
    \norm{D^N D_{t,q+\bn-1}^M \ov\phi_{TN\triangle}^*}_\infty &\leq  \Tau_{q+\bn}^{\Nindt+1} \de_{q+3\bn}^{\sfrac32} \la_{q+\bn}^{N} (\Tau_{q+\bn-1}\Ga_{q+\bn-1})^{-M}\\
    &\leq \Tau_{q+\bn} \de_{q+3\bn}^{\sfrac32} \la_{q+\bn}^{N} \MM{M, \tau_{q-\bn-1}^{-1}, \Tau_{q+\bn-1}^{-1}}\, 
\end{align*}
for $N,M\leq \sfrac{\Nind}4$.

Next, we consider $\triangle = M2$. As we see from \eqref{ER:new:error:moll1}, $S^{m,*}_{M2}$ is non-trivial only when $m=q+\bn$. We first note that when $q+1\leq k <q+\bn$,  
\begin{align*}
    \norm{D^N D_{t,q+\bn-1}^M \hat w_{k}}_\infty
    =\norm{D^N D_{t,k}^M \hat w_{k}}_\infty
    \lec \Ga_{q}^{\sfrac{\badshaq}2+18} r_{q}^{-1} (\la_k \Ga_k)^N (\Tau_{k-1}^{-1}\Ga_{k-1})^M
\end{align*}
for $N+M \leq \sfrac{3\Nfin}2 +1$, from Hypothesis \ref{hyp:dodging1}, \eqref{eq:nasty:Dt:uq:orangutan}, \eqref{eq:imax:old}, and \eqref{ineq:b:second}. Also, applying Lemma \ref{lem:upgrading.material.derivative} to \eqref{eq:nasty:D:vq:old}, we have 
\begin{align*}
    \norm{D^ND_{t,q+\bn-1}^M \na \hat u_q}_\infty
    \lec \Tau_{q+\bn}^{-1}\la_q \Ga_q^{\sfrac{\badshaq}2+18} r_q^{-1} (\la_{q+\bn-1}\Ga_{q+\bn-1})^N (\Tau_{q+\bn-1}^{-1}\Ga_{q+\bn-1})^M
\end{align*}
Combining these with \eqref{est:stress.mollification}, we have from \eqref{eq:vellie:inductive:dtq-1:uniform:upgraded:statement} that
\begin{align*}
\norm{D^ND_{t,q+\bn-1}^M L_{TN} S_{M2}^{q+\bn,*}}_\infty &\leq \norm{D^ND_{t,q+\bn-1}^{M+1} S_{M2}^{q+\bn,*}}_\infty \notag\\
&\qquad \qquad +\norm{D^ND_{t,q+\bn-1}^M [((\hat w_{q+\bn-1}- \hat w_{q})\cdot\na)\tr + \na \hat u_q:] S_{M2}^{q+\bn,*}}_\infty\\
&\leq \Tau_{q+\bn}^{2\Nindt-2} \de_{q+3\bn} 
(\la_{q+\bn}\Ga_{q+\bn})^N
(\Tau_{q+\bn-1}^{-1}\Ga_{q+\bn-1})^M
\end{align*}
for $N+M\leq 2\Nind-1$. Therefore, we apply  Remark~\ref{rem:inverse.div.spcial} of Proposition~\ref{prop:intermittent:inverse:div} by setting $N_*=M_*=\Nind-1$, $M_\circ = N_\circ = \sfrac{\Nind}4$,
\begin{align*}
    &G = -L_{TN}  S_{\triangle}^{m,*}\, ,\quad \const_{G,\infty} = \Tau_{q+\bn}^{2\Nindt-2} \de_{q+3\bn}, \quad \la = \la_{q+\bn}\Ga_{q+\bn},\quad \nu= \nu' = \Tau_{q+\bn-1}^{-1}\Ga_{q+\bn-1}\, ,\\
   & v=\hat u_{q+\bn-1}\,, \quad D_t=D_{t,q+\bn-1}\,, \quad  \lambda'=\la_{q+\bn-1}\Ga_{q+\bn-1}
    \,, \quad\const_v = \Lambda_{q+\bn-1}^{\sfrac12}\, ,
\end{align*}
and choosing a natural number $K_\circ$ so that
\begin{align*}
  \de_{q+3\bn}^{\sfrac32}  \Tau_{q+\bn}^{\Nind }
  \leq (\la_{q+\bn}\Ga_{q+\bn})^{-K_\circ} 
  \leq 
  \de_{q+3\bn}^{\sfrac32}  \Tau_{q+\bn}^{\Nind+1}.
\end{align*}
Then all required assumptions are satisfied as before.  As a result of the remark, we obtain $\ov\phi_{TNM2}^*$ such that $\div \ov\phi_{TNM2}^*
= - \sum_{m=q+1}^{q+\bn} \mathbb{P}_{\neq 0} L_{TN} S^{m,*}_{M2}$, and for $N, M \leq \sfrac{\Nind}4$,
\begin{align*}
\norm{D^N D_{t, q+\bn-1}^M \ov\phi_{TNM2}^*}_\infty
&\leq   \Tau_{q+\bn}^{\Nindt+1} \de_{q+3\bn}^{\sfrac32} 
(\la_{q+\bn}\Ga_{q+\bn})^N
(\Tau_{q+\bn-1}^{-1}\Ga_{q+\bn-1})^M\\
&\leq \Tau_{q+\bn}\de_{q+3\bn}^{\sfrac32} 
(\la_{q+\bn}\Ga_{q+\bn})^N
\MM{M, \Nindt, \tau_{q+\bn-1}^{-1}, \Tau_{q+\bn-1}^{-1}\Ga_{q+\bn-1}}\, .
\end{align*}
Lastly, we set $\ov\phi_{TNS}^*: = \ov\phi_{TNO}^* + \ov\phi_{TNC1}^* + \ov\phi_{TNTN}^*+ \ov\phi_{TNM2}^*$ and collect the properties of $\ov\phi_{TN\triangle}^*$ to conclude \eqref{eq:ctn.nonlocal}. 
\end{proof}
\bigskip

\begin{remark*}[\bf Collecting pressure and current errors from transport-Nash]\label{rem:sat:morn} We now collect all current errors and pressure increments generated by \eqref{eq:ct:split:one}--\eqref{eq:ct:split:three} and set
\begin{align}\label{decomp:phTN}
    \ov\phi_{TN}^m :=  \ov\phi_{TNW}^m +\ov\phi_{TNC}^m + \ov\phi_{TNS}^m \, , \qquad
    \si_{\ov\phi_{TN}^m}
    := \si_{\ov\phi_{TNW}^m}
    +\si_{\ov\phi_{TNC}^m}
    +\si_{\ov\phi_{TNS}^m} \, ,
\end{align}
where the quantities on the right-hand side are constructed in Lemmas~\ref{lem:ct:general:estimate}, \ref{lem:ctdiv:general:estimate}, and \ref{lem:ctosc:general:estimate}. We use a similar notation for the various functions of time $\bmu$, so that recalling \eqref{bmun:def} and \eqref{bmut:def}, we have that $\bmu_T+\bmu_N = \bmu_{TNW}+\bmu_{TNC}+\bmu_{TNS}$.  Then summing over $m$, we have the transport and Nash current error $\ov\phi_{TN}$. We similarly collect the local and nonlocal parts of $\ov\phi_{TN}^m$ and the $\pm$ part of the pressure increments $\si_{\ov\phi_{TN}^m}$.\index{$\ov\phi_{TN}^m$}\index{$\sigma_{\ov\phi_{TN}^m}$}
\end{remark*}

Lastly, we define and analyze the current error associated to the pressure increments $\sigma_{\ov\phi_{TN}^m}$. 

\begin{lemma*}[\bf Pressure current]\label{lem:ctn:pressure:current}
For every $m'\in\{q+\half+1,\dots,q+\bn\}$, there exists a current error $\phi_{{\ov\phi_{TN}^{m'}}}$ associated to the pressure increments $\si_{\ov\phi_{TN}^{m'}}$ and a function $\bmu_{\si_{\ov\phi_{TN}^{m'}}}$ of time
that satisfy the following properties.  
\begin{enumerate}[(i)]
\item\label{i:ctn:pc:2} We have the decompositions and equalities
\begin{subequations}
\begin{align}\label{eq:ctn:desert:decomp}
    \div \phi_{{\ov\phi_{TN}^{m'}}} + \bmu_{\si_{\ov\phi_{TN}^{m'}}}'
    &= D_{t,q}  \si_{\ov\phi_{TN}^{m'}}\, , \\
    \phi_{{\ov\phi_{TN}^{m'}}} = \phi_{{\ov\phi_{TN}^{m'}}}^* + \sum_{m=q+\half+1}^{{m'}} \phi_{{\ov\phi_{TN}^{m'}}}^{m} \, , \qquad 
    \phi_{\ov \phi_{TN}^{m'}}^{m} &= \phi_{\ov \phi_{TN}^{m'}}^{m,l} + \phi_{\ov \phi_{TN}^{m'}}^{m,*} \, .           
\end{align}
\end{subequations}
\item\label{i:pc:3:TN} For $q+\half+1 \leq m \leq m'$ and $N,M\leq  2\Nind$,
\begin{subequations}
\begin{align}
    &\left|\psi_{i,q} D^N \Dtq^M \ov\phi_{\ov \phi_{TN}^{m'}}^{m,l} \right| < \Ga_{m}^{-100} \left(\pi_q^m\right)^{\sfrac 32} r_m^{-1} (\la_m \Ga_m^2)^M \MM{M,\Nindt,\tau_q^{-1}\Ga_q^{i+{18}},\Tau_q^{-1}\Ga_q^9} \, , \label{eq:desert:estimate:1:TN} \\
    &\left\| D^N \Dtq^M \phi_{\ov \phi_{TN}^{m'}}^{m,*} \right\|_\infty  < \Tau_\qbn^{2\Nindt} \delta_{q+3\bn}^{\sfrac 32} (\la_{m'}\Ga_{m'}^2)^N \tau_q^{-M},\label{eq:desert:estimate:21:TN}\\
    &\left\| D^N\Dtq^M \phi_{\ov\phi_{TN}^{m'}}^{*} \right\|_\infty < \Tau_\qbn^{2\Nindt} \delta_{q+3\bn}^{\sfrac 32} (\la_{q+\bn}\Ga_{q+\bn}^2)^N \tau_q^{-M} \label{eq:desert:estimate:22:TN} \, .
\end{align}
\end{subequations}
\item\label{i:pc:4:TN} For all $q+\half+1\leq m \leq m'$ and all $q+1\leq q' \leq m-1$, 
\begin{align}
        B\left( \supp \hat w_{q'}, \sfrac 12 \lambda_{q'}^{-1} \Ga_{q'+1} \right) \cap \supp \left( \phi^{m,l}_{\ov \phi_{TN}^{m'}} \right) = \emptyset \label{eq:desert:dodging:TN} \, .
\end{align}
\item \label{i:pc:5:TN} For $M\leq 2\Nind$, the mean part $\bmu_{\si_{\ov\phi_{TN}^{m'}}}$ satisfies
\begin{align}\label{eq:desert:mean:tn}
    \left|\frac{d^{M+1}}{dt^{M+1}}
    \bmu_{\si_{\ov\phi_{TN}^{m'}}} \right| 
    \leq (\max(1, T))^{-1}\delta_{q+3\bn} \MM{M,\Nindt,\tau_q^{-1},\Tau_{q+1}^{-1}} \, .
\end{align}
\end{enumerate}
\end{lemma*}
\begin{proof}
From \eqref{decomp:phTN}, the pressure increment $\si_{\ov \phi_{TN}^{m'}}$ consists of $\si_{\ov \phi_{TNW}^{m'}}$, $\si_{\ov \phi_{TNC}^{m'}}$, $\si_{\ov \phi_{TNS}^{m'}}$.  Consider the pressure current for the pressure increment $\si_{\ov \phi_{TNW}^{m'}}$ defined in \texttt{Case 3} of the proof of Lemma~\ref{lem:ct:general:estimate}. As a result of the application of Proposition~\ref{lem.pr.invdiv2.c} to $t_{i,j,k,\xi,\vecl,I,\diamond}^{m'}$, from Part 4 of the proposition, we obtain a pressure current  $\phi_{i,j,k,\xi,\vecl,I,\diamond}$ which has a decomposition 
$$
\ov\phi_{i,j,k,\xi,\vecl,I,\diamond}=\ov\phi_{i,j,k,\xi,\vecl,I,\diamond}^*+\sum_{m=0}^{ {\bar m}} \ov\phi^{m}_{i,j,k,\xi,\vecl,I,\diamond}
= \left(\divH + \divR \right) \Dtq \si_{\divH t^{m'}_{i,j,k,\xi,\vecl,I,\diamond}}. $$
Noticing that the estimates for the pressure increment $\sigma_{\ov\phi_{TNW}^{m'}}$ are similar to those of the pressure increments for the Reynolds stress errors, for example those defined in Lemma~\ref{lem:oscillation:pressure}, we can obtain pointwise estimates for $\ov\phi^{m,l}_{\ov\phi_{TNW}^{m'}}$ analogous to those contained in Lemma~\ref{lem:oscillation:pressure:current}. The properties in \eqref{eq:ctn:desert:decomp}--\eqref{eq:desert:mean:tn} follow from similar arguments as before.  We refer also to~\cite{GKN23}, in which a number of error terms are estimate and analyzed using Proposition~\ref{lem.pr.invdiv2.c}.
\end{proof}

\opsubsection{Mollification current error}\label{op:moll:ss}

Similar to the case of the stress mollification errors, we will have to consider various mollification errors that go into the new unresolved current. These are listed below and are estimated in an analogous way to the mollification stress errors.

We recall the operators $\divR$ from \eqref{eq:inverse:div:error:stress} and $L_{TN}$ from \eqref{def:LTN} and regroup the terms by setting
\begin{align*}
    \ov \phi_M^{q+1}&:= \ph_q^q - \ph_\ell \\
    \ov \phi_{M3}^{q+\bn}&:=\frac12 \left( |\hat w_{q+\bn}|^2\hat w_{q+\bn} - |w_{q+1}|^2 w_{q+1}\right)
    \\
    \ov \phi_{M4}^{q+\bn}&:= \divR  \left[L_{TN} \left(\hat w_{q+\bn}\otimes \hat w_{q+\bn} - w_{q+1}\otimes w_{q+1}\right) + (\hat w_{q+\bn} - w_{q+1}) \cdot (\pa_t u_q + (u_q\cdot \na) u_q + \na p_q ) \right]
    \, . 
\end{align*}
We also define
\begin{align}
    \ov \phi_M^{q+\bn} := \ov \phi_{M3}^{q+\bn} +\ov \phi_{M4}^{q+\bn} \, ,
\end{align}
and we set
\begin{equation}
    \bmu_{M4}(t) := \int_0^t \left\langle L_{TN} \left(\hat w_{q+\bn}\otimes \hat w_{q+\bn} - w_{q+1}\otimes w_{q+1}\right) + (\hat w_{q+\bn} - w_{q+1}) \cdot (\pa_t u_q + (u_q\cdot \na) u_q + \na p_q ) \right\rangle(s) \, ds  \, .
\end{equation}
For details on how these error terms appear in the relaxed local energy inequality, we refer to~\cite[subsection~5.1]{GKN23}.

\begin{lemma*}[\bf Basic estimates and applying inverse divergence]\label{lem:moll:curr} For all $N+ M\leq \sfrac{\Nind}4$, the mollification errors $\ov \phi_M^{q+1}$ and $\ov \phi_M^{q+\bn}$ satisfy 
\begin{subequations}
    \begin{align}
        \norm{D^N D_{t,q}^M \ov \phi_M^{q+1}}_\infty
        &\leq \delta_{q+3\bn}^{\sfrac32} \lambda_{q+1}^N \MM{M,\Nindt, \tau_{q}^{-1},\Gamma_{q}^{-1}\Tau_{q}^{-1}} \, , \label{est:curr.mollification1}\\
        \norm{ D^N D_{t,q+\bn-1}^M \ov \phi_M^{q+\bn}}_\infty
        &\leq \Ga_\qbn^9 {\delta_{q+3\bn}^{\sfrac 32}}\Tau_\qbn^{2\Nindt} \left(\lambda_{q+\bn}\Gamma_{q+\bn}\right)^N \MM{M, \NindRt, \tau_{q+\bn-1}^{-1}, \Tau_{q+\bn-1}^{-1}\Ga_{\qbn-1} } \, . 
    \label{est:curr.mollification}
    \end{align}
\end{subequations}
In addition, the mean portion $\bmu_{M4}$ satisfies
\begin{align}\label{eq:sat:evening:6:32}
        \left| \frac{d^{M+1}}{dt^{M+1}}\bmu_{M4} \right|
        \leq (\max(1,T))^{-1} \de_{q+3\bn} \MM{M,\Nindt,\tau_q^{-1},\Tau_{q+1}^{-1}} \quad \text{for }M\leq \sfrac{\Nind}4 \, .
\end{align}

\end{lemma*}
\begin{proof}[Proof of Lemma~\ref{lem:moll:curr}]
We have that \eqref{est:curr.mollification1} follows immediately from \eqref{eq:phicomm:bounds}.  Next, in order to handle $\ov \phi_{M3}^{q+\bn}$, we recall from \eqref{eq:diff:moll:vellie:statement} that
    \begin{align*}
    \norm{D^N D_{t,q+\bn-1}^M \left(w_{q+1}- \hat w_{q+\bn}\right)}_\infty \lec \delta_{q+3\bn}^3 \Tau_\qbn^{25\NindRt} \left(\lambda_{q+\bn}\Gamma_{q+\bn-1}\right)^N \MM{M, \NindRt, \tau_{q+\bn-1}^{-1}, \Tau_{q+\bn-1}^{-1} }  \, .
\end{align*}
for all $N+M\leq \sfrac{\Nfin}{4}$. Using Lemma \ref{lem:dodging}, we note that $D_{t,q+\bn-1}w_{q+1} = D_{t,q}w_{q+1}$ and $D_{t,q+\bn-1}\hat w_{q+\bn} = D_{t,q}\hat w_{q+\bn}$. Then writing 
$$|\hat w_{q+\bn}|^2 \hat w_{q+\bn} - |w_{q+1}|^2 w_{q+1}
=(\hat w_{q+\bn}-w_{q+1}) | \hat w_{q+\bn}|^2 + w_{q+1} (\hat w_{q+\bn}-w_{q+1})\cdot \hat w_{q+\bn} + w_{q+1} w_{q+1} \cdot (\hat w_{q+\bn}-w_{q+1}) $$
and using \eqref{eq:vellie:upgraded:statement}, \eqref{eq:diff:moll:vellie:statement}, and \eqref{eq:moll:vel:threesie},
we have that for all $N+M\leq 2\Nind$,
\begin{align}
    &\norm{
    D^N D_{t,q+\bn-1}^M
    [|\hat w_{q+\bn}|^2 \hat w_{q+\bn} - |w_{q+1}|^2 w_{q+1}]
    }_\infty\notag\\
    &\qquad\leq \delta_{q+3\bn} \Tau_\qbn^{2\NindRt} \left(\lambda_{q+\bn}\Gamma_{q+\bn}\right)^N \MM{M, \NindRt, \tau_{q+\bn-1}^{-1}, \Tau_{q+\bn-1}^{-1}\Ga_{\qbn-1} }\, .
    \label{est:CM21.proof}
\end{align}
As for the remaining term $\ov \phi_{M4}^{q+\bn}$, we first upgrade the material derivative in the estimate for $\hat u_q$.
Applying Lemma \ref{lem:upgrading.material.derivative} to $F^l = 0$, $F^* = \hat u_q$, $k=q+\bn$, $N_\star = \sfrac{3\Nfin}{4}$ with \eqref{eq:bob:Dq':old} and using \eqref{v:global:par:ineq}, we have that
\begin{align*}
    \norm{ D^N D_{t,q+\bn-1}^M \hat u_q}_\infty
    \lec \Tau_{q}^{-1} \lambda_{q+\bn}^N \Tau_{q+\bn-1}^{-M} \, .
\end{align*}
We can now tackle the part of the error term that involves $L_{TN}$. To estimate this, we use Remark \ref{rem:inverse.div.spcial} with \eqref{v:global:par:ineq}, setting 
\begin{align*}
    G &=L_{TN} \left(\hat w_{q+\bn}\otimes \hat w_{q+\bn} - w_{q+1}\otimes w_{q+1}\right)
    , \quad v = \hat u_{q+\bn-1}
\\
    \const_{G, \infty}&= \de_{q+3\bn} \Tau_{q+\bn}^{2\Nindt}, \quad
    \la=\la' = \la_{q+\bn}\Ga_{q+\bn}, \quad 
    M_t = \Nindt,\quad
    \nu =\nu' = \Tau_{q+\bn}^{-1}, \quad \const_v = \La_{q+\bn-1}^{\sfrac12} \\
    &\hspace{2cm}N_* = \sfrac{\Nfin}{9}, \quad
M_* = \sfrac{\Nfin}{10}, \quad
N_\circ = M_\circ = 2\Nind\,  .
\end{align*}
As a result, with a suitable choice of positive integer $K_\circ$ so that
\begin{align*}
    \de_{q+3\bn} \Tau_{q+\bn}^{2\Nindt }\la_{q+\bn}^{5} 2^{2\Nind}
    \leq \la_{q+\bn}^{-K_\circ}
    \leq \de_{q+3\bn} \Tau_\qbn^{\NindRt}\, ,
\end{align*}
we find that for all $N+M\leq 2\Nind$,
\begin{align}
    \norm{D^N D_{t,q+\bn-1}^M\divR \left[L_{TN} \left(\hat w_{q+\bn}\otimes \hat w_{q+\bn} - w_{q+1}\otimes w_{q+1}\right) \right] }_\infty &\lec \de_{q+3\bn} \Tau_\qbn^{\Nindt} (\la_{q+\bn}\Ga_\qbn)^N \Tau_{q+\bn}^{-M} \notag \\
    &\leq \de_{q+3\bn}(\la_{q+\bn}\Ga_\qbn)^N \MM{M, \Nindt,\tau_{q+\bn-1}^{-1} ,\Tau_{q+\bn}^{-1}}\, . \label{est:CM22.proof}
\end{align}
The estimate for the mean portion follows in the usual way from Remark~\ref{rem:est.mean}.

Now we deal with the other part of the error term. Recall from \eqref{eqn:ER} that
$$ \pa_t u_q + (u_q\cdot \na) u_q + \na p_q = \div ( R_q - \pi_q \Id ) \, . $$
We apply Lemma~\ref{rem:no:decoup:inverse:div2} with the following choices:
\begin{align*}
    &G= \div \left( R_q - \pi_q \Id \right)^\bullet \, , \quad \varrho = \vartheta = (\hat w_{q+\bn} - w_{q+1})^\bullet \,, \quad v=\hat u_{q+\bn-1} \, , \quad \lambda' = \lambda_{q+\bn-1} \Ga_{q+\bn-1} \, ,  \\
    &{\nu = \nu' = \Tau_{q+\bn-1}^{-1}\Ga_{q+\bn-1}^2} \, , \quad N_* = \sfrac{\Nind}{2} \, , \quad M_* = \sfrac{\Nind}{2} \, , \quad \dpot=0 \, , \quad \la = \La_{q+\bn} \Ga_{q+\bn} \, , \\
    &\pi' = \const_{*,\infty} = \delta_{q+3\bn}^3 \Tau_\qbn^{25\Nindt} \, ,\quad \Omega= \T^3 \times \mathbb{R} \, , \quad \pi = \Ga_{q+\bn-1} \pi_q \La_{q+\bn-1} \, , \quad M_t = \Nindt \, , \\
    &\Upsilon = \La = \la_{q+\bn}\Ga_{\qbn-1} \, , \quad M_\circ = N_\circ = \sfrac{\Nind}4 \, , \quad K_\circ \textnormal{ such that $\Tau_{\qbn}^{-10\Nindt} \leq \La^{K_\circ} \leq \Tau_{\qbn}^{-10\Nindt-1}$} \, . 
\end{align*}
The analysis here is similar to the analysis for the  nonlocal transport-Nash current errors, and so we omit the details but note that one can easily check that \eqref{eq:inv:div:extra:pointwise:noflow}, \eqref{eq:inv:div:extra:pointwise2:noflow}, and \eqref{parameter:noflow} are satisfied. Since $\dpot=0$, we move straight to the non-local assumptions and output, which again can be easily checked by direct computation or using similar arguments as for other nonlocal error terms. We therefore have from \eqref{eq:inverse:div:error:stress:bound:no:flow} that for $N+M\leq \sfrac{\Nind}{4}$, 
\begin{align}
    &\left\| D^N D_{t,\qbn-1}^M \divR \left( \div \left( R_q - \pi_q \Id \right)^\bullet (\hat w_{q+\bn} - w_{q+1})^\bullet \right) \right\|_\infty \notag\\
    &\qquad\les \Tau_\qbn^{3\Nindt} \delta_{q+3\bn}^3 (\la_\qbn\Ga_\qbn)^N \MM{M,\Nindt, \tau_{\qbn-1}^{-1}, \Tau_{\qbn-1}^{-1}\Ga_{\qbn}}\label{moll:eckel:nonlocal:prep:2}
\end{align}
Promotion of the material derivatives again follows standard arguments and Lemma~\ref{lem:upgrading.material.derivative}, and we omit further details.
\end{proof}

\section{Inductive cutoffs}\label{sec:inductive.cutoffs}

In this section, we define the new partition of unity $\{\psi_{i,\qbn}\}$ and verify the inductive properties from subsection~\ref{sec:cutoff:inductive}.  At the same time, we verify the inductive velocity bounds from subsection~\ref{sec:inductive:secondary:velocity}. The strategy for these proofs follows quite closely the strategy from \cite[subsections~6.1, 6.2]{BMNV21}.  However, the proofs now use $L^3$ inductive information, rather than $L^2$ inductive information.  Thus for the sake of completeness and for the accuracy of the constants chosen in subsection~\ref{sec:not.general} and~\ref{sec:para.q.ind}, which do depend on the computations in this section, we have included full details of all the proofs.\index{velocity cutoffs}

\subsection{New mollified velocity increment and definition of the velocity cutoff functions}
\label{sec:cutoff:velocity:properties}

We first recall the definition of $\hat w_{q+\bn}$ in \eqref{def.w.mollified}. We have that for a mollifier $\mathcal{\tilde P}_{q+\bn,x,t}$ at spatial scale $\lambda_{q+\bn}^{-1}\Gamma_{q+\bn-1}^{-\sfrac 12}$ and temporal scale $\Tau_{q+1}^{-1}$, we have 
\begin{equation}
    \hat w_{q+\bn} = \mathcal{\tilde P}_{q+\bn,x,t} w_{q+1} \, .
\end{equation}

Before defining the velocity cutoff functions, we need the following translations between $\Ga_{q'-1}$ and $\Ga_{q'}$.

\begin{definition}[\bf Translating $\Gamma$'s between $q'-1$ and $q'$]\label{def:istar:j}
Given $i,j,q' \geq 0$, we define
\begin{subequations}
\begin{align*}
i_* = i_*(j,q') &= i_*(j) = \min\{ i \geq 0 \colon \Gamma_{q'}^{i} \geq \Gamma_{q'-1}^{j} \} \\
j_*(i,q') &= \max\{ j: i_*(j) \leq i \}  \, .
\end{align*}
\end{subequations}
A consequence of this definition is the inequality
\begin{equation}\label{eq:ineq:ij}
    \Gamma_{q'}^{i-1} < \Gamma_{q'-1}^{j_*(i,q)} \leq \Gamma_{q'}^i \, .
\end{equation}
\noindent We also note that for $j=0$, we have that $i_*(j)=0$.  Finally, a simple computation shows that $i_*(j)$ has an upper bound which depends on $j$ but not $q$.
\end{definition}

We may now define the velocity cutoff functions using the cutoff functions presented in Lemma~\ref{lem:cutoff:construction:first:statement}, although $\Gamma_q$ will be replaced with $\Gamma_{q+\bn}$ throughout.

\begin{definition}[\bf Intermediate cutoff functions]\label{def:intermediate:cutoffs}
For stage $q+1$ of the iteration where $\qbn \geq 1$, $m\leq\NcutSmall$, and $j_m\geq 0$, we define\index{$\psi_{m,i_m,j_m,q}$}
\begin{align}
    h_{m,j_m,q+\bn}^2(x,t) &= \Gamma_{q+\bn}^{-2i_*(j_m)} \delta_{q+\bn}^{-1} r_q^{\sfrac 23} \left(\tau_{q+\bn-1}^{-1}\Gamma_{q+\bn}^{i_*(j_m)+2}\right)^{-2m} \sum_{N=0}^{ \NcutLarge} \left( \lambda_{q+\bn} \Gamma_{q+\bn} \right)^{-2N} \left| D^N D_{t,q+\bn-1}^m \hat w_{q+\bn} \right|^2 
    \, . 
\label{eq:h:j:q:def}
\end{align}
\noindent We then define $\psi_{m,i_m,j_m,q+\bn}$ by
\begin{align}
 \psi_{m,i_m,j_m,q+\bn}(x,t) 
 &= \gamma_{m,q+\bn}  \left( \Gamma_{q+\bn}^{-2(i_m-i_*(j_m))(m+1)} h_{m,j_m,q+\bn}^2  (x,t) \right) 
\label{eq:psi:i:j:def}
\end{align}
for $i_m> i_*(j_m)$, 
while for $i_m=i_*(j_m)$,
\begin{align}
 \psi_{m,i_*(j_m),j_m,q+\bn}(x,t) 
 &= \tilde \gamma_{m,q+\bn} \left( h_{m,j_m,q+\bn}^2(x,t) \right) \, .
\label{eq:psi:i:i:def}
\end{align}
The intermediate cutoff functions $\psi_{m,i_m,j_m,q+\bn}$ are equal to zero for $ i_m < i_*(j_m)$.  
\end{definition}
The idea of the intermediary cutoffs $\psi_{m,i_m,j_m,\qbn}$ and $i_m$ and $j_m$ is as follows.  First, we use the subscript $m$ to emphasize that $\psi_{m,i_m,j_m,\qbn}$ is using non-negative integers $i_m$ and $j_m$ to quantify the size of $D_{t,\qbn-1}^m \hat w_{\qbn}$, i.e. $m$ material derivatives applied to $\hat w_\qbn$.  Second, all proofs will have to be written using information from the old velocity cutoffs $\psi_{j_m, \qbn-1}$, which we index with $j_m$ (see Definition~\ref{def:psi:m:im:q:def}). Finally, the new velocity cutoffs will be defined in Definition~\ref{def:psi:i:q:def} using the integer $i$, which is equal to the supremum over $0\leq m \leq \NcutSmall$ of the integer $i_m$ being used to quantify the cost of $D_{t,\qbn-1}^m$. Later, $i_m$ which will be shown to take values no larger than $i_{\rm max}$. With these definitions and using \eqref{eq:tilde:partition} and \eqref{eq:psi:support:base:case}, it follows that
\begin{align}
\sum_{i_m\geq0} \psi_{m,i_m,j_m,q+\bn}^{6} = \sum_{i_m\geq i_*(j_m)} \psi_{m,i_m,j_m,q+\bn}^{6} = \sum_{\{ i_m \colon \Gamma_{q+\bn}^{i_m} \geq \Gamma_{q+\bn-1}^{j_m} \}} \psi_{m,i_m,j_m,q+\bn}^{6} \equiv 1
\label{eq:psi:i:j:partition:0}
\end{align}
for any $m$, and for $|i_m-i'_m|\geq 2$, 
\begin{equation}\label{eq:intermediate:overlapping}
  \psi_{m,i_m,j_m,q+\bn}\psi_{m,i_m',j_m,q+\bn}=0 \, .
\end{equation}

\begin{definition}[\bf $m^{\textnormal{th}}$ Velocity Cutoff Function]\label{def:psi:m:im:q:def}
At stage $q+1$ and for $i_m\geq 0$, we inductively define the $m^{\textnormal{th}}$ velocity cutoff function\index{$\psi_{m,i_m,\qbn}$}
\begin{equation}\label{eq:psi:m:im:q:def}
\psi_{m,i_m,q+\bn}^{6} = \sum\limits_{\{j_m\colon i_m\geq i_*(j_m)\}} \psi_{j_m,q+\bn-1}^{6} \psi_{m,i_m,j_m,q+\bn}^{6} \, .
\end{equation}
\end{definition}

We shall employ the notation
\begin{equation}\label{eq:i:tuple:def}
  \Vec{i} =  \{i_m\}_{m=0}^{\NcutSmall} = \left( i_0,...,i_{\NcutSmall} \right) \in \mathbb{N}_0^{\NcutSmall+1}
\end{equation}
to signify a tuple of non-negative integers of length $\NcutSmall+1$.

\begin{definition}[\bf Velocity cutoff function]
\label{def:psi:i:q:def}
At stage $q+1$ and for $0 \leq i \leq i_{\rm max}$,\index{$\psi_{i,q}$}\index{$i$} we define
\begin{align}
    \psi_{i,q+\bn}^{6} = \sum\limits_{\left\{\Vec{i}\colon\max\limits_{0\leq m\leq\NcutSmall} i_m =i\right\}} \prod\limits_{m=0}^{\NcutSmall} \psi_{m,i_m,q+\bn}^{6} \, .
  \label{eq:psi:i:q:recursive}
\end{align}
\end{definition}
\noindent For $\Vec{i}$ as in the sum of \eqref{eq:psi:i:q:recursive}, we shall denote 
\begin{align}
\supp\left( \prod\limits_{m=0}^{\NcutSmall} \psi_{m,i_m,q+\bn}  \right) = \bigcap_{m=0}^{\NcutSmall} \supp(\psi_{m,i_m,q+\bn}) =: \supp (\psi_{\Vec{i},q+\bn} )\,.
\label{eq:new:supp:notation}
\end{align}
This implies that $(x,t) \in \supp(\psi_{i,q+\bn})$ if and only if there exists $\Vec{i}\in \N_0^{\NcutSmall+1}$ such that $\max_{0\leq m\leq\NcutSmall} i_m =i$, and $(x,t) \in \supp(\psi_{\Vec{i},q+\bn})$.

\subsection{Partitions of unity, dodging, and simple bounds on velocity increments}


\begin{lemma}[\bf $\psi_{m,i_m,q+\bn}$ - Partition of unity]
\label{lem:partition:of:unity:psi:m}
For all $m$, we have that 
\begin{align}
\sum_{i_m\geq 0} \psi_{m,i_m,q+\bn}^{6}\equiv 1\,, \qquad \psi_{m,i_m,q+\bn}\psi_{m,i'_m,q+\bn}=0\quad\textnormal{for}\quad|i_m-i'_m|\geq 2 \, . \label{eq:lemma:partition:1}
\end{align}
\end{lemma}
\begin{proof}[Proof of Lemma~\ref{lem:partition:of:unity:psi:m}]
The proof proceeds inductively in a manner very similar to the proof of \cite[Lemma~6.7]{BMNV21}.  
To show the first part of \eqref{eq:lemma:partition:1}, we may use \eqref{eq:psi:i:j:partition:0} and \eqref{eq:psi:m:im:q:def} and reorder the summation to obtain
\begin{align*}
\sum_{i_m\geq 0} \psi_{m,i_m,q}^6 =&   \sum_{i_m\geq 0} \sum_{\{ j_m\colon i_*(j_m) \leq  i_m\}}~
\psi_{j_m,q-1}^6 \psi_{m,i_m,j_m, q}^6(x,t) \notag \\
=&   \sum_{j_m\geq 0} ~
 \psi_{j_m, q-1}^6  \underbrace{\sum_{\{i_m\colon i_m\geq  i_*(j_m)\}} \psi_{m,i_m,j_m, q}^6}_{\equiv 1\mbox{ by } \eqref{eq:psi:i:j:partition:0}}
=   \sum_{j_m\geq 0} ~
 \psi_{j_m,q-1}^6  \equiv 1 \, 
\end{align*}
where the last ineqaulity follows from the inductive assumption \eqref{eq:inductive:partition}.

The proof of the second claim is more involved and will be split into cases.  Using the definition in \eqref{eq:psi:m:im:q:def}, we have that 
\begin{align*}
    \psi_{m,i_m,q+\bn}\psi_{m,i'_m,q+\bn} = \sum\limits_{\{j_m:i_m\geq i_*(j_m)\}} \sum\limits_{\{j_m':i_m'\geq i_*(j_m')\}} \psi_{j_m,q+\bn-1}^6\psi_{j'_m,q+\bn-1}^6 \psi_{m,i_m,j_m,q+\bn}^6\psi_{m,i'_m,j'_m,q+\bn}^6 \, .
\end{align*}
Recalling the inductive assumption \eqref{eq:inductive:partition}, we have that the above sum only includes pairs of indices $j_m$ and $j'_m$ such that $|j_m-j'_m|\leq 1$. So we may assume that 
\begin{equation}\label{eq:overlap:assumption}
(x,t) \in \supp \psi_{m,i_m,j_m,q} \cap \supp \psi_{m,i'_m,j'_m,q},
\end{equation}
where $|j_m-j'_m|\leq 1$. The first and simplest case is the case $j_m=j'_m$. We then appeal to \eqref{eq:intermediate:overlapping} to deduce that it must be the case that $|i_m-i'_m|\leq 1$ in order for
\eqref{eq:overlap:assumption} to be true.

Before moving to the second and third cases, we recall from the proof of \cite[Lemma~6.7]{BMNV21} that by symmetry it will suffice to prove that $\psi_{m,i_m,q+\bn}\psi_{m,i'_m,q+\bn}\equiv 0$ when $i_m'\leq i_m-2$. We then consider the second case in \eqref{eq:overlap:assumption}, in which $j'_m=j_m+1$. When $i_m = i_*(j_m)$, we use that $i_*(j_m)\leq i_*(j_m+1)$ to obtain
\begin{equation*}
i'_m\leq i_m-2=i_*(j_m)-2 < i_*(j_m+1) = i_*(j'_m) \, ,
\end{equation*}
and so by Definition~\ref{def:intermediate:cutoffs}, we have that $\psi_{m,i'_m,j'_m,q+\bn}=0$. Thus we need only now consider $i_m > i_*(j_m)$ in order to finish the proof of the second case from \eqref{eq:overlap:assumption}. From \eqref{eq:overlap:assumption}, items~\eqref{item:cutoff:1}--\eqref{item:cutoff:2} from Lemma~\ref{lem:cutoff:construction:first:statement}, and Definition~\ref{def:intermediate:cutoffs}, we have that
\begin{subequations}
\begin{align}
    h_{m,j_m,q+\bn}(x,t) &\in \left[ \frac{1}{2} \Gamma_{q+\bn}^{(m+1)(i_m-i_*(j_m))}, \Gamma_{q+\bn}^{(m+1)(i_m+1-i_*(j_m))} \right], \label{eq:hm:overlap:0} \\
    h_{m,j_m+1,q+\bn}(x,t) &\leq  \Gamma_{q+\bn}^{(m+1)(i'_m+1-i_*(j_m+1))} . \label{eq:hm:overlap:1}
\end{align}
\end{subequations}
Note that from the definition of $h_{m,j_m,q+\bn}$ in \eqref{eq:h:j:q:def}, we have that 
\begin{equation*}
 \Gamma_{q+\bn}^{(m+1)(i_*\left(j_m+1\right)-i_*\left(j_m\right))}   h_{m,j_m+1,q+\bn} =  h_{m,j_m,q+\bn} \, .
\end{equation*}
Then, since $ i'_m\leq i_m-2$, from \eqref{eq:hm:overlap:1} we have that
\begin{align*}
 \Gamma_{q+\bn}^{-(m+1)(i_m-i_*(j_m))}   h_{m,j_m,q+\bn} 
 &=  \Gamma_{q+\bn}^{-(m+1)(i_m-i_*(j_m))}   h_{m,j_m+1,q+\bn}  \Gamma_{q+\bn}^{(m+1)(i_*\left(j_m+1\right)-i_*\left(j_m\right))}   \notag \\
    &\leq  \Gamma_{q+\bn}^{-(m+1)(i_m-i_*(j_m))}    \Gamma_{q+\bn}^{(m+1)(i'_m+1-i_*(j_m+1))}  \Gamma_{q+\bn}^{(m+1)(i_*\left(j_m+1\right)-i_*\left(j_m\right))}   \notag\\
    &=   \Gamma_{q+\bn}^{(m+1)(i'_m+1 -i_m )}  \notag\\
    &\leq \Gamma_{q+\bn}^{-(m+1)} \,.
\end{align*}
Since $m\geq 0$, the above estimate contradicts the lower bound on  $h_{m,j_m,q+\bn}$ in \eqref{eq:hm:overlap:0} because $\Gamma_{q+\bn}^{-1} \ll \sfrac 12$ 
for $a$ sufficiently large. 

We move to the third and final case, $j'_m=j_m-1$. As before, if $i_m=i_*(j_m)$, then since $i_*(j_m)\leq i_*(j_m-1)+1$, we have that 
\begin{equation*}
    i'_m \leq i_m - 2 = i_*(j_m)-2 \leq i_*(j_m-1) -1 < i_*(j_m-1) = i_*(j'_m)\,,
\end{equation*}
which by Definition~\ref{def:intermediate:cutoffs} implies that $\psi_{m,i'_m,j'_m,q+\bn}=0$, and there is nothing to prove. Thus, we only must consider the case $i_m > i_*(j_m)$. Using the definition \eqref{eq:h:j:q:def} we have that
\begin{equation*}
  h_{m,j_m,q+\bn} = \Gamma_{q+\bn}^{(m+1)(i_*(j_m-1) -i_*(j_m))}    h_{m,j_m-1,q+\bn}
    \,.  
\end{equation*}
On the other hand, for $i'_m \leq i_m-2$ we have from \eqref{eq:hm:overlap:1} that 
\begin{align*}
    h_{m,j_m-1,q+\bn} 
    \leq \Gamma_{q+\bn}^{(m+1)(i'_m+1-i_*(j_m-1))}
    \leq \Gamma_{q+\bn}^{(m+1)(i_m-1-i_*(j_m-1))} \,.
\end{align*}
Therefore, combining the above two displays and the inequality $-i_*(j_m)\geq -i_*(j_m-1)-1$, we obtain the bound
\begin{align*}
  \Gamma_{q+\bn}^{-(m+1)(i_m-i_*(j_m))}  h_{m,j_m,q+\bn} 
  &\leq  \Gamma_{q+\bn}^{-(m+1)(i_m-i_*(j_m))}   \Gamma_{q+\bn}^{(m+1)(i_*(j_m-1) -i_*(j_m))} \Gamma_{q+\bn}^{(m+1)(i_m-1-i_*(j_m-1))}
    \notag\\
      &=     \Gamma_{q+\bn}^{-(m+1)} \,,
\end{align*}
As before, since $m\geq 0$ this produces a contradiction with the lower bound on $h_{m,j_m,q+\bn}$ given in \eqref{eq:hm:overlap:0}, since $\Gamma_{q+\bn}^{-1} \ll \sfrac 12$.
\end{proof}


\begin{lemma}[\bf $\psi_{i,q+\bn}$ - Partition of unity]
\label{lem:partition:of:unity:psi}
We have that
\begin{align}
\sum_{i\geq 0} \psi_{i,q+\bn}^{6}\equiv 1 \,, \qquad \psi_{i,q+\bn}\psi_{i',q+\bn} \equiv 0\quad\textnormal{for}\quad|i-i'|\geq 2 \, . \label{eq:lemma:partition:2}
\end{align}
\end{lemma}

\begin{proof}[Proof of Lemma~\ref{lem:partition:of:unity:psi}]
To prove the first claim for $q+\bn\geq 1$, let us introduce the notation 
\begin{equation}\label{eq:Lambda:i:def}
    \Lambda_i = \left\{ \Vec{i} = (i_0,...,i_{\NcutSmall}) \colon \max_{0\leq m\leq\NcutSmall}i_m=i. \right\}
\end{equation}
Then
\begin{equation*}
    \psi_{i,q+\bn}^6 = \sum_{\Vec{i}\in\Lambda_i} \prod\limits_{m=0}^{\NcutSmall} \psi_{m,i_m,q+\bn}^6 \, ,
\end{equation*}
and thus
\begin{align}
    \sum_{i\geq 0} \psi_{i,q}^6 
    = \sum_{i\geq 0}\sum_{\Vec{i}\in\Lambda_i} \prod\limits_{m=0}^{\NcutSmall} \psi_{m,i_m,q}^6 
    &=\sum_{\Vec{i}\in\mathbb{N}_0^{\NcutSmall+1}} \left( \prod\limits_{m=0}^{\NcutSmall} \psi_{m,i_m,q}^6 \right) \notag\\
    &= \prod\limits_{m=0}^{\NcutSmall} \left( \sum_{i_m\geq 0} \psi_{m,i_m,q}^6 \right)
    = \prod\limits_{m=0}^{\NcutSmall} 1 
    =1 \notag
\end{align}
after using \eqref{eq:lemma:partition:1}.

To prove the second claim, assume towards a contradiction that there exists $|i-i'|\geq 2$ such that $\psi_{i,q}\psi_{i',q}\geq 0$.  Then 
\begin{align}
    0 &\neq \psi_{i,q+\bn}^6 \psi_{i',q+\bn}^6 = \sum_{\Vec{i}\in\Lambda_i} \sum_{\Vec{i}'\in\Lambda_{i'}} \prod\limits_{m=0}^{\NcutSmall} \psi_{m,i_m,q+\bn}^6 \psi_{m,i'_m,q+\bn}^6 \, . \label{eq:lemma:partition:nonzero}
\end{align}
In order for \eqref{eq:lemma:partition:nonzero} to be non-vanishing, by \eqref{eq:lemma:partition:1}, there must exist $\Vec{i}=(i_0,...,i_{\NcutSmall})\in\Lambda_i$ and $\Vec{i}'=(i'_0,...,i'_{\NcutSmall})\in\Lambda_{i'}$ such that $|i_m-i'_m|\leq 1$ for all $0\leq m \leq \NcutSmall$. By the definition of $i$ and $i'$, there exist $m_*$ and $m'_*$ such that
\begin{equation*}
    i_{m_*} = \max_m i_m = i, \qquad     i'_{m'_*} = \max_m i'_m = i'.
\end{equation*}
But then
\begin{align}
    i=i_{m_*} \leq i'_{m_*} + 1 \leq i'_{m'_*} + 1 = i'+1 \, , \qquad \qquad 
    i'=i'_{m'_*} \leq i_{m'_*} + 1 \leq i_{m_*} + 1 = i+1, \notag
\end{align}
implying that $|i-i'|\leq 1$, a contradiction.
\end{proof}


\begin{lemma}[\bf Lower order derivative bounds on $\hat w_{q+\bn}$]
\label{lem:h:j:q:size}
If $(x,t)\in \supp( \psi_{m,i_m,j_m,q+\bn})$ then
\begin{align}\label{eq:h:psi:supp}
h_{m,j_m,q+\bn}\leq  \Gamma_{q+\bn}^{(m+1)\left(i_m+1-i_*(j_m)\right)}.
\end{align}
Moreover, if $i_m>i_*(j_m)$ we have
\begin{align}\label{eq:psi:supp:upper}
h_{m,j_m,q+\bn} \geq  (\sfrac 12) \Gamma_{q+\bn}^{(m+1)(i_m-i_*(j_m))} 
\end{align}
on the support of $\psi_{m,i_m,j_m,q+\bn}$. 
As a consequence, we have that for all $0 \leq m, M \leq \NcutSmall$ and $0 \leq N \leq \NcutLarge$,
\begin{subequations}
\begin{align}
\norm{D^N D_{t,q+\bn-1}^{m} \hat w_{q+\bn}}_{L^\infty( \supp \psi_{m,i_m,q+\bn})}
&\leq {\delta_{q+\bn}^{\sfrac 12}r_{q}^{-\sfrac13}} \Gamma_{q+\bn}^{i_m+1} (\lambda_{q+\bn} \Gamma_{q+\bn})^N (\tau_{q+\bn-1}^{-1} \Gamma_{q+\bn}^{i_m+3})^{m} \label{eq:derivatives:psi:i:m:q} \\
\norm{D^N D_{t,q+\bn-1}^M \hat w_{q+\bn}}_{L^\infty( \supp \psi_{i,q+\bn})}
&\leq {\delta_{q+\bn}^{\sfrac 12}r_{q}^{-\sfrac13}} \Gamma_{q+\bn}^{i+1} (\lambda_{q+\bn} \Gamma_{q+\bn})^N (\tau_{q+\bn-1}^{-1} \Gamma_{q+\bn}^{i+3})^{M}  \, .
\label{eq:derivatives:psi:i:q}
\end{align}
\end{subequations}
\end{lemma}
\begin{proof}[Proof of Lemma~\ref{lem:h:j:q:size}]
Estimates \eqref{eq:h:psi:supp} and \eqref{eq:psi:supp:upper} follow directly from the definitions of $\tilde \gamma_{m,q+\bn}$ and $\gamma_{m,q+\bn}$ in Lemma~\ref{lem:cutoff:construction:first:statement} and the definition of $h_{m,j_m,q+\bn}$ in \eqref{eq:h:j:q:def}. In order to prove \eqref{eq:derivatives:psi:i:m:q}, we note that for $(x,t) \in \supp (\psi_{m,i_m,q+\bn})$, by \eqref{eq:psi:m:im:q:def} there must exist a $j_m$ with $i_*(j_m) \leq i_m$ such that $(x,t) \in \supp (\psi_{m,i_m,j_m,q+\bn})$. Using \eqref{eq:h:psi:supp}, we conclude that
\begin{align}
\norm{ D^N D_{t,q+\bn-1}^m \hat w_{q+\bn}  }_{L^\infty(\supp \psi_{m,i_m,j_m,q+\bn})} 
&\leq \delta_{q+\bn}^{\sfrac 12} r_q^{-\sfrac 13} \Gamma_{q+\bn}^{i_m+1} \left(\lambda_{q+\bn} \Gamma_{q+\bn}\right)^N \left( \tau_{q+\bn-1}^{-1}  \Gamma_{q+\bn}^{i_m+3} \right)^m  
\label{eq:derivatives:psi:i:j:q}
\end{align}
which completes the proof of \eqref{eq:derivatives:psi:i:m:q}. The proof of \eqref{eq:derivatives:psi:i:q} follows from the fact that we have employed the \textit{maximum} over $m$ of $i_m$ to define $\psi_{i,q+\bn}$ in \eqref{def:psi:i:q:def}.
\end{proof}


\begin{corollary}[\bf Higher order derivative bounds on $\hat w_{q+\bn}$]
\label{cor:D:Dt:wq:psi:i:q}
For $N+M \leq 2\Nfin$ and $i \geq 0$, we have the bound 
\begin{align}
&\norm{D^N D_{t,q+\bn-1}^M \hat w_{q+\bn}}_{L^\infty(\supp \psi_{i,q+\bn})} \notag\\
&\qquad \qquad \leq \Gamma_{q+\bn}^{i+1} {\delta_{q+\bn}^{\sfrac 12}r_{q}^{-\sfrac13}}  (\la_{q+\bn}\Ga_{q+\bn})^N \MM{M,\Nindt,\Gamma_{q+\bn}^{i+3} \tau_{q+\bn-1}^{-1},\Tau_{q+\bn-1}^{-1}\Gamma_{q+\bn-1}} 
 \,.
\label{eq:D:Dt:wq:psi:i:q}
\end{align}
\end{corollary}
\begin{proof}[Proof of Corollary~\ref{cor:D:Dt:wq:psi:i:q}]
When $0 \leq N \leq \NcutLarge$ and $0\leq M \leq \NcutSmall \leq \Nindt$, the desired bound was already established in \eqref{eq:derivatives:psi:i:q}. For the remaining cases in which either $N>\NcutLarge$ or $M> \NcutSmall$, note that if $0\leq m \leq \NcutSmall$ and $(x,t) \in \supp \psi_{m,i_m,q+\bn}$, there exists $j_m\geq 0$ with $i_*(j_m) \leq i_m$ such that $(x,t) \in \supp \psi_{j_m,q+\bn-1}$. Thus, we may appeal to \eqref{eq:vellie:inductive:dtq-1:uniform:upgraded:statement}, which gives that for $N+M \leq 2\Nfin$,
\begin{align*}
\abs{D^N D^M_{t,q+\bn-1} \hat w_{q+\bn}(x,t)} \les \Ga_q^{\sfrac{\badshaq}{2}+16} r_q^{-1} (\la_{q+\bn}\Ga_{q+\bn-1})^N \MM{M,\Nindt,\Gamma_{q+\bn-1}^{j_m-1} \tau_{q+\bn-1}^{-1},\Tau_{q+\bn-1}^{-1}\Ga_{q+\bn-1}} \,.
\end{align*}
Since $i_*(j_m) \leq i_m$ implies $\Gamma_{q+\bn-1}^{j_m} \leq \Gamma_{q+\bn}^{i_m}$, we deduce that for $N+M \leq 2\Nfin$, 
\begin{align*}
&\norm{D^N D^M_{t,q+\bn-1} \hat w_{q+\bn}}_{L^\infty(\supp \psi_{m,i_m,q+\bn})} \\
&\qquad \les\Ga_q^{\sfrac{\badshaq}{2}+16} r_q^{-1} (\la_{q+\bn}\Ga_{q+\bn-1})^N \MM{M,\Nindt,\Gamma_{q+\bn}^{i_m} \tau_{q+\bn-1}^{-1},\Tau_{q+\bn-1}^{-1}\Ga_{q+\bn-1}} \\
&\qquad \leq \Ga_{q+\bn}^{i_m+1} \delta_{q+\bn}^{\sfrac 12} r_q^{-\sfrac 13} (\la_{q+\bn}\Ga_{q+\bn})^N \MM{M,\Nindt,\Gamma_{q+\bn}^{i_m+3} \tau_{q+\bn-1}^{-1},\Tau_{q+\bn-1}^{-1}\Gamma_{q+\bn-1}} 
\end{align*}
after using that either $N>\NcutLarge$ or $M> \NcutSmall$, the parameter inequality \eqref{condi.Ncut0.2}, and a large choice of $a$ to absorb the implicit constant in the spare factor of $\Ga_{q+\bn}$. The desired estimate in \eqref{eq:D:Dt:wq:psi:i:q} then follows from taking the maximum over $m$ from Definition~\ref{def:psi:i:q:def}.
\end{proof}

\subsection{Pure spatial derivatives}
In this section we prove that the cutoff functions $\psi_{i,q+\bn}$ satisfy sharp spatial derivative estimates which are consistent with \eqref{eq:sharp:Dt:psi:i:q:old} for $q'=q+\bn$.
\begin{lemma}[\bf Spatial derivatives for the cutoffs]
\label{lem:sharp:D:psi:i:q}
Fix $q+\bn\geq 1$,  $0 \leq m \leq \NcutSmall$, and $i_m\geq0$. For all $j_m \geq 0$ such that $i_m \geq i_*(j_m)$, all $i\geq 0$, and all  $N \leq \Nfin$, we have 
\begin{subequations}
\begin{align}
{\bf 1}_{\supp( \psi_{j_m,q+\bn-1})} \frac{|D^N \psi_{m,i_m,j_m,q+\bn}|}{\psi_{m,i_m,j_m,q+\bn}^{1 -  N/\Nfin}} 
\les (\la_{q+\bn}\Ga_{q+\bn})^N
\,,
\label{eq:sharp:D:psi:i:j:q} \\
\frac{|D^N \psi_{i,q+\bn}|}{\psi_{i,q+\bn}^{1-N/\Nfin}} \les (\la_{q+\bn}\Ga_{q+\bn})^N \, .
\label{eq:sharp:D:psi:i:q}
\end{align}
\end{subequations}
\end{lemma}
\begin{proof}[Proof of Lemma~\ref{lem:sharp:D:psi:i:q}]
\textbf{Step 1: proof of \eqref{eq:sharp:D:psi:i:j:q}.}  We distinguish two cases.  The first case is when $\psi=\tilde\gamma_{m,q}$, or $\psi=\gamma_{m,q}$ \emph{and} we have the lower bound
\begin{equation}\label{eq:clusterf:1}
h_{m,j_m,q+\bn}^2 \Gamma_{q+\bn}^{-2\left(i_m-i_*(j_m)\right)(m+1)} \geq \frac{1}{4} \Gamma_{q+\bn}^{2(m+1)} \, ,
\end{equation}
so that \eqref{eq:DN:psi:q:gain} applies. The goal is then to apply \cite[Lemma~A.5]{BMNV21} to the function $\psi = \tilde \gamma_{m,q}$ or  $\psi = \gamma_{m,q}$ with the choices $\Gamma_\psi = \Gamma_{q+\bn}^{m+1}$, $\Gamma = \Gamma_{q+\bn}^{(m+1)(i_m-i_*(j_m))}$, and $h = h^2_{m,j_m,q+\bn}$. The assumption in \cite[equation~(A.24)]{BMNV21} holds by \eqref{eq:DN:psi:q:0} or \eqref{eq:DN:psi:q:gain} for all $N \leq \Nfin$, and so we need to obtain bounds on the derivatives of $h_{m,j_m,q+\bn}^2$ which are consistent with assumption in \cite[equation~(A.25)]{BMNV21} of \cite[Lemma~A.5]{BMNV21}. For $B\leq \Nfin$, the Leibniz rule gives
\begin{align}
&\left| D^{B} h_{m,j_m,q+\bn}^2\right| \notag\\
&\quad \lesssim (\lambda_{q+\bn} \Gamma_{q+\bn})^B \sum_{B' =0}^{B} \sum_{n=0}^{\NcutLarge}  \Gamma_{{q+\bn}}^{-i_*(j_m)} (\tau_{{q+\bn}-1}^{-1} \Gamma_{{q+\bn}}^{i_*(j_m)+2})^{-m} (\lambda_{q+\bn} \Gamma_{q+\bn})^{-n-B'}  \delta_{q+\bn}^{-\sfrac 12} r_q^{\sfrac 13} | D^{n + B'} D^m_{t,{q+\bn}-1} \hat w_{{q+\bn}}| 
\notag\\
&\quad \qquad \times \Gamma_{{q+\bn}}^{-i_*(j_m)} (\tau_{{q+\bn}-1}^{-1} \Gamma_{{q+\bn}}^{i_*(j_m)+2})^{-m} (\lambda_{q+\bn} \Gamma_{q+\bn})^{-n - B + B'} \delta_{{q+\bn}}^{-\sfrac 12} r_q^{\sfrac 13} | D^{n + B-B'} D^m_{t,{q+\bn}-1} \hat w_{{q+\bn}}| 
\,.
\label{eq:cutoff:spatial:derivatives:000}
\end{align}
For the terms with $L \in \{n+B',n+B-B'\} \leq \NcutLarge$, we may appeal to appeal to estimate \eqref{eq:h:psi:supp}, which gives
\begin{align}
&\Gamma_{{q+\bn}}^{-i_*(j_m)} (\tau_{{q+\bn}-1}^{-1} \Gamma_{{q+\bn}}^{i_*(j_m)+2})^{-m} (\lambda_{q+\bn} \Gamma_{q+\bn})^{-L} \delta_{{q+\bn}}^{-\sfrac 12} r_q^{\sfrac 13} \norm{D^L D_{t,{q+\bn}-1}^m \hat w_{q+\bn}}_{L^\infty(\supp \psi_{m,i_m,j_m,{q+\bn}})} \notag\\
&\qquad \qquad \qquad \leq  \Gamma_{{q+\bn}}^{(m+1)(i_m+1 - i_*(j_m))} \, .
\label{eq:cutoff:spatial:derivatives:00}
\end{align}
On the other hand, for $\NcutLarge < L \in \{n+B',n+B-B'\} \leq \NcutLarge + B \leq 2\Nfin - \Nindt$,  we may appeal to appeal to \eqref{eq:vellie:inductive:dtq-1:uniform:upgraded:statement}, and since $m\leq \NcutSmall < \Nindt$, we deduce that
\begin{align}
&\Gamma_{{q+\bn}}^{-i_*(j_m)} (\tau_{{q+\bn}-1}^{-1} \Gamma_{{q+\bn}}^{i_*(j_m)+2})^{-m} (\lambda_{q+\bn} \Gamma_{q+\bn})^{-L} \delta_{{q+\bn}}^{-\sfrac 12} r_q^{\sfrac 13} \norm{D^L D_{t,{q+\bn}-1}^m \hat w_{q+\bn}}_{L^\infty(\supp \psi_{j_m,{q+\bn}-1})} \notag \\
&\lesssim \Ga_{q+\bn}^{-i_*(j_m)(m+1)-2m} \tau_{q+\bn-1}^{m} (\la_{q+\bn}\Ga_{q+\bn})^{-L} \delta_{q+\bn}^{-\sfrac 12} r_q^{\sfrac 13} \Ga_q^{\sfrac{\badshaq}{2}+16} r_q^{-1} (\la_{q+\bn}\Ga_{q+\bn-1})^{L} (\tau_{q+\bn-1}^{-1}\Ga_{q+\bn-1}^{j_m-1})^m \notag\\
&\lesssim \Ga_{q+\bn}^{-i_*(j_m)(m+1)-2m} \delta_{q+\bn}^{-\sfrac 12} r_q^{\sfrac 13} \Ga_q^{\sfrac{\badshaq}{2}+16} r_q^{-1} \left( \frac{\Ga_{q+\bn-1}}{\Ga_{q+\bn}} \right)^L \Ga_{q+\bn-1}^{m(j_m-1)} \notag\\
&\leq \Ga_{q+\bn}^{(i_m+1-i_*(j_m))(m+1)}  \, .
\label{eq:cutoff:spatial:derivatives:0}
\end{align}
In the last inequality we have used that $i_m \geq i_*(j_m)$ in order to convert $\Ga_{q+\bn-1}^{m(j_m-1)}$ into $\Ga_{q+\bn}^{m i_m}$ and \eqref{condi.Ncut0.3}, which is applicable by the assumption that $L>\NcutLarge$.  Summarizing the bounds \eqref{eq:cutoff:spatial:derivatives:000}--\eqref{eq:cutoff:spatial:derivatives:0}, since $n\leq \NcutLarge$ and $\Nindt\leq \Nfin$, we arrive at
\begin{align*}
{\bf 1}_{ \supp(\psi_{j_m,q+\bn-1} \psi_{m,i_m,j_m,q+\bn}  )} \left| D^{B} h_{m,j_m,q+\bn}^2\right| & \lesssim (\lambda_{q+\bn} \Gamma_{q+\bn})^B \Gamma_{q+\bn}^{2(m+1)(i_m+1 - i_*(j_m))}
\end{align*}
whenever $B\leq \Nfin$. Thus, the assumption in \cite[A.25]{BMNV21} holds with $C_h = \Gamma_{q+\bn}^{2(m+1)(i_m+1 - i_*(j_m))}$, $\lambda = \tilde\lambda = \lambda_{q+\bn}\Ga_{q+\bn}$, $N_* = \infty$, $N=\Nfin$, $M=0$. Note that with these choices of parameters, we have $ C_h \Gamma_\psi^{-2} \Gamma^{-2} =1 $. We may thus apply \cite[Lemma~A.5]{BMNV21} and conclude that 
\begin{align*}
{\bf 1}_{\supp(\psi_{j_m,q+\bn-1})} \frac{\left|D^N \psi_{m,i_m,j_m,q+\bn} \right|}{\psi_{m,i_m,j_m,q+\bn}^{1-N/\Nfin}} \les (\la_{q+\bn}\Ga_{q+\bn})^N
\end{align*}
for all $N \leq \Nfin$, proving \eqref{eq:sharp:D:psi:i:j:q} in the first case.

Recalling the inequality \eqref{eq:clusterf:1}, the second case is when $\psi=\gamma_{m,q}$ and
\begin{equation}\label{eq:clusterf:2}
h_{m,j_m,q+\bn}^2 \Gamma_{q+\bn}^{-2\left(i_m-i_*(j_m)\right)(m+1)} \leq \frac{1}{4} \Gamma_{q+\bn}^{2(m+1)} \, .
\end{equation}
However, since $\gamma_{m,q}$ is uniformly equal to $1$ when the left hand side of the above display takes values in $\left[1,\frac{1}{4}\Gamma_{q}^{2(m+1)}\right]$ from item~\eqref{item:cutoff:2} in Lemma~\ref{lem:cutoff:construction:first:statement}, \eqref{eq:sharp:D:psi:i:j:q} is trivially satisfied in this range of values of the left-hand side.  Thus the analysis of the second case reduces to analyzing the subcase when
\begin{equation}\label{eq:clusterf:3}
h_{m,j_m,q+\bn}^2 \Gamma_{q+\bn}^{-2\left(i_m-i_*(j_m)\right)(m+1)} \leq 1 \, .
\end{equation}
As in the first case, we aim to apply \cite[Lemma~A.5]{BMNV21} with $h=h_{m,j_m,q}^2$, but now with $\Gamma_\psi=1$ and $\Gamma=\Gamma_{q+\bn}^{(m+1)(i_m-i_*(j_m))}$.  From \eqref{eq:DN:psi:q}, the assumption in \cite[(A.24)]{BMNV21} holds.  Towards estimating derivatives of $h$, for the terms with $L\in\{n+B',n+B-B'\}\leq\NcutLarge$, \eqref{eq:clusterf:3} gives immediately that
\begin{align}
&\Gamma_{q+\bn}^{-i_*(j_m)} (\tau_{{q+\bn}-1}^{-1} \Gamma_{q+\bn}^{i_*(j_m)+2})^{-m} (\lambda_{q+\bn} \Gamma_{q+\bn})^{-L} \delta_{q+\bn}^{-\sfrac 12} r_q^{\sfrac 13} \norm{D^L D_{t,{q+\bn}-1}^m \hat w_{q+\bn}}_{L^\infty(\supp \psi_{m,i_m,j_m,{q+\bn}})} \notag\\
&\qquad \qquad \leq  \Gamma_{q+\bn}^{(m+1)(i_m-i_*(j_m))} \, .
\label{eq:cutoff:spatial:derivatives:00:cf}
\end{align}
Conversely, when $\NcutLarge>L$, we may argue as in the estimates which gave \eqref{eq:cutoff:spatial:derivatives:0}, except we achieve the slightly improved bound of $\Ga_{q+\bn}^{(m+1)(i_m-i_*(j_m))}$ as above. We then arrive at
\begin{align*}
{\bf 1}_{ \supp(\psi_{j_m,{q+\bn}-1} \psi_{m,i_m,j_m,{q+\bn}}  )} \left| D^{B} h_{m,j_m,{q+\bn}}^2\right| & \lesssim \Gamma_{q+\bn}^{2(m+1)(i_m-i_*(j_m))} (\lambda_{q+\bn} \Gamma_{q+\bn})^B \notag
\end{align*}
whenever $B\leq \Nfin$. Thus, the assumption in \cite[(A.25)]{BMNV21} now holds with the same choices as before, except now $\const_h = \Gamma_{q+\bn}^{2(m+1)(i_m- i_*(j_m))}$, $\lambda =  \tilde\lambda = \lambda_{q+\bn}\Gamma_{q+\bn}$. Note that with these new choices of parameters, we still have $\const_h \Gamma_\psi^{-2} \Gamma^{-2} =1 $. We may thus apply \cite[Lemma~A.5]{BMNV21} and conclude that 
\begin{align*}
{\bf 1}_{\supp(\psi_{j_m,q+\bn-1})} \frac{\left|D^N \psi_{m,i_m,j_m,q+\bn} \right|}{\psi_{m,i_m,j_m,q+\bn}^{1-N/\Nfin}} \les (\la_{q+\bn}\Ga_{q+\bn})^N
\end{align*}
for all $N \leq \Nfin$, proving \eqref{eq:sharp:D:psi:i:j:q} in the second case.
\smallskip

\noindent\textbf{Step 2: differentiating $\psi_{m,i_m,q}$.}  From the definition \eqref{eq:psi:m:im:q:def} and the bound \eqref{eq:sharp:D:psi:i:j:q}, we next estimate derivatives of the $m^{\rm th}$ velocity cutoff function $\psi_{m,i_m,q}$ and claim that 
 \begin{align}
\frac{|D^N \psi_{m,i_m,q+\bn}|}{\psi_{m,i_m,q+\bn}^{1-N/\Nfin}} \lesssim (\la_{q+\bn}\Ga_{q+\bn})^N
\label{eq:sharp:D:psi:im:q}
\end{align}
for all $i_m\geq 0$ and all $N \leq \Nfin$.  We prove \eqref{eq:sharp:D:psi:im:q} by induction on $N$. When $N=0$ the bound trivially holds, which gives the induction base. For the induction step, assume that \eqref{eq:sharp:D:psi:im:q}  holds  for all $N' \leq N-1$.  By the Leibniz rule from Lemma~\ref{lem:leib} with $p=6$, we obtain
\begin{align}
D^N (\psi_{m,i_m,q+\bn}^6) = 6 \psi_{m,i_m,q+\bn}^5 D^N \psi_{m,i_m,q+\bn} + 
\sum_{\left\{ \substack{\alpha \, : \, \sum_{i=1}^6 \alpha_i = N \, , \\ \alpha_i < N \, \forall \, i }\right\}}
{\binom{N}{\alpha_1,\dots,\alpha_6}} 
\prod_{i=1}^6 D^{\alpha_i} \psi_{m,i_m,q+\bn}
\label{eq:psi:m:i:q:Leibniz}
\end{align}
and thus
\begin{align*}
\frac{D^N \psi_{m,i_m,q+\bn}}{\psi_{m,i_m,q+\bn}^{1-N/\Nfin}} 
&= \frac{D^N(\psi_{m,i_m,q+\bn}^6)}{6 \psi_{m,i_m,q+\bn}^{6-N/\Nfin}} - \frac{1}{6} \sum_{\left\{\substack{\alpha \, : \, \sum_{i=1}^p \alpha_i = N \, , \\ \alpha_i < N \, \forall \, i }\right\}} {\binom{N}{\alpha_1,\dots,\alpha_6}} \prod_{i=1}^6 \frac{D^{\alpha_i} \psi_{m,i_m,q+\bn}}{\psi_{m,i_m,q+\bn}^{1-\alpha_i/\Nfin}} \, .
\end{align*}
Since $\al_i\leq N-1$, by the induction assumption \eqref{eq:sharp:D:psi:im:q} we obtain
\begin{align}
\frac{\left| D^N \psi_{m,i_m,q+\bn} \right|}{\psi_{m,i_m,q+\bn}^{1-N/\Nfin}} 
&\lesssim \frac{ |D^N(\psi_{m,i_m,q+\bn}^6) |}{\psi_{m,i_m,q+\bn}^{6-N/\Nfin}} + (\la_{q+\bn}\Ga_{q+\bn})^N \, .
\label{eq:cutoff:spatial:derivatives:1}
\end{align}
Thus establishing \eqref{eq:sharp:D:psi:im:q}  for the $N$th derivative reduces to bounding the first term on the right side of the above. For this purpose we recall \eqref{eq:psi:m:im:q:def} and \eqref{eq:leibniz:one} and compute 
\begin{align*}
\frac{\left|D^N(\psi_{m,i_m,q+\bn}^6)\right|}{\psi_{m,i_m,q+\bn}^{6-N/\Nfin}}  
&= \frac{1}{\psi_{m,i_m,q+\bn}^{6-N/\Nfin}} \sum_{\{ j_m \colon i_*(j_m) \leq i_m\}} \sum_{K=0}^N {\binom{N}{K}} D^K(\psi_{j_m,q+\bn-1}^6) D^{N-K}(\psi_{m,i_m,j_m,q+\bn}^6) 
\\
&= \frac{\psi_{j_m,q+\bn-1}^{6-K/\Nfin} \psi_{m,i_m,j_m,q+\bn}^{6-(N-K)/\Nfin} }{\psi_{m,i_m,q+\bn}^{6-N/\Nfin}} \sum_{\{ j_m \colon i_*(j_m) \leq i_m\}} \sum_{K=0}^N {\binom{N}{K}} \notag\\
&\qquad \times \sum_{\alpha:\sum_{i=1}^6 \alpha_i=K} {\binom{K}{\alpha_1,\dots,\alpha_6}} \prod_{i=1}^6 \frac{D^{\alpha_i} \psi_{j_m,q+\bn-1}}{\psi_{j_m,q+\bn-1}^{1-\alpha_i/\Nfin}} \notag\\
&\qquad\qquad \times \sum_{\beta:\sum_{i=1}^6 \beta_i=N-K} {\binom{N-K}{\beta_1,\dots,\beta_6}} \prod_{i=1}^6 \frac{D^{\beta_i} \psi_{m,i_m,j_m,q+\bn}}{\psi_{m,i_m,j_m,q+\bn}^{1-\beta_i/\Nfin}} \, .
\end{align*}
Since $K, N-K \leq N$, and $\psi_{j_m,q+\bn-1}, \psi_{m,i_m,j_m,q+\bn} \leq 1$, we have by \eqref{eq:psi:m:im:q:def} that 
\begin{align*}
\frac{\psi_{j_m,q+\bn-1}^{6-K/\Nfin} \psi_{m,i_m,j_m,q+\bn}^{6-(N-K)/\Nfin} }{\psi_{m,i_m,q+\bn}^{6-N/\Nfin}}   \leq \frac{\psi_{j_m,q+\bn-1}^{6-N/\Nfin} \psi_{m,i_m,j_m,q+\bn}^{6-N/\Nfin} }{\psi_{m,i_m,q+\bn}^{6-N/\Nfin}}  \leq 1 \, .
\end{align*}
Then the estimate \eqref{eq:sharp:D:psi:i:j:q} and the inductive assumption \eqref{eq:sharp:Dt:psi:i:q:old} conclude the proof of \eqref{eq:sharp:D:psi:im:q}. In particular, note that this bound is independent of the value of $i_m$. 
\smallskip

\noindent\textbf{Step 3: proof of \eqref{eq:sharp:D:psi:i:q}}
In order to conclude the proof of the Lemma, we must argue that \eqref{eq:sharp:D:psi:im:q} implies \eqref{eq:sharp:D:psi:i:q}. Recalling \eqref{eq:psi:i:q:recursive}, we have that $\psi_{i,q+\bn}^6$ is given as a sum of products of  $\psi_{m,i_m,q+\bn}^6$, for which suitable derivative bounds are available due to \eqref{eq:sharp:D:psi:im:q}. Thus, the proof of \eqref{eq:sharp:D:psi:i:q} is again done by induction on $N$, mutatis mutandi to the proof of \eqref{eq:sharp:D:psi:im:q}. Indeed, we note that $\psi_{m,i_m,q+\bn}^6$ was also given as a sum of squares of cutoff functions for which derivative bounds were available. The proof of the induction step is thus again based on the application of the Leibniz rule for $\psi_{i,q+\bn}^6$; in order to avoid redundancy we omit these details.
\end{proof}

\subsection{Maximal index appearing in the cutoff}

\begin{lemma}[\bf Maximal $i$ index in the definition of $\psi_{i,q+\bn}$]
\label{lem:maximal:i}
There exists $\imax = \imax(q+\bn) \geq 0$, determined by \eqref{eq:imax:def} below, such that if $\lambda_0$ is sufficiently large, then\index{$\imax$}
\begin{subequations}
\begin{align}
\psi_{i,q+\bn} &\equiv 0 \quad \mbox{for all} \quad i > i_{\rm max} \, ,
\label{eq:imax} \\
\Gamma_{q+\bn}^{i_{\rm max}} &\leq  \Ga_{q}^{\sfrac{\badshaq}{2}+18} \delta_{q+\bn}^{-\sfrac 12} r_{q}^{-\sfrac 23} \, ,
\label{eq:imax:bound}\\
\imax(q)  &\leq \frac{\badshaq+12}{(b-1)\varepsilon_\Gamma}
\,.
\label{eq:imax:upper:bound:uniform}
\end{align}
\end{subequations}
\end{lemma}

\begin{proof}[Proof of Lemma~\ref{lem:maximal:i}]
Assume $i\geq 0$ is such that $\supp(\psi_{i,q+\bn}) \neq  \emptyset$. We will prove that 
\begin{equation}\label{eq:imax:proofy}
 \Gamma_{q+\bn}^i \leq \Ga_{q}^{\sfrac{\badshaq}{2}+18} \delta_{q+\bn}^{-\sfrac 12} r_q^{-\sfrac 23} \, . 
\end{equation}
From \eqref{eq:psi:i:q:recursive} it follows that for any $(x,t) \in \supp (\psi_{i,{q+\bn}})$, there must exist at least one $\Vec{i} = (i_0,\ldots,i_{\NcutSmall})$ such that $\max\limits_{0\leq m\leq\NcutSmall} i_m =i$ and $\psi_{m,i_m,{q+\bn}}(x,t) \neq 0$ for all $0\leq m \leq \NcutSmall$. Therefore, in light of \eqref{eq:psi:m:im:q:def}, for each such $m$ there exists a maximal $j_m$ such that $i_*(j_m) \leq i_m$, with $(x,t) \in \supp(\psi_{j_m,{q+\bn}-1}) \cap \supp(\psi_{m,i_m,j_m,{q+\bn}})$. In particular, this holds for any of the indices $m$ such that $i_m = i$. For the remainder of the proof, we fix such an index $0\leq m \leq \NcutSmall$.

If we have $i = i_m = i_{*}(j_m) = i_*(j_m,q)$, then using that $(x,t) \in \supp(\psi_{j_m,{q+\bn}-1})$ and the inductive assumption \eqref{eq:imax:old}, we have that $j_m \leq \imax({q+\bn}-1)$. Now using \eqref{eq:imax:old}, \eqref{ineq.eps.gamma.imax}, and the inequalities $\Gamma_{q+\bn}^{i-1} < \Gamma_{{q+\bn}-1}^{j_m} \leq \Gamma_{{q+\bn}-1}^{\imax({q+\bn}-1)}$, we deduce that
\[
\Gamma_{q+\bn}^{i}
\leq \Gamma_{q+\bn} \Gamma_{{q+\bn}-1}^{\imax({q+\bn}-1)}
\leq \Gamma_{q+\bn}
\Ga_{q-1}^{\sfrac{\badshaq}{2}+18}
 \delta_{q+\bn-1}^{-\sfrac 12} r_{q-1}^{-\sfrac 23} \leq \Ga_{q}^{\sfrac{\badshaq}{2}+18} \delta_{q+\bn}^{-\sfrac 12} r_q^{-\sfrac 23}
\,,
\]
Thus, in this case \eqref{eq:imax:proofy} holds.

On the other hand, if $i = i_m \geq i_{*}(j_m) +1$, then from \eqref{eq:psi:supp:upper} we have that 
$$|h_{m,j_m,q+\bn}(x,t)| \geq (\sfrac{1}{2}) \Gamma_{q+\bn}^{(m+1)(i_m - i_*(j_m))} \, . $$
Now from the pigeonhole principle, there exists $0 \leq n \leq \NcutLarge$ such that 
\begin{align*}
|D^n D_{t,q+\bn-1}^m \hat w_{q+\bn}(x,t)| &\geq \frac{1}{2 \NcutLarge} \Gamma_{q+\bn}^{(m+1)(i_m - i_*(j_m))} {\Gamma_{q+\bn}^{i_*(j_m)}} \delta_{q+\bn}^{\sfrac 12} r_q^{-\sfrac 13} (\lambda_{q+\bn} \Gamma_{q+\bn})^{n} (\tau_{{q+\bn}-1}^{-1} \Gamma_{q+\bn}^{i_*(j_m)+2})^m \notag\\
&\geq   \frac{1}{2 \NcutLarge} \Gamma_{q+\bn}^{i_m} \delta_{q+\bn}^{\sfrac 12} r_q^{-\sfrac 13} (\lambda_{q+\bn} \Ga_{q+\bn})^{n} (\tau_{{q+\bn}-1}^{-1} \Gamma_{q+\bn}^{i_m+2})^m \, ,
\end{align*}
and we also know that $(x,t) \in \supp(\psi_{j_m,q+\bn-1})$. By \eqref{eq:vellie:inductive:dtq-1:uniform:upgraded:statement} and the inequality $\NcutSmall\leq \NindSmall$ from \eqref{condi.Nindt}, we know that
\begin{align*}
|D^n D_{t,q+\bn-1}^m \hat w_{q+\bn}(x,t)|
&\leq \Gamma_{q}^{\sfrac{\badshaq}{2} + 17} r_q^{-1} (\lambda_{q+\bn}\Ga_{q+\bn-1})^{n} (\tau_{{q+\bn}-1}^{-1} \Gamma_{q+\bn-1}^{j_m-1})^m \\
&\leq \Gamma_{q}^{\sfrac{\badshaq}{2} + 17} r_q^{-1} (\lambda_{q+\bn}\Ga_{q+\bn})^{n} (\tau_{{q+\bn}-1}^{-1} \Gamma_{q+\bn}^{i_m})^m \, ,
\end{align*}
where in the last inequality we used the assumption that $i_m \geq i_*(j_m)$ and converted the $\Ga_{q+\bn-1}^{j_m-1}$ into $\Ga_{q+\bn}^{i_m}$. The proof is now completed, since the previous two inequalities and $i_m = i$ imply that 
\begin{equation}\label{eq:imax:contradiction:1}
 \Gamma_{q+\bn}^i \leq 2\NcutLarge \delta_{q+\bn}^{-\sfrac 12} r_q^{-\sfrac 23} \Ga_{q}^{\sfrac{\badshaq}{2}+17} \leq \delta_{q+\bn}^{-\sfrac 12} r_q^{-\sfrac 23} \Ga_{q}^{\sfrac{\badshaq}{2}+18} \, ,
\end{equation}
where in the last inequality we used \eqref{eq:badshaq:is:bad} and a large choice of $a$ to ensure that $\Ga_0 \geq 2\NcutLarge$.

In view of the above inequality, the value of $i_{\rm max}$ is chosen as
\begin{align}
\imax(q) = \sup \{i' \, : \, \Gamma_{q+\bn}^{i'} \leq \Ga_{q}^{\sfrac{\badshaq}{2}+18} r_q^{-\sfrac 23} \delta_{q+\bn}^{-\sfrac 12} \} \, .
\label{eq:imax:def}
\end{align}
With this definition, if $i > \imax(q+\bn)$, then $\supp(\psi_{i,q+\bn})= \emptyset$. To show that $\imax(q+\bn)$ is bounded  independently of $q$, simple (and brutal) computations give that
\begin{align*}
\frac{\log(  \Ga_{q}^{\sfrac{\badshaq}{2}+18}\delta_{q+\bn}^{-\sfrac 12}r_q^{-\sfrac 23})}{\log (\Gamma_{q+\bn})} \leq \frac{\badshaq+12}{(b-1)\varepsilon_\Gamma} \, ,
\end{align*}
verifying that \eqref{eq:imax:upper:bound:uniform} holds.
\end{proof}

\subsection{Mixed derivative estimates}
\newcommand{\Dqp}{D_{q+\bn}}
We will use the notation $D_{q+\bn} = \hat w_{q+\bn} \cdot \nabla$ for the directional derivative in the direction of $\hat w_{q+\bn}$. With this notation we have $D_{t,q+\bn} = D_{t,q+\bn-1} + D_{q+\bn}$. Next, we recall from \cite[equations~(6.54)-(6.55)]{BMNV21} that
\begin{align}
D_{q+\bn}^K = \sum_{j=1}^{K} f_{j,K} D^j \, ,
\label{eq:D:q:K:i}
\end{align}
where 
\begin{align}
f_{j,K} = \sum_{\{ \gamma \in \N^K \colon |\gamma|=K-j \} } c_{j,K,\gamma} \prod_{\ell=1}^K D^{\gamma_{\ell}} \hat w_{q+\bn} \, .
\label{eq:D:q:K:ii}
\end{align}
The  $c_{j,K,\gamma}$'s are explicitly computable coefficients that depend only on $K,j$, and $\gamma$. With the notation in \eqref{eq:D:q:K:ii} we have the following bounds.

\begin{lemma}[\bf Bounds for $D_{q+\bn}^K$]
\label{lem:D:a:fj}
 For $q+\bn\geq 1$ and $1 \leq K \leq  2\Nfin$, the functions $\{ f_{j,K}\}_{j=1}^{K}$ defined in \eqref{eq:D:q:K:ii} obey the estimate
\begin{align}
\norm{D^a f_{j,K}}_{L^\infty (\supp \psi_{i,q+\bn})} &\les  (\Gamma_{q+\bn}^{i+1} \delta_{q+\bn}^{\sfrac 12}r_{q}^{-\sfrac 13})^K (\la_{q+\bn}\Ga_{q+\bn})^{a+K-j}
\label{eq:D:a:fj} 
\end{align}
for any $a \leq 2 \Nfin-K+j$, and any $0\leq i\leq \imax(q+\bn)$.
\end{lemma}
\begin{proof}[Proof of Lemma~\ref{lem:D:a:fj}]
Note that no material derivative appears in  \eqref{eq:D:q:K:ii}, and thus to establish \eqref{eq:D:a:fj} we appeal to  Corollary~\ref{cor:D:Dt:wq:psi:i:q} with $M=0$ and \eqref{eq:vellie:inductive:dtq-1:uniform:upgraded:statement}. From the product rule we obtain that
\begin{align*}
\norm{D^a f_j}_{L^\infty(\supp \psi_{i,q+\bn})} 
&\lesssim \sum_{\{ \gamma \in \N^K \colon |\gamma|=K-j \}} \sum_{\{ \alpha \in \N^k \colon |\alpha| = a\} } \prod_{\ell=1}^K \norm{D^{\alpha_\ell + \gamma_{\ell}} \hat w_{q+\bn}}_{L^\infty(\supp \psi_{i,q+\bn})}  \notag\\
&\lesssim \sum_{\{ \gamma \in \N^K \colon |\gamma|=K-j \}} \sum_{\{ \alpha \in \N^k \colon |\alpha| = a\} } \prod_{\ell=1}^K 
\Gamma_{q+\bn}^{i+1} \delta_{q+\bn}^{\sfrac 12}r_q^{-\sfrac 13} (\la_{q+\bn}\Ga_{\qbn})^{\alpha_\ell+\gamma_\ell}
 \notag\\
&\lesssim (\Gamma_{\qbn}^{i+1} \delta_\qbn^{\sfrac 12}r_q^{-\sfrac 13})^K (\la_{q+\bn}\Ga_{\qbn})^{a+K-j}
\end{align*} 
since $|\gamma|= K-j$.
\end{proof}

\begin{lemma}[\bf Mixed derivatives for $\hat w_\qbn$]
\label{lem:Dt:Dt:wq:psi:i:q}
For $q+\bn\geq 1$ and $0 \leq i \leq \imax$, we have that 
\begin{align*}
&\norm{D^N D_\qbn^K D_{t,\qbn-1}^M \hat w_\qbn}_{L^\infty(\supp \psi_{i,\qbn})} \notag\\
&\qquad \qquad \les  
(\Gamma_{\qbn}^{i+1} \delta_\qbn^{\sfrac 12}r_{q}^{-\sfrac13})^{K+1} 
(\la_{q+\bn}\Ga_\qbn)^{N+K}
\MM{M,\NindSmall,\Gamma_{\qbn}^{i+3}  \tau_{\qbn-1}^{-1},\Tau_{\qbn-1}^{-1}\Ga_{\qbn-1}}\notag\\
&\qquad \qquad \les  
(\Gamma_{\qbn}^{i+1} \delta_\qbn^{\sfrac 12}r_{q}^{-\sfrac13}) 
(\la_\qbn\Ga_\qbn)^N (\Gamma_\qbn^{i-5} \tau_{\qbn}^{-1} )^{K}  \MM{M,\NindSmall,\Gamma_\qbn^{i+3}   \tau_{\qbn-1}^{-1},\Tau_{\qbn-1}^{-1}\Ga_{\qbn-1}}
\end{align*}
holds for $0 \leq K + N  + M \leq 2\Nfin$.
\end{lemma}

\begin{proof}[Proof of Lemma~\ref{lem:Dt:Dt:wq:psi:i:q}]
The second estimate in the Lemma follows from the parameter inequality \eqref{ineq:tau:q}. In order to prove the first estimate, we let $0 \leq a \leq N$ and $1 \leq j \leq K$.
From estimate \eqref{eq:D:Dt:wq:psi:i:q}, we obtain that
\begin{align*}
\norm{D^{N-a+j} D_{t,\qbn-1}^M \hat w_\qbn}_{L^\infty (\supp \psi_{i,\qbn})} 
&\les  \Gamma_\qbn^{i+1} \delta_\qbn^{\sfrac 12} r_q^{-\sfrac 13} (\la_\qbn \Ga_\qbn)^{N-a+j} \\
&\qquad \qquad \times \MM{M,\NindSmall,\Gamma_{\qbn}^{i+3} \tau_{\qbn-1}^{-1},\Tau_{\qbn-1}^{-1}\Ga_{\qbn-1}}
\end{align*}
for $N-a+j+M\leq \Nfin$, which may be combined with \eqref{eq:D:q:K:i}--\eqref{eq:D:a:fj} to obtain that
\begin{align*}
&\norm{D^N D_\qbn^K D_{t,\qbn-1}^M \hat w_\qbn}_{L^\infty(\supp \psi_{i,\qbn})}  \notag\\
&\quad \lesssim \sum_{a = 0}^N \sum_{j=1}^K   \norm{D^a f_{j,K}}_{L^\infty (\supp \psi_{i,\qbn})} \norm{D^{N-a+j} D_{t,\qbn-1}^M \hat w_\qbn}_{L^\infty (\supp \psi_{i,\qbn})} \notag\\
&\quad \lesssim (\Gamma_\qbn^{i+1} \delta_\qbn^{\sfrac 12}r_q^{-\sfrac 13})^{K+1}  (\la_\qbn\Ga_\qbn)^{N+K} \MM{M,\NindSmall,\Gamma_{\qbn}^{i+3} \tau_{\qbn-1}^{-1},\Tau_{\qbn-1}^{-1}\Ga_{\qbn-1}}
\end{align*}
holds for $N+M+K\leq 2\Nfin$, concluding the proof of the lemma. 
\end{proof}

\begin{lemma}[\bf More mixed derivatives for $\hat w_\qbn$ and derivatives for $\hat u_\qbn$]
\label{lem:Dt:Dt:wq:psi:i:q:multi}
For $q+\bn\geq 1$, $k\geq 1$, $\alpha,\beta \in {\mathbb N}^k$ with $|\alpha| = K$, $|\beta| = M$, and $K+M\leq \sfrac{3\Nfin}{2}+1$, we have
\begin{align}
&\norm{ \Big( \prod_{i=1}^k D^{\alpha_i} D_{t,\qbn-1}^{\beta_i} \Big) \hat w_\qbn }_{L^\infty(\supp \psi_{i,\qbn})} \notag\\
& \qquad  \lesssim \Gamma_\qbn^{i+1} \delta_\qbn^{\sfrac 12}r_{q}^{-\sfrac13} (\la_\qbn\Ga_\qbn)^K \MM{M,\NindSmall,\Gamma_\qbn^{i+3}  \tau_{\qbn-1}^{-1}, \Gamma_{\qbn-1} \Tau_{\qbn-1}^{-1}}   \, . \label{eq:nasty:D:wq}
\end{align}
Next, we have that
\begin{subequations}
\begin{align}
&\norm{D^N \Big( \prod_{i=1}^k D_\qbn^{\alpha_i} D_{t,\qbn-1}^{\beta_i} \Big) \hat w_\qbn }_{L^\infty(\supp \psi_{i,\qbn})}\notag\\
&\qquad\lesssim   (\Gamma_\qbn^{i+1} \delta_\qbn^{\sfrac 12}r_{q}^{-\sfrac13})^{K+1} (\la_\qbn\Ga_\qbn)^{N+K} \MM{M,\NindSmall,\Gamma_\qbn^{i+3}   \tau_{\qbn-1}^{-1},\Gamma_{\qbn-1} \Tau_{\qbn-1}^{-1}}  \label{eq:nasty:Dt:wq} \\
&\qquad  
\lesssim  \Gamma_\qbn^{i+1} \delta_\qbn^{\sfrac 12}r_{q}^{-\sfrac13} (\la_\qbn\Ga_\qbn)^N (\Gamma_\qbn^{i-5}  \tau_\qbn^{-1})^{K} \MM{M,\NindSmall,\Gamma_\qbn^{i+3}   \tau_{\qbn-1}^{-1},\Gamma_{\qbn-1} \Tau_{\qbn-1}^{-1}} 
\label{eq:nasty:Dt:wq:WEAK}
\end{align}
\end{subequations}
holds for all  $0 \leq K + M + N \leq \sfrac{3\Nfin}{2}+1$. Lastly, we have the estimate 
\begin{align}
& \norm{ \Big( \prod_{i=1}^k D^{\alpha_i} D_{t,\qbn}^{\beta_i} \Big) D \hat u_\qbn }_{L^\infty(\supp \psi_{i,\qbn})} \notag\\
&\qquad \lesssim \tau_\qbn^{-1}\Ga_\qbn^{i-5} (\la_{\qbn}\Ga_\qbn)^K
\MM{M,\NindSmall,\Gamma_\qbn^{i-5}   \tau_{\qbn}^{-1},\Gamma_{\qbn-1} \Tau_{\qbn-1}^{-1}} 
\label{eq:nasty:D:vq}
\end{align}
for all  $ K + M \leq \sfrac{3\Nfin}{2}$, the estimate 
\begin{align}
& \norm{ \Big( \prod_{i=1}^k D^{\alpha_i} D_{t,\qbn}^{\beta_i} \Big)  \hat u_\qbn }_{L^\infty(\supp \psi_{i,\qbn})} \notag\\
&\qquad \lesssim  \Gamma_\qbn^{i+1} \delta_\qbn^{\sfrac 12}r_{q}^{-\sfrac13}\lambda_\qbn^{2} (\la_\qbn\Ga_\qbn)^K \MM{M,\NindSmall,\Gamma_\qbn^{i-5}   \tau_{\qbn}^{-1},\Gamma_{\qbn-1} \Tau_{\qbn-1}^{-1}} 
\label{eq:nasty:no:D:vq}
\end{align}
for all  $ K + M \leq \sfrac{3\Nfin}{2}+1$, and the estimate
\begin{align}
    \left\| D^K \partial_t^M \hat u_\qbn \right\|_\infty \leq \la_\qbn^{\sfrac 12} (\la_\qbn\Ga_\qbn)^K \Tau_{\qbn}^{-M} \label{eq:bobby:new}
\end{align}
for all $K+M\leq 2\Nfin$.
\end{lemma}


\begin{proof}[Proof of Lemma~\ref{lem:Dt:Dt:wq:psi:i:q:multi}]

We note that \eqref{eq:nasty:Dt:wq:WEAK} follows directly from \eqref{eq:nasty:Dt:wq} by appealing to \eqref{ineq:tau:q}. We first  show that \eqref{eq:nasty:D:wq} holds, then establish \eqref{eq:nasty:Dt:wq}, and lastly, prove the bounds \eqref{eq:nasty:D:vq}--\eqref{eq:bobby:new}.
\smallskip

\noindent{\bf Proof of \eqref{eq:nasty:D:wq}.\,}
The statement is proven by induction on $k$. For $k=1$ the estimate holds for $K+M\leq 2\Nfin$ from  Corollary~\ref{cor:D:Dt:wq:psi:i:q}. For the induction step, assume that \eqref{eq:nasty:D:wq} holds  for any $k' \leq k-1$. We denote
\begin{align}
P_{k'} = \Big( \prod_{i=1}^{k'} D^{\alpha_i} D_{t,\qbn-1}^{\beta_i}\Big) \hat w_\qbn
\label{eq:P:k':def}
\end{align}
and write
\begin{align}
\Big( \prod_{i=1}^k D^{\alpha_i} D_{t,\qbn-1}^{\beta_i} \Big) \hat w_\qbn 
&= (D^{\alpha_k} D_{t,\qbn-1}^{\beta_k}) (D^{\alpha_{k-1}} D_{t,\qbn-1}^{\beta_{k-1}}) P_{k-2} \notag\\
&= (D^{\alpha_k+\alpha_{k-1}} D_{t,\qbn-1}^{\beta_k+\beta_{k-1}}) P_{k-2} + D^{\alpha_k} \left[D_{t,\qbn-1}^{\beta_k}, D^{\alpha_{k-1}}\right] D_{t,\qbn-1}^{\beta_{k-1}}P_{k-2} \, .
\label{eq:Merlot:*}
\end{align}
The first term in \eqref{eq:Merlot:*} already obeys the correct bound, since we know that \eqref{eq:nasty:D:wq}  holds for $k' = k-1$. In order to treat the second term on the right side of \eqref{eq:Merlot:*}, we use \cite[Lemma~A.12]{BMNV21} to write the commutator as\footnote{Following~\cite[subsection~A.7]{BMNV21}, we are using the following notation for iterated commutators.  First, $(\ad D_{t})^0(D) = D$ denotes a spatial derivative, i.e. a zeroth order commutator of $D_t$ and $D$.  Then for $k\geq 1$, we inductively set $(\ad D_{t})^k(D) = [D_t, (\ad D_t)^{k-1}(D)]$.}
\begin{align}
&D^{\alpha_k} \left[D_{t,\qbn-1}^{\beta_k}, D^{\alpha_{k-1}}\right] D_{t,\qbn-1}^{\beta_{k-1}}P_{k-2} \notag\\
&= D^{\alpha_k} \sum_{1 \leq |\gamma| \leq \beta_k} \frac{\beta_k!}{\gamma! (\beta_k - |\gamma|)!} \left(\prod_{\ell=1}^{\alpha_{k-1}} (\ad D_{t,\qbn-1})^{\gamma_\ell}(D) \right) D_{t,\qbn-1}^{\beta_k+\beta_{k-1}-|\gamma|}P_{k-2} \, .
 \label{eq:product:of:ad:Dt:q-1:comm}
\end{align}
From \cite[Lemma~A.13]{BMNV21} and the Leibniz rule we claim that one may expand
\begin{align}
\prod_{\ell=1}^{\alpha_{k-1}} (\ad D_{t,\qbn-1})^{\gamma_\ell}(D) = \sum_{j=1}^{\alpha_{k-1}} g_j D^{j}
\label{eq:product:of:ad:Dt:q-1}
\end{align}
for some explicit functions $g_j$ which obey the estimate
\begin{align}
\norm{ D^a g_j }_{L^\infty(\supp \psi_{i,q})} \lesssim   (\la_{\qbn-1}\Ga_{\qbn-1})^{ a+\alpha_{k-1}-j} 
\MM{|\gamma|,\Nindt,\Gamma_{\qbn}^{i+1}\tau_{\qbn-1}^{-1}, \Gamma_{\qbn-1}^{-1}\Tau_{\qbn-1}^{-1}}
\label{eq:product:of:ad:Dt:q-1:bnd}
\end{align}
for all $a$ such that $a + \alpha_{k-1}- j +|\gamma|\leq \sfrac{3\Nfin}{2}$. The claim  \eqref{eq:product:of:ad:Dt:q-1:bnd} requires a proof, which we sketch next.  Using the definition~\eqref{eq:psi:m:im:q:def} and the inductive estimate \eqref{eq:nasty:D:vq:old} at level $q'=\qbn-1$ and with $k=1$, we have that
\begin{align*}
&\norm{D^a D_{t,\qbn-1}^b D \hat u_{{\qbn-1}}}_{L^\infty(\supp \psi_{m,i_m,\qbn})}  
\notag\\
&\les \sum_{\{ j_m \colon \Gamma_{\qbn-1}^{j_m} \leq \Gamma_\qbn^{i_m}\}} \norm{D^a D_{t,\qbn-1}^b D \hat u_{{\qbn-1}}}_{L^\infty(\supp \psi_{j_m,\qbn-1})} \notag\\
&\les \sum_{\{ j_m \colon \Gamma_{\qbn-1}^{j_m} \leq \Gamma_\qbn^{i_m}\}} 
\tau_{\qbn-1}^{-1}\Ga_{\qbn-1}^{j_m+1} (\la_{\qbn-1}\Ga_{\qbn-1})^a
\MM{b,\Nindt,\Gamma_{\qbn-1}^{j_m+1}\tau_{\qbn-1}^{-1}, \Gamma_{\qbn-1}^{-1}\Tau_{\qbn-1}^{-1}} \notag\\
&\lesssim (\la_{\qbn-1}\Ga_{\qbn-1})^a \MM{b+1,\Nindt,\Gamma_{\qbn}^{i_m+1}\tau_{\qbn-1}^{-1}, \Gamma_{\qbn-1}^{-1}\Tau_{\qbn-1}^{-1}}
\end{align*}
for any $0\leq m \leq \NcutSmall$ and for all $a + b \leq \sfrac{3\Nfin}{2}$. Thus, from the definition \eqref{eq:psi:i:q:recursive} we deduce that 
\begin{align}
\norm{D^a D_{t,\qbn-1}^b D \hat u_{{\qbn-1}}}_{L^\infty(\supp \psi_{i,\qbn})} \lesssim  (\la_{\qbn-1}\Ga_{\qbn-1})^a \MM{b+1,\Nindt,\Gamma_{\qbn}^{i_m+1}\tau_{\qbn-1}^{-1}, \Gamma_{\qbn-1}^{-1}\Tau_{\qbn-1}^{-1}}
\label{eq:ad:Dt:a:D:Merlot}
\end{align}
for all  $a + b \leq \sfrac{3 \Nfin}{2}$.  When combined with the formula in \cite[equation~(A.49)]{BMNV21}, which allows us to write
\begin{align}
(\ad D_{t,\qbn-1})^\gamma(D) = f_{\gamma,\qbn-1} \cdot \nabla
\label{eq:ad:Dt:a:D:Merlot:1}
\end{align}
for an explicit function $f_{\gamma,\qbn-1}$ which is defined in terms of $\hat u_{{\qbn-1}}$, estimate \eqref{eq:ad:Dt:a:D:Merlot} and the Leibniz rule gives the estimate
\begin{align}
\norm{ D^a f_{\gamma,\qbn-1} }_{L^\infty(\supp \psi_{i,q})} \lesssim (\la_{\qbn-1}\Ga_{\qbn-1})^a \MM{\gamma,\Nindt,\Ga_\qbn^{i+1}\tau_{\qbn-1}^{-1},\Ga_{\qbn-1}^{-1}\Tau_{\qbn-1}^{-1}}
\label{eq:ad:Dt:a:D:Merlot:2}
\end{align}
for all $a + \gamma\leq \sfrac{3\Nfin}{2}$. In order to conclude the proof of \eqref{eq:product:of:ad:Dt:q-1}--\eqref{eq:product:of:ad:Dt:q-1:bnd}, we use \eqref{eq:ad:Dt:a:D:Merlot:1} to write
\begin{align*}
\prod_{\ell=1}^{\alpha_{k-1}} (\ad D_{t,\qbn-1})^{\gamma_\ell}(D) = \prod_{\ell=1}^{\alpha_{k-1}} \left( f_{\gamma_\ell,\qbn-1} \cdot \nabla \right)= \sum_{j=1}^{\alpha_{k-1}} g_j D^j \, ,
\end{align*} 
and now the claimed estimate for $g_j$ follows from the previously established bound \eqref{eq:ad:Dt:a:D:Merlot:2} for the $f_{\gamma_\ell,q-1}$'s and their derivatives and the Leibniz rule.

With \eqref{eq:product:of:ad:Dt:q-1}--\eqref{eq:product:of:ad:Dt:q-1:bnd} and \eqref{eq:nasty:D:wq} with $k' = k-1$ in hand, we return to \eqref{eq:product:of:ad:Dt:q-1:comm} and obtain
\begin{align}
&\norm{D^{\alpha_k} \left[D_{t,\qbn-1}^{\beta_k}, D^{\alpha_{k-1}}\right] D_{t,\qbn-1}^{\beta_{k-1}}P_{k-2}}_{L^\infty( \supp \psi_{i,\qbn})} \notag\\
&\lesssim \sum_{j=1}^{\alpha_{k-1}} \sum_{1 \leq |\gamma| \leq \beta_k} \norm{ D^{\alpha_k} \left( g_j \; D^j D_{t,\qbn-1}^{\beta_k+\beta_{k-1}-|\gamma|}P_{k-2} \right) }_{L^\infty( \supp \psi_{i,\qbn})} \notag\\
&\lesssim  \sum_{j=1}^{\alpha_{k-1}} 
\sum_{1 \leq |\gamma| \leq \beta_k}  \sum_{a'=0}^{\alpha_k} \norm{ D^{\alpha_k-a'} g_j}_{L^\infty(\supp \psi_{i,\qbn})} \norm{D^{a'+j} D_{t,\qbn-1}^{\beta_k+\beta_{k-1}-|\gamma|}P_{k-2} }_{L^\infty( \supp \psi_{i,\qbn})}   \notag\\
&\lesssim \sum_{j=1}^{\alpha_{k-1}} \sum_{|\gamma|=1}^{\beta_k} \sum_{a'=0}^{\alpha_k} \Ga_\qbn^{i+1} \delta_\qbn^{\sfrac 12}r_q^{-\sfrac 13} (\la_\qbn\Ga_\qbn)^{\alpha_k-a'+\alpha_{k-1}-j} \MM{|\gamma|,\Nindt,\Ga_\qbn^{i+1}\tau_{\qbn-1}^{-1},\Ga_{\qbn-1}^{-1}\Tau_{\qbn-1}^{-1}} \notag\\ &\qquad \qquad \times (\la_\qbn\Ga_\qbn)^{a'+j+K-\alpha_{k-1}-\alpha_k} \MM{M-|\gamma|,\Nindt,\Ga_\qbn^{i+3}\tau_{\qbn-1}^{-1},\Ga_{\qbn-1}\Tau_{\qbn-1}^{-1}} \notag \\
&\lesssim \Ga_\qbn^{i+1} \delta_\qbn^{\sfrac 12}r_q^{-\sfrac 13} (\la_\qbn\Ga_\qbn)^{K} \MM{M,\Nindt,\Ga_\qbn^{i+3}\tau_{\qbn-1}^{-1},\Ga_{\qbn-1}\Tau_{\qbn-1}^{-1}}
\label{eq:nasty:D:wq:***}
\end{align}
for $K+M\leq \sfrac{3\Nfin}{2}+1$. The $+1$ in the range of derivatives is simply a consequence of the fact that the summand in the third line of the above display starts with $j\geq 1$ and with $|\gamma|\geq 1$, so that only $\sfrac{3\Nfin}{2}$ derivatives may fall on $g_j$, which is the extent of the bounds from \eqref{eq:product:of:ad:Dt:q-1:bnd}. This concludes the proof of the inductive step for \eqref{eq:nasty:D:wq}.
\smallskip

\noindent{\bf Proof of \eqref{eq:nasty:Dt:wq}.\,} 
This estimate follows from Lemma~\ref{lem:cooper:1}.  Indeed,  letting $ v = f = \hat w_\qbn$, $B = D_{t,\qbn-1}$, $\Omega = \supp \psi_{i,\qbn}$, $p=\infty$, the previously established bound \eqref{eq:nasty:D:wq} allows us to verify conditions \eqref{eq:cooper:v}--\eqref{eq:cooper:f}  of Lemma~\ref{lem:cooper:1} with  $N_* = \sfrac{3\Nfin}{2}+1$, $\const_v = \const_f = \Gamma_\qbn^{i+1} \hat \delta_\qbn^{\sfrac 12}r_q^{-\sfrac 13}$, $\lambda_v = \lambda_f = \tilde \lambda_v = \tilde \lambda_f= \Ga_\qbn \la_\qbn$, $N_x = \infty$, $\mu_v =  \mu_f = \Gamma_\qbn^{i+3}  \tau_{\qbn-1}^{-1}$, $\tilde \mu_v = \tilde \mu_f = \Gamma_{\qbn-1} \Tau_{\qbn-1}^{-1}$, and $N_t = \Nindt$. The bound \eqref{eq:nasty:Dt:wq} now is a direct consequence of \eqref{eq:cooper:f:**}.
\smallskip

\noindent{\bf Proof of \eqref{eq:nasty:D:vq}.\,}
First we consider the bound \eqref{eq:nasty:D:vq}, inductively on $k$. For the case $k=1$ we appeal to estimate \eqref{eq:cooper:f:*} in Lemma~\ref{lem:cooper:1} with the operators $A = D_\qbn, B = D_{t,\qbn-1}$ and the functions $v = \hat  w_\qbn$ and $f = D \hat  u_\qbn$, so that $D^n (A+B)^m f  = D^n D_{t,\qbn}^m D \hat  u_\qbn$. As before, the assumption \eqref{eq:cooper:v} holds due to \eqref{eq:nasty:D:wq} with the same parameter choices. Verifying condition \eqref{eq:cooper:f} is this time more involved, and follows by rewriting $f = D \hat u_{q} = D \hat w_q + D \hat u_{{q-1}}$. By using \eqref{eq:nasty:D:wq}, and the parameter inequality \eqref{ineq:tau:q}, we conveniently obtain
\begin{align}
&\norm{ \Big( \prod_{i=1}^k D^{\alpha_i} D_{t,\qbn-1}^{\beta_i} \Big) D \hat  w_\qbn }_{L^\infty(\supp \psi_{i,\qbn})} \notag\\
&\qquad \lesssim  \Gamma_{\qbn}^{i-5} \tau_\qbn^{-1} (\la_\qbn\Ga_\qbn)^K
\MM{M,\Nindt,\Gamma_\qbn^{i-5}  \tau_\qbn^{-1},\Gamma_{\qbn-1}^{} \Tau_{\qbn-1}^{-1}}
\label{eq:cutoff:nam:2}
\end{align}
for all $|\alpha| + |\beta| = K+M  \leq \sfrac{3\Nfin}{2}$ (note that the maximal number of derivatives is not $\sfrac{3\Nfin}{2} +1$ anymore, but instead it is just  $\sfrac{3\Nfin}{2}$; the reason is that we are estimating $D \hat  w_q$ and not $\hat  w_q$). On the other hand, from the inductive assumption \eqref{eq:nasty:D:vq:old} with $q' = \qbn-1$ we obtain that 
\begin{align*}
&\norm{ \Big( \prod_{i=1}^k D^{\alpha_i} D_{t,\qbn-1}^{\beta_i} \Big) D \hat u_{{\qbn-1}} }_{L^\infty(\supp \psi_{j,\qbn-1})} \notag\\
&\qquad \qquad \les \tau_{\qbn-1}^{-1}\Ga_{\qbn-1}^{j-4} (\la_{\qbn-1}\Ga_{\qbn-1})^K \MM{M,\Nindt, \Gamma_{\qbn-1}^{j} \tau_{\qbn-1}^{-1}, \Tau_{\qbn-1}^{-1}\Ga_{q+\bn}}
\end{align*}
for $K+M  \leq \sfrac{3 \Nfin}{2}$. Recalling the definitions \eqref{eq:psi:m:im:q:def}--\eqref{eq:psi:i:q:recursive} and the notation \eqref{eq:new:supp:notation}, we have that $(x,t) \in \supp (\psi_{i,\qbn})$ if and only if $(x,t)\in \supp(\psi_{\Vec{i},\qbn})$, and so for every $m\in \{0,\ldots,\NcutSmall\}$, there exists $j_m$ with $\Gamma_{\qbn-1}^{j_m} \leq \Gamma_{\qbn}^{i_m} \leq \Gamma_\qbn^i$ and $(x,t) \in \supp(\psi_{j_m,\qbn-1})$. Thus, the above stated estimate and \eqref{ineq:tau:q} imply that
\begin{align}
&\norm{ \Big( \prod_{i=1}^k D^{\alpha_i} D_{t,\qbn-1}^{\beta_i} \Big) D \hat u_{{\qbn-1}} }_{L^\infty(\supp \psi_{i,\qbn})} \notag\\
&\qquad \qquad \les \tau_{\qbn}^{-1}\Ga_{\qbn}^{i-10}  (\la_{\qbn-1}\Ga_{\qbn-1})^K \MM{M,\Nindt, \Gamma_{\qbn}^{i-10} \tau_{\qbn}^{-1}, \Tau_{\qbn-1}^{-1}\Ga_{\qbn}}
\label{eq:cutoff:nam:4}
\end{align}
whenever $K+M\leq \sfrac{3\Nfin}{2}$. Combining \eqref{eq:cutoff:nam:2} and \eqref{eq:cutoff:nam:4}, we may now verify condition \eqref{eq:cooper:f} for $f = D \hat  u_\qbn$, with $p = \infty$, $\Omega = \supp (\psi_{i,\qbn})$, $\const_f = \Gamma_{\qbn}^{i-5} \tau_\qbn^{-1}$, $\lambda_f = \tilde \lambda_f =   \lambda_\qbn \Ga_\qbn$, $N_x = \infty$, $\mu_f =   \Gamma_\qbn^{i-5}  \tau_\qbn^{-1}, \tilde \mu_f = \Gamma_{\qbn-1} \Tau_{\qbn-1}^{-1}, N_t =  \Nindt$, and $N_*= \sfrac{3\Nfin}{2}$. We may thus appeal to \eqref{eq:cooper:f:*} and obtain that 
\begin{align*}
\norm{ D^{K} D_{t,\qbn}^{M}  D \hat u_\qbn }_{L^\infty(\supp \psi_{i,\qbn})}
& \les 
\Gamma_\qbn^{i-5} \tau_\qbn^{-1} (\la_\qbn\Ga_\qbn)^K \MM{M,\Nindt,\Gamma_\qbn^{i-5}  \tau_\qbn^{-1},\Gamma_{\qbn-1}^{} \Tau_{\qbn-1}^{-1}}
\end{align*}
whenever $K+M\leq \sfrac{3\Nfin}{2}$, concluding the proof of  \eqref{eq:nasty:D:vq} for $k=1$.

In order to prove \eqref{eq:nasty:D:vq} for a general $k$, we proceed by induction. Assume the estimate holds for every $k' \leq k-1$. Proving \eqref{eq:nasty:D:vq} at level $k$ is done in the same way as we have established the induction step (in $k$) for \eqref{eq:nasty:D:wq}. We let 
\begin{align*}
\tilde P_{k'} = \left( \prod_{i=1}^{k'} D^{\alpha_i} D_{t,\qbn}^{\beta_i} \right) D \hat u_\qbn
\end{align*}
and decompose
\begin{align*}
\left( \prod_{i=1}^k D^{\alpha_i} D_{t,\qbn}^{\beta_i} \right) D \hat u_\qbn 
&= (D^{\alpha_k+\alpha_{k-1}} D_{t,\qbn}^{\beta_k+\beta_{k-1}}) \tilde P_{k-2}  + D^{\alpha_k} \left[D_{t,\qbn}^{\beta_k}, D^{\alpha_{k-1}}\right] D_{t,\qbn}^{\beta_{k-1}}  \tilde P_{k-2} \, .
\end{align*}
Note that the first term is directly bounded using the induction assumption at level $k-1$. To bound the commutator term, similarly to \eqref{eq:product:of:ad:Dt:q-1:comm}--\eqref{eq:product:of:ad:Dt:q-1:bnd}, we obtain that 
\begin{align*}
D^{\alpha_k} \left[D_{t,\qbn}^{\beta_k}, D^{\alpha_{k-1}}\right] D_{t,\qbn}^{\beta_{k-1}} \tilde P_{k-2}  
= D^{\alpha_k} \sum_{1 \leq |\gamma| \leq \beta_k} \frac{\beta_k!}{\gamma! (\beta_k - |\gamma|)!} \left( \sum_{j=1}^{\alpha_{k-1}} \tilde g_j D^j \right) D_{t,\qbn}^{\beta_k+\beta_{k-1}-|\gamma|} \tilde P_{k-2} \, ,
\end{align*}
where one may use the previously established bound \eqref{eq:nasty:D:vq} with $k=1$ (instead of \eqref{eq:ad:Dt:a:D:Merlot}) to estimate $\norm{D^a \tilde g_j}_{L^\infty(\supp \psi_{i,\qbn})}$ 
 The estimate 
\begin{align}
&\norm{D^{\alpha_k} \left[D_{t,\qbn}^{\beta_k}, D^{\alpha_{k-1}}\right] D_{t,\qbn}^{\beta_{k-1}} \tilde P_{k-2}  }_{L^\infty(\supp \psi_{i,\qbn})}  \notag\\
&\quad \lesssim  \tau_\qbn^{-1}\Ga_{\qbn}^{i-5} (\la_\qbn\Ga_\qbn)^K \MM{M,\Nindt,\Gamma_{\qbn}^{i-5}  \tau_\qbn^{-1},\Gamma_{\qbn-1} \Tau_{\qbn-1}^{-1}}
\label{eq:cutoff:nam:5}
\end{align}
follows similarly to \eqref{eq:nasty:D:wq:***}, from the estimate for $\tilde g_j$ and the bound \eqref{eq:nasty:D:vq} with $k-1$ terms in the product. 
This concludes the proof of  estimate \eqref{eq:nasty:D:vq}.
\smallskip

\noindent\textbf{Proof of  \eqref{eq:nasty:no:D:vq}.\,} The proof of this bound is nearly identical to that of \eqref{eq:nasty:D:vq}, as is readily seen for $k=1$:  we just need to replace $D\hat w_\qbn$ estimates with $\hat w_\qbn$ estimates, and $D \hat u_{{\qbn-1}}$ bounds with $\hat u_{{\qbn-1}}$ bounds. For instance, instead of \eqref{eq:cutoff:nam:2}, we appeal to \eqref{eq:nasty:Dt:wq:WEAK} and obtain a bound for $D^K D_{t,\qbn}^M \hat w_\qbn$ which is better than \eqref{eq:cutoff:nam:2} by a factor of $\lambda_\qbn\Ga_\qbn$, and which holds for $K+M \leq \sfrac{3\Nfin}{2}+1$. This estimate is sharper than required by \eqref{eq:nasty:no:D:vq}. The estimate for $D^K D_{t,\qbn}^M \hat u_{{\qbn-1}}$ is obtained similarly to \eqref{eq:cutoff:nam:4}, except that instead of appealing to the induction assumption \eqref{eq:nasty:D:vq:old} at level $q'=\qbn-1$, we use \eqref{eq:bob:Dq':old} with $q'=\qbn-1$. The estimates hold for $K+M \leq \sfrac{3\Nfin}{2}+1$. These arguments establish \eqref{eq:nasty:no:D:vq} with $k=1$. The case of general $k\geq 2$ is treated inductively exactly as before, because the commutator term is bounded in the same way as \eqref{eq:cutoff:nam:5}, except that $K+1$ is replaced by $K$. To avoid redundancy, we omit these details.
\smallskip

\noindent\textbf{Proof of \eqref{eq:bobby:new}.\,} The proof of this bound is immediate from \eqref{eq:vellie:inductive:dtq-1:uniform:upgraded:statement}, the definition of $\hat w_\qbn$ in Lemma~\ref{lem:mollifying:w}, the inductive assumption \eqref{eq:bobby:old}, and the triangle inequality.
\end{proof}

\subsection{Material derivatives}

\begin{remark}[\bf Rewriting $\psi_{i,\qbn}$]
\label{rem:rewrite:cutoffs}
In order to take material derivatives of $\psi_{i,\qbn}$, we need to take advantage of certain cancellations. For this purpose, we introduce the summed cutoff function 
\begin{align}
\Psi_{m,i,\qbn}^6 = \sum_{i_m=0}^i \psi_{m,i_m,\qbn}^6 
\label{eq:fancy:cutoff}
\end{align}
for any given $0\leq m \leq \NcutSmall$ and note via Lemma~\ref{lem:partition:of:unity:psi:m} that
\begin{align}
D(\Psi_{m,i,\qbn}^6) = D(\psi_{m,i,\qbn}^6)  {\bf 1}_{\supp(\psi_{m,i+1,\qbn})} 
\,.
\label{eq:fancy:cutoff:supp}
\end{align}
With the notation \eqref{eq:fancy:cutoff} we return to the definition \eqref{eq:psi:i:q:recursive} and note that 
\begin{align}
\psi_{i,\qbn}^6 
&= \sum_{m=0}^{\NcutSmall} \psi_{m,i,\qbn}^6 \prod_{m'=0}^{m-1}  \Psi_{m',i,\qbn}^6  \prod_{m''=m+1}^{\NcutSmall} ( \Psi_{m'',i,\qbn}^6 - \psi_{m'',i,\qbn}^6) 
\notag\\
&= \sum_{m=0}^{\NcutSmall} \psi_{m,i,\qbn}^6 \prod_{m'=0}^{m-1}  \Psi_{m',i,\qbn}^6  \prod_{m''=m+1}^{\NcutSmall}  \Psi_{m'',i-1,\qbn}^6
\,.
\label{eq:cutoff:resummation} 
\end{align}
\end{remark}

Inspecting   \eqref{eq:cutoff:resummation} and using identity \eqref{eq:fancy:cutoff:supp} and the definitions \eqref{eq:new:supp:notation}, \eqref{eq:fancy:cutoff}, we see that 
\begin{align}
(x,t) \in \supp (D_{t,\qbn-1}  \psi_{i,\qbn}^6) \quad {\implies} \quad 
&\exists \Vec{i}\in \N_0^{\NcutSmall+1} \mbox{ and } \exists 0\leq m \leq\NcutSmall \notag\\
&\mbox{with }i_m \in \{i-1,i\} \mbox{ and }  \max_{0\leq m'\leq \NcutSmall} i_{m'} = i \notag\\
&\mbox{such that } (x,t) \in \supp(\psi_{\Vec{i},\qbn}) \cap \supp(D_{t,\qbn-1} \psi_{m,i_m,\qbn}) \notag\\
&  \mbox{and }  i_{m'} \leq i_{m} \mbox{ whenever }  m < m' \leq \NcutSmall\,.
\label{eq:orangutan}
\end{align}
The generalization of characterization \eqref{eq:orangutan} to higher order material derivatives $D_{t,\qbn-1}^M$ is direct:  $(x,t) \in \supp (D_{t,\qbn-1}^M  \psi_{i,\qbn}^6)$ implies that there exists $ \Vec{i}\in \N_0^{\NcutSmall+1}$ with maximal index equal to $i$,  such that for every $0\leq m \leq \NcutSmall$ for which $(x,t) \in \supp(\psi_{\Vec{i},\qbn}) \cap \supp(D_{t,\qbn-1} \psi_{m,i_m,\qbn}) $, we have $i_{m'} \leq i_m \in \{i-1,i\}$ whenever $m< m'$. Using this characterization, we may prove the following.
\begin{lemma}[\bf Mixed derivatives for intermediate velocity cutoff functions]
\label{lem:sharp:Dt:psi:i:j:q}
Let $\qbn\geq 1$, $0\leq i \leq \imax(\qbn)$, and fix $\Vec{i}\in \N_0^{\NcutSmall+1}$ such that $\max_{0\leq m\leq\NcutSmall} i_m =i$, as in the right side of \eqref{eq:orangutan}.
Fix $0 \leq m \leq \NcutSmall$ such that $ i_m \in \{i-1,i\}$ and such that $i_{m'}\leq i_m$ for all $m \leq m'\leq \NcutSmall$, again as in the right hand side of \eqref{eq:orangutan}. Lastly, fix $j_m$   such that $i_*(j_m) \leq i_m$. For  $N,K,M,k \geq 0$, $\alpha,\beta \in {\mathbb N}^k$ such that $|\alpha| = K$ and $|\beta| = M $, we have
\begin{align}
& \frac{ {\bf 1}_{\supp(\psi_{\Vec{i},\qbn}) } {\bf 1}_{\supp(\psi_{j_m,\qbn-1})}}{\psi_{m,i_m,j_m,\qbn}^{1- (K+M)/\Nfin}} \left| \left( \prod_{l=1}^k D^{\alpha_l} D_{t,\qbn-1}^{\beta_l}\right)  \psi_{m,i_m,j_m,\qbn} \right|  \notag\\
&\quad  \les
(\la_\qbn\Ga_\qbn)^K
\MM{M,\Nindt-\NcutLarge, \Gamma_{\qbn}^{i+3} \tau_{\qbn-1}^{-1}, \Gamma_{\qbn-1} \Tau_{\qbn-1}^{-1} }
\label{eq:sharp:Dt:psi:i:j:q}
\end{align}
for all $K$ such that $0 \leq K+M \leq \Nfin$.
Moreover,
\begin{align}
&  \frac{ {\bf 1}_{\supp(\psi_{\Vec{i},\qbn}) } {\bf 1}_{\supp(\psi_{j_m,\qbn-1})}}{\psi_{m,i_m,j_m,\qbn}^{1- (N+K+M)/\Nfin}} \left| D^N \left( \prod_{l=1}^k D_\qbn^{\alpha_l} D_{t,\qbn-1}^{\beta_l}\right)  \psi_{m,i_m,j_m,\qbn} \right|  \notag\\
&\quad  \lesssim 
(\la_\qbn\Ga_\qbn)^N (\tau_{\qbn}^{-1}\Ga_{\qbn}^{i-5})^K
\MM{M,\Nindt-\NcutLarge, \Gamma_{\qbn}^{i+3} \tau_{\qbn-1}^{-1}, \Gamma_{\qbn-1} \Tau_{\qbn-1}^{-1} }
\label{eq:sharp:Dt:psi:i:j:q:mixed}
\end{align}
holds whenever $0 \leq N+ K+M \leq \Nfin$.
\end{lemma}
\begin{proof}[Proof of Lemma~\ref{lem:sharp:Dt:psi:i:j:q}]
Note that for $M=0$ estimate \eqref{eq:sharp:Dt:psi:i:j:q} was already established in \eqref{eq:sharp:D:psi:i:j:q}. The bound \eqref{eq:sharp:Dt:psi:i:j:q:mixed} with $M=0$, i.e., an estimate for the $D^N D_\qbn^K \psi_{m,i_m,j_m,\qbn}$, holds by appealing to the expansion \eqref{eq:D:q:K:i}--\eqref{eq:D:q:K:ii},  the bound \eqref{eq:D:a:fj} (which is applicable since in the context of estimate  \eqref{eq:sharp:Dt:psi:i:j:q:mixed} we work on the support of $\psi_{\vec i,\qbn}$),   to the bound \eqref{eq:sharp:Dt:psi:i:j:q} with $M=0$, and to \eqref{ineq:tau:q}.
The rest of the proof is dedicated to the case $M \geq 1$. The proofs are very similar to the proof of Lemma~\ref{lem:sharp:D:psi:i:q}, but we additionally need to appeal to bounds and arguments from the proof of Lemma~\ref{lem:Dt:Dt:wq:psi:i:q:multi}.  
\smallskip

\noindent{\bf Proof of \eqref{eq:sharp:Dt:psi:i:j:q}.\,}
We start with the case $k=1$ and estimate $D^K D_{t,\qbn-1}^M \psi_{m,i_m,j_m,\qbn}$ for $K+M\leq \Nfin$ and $M\geq 1$. We note that the operator $D_{t,\qbn-1}$ is a scalar differential operator, and thus the Faa di Bruno argument which was used to  bound \eqref{eq:sharp:D:psi:i:j:q} may be repeated. As was done there, we recall the definitions \eqref{eq:psi:i:j:def}--\eqref{eq:psi:i:i:def} and split the analysis in two cases, according to whether \eqref{eq:clusterf:1} or \eqref{eq:clusterf:3} holds. 

Let us first consider the case \eqref{eq:clusterf:1}. Our goal is to apply  \cite[Lemma~A.5]{BMNV21} to the function $\psi = \gamma_{m,\qbn}$ or $\psi=\tilde \gamma_{m,\qbn}$, with $\Gamma_\psi = \Gamma_{\qbn}^{m+1}$, $\Gamma = \Gamma_{\qbn}^{(m+1)(i_m-i_*(j_m))}$,   $h(x,t) =  h_{m,j_m,\qbn}^2(x,t)$, and $D_t = D_{t,\qbn-1}$. The estimate in \cite[(A.24)]{BMNV21} again holds by \eqref{eq:DN:psi:q:0} and \eqref{eq:DN:psi:q:gain}, and so it remains to obtain a bound on the derivatives of $(h_{m,j_m,\qbn}(x,t))^2$ on the set $\supp(\psi_{\Vec i,q}) \cap \supp(\psi_{j_m,q-1} \psi_{m,i_m,j_m,q})$ in order to satisfy \cite[(A.25)]{BMNV21}. Similarly to \eqref{eq:cutoff:spatial:derivatives:000}, for $K'+M' \leq \Nfin$ the Leibniz rule and definition \eqref{eq:h:j:q:def} gives
\begin{align}
&\left|D^{K'} D_{t,\qbn-1}^{M'} h_{m,j_m,\qbn}^2\right| \notag\\
&\qquad \lesssim (\lambda_\qbn \Gamma_\qbn)^{K'} 
(\tau_{\qbn-1}^{-1} \Gamma_\qbn^2 )^{M'} \Gamma_\qbn^{-2(m+1) i_*(j_m) }   \notag\\
&\qquad \quad \times
\sum_{K''=0}^{K'} \sum_{M'' =0}^{M'} \sum_{n=0}^{\NcutLarge}  
( \tau_{\qbn-1}^{-1} \Gamma_\qbn^2)^{-m-M''} 
(\lambda_\qbn \Gamma_\qbn)^{-n-K''}  \delta_{\qbn}^{-\sfrac 12} r_q^{\sfrac 13}
| D^{n + K''} D^{m+M''}_{t,\qbn-1} \hat w_\qbn| 
\notag\\
&\qquad \qquad     \times 
( \tau_{\qbn-1}^{-1} \Gamma_\qbn^2)^{-m-M'+M''} 
(\lambda_\qbn \Gamma_\qbn)^{-n - K' + K''} \delta_\qbn^{-\sfrac 12} r_q^{\sfrac 13}
| D^{n + K'-K''} D^{m+M'-M''}_{t,\qbn-1} \hat w_\qbn| 
\,.
\label{eq:5:123}
\end{align}
By the characterization \eqref{eq:orangutan}, for every $(x,t)$ in the support described on the left side of \eqref{eq:sharp:Dt:psi:i:j:q}  we have that for every $m \leq R \leq \NcutSmall$, there exists $i_R \leq i_m$ and $j_R$ with $i_*(j_R) \leq i_R$, such that $(x,t) \in \supp  \psi_{j_R,\qbn-1} \psi_{R,i_R,j_R,\qbn}$. As a consequence, for the terms in the sum \eqref{eq:5:123} with $L \in \{n+K'',n+K'-K''\} \leq \NcutLarge$ 
and $R \in \{m+M'',m+M'-M''\} \leq \NcutSmall$, 
we may appeal to estimate \eqref{eq:h:psi:supp} which gives a bound on $h_{R,j_R,\qbn}$, and thus obtain
\begin{align*}
&(\tau_{\qbn-1}^{-1} \Gamma_\qbn^2)^{-R}  (\lambda_\qbn \Gamma_\qbn)^{-L} \delta_\qbn^{-\sfrac 12} r_q^{\sfrac 13} \norm{D^L D_{t,q-1}^R \hat w_\qbn}_{L^\infty( \supp \psi_{R,i_R,j_R,\qbn} )} \notag\\
&\qquad \qquad \leq \Gamma_\qbn^{(R+1)i_*(j_R)}  \Gamma_\qbn^{(R+1) (i_R + 1 - i_*(j_R))} 
\notag\\
&\qquad \qquad \leq  \Gamma_\qbn^{(R+1) (i_m + 1)}  \, .
\end{align*}
On the other hand, if $L > \NcutLarge$, or if $R > \NcutSmall$, then by \eqref{eq:vellie:inductive:dtq-1:uniform:upgraded:statement}, we have that 
\begin{align}
(\tau_{\qbn-1}^{-1} \Gamma_\qbn^2)^{-R}  (\lambda_\qbn &\Gamma_\qbn)^{-L} \delta_\qbn^{-\sfrac 12} r_q^{\sfrac 13} \norm{D^L D_{t,\qbn-1}^R \hat w_\qbn}_{L^\infty(\supp \psi_{j_m,\qbn-1})} 
\notag\\
&  \leq \Ga_q^{\sfrac{\badshaq}{2}+16} r_q^{-1} \Ga_\qbn^{-L} \Ga_{\qbn-1}^L \Ga_\qbn^{-2R}
\MM{R,\NindSmall,\Gamma_{\qbn-1}^{j_m-1}, \tau_{\qbn-1} \Tau_{\qbn-1}^{-1}}
\notag\\
&  \leq  \MM{R,\NindSmall,\Gamma_{\qbn}^{i_m-1}, \tau_{\qbn-1} \Tau_{\qbn-1}^{-1}}
\, .
\label{eq:cutoff:spatial:derivatives:0000}
\end{align}
since $\NcutLarge$ and $\NcutSmall$  were taken sufficiently large to obey \eqref{condi.Ncut0} and $i_m\geq i_*(j_m)$. Combining \eqref{eq:5:123}--\eqref{eq:cutoff:spatial:derivatives:0000}, we have that
\begin{align}
&{\bf 1}_{\supp(\psi_{\Vec{i},\qbn}) } {\bf 1}_{\supp(\psi_{j_m,\qbn-1})} \left|D^{K'} D_{t,\qbn-1}^{M'} h_{m,j_m,\qbn}^2\right|
\notag\\
&\lesssim \Gamma_{\qbn}^{2(m+1)(i_m - i_*(j_m)+1)}    (\la_{\qbn}\Ga_\qbn)^{K'} \MM{ M',\NindSmall-\NcutSmall,\tau_{\qbn-1}^{-1} \Gamma_{\qbn}^{i+3}, \Tau_{\qbn-1}^{-1}}
\label{eq:5:126}
\end{align}
for all $K'+M'\leq \Nfin$. The upshot of \eqref{eq:5:126} is that the condition in \cite[(A.25)]{BMNV21} is now verified, with  $\const_h = \Gamma_{\qbn}^{2(m+1)(i_m - i_*(j_m)+1)}$, and $\lambda = \tilde \lambda =  \Gamma_\qbn \lambda_\qbn$,  $\mu =\tau_{\qbn-1}^{-1} \Gamma_{\qbn}^{i+3}$, $\tilde  \mu = \Tau_{\qbn-1}^{-1}$, and $N_t = \Nindt-\NcutSmall$. We obtain from \cite[(A.26)]{BMNV21} and the fact that $(\Gamma_\psi \Gamma)^{-2} \const_h = 1$ that \eqref{eq:sharp:Dt:psi:i:j:q} holds when $k=1$  for those $(x,t)$ such that $h_{m,j_m,\qbn}(x,t)$ satisfies \eqref{eq:clusterf:1}. The case when $h_{m,j_m,\qbn}(x,t)$ satisfies the bound  \eqref{eq:clusterf:3} is nearly identical, as  was the case in the proof of Lemma~\ref{lem:sharp:D:psi:i:q}. The only changes are that now $\Gamma_\psi = 1$ (according to \eqref{eq:DN:psi:q}), and that the constant $\const_h$ which we read from the right side of \eqref{eq:5:126} is now improved to $\Gamma_{\qbn}^{2(m+1)(i_m - i_*(j_m))}$. These two changes offset each other, resulting in the same exact bound. Thus, we have shown that \eqref{eq:sharp:Dt:psi:i:j:q} holds when $k=1$. 

The general case $k\geq 1$ in \eqref{eq:sharp:Dt:psi:i:j:q} is obtained via induction on $k$, in precisely the same fashion as the proof of estimate \eqref{eq:nasty:D:wq} in Lemma~\ref{lem:Dt:Dt:wq:psi:i:q:multi}. At the heart of the matter lies a commutator bound similar to \eqref{eq:nasty:D:wq:***}, which  is proven  in precisely the same way by appealing to the fact that we work on $\supp (\psi_{\Vec{i},\qbn}) \subset \supp (\psi_{i,\qbn})$, and thus bound \eqref{eq:product:of:ad:Dt:q-1:bnd} is available; in turn, this bound provides sharper space and material estimates than required in \eqref{eq:sharp:Dt:psi:i:j:q}, completing the proof. In order to avoid redundancy we omit further details.
\smallskip

\noindent{\bf Proof of \eqref{eq:sharp:Dt:psi:i:j:q:mixed}.\,}
This estimate follows from Lemma~\ref{lem:cooper:1} in a manner identical to the proof of \cite[(6.77)]{BMNV21}, and we omit the details.
\end{proof}

\begin{lemma}[\bf Mixed spatial and material derivatives for velocity cutoffs]
\label{lem:sharp:Dt:psi:i:q}
Let $q+\bn\geq 1$, $0 \leq i \leq i_{\rm max}(q+\bn)$, $N,K,M,k \geq 0$, and let $\alpha,\beta \in {\mathbb N}^k$ be such that $|\alpha|=K$ and $|\beta|=M$.  Then we have
\begin{align}
&\frac{1}{\psi_{i,q+\bn}^{1- (K+M)/\Nfin}} \left|\left(\prod_{l=1}^k D^{\alpha_l} D_{t,q+\bn-1}^{\beta_l}\right) \psi_{i,q+\bn}\right| \notag\\
&\qquad \les (\la_\qbn\Ga_\qbn)^K
\MM{M,\Nindt-\NcutSmall,\Gamma_{{q+\bn}}^{i+3}  \tau_{q+\bn-1}^{-1}, \Gamma_{\qbn+1} \Tau_{\qbn-1}^{-1}}
\label{eq:sharp:Dt:psi:i:q}
\end{align}
for $K + M \leq \Nfin$, and 
\begin{align}
&\frac{1}{\psi_{i,\qbn}^{1- (N+K+M)/\Nfin}} \left| D^N \left( \prod_{l=1}^k D_\qbn^{\alpha_l} D_{t,\qbn-1}^{\beta_l}\right)  \psi_{i,\qbn} \right| \notag\\
&\qquad \les  (\la_\qbn\Ga_\qbn)^N
(\Gamma_\qbn^{i-5} \tau_\qbn^{-1})^K 
\MM{M,\Nindt-\NcutSmall,\Gamma_{{q+\bn}}^{i+3}  \tau_{q+\bn-1}^{-1}, \Gamma_{\qbn+1} \Tau_{\qbn-1}^{-1}}
\label{eq:sharp:Dt:psi:i:q:mixed}
\end{align}
holds for $N+ K+ M \leq \Nfin$.
\end{lemma}

\begin{proof}[Proof of Lemma~\ref{lem:sharp:Dt:psi:i:q}]
Note that for $M = 0$ estimate \eqref{eq:sharp:Dt:psi:i:q} holds by \eqref{eq:sharp:D:psi:i:q}. The bound \eqref{eq:sharp:Dt:psi:i:q:mixed} holds for $M=0$, due to the expansion \eqref{eq:D:q:K:i}--\eqref{eq:D:q:K:ii},  the bound \eqref{eq:D:a:fj}  on the support of $\psi_{i,\qbn}$, the bound \eqref{eq:sharp:Dt:psi:i:q} with $M=0$, and to the parameter inequality \eqref{ineq:tau:q}. The rest of the proof is dedicated to the cases $M \geq 1$ for both bounds. 

The argument is very similar to the proof of Lemma~\ref{lem:sharp:D:psi:i:q} and so we only emphasize the main differences. We start with the proof of \eqref{eq:sharp:Dt:psi:i:q}. We claim that in a the same way that \eqref{eq:sharp:D:psi:i:j:q} was shown to imply \eqref{eq:sharp:D:psi:im:q}, one may show that estimate \eqref{eq:sharp:Dt:psi:i:j:q} implies that for any $\Vec{i}$ and $0\leq m \leq \NcutSmall$ as on the right side of \eqref{eq:orangutan} (in particular, as in Lemma~\ref{lem:Dt:Dt:wq:psi:i:q:multi}), we have that 
\begin{align}
& \frac{ {\bf 1}_{\supp(\psi_{\Vec{i},\qbn})}}{\psi_{m,i_m,\qbn}^{1- (K+M)/\Nfin}} \left| \left( \prod_{l=1}^k D^{\alpha_l} D_{t,\qbn-1}^{\beta_l}\right)  \psi_{m,i_m,\qbn} \right|  \notag\\
&\qquad  \les
(\la_\qbn\Ga_\qbn)^K
\MM{M,\Nindt-\NcutLarge, \Gamma_{\qbn-1}^{i+3}  \tau_{\qbn-1}^{-1}, \Gamma_{{\qbn}} \Tau_{\qbn-1}^{-1} }
\,.
\label{eq:5:129}
\end{align}
The proof of the above estimate is done by induction on $k$. For $k=1$, the first step in establishing \eqref{eq:5:129} is to use the Leibniz rule and induction on the number of material derivatives to reduce the problem to an estimate for $\psi_{m,i_m,\qbn}^{-6+ (K+M)/\Nfin} D^K D_{t,\qbn-1}^M (\psi_{m,i_m,\qbn}^6)$; this is achieved in precisely the same way that \eqref{eq:cutoff:spatial:derivatives:1} was proven. The derivatives of $\psi_{m,i_m,\qbn}^6$ are now bounded via the Leibniz rule and the definition \eqref{eq:psi:m:im:q:def}. Indeed, when $D^{K'} D_{t,\qbn-1}^{M'}$ derivatives fall on $\psi_{m,i_m,j_m,\qbn}^6$, the required bound is obtained from \eqref{eq:sharp:Dt:psi:i:j:q}, which gives the same upper bound as the one required by \eqref{eq:5:129}. On the other hand, if $D^{K-K'} D_{t,\qbn-1}^{M-M'}$ derivatives fall on $\psi_{j_m,\qbn-1}^6$, the required estimate is provided by \eqref{eq:nasty:Dt:psi:i:q:orangutan} with $q' = \qbn-1$ and $i$ replaced by $j_m$; the resulting estimates are strictly better than what is required by \eqref{eq:5:129}. This shows that estimate \eqref{eq:5:129} holds for $k=1$. We then proceed inductively in $k\geq 1$, in  the same fashion as the proof of estimate \eqref{eq:nasty:D:wq} in Lemma~\ref{lem:Dt:Dt:wq:psi:i:q:multi}; the corresponding commutator bound is applicable because we work on $\supp(\psi_{m,i_m,\qbn}) \cap \supp (\psi_{i,\qbn})$.  In order to avoid redundancy we omit these details, and conclude the proof of \eqref{eq:5:129}. 

As in the proof of Lemma~\ref{lem:sharp:D:psi:i:q}, we are now able to show that \eqref{eq:sharp:Dt:psi:i:q} is a consequence of \eqref{eq:5:129}. As before, by induction on the number of material derivatives and the Leibniz rule we reduce the problem to an estimate for $\psi_{i,\qbn}^{-6 + (K+M)/\Nfin} \prod_{l=1}^k D^{\alpha_l} D_{t,\qbn-1}^{\beta_l} (\psi_{i,\qbn}^6)$; see the proof of \eqref{eq:cutoff:spatial:derivatives:1} for details. In order to estimate derivatives of $\psi_{i,\qbn}^6$, we use identities~\eqref{eq:fancy:cutoff:supp} and~\eqref{eq:cutoff:resummation}, which 
imply upon applying a differential operator, say $D_{t,\qbn-1}$, that 
\begin{align}
D_{t,\qbn-1} &(\psi_{i,\qbn}^6) \notag\\
&= D_{t,\qbn-1} \left( \sum_{m=0}^{\NcutSmall}  \prod_{m'=0}^{m-1} \Psi_{m',i,\qbn}^6 \cdot \psi_{m,i,\qbn}^6 \cdot \prod_{m''=m+1}^{\NcutSmall} \Psi_{m'',i-1,\qbn}^6  \right) \notag\\
&=  \sum_{m=0}^{\NcutSmall} \sum_{\bar m'=0}^{m-1} D_{t,\qbn-1}(\psi_{\bar m',i,\qbn}^6)  \prod_{\substack{0\leq m' \leq m-1 \\ m' \neq \bar m'}}  \Psi_{m',i,\qbn}^6 \cdot \psi_{m,i,\qbn}^6\cdot \prod_{m''=m+1}^{\NcutSmall} \Psi_{m'',i-1,\qbn}^6 \notag\\
&\quad + \sum_{m=0}^{\NcutSmall} \sum_{\bar m''=m+1}^{\NcutSmall} \prod_{m'=0}^{m-1} \Psi_{m',i,\qbn}^6 \cdot \psi_{m,i,\qbn}^6 \cdot D_{t,\qbn-1}(\Psi_{\bar m'',i-1,\qbn}^6)  \prod_{\substack{m+1\leq m'' \leq \NcutSmall \\ m'' \neq \bar m''}} \Psi_{m'',i-1,\qbn}^6 \notag\\
&\quad + \sum_{m=0}^{\NcutSmall}  \prod_{m'=0}^{m-1} \Psi_{m',i,\qbn}^6 \cdot  D_{t,\qbn-1} (\psi_{m,i,\qbn}^6)  \cdot  \prod_{m''=m+1}^{\NcutSmall} \Psi_{m'',i-1,\qbn}^6\,.
\label{eq:orangutan:1}
\end{align}
Higher order material derivatives of $\psi_{i,\qbn}^6$, and mixtures of space and material derivatives are obtained similarly, by an application of the Leibniz rule. 
Equality \eqref{eq:orangutan:1} in particular justifies why we have only proven \eqref{eq:5:129} for $\Vec{i}$ and $0\leq m \leq \NcutSmall$ as on the right side of \eqref{eq:orangutan}! With \eqref{eq:5:129} and \eqref{eq:orangutan:1} in hand, we now repeat the argument from the proof of Lemma~\ref{lem:sharp:D:psi:i:q} (see the two displays below \eqref{eq:cutoff:spatial:derivatives:1}) and conclude that \eqref{eq:sharp:Dt:psi:i:q} holds. 

In order to conclude the proof of the Lemma, it remains to establish \eqref{eq:sharp:Dt:psi:i:q:mixed}. This bound follows now directly from \eqref{eq:sharp:Dt:psi:i:q} and an application of Lemma~\ref{lem:cooper:1} (to be more precise, we need to use the proof of this Lemma), in precisely the same way that \eqref{eq:sharp:Dt:psi:i:j:q} was shown earlier to imply \eqref{eq:sharp:Dt:psi:i:j:q:mixed}. As there are no changes to be made to this argument, we omit these details. 
\end{proof}

\subsection{\texorpdfstring{$L^r$}{Lr} size of the velocity cutoffs}
The purpose of this section is to show that the inductive estimate \eqref{eq:psi:i:q:support:old} holds with $q'=\qbn$. 
\begin{lemma}[\bf Support estimate]
\label{lem:psi:i:q:support}
For all $0 \leq i \leq \imax(\qbn)$ and $1\leq r \leq \infty$, we have that 
\begin{align}
\norm{\psi_{i,\qbn}}_{r}  &\lesssim  \Gamma_\qbn^{\frac{-3i+\CLebesgue}{r}}
\label{eq:psi:i:q:support}
\end{align}
where $ \CLebesgue$ is defined in \eqref{eq:psi:i:q:support:old} and thus depends only on $b$. 
\end{lemma}

\begin{proof}[Proof of Lemma~\ref{lem:psi:i:q:support}]
First, note that the cases $1<r\leq \infty$ follow from the case $r=1$ and interpolation.  Next, observe that if $i \leq \sfrac 13 \CLebesgue $, then \eqref{eq:psi:i:q:support} trivially holds because  $0 \leq \psi_{i,\qbn} \leq 1$ for all $\qbn\geq 1$ once $a$ is chosen to be sufficiently large. Thus, we only consider $i$ such that $\sfrac 13 \CLebesgue < i \leq \imax(\qbn)$.

First, we note that Lemma~\ref{lem:partition:of:unity:psi:m} implies that the functions $\Psi_{m,i',\qbn}$ defined in \eqref{eq:fancy:cutoff} satisfy $0 \leq \Psi_{m,i',q}^2 \leq 1$, and thus \eqref{eq:cutoff:resummation} implies that 
\begin{align}
\norm{\psi_{i,\qbn}}_1 \leq \sum_{m=0}^{\NcutSmall} \norm{\psi_{m,i,\qbn}}_1 \,. 
\label{eq:orange:orangutan:0}
\end{align}
Next, we let $j_*(i) = j_*(i,\qbn)$ be the {\em maximal} index of $j_m$ appearing in \eqref{eq:psi:m:im:q:def}. In particular, recalling also \eqref{eq:ineq:ij}, we have that 
\begin{align}
\Gamma_{\qbn}^{i-1} < \Gamma_{\qbn-1}^{j_*(i)} \leq \Gamma_{\qbn}^i < \Gamma_{\qbn-1}^{j_*(i)+1}  \,.
\label{eq:orange:orangutan:1}
\end{align}
Using \eqref{eq:psi:m:im:q:def}, in which we simply write $j$ instead of $j_m$,   the fact that $0\leq \psi_{j,\qbn-1}^2, \psi_{m,i,j,\qbn}^2 \leq 1$, and the inductive assumption \eqref{eq:psi:i:q:support:old} at level $\qbn-1$, we may  deduce that 
\begin{align}
\norm{\psi_{m,i,\qbn}}_1 
&\leq \norm{\psi_{j_*(i),\qbn-1}}_1 + \norm{\psi_{j_*(i)-1,\qbn-1}}_1 + \sum_{j=0}^{j_*(i)-2}  \norm{\psi_{j,\qbn-1} \psi_{m,i,j,\qbn}}_1 \notag\\
&\leq \Gamma_{\qbn-1}^{-3 j_*(i) + \CLebesgue} + \Gamma_{\qbn-1}^{-3 j_*(i) + 3 + \CLebesgue} + \sum_{j=0}^{j_*(i)-2}  \abs{ \supp(\psi_{j,\qbn-1} \psi_{m,i,j,\qbn})} 
\,. \label{eq:orange:orangutan:2}
\end{align}
The second term on the right side of \eqref{eq:orange:orangutan:2} is estimated using the last inequality in \eqref{eq:orange:orangutan:1} as 
\begin{align}
\Gamma_{\qbn-1}^{-3 j_*(i) + 3 + \CLebesgue} 
\leq  \Gamma_{\qbn}^{-3 i} \Gamma_{\qbn-1}^{6 + \CLebesgue}
\leq \Gamma_{\qbn}^{-3 i + \CLebesgue -1}  \Gamma_{\qbn-1}^{6 + \CLebesgue - b (\CLebesgue -1)} 
= \Gamma_{\qbn}^{-3 i + \CLebesgue -1}
\label{eq:orange:orangutan:3}
\end{align}
where in the last equality we have used the definition of $\CLebesgue$ in \eqref{eq:psi:i:q:support:old}.
Clearly, the first term on the right side of \eqref{eq:orange:orangutan:2} is also bounded by the right side of \eqref{eq:orange:orangutan:3}. We are left to estimate the terms appearing in the sum on the right side of \eqref{eq:orange:orangutan:2}. The key fact is that for any $j \leq j_*(i)-2$ we have that $i \geq  i_*(j)+1$; this can be seen to hold because $b < 2$. Recalling \eqref{eq:psi:supp:upper}, for $j\leq j_*(i)-2$ we have that
\begin{align}
\supp(\psi_{j,\qbn-1} \psi_{m,i,j,\qbn}) 
&\subseteq \left\{ (x,t) \in \supp(\psi_{j,\qbn-1}) \colon h_{m,j,\qbn}^3 \geq \frac 18 \Gamma_\qbn^{ 3(m+1)(i-i_*(j))} \right\} \notag\\
&\subseteq \left\{ (x,t)   \colon \psi_{j\pm,\qbn-1}^6 h_{m,j,\qbn}^3 \geq \frac 18 \Gamma_\qbn^{ 3(m+1)(i-i_*(j))} \right\} \,.
\label{eq:orange:orangutan:4}
\end{align}
Here, $\psi_{j\pm,\qbn-1}$ denotes $\psi_{j\pm,\qbn-1}^6 = \displaystyle\sum_{j'=j-1}^{j+1}\psi_{j',\qbn-1}^6$.
In the second inclusion of \eqref{eq:orange:orangutan:4} we have appealed to \eqref{eq:inductive:partition} at level $\qbn-1$. By Chebyshev's inequality and the definition of $h_{m,j,\qbn}$ in \eqref{eq:h:j:q:def} we deduce that
\begin{align*}
\abs{\supp(\psi_{j,\qbn-1} \psi_{m,i,j,\qbn})}
& \leq (2\NcutLarge)^3 \Gamma_\qbn^{-3(m+1)(i-i_*(j))} \sum_{n=0}^{\NcutLarge} \Gamma_\qbn^{-3i_*(j)} \delta_\qbn^{-\sfrac 32} r_q  (\lambda_\qbn \Gamma_\qbn)^{-3n} \notag\\
&\qquad \qquad \qquad \times \left(\tau_{\qbn-1}^{-1} \Gamma_{\qbn}^{i_*(j)+2}\right)^{-3m} \norm{\psi_{j\pm,\qbn-1} D^n D_{t,\qbn-1}^m \hat w_\qbn}_3^3 \,.
\end{align*}
Since in the above display we have that $m\leq \NcutSmall \leq \Nindt$ from \eqref{condi.Nindt}, we may combine the above estimate with \eqref{eq:vellie:inductive:dtq-1:upgraded:statement} to deduce that 
\begin{align}
 \abs{\supp(\psi_{j,\qbn-1} \psi_{m,i,j,\qbn})} 
&\leq 8 \NcutLarge^4 \Gamma_\qbn^{-3(m+1)(i-i_*(j))}  \Gamma_{\qbn}^{-3i_*(j)} \Ga_q^{60} \left( \Gamma_{\qbn-1}^{j-1} \Gamma_\qbn^{-i_*(j)-2} \right)^{3m} \notag\\
&\leq 8 \NcutLarge^4 \Ga_q^{60} \Ga_\qbn^{-3i} \notag\\
&\leq \Ga_{\qbn}^{-3i+\CLebesgue-1}
\,.
\label{eq:orange:orangutan:5}
\end{align}
We have used here that $\Gamma_{\qbn-1}^j \leq \Gamma_\qbn^{i_*(j)}$, that $m\geq 0$, and that $   \CLebesgue \geq 62$ since $b\leq \sfrac{25}{24}$ from \eqref{ineq:b}. 

Combining \eqref{eq:orange:orangutan:0}, \eqref{eq:orange:orangutan:2}, \eqref{eq:orange:orangutan:3}, and \eqref{eq:orange:orangutan:5} we deduce that 
\begin{align*}
\norm{\psi_{i,\qbn}}_1 \leq \NcutSmall \,  j_*(i) \, \Gamma_{q+1}^{-3 i +\CLebesgue - 1} \,.
\end{align*}
In order to conclude the proof of the Lemma, we use that $\NcutSmall$ is a constant independent of $q$, and that by \eqref{eq:orange:orangutan:2} and \eqref{eq:imax:upper:lower} we have
\begin{align*}
j_*(i) \leq i \frac{\log \Gamma_\qbn}{\log \Gamma_{\qbn-1}} \leq \imax(\qbn-1) b \leq \frac{\badshaq+12}{(b-1)\varepsilon_\Ga} b \,.
\end{align*}
Thus $j_*(i)$ is also bounded from above by a constant independent of $q$, and upon taking $a$ sufficiently large we conclude the proof.
\end{proof}

\subsection{Verifying Eqn.~\texorpdfstring{\eqref{eq:inductive:timescales}}{eqinductive}}\label{ss:dodging:verification}
The following lemma verifies the inductive assumption \eqref{eq:inductive:timescales} at level $q'=q+\bn$.
\begin{lemma}[\bf Overlapping and timescales]\label{lem:overlap:timescales}
Let $q''\in\{q+1,\dots,q+\bn\}$. Assume that $\psi_{i,q+\bn}\psi_{i'',q''}\not\equiv 0$. Then it must be the case that $\tau_{q+\bn}\Gamma_{q+\bn}^{-i}\leq \tau_{q''}\Gamma_{q''}^{-i''-{25}}$.
\end{lemma}
\begin{proof}[Proof of Lemma~\ref{lem:overlap:timescales}]
We split the proof into two steps. In the first step, we prove the claim for $q''=q+\bn-1$, while in the second step we prove the claim for the remaining cases.
\smallskip

\noindent\textbf{Step 1:} We must prove that if $\psi_{i,q+\bn}\psi_{i'',q+\bn-1}\not \equiv 0$, then $\tau_{q+\bn}\Gamma_{q+\bn}^{-i}\leq \tau_{q+\bn-1}\Gamma_{q+\bn-1}^{-i''-{25}}$. By \eqref{eq:psi:i:q:recursive}, if $\psi_{i,q+\bn}(t,x)\neq 0$, then there exists $\vec{i}=(i_0,\dots,i_{\NcutSmall})$ such that $\max_m i_m = i$, and $\psi_{m,i_m,q+\bn}\neq 0$ for all $0\leq i \leq \NcutSmall$.  By \eqref{eq:psi:m:im:q:def} and Definition~\eqref{def:istar:j}, for each $i_m$ there exists a corresponding $j_m$ such that $\psi_{j_m,q+\bn-1}(t,x) \neq 0$ and $\Gamma_{q+\bn}^{i_m} \geq \Gamma_{q+\bn-1}^{j_m}$. From \eqref{eq:inductive:partition} and \eqref{ineq:tau:q}, it then follows that if $\psi_{m,i_m,q+\bn} \psi_{j',q+\bn-1} \neq 0$, then
$$  \tau_{q+\bn}\Gamma_{q+\bn}^{-i_m} \leq \tau_{q+\bn-1}\Gamma_{q+\bn-1}^{-j'-40} \, .  $$
Then \eqref{eq:psi:i:q:recursive} gives that if $\psi_{i,q+\bn}\psi_{i'',q+\bn-1}\not\equiv 0$, 
$$  \tau_{q+\bn}\Gamma_{q+\bn}^{-i} \leq \tau_{q+\bn-1}\Gamma_{q+\bn-1}^{-i''-30} \, .  $$
\smallskip

\noindent\textbf{Step 2:} Suppose that $q'' \leq q+\bn-2$ and that $\psi_{i,q+\bn}(t,x)\psi_{i'',q''}(t,x)\neq 0$.  Then from \eqref{eq:inductive:partition}, there exists $j$ such that $\psi_{i,q+\bn}(t,x) \psi_{j,q+\bn-1}(t,x) \psi_{i'',q''}(t,x)\neq 0$. Applying the result of Step 1 in combination with the inductive assumption \eqref{eq:inductive:timescales} concludes the proof.
\end{proof}

\opsection{Velocity increment potential}\label{opsection:vel:inc:pot}

In order to analyze certain current errors (see for example \cite[Lemma~6.13]{GKN23}), it will be necessary to write the mollified velocity increment $\hat w_\qbn$ as the iterated Laplacian of a potential.  We first carry out this construction for $w_{q+1}$ in the first subsection, as well as construct a pressure increment which dominates the resulting velocity increment potential and analyze its associated pressure current error.  Then in subsection~\ref{op:wed}, we analyze the mollified velocity increment potential, which completes the bulk of the work required to verify the inductive assumptions in subsubsection~\ref{sec:ind.vel.inc.pot}.  Finally, in subsection~\ref{op:wed:2} we prove a lemma which allows us to verify~\eqref{eq:psi:q:q'} at level $\qbn$ in \cite[Lemma~6.8]{GKN23}.\index{velocity increment potential}

\opsubsection{Defining the velocity increment potential}\label{sec:vel.inc}

In this section, we define a potential for $w_{q+1}$ along with an error term, construct its pressure increment and the associated current errors, and investigate their properties.

\begin{lemma*}[\bf Velocity increment potential]\label{lem:rep:vel:inc:potential}
For a given $w_{q+1}^{(\texttt{l})}$, $\texttt{l}=p,c$, as in \eqref{defn:w}, there exists a tensor $\upsilon_{q+1}^{(\texttt{l})}$ and an error $e_{q+1}^{(\texttt{l})}$ such that the following hold.
\begin{enumerate}[(i)] 
    \item Let $\dpot$ be as in \eqref{i:par:10}. Then\index{$\upsilon_{q+1}$}\index{$\upsilon^{(\texttt{l})}_{q+1}$} $w_{q+1}^{(\texttt{l})}$ can be written in terms of $\upsilon_{q+1}^{(\texttt{l})}$ and $e_{q+1}^{(\texttt{l})}$ as
     \begin{equation}\label{exp.w.q+1}
        \begin{split}
w_{q+1}^{({p})} &= \div^{\dpot} \upsilon_{q+1}^{({p})} + e_{q+1}^{({p})} \\
w_{q+1}^{({c})}&= \div^{\dpot}(r_q\Ga_q^{-1} \upsilon_{q+1}^{({c})}) + r_q\Ga_q^{-1}e_{q+1}^{({c})} \, , 
\end{split}
\end{equation}
or equivalently notated component-wise as $(w_{q+1}^{(p)})^\bullet = 
\pa_{i_1} \dots \pa_{i_\dpot} \upsilon_{q+1}^{(p, \bullet, i_1, \dots, i_\dpot)} + e_{q+1}^\bullet
$. 

\item $\upsilon_{q+1}^{(\texttt{l})}$ and $e_{q+1}^{(\texttt{l})}$ have the support property\footnote{For any smooth set $\Omega\subset\mathbb{T}^3$, we use $\Omega\circ\Phiik$ to denote the set $\Phiik^{-1}(\Omega)\subset\mathbb{T}^3\times\R$, i.e. the space-time set whose characteristic function is annihilated by $\Dtq$.}
\begin{align}
\supp(\upsilon_{q+1}^{(\texttt{l})}) &, \, \supp(e_{q+1}^{(\texttt{l})}) \notag\\ &\subseteq \bigcup_{\xi,i,j,k,\vecl,I,\diamond} \supp \left(\chi_{i,k,q} \zeta_{q,\diamond,i,k,\xi,\vecl} \left(\rhob_{(\xi)}^\diamond \zetab_{\xi}^{I,\diamond} \right)\circ \Phiik \right) \cap
    B\left( \supp \varrho^I_{(\xi),\diamond} , 2 \lambda_\qbn^{-1} \right)\circ \Phiik \, .  \label{supp.upsilon.e.lem}
\end{align}

\item\label{item:eckel:1} For $0\leq k \leq \dpot$,  $(\upsilon_{q+1,k}^{(\texttt{l})})^\bullet :=\la_{q+\bn}^{\dpot-k}\partial_{i_1}\cdots \partial_{i_k}
    \upsilon_{q+1}^{(\texttt{l}, \bullet,i_1,\dots,i_\dpot)}$,\footnote{If $k=0$, we adopt the convention that $\partial_{i_1}\cdots\partial_{i_k}$is the identity operator.} satisfies the estimates
\begin{subequations}
\begin{align}
    &\norm{\psi_{{i,q}} D^N {D_{t,q}^{M}} 
   \upsilon_{q+1,k}^{(\texttt{l})}
    }_3\leq \Ga_q^{10} \de_{q+\bn}^\frac12 r_{q}^{-\frac13}
     \la_{q+\bn}^N
     \MM{M, \Nindt, \Ga_{q}^{i+14} \tau_{q}^{-1}, \Ga_{q}^{8}\Tau_{q}^{-1}}
    \label{est.new.upsilon.3}
    \\
    &\norm{\psi_{{i,q}} D^N {D_{t,q}^{M}} 
   \upsilon_{q+1,k}^{(\texttt{l})}
    }_\infty \leq 
    \Ga_q^{\frac{\badshaq}2+ 10} r_q^{-1}
    \la_{q+\bn}^{N}
     \MM{M, \Nindt, \Ga_{q}^{i+14} \tau_{q}^{-1}, \Ga_{q}^{8}\Tau_{q}^{-1}}
         \label{est.new.upsilon.inf}
\end{align}
\end{subequations}
for $N\leq \sfrac{\Nfin}4-2\dpot^2$ and $M\leq \sfrac{\Nfin}5$.

\item $e_{q+1}^{(\texttt{l})}$ satisfies
\begin{align} 
    \norm{ D^N D_{t,q}^{M} e_{q+1}^{(\texttt{l})}}_{\infty} 
    \leq \delta_{q+3\bn}^3\Tau_{\qbn}^{20\Nindt}\la_{q+\bn}^{-10}
    \la_{q+\bn}^{N}
     \MM{M, \Nindt,  {\tau_{q}^{-1}}, \Ga_{q}^{8}\Tau_{q}^{-1}}
    \label{est.new.e.inf}\, .
\end{align}
for $N\leq \sfrac{\Nfin}4-2\dpot^2$ and $M\leq \sfrac{\Nfin}5$.
\end{enumerate}

\end{lemma*}

\begin{remark*}[\bf Notation for cumulative velocity increment potential]\label{rem:rep:vel:inc:potential}
We let $\upsilon_{q+1} : = \upsilon_{q+1}^{(p)} + r_q\Ga_q^{-1} \upsilon_{q+1}^{(c)}$ and $\upsilon_{q+1,k}^\bullet :=\la_{q+\bn}^{\dpot-k}\partial_{i_1}\cdots \partial_{i_k}
    \upsilon_{q+1}^{( \bullet,i_1,\dots,i_\dpot)}$. As a corollary of Lemma \ref{lem:rep:vel:inc:potential}, we have that
\begin{align*}
    w_{q+1} = \div^{\dpot} \upsilon_{q+1} + e_{q+1} \, ,
\end{align*}
where $\upsilon_{q+1}$ and $e_{q+1}$ share the properties \eqref{supp.upsilon.e.lem}--\eqref{est.new.e.inf} with $\upsilon_{q+1}^{(\texttt{l})}$ and $e_{q+1}^{(\texttt{l})}$ after adjusting the inequalities to include an implicit constant.     
\end{remark*}

\begin{proof} Recall from subsection~\ref{ss:corr-ec-tor} that $w_{q+1} = w_{q+1,R} + w_{q+1,\ph}$ where 
\begin{align}
    w_{q+1,\diamond} 
    &= \sum_{i,j,k,\xi,\vecl,I} a_{(\xi),\diamond} \nabla\Phi_{(i,k)}^{-1} (\rhob_{(\xi)}^\diamond \etab^{I,\diamond}_\xi)\circ \Phiik \WW_{(\xi),\diamond}^{I} \circ \Phi_{(i,k)}  \label{exp.w.diamond1}\\ 
    &\qquad + \sum_{i,j,k,\xi,\vecl,I} \nabla \left( (\rhob_{(\xi)}^\diamond \zetab^{I,\diamond}_\xi)\circ \Phiik a_{(\xi),\diamond} \right) \times \left( \nabla\Phiik \UU_{(\xi),\diamond}^{I}\circ \Phiik \right)\label{exp.w.diamond2}
\end{align}
for $\diamond = R, \ph$. To construct $ \upsilon_{q+1}$ and $e_{q+1}$, we will apply Corollary \ref{cor:inverse.div} to the right hand side terms.  We shall adhere to the convention set out in Remark~\ref{rem:scalar:inverse:div} and treat each component separately, so that the resulting tensor potential does not have any special symmetry properties.

Fix values for all indexes $i,j,k,\xi,\vecl,I$, set $\diamond = R$, and consider one component, indexed by $\bullet$, of the vector field in \eqref{exp.w.diamond1}. Set 
\begin{align*}
p&=3,\infty \, , \quad N_* = \sfrac{\Nfin}4,  \quad 
M_*=\sfrac{\Nfin}5, \quad M_t = \Nindt,\quad\\
G&=a_{(\xi),R} \nabla\Phi_{(i,k)}^{-1} (\rhob_{(\xi)}^R \etab^{I,R}_\xi)\circ \Phiik \xi^\bullet r_q^{-\sfrac 13} \, , \quad
\Phi = \Phi_{(i,k)}, 
\quad \pi = \pi_\ell \Ga_q^{30},  \quad   r_G = r_q \\
\const_{G,p} &=\left| \supp \left(\eta_{i,j,k,\xi,\vecl,R}\zetab_\xi^{I,R} \right) \right|^{\sfrac 1p} \de_{q+\bn}^\frac12 \Ga_q^{j+7},\quad
\la = \la_{q+\half}, \quad \la' = {\lambda_q \Ga_q},
\quad \nu = \tau_q^{-1}\Ga_q^{i+13},
\quad  \nu' = \Tau_q^{-1} \Ga_q^8, 
\\
\varrho&=r_q^{\sfrac 13} \varrho_{\xi,\lambda_{q+\bn},\frac{\lambda_{q+\lfloor\sfrac \bn 2 \rfloor}\Gamma_q}{\lambda_{q+\bn}}}, \quad
\td\vartheta=\td\vartheta_{\xi,\lambda_{q+\bn},\frac{\lambda_{q+\lfloor\sfrac \bn 2 \rfloor}\Gamma_q}{\lambda_{q+\bn}},R}\quad\\
\const_{*,3} 
&=1 \, , \quad
\const_{*,\infty}
={r_q^{-\sfrac 23}},
\quad \mu = \la_{q+\half}\Ga_{q},
\quad \Upsilon = \Upsilon' = \Lambda = \la_{q+\bn},
\end{align*}
where $\td\vartheta$ is constructed from Proposition~\ref{prop:pipeconstruction} with $\Dpot=\dpot^2$. Then, all assumptions of Corollary \ref{cor:inverse.div} hold by \eqref{condi.Nfin0}, \eqref{e:a_master_est_p_R:zeta}, \eqref{e:a_master_est_p_R_pointwise}, \eqref{eq:nasty:D:vq:old}, Corollary \ref{cor:deformation}, \eqref{e:Phi}, and Proposition~\ref{prop:pipeconstruction}.  Then from \eqref{divd.expression}, there exist $ R=:\upsilon_{(\xi),I, R}^{(p)}$ and $E=:e_{(\xi),I, R}^{(p)}$ such that
\begin{align*}
    a_{(\xi),R} \nabla\Phi_{(i,k)}^{-1} (\rhob_{(\xi)}^R \etab^{I,R}_\xi)\circ \Phiik \WW_{(\xi),R}^{I} \circ \Phi_{(i,k)}= \div^{\dpot}  \upsilon_{(\xi),I, R}^{(p)}
    + e_{(\xi),I, R}^{(p)} \, .
\end{align*}
From \eqref{eq:inverse:divd:stress:1}, we have that
\begin{subequations}
\begin{align}
    \norm{D^{N} D_{t,q}^M \pa_{i_1}\cdots\pa_{i_l} \upsilon_{(\xi),I, R}^{(p)}}_{3}
&\leq
    \norm{D^{N} D_{t,q}^M D^l \upsilon_{(\xi),I, R}^{(p)}}_{3} \nonumber\\
 &\les  \left| \supp \left(\eta_{(\xi),R}\zetab_\xi^{I,R} \right) \right|^{\sfrac 13} \de_{q+\bn}^\frac12 \Ga_q^{j+7} r_q^{-\sfrac13}   \la_{q+\bn}^{l-\dpot+N} \MM{M,M_{t},\tau_q^{-1}\Ga_q^{i+13},\Tau_q^{-1} \Ga_q^8} \, , 
 \label{new.upsilon.3.piece} \\
\norm{D^{N} D_{t,q}^M \pa_{i_1}\cdots\pa_{i_l} \upsilon_{(\xi),I, R}^{(p)}}_{\infty}
 &\leq
 \norm{D^{N} D_{t,q}^M D^l \upsilon_{(\xi),I, R}^{(p)}}_{\infty}\nonumber\\
 &\les \de_{q+\bn}^\frac12 \Ga_q^{j+7} r_q^{-1}   \la_{q+\bn}^{l-\dpot} \la_{q+\bn}^{N+\alpha} \MM{M,M_{t},\tau_q^{-1}\Ga_q^{i+13},\Tau_q^{-1}\Ga_q^8}\nonumber\\ 
 &\les \Ga_q^{\sfrac{(\badshaq+20)}{2}} r_q^{-1}   \la_{q+\bn}^{l-\dpot} \la_{q+\bn}^{N+\alpha} \MM{M,M_{t},\tau_q^{-1}\Ga_q^{i+3},\Tau_q^{-1}}\,,  \label{new.upsilon.infty.piece} 
\end{align}
\end{subequations}
for $0\leq l \leq \dpot$, $N+l\leq \sfrac{\Nfin}4-\dpot^2$,  and $M\leq \sfrac{\Nfin}5$, where we used \eqref{ineq:jmax:use} in the last inequality.  From \eqref{eq:inverse:divd:error:1}, we have that 
\begin{align}\label{est.e.piece}
    \norm{D^N D_{t,q}^M e_{(\xi),I, R}^{(p)}}_{\infty}  
\les \de_{q+\bn}^\frac12 \Ga_q^{j+7} r_q^{-1}  \left( \sfrac{\la_{q+\half}}{\la_{q+\bn}}\right)^{ \dpot} \la_{q+\bn}^{N+\alpha} \MM{M,M_{t},\tau_q^{-1}\Ga_q^{i+13},\Tau_q^{-1}\Ga_q^8}\, .
\end{align}
for $N\leq \sfrac{\Nfin}4-\dpot^2$,  and $M\leq \sfrac{\Nfin}5$.
Furthermore, from \eqref{eq:inverse:divd:linear} and \eqref{item:pipe:3.5} from  Proposition~\ref{prop:pipeconstruction}, we have that the supports of $\upsilon_{(\xi),I, R}^{(p)}$ and $e_{(\xi),I, R}^{(p)}$ are contained in the set on the right-hand side of \eqref{supp.upsilon.e.lem}.

We now sum over indexes $i,j,k,\xi,\vecl,I$ and set
\begin{align}\label{defn:upsilon.p.R}
 \upsilon_{q+1,R}^{(p)}
&= \sum_{i,j,k,\xi,\vecl,I}
 \upsilon_{(\xi),I, R}^{(p)} \, , \qquad
e_{q+1,R}^{(p)}
= \sum_{i,j,k,\xi,\vecl,I}
e_{(\xi),I, R}^{(p)} \, ,   
\end{align}
which verifies the first equality in \eqref{exp.w.q+1} and \eqref{supp.upsilon.e.lem}. Using \eqref{eq:desert:cowboy:sum} to obtain an $L^\infty$ bound for the sum and Corollary~\ref{rem:summing:partition} with $H_{i,j,k,\xi,\vecl,R} = \upsilon_{(\xi),I, R}^{(p)}, \theta_2 = \theta = 1, p=3 , \const_H = \de_{q+\bn}^\frac12 \Ga_q^{7} r_q^{-1} , N_x = N_\ast = \sfrac{\Nfin}4-\dpot^2$, the obvious choices for the other parameters, \eqref{new.upsilon.3.piece},  \eqref{new.upsilon.infty.piece},  \eqref{est.e.piece}, and \eqref{ineq:dpot:1}, we have that $ \upsilon_{q+1,R}^{(p)}$ and $e_{q+1,R}^{(p)}$ satisfy 
\begin{align*}
\norm{\psi_{i,q}D^N D_{t,q}^M \partial_{i_1}\dots \partial_{i_k}
    ({\upsilon_{q+1, R}^{(p)}})
    ^{(i_1,\dots, i_\dpot)}}_{3}
 &\les \Ga_q^{10} \de_{q+\bn}^\frac12 r_q^{-\sfrac13}   \la_{q+\bn}^{k-\dpot} \la_{q+\bn}^{N+\alpha} \MM{M,M_{t},\tau_q^{-1}\Ga_q^{i+14},\Tau_q^{-1}\Ga_q^8} \\
 \norm{\psi_{i,q}D^N D_{t,q}^M \partial_{i_1}\dots \partial_{i_k}
    ({\upsilon_{q+1, R}^{(p)}})
    ^{(i_1,\dots, i_\dpot)}}_{\infty}
 &\les \Ga_q^{\frac{\badshaq}2+10}r_q^{-1}   \la_{q+\bn}^{k-\dpot} \la_{q+\bn}^{N+\alpha} \MM{M,M_{t},\tau_q^{-1}\Ga_q^{i+14},\Tau_q^{-1}\Ga_q^8}\, ,  \\
   \norm{D^N D_{t,q}^M e_{q+1,R}^{(p)}}_{\infty}  
&\lec \delta_{q+3\bn}^3\Tau_{\qbn}^{2\Nindt}\la_{q+\bn}^{-10} \la_{q+\bn}^{N} \MM{M,M_{t},\tau_q^{-1}\Ga_q^{i+14},\Tau_q^{-1}\Ga_q^8}  
\end{align*}
for $N\leq \sfrac{\Nfin}4-\dpot^2$,  and $M\leq \sfrac{\Nfin}5$.
The first inequality follows from Lemma~\eqref{lemma:cumulative:cutoff:Lp} and Remark \ref{rem:summing:partition}, and the second and the last inequalities use the support property noted earlier.

In a similar way, we work on 
\eqref{exp.w.diamond1} with $\ph$ and
\eqref{exp.w.diamond2} with $R, \ph$ and generate $( \upsilon _{q+1, \ph}^{(p)}, e_{q+1, \ph}^{(p)})$, $( \upsilon _{q+1, R}^{(c)},e_{q+1, R}^{(c)})$, and $( \upsilon _{q+1, \ph}^{(c)}, e_{q+1, \ph}^{(c)})$, respectively. Indeed, for \eqref{exp.w.diamond1} with $\ph$, we set
\begin{align*}
    G= a_{(\xi),\ph} \nabla\Phi_{(i,k)}^{-1} (\rhob_{(\xi)}^\ph \etab^{I,\ph}_\xi)\circ \Phiik \xi , \quad 
    \varrho = \varrho_{\xi,\lambda_{q+\bn},\frac{\lambda_{q+\lfloor\sfrac \bn 2 \rfloor}\Gamma_q}{\lambda_{q+\bn}},{\ph}}, \quad
\td\vartheta=r_q^{-\sfrac13}\td\vartheta_{\xi,\lambda_{q+\bn},\frac{\lambda_{q+\lfloor\sfrac \bn 2 \rfloor}\Gamma_q}{\lambda_{q+\bn}},\ph}
\end{align*}
where $\td\vartheta$ is constructed from Proposition~\ref{prop:pipe.flow.current} with $\Dpot=\dpot^2$, and choose the rest of parameters and functions as in the case $\diamond=R$. The rest of the conclusions follow analogously to the case $\diamond=R$, and we omit further details. In the case of \eqref{exp.w.diamond2}, we write 
\begin{align*}
(w_{(\xi),\diamond}^{(c),I})^{\bullet} = r_q\Ga_q^{-1} G_\diamond (\varrho_\diamond \circ \Phi) \, ,     
\end{align*}
where $G_\diamond$ and $\varrho_\diamond$ are defined by
\begin{align*}
G_R = \la_{q+\half}^{-1}\epsilon_{\bullet pr} \partial_p \left( a_{(\xi),R} \left( \rhob_{(\xi)}^R \zetab_\xi^{I,R} \right)\circ \Phiik\right)  \partial_r \Phi_{(i,k)}^s , \quad
\varrho_R &= \la_{q+\bn}(\mathbb{U}_{(\xi),R}^{I})^s, \quad \Phi = \Phi_{(i,k)}\\
G_\ph = r_q^{\sfrac 13}\la_{q+\half}^{-1}\epsilon_{\bullet pr} \partial_p \left( a_{(\xi),\ph} \left( \rhob_{(\xi)}^\ph \zetab_\xi^{I,\ph} \right)\circ \Phiik\right)  \partial_r \Phi_{(i,k)}^s , \quad
\varrho_\ph &= r_q^{-\sfrac13}\la_{q+\bn}(\mathbb{U}_{(\xi),\ph}^{I})^s, \quad \Phi = \Phi_{(i,k)}  \, .
\end{align*}
Due to the rescaling by $r_q\Ga_q^{-1}$, we may apply Corollary \ref{cor:inverse.div} to $(r_q\Ga_q^{-1})^{-1}(w_{(\xi),\diamond}^{(c),I})^{\bullet}$ with the same choice of parameters as in the case $\texttt{l}=p$. As a consequence, we obtain $( \upsilon _{q+1, \diamond}^{(c)},e_{q+1, \diamond}^{(c)})$, defined as in \eqref{defn:upsilon.p.R}, which enjoy the same properties as $( \upsilon _{q+1, R}^{(p)},e_{q+1, R}^{(p)})$. Note that from the construction, the velocity increment potential associated to the correctors satisfies
\begin{align*}
    w_{q+1, \diamond}^{(c)} =  \div^{\dpot}(r_q\Ga_q^{-1} \upsilon_{q+1, \diamond}^{(c)}) + r_q\Ga_q^{-1}e_{q+1, \diamond}^{(c)} \, . 
\end{align*}
We may now set
\begin{align*}
    \upsilon _{q+1}
 &=\sum_{\diamond= R,\ph} \upsilon _{q+1, \diamond}^{(p)}
 + r_q\Ga_q^{-1}\upsilon _{q+1, \diamond}^{(c)}
 =: \upsilon _{q+1}^{(p)}+ r_q\Ga_q^{-1}\upsilon _{q+1}^{(c)}
 \\
 e_{q+1}
 &=\sum_{\diamond= R,\ph}e _{q+1, \diamond}^{(p)}
 +r_q\Ga_q^{-1}e_{q+1, \diamond}^{(c)}
 =: e _{q+1}^{(p)}+ r_q\Ga_q^{-1}e _{q+1}^{(c)}\,  .
\end{align*}
which leads to \eqref{supp.upsilon.e.lem}, \eqref{est.new.upsilon.3}, \eqref{est.new.upsilon.inf}, and \eqref{est.new.e.inf}.
\end{proof}

\begin{remark*}[\bf Decompositions of potentials into pieces to facilitate pressure creation]\label{rem:pre.pr.vel.inc}
From the proof of Lemma~\ref{lem:rep:vel:inc:potential}, the velocity increment potentials $\upsilon_{q+1,k}^{(\texttt{l})}$, $\texttt{l} = p, c$, $k=0, \cdots, \dpot$, have the additional properties listed below.
\begin{enumerate}[(i)]
    \item Using  Corollary~\ref{cor:inverse.div}, \eqref{item:divd:local:ii}, we have that $\upsilon_{q+1,\dpot}^{(\texttt{l})}=\lambda_\qbn^{\dpot}\upsilon_{q+1}^{(\texttt{l})}$ can be decomposed as
    \begin{align}
    \upsilon_{q+1,\dpot}^{(\texttt{l})}
    &= \la_\qbn^\dpot \sum_{i,j,k,\xi,\vecl,I, \diamond}\sum_{\ovj = 0}^{\ov \const_\divH} H^{\al(\ov j)}_{{(\xi), I,\diamond}}(\rho^{\beta(\ov j)}_{{(\xi), I,\diamond}}  \circ\Phi_{(i,k)})  \notag \\
    &=: \sum_{\pxi,I, \diamond} H_{{(\xi), I,\diamond}} \rho_{{(\xi), I,\diamond}}  \circ\Phi_{(i,k)} \label{rep:vel:inc:pot}
    \end{align}
where we abuse notation slightly by using $\pxi$ to include the indices $i,j,k,\xi,\vecl,\ovj$ as well as the indices in $\alpha(\ov j)$ or $\beta(\ovj)$ in the final expression, which take a finite number of values independent of $q$.
\item  Let $p=3$ or $\infty$. $H_{\pxiid}$ satisfies 
\begin{subequations}
\begin{align}
\supp H_{\pxiid} 
&\subseteq \supp \left((\rhob_{(\xi)}^\diamond \etab^{I,\diamond}_\xi)\circ \Phiik\right)  \, ,  
\label{eq:inverse:divd:linear.ap.H}\\
\left\|\prod_{i=1}^{k}D^{\alpha_i} D_{t,q}^{\be_i} H_{\pxiid} \right\|_{p} &\les \left| \supp \left(\eta_{i,j,k,\xi,\vecl,\diamond}\zetab_\xi^{I,\diamond} \right) \right|^{\sfrac 1p} \de_{q+\bn}^{\sfrac 12} \Ga_q^{j+7} r_q^{-\sfrac 13} \notag\\
&\qquad \qquad \times \la_{q+\half}^{|\al|} \MM{|\be|,\Nindt,\tau_q^{-1}\Ga_q^{i+13},\Tau_q^{-1}\Ga_q^8} \, , \label{eq:inverse:divd:sub:2.ap} \\
\left|\prod_{i=1}^{k}D^{\alpha_i} D_{t,q}^{\be_i} H_{\pxiid} \right|
&\lec (\pi_\ell\Ga_q^{30})^{\sfrac 12} r_q^{-\sfrac 13} \la_{q+\half}^{|\al|} \MM{|\be|,\Nindt,\tau_q^{-1}\Ga_q^{i+13},\Tau_q^{-1}\Ga_q^8} \, , \label{eq:invser:divd:point.ap}
\end{align}
\end{subequations}
for all integer $k\geq 1$ and multi-indices $\al, \be\in \mathbb{N}^k$ with $|\al|\leq \sfrac{\Nfin}4-\dpot^2$ and $|\be|\leq \sfrac{\Nfin}5$.
\item $\rho_{\pxiid}$ is $(\sfrac{\mathbb{T}}{\la_{q+\half}}\Ga_q)^3$-periodic and satisfies 
\begin{subequations}
\begin{align}
\supp \rho_{\pxiid} &\subseteq \supp \left(\td\vartheta_{\xi,\lambda_{q+\bn},\frac{\lambda_{q+\lfloor\sfrac \bn 2 \rfloor}\Gamma_q}{\lambda_{q+\bn}},\diamond}\right) \, \label{eq:inverse:divd:linear.ap}\\
\left\| D^N \rho_{\pxiid} \right\|_{L^p} &\les r_q^{\frac2p-\frac 23} \la_{q+\bn}^{N}  \label{eq:inverse:divd:sub:1:ap}
\end{align}
\end{subequations}
for all $N \leq \sfrac{\Nfin}4-\dpot^2$ and $(\pxiid)$.
\end{enumerate}
These properties of $H_{{(\xi), I, \diamond}}$ and $\rho_{\pxiid}$ follow from items~\eqref{item:divd:local:i}--\eqref{item:divd:local:iv}. 

From the above properties, we may derive similar formulae and properties for \emph{all} of the various velocity increment potentials $\upsilon_{q+1,h}^{\pttl}$ defined in item~\eqref{item:eckel:1} for $0\leq h \leq \dpot$.  Specifically, we have that $\upsilon_{q+1,h}^{\pttl}$ can be decomposed using \eqref{rep:vel:inc:pot} and the Leibniz rule\footnote{We use the notation $$ \partial_{i_1}\cdots \partial_{i_h} (f g ) = \sum_{\substack{\vec{a}_h=(a_1,\dots,a_A),  \\ \vec{b}_h=(b_1,\dots,b_B)}} C_{\vec{a}_h,\vec{b}_h} \partial_{i_{a_1}} \cdots \partial_{i_{a_A}} f \,  \partial_{i_{b_1}} \cdots \partial_{i_{b_A}} g = \sum_{\vec{a}_h,\vec{b}_h} C_{\vec{a}_h,\vec{b}_h} \partial_{\vec{a}_h} f \partial_{\vec{b}_h} g \, , $$ where $\vec{a}_h,\vec{b}_h$ are multi-indices with $A$, respectively $B$ distinct components for which the union of all indices belonging to either $\vec{a}_h$ or $\vec{b}_h$ is $\{i_1,\dots,i_h\}$.} as
\begin{align}
    \upsilon_{q+1,h}^{(\texttt{l},\bullet,i_{h+1}, \cdots, i_\dpot)}
    &= \la_{q+\bn}^{\dpot-h}\partial_{i_1}\cdots \partial_{i_h}
    \upsilon_{q+1}^{(\texttt{l}, \bullet,i_1,\dots,i_\dpot)} \notag \\
    &=\la_{q+\bn}^{\dpot-h} \sum_{\vec{a}_h,\vec{b}_h} \const_{\vec{a}_h,\vec{b}_h} \sum_{i,j,k,\xi,\vecl,I,\diamond} \sum_{\ovj=0}^{\ov \const_{\divH}} \pa_{\vec{a}_h} H^{\al(\ov j)}_{{ (\xi), I,\diamond}} \pa_{\vec{b}_h} \left(\rho^{\beta(\ov j)}_{{ (\xi), I,\diamond}}  \circ\Phi_{(i,k)} \right) \notag \\
    &=: \sum_{\pxiid,h'} H_{\pxiid}^{h,h'} \rho_{\pxiid}^{h,h'} \circ \Phiik \notag \\
    &=: \sum_{\pxiid,h'}  \Upsilon^{h,h'}_{(\xi),I,\diamond} \, , \label{rep:vel:inc:pot:eckel}
\end{align}
where $H_{\pxiid}^{h,h'}$, $\rho_{\pxiid}^{h,h'}$, and $\Upsilon^{h,h'}_{(\xi),I,\diamond}$ satisfy the following, and we again abuse notation slightly by letting $\pxi$ denote all indices $i,j,k,\xi,\vecl,\ovj$, as well as those indices needed for the application of the Faa di Bruno formula from \eqref{eq:Faa:di:Bruno:2} to $\pa_{\vec{b}_h} \left(\rho^{\beta(\ov j)}_{{ (\xi), I,\diamond}}  \circ\Phi_{(i,k)} \right)$. We again have that $\pxi$ includes $i,j,k,\xi,\vecl,\xi$, as well as the a finite, $q$-independent number of indices.  
\begin{enumerate}[(i)]
    \item  Let $p=3$ or $\infty$. $H_{\pxiid}^{h,h'}$ satisfies 
\begin{subequations}
\begin{align}
\supp H^{h,h'}_{\pxiid} 
&\subseteq \supp \left((\rhob_{(\xi)}^\diamond \etab^{I,\diamond}_\xi)\circ \Phiik\right)  \, ,  
\label{eq:inverse:divd:linear.ap.H.eckel}\\
\left\|\prod_{i=1}^{k}D^{\alpha_i} D_{t,q}^{\be_i} H^{h,h'}_{\pxiid} \right\|_{p} &\les \left| \supp \left(\eta_{i,j,k,\xi,\vecl,\diamond}\zetab_\xi^{I,\diamond} \right) \right|^{\sfrac 1p} \de_{q+\bn}^{\sfrac 12} \Ga_q^{j+7} r_q^{-\sfrac 13} \notag\\
&\qquad \qquad \times \la_{q+\half}^{|\al|} \MM{|\be|,\Nindt,\tau_q^{-1}\Ga_q^{i+13},\Tau_q^{-1}\Ga_q^8} \, , \label{eq:inverse:divd:sub:2.ap.eckel} \\
\left|\prod_{i=1}^{k}D^{\alpha_i} D_{t,q}^{\be_i} H^{h,h'}_{\pxiid} \right|
&\lec (\pi_\ell\Ga_q^{30})^{\sfrac 12} \la_{q+\half}^{|\al|} \MM{|\be|,\Nindt,\tau_q^{-1}\Ga_q^{i+13},\Tau_q^{-1}\Ga_q^8} \, , \label{eq:invser:divd:point.ap.eckel}
\end{align}
\end{subequations}
for all integer $k\geq 1$ and multi-indices $\al, \be\in \mathbb{N}^k$ with $|\al|\leq \sfrac{\Nfin}4-2\dpot^2$ and $|\be|\leq \sfrac{\Nfin}5$.
\item $\rho_{\pxiid}^{h,h'}$ is $(\sfrac{\mathbb{T}}{\la_{q+\half}}\Ga_q)^3$-periodic and satisfies 
\begin{subequations}
\begin{align}
\supp \rho_{\pxiid}^{h,h'} &\subseteq \supp \left(\td\vartheta_{\xi,\lambda_{q+\bn},\frac{\lambda_{q+\lfloor\sfrac \bn 2 \rfloor}\Gamma_q}{\lambda_{q+\bn}},\diamond}\right) \, \label{eq:inverse:divd:linear.ap.eckel}\\
\left\| D^N \rho_{\pxiid}^{h,h'} \right\|_{L^p} &\les r_q^{\frac2p-\frac 23} \la_{q+\bn}^{N}  \label{eq:inverse:divd:sub:1.eckel}
\end{align}
\end{subequations}
for all $N \leq \sfrac{\Nfin}4-2\dpot^2$ and $(\pxiid)$.
\item For $p=3,\infty$, we have that
\begin{align}\label{eq:inverse:divd:linear:ap:total:eckel}
    \left\| \Upsilon_\pxiid^\hhp \right\|_p \les \left| \supp \left(\eta_{i,j,k,\xi,\vecl,\diamond}\zetab_\xi^{I,\diamond} \right) \right|^{\sfrac 1p} \de_{q+\bn}^{\sfrac 12} \Ga_q^{j+7} r_q^{\sfrac 2p -1} \, .
\end{align}
\end{enumerate}
The proofs of these properties follows from backwards induction on the index $h$.  Indeed, the case $h=\dpot$ has already been shown in the beginning of the remark.  The subsequent cases follow from application of the Faa di Bruno formula to \eqref{rep:vel:inc:pot} to derive \eqref{rep:vel:inc:pot:eckel}, \eqref{eq:inverse:divd:linear.ap.H}--\eqref{eq:inverse:divd:sub:1:ap}, Corollary~\ref{cor:deformation}, and Lemma~\ref{l:slow_fast}.
\end{remark*}

\begin{lemma*}[\bf Pressure increment] \label{lem:pr.inc.vel.inc.pot}
Define $\upsilon_{q+1,k}^{\pttl}$, $0\leq k \leq \dpot$, $\texttt{l}=p,c$, as in Lemma \ref{lem:rep:vel:inc:potential}.
Then there exists a pressure increment\index{$\sigma_\upsilon$} $\si_{\upsilon^{(\texttt{l})}}=\si_{\upsilon^{(\texttt{l})}}^+ - \si_{\upsilon^{(\texttt{l})}}^-$ associated to the sum $\sum_{k=0}^\dpot \upsilon_{q+1,k}^{(\texttt{l})}$ of velocity increment potentials such that the following properties hold.
\begin{enumerate}[(i)]
    \item We have that for all $k=0,1,\dots,\dpot$,
\begin{align}
\left|\psi_{i,q}D^N D_{t,q}^M \upsilon_{q+1,k}^{(\texttt{l})}\right|
    \lec (\si_{\upsilon^{(\texttt{l})}}^+  + \de_{q+3\bn})^{\sfrac12}r_{q}^{-1} (\la_{q+\bn}\Ga_{q+\bn}^{{\sfrac1{10}}})^N\MM{M,\Nindt,\tau_q^{-1}\Ga_q^{i+16},\Tau_q^{-1}\Ga_{q}^9}
    \label{est.vel.inc.pot.by.pr}
\end{align}
for any $0\leq k \leq \dpot$ and $N,M \leq \sfrac{\Nfin}5$. 
\item Set 
\begin{align}\label{defn:si.upsilon}
 \si_{\upsilon}^\pm:= \si_{\upsilon^{(p)}}^\pm
+\si_{\upsilon^{(c)}}^\pm \, , \qquad \si_{\upsilon}= \si_{\upsilon}^+-\si_{\upsilon}^- \, .
\end{align}
Then we have that
\begin{subequations}
\begin{align}
\left|\psi_{i,q}D^N D_{t,q}^M \si_{\upsilon}^+\right|
    &\lec (\si_{\upsilon}^{+}+ \de_{q+3\bn})
    (\la_{q+\bn}\Ga_{q+\bn}^{\sfrac{1}{10}})^N\MM{M,\Nindt,\tau_q^{-1}\Ga_q^{i+16},\Tau_q^{-1}\Ga_{q}^9}\, , \label{est.pr.vel.inc}\\
\norm{\psi_{i,q}D^N D_{t,q}^M\si_{\upsilon}^+}_{\sfrac32}
    &\leq {\Ga_{q+\bn}^{-9}}\de_{q+2\bn}(\la_{q+\bn}\Ga_{q+\bn}^{\sfrac{1}{10}})^N \MM{M,\Nindt,\tau_q^{-1}\Ga_q^{i+16},\Tau_q^{-1}\Ga_{q}^9}\, , \quad \label{est.pr.vel.inc.32}\\
    \norm{\psi_{i,q}D^N D_{t,q}^M\si_{\upsilon}^+}_{\infty}
    &\leq \Ga_{q+\bn}^{\badshaq-9}(\la_{q+\bn}\Ga_{q+\bn}^{\sfrac{1}{10}})^N \MM{M,\Nindt,\tau_q^{-1}\Ga_q^{i+16},\Tau_q^{-1}\Ga_{q}^9}\, , \quad \label{est.pr.vel.inc.infty}\\
     \norm{\psi_{i,q}D^N D_{t,q}^M\si_{\upsilon}^-}_{\sfrac32}&\leq {\Ga_{q+\bn}^{-9}}\de_{q+2\bn}(\la_{q+\half}\Ga_{q+\half})^N\MM{M,\Nindt,\tau_q^{-1}\Ga_q^{i+16},\Tau_q^{-1}\Ga_{q}^9}\,  , \label{est.pr.vel.inc.-32}\\
     \norm{\psi_{i,q}D^N D_{t,q}^M\si_{\upsilon}^-}_{\infty}&\leq \Ga_{q+\bn}^{\badshaq-9}(\la_{q+\half}\Ga_{q+\half})^N\MM{M,\Nindt,\tau_q^{-1}\Ga_q^{i+16},\Tau_q^{-1}\Ga_{q}^9}\,  , \label{est.pr.vel.inc.-infty}\\
    {\left| \psi_{i,q} D^N D_{t,q}^M \si_{\upsilon}^-\right|}
    &\lec  \pi_\ell \Ga_q^{30}  {r_{q}^{\sfrac 43}} (\la_{q+\half}\Ga_{q+\half})^N\MM{M,\Nindt,\tau_q^{-1}\Ga_q^{i+16},\Tau_q^{-1}\Ga_{q}^9}\, .\label{est.pr.vel.inc-}
\end{align}
\end{subequations}
for all $N \leq \sfrac{\Nfin}5$ and $M\leq \sfrac{\Nfin}5-\NcutSmall$. 
\item We have that
\begin{align}
    \supp(\si_{\upsilon}^+) \cap B(\hat w_{q''},\la_{q''}^{-1}\Ga_{q''{+1}}) \, , \quad
    \supp(\si_{\upsilon}^-)\cap B( \hat{w}_{q'}, \la_{q'}^{-1}\Ga_{q'{+1}}) = \emptyset \label{supp:pr:vel.inc}
\end{align}
for $q+1\leq q'' \leq \qbn-1$ and $q+1\leq q'\leq q+\half$. 
\item\label{i:presh:vel:mean} Define 
\begin{equation}\label{def:bmu:vel:presh}
    \bmu_{\sigma_{\upsilon}}(t) = \int_0^t \left \langle \Dtq \sigma_{\upsilon}  \right \rangle (s) \, ds \, .
\end{equation}
Then we have that
    \begin{align}\label{th:billys:3}
      \left|\frac{d^{M+1}}{dt^{M+1}} \bmu_{\sigma_{\upsilon}} \right| 
      \leq (\max(1, T))^{-1}\delta_{q+3\bn}^2 \MM{M,\Nindt,\tau_q^{-1},\Tau_{q+1}^{-1}}
    \end{align}
    for $0\leq M\leq 2\Nind$. 
\end{enumerate}
\end{lemma*}

\medskip

\begin{remark*}[\bf Pointwise bounds for principal and corrector parts]
From \eqref{exp.w.q+1}--\eqref{est.new.e.inf}, \eqref{est.vel.inc.pot.by.pr}, and \eqref{condi.Nfin0}, we have that
\begin{subequations}\label{sunday:morning:1}
\begin{align}
    \left|\psi_{i,q}D^N D_{t,q}^M w_{q+1}^{(p)}\right|
    &\lec (\si_{\upsilon^{(p)}}^+  + \de_{q+3\bn})^{\sfrac12}r_{q}^{-1} (\la_{q+\bn}\Ga_{q+\bn}^{{\sfrac1{10}}})^N\MM{M,\Nindt,\tau_q^{-1}\Ga_q^{i+16},\Tau_q^{-1}\Ga_{q}^9} \, ,
    \label{est.vel.inc.p.by.pr} \\
    \left|\psi_{i,q}D^N D_{t,q}^M w_{q+1}^{(c)}\right|
    &\lec (\si_{\upsilon^{(c)}}^+  + \de_{q+3\bn})^{\sfrac12}\Ga_{q}^{-1} (\la_{q+\bn}\Ga_{q+\bn}^{{\sfrac1{10}}})^N\MM{M,\Nindt,\tau_q^{-1}\Ga_q^{i+16},\Tau_q^{-1}\Ga_{q}^9}
    \label{est.vel.inc.c.by.pr}
\end{align}
\end{subequations}
for $N, M\leq \sfrac{\Nfin}5$. Note that thanks to the factor $r_q\Ga_q^{-1}$ in \eqref{exp.w.q+1}, the bound in \eqref{est.vel.inc.c.by.pr} has extra gain of $r_q\Ga_q^{-1}$ compared to \eqref{est.vel.inc.p.by.pr}. This gain will be useful when we deal with the divergence corrector stress errors in subsection \ref{sss:dce} and divergence corrector current errors in \cite[subsection~5.5]{GKN23}.  We also record an upgraded version of \eqref{sunday:morning:1}, which states that in the same range of $N$ and $M$, we have that
\begin{subequations}\label{sunday:morning:2}
\begin{align}
    \left|\psi_{i,\qbn-1} D^N D_{t,\qbn-1}^M w_{q+1}^{(p)}\right|
    &\lec (\si_{\upsilon^{(p)}}^+  + \de_{q+3\bn})^{\sfrac12}r_{q}^{-1} (\la_{q+\bn}\Ga_{q+\bn}^{{\sfrac1{10}}})^N \notag \\
    &\qquad\quad \times\MM{M,\Nindt,\tau_{\qbn-1}^{-1}\Ga_{\qbn-1}^{i-5},\Tau_{\qbn-1}^{-1}\Ga_{\qbn}^{-1}} \, ,
    \label{est.vel.inc.p.by.pr.upup}\\
    \left|\psi_{i,\qbn-1} D^N D_{t,\qbn-1}^M w_{q+1}^{(c)}\right|
    &\lec (\si_{\upsilon^{(c)}}^+  + \de_{q+3\bn})^{\sfrac12}\Ga_{q}^{-1} (\la_{q+\bn}\Ga_{q+\bn}^{{\sfrac1{10}}})^N \notag \\
    &\qquad\quad \times\MM{M,\Nindt,\tau_{\qbn-1}^{-1}\Ga_{\qbn-1}^{i-5},\Tau_{\qbn-1}^{-1}\Ga_{\qbn}^{-1}} \, .
    \label{est.vel.inc.c.by.pr.upup}
\end{align}
\end{subequations}
The proof of \eqref{sunday:morning:2} is immediate from Hypothesis~\ref{eq:inductive:timescales} at level $q$ and Remark~\ref{rem:checking:hyp:dodging:1}, which asserts that Hypothesis~\ref{hyp:dodging1} is verified at level $q+1$ with $q'=\qbn$.
\end{remark*}

Before giving the proof of Lemma~\ref{lem:pr.inc.vel.inc.pot}, we record the following lemma, which investigates the current error generated by the pressure increment $\si_{\upsilon}$. The proof of both lemmas will proceed using Proposition~\ref{lem:pr.vel}.

\begin{lemma*}[\bf Current error from  the pressure increment]\label{lem:pr.current.vel.inc}
There exists a current error $\phi_{\upsilon}$ generated by $\sigma_\upsilon$ such that the following hold.
\begin{enumerate}[(i)]
\item\label{item:cpi:1} We have the decomposition and equalities 
    \begin{subequations}
    \begin{align}
\phi_{{\upsilon}} &= \underbrace{\phi_{{\upsilon}}^*}_{\textnormal{nonlocal}} + \underbrace{\sum_{m'=q+\half+1}^{\qbn} \phi_{{\upsilon}}^{m'}}_{\textnormal{local}} 
\label{S:pr:velocity:dec:statement}\\
\div \phi_{\upsilon} &= D_{t,q}\si_{\upsilon} - \bmu_{\si_{\upsilon}},  \notag
\end{align}
where $\bmu_{\si_{\upsilon}}$ is defined as in \eqref{def:bmu:vel:presh}.
\end{subequations}
\item For all $N\leq \sfrac{\Nfin}5$ and $M\leq \sfrac{\Nfin}5-\NcutSmall-1$ and $q+\half+1\leq m' \leq q+\bn$,
\begin{align}
\left|\psi_{i,q} D^ND_{t,q}^M\phi_{{\upsilon}}^{m'}\right|
&\lec \Ga_m^{-100}(\pi_q^{m'})^{\sfrac32}r_m^{-1} (\lambda_{m}\Ga_{m'})^N \MM{M,\Nindt, \tau_q^{-1}\Gamma_{{q}}^{i+16},\Tau_q^{-1}\Ga_q^9}
\label{pr:current:vel:loc:pt:q}
\, .
\end{align}
\item For all $N\leq 3\Nind$ and $M\leq 3\Nind$,
\begin{align}
    \norm{D^N\Dtq^M\phi_{{\upsilon}}^{*}}_{\infty}
    \lec \delta_{q+3\bn}^{\sfrac 32} \Tau_{q+\bn}^{2\Nindt} \lambda_{q+\bn+2}^{-10}
    (\lambda_{q+\bn}\Ga_{q+\bn})^N \MM{M,\Nindt, {\tau_q^{-1}},\Tau_q^{-1}\Ga_q^9} \, .
 \label{pr:current:vel:nonloc:infty1:q}
\end{align}
\item For all $q+1\leq q'\leq q+\half$, $q+\half+2\leq m \leq q+\bn $, and $q+1\leq q''\leq m-1$, we have the support properties
\begin{align}
  \supp(\phi_{{\upsilon}}^{q+\half+1})
    \cap B(\hat{w}_{q'}, \la_{q+1}^{-1}\Ga^2_{q})=\emptyset \, , \qquad
    \supp(\phi_{{\upsilon}}^{m}) \cap \supp \hat w_{q''} = \emptyset  \label{supp:pr:vel:q} \, .
\end{align}
\end{enumerate}
\end{lemma*}

\begin{proof}[Proofs of Lemma~\ref{lem:pr.inc.vel.inc.pot} and Lemma~\ref{lem:pr.current.vel.inc}]
\texttt{Step 1: Setup and Assumptions from Proposition~\ref{lem:pr.vel}.} In order to create a pressure increment which dominates \emph{all} of the various velocity increment potentials $\upsilon_{q+1,h}^{\pttl}$ defined in item~\eqref{item:eckel:1}, we shall create pressure increments which dominate each separate piece, and then sum at the end. We fix all indices $\pxi,I,\diamond,h,h'$ from the formula in \eqref{rep:vel:inc:pot:eckel} and apply Proposition~\ref{lem:pr.vel} with the following choices:
\begin{align*}
    N_*&=\sfrac{\Nfin}4 - 2\dpot^2 , \quad M_* = \sfrac{\Nfin}5, \quad M_t = \Nindt, \quad N_\circ = M_\circ = 3\Nind \, , \\
    \hat \upsilon &= \Upsilon_{\pxiid}^{h,h'} \, , \quad  G = H_\pxiid^{h,h'} \, , \quad \rho = \rho_\pxiid^\hhp \, , \quad \pi = \pi_\ell \Ga_q^{30}\, , \quad K_\circ\textnormal{ as in } \\
  \const_{G,p} &= \left| \supp \left(\eta_{i,j,k,\xi,\vecl,\diamond}\zetab_\xi^{I,\diamond} \right) \right|^{\sfrac 1p}\Ga_q^{j+7} \de_{q+\bn}^\frac12 r_q^{-\sfrac 13} + \la_\qbn^{-10} \, , \quad K_\circ \textnormal{ as in item~\eqref{i:par:9.5}} \\
    \const_{\rho, p}&= r_q^{\frac2p-\frac 23} \, , \quad \la = \la_{q+\half}, \quad \la' = \La_q, \quad \nu = \tau_q^{-1}\Ga_q^{i+13},\quad \nu'= \Tau_q^{-1} \Ga_q^8 \, , \quad \La =\la_{q+\bn}, \\
    r_G &= r_{\hat\upsilon} = r_q \, , \quad  \mu = \la_{q+\half}\Ga_{q} \, , \quad 
    \Ga = \Ga_q^{\sfrac1{10}} \, , \quad \Phi = \Phi_{(i,k)}\, , \quad v = \hat u_q \, , \quad \const_v =  \La_q^{\sfrac12} \, ,  
    \\
    \mu_0 &=\la_{q+\half+1}, \quad
    \mu_1 = \la_{q+\half+\sfrac32}, \quad
    \mu_m = \la_{q+\half+m}, \quad \mu_{\bar m} = \la_{q+\bn+1}, \quad
    \de_{\rm tiny} = \de_{q+3\bn}\, ,
\end{align*}
where $\mu_m=\la_{q+\half+m}$ above is defined for $2\leq m \leq \bm$. Then we have that \eqref{est.G.sample3}--\eqref{eq:desert:def:svi} are verified from \eqref{eq:inverse:divd:linear.ap.H.eckel}--\eqref{eq:inverse:divd:linear:ap:total:eckel}, \eqref{eq:sample:5:decoup:vel} holds by definition and by \eqref{condi.Ndec0}, \eqref{eq:DDpsi:sample3}--\eqref{eq:inverse:div:v:global.sample3} hold from \eqref{eq:nasty:D:vq:old}, Corollary \ref{cor:deformation}, \eqref{eq:bobby:old}, and \eqref{v:global:par:ineq}, \eqref{par.con.sample3.Ncut} holds from \eqref{condi.Ncut0.1}, \eqref{par.con.sample3} holds due to \eqref{condi.Ncut0.2}, \eqref{par.con.sample3.dec} holds due to \eqref{condi.Nfin0}, \eqref{eq:sample:prop:parameters:0:vel} holds from direct computation, and \eqref{eq:sample:riots:4:vel:1}--\eqref{eq:sample:riots:4:vel:3} hold due to \eqref{i:par:10}.
\smallskip

\noindent\texttt{Step 2: Part 2 from Proposition~\ref{lem:pr.vel} and proof of Lemma~\ref{lem:pr.inc.vel.inc.pot}. } We now apply the conclusions from Part 2 of Proposition~\ref{lem:pr.vel}.  We first have from \eqref{eq:desert:some:stuff} and \eqref{est.w.by.pr} the existence of a pressure increment $\sigma_{\Upsilon_{\pxiid}^{\hhp}}=\sigma^+_{\Upsilon_{\pxiid}^{\hhp}}-\sigma^-_{\Upsilon_{\pxiid}^{\hhp}}$ such that
\begin{align}
\left|D^N D_{t,q}^M  \Upsilon_\pxiid^\hhp\right|
    \lec \left(\si_{\upsilon_\pxiid^\hhp}^+  + \de_{q+3\bn}\right)^{\sfrac12}r_{q}^{-1} (\la_{q+\bn}\Ga_{q+\bn}^{\sfrac{1}{10}})^N\MM{M,\Nindt,\tau_q^{-1}\Ga_q^{i+15},\Tau_q^{-1}\Ga_{q}^9} \label{eckel:pressure}
\end{align}
for all $N \leq \sfrac{\Nfin}4-2\dpot^2$ and $M\leq \sfrac{\Nfin}5$.  Then using items~\eqref{sample3:item:3}--\eqref{sample3:item:1} and \eqref{condi.Nindt}, we have that
\begin{subequations}
\begin{align}
\left|D^N D_{t,q}^M \sigma^+_{\Upsilon_{\pxiid}^{\hhp}}\right|&\lec \left(\sigma^+_{\Upsilon_{\pxiid}^{\hhp}} + \de_{q+3\bn}\right) (\la_{q+\bn}\Ga_{q+\bn}^{\sfrac{1}{10}})^N\MM{M,\Nindt,\tau_q^{-1}\Ga_q^{i+15},\Tau_q^{-1}\Ga_{q}^9}\, , \label{est.pr.vel.inc:piece}\\
\norm{D^N D_{t,q}^M \sigma^+_{\Upsilon_{\pxiid}^{\hhp}} }_{\sfrac32} &\lec \left| \supp \left(\eta_{i,j,k,\xi,\vecl,\diamond}\zetab_\xi^{I,\diamond} \right) \right|^{\sfrac 23}\Ga_q^{2j+14}\de_{q+\bn} r_q^{\sfrac 43} \notag\\
&\qquad \qquad \times (\la_{q+\bn}\Ga_{q+\bn}^{\sfrac{1}{10}})^N\MM{M,\Nindt,\tau_q^{-1}\Ga_q^{i+15},\Tau_q^{-1}\Ga_{q}^9}\, , \label{est.pr.vel.inc.32:piece}\\
\norm{D^N D_{t,q}^M \sigma^+_{\Upsilon_{\pxiid}^{\hhp}} }_{\infty} &\lec \Ga_q^{\badshaq+20}(\la_{q+\bn}\Ga_{q+\bn}^{\sfrac{1}{10}})^N\MM{M,\Nindt,\tau_q^{-1}\Ga_q^{i+15},\Tau_q^{-1}\Ga_{q}^9}\, , \label{est.pr.vel.inc.infty:piece}\\
\norm{D^N D_{t,q}^M \sigma^-_{\Upsilon_{\pxiid}^{\hhp}}}_{\sfrac32} &\lec\left| \supp \left(\eta_{i,j,k,\xi,\vecl,\diamond}\zetab_\xi^{I,\diamond} \right) \right|^{\sfrac 23}\Ga_q^{2j+14} \de_{q+\bn} r_q^{\sfrac 43} \notag\\
&\qquad \qquad \times (\la_{q+\half}\Ga_{q+\half})^N\MM{M,\Nindt,\tau_q^{-1}\Ga_q^{i+15},\Tau_q^{-1}\Ga_{q}^9} \,  , \label{est.pr.vel.inc.-32:piece}\\
\norm{D^N D_{t,q}^M \sigma^-_{\Upsilon_{\pxiid}^{\hhp}}}_\infty &\lec \Ga_q^{\badshaq+20}(\la_{q+\half}\Ga_{q+\half})^N\MM{M,\Nindt,\tau_q^{-1}\Ga_q^{i+15},\Tau_q^{-1}\Ga_{q}^9}\, , \quad \label{est.pr.vel.inc.-infty:piece}\\
\left| D^N D_{t,q}^M \sigma^-_{\Upsilon_{\pxiid}^{\hhp}}\right|
&\lec  \pi_\ell\Ga_q^{30} r_{q}^{\sfrac 43} (\la_{q+\half}\Ga_{q+\half})^N\MM{M,\Nindt,\tau_q^{-1}\Ga_q^{i+15},\Tau_q^{-1}\Ga_{q}^9}\, ,
    \label{est.pr.vel.inc-:piece}
\end{align}
\end{subequations}
for all $N \leq \sfrac{\Nfin}4-2\dpot^2-\NcutLarge$ and $M\leq \sfrac{\Nfin}5-\NcutSmall$. In \eqref{est.pr.vel.inc.infty:piece} and \eqref{est.pr.vel.inc.-infty:piece}, we used \eqref{ineq:jmax:use}.  Finally, from \eqref{supp:pr:vecl}, \eqref{eq:inverse:divd:linear.ap.H.eckel}, \eqref{eq:inverse:divd:linear.ap.eckel}, \eqref{supp.upsilon.e.lem}, and Lemma \ref{lem:dodging}, we get the support properties
\begin{align*}
    \supp\left(\si_{\Upsilon_\pxiid^\hhp}^{+}\right)
    &\subseteq \supp\left(\Upsilon_\pxiid^\hhp\right) \\
    &\subseteq \supp \left(\chi_{i,k,q} \zeta_{q,\diamond,i,k,\xi,\vecl} \left(\rhob_{(\xi)}^\diamond \zetab_{\xi}^{I,\diamond} \right)\circ \Phiik \right) \cap
    B\left( \supp \varrho^I_{(\xi),\diamond} , 2 \lambda_\qbn^{-1} \right)\circ \Phiik  \, ,\\
    \supp\left(\si_{\Upsilon_\pxiid^\hhp}^{-}\right) &\cap B(\hat{w}_{q'}, \la_{q'}^{-1}\Ga_{q'})
    \subseteq  \supp\left(\eta_{i,j,k,\xi,\vecl,\diamond}\zetab_\xi^{I,\diamond}\right)  \cap B(\hat{w}_{q'}, \la_{q'}^{-1}\Ga_{q'}) = \emptyset \, ,
\end{align*}
for $q+1\leq q'\leq q+\half$.

We now sum over $h,h',\pxi,i,\diamond$ (while recalling from \eqref{rep:vel:inc:pot:eckel} that summation over $\pxi$ includes summation over $i,j,k,\xi,\vecl,\ovj$ as well as any indices needed for the application of the Faa di Bruno formula) and set
\begin{align*}
    \si_{\upsilon}^{\pm}
     := \sum_{\pxi,I,\diamond,h',h}
    \si_{\Upsilon_\pxiid^\hhp}^{\pm} \, .
\end{align*}
From \eqref{eckel:pressure}, \eqref{rep:vel:inc:pot:eckel}, \eqref{eq:desert:cowboy:sum}, and Corollary~\ref{lem:agg.pt} with $H={\Upsilon_\pxiid^\hhp}$ and $\varpi=\si_{\Upsilon_\pxiid^\hhp}^{+}+\mathbf{1}_{\supp {\Upsilon_\pxiid^\hhp}}\delta_{q+3\bn}$, we have that \eqref{est.vel.inc.pot.by.pr} holds. We have \eqref{defn:si.upsilon} from the formula above.  In order to verify \eqref{est.pr.vel.inc}--\eqref{est.pr.vel.inc-}, we appeal to \eqref{est.pr.vel.inc:piece}--\eqref{est.pr.vel.inc-:piece} and Corollaries~\ref{rem:summing:partition} and \ref{lem:agg.pt}. Specifically, the $L^{\sfrac 32}$ estimates in \eqref{est.pr.vel.inc.32} and \eqref{est.pr.vel.inc.-32} use \eqref{eq:par:div:1} and Corollary~\ref{rem:summing:partition} with $\theta_2=\theta=2$, $H=\sigma_{\Upsilon_\pxiid^\hhp}^{\pm}$, and $\const_H=\delta_{q+\bn}r_q^{\sfrac 43}\Ga_q^{14}$. The $L^\infty$ estimates in \eqref{est.pr.vel.inc.infty} and \eqref{est.pr.vel.inc.-infty} follow from \eqref{eq:desert:cowboy:sum}, \eqref{eq:par:div:2}, and Corollary~\ref{lem:agg.pt} and  with the same choice of $H$ and $\varpi=\Ga_q^{\badshaq+20}\mathbf{1}_{\supp \Upsilon_\pxiid^\hhp}$.  Finally, the pointwise estimates in \eqref{est.pr.vel.inc} and \eqref{est.pr.vel.inc-} follow from Corollary~\ref{lem:agg.pt} in much the same manner as the $L^\infty$ estimates just derived, and we omit further details.

\smallskip

\noindent\texttt{Step 3: Part 3 from Proposition~\ref{lem:pr.vel} and proof of Lemma~\ref{lem:pr.current.vel.inc}.} We now apply the conclusions from Part 3 of Propsition~\ref{lem:pr.vel}. From item~\eqref{sample3:item:decomp}, there exist current errors $\phi_{\Upsilon_\pxiid^\hhp}$ such that we have the decompositions and equalities
\begin{subequations}
\begin{align}
\phi_{\Upsilon_\pxiid^\hhp} &= \underbrace{\phi_{\Upsilon_\pxiid^\hhp}^*}_{\textnormal{nonlocal}} + \underbrace{\sum_{m'=q+\half+1}^{\qbn} \phi_{\Upsilon_\pxiid^\hhp}^{m'}}_{\textnormal{local}} \label{eq:cpi:defs:proof}\\
&= \underbrace{(\divH+\divR)\left(D_t \sigma_{\Upsilon_\pxiid^\hhp}^*\right) + \sum_{m'=q+\half+1}^{\qbn}  \divR\left(D_t \sigma_{\Upsilon_\pxiid^\hhp}^{m'}\right)}_{\textnormal{nonlocal}} + \underbrace{\sum_{m'=q+\half+1}^{\qbn}  \divH\left(D_t \sigma_{\Upsilon_\pxiid^\hhp}^{m'}\right)}_{\textnormal{local}} \, , \notag \\
\div&\left( \phi_{\Upsilon_\pxiid^\hhp}^{m'}(t,x) + \divR\left(D_t \sigma_{\Upsilon_\pxiid^\hhp}^{m'} \right)(t,x) \right) = D_t \sigma_{\Upsilon_\pxiid^\hhp}^{m'}(t,x) - \int_{\T^3} D_t \sigma_{\Upsilon_\pxiid^\hhp}^{m'}(t,x') \, dx' \, , \notag  \\
\div &\left( \phi_{\Upsilon_\pxiid^\hhp}^*(t,x) - \sum_{m=0}^{\bm} \divR \left(D_t \sigma_{\Upsilon_\pxiid^\hhp}^{m'} \right)(t,x) \right) = D_t \sigma_{\Upsilon_\pxiid^\hhp}^*(t,x) - \int_{\T^3} D_t \sigma_{\Upsilon_\pxiid^\hhp}^*(t,x') \, dx' \, . \notag 
\end{align}
\end{subequations}
Next, from \eqref{sample3:item:7} in Proposition~\ref{lem:pr.vel}, \eqref{condi.Nindt}, and \eqref{condi.Nfin0}, we have that for $(p,p')=(3, \sfrac32)$ or $(\infty,\infty)$ and $2\leq m \leq \bm$,
\begin{subequations}
\begin{align}
\norm{D^ND_t^M\phi_{{\Upsilon_\pxiid^\hhp}}^{0}}_{p'} &\lec \tau_q^{-1} \Ga_q^{i+14} \left( \de_{q+\bn} r_q^{-\sfrac 23}
\Ga_q^{2j+14} \left| \supp \left(\eta_{i,j,k,\xi,\vecl,\diamond}\zetab_\xi^{I,\diamond} \right) \right|^{\sfrac 2p} + \la_\qbn^{-20} \right) \notag\\
&\qquad \times \left( \frac{\la_{q+\half+1}}{\la_{q+\half}} \right)^{\frac 43 - \frac{2}{p'}} r_q^2 \la_{q+\half}^{-1} (\la_{q+\half+1})^N \MM{M,\Nindt,\tau_q^{-1}\Ga_q^{i+15},\Tau_q^{-1}\Ga_q^9}\, ,
\label{pr:current:vel:loc:p:0:piece}\\
\left|D^ND_t^M\phi_{{\Upsilon_\pxiid^\hhp}}^{0}\right|&\lec \tau_q^{-1}\Ga_q^{i+50} \pi_\ell  r_{q}^{\sfrac43} \left(\frac{\la_{q+\half+1}}{\la_{q+\half}}\right)^{\sfrac43} \la_{q+\half}^{-1} \notag\\
&\qquad \times (\la_{q+\half+1})^N \MM{M,\Nindt,\tau_q^{-1}\Ga_q^{i+15},\Tau_q^{-1}\Ga_q^9}
\label{pr:current:vel:loc:pt:0:piece}
\, , \\
\norm{D^ND_t^M\phi_{{\Upsilon_\pxiid^\hhp}}^{m}}_{p'}&\lec 
\tau_q^{-1} \Ga_q^{i+16} \left( \de_{q+\bn} r_q^{-\sfrac 23}
\Ga_q^{2j+14} \left| \supp \left(\eta_{i,j,k,\xi,\vecl,\diamond}\zetab_\xi^{I,\diamond} \right) \right|^{\sfrac 2p} + \la_\qbn^{-20} \right) \notag\\
&\qquad \times \left( \frac{\min\left(\la_{q+\half+m},\la_\qbn\right)}{\la_{q+\half}} \right)^{\frac 43 - \frac{2}{p'}} r_q^2 (\la_{q+\half+{m-1}}^{-2}\la_{q+\half+m}) \notag\\
&\qquad \qquad \quad \times (\min(\la_{q+\half+m}, \la_{q+\bn}\Ga_{q+\bn}))^N  \MM{M,\Nindt,\tau_q^{-1}\Ga_q^{i+15},\Tau_q^{-1}\Ga_q^9}\, , 
\label{pr:current:vel:loc:p:piece}
\\
\left|D^ND_t^M\phi_{{\Upsilon_\pxiid^\hhp}}^{m}\right|
&\lec 
\tau_q^{-1}\Ga_q^{i+50} \pi_\ell r_{q}^{\sfrac43} \left(\frac{\min(\la_{q+\half+m},\la_\qbn)}{\la_{q+\half}\Ga_{q}}\right)^{\sfrac43} (\la_{q+\half+{m-1}}^{-2}\la_{q+\half+m})\nonumber\\
&\qquad \times 
(\min(\la_{q+\half+m}, \la_{q+\bn}\Ga_{q+\bn}))^N \MM{M,\Nindt,\tau_q^{-1}\Ga_q^{i+15},\Tau_q^{-1}\Ga_q^9}\, ,
\label{pr:current:vel:loc:pt:piece}
\end{align}
\end{subequations}
for $N\leq \sfrac{\Nfin}{5}$ and $M\leq \sfrac{\Nfin}{5}-\NcutLarge-1$. In the case $m=1$, we have bounds which match the bounds for $m=2$ above, except that the inverse divergence gain of $\la_{q+\half+m-1}^{-2}\la_{q+\half+m}$ is replaced with $\la_{q+\half+\sfrac 32}^{-2} \la_{q+\half+1}$. Furthermore, we have from \eqref{pr:current:vel:nonloc:infty1} and item~\eqref{i:par:9.5} that
\begin{align}
\norm{D^ND_t^M\phi_{\Upsilon_\pxiid^\hhp}^{*}}_{\infty}
&\lec \delta_{q+3\bn}^{\sfrac 32} \Tau_{q+\bn}^{2\Nindt} \la_\qbn^{-60} (\la_\qbn\Ga_\qbn)^N (\tau_q^{-1}\Ga_q^{i+14})^M 
\label{pr:current:vel:nonloc:piece}
\end{align}
for $N, M\leq 3\Nind$. Finally, \eqref{sample3:item:6} from Proposition~\ref{lem:pr.vel}, \eqref{eq:inverse:divd:linear.ap.eckel}, \eqref{eq:inverse:divd:linear.ap.H.eckel}, and Lemma \ref{lem:dodging} give that for each $1\leq m\leq \bar m$ and any $q+1\leq q'\leq q+\half$ and $q+1\leq q'' \leq q+\half+m-1$
\begin{align}
    \supp\left(\phi_{\Upsilon_\pxiid^\hhp}^{0}\right) \cap B(\hat w_{q'}, \la_{q+1}^{-1}\Ga_{q}^2)=\emptyset \, , \qquad
    \supp\left(\phi_{\Upsilon_\pxiid^\hhp}^{m}\right)
    \cap \supp \hat w_{q''} \, , \notag \\
      \supp\left(\phi_{\Upsilon_\pxiid^\hhp}^{0}\right) \, , \,  \supp\left(\phi_{\Upsilon_\pxiid^\hhp}^{m}\right) \subseteq \supp\left( \eta_{i,j,k,\xi,\vecl,\diamond} \zetab_\xi^{I,\diamond} \right) \, . \label{supp:current:pr:piece} 
\end{align}

We now sum over $h,h',\pxi,i,\diamond$ (while recalling from \eqref{rep:vel:inc:pot:eckel} that summation over $\pxi$ includes summation over $i,j,k,\xi,\vecl,\ovj$ as well as any indices needed for the application of the Faa di Bruno formula) and set
\begin{align}
\phi_{{\upsilon}}^{q+\half+1}
&:= \sum_{\pxiid,h',h} \phi_{\Upsilon_\pxiid^\hhp}^{0} \, , \qquad 
\phi_{{\upsilon}}^{q+\half+2}
:= \sum_{\pxiid,h',h} \sum_{m=1}^2
\phi_{\Upsilon_\pxiid^\hhp}^{m} \label{eq:cpi:defs} \\
\phi_{{\upsilon}}^{q+\half+m}
&:= \sum_{\pxiid,h',h} \phi_{\Upsilon_\pxiid^\hhp}^{m} \, , \qquad 
\phi_{{\upsilon}}^{q+\bn}
:= \sum_{\pxiid,h',h} \sum_{m=\bar m -1}^{\bar m}\phi_{\Upsilon_\pxiid^\hhp}^{m} \, , \qquad 
\phi_{{\upsilon}}^{*}
:= \sum_{\pxiid,h',h}
\phi_{\Upsilon_\pxiid^\hhp}^{*} \notag 
\end{align}
for $3\leq m\leq \bar m-2$.

We can now conclude the proof of Lemma~\ref{lem:pr.current.vel.inc}.  First, we have that item~\eqref{item:cpi:1} follows from the definitions in \eqref{eq:cpi:defs} and \eqref{eq:cpi:defs:proof}. Next, we have that \eqref{supp:pr:vel:q} follows from the same definitions, \eqref{supp:current:pr:piece}, and Lemma~\ref{lem:dodging}. We can achieve the nonlocal bounds in \eqref{pr:current:vel:nonloc:infty1:q} from \eqref{pr:current:vel:nonloc:piece} and summation over all indices $\pxi,I,\diamond,h',h$, which from Lemma~\ref{lem.cardinality}, \eqref{eq:imax:upper:lower}, Lemma~\ref{lem:maximal:j}, and the discussion following \eqref{rep:vel:inc:pot:eckel} is bounded by $\la_\qbn^4$.  The bound for $\bmu_{\sigma_\upsilon}$ in item~\eqref{i:presh:vel:mean} follows similarly from \eqref{est:mean.Dtsivel} \eqref{ineq:K_0}, and a large choice of $a_*$ in \eqref{i:choice:of:a} to ensure that we can put the prefactor of $\max(1,T)^{-1}$ in the amplitude. Finally, we may conclude \eqref{pr:current:vel:loc:pt:q} from an application of Corollary~\ref{lem:agg.Dtq} with $H=\phi_{\Upsilon^\hhp_\pxiid}^\bullet$ (with the value of $\bullet$ according to the divisions in \eqref{eq:cpi:defs}) and 
$$  \varpi=\Ga_q^{50}\pi_\ell r_q^{\sfrac 43}\left(\frac{\min(\la_{q+\half+m},\la_{\qbn})}{\la_{q+\half}}\right)^{\sfrac 43} \la_{q+\half+m-1}^{-2}\la_{q+\half+m} \, . $$
Indeed appealing to \eqref{eq:aggDtq:conc:1}, \eqref{ind:pi:lower}, \eqref{eq:ind.pr.anticipated}, \eqref{eq:desert:ineq}, and the fact that 
$$  r_q^{\sfrac 43} \left(\frac{\min(\la_{q+\half+m},\la_{\qbn})}{\la_{q+\half}}\right)^{\sfrac 43} \leq \Ga_q^{10}  $$
from the definition of $r_q$, we conclude the proof.
\end{proof}

\opsubsection{Estimates for the velocity increment potentials}\label{op:wed}
We will now verify the inductive assumptions of subsubsection~\ref{sec:ind.vel.inc.pot} in the following proposition. We first recall the definitions of $\upsilon_{q+1}$ and $e_{q+1}$ from Remark~\ref{rem:rep:vel:inc:potential} and the mollifier $\tilde{\mathcal{P}}_{\qbn,x,t}$ from Definition~\ref{def:wqbn} and define\index{$\hat\upsilon_\qbn$}\index{$\hat e_{q+\bn}$}
\begin{equation}\label{def.vel.inc.pot.mollified}
    \hat \upsilon_{q+\bn} := \mathcal{\tilde P}_{q+\bn,x,t} \upsilon_{q+1}, \qquad \hat e_{q+\bn} := \mathcal{\tilde P}_{q+\bn,x,t} e_{q+1}\,.
\end{equation}

\begin{proposition*}[\bf Verifying \eqref{exp.w.q'}, \eqref{supp.upsilon.e.ind}, and \eqref{est.e.inf} and setting up \eqref{est.upsilon.ptwise} at level $q+1$]\label{eq:prop:saturday:night} The velocity increment and velocity increment potentials satisfy the following.
\begin{enumerate}[(i)] 
\item $\hat w_{q+\bn}$ can be decomposed as 
\begin{align}\label{exp.w.qbn}
\hat w_{q+\bn} = \div^{\dpot} \hat \upsilon_{q+\bn} +\hat e_{q+\bn}  \,, 
\end{align}
which written component-wise gives $
\hat w_{q+\bn}^\bullet = 
\pa_{i_1}\cdots \pa_{i_\dpot} \hat\upsilon_{q+\bn}^{(\bullet, i_1, \cdots, i_\dpot)} + \hat e_{q+\bn}^\bullet
$. 
\item {For all $q+1\leq q' \leq q+\bn-1$, the supports of $\hat\upsilon_{q+\bn}$ and $\hat e_{q+\bn}$ satisfy
\begin{align}\label{supp.upsilon.e.verify}
B\left(\supp(\hat w_{q'}), \frac14 \la_{q'}\Ga_{q'}^2\right) \cap \left(\supp(\hat \upsilon_{q+\bn}) \cup \supp(\hat e_{q+\bn})\right) = \emptyset \, .  
\end{align}}

\item For $N+M \leq \sfrac{3\Nfin}{2}$, we have that $\hat\upsilon_{q+\bn,k}^\bullet :=\la_{q+\bn}^{\dpot-k}\partial_{i_1}\cdots \partial_{i_k}
    \hat\upsilon_{q+\bn}^{(\bullet,i_1,\dots,i_\dpot)}$, $0\leq k \leq \dpot$, satisfies the estimates
\begin{align}
    \left|\psi_{i,q+\bn-1} D^ND_{t,q+\bn-1}^{M} 
   \hat \upsilon_{q+\bn,k}
    \right|&< {\Ga_{q+\bn}} \left(\si_{\upsilon^{(p)}}^+ + \si_{\upsilon^{(c)}}^+ + 2\de_{q+3\bn} \right)^{\sfrac12} r_{q}^{-1}
     (\la_{q+\bn}{\Ga_{q+\bn}})^N \notag \\
     & \qquad \qquad \qquad \times \MM{M, \Nindt, \Ga_{q+\bn-1}^i \tau_{q+\bn-1}^{-1}, \Tau_{q+\bn-1}^{-1}{\Ga_{q+\bn-1}^2}} \, . \label{est.upsilon.ptwise.verify}
\end{align}
\item For $N+M \leq \sfrac{3\Nfin}{2}$, $\hat e_{q+\bn}$ satisfies
\begin{align} 
    \norm{ D^ND_{t,q+\bn-1}^{M} \hat e_{q+\bn}}_{\infty} 
    &\leq \de_{q+3\bn}^3 \Tau_\qbn^{10\Nindt}\la_{q+\bn}^{-10} 
    (\la_{q+\bn}{\Ga_{q+\bn}})^N  \notag \\
     & \qquad \qquad \qquad \times
     \MM{M, \Nindt, \tau_{q+\bn-1}^{-1}, \Tau_{q+\bn-1}^{-1} \Ga_{q+\bn-1}^2}
    \label{est.e.inf.verify}\, .
\end{align}
\end{enumerate}
\end{proposition*}

\begin{proof}[Proof of Proposition~\ref{eq:prop:saturday:night}]
We first note that \eqref{exp.w.qbn} follows immediately from the definition of $\hat \upsilon_{q+\bn}$ and $\hat e_{q+\bn}$ in \eqref{def.vel.inc.pot.mollified} and the identity in Remark~\ref{rem:rep:vel:inc:potential}. 

Next, an immediate consequence of \eqref{supp.upsilon.e.lem} and \eqref{eq:dodging:oldies:prep} is that
\begin{align*}
    B\left(\supp(\hat w_{q'}), \frac12 \la_{q'}\Ga_{q'}^2, 2\Tau_q \right) \cap  \left(\supp( \upsilon_{q+1}) \cup \supp(e_{q+1})\right) = \emptyset \, .
 \end{align*}
for all $q+1\leq q' \leq q+\bn-1$. Now notice that by properties of the mollification, we have that
$$\supp(\hat \upsilon_{q+\bn}) \subseteq B\left(\supp(\upsilon_{q+1}), \left(\la_{q+\bn}\Ga_{q+\bn-1}^{\sfrac12}\right)^{-1}, \Tau_{q+1}^{-1}\right) \, , $$ 
and similarly 
$$\supp(\hat e_{q+\bn}) \subseteq B\left(\supp(e_{q+1}), \left(\la_{q+\bn}\Ga_{q+\bn-1}^{\sfrac12}\right)^{-1}, \Tau_{q+1}^{-1}\right) \, . $$ 
With this we now see that \eqref{supp.upsilon.e.verify} is satisfied.

Note that from \eqref{supp.upsilon.e.verify} and an application of Lemma~\ref{lem:upgrading.material.derivative}, we see that \eqref{est.vel.inc.pot.by.pr} implies that for all $N,M \leq \sfrac{\Nfin}5$, $0\leq k\leq \dpot$ and $\ttl=p,c$,
\begin{align}\label{est.vel.inc.pot.by.pr.upgraded}
    \left|\psi_{i,q+\bn-1}D^N D_{t,q+\bn-1}^M \hat\upsilon_{q+\bn,k}^{(\ttl)}\right|
    &\lec (\si_{\upsilon^{(p)}}^+  + \de_{q+3\bn})^{\sfrac12}r_{q}^{-1} (\la_{q+\bn}\Ga_{q+\bn}^{{\sfrac1{10}}})^N\MM{M,\Nindt,\tau_{q+\bn-1}^{-1}\Ga_{q+\bn-1}^{i-4},\Tau_q^{-1}\Ga_{q}^9} \, . 
\end{align}
Now we apply Proposition~\ref{lem:mollification:general} with the parameter choices
\begin{align*}
    &p=3, \infty \, , \quad N_{\rm g}, N_{\rm c} \textnormal{ as in \eqref{i:par:12}} \, , \quad M_t = \Nindt\,  , \quad N_* = \sfrac{\Nfin}{5} \, , \\
    &N_\gamma = 2\Nfin \, , \quad  \Omega = \supp \psi_{i,q+\bn-1}\, , \quad v = \hat u_{q+\bn-1}\, , \quad i=i \, , \quad c=-1 \, , \\
    &\lambda = \lambda_{q+\bn} \, , \quad \Lambda = \lambda_{q+\bn}\Ga_{q+\bn-1} \,, \quad \Ga = \Ga_{q+\bn-1}, \quad \tau = \tau_{q+\bn-1}\Gamma_{q+\bn-1}^{-2} \, , \quad \Tau = \Tau_{q+\bn-1} \, ,\\
    &f = \upsilon_{q+1,k}^{(l)} \, , \quad \const_{f,3} = \Ga_q^{20} \de_{q+\bn}^{\sfrac 12}r_q^{-\sfrac 13} \, , \quad \const_{f,\infty} = \tilde \const_f =  \Ga_q^{\sfrac{\badshaq}{2} + 16} r_q^{-1} \, , \quad \const_v = \Lambda_{q+\bn-1}^{\sfrac 12} \,  .
\end{align*}
In a similar way to the proof of Lemma~\ref{lem:mollifying:w}, we see that all the assumptions of the proposition are satisfied. Therefore, conclusion \eqref{eq:moll:conc:2} implies that $N,M \leq \sfrac{\Nfin}5$, $0\leq k\leq \dpot$ and $l=p,c$,
\begin{align*}
    \left\|D^N D_{t,q+\bn-1}^M \left(\hat\upsilon_{q+\bn,k}^{(l)} - \upsilon_{q+1,k}^{(l)}\right) \right\|_\infty
    &\lec \de_{q+3\bn}^{3}\Tau_{q+\bn}^{25\Nindt} (\la_{q+\bn}\Ga_{q+\bn-1})^N\MM{M,\Nindt,\tau_{q+\bn-1}^{-1},\Tau_{q+\bn-1}^{-1}\Ga_{q+\bn-1}} \, .
\end{align*}
Combining this estimate with the pointwise estimate \eqref{est.vel.inc.pot.by.pr.upgraded} implies \eqref{est.upsilon.ptwise.verify} for $N,M \leq \sfrac{\Nfin}5$. The case when $\sfrac{\Nfin}5 \leq N+M \leq \sfrac{3\Nfin}2$ follows from first noticing that conclusion \eqref{eq:moll:conc:1} implies that for all $N,M \leq 2\Nfin$, $0\leq k\leq \dpot$ and $\ttl=p,c$, we have
\begin{align*}
    \left\|\psi_{i,q+\bn-1}D^N D_{t,q+\bn-1}^M \hat\upsilon_{q+\bn,k}^{(l)}\right\|_\infty
    &\lec \Ga_q^{\sfrac{\badshaq}2 +16} r_{q}^{-1} (\la_{q+\bn}\Ga_{q+\bn}^{{\sfrac1{10}}})^N\MM{M,\Nindt,\tau_{q+\bn-1}^{-1}\Ga_{q+\bn-1}^{i-4},\Tau_{q+\bn-1}^{-1}\Ga_{q+\bn-1}} \, .
\end{align*}
Then combining this estimate with \eqref{eq:Nind:darnit} implies estimate \eqref{est.upsilon.ptwise.verify} in this case.

Finally, to prove \eqref{est.e.inf.verify}, we must upgrade the nonlocal derivative bound in \eqref{est.new.e.inf}.  This is trivial using the extra prefactors of $\Tau_\qbn^{20\Nindt}$, and so we omit the details. 
\end{proof}

\opsubsection{New inductive cutoffs are dominated by the pressure increment}\label{op:wed:2}

We conclude this section with a lemma which shows that a rescaled combination of the intermittent pressure and the velocity pressure increment can be used to dominate a weighted sum of the velocity cutoff functions.

\begin{lemma*}\label{lem:pr.vel.dom.cutoff}
    The new velocity cutoff functions $\psi_{i,q+\bn}$ satisfy
    \begin{align}
    \sum_{i=0}^{\imax} \psi_{i,q+\bn}^2 \delta_{q+\bn} r_{q}^{-\sfrac 23} \Gamma_{q+\bn}^{2i} &\les  {r_{q}^{-2}} \left(\pi_{q}^{q+\bn} + \si_\upsilon^+ + \de_{q+3\bn}\right) \,   \label{eq:psi:q:qplusbn:ineq:0:recall}
\end{align}
for a $q$-independent implicit constant.
\end{lemma*}

\begin{proof}
    From \eqref{eq:psi:i:q:recursive} and the fact that all cutoff functions are bounded in between $0$ and $1$, we have that
\begin{align}
    \sum_{i=0}^{\imax} \psi_{i,q+\bn}^2 \delta_{q+\bn} r_{q}^{-\sfrac 23} \Gamma_{q+\bn}^{2i} &\lec \delta_{q+\bn} r_{q}^{-\sfrac 23} \sum_{i=0}^{\imax} \Gamma_{q+\bn}^{2i} \sum\limits_{\left\{\Vec{i}\colon\max\limits_{0\leq m\leq\NcutSmall} i_m =i\right\}} \prod\limits_{m=0}^{\NcutSmall} \psi_{m,i_m,q+\bn}^{2} \notag\\
    &\leq \sum_{m=0}^{\NcutSmall} \delta_{q+\bn} r_{q}^{-\sfrac 23}  \sum_{i_m \geq 0} \psi_{m,i_m,q+\bn}^2 \Gamma_{q+\bn}^{2i_m} \, . \label{eq:dd:one}
\end{align}
Therefore it will suffice to show that the right-hand side of \eqref{eq:psi:q:qplusbn:ineq:0:recall} dominates the double sum above. We will in fact fix $m$, take the sum over $i_m\geq 0$, multiply by $\Gamma_{q+\bn}$, and show that this is dominated by the right-hand side of \eqref{eq:psi:q:qplusbn:ineq:0:recall}.  Using that $m$ is bounded by $\NcutSmall$ and choosing $a$ large enough will then conclude the proof.

From the definition of $\psi_{m,i_m,q+\bn}$ in \eqref{eq:psi:m:im:q:def}, we have that
\begin{align}
    \Gamma_{q+\bn}^{2i_m} \psi_{m,i_m,q+\bn}^2 &\lec \Gamma_{q+\bn}^{2i_m} \sum_{\{j_m: i_*(j_m)\leq i_m\}} \psi_{j_m,q+\bn-1}^2 \psi_{m,i_m,j_m,q+\bn}^2 \notag \\
    &= \Gamma_{q+\bn}^{2i_*(j_m)} \psi_{j_m,q+\bn-1}^2 \psi_{m,i_*(j_m),j_m,q+\bn}^2 + \Gamma_{q+\bn}^{2i_m} \sum_{\{j_m: i_*(j_m) < i_m\}} \psi_{j_m,q+\bn-1}^2 \psi_{m,i_m,j_m,q+\bn}^2 \, . \label{eq:dom:proof:1}
\end{align}
From \eqref{eq:ineq:ij}, we know that the first term above is dominated by
\begin{align}
    \Gamma_{q+\bn-1}^{2j_m+4} \psi_{j_m,q+\bn-1}^2 \, . \notag
\end{align}
Since $m$ and $i_m$ only take finitely many values, we may bound the contribution to the right-hand sides of \eqref{eq:dd:one} and \eqref{eq:dom:proof:1} from the terms with $j_m$ such that $i_*(j_m)=i_m$ by an implicit constant multiplied by
\begin{align}
    \sum_{j_m \geq 0} \Gamma_{q+\bn-1}^{2j_m+4} \psi_{j_m,q+\bn-1}^2 \delta_{q+\bn}  r_{q}^{-\sfrac 23} 
    &\leq r_{q-1}^{-2}\pi_q^{q+\bn-1}  \Gamma_{q+\bn-1}^5 \frac{\delta_{q+\bn}  r_{q}^{-\sfrac 23}  }{\delta_{q+\bn-1}  r_{q-1}^{-\sfrac 23} }\notag\\
    &\leq \Ga_\qbn^{-2} r_{q}^{-2} \pi_{q}^{q+\bn}  \, . \notag
\end{align}
Here we have used the inductive assumption \eqref{eq:psi:q:q'} to achieve the first inequality above and the inequalities \eqref{eq:pr.subtle} and \eqref{eq:ind.pr.anticipated} to achieve the second inequality. We have thus concluded that the \emph{lowest terms} with $i_m=i_*(j_m)$ from \eqref{eq:dom:proof:1}, summed over $i_m$ and appropriately weighted, are indeed dominated by the right-hand side of \eqref{eq:psi:q:qplusbn:ineq:0:recall}.

We now must consider the rest of the terms in \eqref{eq:dom:proof:1}, for which $i_*(j_m)<i_m$. Assume that $(t,x)\in \supp (\psi_{j_m,q+\bn-1}^2 \psi_{m,i_m,j_m,q+\bn}^2)$. By \eqref{eq:h:j:q:def} and Lemma~\ref{lem:cutoff:construction:first:statement}, item~\eqref{item:cutoff:2}, and there exists $n\leq \NcutLarge$ such that
\begin{align}
    \frac{1}{4\NcutLarge} &\leq \Gamma_{q+\bn}^{-2i_m(m+1)} \delta_{q+\bn}^{-1} r_q^{\sfrac 23} (\lambda_{q+\bn} \Gamma_{q+\bn})^{-2n} (\tau_{q+\bn-1}^{-1}\Gamma_{q+\bn}^2)^{-2m} |D^n D_{t,q+\bn-1}^m \hat w_{\qbn}|^2 \, .\notag
\end{align}
Note that due to Definition~\ref{def:istar:j}, the fact that we consider $(t,x)\in \supp (\psi_{j_m,q+\bn-1}^2 \psi_{m,i_m,j_m,q+\bn}^2)$, and \eqref{eq:psi:i:j:partition:0}, which gives $i_m \geq i_*(j_m)$, we have that $\Ga_\qbn^{-i_m}\Ga_{\qbn-1}^{j_m} \leq 1$.
Now using \eqref{condi.Nindt} and that we are on the support of $\psi_{j,\qbn-1}$ by assumption so that  we may appeal to \eqref{sunday:morning:2}, we have that
\begin{align}
    \Gamma_{q+\bn}^{2 i_m} \delta_{q+\bn} r_q^{-\sfrac 23} &\lesssim (\lambda_{q+\bn} \Gamma_{q+\bn})^{-2n} (\tau_{q+\bn-1}^{-1}\Gamma_{q+\bn}^{2+i_m})^{-2m} \left(\si_{\upsilon}^+ + \de_{q+3\bn}\right) r_q^{-2} (\la_{q+\bn}\Ga_\qbn)^{2n} (\tau_{q+\bn-1}^{-1}\Ga_{q+\bn-1}^{j_m-5})^{2m}  \notag \\
    &\leq \left(\si_{\upsilon}^+ + \de_{q+3\bn}\right) r_q^{-2} \, . \label{eq:showing:mess:2}
\end{align}
Thus, \eqref{eq:psi:q:qplusbn:ineq:0:recall} follows from summing \eqref{eq:showing:mess:2} over $i_m\geq 0$, from which we find that
$$  \sum_{i_m\geq 0} \sum_{\{j_m: i_*(j_m) < i_m\}} \psi_{j_m,q+\bn-1}^2 \psi_{m,i_m,j_m,q+\bn}^2 \Gamma_{q+\bn}^{2i_m} \delta_{q+\bn} r_q^{-\sfrac 23} \les r_q^{-2} \left(\pi_{q}^{q+\bn} + \si_\upsilon^+ + \de_{q+3\bn}\right) \, . $$
Now summing over $0\leq m\leq \NcutSmall$ concludes the proof of \eqref{eq:psi:q:qplusbn:ineq:0:recall}.
\end{proof}

\section{Parameters}\label{sec:parm}
\subsection{Definitions and inequalities}\label{sec:para.q.ind}
In this section, we choose the values of the parameters and list important consequences. The choices in items~\eqref{i:par:1}--\eqref{item:choice:of:alpha} are rather delicate, while all the choices in items~\eqref{i:par:tau}--\eqref{i:choice:of:a} follow the plan of ``choosing a giant parameter which dwarfs all the preceding parameters."  It is imperative that each inequality below depends \emph{only} on parameters which have already been chosen, and that none depend on $q$.  We point out that in item~\eqref{i:par:2.5}, we define two parameters $\la_q$ and $\delta_q$ in terms of an undetermined large natural number $a$.  This is merely for ease of notation and computation.  Indeed one can check that none of the inequalities below require a precise choice of $a$, nor depend on $q$; rather, any sufficiently large choice of $a$ which may be used to absorb implicit constants will do. Therefore the precise choice of $a$ is made at the very end in item~\eqref{i:choice:of:a}.
\begin{enumerate}[(i)]
    \item\label{i:par:1} Choose $\beta\in(0,\sfrac 13)$ and $\bn$ a large positive multiple of $6$ as in \eqref{eq:choice:of:bn}.
    \item\label{i:par:beeeeeee} Choose $b\in (1,\sfrac{25}{24})$ as in \eqref{ineq:b}.
    \item\label{i:par:2.5} For an undetermined natural number $a$, define $\la_q$ and $\de_q$ as in \eqref{eq:def:la:de}.  Note that with this definition of $\la_q$, we have that
    \begin{equation}
        a^{(b^q)} \leq \la_q \leq 2a^{(b^q)} \qquad \textnormal{and} \qquad  \frac{1}{3} \la_q^b \leq \la_{q+1} \leq 2\la_q^b \, . \label{eq:la:comparability}
    \end{equation}
    As a consequence of these definitions, we shall deduce a number of inequalities, each of which is independent of the choice of $a$ and of $q$ once $a$ is sufficiently large.  At the end we will thus choose $a$ sufficiently large to absorb a number of implicit constants, including those in \eqref{eq:la:comparability}.  Therefore, in many of the following computations, we may make the slightly incorrect assumption that $\la_q$ is \emph{actually equal} to $a^{(b^q)}$ in order to streamline the arithmetic.
    \begin{enumerate}
        \item An immediate consequence of these definitions and of the first inequality in \eqref{ineq:b:first} is that
    \begin{align*}
        \delta_\qbn &\left( \la_q \la_{q+\sfrac \bn 3}^{-1} \right)^{\sfrac 23} \la_{\qbn+1}^4 \la_\qbn^{-4} \frac{\la_q^4\la_\qbn^4}{\la^8_{q+\half}} < \delta_{q+\sfrac{4\bn}{3}+2}\\
        &\iff 2\beta b^{\sfrac{4\bn}{3}+2} - 2\beta b^{\bn} < \frac 23 b^{\sfrac \bn 3} - \frac 23 - 4b^{\bn+1} + 4b^{\bn} - 4b^{\bn} + 8 b^{\sfrac \bn 2} - 4 \\
        &\iff 2\beta b^{\bn} (b-1) (1+b+\dots+b^{\sfrac \bn 3 +1}) < \frac 23 (b-1) (1+b+\dots+b^{\sfrac \bn 3 -1}) \\
        &\qquad \qquad \qquad \qquad \qquad \qquad \qquad \qquad \qquad \qquad - 4 b^{\bn}(b-1) - 4(1+\dots+b^{\sfrac \bn 2-1})^2(b-1)^2 \\
        &\iff \beta < \frac{1}{3b^{\bn}} \cdot \frac{1+b+\dots+b^{\sfrac \bn 3 -1}}{1+b+\dots+b^{\sfrac \bn 3 +1}} - \frac{2}{1+b+\dots+b^{\sfrac \bn 3 +1}} - \frac{2(b-1)(1+\dots+b^{\sfrac \bn 2-1})^2}{(1+b+\dots+b^{\sfrac \bn 3 +1})b^{\bn}} \, ,
    \end{align*}
    where we have written out the quantity at the beginning in terms of $\la_q\approx a^{(b^q)}$ and then compared exponents on both sides. It is easy to generalize the above to
    \begin{equation}\label{delta:useful}
        \delta_\qbn \left( \la_q \la_{q+k}^{-1}  \right)^{\sfrac 23} \la_{\qbn+1}^4 \la_\qbn^{-4} \frac{\la_q^4\la_\qbn^4}{\la^8_{q+\half}} < \delta_{q+\bn + k + 2} \qquad \forall k\geq \sfrac{\bn}{3} \, .
    \end{equation}
    
    \item\label{item:ugly:1} A consequence of the second inequality in \eqref{ineq:b:first} is that
    \begin{align*}
      &\frac{\de_{q+\bn}}{\de_{q+\bn-1}} \left(\frac{\sfrac{\la_{q+\half}}{\la_\qbn}}{\sfrac{\la_{q+\half-1}}{\la_{\qbn-1}}}\right)^{\sfrac43} < \frac{\de_{q+2\bn}}{\de_{q+2\bn-1}} \\
      &\qquad \iff -2\beta b^{\bn} + 2\beta b^{\bn-1} + (b^{\sfrac{\bn}{2}}-b^{\bn}) (b-1)\frac{4}{3b} < -2\beta b^{2\bn} + 2\beta b^{2\bn-1} \\
      &\qquad \iff 2\beta b^{\bn-1} \left( b^{\bn+1} - b - b^{\bn} + 1 \right) < (b^{\bn}-b^{\sfrac \bn 2})(b-1)\frac{4}{3b} \\
      &\qquad \iff \beta < \frac{2}{3b^{\sfrac \bn 2}} \frac{1+\dots+b^{\sfrac \bn 2 -1}}{1+\dots+b^{\bn-1}} \, .
    \end{align*}
    \item\label{item:ugly:2} A consequence of the definition of $\la_q$ is that for $q'\geq q-\sfrac \bn 2+1$, 
    \begin{equation}
        \frac{\la_{q'+\half}\la_{q+\half}}{\la_q\la_{q'+\bn}} < 1 \, . \label{eq:sup:growth}
    \end{equation}
    Indeed when $q'=q-\sfrac \bn 2 +1$, the inequality reduces to $\la_{q+1}\la_q^{-1} \la_{q+\half} \la_{q+\half+1}^{-1}<1$, which is an immediate consequence of the super-exponential growth; larger $q'$ are similar.
    \item\label{item:ugly:3} We have that $\delta_q \la_q^{\sfrac 23} < \delta_{q'}\la_{q'}^{\sfrac 23}$ for all $q'>q$. A stronger inequality is that for all $k\geq 1$, $\delta_{q+\bn}\la_q^{\sfrac 23} < \delta_{q+k+\bn} \la_{q+k}^{\sfrac 23}$, which is in fact equivalent to $\beta < \sfrac{1}{3b^{\bn}}$, which is implied by the first inequality in \eqref{ineq:b}. A final consequence of both inequalities is \begin{equation}
        \delta_\qbn \frac{\la_q^{\sfrac 23}}{\la_\qbn^{\sfrac 23}} < \delta_{q+2\bn} \implies \delta_\qbn^{\sfrac 12} \delta_q^{\sfrac 12} \frac{\la_q}{\la_\qbn} < \delta_{q+2\bn} \implies \delta_\qbn^{\sfrac 12} \delta_q^{\sfrac 12} \frac{\la_q}{\la_\qbn} \frac{\la_{q+\half}^{\sfrac 13}\la_\qbn^{-\sfrac 13}}{\la_{q+\half-1}^{\sfrac 13}\la_{\qbn-1}^{-\sfrac 13}} < \delta_{q+2\bn} \, .
    \end{equation}
    \item From the second inequality in \eqref{ineq:b:first}, we have that 
    \begin{align*}
       \beta < \frac{2}{3b^{\sfrac{\bn}{2}}} \cdot \frac{1+\dots+b^{\sfrac{\bn}{2}-1}}{1+\dots+b^{\bn-1}} \implies \delta_{q+\bn} \la_{q+\half}^{\sfrac 43} < \delta_{q+2\bn} \la_{\qbn}^{\sfrac 43} \, .
    \end{align*}
    \end{enumerate}
    \item\label{i:par:3} Choose $\CLebesgue=\frac{6+b}{b-1}$.
    \item\label{i:par:4} Define $\Ga_q$, $r_q$, $\tau_q$, and $\La_q$ by\footnote{The same type of comparability that we have in \eqref{eq:la:comparability} holds for $\Ga_q$ as defined in \eqref{eq:deffy:of:gamma}.}
    \begin{align}
        \Ga_q &= 2^{\left \lceil \varepsilon_\Gamma \log_2 \left( \frac{\la_{q+1}}{\la_q} \right) \right \rceil} \approx \left( \frac{\la_{q+1}}{\la_q} \right)^{\varepsilon_\Gamma} \approx \la_{q}^{(b-1)\varepsilon_\Gamma} \, , \qquad r_q = \frac{\la_{q+\half}\Ga_q}{\la_\qbn} \label{eq:deffy:of:gamma} \\
        \tau_q^{-1} &= \delta_q^{\sfrac 12} \la_q r_{q-\bn}^{-\sfrac 13} \Ga_q^{35} \, , \qquad \qquad \qquad \qquad \qquad \quad \qquad  \La_q = {\la_q \Ga_q^{10} } \, , 
    \end{align}
    where we choose $0<\varepsilon_\Gamma\ll (b-1)^2<1$ such that
\begin{subequations}
    \begin{align}
    (\delta_{q-\bn} \delta_{q-\bn-1}^{-1})^{\sfrac{1}{10}} \Ga_\qbn^{1000} &\leq 1 \, ,  \label{eq:delta:helping:matt}\\
    \label{ineq:tau:q}
    \Gamma_q^{{25}} \la_q \delta_q^{\sfrac 12}r_{q-\bn}^{-\sfrac 13} &\leq \tau_q^{-1} {\leq  \Gamma_q^{{50}} \la_q \delta_q^{\sfrac 12}r_{q-\bn}^{-\sfrac 13}} \, , \qquad \Gamma_\qbn^{300} \tau_{\qbn-1}^{-1} \leq \tau_\qbn^{-1} \, , \\
    \label{eq:pr.subtle}
    \Ga_{q+\bn}^{25} \frac{\de_{q+\bn}}{\de_{q+\bn-1}} \left(\frac{r_q}{r_{q-1}}\right)^{\sfrac43} &< \frac{\de_{q+2\bn}}{\de_{q+2\bn-1}} \\
      \lambda_{n-1}^{-2} \lambda_n \lhalf &\leq \Gamma_q^{-1} \quad \textnormal{for $q+\half+3\leq n \leq q+\bn+2$} \, , \label{ineq:stinky:3} \\
    \label{eq:dodging:parameterz}
     \Gamma_{q+\bn}^3 \Gamma_q^{-2} \frac{\lambda_{q'+\half}\lambda_{q+\half}}{\lambda_q \lambda_{q'+\bn}} &\leq 1 \qquad \textnormal{for all $q'$ such that $q+\half +1 - \bn \leq q' \leq q$} \, , \\
    \label{la.beats.de}
    \left(\frac{\la_q}{\la_{q'}}\right)^{\sfrac23} \Ga_{q+\bn}^{2000+10\CLebesgue} &< \left(\frac{\de_q}{\de_{q'}}\right)^{-1} \\
    \label{eq:par:div:1}
    r_q^{\sfrac 43} \delta_{q+\bn} \Gamma_q^{600} &\leq \delta_{q+2\bn} \implies r_q^{\sfrac 43} \Ga_q^{600} \delta_{q} \leq \delta_{q+\bn} \\
    \left( \frac{r_{q+1}}{r_q} \right) \Ga_\qbn^{1000+10\CLebesgue} &\leq 1 \label{ineq:r's:eat:Gammas} \\
     \Gamma_q^{5\CLebesgue+300} \delta_{q+\bn}^{\sfrac 12} r_q^{\sfrac 13} \lambda_{q+\bn}^{-1} \tau_q^{-1} &\leq \Ga_{q+\bn}^{-10} \delta_{q+2\bn} \, , \label{ineq:transport:basic}\\
    \Ga_{q+\bn} \delta_{q+\bn-1}^{-\sfrac 12} r_{q-1}^{-\sfrac 23} &\leq \de_{q+\bn}^{-\sfrac 12}r_q^{-\sfrac 23} \label{ineq.eps.gamma.imax} \, , \\
    \Ga_\qbn^{1000} &< \min\left(\la_q \la_\qbn^{-1}  r_q^{-2}, \la_q^{-\sfrac{1}{10}} \la_{q+1}^{\sfrac{1}{10}}, \delta_{q}^{\sfrac{1}{10}} \delta_{q+1}^{-\sfrac{1}{10}} \right) \label{eq:prepping:badshaq} \\
    \left \lceil \frac{(b^{\sfrac \bn 2-1}+\dots+b+1)^2}{\varepsilon_\Gamma(b^{\bn-1}+\dots+b+1)} \right \rceil &\geq 20 \, , \qquad 2000 \varepsilon_\Gamma b^{\bn} < 1 \, . \label{eq:prepping:badshaq:2}
    \end{align}
    \end{subequations}
    Indeed we have that the first inequality in \eqref{ineq:tau:q} is immediate, the second is possible since $\tau_q^{-1}$ is increasing in $q$, \eqref{eq:pr.subtle} is possible due to item~\eqref{item:ugly:1}, \eqref{ineq:stinky:3} and \eqref{eq:prepping:badshaq:2} are possible from immediate computation, \eqref{eq:dodging:parameterz} is possible due to item~\eqref{item:ugly:2}, \eqref{la.beats.de}, \eqref{eq:par:div:1}, and \eqref{ineq:transport:basic} are possible due to item~\eqref{item:ugly:3}, \eqref{ineq:r's:eat:Gammas}, \eqref{ineq.eps.gamma.imax}, and \eqref{eq:delta:helping:matt} are possible since $r_{q}$ and $\delta_q$ are decreasing in $q$, and \eqref{eq:prepping:badshaq} is possible due to \eqref{eq:deffy:of:gamma} and the super-exponential growth, which shows that $\la_q \la_\qbn^{-1} \la_{q+\half}^{-2} \la_{\qbn}^2, \la_q \la_{q+2} \la_{q+1}^{-2} > 1$.
    \item\label{i:par:5} Choose $\badshaq$ as\index{$\badshaq$}
    \begin{equation}\label{eq:badshaq:choice}
        \badshaq = 3\left \lceil \frac{(b^{\sfrac \bn 2}-1)^2}{(b-1)^2\varepsilon_\Gamma(b^{\sfrac \bn 2 -1}+\dots+b+1)}  + \frac{2000 b^{\bn}}{b^{\sfrac \bn 2}-1} + \frac{4b^{\bn-1}}{(b-1)\varepsilon_\Gamma(1+\dots+b^{\sfrac \bn 2 -1})} \right \rceil \, .
    \end{equation}
    As a consequence of this definition and \eqref{eq:prepping:badshaq:2}, we have that
    \begin{equation}
         10\leq \badshaq \, . \label{eq:badshaq:is:bad}
    \end{equation}
    We furthermore have that for all $\sfrac \bn 2 \leq k \leq \bn$,
    \begin{align*}
        \Ga_q^{\badshaq} &\la_q^2 \la_{q+k}^4 \la_{q+\half}^{-4} \la_{q+k}^{2} \la_{q+k-1}^{-4} < \Ga_{q+\half}^{\badshaq} \\
        &\iff 2\left( 1 - 2b^{\sfrac \bn 2} + b^{k} +2 b^{k} - 2b^{k-1} \right) < \badshaq (b-1)\varepsilon_\Gamma (b^{\sfrac{\bn}{2}}-1) \\
        &\impliedby 2\left( 1 - 2b^{\sfrac k 2} + b^{k} +2 b^{\bn} - 2b^{\bn-1} \right) < \badshaq (b-1)\varepsilon_\Gamma (b^{\sfrac{\bn}{2}}-1) \\
        &\iff 2 \left( b^{\sfrac k 2}-1 \right)^2 + 4b^{\bn-1}(b-1) < \badshaq (b-1)^2 \varepsilon_\Gamma (1+\dots+b^{\sfrac \bn 2 -1}) \\
        &\iff \frac{2\left(b^{\sfrac k 2}-1\right)^2}{(b-1)^2\varepsilon_\Gamma(1+\dots+b^{\sfrac \bn 2 -1})} + \frac{4b^{\bn-1}}{(b-1)\varepsilon_\Gamma (1+\dots+b^{\sfrac \bn 2 -1})} < \badshaq \, ,
    \end{align*}
    which is implied by \eqref{eq:badshaq:choice}. As a consequence of the above inequality, \eqref{eq:prepping:badshaq:2}, \eqref{eq:prepping:badshaq}, and \eqref{eq:badshaq:choice}, we have that for all $\sfrac \bn 2 \leq k \leq \bn$,
      \begin{subequations}
    \begin{align}
            \label{eq:par:div:2}
        \Gamma_q^{\badshaq} \leq \Gamma_{q+\half}^{\badshaq} \Ga_\qbn^{-2000} \, , \qquad \Gamma_q^{\badshaq+500} \Lambda_q \left( \frac{\la_{q+k}}{\la_{q+\half}} \right)^2 \lambda_{q+k-1}^{-2}\la_{q+k} &\leq \Gamma_{q+\half}^{\badshaq} \Ga_\qbn^{-200}  \, .
    \end{align}
    \end{subequations}
    \item\label{item:choice:of:alpha} Choose $\alpha=\alpha(q)\in (0,1)$ such that
    \begin{align}
     \lambda_{q+\bn}^\alpha = \Gamma_q^{\sfrac 1{10}} \, . \label{eq:choice:of:alpha}
    \end{align}
    \item\label{i:par:tau} Choose $\Tau_q$ according to the formula
    \begin{align}
    \frac12 \Tau_{q-1}^{-1} = \tau_q^{-1} \Gamma_q^{\badshaq+100} \delta_q^{-\sfrac 12} r_q^{-\sfrac 23}  + \Gamma_q^{\badshaq+100}\delta_q^{-\sfrac 12}r_q^{-1} \Lambda_q^3 \, . \label{v:global:par:ineq}
    \end{align}
    \item\label{i:par:4.5} Choose $\Npr$ such that\index{$\Npr$}
    \begin{align}\label{defn:Npr}
    \Ga_{q+\Npr}\La_{q+\bn}^4 \leq \Ga_{q+\Npr+1} \, .
    \end{align}
    \item\label{i:par:6} Choose $\NcutSmall
    $ and $\NcutLarge$ such that\index{$\NcutLarge$}\index{$\NcutLarge$}
    \begin{subequations}\label{condi.Ncut0}
    \begin{align}
     \NcutSmall &\leq \NcutLarge \, , \label{condi.Ncut0.1} \\
     \la_\qbn^{200}\left(\frac{\Ga_{q-1}}{\Ga_q}\right)^{\frac{\NcutSmall}5}
    &\leq \min\left(\la_{q+\bn}^{-4} \de_{q+3\bn}^2, \Ga_{q+\bn}^{-\badshaq-17-\CLebesgue} \delta_{q+3\bn}^2r_q \right) \, , \label{condi.Ncut0.2} \\
    \delta_{q+\bn}^{-\sfrac 12} r_q^{-1} \Ga_{q+\bn}^{\sfrac{\badshaq}{2}+16+\CLebesgue} \left( \frac{\Ga_{q+\bn-1}}{\Ga_{q+\bn}}\right)^{\NcutLarge} &\leq \Ga_{q+\bn}^{-1} \, . \label{condi.Ncut0.3}
    \end{align}
    \end{subequations}
    \item\label{i:par:7} Choose $\Nindt$ such that
    \begin{align}\label{condi.Nindt}
        \Nindt \geq \NcutSmall, \quad       
        \Ga_q^{-\Nindt} (\tau_q^{-1}\Ga_q^{i+40})^{-\NcutSmall-1} 
        (\Tau_q^{-1}\Ga_q)^{\NcutSmall+1}\leq 1 \, .
    \end{align}
    \item\label{i:par:12} Choose $N_{\rm g}, N_{\rm c}$ so that
\begin{subequations}\label{eq:darnit}
\begin{align}
    \Gamma_{q-1}^{-N_{\rm g}} \Gamma_q^2 &\leq \Gamma_{q+1} \Tau_{q+1}^{50\Nindt} \delta_{q+3\bn}^3 \, , \label{eq:darnit:2} \\
    2(\Tau_{q+\bn-1}^{-1}{\Ga_{q+\bn-1}^{10}})^{5\Nindt} \Gamma_{q+\bn}^{2\badshaq+\CLebesgue+100} r_q^{-2} \Gamma_{q-1}^{-\sfrac{N_{\rm c}}{2}} &\leq \Gamma_{q+\bn}^{-N_{\rm g}}
    \delta_{q+3\bn}^3\tau_{q+\bn-1}^{50\Nindt} \, , \label{eq:darnit:3} \\
    N_{\rm g} &\leq N_{\rm c} \leq \frac{\Nind}{40} \, . \label{eq:darnit:1} \footnotemark
\end{align}
\end{subequations}\footnotetext{This inequality is independent from the first two, and can be ensured by a large choice of $\Nind$ in the next step.  Since all the inequalities in \eqref{eq:darnit} are used together, we break the order slightly and include \eqref{eq:darnit:1} in this bullet point.}
    \item\label{i:par:8} Choose $\Nind$ such that \eqref{eq:darnit:1} is satisfied and\index{$\Nind$}
    \begin{subequations}
    \begin{align}
    \Nindt &\leq \Nind \, , \label{eq:Nindsy:1} \\
    \left(\Ga_{q-1}^{\Nind} \Ga_q^{-\Nind}\right)^{\sfrac{1}{10}} &\leq \delta_{q+5\bn}^3 \Gamma_q^{-2\badshaq-3} r_q \, . \label{eq:Nind:darnit}
    \end{align}
    \end{subequations}
    \item\label{i:par:9} Choose $\Ndec$ such that\index{$\Ndec$}
    \begin{align}\label{condi.Ndec0}
        (\la_{q+\bn+2}\Ga_q)^4 
        &\leq \left(
        \frac{\Ga_q^{\sfrac{1}{10}}}{{4\pi}}
        \right)^\Ndec \, , \qquad \qquad 
        \Nind \leq \Ndec \, .
    \end{align}
    \item\label{i:par:9.5} Choose $K_\circ$ large enough so that\index{$K_\circ$}
\begin{equation}\label{ineq:K_0}
    \lambda_q^{-K_\circ} \leq \delta_{q+3\bn}^{3} \Tau_{q+\bn}^{{5\Nind}} \lambda_{q+\bn+2}^{-100} \, .
\end{equation}
\item\label{i:par:10} Choose $\dpot$ and $N_{**}$ such that\index{$\dpot$}\index{$N_{**}$}
\begin{subequations}
    \begin{align}
     2\dpot + 3 &\leq N_{**} \, , \label{ineq:Nstarstar:dpot} \\
     \lambda_{q+\bn}^{100} \Gamma_q^{-\sfrac{\dpot}{200}} \Lambda_{q+\bn+2}^{5+K_\circ} \left( 1 + \frac{\max(\la_\qbn^2\Tau_q^{-1},\Lambda_q^{\sfrac 12}\Lambda_{q+\bn})}{\tau_q^{-1}} \right)^{20\Nind} &\leq \Tau_\qbn^{200\Nindt} \, , \label{ineq:dpot:1} \\
     \lambda_{q+\bn}^{100} \Gamma_q^{-\sfrac{N_{**}}{20}} \Lambda_{q+\bn+2}^{5+K_\circ} \left( 1 + \frac{\max(\la_\qbn^2\Tau_q^{-1},\Lambda_q^{\sfrac 12}\Lambda_{q+\bn})}{\tau_q^{-1}} \right)^{20\Nind} &\leq \Tau_\qbn^{20\Nindt} \, . \label{ineq:Nstarz:1}
\end{align}
\end{subequations}
\item\label{i:par:11} Choose $\Nfin$ such that\index{$\Nfin$}
\begin{subequations}
    \begin{align}
        2\Ndec + 4 + 10\Nind &\leq \sfrac{\Nfin}{40000} - \dpot^2 - 10\NcutLarge - 10\NcutSmall - N_{**} - 300 \, . \label{condi.Nfin0}
        \end{align}
    \end{subequations}
\item\label{i:choice:of:a} Having chosen all the parameters mentioned in items~\eqref{i:par:1}--\eqref{i:par:11} except for $a$, there exists a sufficiently large parameter $a_*$ such that $a_*^{(b-1)\varepsilon_\Gamma b^{-2\bn}}$ is at least fives times larger than \emph{all} the implicit constants throughout the paper, as well as those which have been suppressed in the computations in this section. Choose $a$ to be any natural number larger than $a_*$.\index{$a_*$}
\end{enumerate}

\subsection{A few more inequalities}
For all $q+\half-1 \leq m \leq m' \leq  \qbn$, we have that
\begin{align}
 &\Ga_q^{500+5\CLebesgue} \la_q \left( \frac{\delta_{q+\bn}}{\delta_{m+\bn}} \right)^{\sfrac 32} \La_q^{\sfrac 23} \left(\la_{m'-1}^{-2} \la_{m'}\right)^{\sfrac 23} \left( \frac{\min(\la_{m},\la_\qbn)\Ga_q}{\la_{q+\half}} \right)^{\sfrac 43}  \la_{m-1}^{-2} \la_{m} \leq  \Ga_q^{-250} \, , \label{ineq:in:the:afternoon}
\end{align}
and
\begin{align}
 &\Ga_q^{500+5\CLebesgue} \La_q \left( \frac{\min(\la_{m'},\la_\qbn)}{\la_{q+\half}} \right)^{\sfrac 23} \left( \frac{\delta_{q+\bn}}{\delta_{m+\bn}} \right)^{\sfrac 32} \La_q \la_{m'-1}^{-2} \la_{m'} \left( \frac{\min(\la_{m},\la_\qbn)\Ga_q}{\la_{q+\half}} \right)^{\sfrac 43}  \la_{m-1}^{-2} \la_{m} \leq  \Ga_q^{-250} \, . \label{ineq:in:the:afternoon:ER}
\end{align}
We claim the first inequality is morally equivalent to 
$$ \la_q \left( \frac{\delta_{q+\bn}}{\delta_{m+\bn}} \right)^{\sfrac 32} \la_q^{\sfrac 23}  \left( \min(\la_{m}, \la_\qbn ) \right)^{\sfrac 23} \la_{q+\half}^{-\sfrac 43}  \la_{m}^{-1} \leq 1 \, . $$
This equivalence is due to \eqref{delta:useful} (used to absorb a feq meaningless losses of $\la_k \la_{k-1}^{-1}$) and \eqref{la.beats.de} (used to absorb $\Ga_{\qbn}^{2000+10\CLebesgue}$, which itself can be absorbed in on meaningless loss of $\la_k \la_{k-1}^{-1}$ from \eqref{eq:prepping:badshaq}). Checking the simplified inequality then boils down to applying \eqref{delta:useful}. We leave further details to the reader. The second inequality is morally equivalent to 
$$  \la_q \left( \frac{\la_{m'}}{\la_{q+\half}} \right)^{\sfrac 23} \left( \frac{\delta_{\qbn}}{\delta_{m+\bn}} \right)^{\sfrac 32} \la_q \la_{m'}^{-1} \la_m^{-1} \left( \frac{\la_m}{\la_{q+\half}} \right)^{\sfrac 43} \leq 1 \, ,  $$
which can be checked by again using similar reasoning.

At this point, we list a number of additional inequalities, each of which can be checked by similar reasoning as the two inequalities above.  We leave further details to the reader.
\begin{subequations}
\begin{align}
    \label{ineq:in:the:morning}
    \la_q \Ga_q^{250} \La_q^{\sfrac 23} \left( \frac{r_{q+\half+1}}{r_q} \right)^{\sfrac 23} \la_{q+\half}^{-\sfrac 23} \left( \frac{\la_{q+\half+1}\Ga_q}{\la_{q+\half}} \right)^{\sfrac 43} \la_{q+\half}^{-1} \delta_\qbn^{{\sfrac 32}} \leq \delta^{{\sfrac 32}}_{q+\bn+\half+1} \, , \\
    \label{ineq:in:the:morning:ER}
    \la_q \Ga_q^{250+5\CLebesgue} \La_q \la_{q+\half}^{{-1}} \left( \frac{\la_{q+\half+1}\Ga_q}{\la_{q+\half}} \right)^{2} \la_{q+\half}^{-1} \delta_\qbn^{{\sfrac 32}} \leq \delta_{q+\bn+\half+1}^{{\sfrac 32}} \, , \\
\label{eq:desert:ineq}
   \delta_\qbn  \Ga_q^{500} \La_q^{\sfrac 23} \left(\la_{m-1}^2\lambda_m^{-1}\right)^{-\sfrac 23} \leq \delta_{m+\bn} \quad \textnormal{for} \quad q+\half-5 \leq m \leq q+\bn+5 \, , \\
\label{ineq:osc:general}
    \delta_{q+\bn} \La_q \Gamma_q^{400+5\CLebesgue} \left(\frac{\la_{m}}{\la_{q+\bn}r_q}\right)^{2/3} \lambda_{m-1}^{-2} \lambda_{m} \leq \Ga_m^{-9} \de_{m+\bn} \, , \\
\label{o:5}
\frac{\de_{q+\bn}}{\de_{{m}+\bn}}\Gamma_q^{200+5\CLebesgue} \left(\frac{\min(\la_m, \la_{q+\bn})}{\la_{q+\bn}r_q} \right)^{\sfrac 23} \Lambda_q \lambda_{m-1}^{-2} \min(\la_m, \la_{q+\bn})
\leq \Ga_{q+\half}^{-100} \, .
\end{align}
\end{subequations}

\appendix

\section{Appendix and toolkit}\label{sec:app:tools}

The appendix serves a number of purposes.  First, we prove general $L^p$ decoupling lemmas in subsection~\ref{appetizer:1}.  Then in subsection~\ref{sec:operator:iterates}, we recall a number of lemmas from \cite{BMNV21, NV22} which handle sums, iterates, and commutators of different differential operators.  Then in subsection~\ref{appetizer:3}, we construct and prove estimates for the various inverse divergence operators used throughout the proofs of Theorems~\ref{thm:main} and~\ref{thm:main:ER}.  Subsection~\ref{appetizer:4} contains a general lemma which allows us to upgrade material derivative estimates from $\Dtq$ to $D_{t,k}$ for $k>q$.  Finally, subsection~\ref{appetizer:5} contains a general mollification lemma which we apply whenever we need to estimate a mollified function and its difference with the original function.

\subsection{Decoupling lemmas and consequences of the Fa\`a di Bruno formula}\label{appetizer:1}

We begin with an $L^p$ decoupling\index{decoupling} lemma in the spirit of that from \cite{BMNV21}.  Some adjustments to the proof are required to treat the cases $p\neq 1,2,\infty$ and $d\neq 3$, as well as the slight adjustment to the assumption \eqref{eq:decoup:cond:alt} on the high-frequency function, which provides a slight increase in generality.  Note that the first inequality in \eqref{e:decoup:con:1} is implied by the second and the assumption that $\lambda \geq 2$, and so in practice we shall only check the second inequality.
\begin{lemma}[\bf $L^p$ decoupling]\label{lem:decoup}
Let $\Ndec,\kappa,\lambda\geq 1$ be such that
\begin{equation}\label{e:decoup:con:1}
    \left( 2 \cdot \frac{2\pi \sqrt{d}}{\kappa} \right) \cdot \lambda \leq \frac{2}{3} \, , \qquad \lambda^{\Ndec+d+1} \cdot \left( 2 \cdot \frac{2\pi \sqrt{d}}{\kappa} \right)^\Ndec \leq 1 \, .
\end{equation}
Let $p\in [1,\infty)$, and for $d\geq 1$, let $f$ be a $\T^d$-periodic function such that there exists $\const_f$ such that for all $0\leq j \leq \Ndec+d+1$,
\begin{equation}\label{e:decoup:cond:2}
    \left\| D^j f \right\|_{L^p} \leq \const_f \lambda^j \, .
\end{equation}
Let $g$ be a $\T^d$-periodic function and $\const_g > 0$ a constant such that for any cube $T$ of side-length $\sfrac{2\pi}{\kappa}$,
\begin{align}\label{eq:decoup:cond:alt}
{\kappa^{\sfrac dp}} \left\| g \right\|_{L^p(T)} \leq \const_g \, .
\end{align}
Then there exists a dimensional constant $C=C(p,d)$ which is independent of $f$ and $g$ such that
\begin{equation}\label{e:decoup:conc}
    \left\| f g \right\|_{L^p(\T^d)} \leq C(p,d) \const_f \const_g \, .
\end{equation}
\end{lemma}
\begin{proof}[Proof of Lemma~\ref{lem:decoup}]
Let $\{T_j\}_j$ be disjoint cubes of side-length $\sfrac{2\pi}{\kappa}$ such that $$ \bigcup_j T_j = \T^d \, . $$
For any Lebesgue integrable function $h$, let 
$$  \Bar{h}_j := \dashint_{T_j} h(x) \, dx \, .  $$
Note that from Jensen's inequality, we have that
\begin{equation}\label{e:convexity}
    | \Bar{h}_j |^p = \left| \dashint_{T_j} h(x) \, dx \right|^p \leq \dashint_{T_j} |h(x)|^p \, dx = \overline{|h|^p}_j \, .
\end{equation}
For any $x\in T_j$, we have that
\newcommand{\barf}{\Bar{f}}
\begin{align}
    |f(x)|^p &= \left( |\Bar{f}_j| + |f(x) - \Bar{f}_j | \right)^p \notag\\
    &\leq 2^p \left( |\Bar{f}_j|^p + |f(x) - \Bar{f}_j|^p \right) \notag\\
    &\leq 2^p \left( |\Bar{f}_j|^p + \left(\sup_{x\in T_j} |f(x) - \barf_j |\right)^p \right) \notag\\
    &\leq 2^p \left( |\barf_j|^p + \left( \frac{2\pi\sqrt{d}}{\kappa} \sup_{T_j} |Df| \right)^p \right) \notag\\
    &\leq 2^p \overline{|f|^p}_j + 2^p \left( \frac{2\pi\sqrt{d}}{\kappa} \right)^p \sup_{T_j} |Df|^p \, ,
\end{align}
where in the last line we have used \eqref{e:convexity}. Iterating, we obtain
\begin{align}
    |f(x)|^p &\leq 2^p \overline{|f|^p}_j + 2^p \left( \frac{2\pi\sqrt{d}}{\kappa}\right)^p \left( 2^p \overline{|Df|^p}_j + 2^p \left( \frac{2\pi\sqrt{d}}{\kappa} \right)^p \sup_{T_j} |D^2 f|^p \right) \notag\\
    &\leq \sum_{m=0}^{\Ndec-1} 2^{(m+1)p} \left( \frac{2\pi\sqrt{d}}{\kappa}\right)^{mp} \overline{|D^m f |^p}_j + \left( 2 \cdot \frac{2\pi\sqrt{d}}{\kappa} \right)^{\Ndec p} \left\| D^\Ndec f \right\|_{L^\infty}^p \, . \notag
\end{align}
Multiplying by $g$, integrating over $T_j$, and using \eqref{eq:decoup:cond:alt}, we obtain\footnote{Note that in the third line, we move the average from $|D^m f|^p$ to $|g|^p$. In the fourth line, we used the assumption \eqref{eq:decoup:cond:alt} on $g$. In the second to last line, we used the assumption \eqref{e:decoup:con:1}.}
\begin{align}
\left\| fg \right\|_{L^p}^p &= \sum_j \int_{T^j} |fg|^p \notag\\
&\leq \sum_j \int_{T_j} |g|^p \sum_{m=0}^{\Ndec-1} 2^{(m+1)p} \left( \frac{2\pi\sqrt{d}}{\kappa}\right)^{mp} \overline{|D^m f |^p}_j + \left( 2 \cdot \frac{2\pi\sqrt{d}}{\kappa} \right)^{\Ndec p} \left\| D^\Ndec f \right\|_{L^\infty}^p \const_g^p \notag\\
&= \sum_j \dashint_{T_j} |g|^p \sum_{m=0}^{\Ndec-1} 2^{(m+1)p} \left( \frac{2\pi\sqrt{d}}{\kappa}\right)^{mp} \left\| D^m f \right\|_{L^p(T_j)}^p + \left( 2 \cdot \frac{2\pi\sqrt{d}}{\kappa} \right)^{\Ndec p} \left\| D^\Ndec f \right\|_{L^\infty}^p \const_g^p \notag\\
&\leq (C(d))^p \const_g^p \sum_{m=0}^{\Ndec-1} 2^{(m+1)p} \left( \frac{2\pi\sqrt{d}}{\kappa}\right)^{mp} \const_f^p \lambda^{mp} + \left( 2 \cdot \frac{2\pi\sqrt{d}}{\kappa} \right)^{\Ndec p} \left( C'(d) \const_f \lambda^{\Ndec+d+1} \const_g \right)^p  \notag\\
&\leq (C(d))^p \const_g^p 2^p \cdot 3 \cdot \const_f ^p + (C'(d))^p \const_f^p \const_g^p  \notag\\
&=: (C(p,d))^p \const_f^p \const_g^p \, .
\end{align}
Taking $p^{\rm th}$ roots on both sides concludes the proof.
\end{proof}

We now recall the multivariable Fa\`a di Bruno formula (see for example the appendix in \cite{BMNV21}). Let $g= g(x_1,\ldots,x_d) = f ( h(x_1,\ldots,x_d))$, where $f \colon \R^m \to \R$, and $h \colon \R^d \to \R^m$ are $C^n$ functions. Let $\alpha \in {\mathbb N}_0^d$ be such that $|\alpha|=n$, and let $\beta \in \N_0^m$ be such that $1\leq |\beta|\leq n$. We then define  
\begin{align} 
p(\alpha,\beta) 
&= \Bigg\{ (k_1,\ldots,k_n; \ell_1,\ldots,\ell_n) \in (\N_0^m)^n \times (\N_0^d)^n \colon \exists s \mbox{ with }1\leq s \leq n \mbox{ s.t. } \notag \\
&\qquad   \qquad  |k_j|, |\ell_j| > 0 \Leftrightarrow 1 \leq j \leq s, \, 0 \prec \ell_{1} \prec \ldots \prec \ell_s,  \sum_{j=1}^s k_j =  \beta, \sum_{j=1}^s |k_j| \ell_j = \alpha \Bigg\} \, .
\label{eq:Faa:di:Bruno:1}
\end{align} 
The multivariable Fa\`a di Bruno formula states that
\begin{align}
\partial^\alpha g (x) = \alpha! \sum_{|\beta|=1}^n (\partial^\beta f)(h(x)) \sum_{p(\alpha,\beta)} \prod_{j=1}^n \frac{(\partial^{\ell_j} h(x))^{k_j}}{k_j! (\ell_j!)^{k_j}}\, .
\label{eq:Faa:di:Bruno:2}
\end{align}


Throughout this manuscript, we must estimate only finitely many derivatives.  Therefore we ignore the factorials in \eqref{eq:Faa:di:Bruno:2} and absorb them into the implicit constant of the symbol ``$\les$." We now recall the following lemma from \cite{BMNV21}, which gives a useful consequence of the Fa\`a di Bruno formula.

\begin{lemma}[\bf Compositions with flow maps]\label{l:useful_ests}
Given a smooth function $f \colon \R^d \times \R \to \R$, suppose that for $\lambda \geq 1$ the vector field $\Phi \colon \R^d \times \R \to \R^d$ satisfies the estimate
\begin{align}
\norm{D^{N+1} \Phi}_{L^{\infty}(\supp f)}&\les \lambda^{N}
\label{e:Faa_di_Bruno_0}
\end{align}
for $0 \leq N \leq N_*$.  Then for any $1 \leq N\leq N_*$ we have
\begin{align}
\abs{D^N\left( f \circ \Phi \right)(x,t)} \lesssim& \sum_{m=1}^N \lambda^{N-m}\abs{(D^m f)\circ \Phi (x,t)}
\label{e:Faa_di_Bruno_1} 
\end{align}
and thus trivially we obtain
\begin{align*}
\abs{D^N\left( f \circ \Phi \right)(x,t)} \lesssim& \sum_{m=0}^N \lambda^{N-m}\abs{(D^m f)\circ \Phi (x,t)}.
\end{align*}
for any $0 \leq N\leq N_*$.
\end{lemma}

Many estimates will require estimates for derivatives of products of functions which decouple and which are composed with a diffeomorphism.  The proof is a minor variation on \cite[Lemma~A.7]{BMNV21}.

\begin{lemma}[\bf Decoupling with flow maps]
\label{l:slow_fast}
Let ${p\in[1,\infty]}$, and fix integers $N_* \geq M_* \geq \Ndec\geq 1$. Fix $d\geq 2$ and $f \colon \R^d \times \R \to \R$, and let $\Phi \colon \R^d\times \R \to \R^d$ be a vector field satisfying $D_{t}\Phi = (\partial_t + v\cdot \nabla) \Phi=0$. Denote by $\Phi^{-1}$ the inverse of the flow $\Phi$,  which is the identity at a time slice which intersects the support of $f$. Assume that for some $\lambda , \tau^{-1}, \Tau^{-1} \geq 1$ and $\const_f>0$  the function $f$ satisfies the estimates
\begin{align}
\norm{D^N D_{t}^M f}_{L^{p}}&\les \const_f \lambda^N\MM{M,N_{t},\tau^{-1},\Tau^{-1}} \label{eq:slow_fast_0}
\end{align}
for all $N\leq N_*$ and $M\leq M_*$, and that $\Phi$ and $\Phi^{-1}$ are bounded for all $N\leq N_*$ by
\begin{align}
\norm{D^{N+1} \Phi}_{L^{\infty}(\supp f)}&\les \lambda^{N}\label{eq:slow_fast_1} \\
\norm{D^{N+1} \Phi^{-1}}_{L^{\infty}(\supp f)} &\les  \lambda^{N} \, .
\label{eq:slow_fast_2}
\end{align}
Lastly, suppose that there exist $\varrho:\T^d\rightarrow \R$ and parameters  $\Lambda \geq \Upsilon \geq \mu$  and $\const_\varrho>0$ such that for any cube $T$ of side length $\mu^{-1}$,
\begin{align}
\frac{1}{\mu^{\sfrac dp}} \left\| D^N \varrho \right\|_{L^p(T)} + \norm{D^N \varrho}_{L^p(\T^d)} \les \const_\varrho \MM{N, N_x, \Upsilon, \Lambda} \, \label{eq:slow_fast_4}
\end{align}
for all $0\leq N \leq N_*$.  If the parameters $$\lambda\leq \mu \leq \Upsilon \leq \Lambda$$ satisfy
\begin{align}
 { \Lambda^{d+1}  \leq  \left( \frac{\mu}{4 \pi \sqrt{3} \lambda}\right)^{\Ndec}   } \label{eq:slow_fast_3}
\,,
\end{align}
and we have
\begin{equation}\label{eq:slow_fast_3_a}
 {2 \Ndec + d + 1 \leq N_*}  \,,
\end{equation}
then for $N\leq N_*$ and $M\leq M_*$ we have the bound
\begin{align}
\norm{D^N  D_{t }^M \left( f\; \varrho\circ \Phi
\right)}_{L^p}
&\les \const_f \const_\varrho \MM{N,  {N_x, \Upsilon}, \Lambda}  \MM{M, {N_{t}},\tau^{-1},\Tau^{-1}} \, .
\label{eq:slow_fast_5} 
\end{align}
\end{lemma} 

\begin{remark}
\label{rem:slow:fast}
We note that if estimate \eqref{eq:slow_fast_0} is known to hold for $N+M \leq N_\circ$ for some $N_\circ \geq 2\Ndec + d+1$ (instead of $N\leq N_*$ and $M\leq M_*$), and if the bounds \eqref{eq:slow_fast_1}--\eqref{eq:slow_fast_2} hold for all $N\leq N_\circ$, then it follows from the method of proof that the bound \eqref{eq:slow_fast_5} holds for $N+ M \leq N_\circ$ and $M \leq N_\circ - 2 \Ndec - d-1$. The only modification required is that instead of considering the cases $N' \leq N_* - \Ndec - d-1$ and $N' > N_* - \Ndec - d-1$, we now split into $N' + M \leq N_\circ - \Ndec - d-1$ and $N' + M > N_\circ  - \Ndec - d-1$. In the second case we use that $N-N'' \geq N_0 - M - \Ndec - d-1 \geq \Ndec$, where the last inequality holds precisely because $M \leq   N_\circ - 2 \Ndec - d-1$.
\end{remark}

\begin{proof}[Proof of Lemma~\ref{l:slow_fast}]
Since $D_{t}\Phi = 0$ we have $D_{t}^M (\varrho\circ \Phi) = 0$. 
Furthermore, since $\div\, v\equiv 0$, we have that $\Phi$ and $\Phi^{-1}$ preserve  volume.  Then using Lemma~\ref{l:useful_ests}, which we may apply due to \eqref{eq:slow_fast_1}, we have
\begin{align}
\norm{D^ND_{t}^M\left( f \; \varrho\circ \Phi \right)}_{L^p} 
& \lesssim \sum_{N'=0}^N \norm{D^{N'} D_{t}^M f \; D^{N-N'}\left(  \varrho \circ \Phi \right)}_{L^p}
\notag \\
&
\lesssim \sum_{N'=0}^N \sum_{N''=0}^{N-N'}\lambda^{N-N'-N''}\norm{D^{N'} D_{t}^Mf \;(D^{N''}\varrho)\circ \Phi}_{L^p}
\notag \\
&
\lesssim \sum_{N'=0}^N \sum_{N''=0}^{N-N'}\lambda^{N-N'-N''} \norm{\left(D^{N'} D_{t}^Mf\right) \circ \Phi^{-1} D^{N''}\varrho}_{L^p}.
\label{eq:slow_fast_temp_1}
\end{align}
In \eqref{eq:slow_fast_temp_1} let us first consider the case $N' \leq N_* - \Ndec - d - 1$.  Due to assumption \eqref{eq:slow_fast_2}, we may apply Lemma \ref{l:useful_ests}, and appealing to \eqref{eq:slow_fast_0} we have that
\begin{align}
\norm{D^{n}\left( (D^{N'} D_{t}^Mf) \circ (\Phi^{-1},t)\right)}_{L^{p}}&\les \sum_{n'=0}^{n}\lambda^{n-n'}\norm{(D^{n'+N'} D_{t}^M f)\circ \Phi^{-1}  }_{L^{p}} \notag \\
&\les  \const_f \sum_{n'=0}^{n}\lambda^{n-n'}\lambda^{n'+N'}\MM{M,N_{t},\tau^{-1},\Tau^{-1}} \notag \\
&\les  \left(\const_f \lambda^{N'}\MM{M,N_{t},\tau^{-1},\Tau^{-1}}\right) \lambda^n
\label{eq:slow_fast_temp_2}
\end{align}
for all $n\leq \Ndec+d+1$. This bound matches \eqref{e:decoup:cond:2}, with $\const_f$ replaced by $\const_f \lambda^{N'} \MM{M,N_t,\tau^{-1},\Tau^{-1}}$. 
Since the function $D^{N''}\varrho$ satisfies \eqref{eq:slow_fast_4}, we may apply \eqref{eq:slow_fast_temp_2}, the fact that $\lambda \leq \Upsilon\leq \Lambda$, assumption \eqref{eq:slow_fast_3}, and Lemma~\ref{lem:decoup} to conclude that
\begin{align*}
\norm{\left(D^{N'} D_{t}^M f\right) \circ \Phi^{-1} D^{N''}\varrho}_{L^p}\les \const_f \lambda^{N'} \MM{M,N_t,\tau^{-1},\Tau^{-1}} \const_\varrho \MM{N'',N_x,\Upsilon,\Lambda}
\, .
\end{align*}
Inserting this bound back into \eqref{eq:slow_fast_temp_1} concludes the proof of \eqref{eq:slow_fast_5} for $N'\leq N_*-\Ndec - d - 1$ as considered in this case.
 
Next, let us consider the case  $N' > N_* - \Ndec - d - 1$. Since $0\leq N' \leq N$, in particular this implies that $N>N_* - \Ndec-d-1$. Using furthermore that $N'' \leq N - N'$ and \eqref{eq:slow_fast_3_a}, we also obtain that $N-N'' \geq N' > N_* - \Ndec-d-1 \geq \Ndec$. Then H\"older's inequality, the fact that $\Phi^{-1}$ is volume preserving, the Sobolev embedding $W^{d+1,1} \subset L^\infty$, the ordering $\Lambda\geq\Upsilon\geq \mu \geq 1$, and assumption \eqref{eq:slow_fast_3} implies that
\begin{align*}
\lambda^{N - N'-N''} \norm{\left(D^{N'} D_{t}^Mf\right) \circ \Phi^{-1} D^{N''}\varrho}_{L^p} 
&\les \lambda^{N - N'-N''} \norm{D^{N'} D_{t}^M f }_{L^p}  \norm{D^{N''}\varrho}_{L^\infty} 
\notag\\
&\les \lambda^{N - N'-N''}  \const_f \lambda^{N'}  \MM{M,N_{t},\tau^{-1},\Tau^{-1}} \const_\varrho  \MM{N'' + d+1 ,N_x,\Upsilon,\Lambda} 
\notag\\
&\les \const_f \const_\varrho \MM{N,N_x,\Upsilon,\Lambda} \MM{M,N_{t},\tau^{-1},\Tau^{-1}} \Lambda^{d+1}  \left(\frac{\lambda}{\Upsilon}\right)^{N-N''}  
\notag\\ 
&\les \const_f \const_\varrho \MM{N,N_x,\Upsilon,\Lambda} \MM{M,N_{t},\tau^{-1},\Tau^{-1}} \Lambda^{d+1}  \left(\frac{\lambda}{\mu}\right)^{\Ndec}   
\notag\\ 
&\les \const_f \const_\varrho \MM{N,N_x,\Upsilon,\Lambda} \MM{M,N_{t},\tau^{-1},\Tau^{-1}} \,.
\end{align*}
Combining the above estimate with \eqref{eq:slow_fast_temp_1}, we deduce that the bound \eqref{eq:slow_fast_5} holds also for $N' > N_* - \Ndec-d-1$, concluding the proof of the lemma.
\end{proof}

\subsection{Sums and iterates of operators and commutators with material derivatives}
\label{sec:operator:iterates}

We first record the following identity for material and spatial derivatives applied to functions raised to a positive integer power. 

\begin{lemma}[\bf Leibniz rule with material and spatial derivatives] \label{lem:leib}
Let $d\geq 2$ be given, $g:\T^d\rightarrow \R$ be a smooth function, $v:\T^d\times\R\rightarrow\R^d$ a divergence-free vector field, and set $D_t=\partial_t+v\cdot\nabla$, and $p\in\mathbb{N}$.  Fix $M,N\in \N$, and use $\alpha=(\alpha_1,\alpha_2,\dots,\alpha_p)$ and $\beta=(\beta_1,\beta_2,\dots,\beta_p)$ to denote multi-indices with $|\alpha|=N,|\beta|=M$. Then we have the identities
\begin{subequations}
\begin{align}
D^{N} D_{t}^{M} g^p 
&= \sum_{\left\{\substack{\alpha , \beta \, : \, \sum_{i=1}^p \alpha_i = N \, , \\ \sum_{i=1}^p \beta_i = M }\right\}} {\binom{N}{\alpha_1,\dots,\alpha_p}} {\binom{M}{\beta_1,\dots,\beta_p}} \prod_{i=1}^p D^{\alpha_i} D_t^{\beta_i} g \label{eq:leibniz:one} \\
p g^{p-1} D^{N} D_{t}^{M} g 
&=  D^N D_t^M g^p -  \sum_{\left\{\substack{\alpha , \beta  \, : \, \sum_{i=1}^p \alpha_i = N \, , \\ \sum_{i=1}^p \beta_i = M \, , \\  \alpha_i+\beta_i < N+M \, \forall \, i }\right\}} {\binom{N}{\alpha_1,\dots,\alpha_p}} {\binom{M}{\beta_1,\dots,\beta_p}} \prod_{i=1}^p D^{\alpha_i} D_t^{\beta_i} g \, . \label{eq:leibniz:two} 
\end{align}
\end{subequations}
\end{lemma}

We recall \cite[Lemma~A.10]{BMNV21}.  We have generalized the statement slightly so that it applies in $\T^d$ rather than just $\T^3$; in fact the statement and proof have nothing to do with the dimension.
\begin{lemma}
\label{lem:cooper:1}
Fix $N_x,N_t,N_*  \in \N$, $\Omega \in \T^d \times \R$ a space-time domain, and let $v$ be a vector field and $B$ a differential operator. For  $k\geq 1$ and $\alpha,\beta \in \N^k$ such that $|\alpha|+ |\beta| \leq N_*$, we assume that we have the bounds
\begin{align}
\norm{ \left(\prod_{i=1}^k D^{\alpha_i} B^{\beta_i} \right) v}_{L^\infty(\Omega)} \les \const_v \MM{|\alpha|,N_x,\lambda_v,\tilde \lambda_v} \MM{|\beta|,N_t,\mu_v,\tilde \mu_v}
\label{eq:cooper:v}
\end{align}
for some $\const_v\geq 0$, $1\leq \lambda_v \leq \tilde \lambda_v$, and $1\leq \mu_v \leq \tilde \mu_v$.
With the same notation and restrictions on $|\alpha|,|\beta|$, let $f$ be a function which for some  $p \in [1,\infty]$  obeys 
\begin{align}
\norm{ \left(\prod_{i=1}^k D^{\alpha_i} B^{\beta_i} \right) f}_{ L^p (\Omega)} \les \const_f \MM{|\alpha|,N_x,\lambda_f,\tilde \lambda_f}  \MM{|\beta|,N_t,\mu_f,\tilde \mu_f}
\label{eq:cooper:f}
\end{align}
for some $\const_f\geq 0$, $1\leq \lambda_f \leq \tilde \lambda_f$, and $1\leq \mu_f \leq \tilde \mu_f$. 
Denote
\begin{align*}
\lambda = \max\{ \lambda_f,\lambda_v\}, \quad \tilde \lambda= \max\{\tilde \lambda_f,\tilde \lambda_v\}, \quad \mu = \max\{\mu_f,\mu_v\}, \quad \tilde \mu = \max\{\tilde \mu_f,\tilde \mu_v\}.
\end{align*}
Then, for \[ A = v\cdot \nabla \] we have the bounds 
\begin{align}
\norm{D^n \left( \prod_{i=1}^{k} A^{\alpha_i} B^{\beta_i} \right) f}_{ L^p (\Omega)} 
&\les \const_f \const_v^{|\alpha|} \MM{n+|\alpha|,N_x,\lambda,\tilde \lambda}  \MM{|\beta|,N_t,\mu,\tilde \mu} \label{eq:cooper:f:**} \\
&\les \const_f  \MM{n,N_x,\lambda,\tilde \lambda} (\const_v \tilde \lambda)^{|\alpha|} \MM{|\beta|,N_t,\mu,\tilde \mu} 
\notag\\ 
&\les \const_f  \MM{n,N_x,\lambda,\tilde \lambda}  \MM{|\alpha|+|\beta|,N_t, \max\{\mu, \const_v \tilde \lambda\}, \max\{\tilde \mu,\const_v \tilde \lambda\} } 
\label{eq:cooper:f:***}
\end{align}
as long as $n+|\alpha|+|\beta|\leq N_*$.
As a consequence,  if $k=m$ then \eqref{eq:cooper:f:***} and an expansion of the operator $(A+B)^M$ imply that for all $n+m \leq N_*$, 
\begin{align}
\norm{D^n (A + B)^m f}_{ L^p (\Omega)} 
\les \const_f  \MM{n,N_x,\lambda,\tilde \lambda}  \MM{m,N_t, \max\{\mu, \const_v \tilde \lambda\}, \max\{\tilde \mu,\const_v \tilde \lambda\} } \, .
\label{eq:cooper:f:*}
\end{align}
\end{lemma}

A corollary of the previous lemma is the commutator lemma \cite[Lemma~A.14]{BMNV21}, which we now record along with several useful remarks.
\begin{lemma}
\label{lem:cooper:2}
 Let $p\in [1,\infty]$. 
Fix $N_x,N_t,N_*,M_* \in \N$, let $v$ be a vector field, let $D_t = \partial_t + v\cdot \nabla$ be the associated material derivative, and let $\Omega$ be a space-time domain. Assume that the vector field $v$ obeys 
\begin{align}
\norm{D^N D_t^M D v}_{L^\infty(\Omega)} \les \const_v \MM{N+1,N_x,\lambda_v,\tilde \lambda_v} \MM{M,N_t,\mu_v,\tilde \mu_v}
\label{eq:cooper:2:v}
\end{align}
for $N \leq N_*$  and $M \leq M_*$.
Moreover, let $f$ be a function which obeys
\begin{align}
\norm{D^N D_t^M f}_{L^p(\Omega)} \les \const_f \MM{N,N_x,\lambda_f,\tilde \lambda_f} \MM{M,N_t,\mu_f,\tilde \mu_f}
\label{eq:cooper:2:f}
\end{align}
for all $N\leq N_*$ and $M \leq M_*$. 
Denote
\begin{align*}
\lambda = \max\{ \lambda_f,\lambda_v\}, \quad \tilde \lambda= \max\{\tilde \lambda_f,\tilde \lambda_v\}, \quad \mu = \max\{\mu_f,\mu_v\}, \quad \tilde \mu = \max\{\tilde \mu_f,\tilde \mu_v\}.
\end{align*}
Let $m,n,\ell \geq 0$ be such that $n+\ell \leq N_*$ and $m\leq M_*$. 
Then, we have that the commutator $[D_{t}^m,D^n]$ is bounded  as
\begin{align}
\norm{D^\ell \left[ D_t^m,D^n \right] f}_{L^{p}(\Omega)} 
&\les \const_f  \const_v \tilde \lambda_v \MM{\ell+n,N_x,\lambda,\tilde \lambda}   \MM{m-1,N_t,\max\{\mu,\const_v \tilde \lambda_v\},\max\{\tilde \mu,\const_v \tilde \lambda_v\}} 
\label{eq:cooper:2:f:1:*}
\\
&\les \const_f \MM{\ell+n,N_x,\lambda,\tilde \lambda}   \MM{m ,N_t,\max\{\mu,\const_v \tilde \lambda_v\},\max\{\tilde \mu,\const_v \tilde \lambda_v\}}.
\label{eq:cooper:2:f:1}
\end{align}
Moreover, we have that for $k\geq 2$, and any $\alpha,\beta\in \N^k$  with $|\alpha|\leq N_*$ and $|\beta|\leq M_*$, the estimate
\begin{align}
\norm{\left( \prod_{i=1}^{k} D^{\alpha_i} D_t^{\beta_i}\right) f}_{L^{p}(\Omega)} 
&\les \const_f \MM{|\alpha|,N_x,\lambda,\tilde \lambda}   \MM{|\beta|,N_t,\max\{\mu,\const_v \tilde \lambda_v\},\max\{\tilde \mu,\const_v \tilde \lambda_v\}}
\label{eq:cooper:2:f:2}
\end{align}
holds.
\end{lemma}

\begin{remark}
\label{rem:cooper:2:sum}
If instead of \eqref{eq:cooper:2:v} and \eqref{eq:cooper:2:f} holding for $N\leq N_*$ and $M\leq M_*$, we know that both of these inequalities hold for all $N+M \leq N_\circ$ for some $N_\circ \geq 1$, then the conclusions of the Lemma hold as follows: the bounds \eqref{eq:cooper:2:f:1:*} and \eqref{eq:cooper:2:f:1} hold for $\ell+n+m\leq N_\circ$, while \eqref{eq:cooper:2:f:2} holds for $|\alpha|+|\beta| \leq N_\circ$. We refer to \cite{BMNV21} for further discussion.
\end{remark}

\begin{remark}
\label{rem:cooper:2}
If the assumption \eqref{eq:cooper:2:f} is replaced by 
\begin{align}
\norm{D^{N} D_t^M f}_{L^{p}(\Omega)} \les \const_f \MM{N-1,N_x,\lambda_f,\tilde \lambda_f} \MM{M,N_t,\mu_f,\tilde \mu_f}
\,,
\label{eq:cooper:2:f:alt}
\end{align}
whenever $1 \leq N\leq N_*$, then the conclusion \eqref{eq:cooper:2:f:2} instead becomes
\begin{align}
\norm{\left( \prod_{i=1}^{k} D^{\alpha_i} D_t^{\beta_i}\right) f}_{L^{p}(\Omega)} 
&\les \const_f \MM{|\alpha|-1,N_x,\lambda,\tilde \lambda}   \MM{|\beta|,N_t,\max\{\mu,\const_v \tilde \lambda_v\},\max\{\tilde \mu,\const_v \tilde \lambda_v\}}
\label{eq:cooper:2:f:2:alt}
\end{align}
whenever $|\alpha|\geq 1$. We again refer to \cite{BMNV21} for further discussion.
\end{remark}

\begin{remark}\label{rem:upgrade.material.derivative.end}
Fix $p\in [1,\infty]$, $N_x,N_t,N_*  \in \N$, and a space-time domain $\Omega \in \T^d \times \R$.
Define $D_t = \pa_t + (v\cdot \na)$ as in Lemma \ref{lem:cooper:2}.
Suppose that for  $k\geq 1$ and $\alpha,\beta \in \N^k$ such that $|\alpha|+ |\beta| \leq N_*$, we have the bounds
\begin{align}
\norm{ \left(\prod_{i=1}^k D^{\alpha_i} {D_t^{\beta_i}} \right) w}_{L^\infty(\Omega)} \les \const_w \MM{|\alpha|,N_x,\lambda_w,\tilde \lambda_w} \MM{|\beta|,N_t,\mu_w,\tilde \mu_w}
\label{eq:cooper:w}
\end{align}
for some $\const_w\geq 0$, $1\leq \lambda_w \leq \tilde \lambda_w$, and $1\leq \mu_w \leq \tilde \mu_w$.
Then, under the assumption \eqref{eq:cooper:2:v} and \eqref{eq:cooper:2:f} in Lemma \ref{lem:cooper:2} with $M_*=N_*$, we have that for all $N,M\leq N_*$, 
\begin{align}
\norm{D^N (D_t + (w\cdot \na))^M f}_{ L^p (\Omega)} 
\les \const_f  \MM{n,N_x,\lambda,\tilde \lambda}  \MM{m,N_t, \mu, \tilde \mu } \, 
\label{eq:cooper:f:mat}
\end{align}
where 
\begin{align*}
\lambda = \max\{ \lambda_f,\lambda_v, \lambda_w\}, \quad \tilde \lambda&= \max\{\tilde \lambda_f,\tilde \lambda_v,  \tilde \lambda_w\}, \quad \mu = \max\{\mu_f,\mu_v, \mu_w, 
\const_v\td\la_v,
\const_w \td\la_w \}, \\
\tilde \mu &= \max\{\tilde \mu_f,\tilde \mu_v, \tilde \mu_w, \const_v\td\la_v, \const_w \td\la_w\} \, .
\end{align*}
If \eqref{eq:cooper:2:v} and \eqref{eq:cooper:2:f} hold for $N+M\leq N_*$, as in Remark~\ref{rem:cooper:2:sum}, then \eqref{eq:cooper:f:mat} holds also for $N+M\leq N_*$.
\end{remark}

\subsection{Inversion of the divergence}\label{appetizer:3}

\begin{proposition}[\bf Inverse divergence iteration step]
\label{prop:Celtics:suck}
Let $n\geq 2$ be given. Fix a zero-mean $\T^n$-periodic function $\varrho$ and a zero-mean $\T^n$-periodic symmetric tensor field $\vartheta^{(i,j)}$ which are related by {$\varrho  =  \partial_{ij} \vartheta^{(i,j)}$}. Let $\Phi$ be a volume preserving diffeomorphism of $\T^n$. Define the matrix $A = (\nabla \Phi)^{-1}$. Given a vector field $G^k$, we have
\begin{align}
G^k (\varrho \circ \Phi)
= \partial_\ell R^{k \ell} + E^k
\label{eq:Celtics:suck:total}
\end{align}
where the symmetric stress $R^{k\ell}$ is given by
\begin{align}
R^{k\ell}
&= 
G^k A_i^\ell(\partial_j\vartheta^{(i,j)} \circ \Phi)
    +G^\ell A_i^k (\partial_j\vartheta^{(i,j)} \circ \Phi)
- G^n\pa_n\Phi^m A_i^k   A_j^\ell(\pa_m\vartheta^{(i,j)} \circ \Phi)
\, ,
\label{eq:Celtics:suck:stress}
\end{align}
and the error term $E^k$ is given by 
\begin{align}
E^k
&= -\pa_\ell (G^\ell A_i^k) (\partial_j\vartheta^{(i,j)} \circ \Phi)
    - (\pa_\ell G^k)A_i^\ell(\partial_j\vartheta^{(i,j)} \circ \Phi)
    + \pa_n (G^\ell A_i^k \pa_\ell\Phi^m) A_j^n(\pa_m\vartheta^{(i,j)} \circ \Phi)
\, .
\label{eq:Celtics:suck:error}
\end{align}
\end{proposition}

\begin{remark}[\bf Linearity with respect to $G$]\label{rem:linear:inverse:div}
From \eqref{eq:Celtics:suck:stress} and \eqref{eq:Celtics:suck:error}, it is clear that the symmetric stress and error term are \emph{linear} in $G$; more precisely, each term of the symmetric stress and error may be written as a product of flow maps, high frequency functions, and a single component of either $G$ or $\nabla G$. This will be a useful observation when determining the support properties of the symmetric stresses and error terms.
\end{remark}

\begin{proof}[Proof of Proposition~\ref{prop:Celtics:suck}]  

By the definition of $A$, we have $A^n_\ell \pa_k\Phi^\ell = \de_{nk}$, and the volume-preserving property of $\Phi$ gives the Piola identity $\pa_n A^n_\ell =0$. These then imply a useful identity $(\pa_\ell \ph) \circ \Phi
= \pa_n (A^n_\ell (\ph\circ\Phi))$. Using this, we first get
\begin{align*}
    G^k (\varrho \circ \Phi)
    &=G^k (\partial_i \partial_j \vartheta^{(i,j)} \circ \Phi)
    =G^k \pa_\ell( A_i^\ell(\partial_j\vartheta^{(i,j)} \circ \Phi))
    =\pa_\ell (G^k A_i^\ell(\partial_j\vartheta^{(i,j)} \circ \Phi))
    - (\pa_\ell G^k)A_i^\ell(\partial_j\vartheta^{(i,j)} \circ \Phi)\\
    &=\pa_\ell (G^k A_i^\ell(\partial_j\vartheta^{(i,j)} \circ \Phi)
    +G^\ell A_i^k (\partial_j\vartheta^{(i,j)} \circ \Phi))
    -G^\ell A_i^k \pa_\ell\Phi^m (\pa_m\partial_j\vartheta^{(i,j)}) \circ \Phi\\
    &\quad-\pa_\ell (G^\ell A_i^k) (\partial_j\vartheta^{(i,j)} \circ \Phi)
    - (\pa_\ell G^k)A_i^\ell(\partial_j\vartheta^{(i,j)} \circ \Phi) \, .
\end{align*}
In the last equality,the first two terms match the first two terms in $\partial_\ell R^{k\ell}$, while the last two terms will go into the error term $E^k$. To deal with the remaining term, we use
\begin{align*}
 G^\ell A_i^k \pa_\ell\Phi^m (\pa_m\partial_j\vartheta^{(i,j)}) \circ \Phi
 &=G^\ell A_i^k \pa_\ell\Phi^m \pa_n (A_j^n(\pa_m\vartheta^{(i,j)} \circ \Phi))\\
 &=\pa_n (
 G^\ell\pa_\ell\Phi^m A_i^k   A_j^n(\pa_m\vartheta^{(i,j)} \circ \Phi))
 - \pa_n (G^\ell A_i^k \pa_\ell\Phi^m) A_j^n(\pa_m\vartheta^{(i,j)} \circ \Phi) \, .
\end{align*}
Indeed, plugging this identity into the second term, we obtain the symmetric stress $R^{k\ell}$ and error term $E^k$. Note that the first term above is symmetric due to the assumed symmetry of $\vartheta^{(i,j)}$. 
\end{proof}

With the iterative step in hand, we can now state the proposition which contains our main inverse divergence algorithm.  The spirit of the statement and proof is similar to the corresponding statements and proofs in \cite{BMNV21,NV22}, modulo minor adjustments.  After stating the main proposition, we record a number of useful remarks which follow from the proof.\index{inverse divergence}

\begin{proposition}[\bf Main inverse divergence operator]
\label{prop:intermittent:inverse:div}
Let dimension $n\geq 2$ and Lebesgue exponent $p\in[1,\infty]$ be free parameters. The remainder of the proposition is composed first of \emph{low and high-frequency assumptions}, which then produce a \emph{localized output} satisfying a number of properties.  Finally, the proposition concludes with \emph{nonlocal assumptions and output}. 
\smallskip

\noindent\textbf{Part 1: Low-frequency assumptions}
\begin{enumerate}[(i)]
\item\label{item:cond.G.inverse.div} Let $G$ be a vector field and assume there exist a constant $\const_{G,p} > 0$ and parameters 
\begin{equation}\label{eq:inv:div:NM}
N_*\geq M_*\geq 1 \, ,
\end{equation}
$M_t$, and $\lambda, \nu,\nu' \geq 1$ such that 
\begin{align}
\norm{D^N D_{t}^M G}_{L^p}&\lesssim \const_{G,p} \lambda^N\MM{M,M_{t},\nu,\nu'}
\label{eq:inverse:div:DN:G}
\end{align}
for all $N \leq N_*$ and $M \leq M_*$.
\item Fix an incompressible vector field $v(t,x):\R\times\T^n\rightarrow \R^n$ and denote its material derivative by $D_t = \partial_t + v\cdot\nabla$. Let $\Phi$ be a volume preserving diffeomorphism of $\T^n$ such that 
\begin{align}
D_t \Phi = 0 \,
\qquad \mbox{and} \qquad
\norm{\nabla \Phi - \Id}_{L^\infty(\supp G)} \leq \sfrac 12 \,. \label{eq:DDpsi2}
\end{align}
Denote by $\Phi^{-1}$ the inverse of the flow $\Phi$,  which is the identity at a time slice which intersects the support of $G$.
Assume that  the velocity field $v$ and the flow functions $\Phi$ and $\Phi^{-1}$ satisfy the bounds 
\begin{subequations}
\begin{align}
\norm{D^{N+1}   \Phi}_{L^{\infty}(\supp G)} + \norm{D^{N+1}   \Phi^{-1}}_{L^{\infty}(\supp G)} 
&\les \lambda'^{N}
\label{eq:DDpsi}\\
\norm{D^ND_t^M D v}_{L^{\infty}(\supp G)}
&\les \nu \lambda'^{N}\MM{M,M_{t},\nu,\nu'}
\label{eq:DDv}
\,,
\end{align}
\end{subequations}
for all $N \leq N_*$, $M\leq M_*$, and some $\lambda'>0$. 
\end{enumerate}
\smallskip

\noindent\textbf{Part 2: High-frequency assumptions}
\begin{enumerate}[(i)]
\item\label{item:inverse:i} Let $\varrho \colon \T^n \to \R$ be a zero mean scalar function such that there exists a large positive even integer $\dpot \gg 1$ and a smooth, mean-zero, adjacent-pairwise symmetric tensor potential\footnote{We use $i_j$ for $1\leq j \leq \dpot$ to denote any number in the set $\{1,\dots,n\}$. We refer to Lemma~\ref{lem:special:cases} for the meaning of adjacent-pairwise symmetric.} $\vartheta^{(i_1,\dots, i_\dpot)}:\T^n \rightarrow \R^{\left(n^\dpot\right)}$ such that $\varrho(x) = \partial_{i_1}\dots\partial_{i_\dpot} \vartheta^{(i_1 \dots i_\dpot)}(x)$.
\item \label{item:inverse:ii} There exists a parameter $\mu\geq 1$ such that $\varrho$ and $\vartheta$ are $(\sfrac{\T}{\mu})^n$-periodic.
\item \label{item:inverse:iii} There exist parameters $1 \ll \Upsilon \leq \Upsilon' \leq \Lambda$, $\const_{*,p}>0$ such that for all $0\leq N \leq {N_*}$ and all $0\leq k \leq \dpot$,
\begin{align}
\norm{D^N \partial_{i_1}\dots \partial_{i_k} \vartheta^{(i_1,\dots, i_\dpot)}}_{L^p} \les \const_{*,p} \Upsilon^{k-\dpot} \MM{N, \dpot - k , \Upsilon', \Lambda} \, .
\label{eq:DN:Mikado:density}
\end{align} 
\item\label{item:inverse:iv} There exists $\Ndec$ such that the above parameters satisfy
\begin{align}
 \lambda', \lambda \ll \mu \leq \Upsilon \leq \Upsilon' \leq \Lambda  \,, \qquad \max(\lambda,\lambda') \Upsilon^{-2} \Upsilon' \leq 1 \, , \qquad N_*-\dpot \geq 2\Ndec + n+1 \, , 
 \label{eq:inverse:div:parameters:0}
\end{align}
where by in the first inequality in \eqref{eq:inverse:div:parameters:0} we mean that 
\begin{align}
\Lambda^{n+1} \left(\frac{\mu}{2\pi \sqrt{3} \max(\lambda,\lambda')}\right)^{-\Ndec} \leq 1
\,.
\label{eq:inverse:div:parameters:1}
\end{align}
\end{enumerate}
\smallskip

\noindent\textbf{Part 3: Localized output}
\begin{enumerate}[(i)]
\item\label{item:div:local:0} There exists a symmetric tensor $R$ and a vector field $E$ such that\index{$\divH$}
\begin{align}
G \; \varrho\circ \Phi  &=  \div R + E  
=: \div\left( \divH \left( G \varrho \circ \Phi \right) \right) + E \, . \label{eq:inverse:div}
\end{align}
We use the notation $R=\divH( G \varrho \circ \Phi)$ for the symmetric stress.
\item\label{item:div:local:i} The support of $R$ is a subset of $\supp G \cap \supp \vartheta$.
\item\label{item:div:local:ii} There exists an explicitly computable positive integer $\const_\divH$, an explicitly computable function $r(j):\{0,1,\dots,\const_\divH\}\rightarrow \mathbb{N}$ and explicitly computable tensors
\begin{align*}
    &\rho^{\beta(j)} \, , \qquad \beta(j)=(\beta_1,\beta_2,\dots,\beta_{r(j)})\in  \{1,\dots,n\}^{r(j)} \, , \\
    &H^{\alpha(j)} \, , \qquad \alpha(j)=(\alpha_1,\alpha_2,\dots,\alpha_{r(j)},k,\ell) \in \{1,\dots,n\}^{r(j)+2} \, 
\end{align*}
of rank $r(j)$ and $r(j)+2$, respectively, all of which depend only on $G$, $\varrho$, $\Phi$, $n$, $\dpot$, such that the following holds.  The symmetric, localized stress $R$ can be decomposed into a sum of symmetric, localized stresses as\footnote{The contraction is on the first $r(j)$ indices, and the resulting rank two tensor is symmetric.}
\begin{align}\label{eq:divH:formula}
    \divH^{k\ell} (G \varrho \circ \Phi) = R^{k\ell} = \sum_{j=0}^{\const_\divH}  H^{\alpha(j)} \rho^{\beta(j)} \circ \Phi \, .
\end{align}
Furthermore, we have that
\begin{align}
\supp H^{\alpha(j)} \subseteq \supp G  \, , \qquad \supp \rho^{\beta(j)} \subseteq \supp \vartheta \,  . \label{eq:inverse:div:linear}
\end{align}
\item\label{item:div:local:iii} For all $N \leq N_* - \sfrac \dpot 2$, $M\leq M_*$, and $j\leq \const_\divH$, we have the subsidiary estimates\footnote{In fact it is clear from the algorithm that as $j$ increases, the estimates become much stronger.  For simplicity's sake we simply record identical estimates for each term which are sufficient for our aims.}
\begin{subequations}\label{eq:inverse:div:sub:main}
\begin{align}
    \left\| D^N \rho^{\beta{(j)}} \right\|_{L^p} &\les \const_{*,p} \Upsilon^{-2} \Upsilon' \MM{N,1,\Upsilon',\Lambda} \label{eq:inverse:div:sub:1} \\
    \left\|D^N D_{t}^M H^{\alpha{(j)}}\right\|_{L^p} &\lesssim \const_{G,p} \left(\max(\lambda,\lambda')\right)^N \MM{M,M_{t},\nu,\nu'} \, . \label{eq:inverse:div:sub:2} 
\end{align}
\end{subequations}
\item\label{item:div:local:iv} For all $N \leq N_* - \sfrac \dpot 2$ and $M\leq M_*$, we have the main estimate
\begin{align}
\norm{D^N D_{t}^M R}_{L^p}
 &\les  \const_{G,p} \const_{*,p}  \Upsilon' \Upsilon^{-2} \MM{N,1,\Upsilon',\Lambda} \MM{M,M_{t},\nu,\nu'} 
\label{eq:inverse:div:stress:1}
\end{align}
\item\label{item:div:nonlocal} For $N \leq N_* - \sfrac \dpot 2 $ and $M\leq M_*$ the error term $E$  in \eqref{eq:inverse:div} satisfies
\begin{align}
\norm{D^N D_{t}^M E}_{L^p}  
\les \const_{G,p} \const_{*,p}   \max(\lambda,\lambda')^{\sfrac \dpot 2} \left( \Upsilon' \Upsilon^{-2}\right)^{\sfrac \dpot 2} \Lambda^{N} \MM{M,M_{t},\nu,\nu'} 
\,.
\label{eq:inverse:div:error:1}
\end{align}
\end{enumerate}
\smallskip

\noindent\textbf{Part 4: Nonlocal assumptions and output}
\begin{enumerate}[(i)]
\item 
\label{item:nonlocal:v}
Let $N_\circ, M_\circ$ be integers such that 
\begin{equation}\label{eq:inv:div:wut}
1 \leq M_\circ \leq N_\circ \leq \sfrac{M_*}{2} \, ,
\end{equation}
and let $K_\circ$ be a positive integer.\footnote{{$K_\circ$ serves as an extra amplitude gain which will be used later to eat some material derivative losses.}}  Assume that in addition to the bound \eqref{eq:DDv} we have the following global lossy estimates
\begin{align}
\norm{D^N \partial_t^M v}_{L^\infty}\les  \const_v \lambda'^N \nu'^M
\label{eq:inverse:div:v:global}
\end{align}
for all  $M \leq M_\circ$ and $N+M \leq N_\circ + M_\circ$, where 
\begin{align}
\const_v \lambda' \les \nu' 
\,.
\label{eq:inverse:div:v:global:parameters}
\end{align}
\item Assume that $\dpot $ is large enough so that\index{$\dpot$}\index{$\divR$}
\begin{align}
\const_{G,p} \const_{*,p} \max(\lambda,\lambda')^{\sfrac \dpot 4} (\Upsilon' \Upsilon^{-2})^{\sfrac \dpot 4}  \Lambda^{{n}+2+K_\circ} \left(1 + \frac{\max\{ \nu', \const_v \Lambda \}}{\nu 
}\right)^{M_\circ}
\leq 1
\, .
\label{eq:riots:4}
\end{align}
\end{enumerate}

Then we may write  
\begin{align}
E = \div \RR_{\rm nonlocal} + \fint_{\T^3} G \varrho \circ \Phi \, dx =: \div \left(\divR(G \varrho \circ \Phi)\right) + \fint_{\T^3} G \varrho \circ \Phi \,  dx \, ,
\label{eq:inverse:div:error:stress}
\end{align}
where $\RR_{\rm nonlocal} = \divR(G \varrho \circ \Phi)$ is a traceless symmetric stress which satisfies
\begin{align}
\norm{D^N D_{t}^M \RR_{\rm nonlocal} }_{L^\infty}  
\leq  \frac{1}{\Lambda^{K_\circ}} \max(\lambda,\lambda')^{\sfrac \dpot 4} (\Upsilon' \Upsilon^{-2})^{\sfrac \dpot 4} \Lambda^N \nu^M
\label{eq:inverse:div:error:stress:bound}
\end{align}
for  $N \leq N_\circ$ and $M\leq M_\circ$.
\end{proposition}
\begin{remark}[\bf Lossy derivatives on $v$ and estimates for $\overline{R}_{\rm nonlocal}$]\label{rem:lossy:choices}
Let us specify the estimates we expect to obtain from \eqref{eq:inverse:div:error:stress:bound} for the nonlocal error term $\overline{R}_{\rm nonlocal}$. For our applications, we need to choose parameters so that the estimate reads
\begin{equation}\label{eq:R:nonlocal:in:practice}
   \left\| D^N D_t^M \overline{R}_{\rm nonlocal} \right\|_{L^\infty} \leq \lambda_{q+\bn}^{-10} \delta_{q+3\bn}^2 {\Tau_{q+\bn}^{4\Nindt}} \la_\qbn^N \tau_q^{-M} 
\end{equation}
for $N,M\leq 2\Nind$. We therefore choose $N_\circ=M_\circ=2\Nind$, and since in applications $M_*$ will be at least $\sfrac{\Nfin}{10000}$, we have from \eqref{condi.Nfin0} that $M_\circ \leq N_\circ \leq \sfrac{M_*}{2}$.  Next, we choose $K_\circ$ large enough so that $\lambda_q^{-K_\circ} \leq \delta_{q+3\bn}^{2} \Tau_{q+\bn}^{4\Nindt}\la_\qbn^{-100}$, which follows from \eqref{ineq:K_0}. The lossy estimates in \eqref{eq:inverse:div:v:global} follow from the inductive assumption \eqref{eq:bobby:old} with $\const_v=\Lambda_q^{\sfrac 12}$; note that \eqref{eq:inverse:div:v:global:parameters} is precisely \eqref{v:global:par:ineq}. Finally, the inequality in \eqref{eq:riots:4} will be a consequence of our choices of $\lambda,\lambda',\Upsilon',\Upsilon$, which from \eqref{ineq:stinky:3} give a gain of at least $\Gamma_q^{-\lfloor \sfrac \dpot {40} \rfloor}$, and \eqref{ineq:dpot:1}.

\end{remark}

\begin{remark}[\bf Special case for negligible error terms] \label{rem:inverse.div.spcial}
The inverse divergence operator defined in the proposition can be applied to an input without the structure of low and high frequency parts when $\varrho=1$ and $\const_{G,p}$ are sufficiently small. More precisely, we 
keep the low-frequency assumption (Part 1), replace the high-frequency assumptions (Part 2) with $\varrho =1$, and set $\Upsilon= \Upsilon'= \Lambda = \max(\la,\la')$, $\const_{*,p}=1$, $\dpot =0$. \footnote{Since we do not need decoupling, $\mu$ does not need to be specified.} Then, as long as $\const_{G,p}$ is small enough to satisfy \eqref{eq:riots:4}, the conclusions in Part 4 hold. In particular, we have that
\begin{align*}
    G = \div \divR G + \fint_{T^3} G \, dx \, . 
\end{align*}
Note that $\divR G = \mathcal{R}G$ in the special case, where $\mathcal{R}$ is the usual inverse divergence operator defined in \eqref{eq:RSZ}. 
\end{remark}

\begin{remark}[\bf High frequency part of the output as a potential]\label{rem:potential.rho}
In order to obtain the conclusions in Remarks \ref{rem.ct.osc}, \ref{rem.ct.tn}, and \ref{rem.ct.divcorr}, we need to write $\rho^{\be(j)}$ as a potential. This can be done if 
the potentials $\vartheta^{(i_1, \cdots, i_{\dpot})}$ used in the application of the inverse divergence in Section~\ref{sec:ER:main} can be written as 
$\vartheta^{(i_1, \cdots, i_{\dpot})} = \pa_{i_{\dpot+1} \cdots i_{2\dpot}} \theta^{(i_{\dpot+1}, \cdots, i_{2\dpot})}$, where $ \theta$ satisfies
\begin{align*}
    \norm{D^N \div^k \theta^{(i_1, \cdots, i_{2\dpot})}}_{L^p} 
    \lec \const_{*,p} \Upsilon^{k-2\dpot}\MM{N, 2\dpot-k, \Upsilon', \La}
\end{align*}
for $0\leq k \leq 2\dpot$ and $N\leq N_*$. This is easily ensured by \emph{initially} choosing $\varrho$ as $\varrho= \pa_{i_1 \cdots i_{2\dpot}} \theta^{(i_1, \cdots, i_{2\dpot})}$, where we save half of the divergences for later to enable the application of the inverse divergence algorithm a second time, as will be done in for the transport/Nash current errors in~\ref{op:tnce}. Since the inverse divergence algorithm shows that $\rho^{\be(j)}$ consists of spatial derivatives and divergences of $\vartheta$, it is clear that $\rho^{\be(j)}$ can be written in potential form as $\rho^{\be(j)} = \pa_{i_{\dpot+1}\cdots i_{\dpot+k}}\ov\theta^{(i_{\dpot+1},\cdots, i_{k}, \be(j))}$ for some potential
$\ov\theta^{(i_{\dpot+1},\cdots, i_{k}, \be(j))}$. Furthermore, we have
\begin{align*}
    \norm{D^N \pa_{i_{\dpot+1}\cdots i_{\dpot+k}}\ov\theta^{(i_{\dpot+1},\cdots, i_{\dpot+ k}, \be(j))}}_{L^p}
    \lec \const_{*,p} (\Upsilon^{-2} \Upsilon') \Upsilon^{k-\dpot}\MM{N, \dpot-k+1, \Upsilon', \La}
\end{align*}
for $0\leq k \leq \dpot$ and $N\leq N_* - \sfrac{\dpot}2$.
\end{remark}

\begin{remark}[\bf Mean of the error term]\label{rem:est.mean} We claim that the mean $\langle  G(\varrho\circ\Phi) \rangle$ satisfies
\begin{align}\notag
    \bigg|\frac{d^M}{dt^M} \langle  G(\varrho\circ\Phi) \rangle
    \bigg|\leq \La^{-K_\circ}{(\max(\la, \la')\Upsilon^{-1})^{\frac34\dpot} }\MM{M,M_t, \nu, \nu'}
\end{align}
for $M\leq M_\circ$. To see this, first note that since $v$ is incompressible, $\frac{d^M}{dt^M} \langle  G(\varrho\circ\Phi) \rangle = \langle (D_t^{M} G)(\varrho\circ\Phi) \rangle$.
Then using Lemma \ref{lem:decoup} with \eqref{eq:inverse:div:parameters:1}, \eqref{eq:inverse:div:DN:G}, \eqref{eq:DDpsi}, \eqref{eq:DN:Mikado:density}, and \eqref{eq:riots:4}, we have the desired estimate
\begin{align*}
\bigg|\int_{\T^3}(D_t^{M} G)(\varrho\circ\Phi) dx\bigg|
&=\bigg|\int_{\T^3}(D_t^{M} G)\circ \Phi^{-1} \div^\dpot \vartheta dx\bigg|
=\bigg|\int_{\T^3}\pa_{(i_1,\cdots, i_\dpot)}((D_t^{M} G)\circ \Phi^{-1}) \vartheta^{(i_1,\cdots, i_\dpot)} dx\bigg|\\
&\lec \norm{\pa_{(i_1,\cdots, i_\dpot)}((D_t^{M} G)\circ \Phi^{-1})}_1\norm{\vartheta^{(i_1,\cdots, i_\dpot)}}_1\\
&\lec \const_{G,p}  \const_{*,p}(\max(\la, \la'))^{\dpot} 
{\Upsilon^{-\dpot}}
\MM{M,M_t, \nu, \nu'}
\\
&\leq \La^{-K_\circ}
{(\max(\la, \la')\Upsilon^{-1})^{\frac34\dpot} }
\MM{M,M_t, \nu, \nu'}\, . 
\end{align*}
Inn particular, under the same choice of parameters suggested in Remark \ref{rem:lossy:choices}, we have
\begin{align*}
    \bigg|\frac{d^M}{dt^M} \langle  G(\varrho\circ\Phi) \rangle
    \bigg|\leq \la_{q+\bn}^{-10}\de_{q+3\bn}^2 \Tau_{q+\bn}^{4\Nindt}\tau_q^{-M}
\end{align*}
for $M\leq 2\Nind$. 
\end{remark}

\begin{remark}[\bf Inverse divergence for scalar fields]\label{rem:scalar:inverse:div}
Adjusting the above proposition so that $G$ is a scalar field and the output is a vector field is simple; one can make the substitution $\displaystyle G\rightarrow \left(G,\underbrace{0, \dots , 0}_{n-1 \,  \, 0\rm{'s}}\right)$, apply the Proposition to the newly constructed vector field, and take the first row or column of the symmetric stress as the output.
\end{remark}

\begin{remark}[\bf Inverse divergence with pointwise bounds]\label{rem:pointwise:inverse:div}
Let us consider the setting in which all the inductive assumptions from the proposition hold, or are adjusted according to Remark~\ref{rem:scalar:inverse:div}, but we know in addition that there exists a smooth, non-negative function $\pi$ such that
\begin{align}\label{eq:inv:div:extra:pointwise}
\left| D^N D_{t}^M G \right| &\lesssim \pi \lambda^N\MM{M,M_{t},\nu,\nu'} \, .
\end{align}
for $N\leq N_*$ and $M\leq M_*$. Then it is clear from the algorithm utilized in the proof that we may additionally conclude that
\begin{align}\label{eq:inv:div:extra:conc}
       \left|D^N D_{t}^M H^{\alpha_{(j)}}\right| &\lesssim \pi \left(\max(\lambda,\lambda')\right)^N \MM{M,M_{t},\nu,\nu'} 
\end{align}
for $N\leq N_*- \lfloor\sfrac \dpot 2 \rfloor$ and $M\leq M_*$.
\end{remark}

\begin{remark}[\bf Avoiding abuses of notations]\label{rem:bar:variables}
Proposition~\ref{prop:intermittent:inverse:div}, and indeed many of the other ``abstract nonsense" lemmas and propositions in the manuscript, are written using generic notations such as $\lambda,\const_{G,\sfrac 32}$, etc. Application of the lemma or proposition then requires specification of values for these various inputs.  Occasionally several such lemmas or propositions will be applied in succession; for example, repeated applications of the inverse divergence as in Corollary~\ref{cor:inverse.div}. In such situations, we shall add bars above all symbols in the statements of the ``abstract nonsense" lemmas, and then specify an input for the ``bar variable."  For example, applying Proposition~\ref{prop:intermittent:inverse:div} to a term from the sum in \eqref{eq:divH:formula} (which has the same form as the input of the inverse divergence, just with different parameters!) would be done using the parameter choices $\overline{\const}_{G,p}=\const_{G,p}$, $\overline{\lambda}=\max(\lambda,\lambda')$, $\overline{\const}_{*,p}=\const_{*,p}\Upsilon^{-2}\Upsilon'$, and $\overline{N}_*=N_*-\lfloor\sfrac \dpot 2 \rfloor$, which are all valid choices due to \eqref{eq:inverse:div:sub:main}.
\end{remark}

\begin{proof}[Proof of Proposition~\ref{prop:intermittent:inverse:div}]
We divide the proof into four steps. First, we collect some simple preliminary bounds.  Next, we apply Proposition~\ref{prop:Celtics:suck} the first time and show that an error term is produced which obeys the estimates required in \eqref{eq:inverse:div:stress:1}.  Afterwards we indicate how to apply the algorithm $\lfloor \sfrac \dpot 2 \rfloor-1$ more times to produce $R$ and $E$ obeying \eqref{eq:inverse:div:stress:1} and \eqref{eq:inverse:div:error:1}, respectively. By construction, both $R$ and $E$ will be supported in $\supp G \cap \supp \vartheta\circ \Phi$.  The support property for $R$ and the conclusions in \eqref{eq:divH:formula}, \eqref{eq:inverse:div:sub:main}, \eqref{eq:inverse:div:stress:1}, and \eqref{eq:inverse:div:error:1} will be proven along the way.  Finally, we outline how to obtain the bounds in \eqref{eq:inverse:div:error:stress:bound} for the nonlocal portion of the inverse divergence.  The entire proof follows closely the method of proof of \cite[Proposition~A.18]{BMNV21}, the main differences being the slight adjustment to the iteration step due to the difference between Proposition~\ref{prop:Celtics:suck} and \cite[Proposition~A.17]{BMNV21}, and the slightly more general assumption in~\eqref{eq:DN:Mikado:density} compared to \cite[A.69]{BMNV21}.  The only significant difference to the conclusion is that the amplitude gain is $\Upsilon' \Upsilon^{-2}$, cf.~\eqref{eq:inverse:div:stress:1} compared to \cite[A.73]{BMNV21}.
\smallskip

\noindent\texttt{Step 1: } An application of Lemma~\ref{lem:cooper:2}, or more precisely Remark~\ref{rem:cooper:2}, yields
\begin{align}
 \norm{D^{N''} D_t^M D^{N' }D \Phi}_{L^\infty(\supp G)} \les \lambda'^{N'+N''} \MM{M,M_t,\nu, \nu'}
 \label{eq:derivatives:phase:gradient}
\end{align}
whenever  $N'+N'' \leq N_*$  and $M\leq M_*$. We similarly obtain
\begin{align}
 \norm{D^{N''} D_t^M D^{N' } (D \Phi)^{-1}}_{L^\infty(\supp G)} \les \lambda'^{N'+N''} \MM{M,M_t,\nu, \nu'}
 \label{eq:derivatives:phase:gradient:inverse}
\end{align}
from the Fa'a di Bruno formula \eqref{eq:Faa:di:Bruno:2}, the formula for matrix inversion in $B_{\sfrac 12}(\Id)$, the Liebniz rule, and \eqref{eq:derivatives:phase:gradient}.  Another application of Lemma~\ref{lem:cooper:2} yields
\begin{align}
 \norm{D^{N''} D_t^M D^{N'}G}_{L^p} \les \const_{G,p} \lambda^{N'+N''} \MM{M,M_t,\nu, \nu'}
 \label{eq:derivatives:mixed:G}
\end{align}
whenever $N'+N''\leq N_*$ and $M \leq M_*$. These preliminary bounds are similar to those from the beginning of the proof of \cite[Proposition~A.18]{BMNV21}, and we refer there for further details. 
\smallskip

\noindent\texttt{Step 2: } For notational purposes, let $\varrho_{(0)} = \varrho$ and $\varrho_{(\dpot)}^{(i_1,\dots,i_\dpot)}=\vartheta^{(i_1,\dots,i_\dpot)}$, and for $1 \leq k < \dpot $ let $\varrho_{(k)}^{i_{\dpot-k+1},\dots, i_{\dpot}} = \partial_{i_1}\dots \partial_{i_{\dpot-k}} \vartheta^{(i_1,\dots,i_\dpot)}$.  Then $\varrho_{(k-1)} = \div \varrho_{(k)}$ (assuming contraction along the proper index, which we omit in a slight abuse of notation), and for any ``pairwise permutation"\footnote{We refer again to Lemma~\ref{lem:special:cases} for the meaning of this.} $\sigma:\{\dpot-k+1,\dots,\dpot\}\rightarrow\{\dpot-k+1,\dots,\dpot\}$, $\varrho_{(k)}^{i_{\dpot-k+1},\dots,i_\dpot}=\varrho_{(k)}^{i_{\sigma(\dpot-k+1)},\dots,i_{\sigma(\dpot)}}$, so that $\varrho_{(k)}$ is pairwise symmetric. We also define $G_{(0)} = G$. Since $\rho_{(0)} = \div \div \rho_{(2)}$ where $\rho_{(2)}$ is pairwise symmetric, we deduce from Proposition~\ref{prop:Celtics:suck}, identities \eqref{eq:Celtics:suck:total}--\eqref{eq:Celtics:suck:error} that 
\begin{align}
G_{(0)}^k \varrho_{(0)}\circ \Phi =  \partial_\ell R_{(0)}^{k\ell} + G_{(1)}^{ijkm} \partial_m \varrho_{(2)}^{(i,j)}\circ \Phi \, .
\label{eq:inverse:div:G0}
\end{align}
The symmetric stress $R_{(0)}$ is given by
\begin{align}
R_{(0)}^{k\ell} = \underbrace{\left(  G_{(0)}^k A^\ell_i \delta_{mj}  + G_{(0)}^\ell A^k_i \delta_{mj} - G^n \partial_n \Phi^m A^k_i A^\ell_j \right)}_{=:S_{(0)}^{i j k \ell m}}   (\partial_m \varrho_{(2)}^{(i,j)})\circ \Phi  \,,
\label{eq:inverse:div:R0}
\end{align}
and the error terms are given by
\begin{align}
G_{(1)}^{ijkm} =  -\partial_\ell(G_{(0)}^\ell A^k_i) \delta_{jm} - \partial_\ell G_{(0)}^k A_i^\ell \delta_{jm} + \partial_n (G_{(0)}^\ell A_i^k \partial_\ell \Phi^m ) A^n_j
\label{eq:inverse:div:G1}
\,,
\end{align}
where as before we denote $(\nabla \Phi)^{-1} = A$. We first show that the symmetric stress $R_{(0)}^{k\ell}$ defined in \eqref{eq:inverse:div:R0} satisfies the estimate \eqref{eq:inverse:div:stress:1}. First, we note that from \eqref{item:inverse:i} and \eqref{item:inverse:ii}, 
the function $\partial_m \varrho_{(2)}^{(i,j)}$ has zero mean, is $(\sfrac{\T}\mu)^3$ periodic, and satisfies  
\begin{align}
\norm{D^N \partial_m \varrho_{(2)}^{(i,j)} }_{L^p} \les \const_{*,p} \Upsilon^{-2} \Upsilon' \MM{N,1,\Upsilon',\Lambda} 
\label{eq:zero:order:stress:1}
\end{align} 
for $N\leq N_*-1$, in view of \eqref{eq:DN:Mikado:density}. Second, we note that since $D_t \Phi =0$, material derivatives may only land on the components of the $5$-tensor $S_{(0)}$. Third, the components of the $5$-tensor $S_{(0)}$ are sums of terms which are linear in $G_{(0)}$ and multilinear in $A$ and $D\Phi$. In particular, due to our assumption \eqref{eq:inverse:div:DN:G} and the previously established bounds in  \eqref{eq:derivatives:phase:gradient} and \eqref{eq:derivatives:phase:gradient:inverse}, upon applying the Leibniz rule, we obtain that 
\begin{align}
\norm{D^N D_t^M S_{(0)}}_{L^p} \les \const_{G,p} \max(\lambda,\lambda')^N  \MM{M,M_t,\nu, \nu'} 
\label{eq:zero:order:stress:2}
\end{align} 
for $N\leq N_*$ and $M\leq M_*$.
Having collected these estimates, the $L^p$ norm of the space-material derivatives of $R_{(0)}$ is obtained from Lemma~\ref{l:slow_fast}. As dictated by \eqref{eq:inverse:div:R0} we apply this lemma with $f = S_{(0)}$ and $\varphi = \partial_m \varrho_{(2)}^{(i,j)}$. Due to \eqref{eq:zero:order:stress:2}, the bound \eqref{eq:slow_fast_0} holds with $\const_f = \const_G$ and a spatial derivative cost of $\max(\lambda,\lambda')$. Due to \eqref{eq:DDpsi}, the assumptions \eqref{eq:slow_fast_1} and \eqref{eq:slow_fast_2} are verified. Next, due to \eqref{eq:zero:order:stress:1}, the assumption \eqref{eq:slow_fast_4} is verified, with $N_x=1$ and $\const_{\varphi}=\const_{*,p}\Upsilon^{-2}\Upsilon' \Lambda^\alpha$. Lastly, assumption \eqref{eq:inverse:div:parameters:1} verifies the condition \eqref{eq:slow_fast_3} of Lemma~\ref{l:slow_fast}. 
Thus, applying estimate \eqref{eq:slow_fast_5} we deduce that 
\begin{align}
\norm{D^N D_t^M R_{(0)}}_{L^p} \les \const_{G,p} \const_{*,p} \Upsilon^{-2} \Upsilon' \MM{N,1,\Upsilon',\Lambda}  \MM{M,M_t,\nu,\nu'}
\label{eq:zero:order:stress:2a}
\end{align}
for all $N\leq N_*-1$ and $M\leq M_*$, which is precisely the bound stated in \eqref{eq:inverse:div:stress:1}.   Here we have used that $N_* \geq 2 \Ndec + n + 1$, which gives that \eqref{eq:slow_fast_3_a} is satisfied.
\smallskip

\noindent\texttt{Step 3: } To continue the iteration, we first analyze the second term in \eqref{eq:inverse:div:G0}. The point is that this term has the {\em same structure} as what we started with; for every fixed $i,j,m$, we may replace $G_{(0)}^k$ by $G_{(1)}^{ijkm}$, and we replace $\varrho_{(0)}$ with $\partial_m \varrho_{(2)}^{(i,j)}$; the only difference is that the bounds for this term are \emph{better}. Indeed, from \eqref{eq:inverse:div:G1} we see that the $4$-tensor $G_{(1)}$ is the sum of various entries from the tensors $ DG_{(0)} \otimes A$ and $DG_{(0)} \otimes A \otimes A \otimes D\Phi$. Recalling \eqref{eq:derivatives:phase:gradient}, \eqref{eq:derivatives:phase:gradient:inverse}, and \eqref{eq:derivatives:mixed:G} and using the Leibniz rule, we deduce that
\begin{align}
 \norm{D^{N''} D_t^M D^{N'}G_{(1)}^{ijkm}}_{L^{{p}}} \les \const_{G,p} \max(\lambda,\lambda')^{N'+N''+1} \MM{M,M_t,\nu, \nu'}
 \label{eq:derivatives:mixed:G:1}
\end{align}
for $N' + N'' \leq N_*-1$ and $M\leq M_*$. The only caveat is that the bounds hold for one fewer spatial derivative. In order to iterate  Proposition~\ref{prop:Celtics:suck}, for simplicity we ignore the $i,j,k,m$ indices, since the argument works in exactly the same way in each case.  Specifically, we write $G_{(1)}^{ijkm}$ simply as $G_{(1)}^{k}$, and for the sake of convenience we suppress indices on the tensors $D\varrho_{(k)}$ and use $D$ as a stand-in for $\partial_m$.  We first note that $D \varrho_{(2)}= \div \div \left( D \varrho_{(4)} \right)$, where $D \varrho_{(4)}$ is a symmetric $2$-tensor once both indices have been specified on the left-hand side of the equality for $D\varrho_{(2)}$. Thus, using identities \eqref{eq:Celtics:suck:total}--\eqref{eq:Celtics:suck:error} and (in a slight abuse of notation) reusing the indices we previously tossed away, we obtain that the second term in \eqref{eq:inverse:div:G0} may be written as 
\begin{align}
G_{(1)}^k (D \varrho_{(2)}) \circ \Phi =  \partial_\ell R_{(1)}^{k\ell} +  G_{(2)}^{ijkm} (  \partial_m D \varrho_{(4)}^{(i,j)}) \circ \Phi
\label{eq:inverse:div:G0:new}
\end{align}
where the symmetric stress $R_{(1)}$ is given by
\begin{align}
R_{(1)}^{k\ell} =   \underbrace{\left(  G_{(1)}^k A^\ell_i \delta_{mj} + G_{(1)}^\ell A_i^k \delta_{mj} - G_{(1)}^n \partial_n \Phi^m A_i^k A_j^\ell \right)}_{=:S_{(1)}^{ijk\ell m} } (\partial_m D \varrho_{(4)}^{(i,j)}) \circ \Phi  \,,
\label{eq:inverse:div:R1}
\end{align}
the error terms are computed as
\begin{align}
G_{(2)}^{ijkm} = -\partial_\ell (G_{(1)}^\ell A_i^k) \delta_{jm} - \partial_\ell G_{(1)}^k A_i^\ell \delta_{jm} + \partial_n ( G_{(1)}^\ell A_i^k \partial_\ell \Phi^m ) A_j^n \, .
\label{eq:inverse:div:G2}
\end{align}
We emphasize that by combining \eqref{eq:inverse:div:R0} and \eqref{eq:inverse:div:G1} with \eqref{eq:inverse:div:R1} and \eqref{eq:inverse:div:G2}, we may compute the tensors $S_{(1)}$ and $G_{(2)}$ {\em explicitly in terms of just space derivatives} of $G$, $D\Phi$, and $A$. Using a similar argument to the one which was used to prove \eqref{eq:zero:order:stress:2}, but by appealing to \eqref{eq:derivatives:mixed:G:1} instead of \eqref{eq:derivatives:mixed:G}, we deduce that  for $N\leq N_*-1$ and $M\leq M_*$, 
\begin{align}
\norm{D^N D_t^M S_{(1)}}_{L^p} \les \const_{G,p} \max(\lambda,\lambda')^{N+1} \MM{M,M_t,\nu, \nu'} \,.
\label{eq:zero:order:stress:3}
\end{align} 
Using the bound \eqref{eq:zero:order:stress:3} and the estimate 
\begin{align*}
\norm{D^N ( \partial_m D \varrho_{(4)}) }_{L^p} \les \const_{*,p} \Upsilon^{-4} \Upsilon'^2 \MM{N,2,\Upsilon',\Lambda}
\,,
\end{align*} 
which is a consequence of \eqref{eq:DN:Mikado:density}, we may deduce from Lemma~\ref{l:slow_fast} that  
\begin{align}
\norm{D^N D_t^M R_{(1)}}_{L^p} \les \const_{G,p} \const_{*,p} \max(\lambda,\lambda') (\Upsilon^{-2}\Upsilon')^2 \MM{N,2,\Upsilon',\Lambda}  \MM{M,M_t,\nu,\nu'}
\label{eq:zero:order:stress:2b}
\end{align}
for $N\leq N_*-2$ and $M\leq M_*$, which is an estimate that is even better than \eqref{eq:zero:order:stress:2a}, aside from the fact that we have lost a spatial derivative. This shows that the first term in \eqref{eq:inverse:div:G0:new} satisfies the expected bound. The low-frequency portion of the second term in \eqref{eq:inverse:div:G0:new} may in turn be shown to satisfy
\begin{align}
 \norm{D^{N''} D_t^M D^{N'}G_{(2)}^{ijkm}}_{L^p} \les \const_{G,p} \max(\lambda,\lambda')^{2+N'+N''} \MM{M,M_t,\nu, \nu'}
 \label{eq:derivatives:mixed:G:2}
\end{align}
for $N'+N''\leq N_*-2$ and $M\leq M_*$.

At this point there is a clear roadmap for iterating this procedure $\lfloor \sfrac \dpot 2 \rfloor$ times, where the limit on the number of steps comes from that fact that $\varrho_{(k)}$ is only defined for $0\leq k \leq \dpot$, and each step in the iteration increases the value of $k$ by 2. Without spelling out these details, the iteration procedure described above produces
\begin{align}
G_{(0)} \varrho_{(0)}\circ \Phi =  \sum_{k=0}^{\lfloor \sfrac \dpot 2 \rfloor-1} \div R_{(k)} +  \underbrace{G_{(\lfloor \sfrac \dpot 2 \rfloor)} \, : \,  \left( D^{\lfloor \sfrac \dpot 2 \rfloor} \varrho_{\left(2\lfloor \sfrac \dpot 2 \rfloor\right)}\right) \circ \Phi}_{=:E}
\label{eq:inverse:div:final}
\end{align}
where each of the $\lfloor \sfrac \dpot 2 \rfloor$ symmetric stresses satisfies 
\begin{align}
\norm{D^N D_t^M R_{(k)}}_{L^p} \les \const_{G,p} \const_{*,p} \max(\lambda,\lambda')^{k} \left(\Upsilon^{-2} \Upsilon' \right)^{k+1} \Lambda^{N}  \MM{M,M_t,\nu,\nu'}
\label{eq:zero:order:stress:final}
\end{align}
for $N\leq N_*-k-1$ and $M\leq M_*$.  Furthermore, the formulae in \eqref{eq:divH:formula} and \eqref{eq:inverse:div:linear} can be computed explicitly from the algorithm already detailed above by keeping track of the high-low product structure of each term in each $R_{(k)}$ and Remark~\ref{rem:linear:inverse:div}, although we forego the details.  The subsidiary estimates are precisely those from \eqref{eq:zero:order:stress:1} and \eqref{eq:zero:order:stress:2}, which are immediate for the terms from the first step of the parametrix expansion, and which follow for the higher order terms by transferring the amplitude gains from the high-frequency function onto the low-frequency function, and using \eqref{eq:inverse:div:parameters:0}.  Each component of the the error tensor $G_{(\lfloor \sfrac \dpot 2 \rfloor)}$ in \eqref{eq:inverse:div:final} is recursively computable solely in terms of $G$, $D\Phi$, and $A$ and their spatial derivatives and satisfies 
\begin{align}
 \norm{D^{N''} D_t^M D^{N'}G_{(\lfloor \sfrac \dpot 2 \rfloor)}}_{L^p} \les \const_{G,p} \max(\lambda,\lambda')^{\lfloor \sfrac \dpot 2 \rfloor+N'+N''} \MM{M,M_t,\nu, \nu'}
 \label{eq:derivatives:mixed:ginal}
\end{align}
for $N'+N''\leq N_*-\lfloor \sfrac \dpot 2 \rfloor$ and $M\leq M_*$. Lastly, a final application of  Lemma~\ref{l:slow_fast}, which is valid due to with \eqref{eq:derivatives:mixed:ginal} and the assumption $N_* -\dpot \geq 2\Ndec + n + 1$, shows that estimate \eqref{eq:inverse:div:error:1} holds. \smallskip

\noindent\texttt{Step 4: } Finally, we turn to the proof of \eqref{eq:inverse:div:error:stress} and \eqref{eq:inverse:div:error:stress:bound}. Recall that $E$ is defined by the second term in \eqref{eq:inverse:div:final}, and thus $\fint_{\T^n} G \varrho \circ \Phi dx = \fint_{\T^n} E dx$.
Using the standard nonlocal inverse-divergence operator
\begin{align}
 (\RSZ f)^{ij} = -\frac{1}{2} \Delta^{-2} \partial_i \partial_j \partial_k f^k - \frac{1}{2} \Delta^{-1} \partial_k \delta_{ij} f^k + \Delta^{-1} \partial_i \delta_{jk} f^k + \Delta^{-1} \partial_j \delta_{ik} f^k
\label{eq:RSZ}
\end{align} 
we may define 
\begin{align*}
\RR_{\rm nonlocal} = \RSZ E  \,.
\end{align*}
By the definition of $\RSZ$ we have that $\RR_{\rm nonlocal}$ is traceless, symmetric, and satisfies $\div \RR_{\rm nonlocal} = E - \fint_{\T^n} E dx$, i.e. \eqref{eq:inverse:div:error:stress} holds. 

Using the formulas in \eqref{eq:RNC}, \eqref{eq:DNC}, the assumption \eqref{eq:inverse:div:v:global}, and the fact that $D$ and $\partial_t$ commute with $\RSZ$, we deduce that for every $N\leq N_\circ$ and $M\leq M_\circ$ we have
\begin{align}
\norm{D^N D_t^M \RR_{\rm nonlocal}}_{L^\infty} 
&\les \sum_{\substack{M'\leq M \\ N'+M' \leq N+M}} \sum_{K=0}^{M-M'}  \const_v^K  (\lambda')^{N-N'+K} \nu'^{-(M-M'-K)} \norm{D^{N'} \partial_t^{M'} \RSZ E}_{L^\infty}
\notag \\
&\les \sum_{\substack{M'\leq M \\ N'+M' \leq N+M}} (\lambda')^{N-N'} \nu'^{-(M-M')} \norm{D^{N'} \partial_t^{M'}  E}_{L^\infty}
\label{eq:riots:1}
\end{align}
where in the last inequality we have used that by assumption $\const_v \lambda' \les \nu'^{-1}$, and that  $\RSZ \colon L^p(\T^n) \to L^p(\T^n)$ is a bounded operator.

Our goal is to appeal to estimate~\eqref{eq:cooper:f:*} in Lemma~\ref{lem:cooper:1}, with $A = - v\cdot \nabla$, $B = D_t$ and $f=E$ in order to estimate the $L^\infty$ norm of $D^{N'} \partial_t^{M'}  E = D^{N'} (A+B)^{M'} E$. First, we claim that $v$ satisfies the lossy estimate
\begin{align}
\norm{D^ND_t^M v}_{L^{\infty}}
\les \const_v \lambda'^{N} \nu'^{-M}
\label{eq:DDv2} 
\end{align}
for $M\leq M_\circ$ and $N+M\leq N_\circ + M_\circ$. This estimate does not follow immediately from either \eqref{eq:DDv} or \eqref{eq:inverse:div:v:global}. For this purpose, we apply Lemma~\ref{lem:cooper:1} with $f = v$, $B= \partial_t$, $A = v\cdot \nabla$, and $p=\infty$.  Using \eqref{eq:inverse:div:v:global}, and the fact that $B = \partial_t$ and $D$ commute, we obtain that bounds \eqref{eq:cooper:v} and \eqref{eq:cooper:f} hold with $\const_f = \const_v$, $\lambda_v = \tilde \lambda_v = \lambda_f = \tilde \lambda_f = \lambda'$, and $\mu_v = \tilde \mu_v = \mu_f = \tilde \mu_f = \nu'^{-1}$. Since $A+ B = D_t$, we obtain from the bound \eqref{eq:cooper:f:*} and the assumption $\const_v \lambda' \les \nu'^{-1}$ that \eqref{eq:DDv2} holds.

Second, we claim that for any $k\geq 1$ we have
\begin{align}
\norm{ \left(\prod_{i=1}^k D^{\alpha_i} D_{t}^{\beta_i}\right) v}_{L^\infty(\supp G)} \les \const_v \lambda'^{|\alpha|} \nu'^{|\beta|}
\label{eq:riots:5}
\end{align}
whenever $|\beta|\leq M_\circ$ and $|\alpha| +|\beta| \leq N_\circ + M_\circ$.
To see this, we use Lemma~\ref{lem:cooper:2} with $f =v$, $p=\infty$, and $\Omega = \supp G$. From \eqref{eq:DDv} we have that \eqref{eq:cooper:2:v} holds with  $\const_v = \nu/ \lambda' $, $\lambda_v = \tilde \lambda_v = \lambda'$, $\mu_v=\nu$, and $\tilde \mu_v = \nu'$. On the other hand, from \eqref{eq:DDv2} we have that \eqref{eq:cooper:2:f} holds with $\const_f = \const_v$, $\lambda_f = \tilde \lambda_f = \lambda'$, and $\mu_f =\tilde \mu_f = \nu'^{-1}$. We then deduce from \eqref{eq:cooper:2:f:2} that \eqref{eq:riots:5} holds. 

Third, we claim that 
\begin{align}
\norm{ \left(\prod_{i=1}^k D^{\gamma_i} D_{t}^{\beta_i}\right) E}_{L^\infty(\supp G)} \les \const_{G,p} \const_{*,p} \max(\lambda,\lambda')^{\lfloor \sfrac \dpot 2 \rfloor} (\Upsilon' \Upsilon^{-2})^{\lfloor \sfrac \dpot 2 \rfloor} \Lambda^{|\gamma|+n+1} \MM{|\beta|,M_t,\nu, \nu'}
\label{eq:riots:6}
\end{align}
holds  whenever $|\gamma|\leq N_* - \lfloor \sfrac \dpot 2 \rfloor - n - 1$ and $|\beta|\leq M_*$. This estimate again follows from Lemma~\ref{lem:cooper:2}, this time with $f= E$, by appealing to the previously established bound \eqref{eq:inverse:div:error:1} and the Sobolev embedding $W^{n+1,1}(\T^n) \hookrightarrow L^\infty(\T^n) $.

At last, we are in the position to apply Lemma~\ref{lem:cooper:1}. The bound 
\eqref{eq:riots:5} implies that assumption \eqref{eq:cooper:v} holds with $B= D_t$, $\lambda_v = \tilde \lambda_v = \lambda'$, and $\mu_v = \tilde \mu_v = \nu'$. The bound
\eqref{eq:riots:6} implies that assumption \eqref{eq:cooper:f} of Lemma~\ref{lem:cooper:1} holds with $\const_f = \const_{G,p} \const_{*,p} \max(\lambda,\lambda')^{\lfloor \sfrac \dpot 2 \rfloor} (\Upsilon' \Upsilon^{-2})^{\lfloor \sfrac \dpot 2 \rfloor} \Lambda^{n+1}$, $\lambda_f = \tilde \lambda_f = \Lambda$, $\mu_f=\nu$, and $\tilde \mu_f = \nu'$. We may now use estimate \eqref{eq:cooper:f:*}, and the assumption that $\Lambda \geq \lambda,\lambda'$ to deduce that 
\begin{align}
 \norm{D^{N'} \partial_t^{M'}  E}_{L^\infty} \les \const_{G,p} \const_{*,p} \max(\lambda,\lambda')^{\lfloor \sfrac \dpot 2 \rfloor} (\Upsilon' \Upsilon^{-2})^{\lfloor \sfrac \dpot 2 \rfloor} \Lambda^{N'+n+1} (\max\{\const_v \Lambda, \nu' \})^{M'}
 \label{eq:riots:2}
\end{align}
holds whenever $M' \leq M_\circ$ and $N'+M' \leq N_\circ + M_\circ$.
Combining \eqref{eq:riots:1} and \eqref{eq:riots:2} we deduce that 
\begin{align}
\norm{D^N D_t^M \RR_{\rm nonlocal}}_{L^\infty} 
&\les \const_{G,p} \const_{*,p} \max(\lambda,\lambda')^{\lfloor \sfrac \dpot 2 \rfloor} (\Upsilon' \Upsilon^{-2})^{\lfloor \sfrac \dpot 2 \rfloor} \Lambda^{n+1} \notag\\
&\qquad \qquad \times 
\sum_{\substack{M'\leq M \\ N'+M' \leq N+M}} \lambda'^{N-N'} \nu'^{-(M-M')} \Lambda^{N'} (\max\{\const_v \Lambda, \nu' \})^{M'}
\notag\\
&\les \const_{G,p} \const_{*,p} \max(\lambda,\lambda')^{\lfloor \sfrac \dpot 2 \rfloor} (\Upsilon' \Upsilon^{-2})^{\lfloor \sfrac \dpot 2 \rfloor} \Lambda^{N+n+1} (\max\{\const_v \Lambda, \nu' \})^{M}
\label{eq:riots:3}
\end{align}
whenever $N\leq N_\circ$ and $M\leq M_\circ$.  Estimate \eqref{eq:inverse:div:error:stress:bound} follows by appealing to the assumption \eqref{eq:riots:4}.
\end{proof}

Observe that in the proof of Proposition \ref{prop:intermittent:inverse:div}, $\rho^{\be{(j)}}$ consists of $\na \varrho_{(2)}$, $\na^2 \varrho_{(4)}, \cdots, \na^{\floor{\sfrac{\dpot}2}} \varrho_{2\floor{\sfrac{\dpot}{2}}}$; recall that $\varrho_{(0)}= \varrho = \div^\dpot \vartheta$ and $\varrho_{(k-1)} = \div \varrho_{(k)} = \div^{\dpot-(k-1)} \vartheta$. Keeping this in mind, when $\varrho$ is given as $\div^{(2\dpot)^2} \vartheta$, we can apply the proposition iteratively to get 
\begin{align*}
    G(\rho\circ\Phi) =\div^{\dpot} R + E.
\end{align*}
The details are described in the following corollary. Since this operator will be applied to velocity increments, some of the adjustments are specified for this particular application.

\begin{corollary}[\bf Iterated inverse divergence for scalar fields]\label{cor:inverse.div}
We suppose that the same assumptions hold as in Proposition \ref{prop:intermittent:inverse:div} together with Remark~\ref{rem:scalar:inverse:div} except for the following substitutions.
\begin{enumerate}[(i)]
    \item  Fix $\Ndec,N_*,M_*\dpot \geq 1$ such that $\dpot$ is even and $N_* - {\dpot^2} \geq 2\Ndec +  n + 1+M_*$ (replacing \eqref{eq:inv:div:NM} and the last inequality in \eqref{eq:inverse:div:parameters:0}).
    \item $\varrho$ is given as an iterated divergence $\varrho = \div^{{(\dpot^2)}} \td\vartheta$ (replacing \eqref{item:inverse:i}).
    \item \label{item:inverse.divd:ii} There exist parameters $1 \ll \Upsilon \leq \Upsilon' = \Lambda$ and $\const_{*,p}>0$ such that for all $0\leq N \leq {N_*}$ and all $0\leq k \leq {\dpot^2}$, \eqref{eq:DN:Mikado:density} is replaced with
\begin{align}
\norm{D^N \partial_{i_1}\dots \partial_{i_k} \td\vartheta^{(i_1,\dots, i_{\dpot^2})}}_{L^p} \les \const_{*,p} \Upsilon^{k-\dpot^2} \Upsilon'^N \, .
\label{eq:DN:Mikado:density.divd}
\end{align} 
\end{enumerate}
Additionally, we assume that there exists a smooth, non-negative function $\pi$ such that
\begin{align}\label{eq:inv:div:extra:pointwise:iterated}
\left| D^N D_{t}^M G \right| &\lesssim \pi^\frac12 r_G^{-\frac13} \lambda^N\MM{M,M_{t},\nu,\nu'} \, 
\end{align}
for $N\leq N_*$ and $M\leq M_*$. Then, we have that
\begin{align}\label{divd.expression}
G(\varrho\circ \Phi)
= \div^{\dpot} R + E
\end{align}
for a rank $dpot$ tensor $R$ and error $E$ satisfying the following properties.
\begin{enumerate}[(i)]
    \item\label{item:divd:local:i} The support of $R$ is a subset of $\supp G \cap \supp (\td\vartheta\circ\Phi)$, and hence so is the support of $E$.
    \item \label{item:divd:local:ii} There exists an explicitly computable positive integer $\ov\const_H$, an explicitly computable function $r(j):\{0,1,\dots,\ov\const_H\}$ and explicitly computable tensors
    \begin{align*}
        &\rho^{\beta(j)} \, , \qquad \beta(j) = (\beta_1,\beta_2,\dots,\beta_{r(j)}) \in \{1,\dots,n\}^{r(j)} \, , \\
        &H^{\alpha(j)} \, , \qquad \alpha(j) = (\alpha_1,\alpha_2,\dots,\alpha_{r(j)}) \in \{1,\dots,n\}^{r(j)+\dpot} \, ,
    \end{align*}
    of rank $r(j)$ and $r(j)+\dpot$, respectively, all of which depend only on $G,\varrho,\Phi,n,\dpot$ such that the following holds.  The localized stress
    $R$ can be decomposed into a sum of localized stresses as
    \begin{align*}
    R &= \sum_{j=0}^{\ov\const_\divH}
H^{\alpha{(j)}}
(\rho^{\beta{(j)}}  \circ \Phi) \, . 
\end{align*}
Furthermore, we have that
\begin{align}
\supp H^{\alpha{(j)}} \subseteq \supp G  \, , \qquad \supp \rho^{\beta{(j)}} \subseteq \supp \td\vartheta \,  . \label{eq:inverse:divd:linear}
\end{align}
\item\label{item:divd:local:iii} We have the subsidiary estimates
\begin{subequations}\label{eq:inverse:divd:sub:main}
\begin{align}
    \left\| D^N \rho^{\beta{(j)}} \right\|_{L^p} &\les \const_{*,p} (\Upsilon^{-2} \Upsilon' )^{{\dpot}} {\Lambda}^{N} \label{eq:inverse:divd:sub:1}
    \end{align}
for all $N \leq N_*-\dpot^2$ and $j\leq \ov\const_\divH$, and
    \begin{align}
    \left\|\prod_{i=1}^{k}D^{\alpha_i} D_{t}^{\be_i} H^{\alpha{(j)}}\right\|_{L^p} &\lesssim \const_{G,p} \left(\max(\lambda,\lambda')\right)^{|\al|} \MM{|\be|,M_{t},\nu,\nu'} \,  \label{eq:inverse:divd:sub:2} \\
      \left|\prod_{i=1}^{k}D^{\alpha_i} D_{t}^{\be_i} H^{\alpha{(j)}}\right|
    &\lec \pi^{\frac12}r^{-\frac13} (\max(\la, \la'))^{|\al|}\MM{|\be|, M_t, \nu, \td\nu}. \label{eq:invser:divd:point}
\end{align}
\end{subequations}
for all integer $k\geq 1$, multi-indices $\al, \be\in \mathbb{N}^k$ with $|\al|\leq N_*-\dpot^2$ and $|\be|\leq M_*$, and $j\leq\ov\const_{\divH}$.

\item\label{item:divd:local:iv} We have the main estimate
\begin{align}
\norm{\prod_{i=1}^{k}D^{\alpha_i} D_{t}^{\be_i} R}_{L^p}
 &\les  \const_{G,p} \const_{*,p}  (\Upsilon' \Upsilon^{-2} )^{{\dpot}} {\Upsilon'}^{|\al|} \MM{|\be|,M_{t},\nu,\nu'} 
\label{eq:inverse:divd:stress:1}
\end{align}
for all integer $k\geq 1$, multi-indices $\al, \be\in \mathbb{N}^k$ with $|\al|\leq N_*-\dpot^2$ and $|\be|\leq M_*$, and $j\leq\ov\const_{\divH}$.
\item\label{item:divd:nonlocal} For $N \leq N_* - \dpot^2$ and $M\leq M_*$ the error term $E$ in \eqref{divd.expression} satisfies\footnote{In our applications, $\Upsilon=\Upsilon'$, so the sum of loss factors is irrelevant.  If one wanted to be more precise, this loss could be eliminated using a more careful algorithm and a few more conditions on the relative sizes of all the frequencies.}
\begin{align}
\norm{D^N D_{t}^M E}_{L^p}  
\les \const_{G,p} \const_{*,p}   \max(\lambda,\lambda')^{{\sfrac \dpot 2}} \left( \Upsilon' \Upsilon^{-2}\right)^{{\sfrac \dpot 2}} {\Lambda}^{N} \MM{M,M_{t},\nu,\nu'} {\sum_{k=0}^{\dpot-1} \left(\frac{\Upsilon'}{\Upsilon}\right)^{2k}}
\,.
\label{eq:inverse:divd:error:1}
\end{align}
\end{enumerate}

\end{corollary}
\begin{proof}
The proof is based on applying Proposition \ref{prop:intermittent:inverse:div} $\dpot$ times.  In the first iteration, we get
\begin{align*}
    G(\varrho\circ \Phi)
    = \sum_{j_1=0}^{\const_{\divH}}
    \div \left(H^{\alpha{(j_1)}}(\rho^{\beta{(j_1)}}\circ\Phi) \right) + E_{(1)}
\end{align*}
where $H^{\alpha{(j_1)}}$ satisfies \eqref{eq:inverse:div:sub:2} and \eqref{eq:inverse:div:v:global}. From \eqref{eq:divH:formula} and Remark~\ref{rem:scalar:inverse:div}, we have that the rank of $H^{\alpha(j_1)}$ is one larger than the rank of $\rho^{\beta(j_1)}$. Also, replacing $\pi$ by $\pi^{\sfrac12} r^{-\sfrac13}$ in Remark \ref{rem:pointwise:inverse:div}, we get 
\begin{align*}
    |D^N D_t^M H^{\alpha{(j_1)}}|
    \lec \pi^\frac12 r^{-\frac13}\la^N \MM{M, M_t, \nu, \td \nu} 
\end{align*}
for $N\leq N_*-{\halfd}$ and $M\leq M_*$. In addition, $E_{(1)}$ satisfies \eqref{eq:inverse:divd:error:1}. Since we use the same $\Phi$, all assumptions on $G$ and $\Phi$ in the proposition holds for $N_*$ replaced with $N_*-\halfd$. From the proof of Proposition \ref{prop:intermittent:inverse:div} we note that $\rho^{\be{(j)}}$ consists of $\na^k \varrho_{(2k)}$, $1\leq k \leq \halfd$, which can be written as $\na^k \div^{\dpot^2 - 2k}\td\vartheta=\div^{\dpot} (\na^k \div^{\dpot^2 - 2k-\dpot}\td\vartheta)$.  Then, $\na^k \varrho_{(2k)}$ and its potential $\na^k \div^{\dpot^2 - 2k-\dpot}\td\vartheta$ satisfy \eqref{item:inverse:i}, \eqref{item:inverse:ii} in the assumption of Proposition \ref{prop:intermittent:inverse:div} and 
\begin{align*}
    \norm{D^N \pa_{i_1}\cdots\pa_{i_{k'}} (\na^k \div^{\dpot^2 - 2k-\dpot}\td\vartheta)}
    \les \const_{*,p} \Upsilon^{-2k-\dpot+k'}\Upsilon'^{N+k}
\end{align*}
for any $N\leq {N_*-k}$ and $0\leq k'\leq \dpot$. In particular, we have
\begin{align}\label{est.rhoj1}
     \left\| D^N \rho^{\beta{(j_1)}} \right\|_{L^p} &\les \const_{*,p} \Upsilon^{-2} \Upsilon' {\Upsilon'}^{N}
\end{align}
for $N \leq N_* - \halfd$ and $j_1 \leq \const_\divH$. This implies that \eqref{eq:DN:Mikado:density} holds for $C_{*,p}$ replaced with $C_{*,p} \Upsilon'\Upsilon^{-2}$ and $N_*$ with $N_*-\halfd$ and $\vartheta$ with the potential of $\rho^{\be{(j)}}$, respectively. Furthermore, from the construction it is easy to see that 
\begin{align*}
    \supp\left(\rho^{\be{(j)}}\right)\subset \supp(\td\vartheta) \, . 
\end{align*}

Iterating this process $\dpot$ times, we get 
\begin{align*}
    G(\varrho\circ \Phi)
    &= \sum_{j_1=0}^{\const_\divH}
    \div \left(H^{\alpha{(j_1)}}(\rho^{\beta{(j_1)}}\circ\Phi) \right) + E_{(1)}
    =\sum_{j_1, j_2=0}^{\const_\divH}
    \div^{2} \left(H^{\alpha{(j_1,j_2)}}(\rho^{\beta{(j_1,j_2)}}\circ\Phi) \right) 
    + \div E_{(2)}
    + E_{(1)}\\
    &=:\sum_{j=0}^{\ov\const_\divH}
\div^{\dpot}\left(H^{\alpha{(j)}}
(\rho^{\beta{(j)}}  \circ \Phi)\right)
+\sum_{k=1}^{\dpot} 
    \div^{k-1} E_{(k)}\, . 
\end{align*}
As a result, we get \eqref{divd.expression}, where $E$ is defined by
\begin{align*}
     E &:= \sum_{k=1}^{\dpot} 
    \div^{k-1} E_{(k)} \, .
\end{align*}
Since we have
\begin{align*}
    \supp H^{\alpha{(j)}}\subset 
\cdots \subset \supp(H^{\al{(j_1)}}) \subset \supp(G), \quad
\supp \rho^{\beta{(j_1)}}\subset \supp(\td\vartheta) \, ,
\end{align*}
\eqref{eq:inverse:divd:linear} holds. Therefore, \eqref{item:divd:local:i} and \eqref{item:divd:local:ii} have been verified, as has \eqref{est.rhoj1} and \eqref{eq:inverse:divd:sub:1}. Furthermore, we have
    \begin{align*}
    \left\|D^N D_{t}^M H^{\alpha{(j)}}\right\|_{L^p} &\lesssim \const_{G,p} \left(\max(\lambda,\lambda')\right)^{N} \MM{M,M_{t},\nu,\nu'} \,  \\
      \left|D^N D_{t}^M H^{\alpha{(j)}}\right|
    &\lec \pi^{\frac12}r^{-\frac13} (\max(\la, \la'))^{N}\MM{M, M_t, \nu, \td\nu}.\\
    \norm{D^N D_{t}^M  R}_{L^p}
 &\les  \const_{G,p} \const_{*,p}  (\Upsilon' \Upsilon^{-2} )^\dpot {\Lambda}^{N} \MM{M,M_{t},\nu,\nu'} 
\end{align*}
for all integers $N\leq N_*-\dpot^2$ and $M\leq M_*$. Also, $E^{(k)}$ satisfies
\begin{align*}
 \norm{D^N D_{t}^M E^{(k)}}_{L^p}  
\les \const_{G,p} \const_{*,p} (\Upsilon'\Upsilon^{-2})^{k-1}   \max(\lambda,\lambda')^{\halfd} \left( \Upsilon' \Upsilon^{-2}\right)^{\halfd} \Upsilon'^{N} \MM{M,M_{t},\nu,\nu'} 
\end{align*}
for $1\leq k\leq\dpot$, $N\leq N_*-k\cdot \halfd$, and $M\leq M_*$. 

Finally, we apply Lemma \ref{lem:cooper:2} to upgrade these estimates to the one with commutations of the operators, \eqref{eq:inverse:divd:sub:2}, \eqref{eq:invser:divd:point}, \eqref{eq:inverse:divd:stress:1}, and \eqref{eq:inverse:divd:error:1}. We will work only for \eqref{eq:inverse:divd:sub:2}, then the last will follow by a similar argument. 
To avoid confusion in the notations, we rewrite some repeated symbols from Lemma~\ref{lem:cooper:2} with bars above on the left-hand side of the equalities below, while the right-hand side are parameters given in the assumptions of the Corollary. Set $\overline p=p$, $\overline N_t=M_t$, $\overline N_*=N_*-\dpot\halfd$, $\overline M_*=M_*$, $\overline v=v$, $\overline \Omega= \supp{G}$, $\const_v = \nu (\la')^{-1}$, $\la_v= \td\la_v = \la'$, $\mu_v = \mu_f = \nu$, $\td\mu_v = \td \mu_ f= \td \nu$, $ f = H^{\al(j)}$, and $\la_f = \td \la_f = \max(\la, \la')$. Then, as a consequence of the lemma, we have $\eqref{eq:inverse:divd:sub:2}$. For \eqref{eq:invser:divd:point}, we work at each point $x$ in a similar way, but set $\overline\Omega = \Omega(x)$ as a small closed neighborhood of $x$ contained in $\supp(G)$ and use the continuity of $\pi$ so that $\sup_{\Omega(x)}\pi \leq 2\pi(x)$. 

\end{proof}

Finally, we shall need a simpler case of the inverse divergence, when the density is not flowed and the input is a scalar field. 

\begin{lemma}[\bf Inverse divergence without flow map]\label{rem:no:decoup:inverse:div2}
Fix dimension $n\geq 2$. Let $G$ be a smooth scalar field and let $\dpot$ be a non-negative integer such that the smooth scalar field $\varrho$ and tensor field $\vartheta$ defined on $\R\times \T^n$ satisfy $\varrho = \partial_{i_1}\dots\partial_{i_\dpot}\vartheta^{(i_1\dots i_\dpot)}(x)$ (note that no symmetry assumptions needed). 
\smallskip

\noindent\textbf{Part 1: Algorithm for inverse divergence}\\
We have a decomposition
\begin{align}\label{decomp:noflow}
    G \varrho =: \div (\divH(G\varrho)) + E
\end{align}
where the vector field $\divH(G\varrho)$ and scalar field $E$ are defined by
\begin{align}\label{defn.divHR:noflow}
    \divH(G\varrho)^\bullet
    &:= \sum_{k=0}^{\dpot-1}
    (-1)^{\dpot-k+1}
     \pa_{i_{k+2}} \dots\pa_{i_{\dpot}} G \, \underbrace{\div^{(k)}}_{{\pa_{i_1},\dots,\pa_{i_k}}}\vartheta^{(i_1,\dots,i_k,\bullet, i_{k+2},\dots, i_\dpot)}, \quad
    E = (-1)^\dpot \na^{\dpot}G: \vartheta\, ,
\end{align}
where we use the convention $\pa_{i_{k+2}} \cdots\pa_{i_{\dpot}} G=G$ and $\vartheta^{(i_1,\dots,i_k,\bullet, i_{k+2},\dots, i_\dpot)}=\vartheta^{(i_1,\dots,i_{\dpot-1},\bullet)}
$ when $k=\dpot-1$.
\smallskip

\noindent\textbf{Part 2: Localized assumptions and output}\\
Fix a set $\Omega\subset \R\times \T^n$.  Let parameters $N_*\geq M_*\geq 1$ be given. Define $v$ and $D_t$ as in Part 1 of Proposition~\ref{prop:intermittent:inverse:div}, where $v$ satisfies \eqref{eq:DDv} with $\lambda',\nu,\nu',N_*,M_*$ and $L^\infty(\supp G)$ replaced with $L^\infty(\Omega)$. Let smooth, non-negative functions $\pi$ and $\pi'$ be given such that
\begin{subequations}
\begin{align}
\left| D^N D_{t}^M G \right| &\lesssim \pi \lambda^N\MM{M,M_{t},\nu,\nu'} \qquad \textnormal{on }\Omega \, 
\label{eq:inv:div:extra:pointwise:noflow}\\
\Upsilon^{\dpot-k} \left| D^N D_{t}^M \partial_{i_1}\dots \partial_{i_k} \vartheta^{(i_1,\dots, i_\dpot)}\right| &\lesssim \pi' \Lambda^N\MM{M,M_{t},\nu,\nu'} \qquad \textnormal{on }\Omega
\label{eq:inv:div:extra:pointwise2:noflow}
\end{align}
\end{subequations}
for $N\leq N_*$ and $M\leq M_*$, where the parameters satisfy
\begin{align}\label{parameter:noflow}
    \la', \la \leq \Upsilon\leq \La, \quad
    \max(\la, \la') \Upsilon^{-1}\leq 1, \quad
    N_* \geq \dpot, \quad \la, \nu, \nu' \geq 1 \, .
\end{align}
Then $\divH(G\varrho)$ satisfies 
\begin{equation}\label{eq:div:no:flow:support}
\supp(\divH(G\varrho))\subseteq \supp(G\vartheta) \, ,
\end{equation}
and for $N\leq N_*-\dpot$ and $M\leq M_*$,
\begin{align}
\left| D^N D_{t}^M \divH(G\varrho) \right| &\lesssim 
\pi \pi' \Upsilon^{-1}\Lambda^N\MM{M,M_{t},\nu,\nu'} \,  \quad \textnormal{on }\Omega\, .
\label{eq:inv:div:pointwise:local}
\end{align}
\smallskip

\noindent\textbf{Part 3: Nonlocal assumptions and output}\\
Finally, we assume that all assumptions from \eqref{item:nonlocal:v} in Part 4 in Proposition~\ref{prop:intermittent:inverse:div} hold.  Next, we assume that for $N\leq N_*$ and $M\leq M_*$,
\begin{subequations}
\begin{align}
\norm{D^N D_{t}^M G}_{L^{\infty}}&\lesssim \const_{G,\infty} \lambda^N (\nu')^M \, ,
\label{eq:inverse:div:DN:G:noflow}\\
\norm{D^ND_t^M \partial_{i_1}\dots \partial_{i_k} \vartheta^{(i_1,\dots, i_\dpot)}}_{L^\infty} &\les \const_{*,\infty} \Upsilon^{k-\dpot} \Lambda^N(\nu')^M \, .
\label{eq:DN:Mikado:density:noflow}
\end{align}
\end{subequations}
Also, we choose $\dpot$ large enough to satisfy
\begin{align}\label{dpot:noflow}
    \const_{G,\infty} \const_{*,\infty} (\max(\lambda,\lambda')\Upsilon^{-1})^{{\halfd}}   \Lambda^{K_\circ} \left(1 + \frac{\max\{ \nu', \const_v \Lambda \}}{\nu
}\right)^{M_\circ}\
\leq 1
\, .
\end{align}
Then we may write  
\begin{align}
E =: \div \left(\divR(G \varrho )\right) + \fint_{\T^3} G \varrho \, dx \, ,
\label{eq:inverse:div:error:stress:no:flow}
\end{align}
where $\divR(G \varrho )$ is a vector field which satisfies
\begin{align}
\norm{D^N D_{t}^M \divR(G \varrho ) }_{L^\infty}  
\lec  \frac{1}{\Lambda^{K_\circ}} (\max(\lambda,\lambda')\Upsilon^{-1})^{{\halfd}}\Lambda^N \nu^M
\label{eq:inverse:div:error:stress:bound:no:flow}
\end{align}
for  $N \leq N_\circ$ and $M\leq M_\circ$.
\end{lemma}
\begin{proof}[Proof of Lemma~\ref{rem:no:decoup:inverse:div2}]
With the definition \eqref{defn.divHR:noflow} in hand, we can easily check \eqref{decomp:noflow}--\eqref{eq:inv:div:pointwise:local}. To define $\divR(G \varrho )$, we use the standard operator $(\mathcal{R}f)^i=\Delta^{-1}\partial_i$ and let $\divR(G \varrho ) = \RSZ E$. The desired estimate for $\divR(G \varrho)$ follows as in the Proof of Proposition~\ref{prop:intermittent:inverse:div} with minor modifications, and we leave the details to the reader.
\end{proof}

\subsection{Upgrading material derivatives}\label{appetizer:4}
\begin{lemma}[\bf Upgrading material derivatives]\label{lem:upgrading.material.derivative}
Fix $p\in [1, \infty]$ and a positive integer $N_\star\leq \sfrac{3\Nfin}4$. Assume that a tensor $F$ is given with a decomposition $F = F^l + F^*$ which satisfy
\begin{subequations}
\begin{align} 
\left\| \psi_{i,q} D^N D_{t,q}^M F^l \right\|_{p}
&\lesssim 
\const_{p,F}  \lambda_{F}^N  \MM{M,\Nindt,\Gamma_{q}^{i+c} \tau_{q}^{-1}, \Gamma_{q}^{-1}\Tau_{q}^{-1} }
\label{eq:before.upgraded}\\
\left\| D^N D_{t,q}^M F^* \right\|_{\infty} 
&\lesssim \const_{*,F}\Tau_{q+\bn}^{\Nindt}\lambda_{F}^N\tau_q^{-M} 
\label{eq:before.upgraded.uniform}
\end{align}
\end{subequations}
for all $M+N\leq N_{\star}$, an absolute constant {$c\leq 20$}, and constants $\const_{p,F}$ and $\const_{*,F}$. Assume furthermore that there exists $k$ such that $q+1<k\leq q+\bn$ and
\begin{align}\label{supp.F}
    \supp(\hat w_{q'}, \la_{q'}^{-1}{\Ga_{q'}}) \cap \supp(F^l) =\emptyset \quad \forall q+1\leq q' < k \, .
\end{align}
Finally, assume that
\begin{equation}\label{eq:timescale:upgrading}
\la_F \Ga_{q+\bn}^{\imax+2}\de_{q+\bn}^{\frac12}r_q^{-\frac13}\leq \Tau_{q+\bn}^{-1} \, .
\end{equation}
Then $F$ obeys the following estimate with an upgraded material derivative for all $M+N\leq N_{\star}$;
\begin{align}
&\left\| \psi_{i,k-1} D^N D_{t,k-1}^M F \right\|_{p} \lesssim \
(\const_{p,F} + \const_{*,F}) \max(\lambda_{F}, \La_{k-1})^N
  \MM{M,\Nindt,\Gamma_{k-1}^{i} \tau_{k-1}^{-1}, \Gamma_{k-1}^{-1}\Tau_{k-1}^{-1} } \, .
\label{eq:after.upgraded}
\end{align}
In particular, the nonlocal part $F^*$ obeys better estimate
\begin{align}\label{est:nonlocal:upgrade}
\norm{D^N D_{t,k-1}^M F^*}_\infty
    \lec \const_{*,F}\max(\lambda_{F}, \la_{k-1}\Ga_{k-1})^N 
    \MM{M, \Nindt, \tau_{k-1}^{-1},  \Tau_{k-1}^{-1}\Gamma_{k-1}^{-1}}
\end{align}
for $N+M\leq N_\star$.

Similarly, if instead of \eqref{eq:before.upgraded}, $F^l$ satisfies 
\begin{align}
    \left|\psi_{i,q} D^N D_{t,q}^M F^l \right| 
&\lesssim 
\pi_{F}  \lambda_{F}^N  \MM{M,\Nindt,\Gamma_{q}^{i+c} \tau_{q}^{-1}, \Gamma_{q}^{-1}\Tau_{q}^{-1} }
\label{eq:before.upgraded.pt}
\end{align}
for all $M+N\leq N_{\star}$, an absolute constant ${c\leq 24}$, and a positive function $\pi_F$ with $\pi_F \geq \const_{*,F}$, we have
\begin{align}
&\left| \psi_{i,k-1} D^N D_{t,k-1}^M F \right| \lesssim \
\pi_{F}\max(\lambda_{F}, \La_{k-1})^N
  \MM{M,\Nindt,\Gamma_{k-1}^{i} \tau_{k-1}^{-1}, \Gamma_{k-1}^{-1}\Tau_{k-1}^{-1} }
\label{eq:after.upgraded.pt}
\end{align}
for all $M+N\leq N_{\star}$, provided that \eqref{eq:timescale:upgrading} holds.
\end{lemma}

\begin{proof} 
We first handle the local portion $F^l$ by upgrading $\psi_{i,q}$ in \eqref{eq:before.upgraded} to the one with $\psi_{i,k-1}$, and then upgrading $D_{t,q}$ to $D_{t,k-1}$. Since $\psi_{i',q}^6$ forms a partition of unity from \eqref{eq:inductive:partition} and we have $\tau_q^{-1}\Ga_q^{i'+{24}} \leq \tau_{k-1}^{-1}\Ga_{k-1}^i$ when $\psi_{i',q}\psi_{i,k-1}\neq 0$ by \eqref{eq:inductive:timescales}, we obtain that
\begin{align}
    \left\| \psi_{i,k-1} D^N D_{t, q}^M F^l \right\|_{p} 
    &=
    \left\| \psi_{i,k-1}\sum_{i'=0}^{\imax} \psi_{i',q}^6 D^N D_{t, q}^M F^l \right\|_{p} \notag\\
    &\lec \sum_{i':\psi_{i',q}\psi_{i,k-1}\neq 0 }\left\| \psi_{i',q} D^N D_{t, q}^M F^l\right\|_{p}\notag\\
    &\lec \const_{p,F} \lambda_{F}^N \MM{M, \Nindt, \tau_{k-1}^{-1}\Gamma_{k-1}^{i}, \Gamma_{k-1}^{-1}\Tau_{k-1}^{-1}}. \label{eq:upgrade.cutoff}
\end{align}
Here we used the maximal cardinality of $i'$ is $\imax$. Then, using \eqref{supp.F}, we have $D_{t,k-1}^M F^l =D_{t,q}^M F^l$ and the desired inequality \eqref{eq:after.upgraded} for $F^l$ follows. In a similar way, we can also get \eqref{eq:after.upgraded.pt} for $F^l$.

On the other hand, we handle the nonlocal portion $F^*$ by claiming that for each $q\leq k' \leq k-1$, we have 
\begin{align}\label{ind.mat}
\norm{D^N D_{t,k'}^M F^*}_\infty
    \lec \const_{*,F}\Tau_{q+\bn}^{\Nindt} \max(\lambda_{F}, \la_{k'}\Ga_{k'})^N 
  (\Tau_{k-1}^{-1}\Gamma_{k-1}^{-1})^{M},
\end{align}
for all $N+ M\leq N_\star$. In particular, this implies that
\begin{align*}
\norm{D^N D_{t,k-1}^M F^*}_\infty
    \lec \const_{*,F}\max(\lambda_{F}, \la_{k-1}\Ga_{k-1})^N 
    \MM{M, \Nindt, \tau_{k-1}^{-1},  \Tau_{k-1}^{-1}\Gamma_{k-1}^{-1}}
\end{align*}
for $N+M\leq N_\star$, which yields \eqref{eq:after.upgraded} and \eqref{eq:after.upgraded.pt}. The proof of the claim is then given by an inductive argument on $k'$. When $k'=q$, it easily follows from \eqref{eq:before.upgraded.uniform}. Next, suppose that \eqref{ind.mat} holds for some $k'<k-1$, and we apply Remark \ref{rem:upgrade.material.derivative.end} to $v= \hat u_{k'}$, $w= \hat w_{k'+1}$, $f=F$, $\Omega = \T^3$, $N_*=N_\star$, $N_t = \Nindt$. Then \eqref{ind.mat} holds for $k'+1$, using \eqref{eq:nasty:D:wq:old}, \eqref{eq:nasty:D:vq:old}, the inductive assumption \eqref{ind.mat} for $k'$, and \eqref{eq:timescale:upgrading}. 
\end{proof}

\subsection{Mollification estimates}\label{appetizer:5}

In this subsection, we require two algebraic identities originally stated in \cite[(5.17a)--(5.17b)]{BMNV21}, which we now recall. Let $v$ be a sufficiently smooth divergence-free vector field and let $D_t = \pa_t + v\cdot \na$ be the material derivative operator associated to $v$. For any sufficiently smooth function $F = F(x,t)$ and any $n,m\geq 0$, the Leibniz rule implies that
\begin{subequations}
\begin{align}
D^n D_t^m F
&=D^n (\partial_t + v \cdot \nabla_x)^m F
= \sum_{\substack{m'\leq m \\ n'+m' \leq n+m}} d_{n,m,n',m'}(v)(x,t) D^{n'} \partial_t^{m'} F \, ,
\label{eq:RNC} \\
d_{n,m,n',m'}(v)
&= \sum_{k = 0}^{m-m'} \sum_{\substack{\{ \gamma \in {\mathbb N}^k \colon |\gamma| = n-n'+k,\\ \beta \in {\mathbb N}^k \colon |\beta| = m-m'-k\}}} c(m,n, k, \gamma, \beta)  \prod_{\ell=1}^{k} \left( D^{\gamma_\ell} \partial_t^{\beta_\ell}  v(x,t) \right) \, ,
\label{eq:DNC}
\end{align}
\end{subequations}
where $c(m,n, k, \gamma, \beta)$ denotes an explicitly computable combinatorial coefficient  which depends only on the factors inside the parentheses. Identities \eqref{eq:RNC}--\eqref{eq:DNC} hold because $D$ and $\partial_t$ commute; the proof is based on induction on $n$ and $m$ and is left to the reader. 
\begin{proposition}[\bf Mollification with spatial and material derivatives]\label{lem:mollification:general}
Let $p\in[1,\infty]$, $N_{\rm g}$, $N_{\rm c}$, $M_t$, $N_\ast$, and $N_\gamma$ be positive integers, $v$ be a divergence-free vector field, and $D_t = \partial_t + v\cdot\nabla$.  Fix parameters $\lambda$, $\Lambda$, $\tau$, $\Tau$, $\Gamma\geq 1$, $i$, $\const_{f,p} \leq \tilde\const_{f}$, $\const_v$, and $c\in [0, 30]$ such that
\begin{subequations}\label{eq:moll:assumps:1}
\begin{align}
N_{\rm g} \leq N_{\rm c} \leq N_* / 4 \, , \quad  M_t \leq N_* \leq N_\gamma\, , \quad 
\lambda\Gamma &\leq \Lambda \, , \quad \tau^{-1}\Gamma^{i+ c}\leq\,  \Tau^{-1} \, , \quad \const_v \la \leq \Tau^{-1}\, ,
\label{eq:moll:assumps:1:1}\\
(\Tau^{-1}{\Ga})^{M_t} \tilde\const_{f} \Ga^{-\sfrac{N_\textnormal{c}}{2}} &\leq  {{\Ga^{-N_{ \rm g}}}} \const_{f,p} \tau^{-M_t} \label{eq:moll:assumps:1:3}
\, .  
\end{align}
\end{subequations}
Let $(a,b)+\Tau$ be a time domain and $\Omega\subset (a,b)+\Tau \times\T^d$ be a subset in the space-time domain. Assume that $v$ satisfies
\begin{equation}\label{moll.assum.v.est}
    \norm{ D^{N} \partial_t^{M}  v(x,t) }_{L^\infty((a,b)+\Tau \times \T^3)}  \lec \const_v \la^N \Tau^{-M}
\end{equation}
for all $N+M \leq N_\gamma$.
Assume that $f:(a,b)+\Tau \times\T^d\rightarrow \R$ satisfies the estimates\footnote{By $L^p(\Omega)$, we mean $L^p$ for each fixed timeslice $\Omega\cap \{t=t_0\}$, continuously in time which is non-empty.}
\begin{subequations}
\begin{align}
    \left\| D^N D_t^M f \right\|_{L^p(\Omega)} &\les \const_{f,p} \lambda^N \MM{M,M_t,\tau^{-1}\Gamma^{i+ c},\Tau^{-1}} \label{eq:moll:f:1}\\
    \left\| D^N \partial_t^M f \right\|_{L^\infty((a,b)+\Tau \times\T^d)} &\les \tilde\const_{f} \lambda^N \Tau^{-M} \label{eq:moll:f:2}
\end{align}
\end{subequations}
for $N+M\leq N_*$. Let $\gamma_x$ be a compactly supported mollifier in space at scale $(\lambda^{-1}\Lambda^{-1})^{\sfrac 12}$, $\gamma_t$ be a compactly supported mollifier in time at scale $\Tau\Ga^{-\sfrac 12}$, and assume that the kernels for both mollifiers have vanishing moments up to $N_{\rm c}$ and are $C^{N_\gamma}$ differentiable.

Set $f_\gamma = \gamma_t \ast \gamma_x \ast f $.  Then for $N+M\leq N_\gamma$, we have that
\begin{align}\label{eq:moll:conc:1}
    \left\| D^N D_t^M f_\gamma \right\|_{L^p(\Omega\cap (a,b)\times\T^d)} \lesssim \const_{f,p} \Lambda^N \MM{M,M_t,\tau^{-1}\Gamma^{i+{c}+1},\Tau^{-1}\Gamma} \, ,
\end{align}
while for $N+M\leq N_\ast$, we have that
\begin{align}\label{eq:moll:conc:2}
    \left\| D^N D_t^M (f-f_\gamma) \right\|_{L^p(\Omega\cap (a,b)\times\T^d)} \les \Gamma^{-N_{\rm g}} \const_{f,p} \Lambda^N \MM{M,M_t,\tau^{-1},\Tau^{-1}\Gamma} \, .
\end{align}
\end{proposition}
\begin{proof}
We split the proof into steps. We first set up the Taylor expansion which allows us to take advantage of the vanishing moments.  Next, we prove \eqref{eq:moll:conc:1} and \eqref{eq:moll:conc:2} for $N,M\leq \sfrac{N_*}{4}$.  Finally, we prove \eqref{eq:moll:conc:1} and \eqref{eq:moll:conc:2} in the remaining cases where either $N>\sfrac{N_*}{4}$ or $M>\sfrac{N_*}{4}$. Note that since $\gamma_t$ has a compact support in time at scale $\Tau \Ga^{-\sfrac12}$, $f_\gamma$ is well-defined in the domain $(a,b)\times \T^d$. 
\smallskip

\noindent\texttt{Step 1:} Let us denote by $K_t$ the kernel for $\gamma_t$ and $K_x$ the kernel for $\gamma_x$ so that $K:=K_t K_x$ is the space-time kernel for $\gamma_t \ast \gamma_x$. We denote space-time points $(t,x) \in (a,b)\times\T^d$ and $(s,y)\in(a,b)+\Tau \times\T^d$ by
\begin{align}
\label{eq:strange:notation} 
(t,x) = \theta, \qquad (s,y) = \kappa \, .
\end{align}
Using this notation we may write out $f_\gamma$ explicitly as
\begin{align}\label{eq:commutator:stress:1}
f_\gamma (\theta) = \int_{\mathbb{T}^d\times\mathbb{R}}  f (\theta-\kappa) K (\kappa) \,d\kappa \,.
\end{align}
Expanding $f$ in a Taylor series in space and time around $\theta$ yields the formula
\begin{align}
    f(\theta-\kappa) = f(\theta) + \sum_{|\alpha|+m=1}^{{N_{\textnormal{c}}-1}} \frac{1}{\alpha! m!}D^\alpha \partial_t^m f(\theta) (-\kappa)^{(\alpha,m)} + R_{N_\textnormal{c}}(\theta,\kappa)
\end{align}
where
\begin{align}
    R_{N_{\textnormal{c}}}(\theta,\kappa) = \sum_{|\alpha| + m =N_\textnormal{c}} \frac{N_{\textnormal{c}}}{\alpha! m!} (-\kappa)^{(\alpha,m)} \int_0^1 (1-\eta)^{N_\textnormal{c}-1} D^\alpha \partial_t^m f (\theta-\eta\kappa) \,d\eta \, .
    \label{eq:R:q:comm:remainder}
\end{align}
\smallskip

\noindent\texttt{Step 2:} Assume that $N,M\leq \sfrac{N_\ast}{4}$. Here we note that because of the vanishing moments of $K$,
\begin{align}
f_\gamma(\theta) - f(\theta)
&= \sum_{|\alpha| + m'' = N_{\textnormal{c}}} \frac{N_{\textnormal{c}}}{\alpha! m''!}  \int_{\mathbb{T}^d\times\mathbb{R}} K(\kappa) (-\kappa)^{(\alpha,m'')} \int_0^1 (1-\eta)^{N_{\textnormal{c}}-1} D^\alpha \partial_t^{m''} f (\theta-\eta\kappa) \,d\eta \, d\kappa \, .
\label{eq:mollify:minus:identity}
\end{align}
Now we appeal to the identity \eqref{eq:RNC} with $F = f_\gamma - f$ to obtain
\begin{align}
\norm{D^n D_t^m (f_\gamma - f)}_{L^\infty((a,b) \times\T^d)}  
\les \sum_{\substack{m'\leq m \\ n'+m' \leq n+m}} \norm{d_{n,m,n',m'}(v)}_{L^\infty} \norm{ D^{n'} \partial_t^{m'} (f_\gamma - f)}_{L^\infty((a,b) \times\T^d)} \, . 
\label{eq:MNF:1}
\end{align}
From assumptions \eqref{eq:moll:assumps:1} and \eqref{moll.assum.v.est} and the formula \eqref{eq:DNC}, we have that 
\begin{align}
\norm{d_{n,m,n',m'}(v)}_{L^\infty} \les \sum_{k=0}^{m-m'} \const_v^k \lambda^{ n-n' +k} (\Tau^{-1})^{m-m' -k} \lec \lambda^{ n-n'} (\Tau^{-1})^{m-m'}
\,.
\label{eq:MNF:2}
\end{align} 
Combining this estimate with the bound  \eqref{eq:moll:f:2}, we deduce that
\begin{align}
&\left\|D^N D_{t}^M (f_\gamma - f)\right\|_{L^\infty((a,b) \times\T^d)} \notag\\
&\qquad\les \sum_{\substack{m'\leq M \\ n'+m' \leq N+M}} \lambda^{ N-n'} (\Tau^{-1} )^{M-m'} \norm{ D^{n'} \partial_t^{m'} (f_\gamma - f)}_{L^\infty((a,b) \times\T^d)} 
\notag\\
&\qquad\les \sum_{\substack{m'\leq M \\ n'+m' \leq N+M}} \sum_{|\alpha| + m'' = N_{\textnormal{c}}} \lambda^{ N-n'} (\Tau^{-1})^{M-m'} \times \tilde\const_{f} \lambda^{n' + |\alpha| }  (\Tau^{-1})^{m'+ m''}
\int_{\T^3\times \R} \abs{\kappa^{(\alpha,m'')}} |K(\kappa)| d\kappa
\notag\\
&\qquad\les \tilde\const_{f} \sum_{|\alpha| + m'' = N_{\textnormal{c}}} \lambda^{ N + |\alpha|} (\Tau^{-1})^{M+m''} 
(\Lambda\lambda)^{-\sfrac{|\alpha|}{2}} (\Tau \Ga^{-\sfrac 12})^{m''}
\notag\\
&\qquad\les \tilde\const_{f} \lambda^{ N } \Tau^{-M}  
\Gamma^{-\sfrac{N_{\textnormal{c}}}{2}}\les {\Gamma^{-{N_{\rm g}}}} \const_{f,p} \Lambda^N \MM{M,M_t,\tau^{-1},\Tau^{-1}\Gamma}\,,
\label{eq:dex:1}
\end{align}
where the last inequality follows from \eqref{eq:moll:assumps:1} and holds for $N,M\leq \sfrac{N_*}{4}$. This establishes \eqref{eq:moll:conc:2} in this range of $N,M$, and by the triangle inequality for $f_\gamma = f_\gamma - f +f$ establishes \eqref{eq:moll:conc:1} in the same range of $N,M$.
\smallskip

\noindent\texttt{Step 3:} We now consider \eqref{eq:moll:conc:1} in the case that either $M\geq \sfrac{N_*}{4}$ or $N\geq \sfrac{N_*}{4}$, and $N+M\leq N_\gamma$. We first note that when $N_*\leq N+M\leq  N_\gamma$, applying the differential operator to the kernels for the mollifiers, we get
\begin{align}
    \norm{D^N \pa_t^M f_\gamma}_{L^\infty((a,b) \times\T^d)}
    \lec \td \const_f \min_{\substack{n+m=N_*\\ n\leq N, \, m\leq M}} \la^n \Tau^{-m} (\la\Lambda)^{\frac12(N-n)} (\Tau^{-1}\Gamma^{\sfrac12})^{M-m}
    \label{eq:moll:f:3}
\end{align}
This implies that when either $N$ or $M$ exceeds $\sfrac{N_*}{4}$ but $N+M\leq N_\gamma$, we have
\begin{align}
    \left\| D^N D_t^M f_\gamma \right\|_{L^\infty((a,b) \times\T^d)} &\lesssim \sum_{\substack{m\leq M \\ n+m\leq N+M}} \left\| d_{N,M,n,m}(v) \right\|_{L^\infty} \left\| D^n \partial_t^m f_\gamma \right\|_{L^\infty} \notag\\
   & \lesssim 
   \tilde\const_{f} \Ga^{-\frac{N_*}{8}} \Lambda^N (\Tau^{-1}\Gamma)^M
   \lesssim 
   \tilde\const_{f} \Ga^{-\frac{N_c}{2}} \Lambda^N (\Tau^{-1}\Gamma)^M
   \\
    &\lesssim \Gamma^{-N_{\rm g}} \const_{f,p} \Lambda^N \MM{M,M_t,\tau^{-1},\Tau^{-1}\Gamma} \notag 
\end{align}
where we have used \eqref{eq:MNF:2},  \eqref{eq:moll:f:2}, \eqref{eq:moll:f:3},
\eqref{eq:moll:assumps:1}, and \eqref{eq:moll:assumps:1:3}
. In the second inequality, the factor $\Ga^{-\frac{N_*}{8}}$ gain has been obtained by paying lossy derivative costs. This completes the proof of \eqref{eq:moll:conc:1} when either $N$ or $M$ exceeds $\sfrac{N_*}{4}$ and $N+M\leq N_\gamma$.

Finally, in order to prove \eqref{eq:moll:conc:2} when either $N$ or $M$ exceeds $\sfrac{N_*}{4}$ and $N+M\leq N_*$, we use the triangle inequality as in the previous step, the estimate just shown, and the estimate
\begin{align*}
    \left\| D^N D_t^M f \right\|_{L^p(\Omega\cap (a,b) \times\T^d)} &\lesssim \const_{f,p} \Gamma^{-(M+N)} \Lambda^N \MM{M,M_t,\tau^{-1}\Gamma^{i+c+1},\Tau^{-1}\Gamma} \\
    &\lesssim \Gamma^{-N_{\rm g}} \const_{f,p} \Lambda^N \MM{M,M_t,\tau^{-1},\Tau^{-1}\Gamma} \, ,
\end{align*}
which follows from \eqref{eq:moll:f:1} and \eqref{eq:moll:assumps:1}.
\end{proof}

\printindex

\medskip

\noindent\textsc{Department of Mathematics, Princeton University, Princeton, NJ, USA.}
\vspace{.03in}
\newline\noindent\textit{Email address}: \href{mailto:vgiri@math.princeton.edu}{vgiri@math.princeton.edu}.
\smallskip

\noindent\textsc{Department of Mathematics, ETH Z\"urich, Z\"urich, Switzerland.}
\vspace{.03in}
\newline\noindent\textit{Email address}: \href{hyunju.kwon@math.ethz.ch}{hyunju.kwon@math.ethz.ch}.
\smallskip

\noindent\textsc{Department of Mathematics, Purdue University, West Lafayette, IN, USA.}
\vspace{.03in}
\newline\noindent\textit{Email address}: \href{mailto:mdnovack@purdue.edu}{mdnovack@purdue.edu}.

\end{document}